\documentclass[10pt]{book}
\usepackage[french]{babel}
\usepackage{amssymb, amsmath, amscd}

\numberwithin{equation}{section}

\usepackage{makeidx}
\newtheorem{thm}{Th\'eor\`eme}[section]
\newtheorem{lem}[thm]{Lemme}
\newtheorem{prop}[thm]{Proposition}
\newtheorem{conj}[thm]{Conjecture}
\newtheorem{cor}[thm]{Corollaire}
\newtheorem{qu}[thm]{Question}
\newtheorem{defi}[thm]{D\'efinition}
\newtheorem{prob}[thm]{Probl\`eme}
\newcommand{\LG}{\mathfrak{g}}
\newcommand{\LK}{\mathfrak{k}}
\newcommand{\LM}{\mathfrak{p}}
\newcommand{\la}{\mathfrak{a}}
\newcommand{\lln}{\mathfrak{n}}

\def\pa{\S\kern.15em}
\def\adots{\mathinner{\mkern2mu\raise1pt\hbox{.}
\mkern3mu\raise4pt\hbox{.}\mkern1mu\raise7pt\hbox{.}}}

\def\ba{\backslash}
\def\1{\hbox{\small 1}\!\!\hbox{\normalsize 1}}
\def\SL{\textrm{SL}}
\def\Gal{\textrm{Gal}}
\def\ind{\textrm{ind}}
\def\GL{\textrm{GL}}

\def\({{\rm (}}
\def\){{\rm )}}
\def\[{{\rm [}}
\def\]{{\rm ]}}
\def\ba{\backslash}
\def\1{\hbox{\small 1}\!\!\hbox{\normalsize 1}}
\def\SL{\textrm{SL}}
\def\Gal{\textrm{Gal}}
\def\ind{\textrm{ind}}
\def\GL{\textrm{GL}}

\title{Spectre automorphe des vari\'et\'es hyperboliques et applications topologiques}

\author{N. Bergeron et L. Clozel}
 
\date{} 
 
\makeindex

\begin{document}

\maketitle

\pagestyle{empty}

\begin{quotation}
\`A Natan, C\'eline et \`a toute l'\'equipe de la cr\`eche parentale ``Le Moulin bleu'' sans qui ce 
livre serait d\'ej\`a fini depuis quelques mois.

N.B.
\end{quotation}

\newpage

\pagestyle{headings}

\pagenumbering{roman}

\chapter*{Introduction}

\section*{Un programme conjectural}

Soit $G$ un groupe alg\'ebrique connexe et presque simple sur ${\Bbb Q}$ tel que le groupe
$G ({\Bbb R} )$ des points r\'eels soit le produit (avec intersection finie) d'un groupe compact et d'un groupe r\'eel 
non compact presque simple. Nous noterons ce dernier groupe $G^{{\rm nc}}$ (nc signifie ici non compact).
Dans cette introduction (comme dans la majeure partie du texte) nous supposons $G^{{\rm nc}}$ isomorphe au groupe 
$SO(n,1)$ (resp. $SU(n,1)$).

L'espace sym\'etrique associ\'e au groupe $G$ (plus pr\'ecisemment au groupe $G ({\Bbb R} )$) est alors l'espace hyperbolique r\'eel (resp. complexe) de 
dimension r\'eelle (resp. complexe) $n$, nous le noterons $X_G$, $d_G$ sa dimension r\'eelle et nous supposerons la m\'etrique 
normalis\'ee de fa\c{c}on \`a ce que les courbures sectionnelles r\'eelles soient comprises entre $-4$ et $-1$ ($-1$ pour les espaces hyperboliques r\'eels associ\'es 
aux groupes $SO(n,1)$ et $-4$ pour le plan hyperbolique r\'eel associ\'e au groupe $SU(1,1)$).  Nous 
noterons $d=1$ dans le cas r\'eel et $d=2$ dans le cas complexe; on a donc $d_G = dn$. Posons enfin 
$\rho_G = \frac{d_G -2 +d}{2}$. 

Soit $K_f$ un sous-groupe compact-ouvert de $G({\Bbb A}_f )$ (o\`u ${\Bbb A}_f$ d\'esigne l'anneau des 
ad\`eles finis\index{ad\`eles, finis} sur ${\Bbb Q}$) tel que le groupe $\Gamma = G ({\Bbb Q} ) \cap K_f$ soit sans torsion.
Un tel groupe $\Gamma$ est appel\'e {\it sous-groupe de congruence (sans torsion) de 
$G$}\,~\footnote{Lorsque $G$ est d\'efini sur ${\Bbb Z}$, on peut pr\'ef\'erer consid\'erer les sous-groupes 
$\Gamma$ de $G(\mathbb{Z} )$ sans torsion et contenant un sous-groupe de la forme ker$(G(\mathbb{Z} ) \rightarrow G(\mathbb{Z} / N\mathbb{Z} ))$ pour 
un certain entier $N\geq 1$.} \index{congruence, sous-groupe de}.
On appellera {\it vari\'et\'e hyperbolique r\'eelle (resp. complexe) de congruence} \index{congruence, vari\'et\'e hyperbolique de} tout quotient de 
l'espace hyperbolique r\'eel (resp. complexe) $X_G$ 
par un sous-groupe de congruence (sans torsion) $\Gamma$ de $G$ comme ci-dessus. 

D'apr\`es un th\'eor\`eme de Borel et Harish-Chandra \cite{BorelHarishChandra}, une vari\'et\'e hyperbolique 
r\'eelle (resp. complexe) de congruence est toujours de volume fini et est compacte si et seulement si le groupe $G$  
est anisotrope sur ${\Bbb Q}$.

En 1965, Selberg \cite{Selberg} a d\'emontr\'e que si $\Gamma$ est un sous-groupe de congruence du ${\Bbb Q}$-groupe $SL(2)$, 
alors la premi\`ere valeur propre non nulle $\lambda_1 $ du laplacien sur les fonctions 
$L^2$ de $\Gamma \backslash {\Bbb H}^2$ est sup\'erieure ou \'egale \`a $\frac{3}{16}$. Il conjecture de plus que 
$\lambda_1 \geq \frac{1}{4}$, la valeur apparaissant d'ailleurs dans le spectre $L^2$ d'un quotient $\Gamma \backslash {\Bbb H}^2$
pour un certain groupe de congruence $\Gamma$. 

Apr\`es des travaux de Selberg \cite{Selberg}, Jacquet et Langlands \cite{JacquetLanglands}, Gelbart et Jacquet \cite{GelbartJacquet},
Sarnak \cite{Sarnak}, Elstrodt, Grunewald et Mennicke \cite{EGM}, Li, Piatetski-Shapiro et Sarnak \cite{LPSS} et Burger et 
Sarnak \cite{BurgerSarnak}, le second auteur a r\'ecemment \'etendu, dans \cite{Clozel}, le th\'eor\`eme de Selberg \`a toutes les vari\'et\'es hyperboliques
de congruence~: si $G$ est un groupe ${\Bbb Q}$-alg\'ebrique comme ci-dessus et $\Gamma$ un sous-groupe de 
congruence de $G$, alors la premi\`ere valeur propre non nulle du laplacien sur les fonctions $L^2$ de 
$\Gamma \backslash X_G$ est minor\'ee par une constante strictement positive ne d\'ependant que de la dimension 
r\'eelle $d_G$ de $X_G$ (``Conjecture $\tau$''). 

Dans la premi\`ere partie de cet article on commence l'\'etude du spectre sur les formes diff\'erentielles. Le fait que 
l'on se restreigne aux espaces hyperboliques r\'eels et complexes sera motiv\'e au chapitre 3 o\`u l'on 
mentionera des r\'esultats concernant d'autres espaces sym\'etriques. 
Il est d'ailleurs int\'eressant de noter que le rang sup\'erieur, $\geq 2$, joue encore un r\^ole sp\'ecial dans l'\'etude du spectre
sur les formes diff\'erentielles (et non seulement sur celui sur les fonctions); on explicitera ce ph\'enom\`ene sur les
groupes unitaires $U(p,q)$. Revenons maintenant \`a notre groupe $G$. Dans notre \'etude nous sommes guid\'es par 
la conjecture suivante, qui g\'en\'eralise la Conjecture de Selberg.

\medskip

\noindent
{\bf Conjecture A} {\it Soit $G$ un groupe ${\Bbb Q}$-alg\'ebrique comme ci-dessus. 
\begin{enumerate}
\item Si $G^{\rm nc}  \cong SO(2n,1)$, alors pour tout entier naturel $i\leq n$, 
$$\lambda_1^i  ( \Gamma \backslash X_G ) \geq \max (2n-2i-2 , \frac{1}{4} ) >0 .$$
\item Si $G^{\rm nc} \cong SO(2n+1,1)$, alors pour tout entier naturel $i\leq n-1$, 
$$\lambda_1^i  ( \Gamma \backslash X_G ) \geq 2n -2i -1 >0 .$$
\item Si $G^{\rm nc} \cong SU(n,1) $, alors pour tout entier naturel $i\leq n$, 
$$\lambda_1^i  ( \Gamma \backslash X_G ) \geq \max (4(n-i-1) ,1) >0 .$$
\end{enumerate}
Ici $\lambda_1^i$ d\'esigne la premi\`ere valeur propre non nulle du laplacien de Hodge sur les formes diff\'erentielles
de degr\'e $i$ et $\Gamma$ est un sous-groupe de congruence quelconque de $G$.}

\medskip \index{Conjecture A}

La Conjecture A est tr\`es profonde (elle contient en particulier la Conjecture de Selberg) et de nature arithm\'etique, la
pr\'ediction que les constantes explicites de minoration du spectre sont toujours des demi-entiers est une ``conjecture 
de puret\'e'' qui renvoie au th\'eor\`eme de puret\'e de Deligne. On peut relaxer un peu cette conjecture en un \'enonc\'e
plus g\'eom\'etrique.

\medskip

\noindent
{\bf Conjecture A$^-$} {\it Soit $G$ comme dans la Conjecture A. Alors, pour tout entier naturel $i\leq 
\frac{d_G}{2} -1$ il existe une constante strictement positive $\varepsilon (G , i)$ telle que pour tout 
sous-groupe de congruence $\Gamma$ dans $G$, 
$$\lambda_1^i (\Gamma \backslash X_G ) \geq \varepsilon (G ,i) .$$}

\medskip \index{Conjecture A$^-$}

Soit $i$ un entier naturel et $G$ un groupe ${\Bbb Q}$-alg\'ebrique comme dans la Conjecture A. On dira 
que la {\it Conjecture A$^-$($i$) est v\'erifi\'ee} \index{Conjecture A$^-$($i$)} par $G$ s'il existe une constante strictement positive 
$\varepsilon (G ,i)$ telle que pour tout sous-groupe de congruence $\Gamma$ de $ G$,
$$\lambda_1^i (\Gamma \backslash X_G ) \geq \varepsilon (G ,i),$$
o\`u $\lambda_1^i$ d\'esigne encore la premi\`ere valeur propre non nulle du laplacien de Hodge sur les 
formes diff\'erentielles de degr\'e $i$.

\medskip

La solution mentionn\'ee ci-dessus de la Conjecture $\tau$ peut donc s'\'enoncer de la fa\c{c}on suivante~:
{\it la Conjecture A$^- (0)$ est v\'erifi\'ee par tout groupe $G$ comme dans la Conjecture A.}

\medskip

La Conjecture A$^-$ pour les groupes $G$ tels que $G^{\rm nc} \cong SU(n,1)$ implique qu'il existe une constante 
strictement positive $\varepsilon (G)$ telle que la premi\`ere valeur propre du laplacien de Hodge sur l'ensemble des 
formes diff\'erentielles (de tous les degr\'es) sur $\Gamma \backslash X_G$, pour n'importe quel sous-groupe de 
congruence $\Gamma$ de $G$, soit sup\'erieure ou \'egale \`a $\varepsilon (G)$.

Cette propri\'et\'e est fausse pour les groupes $G$ tels que  $G^{{\rm nc }} \cong SO(n,1)$ lorsque $n$ est impair. En effet, 
lorsque $G^{{\rm nc }} \cong  SO(2n+1,1)$, $\lambda_1^n (\Gamma \backslash X_G )$, la premi\`ere valeur propre non nulle du laplacien de 
Hodge sur les formes diff\'erentielles de degr\'e $n$, peut \^etre rendue arbitrairement petite quitte \`a changer de 
sous-groupe de congruence $\Gamma$ de $G$. N\'eanmoins si l'on se restreint au spectre des $n$-formes 
diff\'erentielles {\bf ferm\'ees} la Conjecture A$^-$ implique encore une minoration uniforme du spectre sur les $n$-formes
diff\'erentielles. On verra que cette propri\'et\'e est utile g\'eom\'etriquement. \'Etant donn\'e un entier naturel $i$ et 
un groupe $ G$ comme dans la Conjecture A, on dira que la 
{\it Conjecture A$^-_{d=0}$($i$) est v\'erifi\'ee} \index{Conjecture A$^-_{d=0}$($i$)} par $G$ s'il existe une constante strictement positive 
$\varepsilon (G ,i)$ telle que la premi\`ere valeur propre du laplacien de Hodge sur les
formes diff\'erentielles ferm\'ees de degr\'e $i$ soit $\geq \varepsilon (G,i)$.
 
\medskip

Toute la premi\`ere partie du texte sera motiv\'ee par la double Conjecture A, A$^-$. La deuxi\`eme partie
est quant \`a elle motiv\'ee par une seconde conjecture de nature topologique$/$g\'eom\'etrique. Cette conjecture
concerne l'homologie des vari\'et\'es hyperboliques r\'eelles ou complexes de congruence. Pr\'ecisons que tous
les groupes de (co-)homologie que nous consid\`ererons dans cet article seront \`a coefficients rationnels.

\medskip

\noindent
{\bf Conjecture B} {\it Soit $M$ une vari\'et\'e hyperbolique r\'eelle (resp. complexe) de congruence de dimension $m$.
Soit $F$ une sous-vari\'et\'e r\'eelle (resp. complexe) compacte connexe totalement g\'eod\'esique immerg\'ee dans 
$M$. Alors, il existe un rev\^etement fini de congruence $\hat{M}$ de $M$ tel que~:
\begin{enumerate}
\item l'immersion de $F$ dans $M$ se rel\`eve en un plongement de $F$ dans $\hat{M}$,
\item l'application induite~: 
$$H_p (F ) \rightarrow H_p (\hat{M} ) $$
est injective pour tout entier $p\geq \frac{m}{2}$.
\end{enumerate}
}

\medskip \index{Conjecture B}

Nous dirons qu'une vari\'et\'e hyperbolique r\'eelle (resp. complexe) de congruence de dimension $m$ v\'erifie la Conjecture 
B($i$) \index{Conjecture B($i$)} pour un certain entier $i\geq \frac{m}{2}$, si pour toute sous-vari\'et\'e r\'eelle (resp. complexe) compacte connexe
totalement g\'eod\'esique immerg\'ee dans $M$, il existe un rev\^etement fini de congruence $\hat{M}$ de 
$M$ tel que~:
\begin{enumerate}
\item l'immersion de $F$ dans $M$ se rel\`eve en un plongement de $F$ dans $\hat{M}$,
\item l'application induite~: 
$$H_i (F ) \rightarrow H_i (\hat{M} ) $$
est injective.
\end{enumerate}
  
\medskip

La Conjecture B est \`a comparer avec les propri\'et\'es de Lefschetz pour les vari\'et\'es hyperboliques 
complexes de congruence mises en \'evidence par Oda \cite{Oda}, Harris et Li \cite{HarrisLi} et 
Venkataramana \cite{Venky}. Ces propri\'et\'es s'\'enoncent plus facilement dans le langage de Harris et Li et 
Venkataramana que l'on rappelle ci-dessous.

\medskip

Soient $H < G$ deux groupes alg\'ebriques sur ${\Bbb Q}$ comme plus haut. Supposons de plus $H$ et $G$ 
{\bf anisotropes} sur ${\Bbb Q}$ et du m\^eme type, orthogonal ou unitaire.
L'inclusion $H \subset  G$ induit un plongement totalement g\'eod\'esique naturel 
$X_H \rightarrow X_G$ entre les espaces sym\'etriques associ\'es. 

Soit $K_f$ un sous-groupe compact-ouvert de $G ({\Bbb A}_f )$ tel que le groupe $\Gamma = G
({\Bbb Q} ) \cap K_f$ soit sans torsion. D\'esignons par $\Gamma_H$ l'intersection de $\Gamma$ avec $H ({\Bbb Q})$.
On obtient alors une application lisse~: 
$$j= j(\Gamma) : \Gamma_H \backslash X_H \rightarrow \Gamma \backslash X_G$$
\index{$j$} pour chaque sous-groupe de congruence sans torsion de $G ({\Bbb Q} )$. 
Celle-ci induit une application naturelle 
\begin{eqnarray} \label{appl en homologie}
j_* : H_* (\Gamma_H \backslash X_H ) \rightarrow H_* (\Gamma \backslash X_G ).
\end{eqnarray}
Remarquons que si $g$ est un \'el\'ement quelconque de $G ({\Bbb Q})$, on peut translater l'immersion $j$ en 
une application lisse
$$j_g = j_g (\Gamma ) : (H \cap g^{-1} \Gamma g) \backslash X_H \rightarrow \Gamma \backslash X_G .$$ \index{$j_g$}
On appelle {\it application de restriction virtuelle} \index{restriction virtuelle} \index{virtuelle, restriction}, l'application 
\begin{eqnarray} \label{appl en cohomologie}
H^* (\Gamma \backslash X_G ) \rightarrow \prod_{g\in G ({\Bbb Q})} H^* ((H \cap g^{-1} \Gamma g)
\backslash X_H)
\end{eqnarray}
induite en cohomologie par la famille d'applications $(j_g )_{g\in G ({\Bbb Q})}$. 

On s'int\'eresse aux applications (\ref{appl en homologie}) et (\ref{appl en cohomologie}) {\it virtuellement} \index{virtuellement} en  un deuxi\`eme sens~: 
\`a rev\^etements fini (de congruence) pr\`es. Soient donc $\Gamma ' \subset \Gamma$ deux sous-groupes de 
congruence sans torsion de $G ({\Bbb Q})$. Alors le rev\^etement fini (de vari\'et\'es compactes)
$$\Gamma ' \backslash X_G \rightarrow \Gamma \backslash X_G$$
induit, en homologie, une application surjective
\begin{eqnarray} \label{surj en homologie}
H_* (\Gamma ' \backslash X_G ) \rightarrow H_* (\Gamma \backslash X_G ) 
\end{eqnarray}
et, en cohomologie, une application injective
\begin{eqnarray} \label{inj en cohom}
H^* (\Gamma \backslash X_G ) \rightarrow H^* (\Gamma ' \backslash X_G ).
\end{eqnarray}
On obtient ainisi un syst\`eme projectif (resp. inductif) de groupes d'homologie (resp. de cohomologie) index\'e par 
les sous-groupes de congruence de $G ({\Bbb Q})$. On d\'esigne par $H_* (Sh^0 G)$ \index{$H_* (Sh^0 G)$} la limite projective du 
syst\`eme (\ref{surj en homologie}) et par $H^* (Sh^0 G)$ \index{$H^* (Sh^0 G)$} la limite inductive du syst\`eme (\ref{inj en cohom}). En passant
\`a la limite, lorsque $\Gamma$ varie, dans (\ref{appl en homologie}) et (\ref{appl en cohomologie}), on obtient les 
applications naturelles~:
$$H_* (Sh^0 H) \rightarrow H_* (Sh^0 G)$$
et 
$$H^* (Sh^0 G ) \rightarrow \prod_{g \in G ( {\Bbb Q})} H^* (Sh^0 H) .$$

\medskip

Rappelons que les th\'eor\`emes de Lefschetz permettent d'\'etudier la topologie d'une vari\'et\'e projective complexe
par r\'ecurrence, en la coupant par un hyperplan de l'espace projectif. Un yoga semblable, mais moins fort, existe pour
les vari\'et\'es hyperboliques complexes de congruence. Il a \'et\'e d\'egag\'e par Harris et Li \cite{HarrisLi} et 
prouv\'e par Venkataramana \cite{Venky}~: {\it si $H \subset G$ sont deux groupes comme ci-dessus de type 
unitaire alors l'application naturelle
$$H^i (Sh^0 G ) \rightarrow \prod_{g \in G ( {\Bbb Q})} H^i (Sh^0 H) $$
est injective pour tout entier naturel $i\leq \frac{d_H}{2}$}.

La Conjecture B ci-dessus implique une conjecture plus faible renfor\c{c}ant le yoga de type Lefschetz~:

\medskip

\noindent
{\bf Conjecture B$^-$} {\it Soient $H \subset G$ comme plus haut. Alors pour tout entier $i\geq \frac{d_G}{2}$,
l'application naturelle
$$H_i (Sh^0 H)  \rightarrow H_i (Sh^0 G)$$
est injective.}
\index{Conjecture B$^-$}

\medskip

Remarquons que par dualit\'e, la Conjecture B$^-$ implique que l'application naturelle de restriction 
$$H^i (Sh^0 G) \rightarrow H^i (Sh^0 H)$$
est surjective. 

\medskip

Dans la deuxi\`eme partie de cet article, nous relions la Conjecture B \`a la Conjecture A$^-$ \'enonc\'ee plus haut. Les 
r\'esultats de la premi\`ere partie concernant la Conjecture A vont nous permettre de d\'emontrer des r\'esultats partiels
concernant la Conjecture B. Avant de passer \`a la description de nos r\'esultats, remarquons que  
d'autres propri\'et\'es de type Lefschetz sont conjectur\'ees et \'etudi\'ees par le premier auteur dans 
\cite{IRMN}, ces propri\'et\'es sont elles aussi li\'ees \`a la Conjecture A$^-$.

Enfin soulignons que la plupart des m\'ethodes d\'evelopp\'ees dans cet article devraient se g\'en\'eraliser 
aux groupes alg\'ebriques semi-simples g\'en\'eraux. En ce qui concerne les r\'esultats spectraux, on l'a dit, le 
cas que nous consid\'erons est le plus int\'eressant. N\'eanmoins ceci n'est plus vrai en ce qui concerne les
r\'esultats topologiques. \`A la fin du texte nous d\'ecrivons d'ailleurs un th\'eor\`eme de la veine de la Conjecture B
mais pour les vari\'et\'es de congruence associ\'ees aux groupes unitaires $U(p,q)$. C'est un premier r\'esultat, 
le premier auteur esp\`ere revenir \`a l'\'etude de ces vari\'et\'es dans un autre texte, en explicitant la 
combinatoire (plus riche) des propri\'et\'es du type Lefschetz auxquelles on doit s'attendre.

\section*{Description des r\'esultats}

On l'a rappel\'e, le second auteur a r\'ecemment d\'emontr\'e la Conjecture A$^- (0)$ pour tout groupe $G$ 
comme dans la Conjecture A. M\^eme la seule Conjecture A$^- (1)$ semble difficile. Il est naturel de commencer 
\`a \'etudier celle-ci sur de petits groupes. Le groupe $SL(2)$ est trop petit~: la Conjecture A$^- (1)$ d\'ecoule trivialement de la Conjecture A$^- (0)$ (la Conjecture
A$^- (i)$ n'est d'ailleurs formul\'e ci-dessus que pour $i \leq \frac{d_G}{2} -1 $ \'egal \`a $0$ dans notre cas).  
Le plus petit groupe pour lequel la Conjecture A$^- (1)$ ne d\'ecoule pas trivialement de A$^- (0)$ est le groupe special unitaire en 
$3$ variables $SU(2,1)$. C'est donc sur ce groupe que porte une bonne partie de notre \'etude. Nous d\'emontrons
notamment~:

\medskip

\noindent
{\bf Th\'eor\`eme 1} {\it Si $G$ est un groupe comme dans la Conjecture A avec $G^{{\rm nc}} \cong SU(2,1)$,
alors la Conjecture A$^-$ est vraie pour $G$, avec $\varepsilon ( G , 0)= \frac{84}{25}$ et $\varepsilon ( G , 1)
= \frac{9}{25}$.} \footnote{Dans \cite{BergeronClozel} on a annonc\'e \`a tort que $\varepsilon (G , 0) =3$ est optimal, nous verrons qu'en fait la 
valeur propre $3$ n'intervient pas dans le spectre automorphe de $G$. Ceci explique \'egalement que la Conjecture A ci-dessus est plus forte que 
celle figurant dans \cite{BergeronClozel}.}

\medskip

En fait ce r\'esultat est essentiellement contenu dans l'article \cite{HarrisLi} de Harris et Li. On en donne ici une preuve
compl\`ete. Celle-ci repose sur la th\'eorie des repr\'esentations, on profite de cette occasion pour d\'etailler les liens
entre le spectre du laplacien de Hodge sur les espaces localement sym\'etriques et la th\'eorie des repr\'esentations 
des groupes semi-simples.

\medskip

Le chapitre 1 rappelle ainsi la formule de Matsushima telle que nous en aurons besoin, on en donne 
une preuve compl\`ete afin de faciliter la lecture. On rappelle \'egalement dans cette section les d\'efinitions de base
de th\'eorie des repr\'esentations que nous utiliserons dans le texte.

\medskip

Le chapitre 2 aborde l'\'etude du spectre du laplacien sur les quotients arithm\'etiques via la th\'eorie des 
repr\'esentations. On rappelle beaucoup de r\'esultats de base de th\'eorie des repr\'esentations. La d\'emonstration 
de la Conjecture $\tau$ y est esquiss\'ee. Enfin on montre que l'extension na\"{\i}ve de 
la Conjecture A \`a des groupes plus g\'en\'eraux est fausse. Il serait d'ailleurs int\'eressant de d\'egager une 
conjecture g\'en\'erale pour tous les espaces localement symm\'etriques.

\medskip

Le chapitre 3 est quant \`a lui consacr\'e \`a la th\'eorie des repr\'esentations de $GL(n)$. On se contente le 
plus souvent de citer les r\'esultats dont nous aurons besoin, on a n\'eanmoins fait un effort pour rendre les r\'esultats
maniables \`a un lecteur pr\^et \`a les adopter comme ``bo\^{i}te noire''.

\medskip

Le chapitre 4 est l'analogue du chapitre pr\'ec\'edent pour le groupe $U(n,1)$, il est n\'eanmoins beaucoup 
plus d\'etaill\'e, le groupe $U(n,1)$ jouant le r\^ole principal dans le texte. Un certain nombre de r\'esultat
classiques rappel\'es ici sont superflus pour la suite, on a n\'eanmoins tenu \`a inclure une description d\'etaill\'ee de ceux-ci afin 
de pouvoir utiliser ce texte dans des travaux ult\'erieurs. On esp\`ere que l'effort fait de collecter ici des r\'esultats \'epars
dans la litt\'erature pourra \'egalement \^etre utile par ailleurs.

\medskip

Le chapitre 5 g\'en\'eralise certains points du chapitre pr\'ec\'edent aux groupes unitaires plus 
g\'en\'eraux. Il est bien connu qu'en rang sup\'erieur la propri\'et\'e (T) de Kazhdan assure l'existence d'un 
trou spectral pour le laplacien sur les fonctions sur les vari\'et\'es localement sym\'etriques. Dans ce chapitre nous 
montrons que pour les groupes unitaires $SU(p,q)$ d\`es que $p\geq q\geq 2$ (donc d\`es que le rang est sup\'erieur ou 
\'egal \`a $2$) un ph\'enom\`ene analogue appara\^{\i}t dans l'\'etude du spectre du laplacien sur les formes 
diff\'erentielles de suffisamment petit degr\'e. Nous d\'eduisons ce Th\'eor\`eme d'un r\'esultat plus g\'en\'eral, \`a para\^{\i}tre (bien que datant
du d\'ebut des ann\'ees 90), de Vogan. Ce r\'esultat illustre \'egalement le r\^ole sp\'ecial jou\'e par les groupes 
$U(n,1)$ et $O(n,1)$ dans ces questions. Groupes pour lesquels on ne peut esp\'erer d'isolation spectrale structurelle
mais pour lesquels il faut \'etudier le spectre automorphe. 

\medskip

Motiv\'e par le chapitre pr\'ec\'edent, le chapitre 6 est consacr\'e \`a la description du yoga d'Arthur pour les
groupes $U(n,1)$ et $O(n,1)$. Nous pr\'ecisons notre interpr\'etation des Conjectures d'Arthur pour ces deux groupes, ce 
qui permet de conforter la Conjecture A.

\medskip

La suite de la premi\`ere partie du texte vise \`a d\'emontrer le Th\'eor\`eme 1 et les r\'esultats que nous d\'ecrivons ci-dessous.
Le Th\'eor\`eme 1 d\'ecoule des travaux de Rogawski permettant de r\'eduire l'\'etude du spectre automorphe de $U(2,1)$
au spectre automorphe de $GL(3)$ par changement de base. Il reste alors \`a contr\^oler le spectre automorphe de $GL(n)$.
La conjecture standard dans ce cas est la Conjecture de Ramanujan. On a besoin d'une approximation de celle-ci.
Une telle approximation nous est fournie par un th\'eor\`eme de Luo, Rudnick et Sarnak.  On a en fait besoin d'une g\'en\'eralisation 
de celui-ci. La d\'emonstration de cette g\'en\'eralisation fait l'objet du chapitre 7.

\medskip

On peut alors dans le chapitre 8 donner une d\'emonstration compl\`ete du Th\'eor\`eme 1. Nous profitons 
\'egalement de ce chapitre pour revenir sur la Conjecture de changement de base faite par Harris et Li dans \cite{HarrisLi}, que 
l'on \'etend l\'eg\`erement et que l'on compare aux Conjectures d'Arthur. Nous nous servons de cette conjecture dans la 
deuxi\`eme partie du texte pour \'etendre la Conjecture B \`a d'autres groupes semi-simples.

\medskip

Pour chaque r\'eel strictement positif $\varepsilon$, notons A$_{\varepsilon}^p$ \index{hypoth\`ese A$_{\varepsilon}^p$} l'hypoth\`ese spectrale suivante sur le 
groupe $G$.
\begin{eqnarray*}
{\rm A}_{\varepsilon}^p \; : \; 
\left\{ 
\begin{array}{l}
\mbox{si } \lambda \mbox{ est dans le } p-\mbox{spectre automorphe de } G, \mbox{ alors :} \\
\mbox{1. soit } \lambda = (\rho_G - k)^2 - (\rho_G - i)^2 \mbox{ avec } k=0, \ldots ,p \mbox{ et } i=k, \ldots , [ \rho_G ]  , \\
\mbox{2. soit } \lambda \geq (\rho_G -p)^2 - \varepsilon^2 .
\end{array}
\right.
\end{eqnarray*}
Nous expliquerons que les Conjectures d'Arthur semblent pr\'evoir que, pour un groupe $G$ comme
ci-dessus, chaque hypoth\`ese A$_0^p$ est vraie pour $p \leq [\rho_G ]$. Les hypoth\`eses A$_{\varepsilon}^p$
sont donc des approximations aux Conjectures d'Arthur pour le groupe $G$, lorsque 
$\varepsilon$ est suffisamment petit, ces approximations sont en g\'en\'eral plus pr\'ecises que la Conjecture A$^-$.
Le but du chapitre 9 est la d\'emonstration du th\'eor\`eme suivant qui permet de relever les approximations
aux Conjectures d'Arthur.

\medskip

\noindent
{\bf Th\'eor\`eme 2} {\it Soient $H  < G$ deux groupes ${\Bbb Q}$-alg\'ebriques comme dans la Conjecture B.
Supposons que $H$ v\'erifie l'hypoth\`ese A$_{\varepsilon}^p$ pour un certain r\'eel 
$\varepsilon >0$ et pour tout entier naturel $p \leq [\rho_H ]$.  

Alors le groupe $ G$ v\'erifie les hypoth\`eses A$_{\rho_G -\rho_H +\varepsilon}^p$  
pour tout entier naturel $p\leq [\rho_H ]$.
}

\medskip

\`A l'aide des Th\'eor\`emes 1 et 2 et de l'analyse de certains groupes unitaires sp\'eciaux (comme ceux intervenant 
dans la preuve de la Conjecture $\tau$) nous montrons dans le chapitre 10 le th\'eor\`eme suivant.

\medskip

\noindent
{\bf Th\'eor\`eme 3} {\it Si $G$ est un groupe comme dans la Conjecture A avec $G^{{\rm nc}} \cong SU(n,1) $ (et quelques conditions techniques 
suppl\'ementaires si $n+1$ est un puissances de $2$),
alors les Conjectures A$^- (0)$ et A$^- (1)$ sont v\'erifi\'ees par $G$, avec $\varepsilon (G ,0)=2n-1$ et 
$\varepsilon (G ,1) = \frac{10n-11}{25}$.} 

\medskip

Ceci conclut la premi\`ere partie, partie spectrale, du texte. Dans la deuxi\`eme partie nous nous int\'eressons \`a l'homologie des
vari\'et\'es de congruence. Le r\'esultat principal est le th\'eor\`eme suivant.

\medskip

\noindent
{\bf Th\'eor\`eme 4} {\it La Conjecture A$^-$ implique la Conjecture B. 

Plus pr\'ecisemment, soit $G$ un groupe comme dans la Conjecture A et $i$ un entier naturel $\leq \frac{d_G}{2}$.
Si la Conjecture A$^-_{d=0} (i)$ est v\'erifi\'ee par $G$, alors quel que soit $\Gamma$ sous-groupe de congruence 
sans torsion de $G$, la vari\'et\'e $\Gamma \backslash X_G$ v\'erifie la Conjecture B$(d_G -i )$.}

\medskip

Nous ne prouverons le Th\'eor\`eme 4 que pour $G^{{\rm nc }} \cong SU(n,1)$. Le cas o\`u $G^{{\rm nc }} \cong SO(n,1)$,
similaire, est plus simple et compl\`etement trait\'e dans \cite{MathZ}. Remarquons n\'eanmoins que puisque 
la Conjecture A$^-_{d=0} (1)$ est d\'emontr\'ee, le Th\'eor\`eme 4 implique le c\'el\`ebre Th\'eor\`eme de Millson selon
lequel toute vari\'et\'e hyperbolique (r\'eelle) arithm\'etique standard admet un rev\^etement fini dont le premier nombre
de Betti est non nul.

\medskip

Nous d\'emontrons en fait un r\'esultat plus g\'en\'eral que le Th\'eor\`eme 4, voir chapitre 15 Th\'eor\`eme \ref{sur l'homologie}.
Pour cela, apr\`es un premier chapitre consacr\'e \`a l'espace hyperbolique complexe principal objet de cette partie, nous consacrons un long chapitre 
aux espaces sym\'etriques associ\'es aux groupes unitaires, nous suivons pour l'essentiel les articles de Tong et Wang cit\'es
dans la bibliographie. C'est \'egalement en suivant les travaux de Tong et Wang que nous construisons dans le chapitre
13 la forme duale aux sous-vari\'et\'es totalement g\'eod\'esiques que nous consid\`erons dans le texte.

\medskip

Alli\'e au Th\'eor\`eme 3, le Th\'eor\`eme 4 nous permet en fin de deuxi\`eme partie de d\'emontrer le th\'eor\`eme 
suivant qui est la premi\`ere avanc\'ee significative vers la Conjecture B.

\medskip

\noindent
{\bf Th\'eor\`eme 5} {\it Soit $M$ une vari\'et\'e hyperbolique complexe compacte de congruence de dimension $d+2$.
Soit $F$ une sous-vari\'et\'e complexe compacte connexe totalement g\'eod\'esique de dimension $d$ et immerg\'ee dans $M$.
Alors, il existe un rev\^etement fini $\hat{M}$ de $M$ tel que~:
\begin{enumerate}
\item l'immersion de $F$ dans $M$ se rel\`eve en un plongement de $F$ dans $\hat{M}$,
\item l'application induite~:
$$H_{d-1} (F ) \rightarrow H_{d-1} (\hat{M} )$$
est injective.
\end{enumerate}

De plus, pour tout entier $N$ et tout cycle $c$ dans $H_{d-1} (\hat{F} )$, il existe un rev\^etement fini 
$M_N$ de $M$ contenant $N$ relev\'es de $c$ lin\'eairement ind\'ependants dans $H_{d-1} (M_N )$.}

\medskip

On en d\'eduit alors facilement le corollaire suivant.

\medskip

\noindent
{\bf Corollaire 6} {\it Soient $H  < G$ deux groupes comme dans la Conjecture B$^-$. Supposons que 
$$\left\{ 
\begin{array}{c}
H^{{\rm nc}} \cong SU(n-1,1)  \\
G^{{\rm nc}} \cong SU(n,1)  , 
\end{array} \right. $$
alors l'application naturelle 
$$H_{2n-3} (Sh^0 H) \rightarrow H_{2n-3} (Sh^0 G)$$
est injective.}

\bigskip

Remarquons que la d\'emonstration des th\'eor\`emes ci-dessus repose de mani\`ere essentielle sur l'\'etude de la cohomologie
$L^2$ de quotients $V= \Gamma \backslash X$, o\`u $X$ est l'espace hyperbolique complexe de dimension complexe $n$
et $\Gamma$ un groupe discret pr\'eservant un sous-espace totalement g\'eod\'esique de $X$ et agissant cocompactement
sur celui-ci. Un cas particulier de ce que nous d\'emontrons au chapitre 14 est le th\'eor\`eme suivant.

\medskip

\noindent
{\bf Th\'eor\`eme 7} {\it Pour tout entier $k<n$, il existe un isomorphisme naturel~:
$${\cal H}_{(2)}^k (V) \stackrel{\backsimeq}{\rightarrow} H^k_c (V ) .$$
De plus, l'espace ${\cal H}_{(2)}^n$ est de dimension infinie.}

\medskip

Le th\'eor\`eme que nous d\'emontrons est plus g\'en\'eral; il repose de mani\`ere essentielle sur un th\'eor\`eme de
Ohsawa et Tanigushi. 

\medskip

Les d\'emonstrations des Th\'eor\`emes 4 et 5 occupent le chapitre 15. Bon nombre des techniques 
d\'evelopp\'ees dans cette deuxi\`eme partie s'\'etendent \`a des espaces sym\'etriques plus g\'en\'eraux. 
Gr\^ace aux r\'esultats spectraux obtenus au chapitre 5 concernant les groupes unitaires $U(p,q)$, nous 
d\'emontrons \'egalement dans ce chapitre le th\'eor\`eme suivant.

\medskip

\noindent
{\bf Th\'eor\`eme 8} {\it Soient $p,q,r$ trois entiers naturels v\'erifiant $p\geq q \geq 2$ et $p\geq r$.
Soient $H$, $G$ deux groupes alg\'ebriques sur ${\Bbb Q}$ tels que 
\begin{itemize}
\item $H$ soit un ${\Bbb Q}$-sous-groupe de $G$,
\item $H^{{\rm nc}} \cong SU(p-r,q)$,
\item $G^{{\rm nc }} \cong SU(p,q)$. 
\end{itemize}
Alors, l'application naturelle 
$$H_{i} (Sh^0 H ) \rightarrow H_i (Sh^0 G )$$
est injective pour tout entier $i > 2pq-p-q+1$.}

\medskip

Il serait \'evidemment int\'eressant de formuler une conjecture g\'en\'erale suivant le principe de la Conjecture B
mais pour tout groupe semi-simple. La combinatoire devrait
\^etre plus compliqu\'ee. Dans un article en pr\'eparation le premier des deux auteurs \'etudie ce 
probl\`eme en portant une attention particuli\`ere au cas des groupes unitaires et orthogonaux.

\newpage

\tableofcontents

\newpage

\thispagestyle{empty}

$\mbox{}$

\newpage

\pagenumbering{arabic}

\part{Spectre des vari\'et\'es hyperboliques}

\chapter{Th\'eor\`eme de Matsushima}

Soit $X$ un espace sym\'etrique simplement connexe de courbure n\'egative sans facteurs euclidiens 
et soit $G$ un groupe de 
Lie semi-simple r\'eel connexe agissant transitivement sur $X$ par isom\'etries.
Nous supposons que l'application de $G$ vers la composante connexe de l'identit\'e du groupe des 
isom\'etries de $X$ est un rev\^etement fini.
Dans la suite nous supposons que la m\'etrique riemannienne de $X$ est identique \`a celle induite 
par la forme de Killing de $G$. Nous 
notons $K$ le groupe d'isotropie dans $G$ d'un point fix\'e $p$ de $X$.
Puisque $X$ est de courbure n\'egative, le groupe $G$ est de centre fini et sans 
facteur compact. 

Soit $\Gamma$ un sous-groupe discret de $G$ tel que $\Gamma \backslash X$ soit 
compact. Soit $\omega$ une $p$-forme diff\'erentielle sur $X$ invariante par $\Gamma$
{\it i.e.} 
$$\omega (x;v_1 , \ldots , v_p ) = \omega (\gamma x ;\gamma v_1 , \ldots , \gamma v_p )
\mbox{  pour tout } \gamma \in \Gamma .$$
Notons $\pi : G \rightarrow G/K =X$ la projection canonique. 
Soit $\tilde{\omega} = \pi^* \omega$. La forme $\tilde{\omega}$ est une $p$-forme
sur $G$.

Nous adoptons la convention classique de noter avec un indice $0$ les alg\`ebres de Lie r\'elles et si $\mathfrak{l}_0$ est une alg\`ebre de Lie r\'elle, 
$\mathfrak{l}$ l'alg\`ebre de Lie complexe $\mathfrak{l} = \mathfrak{l}_0 \otimes {\Bbb C}$.
Pour tout $g\in G$, la translation \`a gauche par $g^{-1}$ fournit un isomorphisme canonique 
de $T_g (G)$ vers $\mathfrak{g}_0 = T_e (G)$ et donc une identification 
\begin{eqnarray} \label{101}
\Omega^p (G; {\Bbb C})={\rm Hom} (\Lambda^p (\LG ) , C^{\infty} (G;{\Bbb C})).
\end{eqnarray}
Fixons $\LG = \LK \oplus \LM$ une d\'ecomposition de Cartan de $\LG_0$. 

La forme $\tilde{\omega}$ provient d'une forme $\omega$ sur $G/K$ invariante
sous l'action (\`a gauche de $\Gamma$) et via l'identification (\ref{101}) appartient \`a :
$$C^p (\LG ,K;C^{\infty} (\Gamma \backslash G ; {\Bbb C}) ) :=
{\rm Hom}_K (\Lambda^p (\LM ) , C^{\infty} (\Gamma \backslash G;{\Bbb C})),$$
o\`u $K$ agit sur $\mathfrak{p} \cong \mathfrak{g} /\mathfrak{k}$ via la repr\'esentation
adjointe et sur $C^{\infty} (\Gamma \backslash G ; {\Bbb C})$ via la repr\'esentation r\'eguli\`ere 
\`a droite. R\'eciproquement, il n'est pas difficile de voir qu'\`a un \'el\'ement de 
$ C^p (\LG ,K;C^{\infty} (\Gamma \backslash G ; {\Bbb C}))$ correspond une $p$-forme 
sur $G$ provenant d'une $p$-forme sur $G/K$ et invariante sous l'action de $\Gamma$.
Nous r\'esumons tout ceci dans la proposition suivante qui semble avoir \'et\'e pour la premi\`ere fois
constat\'ee par Gelfand et Fomin \cite{GelfandFomin} puis pr\'ecis\'ee par Matsushima \cite{Matsushima}.

\begin{prop} \label{P101}
L'application $\omega \mapsto \omega \circ \pi$ induit un isomorphisme de 
complexes gradu\'es de $\Omega^* (\Gamma \backslash X; {\Bbb C})$
sur $C^* (\LG ,K;C^{\infty} (\Gamma \backslash G ; {\Bbb C}))$.
\end{prop}
\index{Proposition de Gelfand et Fomin}

Notons $\Delta$ le laplacien de Hodge sur les formes diff\'erentielles sur $X$ (pour sa structure riemannienne).
Soit $\omega$ une $p$-forme sur $X$ telle que $\Delta \omega = \lambda \omega$
pour un certain $\lambda \in {\Bbb R}$. D'apr\`es la Proposition \ref{P101}, il correspond \`a 
$\omega$ un \'el\'ement (que l'on note toujours) $\omega$ dans $C^p (\LG ,K;C^{\infty} (\Gamma \backslash G ; {\Bbb C}))$ 

D'un autre c\^ot\'e, il est bien connu (cf. \cite{GGPP}) que la repr\'esentation r\'eguli\`ere
droite dans $L^2 (\Gamma \backslash G )$ se d\'ecompose
en :
\begin{eqnarray} \label{GPP}
L^2 (\Gamma \backslash G) = \bigoplus_{\pi \in \hat{G}} m(\pi , \Gamma ) {\cal H}_{\pi } ,
\end{eqnarray}
somme discr\`ete de sous-espaces $G$-invariants irr\'eductibles, 
index\'ee par le dual unitaire $\hat{G}$ de $G$, et o\`u chaque
$m(\pi , \Gamma )$ est fini. 

On voit ainsi se dessiner une correspondance entre certaines repr\'esentations de $G$
et l'espace des formes diff\'erentielles $\lambda$-propres pour le laplacien. Avant de 
donner un \'enonc\'e pr\'ecis, rappelons que l'{\it op\'erateur de Casimir} \index{Casimir, op\'erateur de} est 
$$C= \sum_{1\leq s \leq n} y_s . y_s '$$
o\`u $(y_s)$ est une base de $\LG$ et $(y_s ')$ la base duale de $\LG$ par rapport 
\` a la forme de Killing $B$. 

La (tr\`es) l\'eg\`ere modification du th\'eor\`eme de Matsushima  \cite{Matsushima} \index{Matsushima, Th\'eor\`eme de} dont nous aurons besoin
s'\'enonce alors :

\begin{thm} \label{mat}
Soit $E_{\lambda}^p$ l'espace des $p$-formes diff\'erentielles sur $\Gamma \backslash
X$, $\lambda$-propres pour le laplacien de $X$.
Alors :
$${\rm dim}(E_{\lambda }^p ) = \sum_{
\begin{array}{c}
\pi \in \hat{G} \\
\pi (C) = -\lambda
\end{array}
} 
m(\pi , \Gamma ) {\rm dim}( {\rm Hom}_{K} (\Lambda^p \LM , {\cal H}_{\pi} )) ,$$
la somme \'etant finie.
\end{thm}
\index{Th\'eor\`eme de Matsushima}
 
\begin{cor}
Soit $G= SO(n,1)$ ou $SU(n,1)$.
Si pour toute repr\'esentation $\pi$ apparaissant dans $L^2 (\Gamma \backslash G)$ et dont le $K$-type
rencontre $\Lambda^p \mathfrak{p}$, on a $-\pi (C) > \epsilon$ ou $\pi (C)=0$, la 
Conjecture A$^-$ est v\'erifi\'ee pour les $p$-formes.
\end{cor}

Avant de d\'emontrer le Th\'eor\`eme \ref{mat}, revenons sur l'op\'erateur de Casimir.

\section{Sur l'op\'erateur de Casimir}

Identifions dor\'enavant l'espace $\LM$ au suppl\'ementaire orthogonal de $\LK$ dans $\LG$
par rapport \`a la forme de Killing $B$.
Soit $\{  X_u \}$ (resp. $\{ X_t \}$) une base orthonormale de $\LM$ (resp. $\LK$)
pour la forme $B$ (resp. $-B$). Alors :
$$C=\sum X_u^2 - \sum X_t^2 .$$ \index{Casimir, op\'erateur de}

Soit $T$ est une repr\'esentation unitaire de $G$ dans un Hilbert ${\cal H}$.
Un vecteur $\varphi \in {\cal H}$ est dit {\it lisse} \index{lisse} ou {\it $C^{\infty}$} \index{$C^{\infty}$} si la fonction $g \mapsto T(g) \varphi$ est $C^{\infty}$.
Soit ${\cal H}^{\infty}$ le sous-espace des vecteurs lisses de ${\cal H}$.

Si $X\in \LG$ et $\varphi \in {\cal H}^{\infty}$ nous posons~:
$$T(X)\varphi = \left[ \frac{d}{dt} T(\exp (tX)) \varphi \right] _{|t=0} .$$
Alors :
$$T(C) := \sum T(X_u )^2 -\sum T(X_t )^2 .$$ \index{Casimir, op\'erateur de}

Remarquons que si la repr\'esentation $T$ est irr\'eductible, alors $T(C)$ est scalaire. Cela d\'ecoule 
du Lemme de Schur et du fait que $C$ est dans le centre de l'alg\`ebre enveloppante de
$\LG$.

Le lemme suivant, d\^u \`a Kuga \cite{BorelWallach}, fait le lien entre le laplacien et l'op\'erateur de Casimir.

\begin{lem} \label{L111}
Via l'application de la Proposition \ref{P101}, les $p$-formes $\lambda$-propres pour le laplacien sur 
$\Gamma \backslash X$ correspondent aux \'el\'ements $\omega$ de 
$C^p (\LG , K; C^{\infty} (\Gamma \backslash G))$ tels que 
$$R(C) \omega = -\lambda \omega ,$$
o\`u $R$ est la repr\'esentation r\'eguli\`ere droite de $G$ dans $C^{\infty} (\Gamma \backslash G)$.
\end{lem} \index{Lemme de Kuga}
{\it D\'emonstration.} Soient $\eta^u$ les $1$-formes invariantes \`a gauche sur $G$ telles que 
$$\eta^u (X_v ) = \delta_u^v .$$
Soit 
$$\omega \in C^p (\LG , K;C^{\infty} (\Gamma \backslash G) )={\rm Hom}_K (\Lambda^* \LM ,C^{\infty} (\Gamma \backslash G))$$
telle que (vue comme forme diff\'erentielle) 
$\Delta \omega = \lambda \omega $.

On \'ecrit 
$$\omega = \sum_U \omega_U \eta^U ,$$
o\`u  si $U=(u_1 , \ldots , u_p )$, $\eta^U = \eta^{u_1} \wedge \ldots \wedge \eta^{u_p} $.
On a imm\'ediatement :
$$(d \omega )_U = -\sum_{1\leq v \leq p+1} (-1)^{v-1} R(X_v ) \omega_{U(v)} $$
o\`u $U(v)$ est obtenu en enlevant aux \'el\'ements de $U$ la $v$-i\`eme coordonn\'ee
(les termes en crochets de Lie dans la diff\'erentielle \'etant nuls car $\LG = \LK \oplus \LM$ avec
$[\LM , \LM] \subset \LM$).

De plus, puisque la m\'etrique riemannienne sur $X$ est induite par la forme de Killing $B$
sur $\LM$, on a :
$$<\omega , \omega '> = \sum_U \int_{\Gamma \backslash G} \omega_U \overline{\omega}_U ' .$$

Un calcul simple implique alors le fait suivant.

\medskip

\noindent
{\bf Fait.} L'adjoint de $d$ est l'op\'erateur $\partial$ donn\'ee par :
$$(\partial \omega )_U = \sum_v R(X_v )^* \omega_{\{ v\} \cup U} $$
o\`u $U$ est cette fois de cardinal $p-1$ et o\`u $R(X_v )^*$ est l'adjoint de $R(X_v )$ (agissant sur $ C^{\infty} (\Gamma \backslash G) $
muni de la norme $L^2$) {\it i.e.} $R(X_v )^*= -R(X_v )$.

\medskip

Finalement :
\begin{eqnarray*}
(\partial d \omega )_U & = & -\sum R(X_v ) (d\omega )_{\{ v\} \cup U} \\
                                    & = & -\sum_v R(X_v ) \left( R(X_v ) \omega_U - \sum_w (-1)^w R(X_{v_w} ) \omega_{\{ v\} \cup U(w)} \right) \\
                                    & = & -\sum_v R(X_v )^2 \omega_U + \sum_{v,w} (-1)^u R(X_v ) R(X_{v_w} ) \omega_{\{ v \} \cup U(w )}  
\end{eqnarray*} 
et
\begin{eqnarray*}
(d \partial \omega )_U & = & - \sum_v (-1)^{v-1} R(X_{w_v} ) (\partial \omega )_{U(v)} \\
                                    & = & \sum_{v,w} (-1)^{v-1} R(X_{w_v} ) R(X_w ) \omega_{\{ w\} \cup U(v) } .
\end{eqnarray*} 
Or, 
$$R(X_v ) R(X_{v_w} ) -R(X_{v_w} ) R(X_v ) = R([X_v , X_{v_w} ]) \in R ( \LK ) $$
agit trivialement sur $\omega$ qui est une forme diff\'erentielle sur $X$. Donc :
\begin{eqnarray*}
\Delta \omega & = & - \left( \sum_v R(X_v )^2 \right) \omega \\
                        & = & - R(C) \omega ,
\end{eqnarray*}
car, l\`a encore, $R(X)$ agit trivialement sur $\omega$ si $X \in \LK$. 
Nous avons donc v\'erifi\'e qu'\`a travers la correspondance de la Proposition \ref{P101}, les formes
propres pour le laplacien correspondent \`a des formes propres pour l'action de $R(C)$, les valeurs propres
\'etant de signes oppos\'es.

\section{D\'emonstration du Th\'eor\`eme \ref{mat}}

Soit $E_{\lambda}^p$ l'espace des $p$-formes de carr\'e int\'egrable et $\lambda$-propres.
C'est un espace de dimension finie de formes $C^{\infty}$, le laplacien \'etant \'elliptique. De la 
d\'ecomposition $L^2$ donn\'ee par (\ref{GPP}), on d\'eduit~:
\begin{eqnarray} \label{121}
E_{\lambda}^p & = & \oplus_{\pi} m( \pi , \Gamma ) {\rm Hom}_{K} ( \Lambda^p \LM , H_{\pi} )_{\lambda} ,
\end{eqnarray}
$E_{\lambda}^p$ \'etant vu comme un espace de formes $L^2$, et l'indice $\lambda$ d\'esignant le sous-espace
sur lequel $C$ op\`ere par $-\lambda$. Puisque $C$ peut s'identifier \`a un op\'erateur du centre de l'alg\`ebre 
enveloppante de $\LG$, cet espace n'est non nul que si $\pi (C)=-\lambda$.

Par cons\'equent dim$(E_{\lambda}^p )$ est donn\'ee par la formule du Th\'eor\`eme \ref{mat}, et la somme est finie.
On a en fait canoniquement :
$$E_{\lambda}^p \cong \bigoplus_{\pi \, : \, \pi (C) = -\lambda} {\rm Hom} (\pi , L^2 (\Gamma \backslash G) ) 
\otimes {\rm Hom}_{K} ( \Lambda^p \LM , {\cal H}_{\pi} ) , $$
o\`u toutes les formes (a priori $L^2$) dans cette somme sont $C^{\infty}$ puisque $C$ est elliptique.

\pagestyle{myheadings}

\markright{1.3. REPR\'ESENTATIONS ADMISSIBLES}

\section{Repr\'esentations admissibles et spectre automorphe}

On dit qu'une repr\'esentation $\pi$ de $G$ dans un espace de Hilbert ${\cal H}$ \footnote{Attention : la repr\'esentation 
$\pi$ n'a pas de raison d'\^etre unitaire.} est {\it admissible} \index{admissible} si $\pi (K)$ op\`ere par op\'erateurs unitaires et si 
chaque repr\'esentation irr\'eductible $\tau$ de $K$ n'intervient qu'avec une multiplicit\'e finie (peut-\^etre nulle) dans la restriction 
$\pi_{|K}$ de $\pi$ \`a $K$.

Un th\'eor\`eme de Harish-Chandra \cite{HarishChandra} (voir aussi \cite{Knapp}) affirme que toute repr\'esentation irr\'eductible unitaire 
$(\pi , {\cal H}_{\pi})$ de $G$ est admissible.  Rappelons que l'espace ${\cal H}_{\pi}^{\infty}$ des vecteurs lisses
de $\pi$ est un sous-espace dense de ${\cal H}_{\pi}$ invariant sous l'action de $G$, et qu'il existe une action 
naturelle de $\mathfrak{g}$ sur ${\cal H}_{\pi}^{\infty}$. Posons 
$${\cal H}_{\pi}^K = \{ v\in {\cal H}_{\pi}^{\infty} \; : \; {\rm dim} \langle \pi (K) v \rangle < +\infty \},$$
o\`u $\langle \pi (K) v \rangle$ d\'esigne le sous-espace engendr\'e par tous les vecteurs de la forme $\pi (k) v$, avec $k\in K$.
D'apr\`es Harish-Chandra \cite{HarishChandra} ${\cal H}_{\pi}^{K}$ est stable sous les actions de $K$ et $\mathfrak{g}$
et comme $\mathfrak{g}$-module il est irr\'eductible, et d\'etermine $\pi$ \`a \'equivalence unitaire pr\`es. En fait 
${\cal H}_{\pi}^{K}$ est un $(\mathfrak{g}, K)$-module, voir \cite{BorelWallach}. \index{($\mathfrak{g} , K)$-module}

\'Etant donn\'e un entier $i$, nous notons 
$$C^i (\pi ) = {\rm Hom}_K (\Lambda^i \mathfrak{p} , {\cal H}_{\pi}^K ) .$$
On a vu dans les sections pr\'ec\'edentes qu'\`a une repr\'esentation $\pi \subset L^2 (\Gamma \backslash G)$ telle que
$C^i (\pi ) \neq 0$ on peut associer une forme diff\'erentielle de degr\'ee $i$ sur $\Gamma \backslash X$. On les 
obtient de la fa\c{c}on suivante. Soit $\varphi : {\cal H}_{\pi}^K \rightarrow C^{\infty} (\Gamma \backslash G)$ l'application
$G$-\'equivariante induite par l'inclusion $\pi \subset L^2 (\Gamma \backslash G)$. 
Fixons une application $K$-\'equivariante non nulle $\omega : \Lambda^i \mathfrak{p} \rightarrow {\cal H}_{\pi}^K$ 
(dans $C^i (\pi )$). La forme diff\'erentielle associ\'ee $\omega_{\varphi}$ est d\'efinie par 
$$\omega_{\varphi} (g. \lambda ) = \varphi(\omega (\lambda ))(g) \; \; (\lambda \in \Lambda^i \mathfrak{p}, \; g\in G).$$
Remarquons que le Th\'eor\`eme 1 implique que si $\pi (C)=0$, la forme $\omega_{\varphi}$ est harmonique de 
degr\'e $i$ et contribue donc \`a la cohomologie de degr\'e $i$ de $\Gamma \backslash X$. On dira donc d'une 
repr\'esentation irr\'eductible unitaire de $G$ qu'elle est {\it cohomologique de degr\'e $i$} \index{repr\'esentation cohomologique} si elle v\'erifie que
\begin{enumerate}
\item $\pi (C) = 0$, et 
\item $C^i (\pi ) \neq 0$.
\end{enumerate}

\medskip

On notera $U (\mathfrak{g} )$ l'alg\`ebre enveloppante sur ${\Bbb C}$ de l'alg\`ebre de Lie $\mathfrak{g}$.
Son centre sera not\'e ${\cal Z} (\mathfrak{g} )$.
 
Rappelons que l'on dit de la repr\'esentation $(\pi , {\cal H}_{\pi} )$ qu'elle admet un {\it caract\`ere infinit\'esimal 
$\chi$} \index{caract\`ere infinit\'esimal} si $\chi$ est un homomorphisme ${\cal Z} (\mathfrak{g}) \rightarrow {\Bbb C}$ tel que 
$\pi (z) = \chi (z) id$ pour tout $z \in {\cal Z} (\mathfrak{g} )$. C'est en particulier toujours le cas si $(\pi , {\cal H}_{\pi} )$
est irr\'eductible et admissible. 

Si $\mathfrak{h}$ est une sous-alg\`ebre de Cartan de l'alg\`ebre de Lie complexe $\mathfrak{g}$, 
l'homomorphisme d'Harish-Chandra $\gamma$, cf. \cite{Knapp}, r\'ealise un isomorphisme d'alg\`ebre de ${\cal Z} (\mathfrak{g} )$
sur la sous-alg\`ebre de $U( \mathfrak{h} )$ constitu\'ee des \'el\'ements fix\'es par le groupe de Weyl $W=W(\mathfrak{h} ,
\mathfrak{g})$. Il correspond donc \`a chaque forme lin\'eaire $\Lambda$ dans $\mathfrak{h}'$ un caract\`ere de ${\cal Z} (\mathfrak{g} )$~:
\begin{eqnarray} \label{chi lambda}
\chi_{\Lambda} (z) = \Lambda (\gamma (z) ) \mbox{  pour  } z \in {\cal Z} (\mathfrak{g} ),
\end{eqnarray}
o\`u l'on a \'etendu $\Lambda$ en un homomorphisme d'alg\`ebre de $U ( \mathfrak{h})$.

Le morphisme $\chi_{\Lambda}$ v\'erifie $\chi_{w \Lambda } = \chi_{\Lambda }$ pour tout $w \in W$ et r\'eciproquement,
si $\chi_{\Lambda '} = \chi_{\Lambda}$ alors $\Lambda ' = w \Lambda $ pour un \'el\'ement $w \in W$ (voir \cite[Proposition 8.21]{Knapp}).

On montre de plus, cf. \cite[Proposition 8.21]{Knapp}, que tout homomorphisme de ${\cal Z} (\mathfrak{g})$ dans 
${\Bbb C}$ est de la forme $\chi_{\Lambda}$, comme dans (\ref{chi lambda}), pour un certain $\Lambda$
dans $\mathfrak{h}'$. Si $\pi$ admet un caract\`ere infinit\'esimal $\chi = \chi_{\Lambda}$, on dira, par 
abus de notation, que $\pi$ a pour caract\`ere infinit\'esimal $\Lambda$; dans ce cas $\Lambda$ n'est d\'etermin\'e
qu'\`a un membre du groupe de Weyl $W$ pr\`es. 

\bigskip

On a donc ramen\'e l'\'etude du spectre des quotients $\Gamma \backslash X$ \`a celle de la d\'ecomposition en 
irr\'eductibles des repr\'esentations $L^2 (\Gamma \backslash G)$ de $G$. Comme annonc\'e dans l'introduction 
nous n'allons consid\'erer que des sous-groupes $\Gamma$ de congruence. 

Supposons maintenant que le groupe r\'eel $G$ est \'egal au groupe des points r\'eels $G({\Bbb R})$ d'un groupe alg\'ebrique (toujours not\'e) $G$ connexe et presque simple sur ${\Bbb Q}$. Nous adoptons la m\^eme notation ($G$) pour d\'esigner le groupe alg\'ebrique et le groupe r\'eel. Le contexte permettra d'\'eviter toute confusion.

On appelle {\it sous-groupe de congruence de $G$} \index{congruence, sous-groupe de} tout sous-groupe $\Gamma = G ({\Bbb Q}) \cap K_f$ de 
$G$ o\`u $K_f$ est un sous-groupe compact-ouvert du groupe $ G ({\Bbb A}_f )$ sur les 
ad\`eles finis. Remarquons que si $G$ est d\'efini sur ${\Bbb Z}$, un sous-groupe $\Gamma$ de 
$G ({\Bbb Z} )$ qui contient un sous-groupe de la forme 
$$\Gamma_N = \mbox{ker} (G ({\Bbb Z} ) \rightarrow G ({\Bbb Z} / N{\Bbb Z} ) ),$$
pour un certain entier $N \geq 1$ est un sous-groupe de congruence. 
Le groupe $\Gamma$ agit alors proprement sur $X$ (via la projection de $\Gamma$ dans $G$)
et d'apr\`es un th\'eor\`eme de Borel et Harish-Chandra 
\cite{BorelHarishChandra} le quotient $\Gamma \backslash X$ est de volume fini. De plus, ce quotient est compact 
si et seulement si le groupe $G$ est anisotrope sur ${\Bbb Q}$. 

Notons $\widehat{G}$ le dual unitaire de $G$ que l'on munit de la topologie de Fell.
On rappellera la d\'efinition et les propri\'et\'es standard de la topologie de Fell au \S 2.2. Suivant Burger et Sarnak 
\cite{BurgerSarnak}, on appelle {\it spectre} \index{spectre} de $\Gamma \backslash G$ - o\`u $\Gamma$ est un 
sous-groupe de congruence de $G$ - l'ensemble des repr\'esentations irr\'eductibles unitaires $\pi \in \widehat{G}$ 
apparaissant (faiblement si $G$ est isotrope sur ${\Bbb Q}$) dans la repr\'esentation r\'eguli\`ere 
de $G$ dans $L^2 (\Gamma \backslash G )$~:
\begin{eqnarray} \label{spectre}
\sigma ( \Gamma \backslash G({\Bbb R}) ) = \{ \pi \in \widehat{G ({\Bbb R})} \; : \; \pi \propto L^2 (\Gamma 
\backslash  G ) \}.
\end{eqnarray}
On rappelle \'egalement la d\'efinition du dual automorphe $\widehat{G}_{{\rm Aut}}$ de $G$ donn\'ee dans \cite{BurgerSarnak}~:
\begin{eqnarray} \label{dual automorphe}
\widehat{G}_{{\rm Aut}} = \overline{ \bigcup_{
\Gamma {\rm cong}
} \sigma ( \Gamma \backslash G  ) }.
\end{eqnarray}
L'adh\'erence dans \ref{dual automorphe} est prise par rapport \`a la topologie de Fell.
D'apr\`es le Th\'eor\`eme \ref{mat} les Conjectures A$^- (i)$ se ram\`enent \`a l'\'etude des spectres (\ref{spectre}) et en fait
comme nous le verrons au Chapitre 2 de tout le dual automorphe (\ref{dual automorphe}). Il s'agit de savoir si 
chaque repr\'esentation cohomologique $\pi$ de degr\'e $i$ est isol\'ee dans la r\'eunion $\{ \pi \} \cup   
\widehat{G}_{{\rm Aut}}$.

\bigskip

Concluons ce premier chapitre en remarquant que l'hypoth\`ese de compacit\'e sur le quotient $\Gamma \backslash X$
n'\'etait pas r\'eellement n\'ecessaire. Il suffit en g\'en\'eral de consid\'erer des formes diff\'erentielles de carr\'e int\'egrable, cf. \cite{BorelWallach}.
On parlera alors de spectre $L^2$.

\newpage

\thispagestyle{empty}

\newpage

\markboth{CHAPITRE 2. SPECTRE DU LAPLACIEN}{2.2. TH\'EORIE DES REPR\'ESENTATIONS}

\chapter{Spectre du laplacien sur les quotients arithm\'etiques}

Dans ce chapitre nous exposons les extensions 
plausibles, aux groupes r\'eductifs g\'en\'eraux et
aux formes diff\'erentielles, de la Conjecture de Selberg relative 
aux valeurs propres du laplacien op\'erant sur
les courbes modulaires.

Rappelons d'abord la Conjecture de Selberg. Soit 
$\Gamma(1)=SL(2,{\Bbb Z})$, et $\Gamma\subset \Gamma(1)$
un sous-groupe de congruence, {\it i.e.}, un sous-groupe contenant, pour 
quelque $N\geq 1$, le sous-groupe
$$
\Gamma(N)=\left\{ \left(
\begin{array}{cc}
a&b \\
c&b
\end{array}
\right) \equiv
\left( 
\begin{array}{cc}
1&0 \\ 
0&1
\end{array}
\right) (\mbox{mod } N) \right\}
$$
de $\Gamma(1)$. Alors $\Gamma$ op\`ere sur le demi-plan de Poincar\'e 
${\cal H} =\{ z\in{\Bbb C} \; : \; {\rm Im}\ z>0\}$ ;
soit $\Delta$ le laplacien invariant (positif) sur ${\cal H}$. Alors 
$L^2(\Gamma\backslash {\cal H})$ se d\'ecompose, selon la
th\'eorie spectrale, en sous-espaces propres de $\Delta$. Le spectre 
{\bf discret} de $\Delta$ est form\'e
des valeurs propres $\{0,(\lambda_n)_{n\geq 1}\}$ o\`u les 
$\lambda_n>0$ sont associ\'es aux
formes paraboliques -- cf. Iwanie\v c \cite[p. 76]{Iwaniec}. 

\begin{conj}[Selberg] \label{Selberg}
$\lambda_n\geq{1\over 4}$.
\end{conj}
\index{Conjecture de Selberg}

Quelques remarques. Il est facile de calculer le spectre (continu !) 
de $\Delta$ dans $L^2({\cal H})$~: il est \'egal
\`a $[\frac{1}{4}, +\infty[$. La conjecture est donc que le spectre 
pour les formes paraboliques est contenu dans
le spectre limite. Par ailleurs l'estim\'ee $\lambda_n\geq \frac{1}{4}$ 
est en g\'en\'eral fausse si
$\Gamma\subset \Gamma(1)$ est un sous-groupe, m\^eme arithm\'etique, 
qui n'est pas de congruence.

Rappelons la minoration connue \`a la suite des travaux de Kim, 
Shahidi et Sarnak :

\begin{thm}[\cite{KimSarnak}] \label{KS}
$\lambda_n\geq \frac{1}{4}-\left( \frac{7}{64}\right)^2$.
\end{thm}
\index{Th\'eor\`eme de Kim, Shahidi et Sarnak}

Dans la suite de ce m\'emoire nous consid\'erons un groupe simple 
$G$ d\'efini sur ${\Bbb Q}$.
Soit
$G ({\Bbb R})^0$ la composante neutre de $G ({\Bbb R})$~:  c'est 
un produit de groupes semi-simples
r\'eels. On supposera par commodit\'e que
$$
G ({\Bbb R} )^0=U\times G^{{\rm nc}}
$$
o\`u $U$ est compact et $G$ est un groupe simple r\'eel non compact. 
Soit $K$ un sous-groupe compact
maximal de $G^{{\rm nc}}$ et $X=G^{{\rm nc}} /K$. Soit $\Gamma\subset G ({\Bbb Q})$ un 
sous-groupe de congruence.
L'analogue de $\Gamma\backslash {\cal H}$ est alors le quotient $\Gamma\backslash 
X=\Gamma\backslash G^{{\rm nc}} /K$. ($\Gamma\subset
G ({\Bbb Q})$ se plonge naturellement dans $U\times G^{{\rm nc}}$, donc dans $G^{{\rm nc}}$ ; 
son image dans $G^{{\rm nc}}$ est
discr\`ete car $U$ est compact).

Soit $A_\infty^i=A_\infty^i(\Gamma\backslash X)$ l'espace des formes 
diff\'erentielles lisses de degr\'e $i$ sur
$\Gamma\backslash X$. On sait d'apr\`es le chapitre~1 que
$$
A_\infty^i(\Gamma\backslash 
X)\cong \mbox{Hom}_{K_\infty}(\Lambda^i \mathfrak{p},C^\infty(\Gamma\backslash X)).
$$

Supposons d'abord $G$ anisotrope, {\it i.e.}, $\Gamma\backslash X$ {\bf 
compact}. On dispose alors sur
$A_\infty^i$ du laplacien de Hodge (positif) $\Delta^i$ (que l'on notera $\Delta$, lorsqu'il n'y aura pas 
d'ambiguit\'e), identifi\'e 
\`a l'op\'erateur de Casimir (Ch.~I). On peut
le consid\'erer comme un op\'erateur non born\'e sur l'espace 
$A_{(2)}^i$ des formes $L^2$. Son noyau est
compos\'e des formes harmoniques. C'est essentiellement ce cas qui 
nous int\'eressera dans la partie
g\'eom\'etrique du volume.

En g\'en\'eral (si $\Gamma\backslash X$ n'est pas compact), $\Delta$ 
d\'efinit encore un op\'erateur auto-adjoint
dans $A_{(2)}^i$, dont le noyau est donn\'e par les formes 
harmoniques \cite{BorelWallach}, mais qui va poss\'eder un spectre
continu. Rappelons que l'espace des ${\cal H}_{(2)}^i$ des formes 
harmoniques est de dimension finie et {\bf a fortiori}
ferm\'e.

\begin{qu} \label{quest}
Fixons $G /{\Bbb Q}$, et $i\in[0,\mbox{dim } X]$. Si $\Gamma$ parcourt 
l'ensemble des sous-groupes de congruence de
$G ({\Bbb Q})$, existe-t-il une minoration uniforme $\varepsilon(
G,i)>0$ du spectre de $\Delta^i$
dans l'orthogonal de ${\cal H}_{(2)}^i(\Gamma\backslash X)$ ?
\end{qu}
\index{Question de la minoration spectrale}

Dans le cas cocompact, on cherche donc une minoration uniforme
\begin{eqnarray} \label{mino lambda}
\lambda\geq\varepsilon(G,i)
\end{eqnarray}
sur les valeurs propres $\neq 0$ de $\Delta^i$. En g\'en\'eral on 
veut que le spectre de $\Delta^i$ dans
$({\cal H}_{(2)}^i)^\perp$ soit contenu dans $[\varepsilon(G,i),+\infty[$.

\bigskip

\section{Le cas des fonctions}

Nous consid\'erons d'abord le cas o\`u $i=0$ ; on s'int\'eresse donc 
au spectre du laplacien sur l'espace
$L_0^2(\Gamma\backslash X)$ des fonstions d'int\'egrale nulle sur $\Gamma\backslash 
X$, $\Gamma$ \'etant un groupe de
congruence. Si $G=SL (2)$, une minoration est donn\'ee par le 
Th\'eor\`eme \ref{KS} (Selberg avait d\'ej\`a
obtenu la minoration $\lambda_n\geq\frac{3}{16}$). Noter que dans ce 
cas $\Delta$ a aussi un spectre continu, mais
l'on sait (inconditionnellement) que celui ci est contenu dans 
$[\frac{1}{4},\infty[$ : cf. \cite[p. 112]{Iwaniec}.

\medskip

Pour des groupes plus g\'en\'eraux, nous nous contenterons en 
g\'en\'eral d'obtenir des r\'esultats qualitatifs, sans
pr\'eciser la borne $\varepsilon$. Ceci suffit pour les applications 
g\'eom\'etriques. N\'eanmoins, il est important de
remarquer que les Conjectures d'Arthur sur le spectre automorphe 
permettent (convenablement interpr\'et\'ees\dots)
d'obtenir pour $G$ donn\'e les valeurs optimales (conjecturales) 
des bornes inf\'erieures $\varepsilon(G,i)$. Nous
allons les pr\'eciser, dans les chapitres suivants, pour les groupes 
associ\'es aux vari\'et\'es hyperboliques r\'eelles et complexes.

\medskip

Revenons \`a un groupe $G$ arbitraire. Pour $i=0$, on peut 
r\'epondre positivement \`a la Question \ref{quest} en toute
g\'en\'eralit\'e. Soit $\Gamma$ un groupe de congruence, et soit 
$L_0^2(\Gamma\backslash X)$ l'espace des fonctions
d'int\'egrale nulle, {\it i.e.}, l'orthogonal de ${\cal H}_{(2)}^0 (\Gamma\backslash X)={\Bbb C}$.

\begin{thm} \label{tau}
Fixons $G/{\Bbb Q}$. Il existe alors $\varepsilon=\varepsilon(G)>0$ tel que le spectre de $\Delta$ dans
$L_0^2(\Gamma\backslash X)$, pour {\bf tout} sous-groupe de congruence 
$\Gamma$ relatif \`a $G$, soit
contenu dans $[\varepsilon,+\infty[$.
\end{thm}
\index{Th\'eor\`eme de Clozel}

Nous allons esquisser la d\'emonstration \cite{Clozel}, puisque les outils de 
celle-ci seront utilis\'es plus loin pour les formes de degr\'e
sup\'erieur. Notons $L_0^2(\Gamma\backslash G)$ l'espace des fonctions $L^2$ 
d'int\'egrale nulle sur $\Gamma\backslash G$. C'est une
repr\'esentation de $G$ (Ch.~I) ; un argument \'evident montre 
qu'elle ne contient aucun sous-espace isomorphe \`a la
repr\'esentation triviale.

L'action triviale de $G$ sur les constantes correspond \`a la valeur 
propre $\lambda=0$ de $\Delta$. Nous voulons montrer
qu'il existe un voisinage fixe $(|\lambda|<\varepsilon)$ de $0$ ne 
rencontrant pas le spectre de $\Delta$ -- m\^eme pour
$\Gamma$ variable. Nous allons traduire ceci en termes de 
repr\'esentations unitaires de $G$ ; cela n\'ecessite quelques
pr\'eliminaires.

\bigskip

\section{Th\'eorie des repr\'esentations}

Soit $G$ un groupe de Lie connexe. On note $\widehat G$ l'ensemble \index{$\widehat{G}$} \index{dual unitaire}
des classes d'\'equivalence de repr\'esentations irr\'eductibles
unitaires de $G$. Rappelons que $\widehat G$ est muni d'une topologie 
naturelle, la topologie de Fell. Soit $\pi \in\widehat
G$ une repr\'esentation d'espace ${\cal H}$. Si $v\in{\cal H}$, le {\bf 
coefficient} associ\'e est la fonction sur $G$~:
\begin{eqnarray} \label{coeff}
c_v(g)=(gv,v) \; \; (g\in G).
\end{eqnarray}
\index{topologie de Fell}
\index{coefficient}

Soit $\Omega\subset G$ un sous-ensemble compact, $v_1,\ldots 
v_d\in{\cal H}$, $\varepsilon>0$. On note
$V(\Omega,v_i,\varepsilon)$ l'ensemble des $\rho\in\widehat G$ telles 
qu'il existe $w_1,\ldots w_d\in{\cal H}_\rho$ v\'erifiant
\begin{eqnarray} \label{fell}
|c_{w_i}(g)-c_{v_i}(g)|<\varepsilon   \; \; (g\in\Omega).
\end{eqnarray}

On obtient ainsi une base de voisinages de $\pi $ ; celle-ci 
d\'efinit la topologie sur $\widehat G$. Cette topologie n'est pas
s\'epar\'ee~: par exemple si $G=SL(2,{\Bbb R})$ et si $s\in[0,\frac{1}{2}[$
et $s\rightarrow \frac{1}{2}$, la limite de la repr\'esentation
irr\'eductible
$$
\pi _s=\mbox{ind}_B^G(|x|^s )
$$
o\`u $B=\left\{
\left(
\begin{array}{cc}
x&* \\
0&x^{-1}
\end{array}
\right)
\right\}\subset G$, et l'induction est l'induction unitaire, est 
form\'ee de la repr\'esentation triviale ${\Bbb C}$ et des deux
s\'eries discr\`etes $\delta_2$, $\delta_{-2}$ (\cite{Lang}). Elle est 
n\'eanmoins s\'eparable (tout point a une base
d\'enombrable de voisinages) de sorte que nous pourrons nous limiter, 
dans les arguments de convergence, \`a consid\'erer
des suites.

Par ailleurs, on peut d\'ecrire la topologie de $\widehat G$ \`a 
l'aide d'un seul coefficient~:

\begin{prop} \label{cv}
Soit $(\pi _n)_{n\geq 1}$, $\pi \in\widehat G$. Alors $\pi _n\rightarrow \pi $ 
si et seulement si un coefficient non nul $c_u$ de $\pi
$ $(u\in {\cal H}_\pi )$ est limite uniforme sur tout compact de 
coefficients $c_{u_n}$ $(u_n\in {\cal H}_{\pi _n})$.
\end{prop}
\noindent \index{topologie de Fell}
Voir Dixmier \cite[Prop. 18.1.5]{Dixmier}.

Dans le reste de ce paragraphe on suppose $G$ simple et non compact. 
Soit $K$ un sous-groupe compact maximal. Soit
$\widehat G_s$ \index{dual sph\'erique} \index{$\widehat G_s$} le sous-ensemble de $G$ form\'e des repr\'esentations 
sph\'eriques, \index{repr\'esentation sph\'erique} {\it i.e.}, ayant un vecteur $\neq 0$ fix\'e
par $K$. On sait que celui-ci est alors unique \`a un scalaire pr\`es \cite{Knapp}.

Le lemme suivant r\'esulte trivialement de la d\'efinition du dual~:

\begin{lem} \label{ouvert}
$\widehat G_s$ est ouvert dans $\widehat G$.
\end{lem}

Noter que le contre-exemple ci-dessus, montrant que $\widehat G$ 
n'est pas s\'epar\'e, montre aussi que $\widehat G_s$
n'est pas ferm\'e.

Puisque $\widehat G$ est muni d'une topologie, il est muni d'une 
structure bor\'elienne. On sait que $\widehat G$ est
bor\'elien standard, c'est-\`a-dire isomorphe \`a $[0,1]$ muni de la 
structure bor\'elienne usuelle \cite[4.6.1]{Dixmier}.

Nous devrons utiliser la d\'ecomposition selon $\widehat G$ de 
repr\'esentations unitaires {\bf r\'eductibles} de $G$.
Celle-ci utilise la notion d'int\'egrale hilbertienne, pour laquelle 
on renvoie le lecteur \`a \cite[A 69]{Dixmier} et \`a \cite{Zimmer}. Dans notre cas, $G$ est
``de type I'' au sens de la th\'eorie des repr\'esentations \cite{Dixmier} et 
nous n'avons besoin que d'un cas tr\`es simple.

Toutes les repr\'esentations de $G$, sauf la repr\'esentation 
triviale, se r\'ealisent dans un espace fixe ${\cal H}$ de dimension
d\'enombrable ; les op\'erateurs $\pi (g)$ $(\pi \in\widehat G)$ sont 
alors fonctions mesurables de $\pi $ pour tout $g$.
Posons, pour $\pi =1_G$ (la repr\'esentation triviale) ${\cal H}_{\pi} 
={\Bbb C}$, et ${\cal H}_{\pi} ={\cal H}$ pour les autres repr\'esentations
de $\widehat G$.

Si $\mu$ est une mesure positive sur $\widehat G$, la repr\'esentation
$$
\rho=\int_{\widehat G} {\cal H}_\pi \ d\mu(\pi )
$$
(sur l'espace des fonctions mesurables $\varphi:\widehat 
G\rightarrow {\cal H} \oplus {\Bbb C}$ telles que $\varphi(\pi )\in {\cal H}_\pi $ et que
$\int ||\varphi(\pi ) ||^2 d\mu(\pi )<\infty$) se 
d\'efinit de fa\c con \'evidente. On d\'efinit de m\^eme
\begin{eqnarray} \label{rho}
\rho=\int_{\widehat G}\widehat \bigoplus^{m(\pi )} {\cal H}_\pi \ d\mu(\pi )
\end{eqnarray}
o\`u $\pi \mapsto m(\pi )\in \{1,2,\ldots ,\infty\}$ est une fonction bor\'elienne.

\begin{thm} \label{decomposition}
Soit $\rho$ une repr\'esentation unitaire de $G$ sur un espace de 
Hilbert ${\cal H}$.
\begin{description}
\item[{\rm (i)}] Il existe une fonction bor\'elienne $\pi \mapsto 
m(\pi )$, $\widehat G\rightarrow \{1,2,\ldots ,\infty\}$ et une mesure
positive $\mu$ sur $\widehat G$ telles que
\begin{eqnarray} \label{rho2}
\rho\cong\int_{\widehat G}\widehat \bigoplus^{m(\pi )} {\cal H}_\pi \ d\mu(\pi ).
\end{eqnarray}
\item[{\rm (ii)}] Si $\rho$ admet deux repr\'esentations (\ref{rho2}), 
associ\'ees \`a $(m,\mu)$ et $(m',\mu')$, les mesures $\mu$ et
$\mu'$ sont \'equivalentes et les fonctions $m$ et $m'$ \'egales p.p. 
(pour $\mu$ ou $\mu'$).
\end{description}
\end{thm}
\index{d\'ecomposition en irr\'eductibles, Th\'eor\`eme de}

Pour les d\'emonstrations, on renvoie \`a Dixmier \cite[8.6.5, 8.6.6]{Dixmier}.

Si $(\rho,{\cal H})$ est donn\'ee et d\'ecompos\'ee selon (\ref{rho2}), le {\bf 
support} de $\rho$ est le support de $\mu$, contenu dans \index{support (d'une repr\'esentation)}
$\widehat G$. Si $\rho,\rho'$ sont deux repr\'esentations, on dit que 
$\rho$ est faiblement contenue dans $\rho'$ si 
Supp$(\rho)\subset$Supp$(\rho') $. (En particulier si $\pi 
\in$Supp$(\rho)$ est irr\'eductible, $\pi $ est faiblement
contenue dans $\rho$.) Pour les applications, il est important 
d'avoir une caract\'erisation intrins\`eque de cette relation.

\begin{thm} \label{faiblt contenue}
Soient $\rho,\rho'$ deux repr\'esentations unitaires de $G$. Les deux 
propri\'et\'es suivantes sont \'equivalentes~:
\begin{description}
\item[{\rm (i)}] $\rho '$ est faiblement contenue dans $\rho$
\item[{\rm (ii)}] Tout coefficient $c$ de $\rho'$ est limite, 
uniform\'ement sur tout compact, d'une suite de combinaisons
lin\'eaires positives de coefficients de $\rho$.
\end{description}
\end{thm}
\index{faiblement contenue}

Nous aurons besoin de deux r\'esultats suppl\'ementaires.

Le premier est une propri\'et\'e de continuit\'e de la restriction. 
Soit $H\subset G$ un sous-groupe semi-simple.

\begin{lem} \label{res} 
Soit $\pi _n\in\widehat G$ et supposons que $\pi _n\rightarrow 1 _G$. Soit 
$S_n\subset \widehat H$ le support de $\pi _{n_{|_H}}$.
Alors
$1_H$ appartient \`a l'adh\'erence de $\displaystyle \sqcup_nS_n$.
\end{lem}

Ceci se d\'eduit ais\'ement de \cite[3.4.2]{Dixmier}.

Le second est un th\'eor\`eme de Howe et Moore \cite{HoweMoore}.

\begin{thm} \label{howe-moore}
Supposons $G$ simple, non compact et connexe. Si $\pi $ est une repr\'esentation 
irr\'eductible non triviale de $G$, les coefficients de $\pi $
tendent vers $0$ \`a l'infini.
\end{thm}
\index{Th\'eor\`eme de Howe-Moore}

\bigskip

\markboth{CHAPITRE 2. SPECTRE DU LAPLACIEN}{2.3. PRINCIPE DE RESTRICTION}

\section{Principe de restriction et d\'emonstration du Th\'eor\`eme \ref{tau}}

La d\'emonstration du Th\'eor\`eme \ref{tau} repose sur une m\'ethode 
introduite par Burger et Sarnak qui permet de d\'emontrer,
pour un groupe $G$, une propri\'et\'e telle que le Th\'eor\`eme 
\ref{tau} par r\'eduction \`a un sous-groupe plus petit pour
lequel le Th\'eor\`eme est d\'ej\`a connu.

Comme dans tout ce chapitre $G$ est un groupe d\'efini sur ${\Bbb Q}$, 
nous le supposerons de plus simple (comme ${\Bbb Q}$-groupe) dans
ce paragraphe. Soit $\widehat G$ le dual du groupe r\'eel $G(=G ({\Bbb R}))$, 
au sens du \S 2.2. On suppose que $G $ n'est pas
compact. Si $\Gamma\subset G({\Bbb Q})$ est un sous-groupe de 
congruence, on peut consid\'erer dans $\widehat G$ le
support de la repr\'esentation $L^2(\Gamma\backslash G)$. On note $\widehat 
G_{\rm Aut}$ \index{dual automorphe} \index{$\widehat G_{\rm Aut}$} la cl\^oture dans $\widehat G$ de la r\'eunion des
supports de $L^2(\Gamma\backslash G )$ quand $\Gamma$ parcourt tous les 
sous-groupes de congruence de $G({\Bbb Q})$. Noter que
$\widehat G_{\rm Aut}$ {\bf d\'epend donc de la ${\Bbb Q}$-forme} de $G$. 

Soit $H\subset G$ un sous-groupe, d\'efini sur ${\Bbb Q}$, et 
semi-simple. Les m\^emes notions s'appliquent alors \`a
$H$.

Par ailleurs, si $\pi $ est une repr\'esentation irr\'eductible de 
$G$, on peut la restreindre \`a $H$. Selon le
Th\'eor\`eme \ref{tau}, le {\bf support} de $\pi _{|_{H}}$ dans $\widehat H$ 
est bien d\'efini. Burger et Sarnak \cite{BurgerSarnak} d\'emontrent le
th\'eor\`eme suivant~:

\begin{thm}[Burger, Sarnak]
Soient $\pi \in\widehat G_{\rm Aut}$ et $\tau\in\widehat H$. Si 
$\tau$ appartient au support de $\pi _{|_{H}}$, $\tau\in\widehat
H_{\rm Aut}$.
\end{thm}

Nous renvoyons \`a \cite{BurgerSarnak} pour la d\'emonstration. Une variante de 
celle-ci (permettant de contr\^oler le laplacien sur les formes
diff\'erentielles plut\^ot que sur les fonctions) sera expliqu\'ee 
dans le Chapitre 9.

\bigskip

Pour d\'emontrer le Th\'eor\`eme \ref{tau}, nous commen\c cons par le 
reformuler en termes de th\'eorie des repr\'esentations.

\begin{thm} \label{tau2}
(Hypoth\`eses du Th\'eor\`eme \ref{tau})
Il existe un voisinage ${\cal V}$ de la repr\'esentation triviale 
dans $\widehat G$ tel que, pour tout sous-groupe de congruence
$\Gamma$, le support de $L_0^2(\Gamma\backslash X )$ soit disjoint de ${\cal V}$.
\end{thm}

Puisque la repr\'esentation triviale $1_G$ de $G$ n'appara\^{\i}t 
pas dans $L_0^2(\Gamma\backslash X)$, ceci \'equivaut par
d\'efinition de $\widehat G_{\rm Aut}$ \`a :
\begin{eqnarray} \label{isolee}
1_G\ \hbox{\bf est\ isol\'ee\ dans}\ \widehat G_{\rm Aut}.
\end{eqnarray}

Les assertions des Th\'eor\`emes \ref{tau} et \ref{tau2} sont \'equivalentes. 
V\'erifions-le en supposant pour simplifier que $G$ est
anisotrope. D'apr\`es le Lemme \ref{ouvert}, le dual sph\'erique $\widehat 
G_s$ est ouvert dans $\widehat G$. Consid\'erons, dans
$L_0^2(\Gamma\backslash G)=\widehat \bigoplus_\pi {\cal H}_\pi  $ (somme directe 
hilbertienne de repr\'esentations irr\'eductibles), la
partie sph\'erique $L_{0,s}^2=\widehat \bigoplus_{\pi \in\widehat 
G_s} {\cal H}_\pi$. Le Th\'eor\`eme \ref{tau2} est donc \'equivalent au
fait que les supports de $L_{0,s}^2$ (pour $\Gamma$ variable) soient 
s\'epar\'es de la repr\'esentation triviale.

Soit $P=MAN$ un sous-groupe parabolique minimal de $G$. Une 
repr\'esentation sph\'erique de $G$ est isomorphe \`a
l'unique sous-quotient sph\'erique $\pi $ d'une induite
$$
\rho=\mbox{ind}_P^{G} (1\otimes e^\nu \otimes 1)
$$
o\`u $\nu$ appartient au dual complexe de $\mathfrak{a}=$Lie$(A)$. Si $\pi $ 
est unitaire, $\nu$ doit \^etre {\bf hermitien}, \index{hermitien} {\it i.e}~:
$\overline \nu=-w\nu$ pour un $w$ dans le groupe de Weyl $W(G,A)$.

Si $\pi $ est sph\'erique et appartient \`a $L_0^2(\Gamma\backslash G)$, le 
sous-espace $\pi ^K$ est un sous-espace de
$L_0^2(\Gamma\backslash X)$ et la valeur propre associ\'ee $\lambda$ de 
$\omega$ est \'egale \`a celle de l'op\'erateur de Casimir,
qui op\`ere par $(\delta,\delta)-(\nu,\nu)$, $\delta\in \mathfrak{a}^*$ 
\'etant la demi-somme des racines de $A$ dans $N$ (cf. \cite[Prop. 8.2.2]{Knapp}).

Ecrivons $\nu=\nu_r+\sqrt{-1}\nu_i$, avec $\nu_r,\nu_i\in \mathfrak{a}^*$. Alors
\begin{eqnarray} \label{nur}
\nu_r=-w\nu_r\ ,\ \nu_i=w\nu_i
\end{eqnarray}
donc $\nu_r$ et $\nu_i$ sont orthogonaux et
\begin{eqnarray} \label{sf}
(\delta,\delta)-(\nu,\nu)=(\delta,\delta)+(\nu_i,\nu_i)-(\nu_r,\nu_r).
\end{eqnarray}

A conjugaison pr\`es par $W$, on peut supposer $\nu_r$ contenu dans 
la chambre de Weyl aigu\"e et ferm\'ee contenant
$\delta$. Si $\pi $ est unitaire on sait alors que $\delta-\nu_r$ est 
dans la chambre de Weyl ferm\'ee obtuse associ\'ee \cite[IV,
5.2]{BorelWallach}.

Par ailleurs la topologie de $\widehat G_s$ d\'eduite de celle de 
$\widehat G$ co\"{\i}ncide avec la restriction aux param\`etres
$\nu$ hermitiens de la topologie de ${\mathfrak{a}}^*\otimes {\Bbb C}$. D'apr\`es 
le Th\'eor\`eme \ref{tau2}, on a donc
$\|\delta-\nu\|=\|\delta-\nu_r\|+\|\nu_i\|\geq\varepsilon_1$ o\`u 
$\varepsilon_1>0$ est d\'etermin\'e par $G$. Pour
$\nu_r$ proche de $\delta$, $\nu_r$ est r\'egulier et (\ref{nur}) implique 
$w=-1$ soit $\nu_i=0$. Il existe donc $\varepsilon_2>0$
tel que $\|\delta-\nu_r\|\geq\varepsilon_2$. Enfin, les 
propri\'et\'es de convexit\'e des chambres de Weyl impliquent alors 
que
$$
(\nu_r,\nu_r)\leq(\delta,\delta)-\varepsilon
$$
o\`u $\varepsilon$ est d\'etermin\'e par $\varepsilon_2$. On en 
d\'eduit d'apr\`es (\ref{sf})~:
\begin{eqnarray} \label{label}
\lambda=(\delta,\delta)-(\nu,\nu)\geq\varepsilon\ .
\end{eqnarray}

R\'eciproquement, (\ref{label}) implique que la repr\'esentation triviale est 
isol\'ee de $L_{0,s}^2$ car $\lambda=0$ pour la
repr\'esentation triviale.

Enfin, si $G$ n'est pas anisotrope, un argument similaire 
s'applique \`a la partie continue de $L_0^2(\Gamma\backslash G)$.

La d\'emonstration du Th\'eor\`eme \ref{tau2} va reposer sur la m\'ethode 
de Burger-Sarnak. Nous utilisons la forme (\ref{isolee}) du
Th\'eor\`eme. Supposons donn\'e un sous-groupe simple $H$ de $
G$ tel que $H( = H({\Bbb R}))$ soit non compact et que le
Th\'eor\`eme soit vrai pour $H$. Si celui-ci est faux pour 
$G$, on peut trouver une suite $\pi _n\in\widehat
G_{\rm Aut}$ $(\pi _n\not\cong 1_G)$ telle que $\pi _n\rightarrow 1_G$. 
D'apr\`es le Lemme \ref{res}, $1_H$ appartient \`a l'adh\'erence
de $\displaystyle \bigcup_n{\rm Supp}(\pi _n{_{|_H}})$. Si $1_H$ est 
isol\'ee dans $\widehat H_{\rm Aut}$, le Th\'eor\`eme
\ref{decomposition} implique que $1_H$ est contenue discr\`etement dans $\pi 
_{n_{|_H}}$ pour $n>>0$. Ceci est impossible d'apr\`es le
Th\'eor\`eme \ref{howe-moore} si $H({\Bbb R})$ n'est pas compact.

Nous devons enfin trouver -- dans tout groupe ${\Bbb Q}$-simple $G$ -- 
un sous-groupe $H$ pour lequel le Th\'eor\`eme
soit d\'ej\`a connu (et $H({\Bbb R})$ non compact).

Nous renvoyons le lecteur \`a \cite{Clozel}, puisque cette partie de l'argument 
n'a rien \`a voir avec les questions abord\'ees dans ce
volume.

\medskip

Le Th\'eor\`eme \ref{tau} a la cons\'equence suivante. Soit $G$ 
anisotrope sur ${\Bbb Q}$ tel que $G({\Bbb R})^0=G^{{\rm nc}} \times U$ (comme au d\'ebut de ce chapitre) 
avec $G^{{\rm nc}}$ simple non compact et $U$ compact. Avec la terminologie 
de l'Introduction~:

\begin{thm} \label{forme fermee}
La Conjecture $A_{d=0}^-(1)$ est v\'erifi\'ee pour $G$.
\end{thm}

Ecrivons en effet, $\Gamma$ \'etant un sous-groupe de congruence~:
$$
L^2(\Gamma\backslash G^{{\rm nc}})={\Bbb C} \oplus \widehat{ 
\bigoplus}_{\pi} \pi ={\Bbb C} \oplus L_0^2(\Gamma\backslash G^{{\rm nc}})
$$
o\`u les repr\'esentations $\pi $ sont irr\'eductibles non triviales. 
Soit $V$ l'espace des vecteurs $K$-finis d'une repr\'esentation
$\pi $ et consid\'erons le complexe calculant la $(\mathfrak{g},K)$-cohomologie de $\pi $ (Ch.~I)~:
\begin{eqnarray} \label{complexe}
\cdots\rightarrow C^i(\pi )=\mbox{Hom}_K (\Lambda^i \mathfrak{p},V)\rightarrow C^{i+1}(\pi )\rightarrow \cdots.
\end{eqnarray}

Si le rang de $G^{{\rm nc}}$ est $>1$, on sait d'apr\`es un th\'eor\`eme de 
Zuckerman et Borel-Wallach \cite[Cor. V.3.4]{BorelWallach}, que
$H^1(\mathfrak{g},K,V)=\{0\}$ si $\pi $ n'est pas triviale. Il en est de 
m\^eme si $\pi $ est triviale car $\mbox{Hom}_K(\mathfrak{p} ,{\Bbb C})=\{0\}$, et ceci
pour tout $G^{{\rm nc}}$.

La suite (\ref{complexe}) est donc exacte en $i=1$ ; puisque le laplacien 
op\`ere par l'op\'erateur de Casimir, dont les valeurs propres sur
$A^0=L^2(\Gamma\backslash G^{{\rm nc }}/K)$ sont minor\'ees hors les constantes, on en 
d\'eduit que les valeurs propres de $\Delta$ sur
$A_{d=0}^1$ sont minor\'ees.

Si $G^{{\rm nc}}$ est de rang~1 (donc localement isomorphe \`a $SU(n,1)$ ou 
$SO(n,1)$), ceci reste vrai en dehors des repr\'esentations
$\pi $ telles que $H^1(\mathfrak{g},K,V)\neq 0$. Mais celles-ci 
repr\'esentent exactement les formes harmoniques, pour lesquelles
$\Delta \alpha=0$.

\bigskip

\markboth{CHAPITRE 2. SPECTRE DU LAPLACIEN}{2.4. REPR\'ESENTATIONS NON ISOL\'EES}

\section{Repr\'esentations non isol\'ees et contre-exemples \`a $A^-$}

Les derniers paragraphes de ce chapitre sont destin\'es \`a mettre en 
valeur le r\^ole particulier des groupes $SU(n,1)$ et
$SO(n,1)$ relativement \`a ces questions. Nous voulons tout d'abord 
expliquer que la Conjecture $A^-(i)$ ne peut \^etre vraie
en g\'en\'eral.

Consid\'erons un groupe $G$ anisotrope tel que $G^{{\rm nc}}$ 
est simple. Soit ${\cal H}^i (\Gamma\backslash X)$ l'espace des
$i$-formes harmoniques. D'apr\`es le Chapitre I, ${\cal H}^i (\Gamma\backslash X)$ 
se d\'ecompose en sous-espaces relatifs aux
repr\'esentations irr\'eductibles $\pi $ de $G$ telles que~:
\begin{description}
\item[(i)] $\pi \subset L^2(\Gamma\backslash G)$
\item[(ii)] Hom$_K (\Lambda^i\mathfrak{p},V(\pi ))\not =0$
\item[(iii)] L'op\'erateur de Casimir op\`ere trivialement sur $V(\pi )$.
\end{description}
Rappelons (Ch.~I) que (ii) et (iii) sont \'equivalents \`a
\begin{description}
\item[(iv)] $H^i(\mathfrak{g},K;V(\pi ))\not =0$.
\end{description}

Soit $\pi \in\widehat G$ une repr\'esentation contenue dans 
$L^2(\Gamma\backslash G)$, et v\'erifiant (iv). Supposons de plus que
$\pi $ n'est pas isol\'ee dans $\widehat G$. Soit donc $\pi 
_n\in\widehat G$ une suite tendant vers $\pi $. La condition (ii) est
alors v\'erifi\'ee par $\pi _n$ pour $n>>0$, par un analogue 
\'evident du Lemme \ref{ouvert}. Par ailleurs, si $\rho\in\widehat G$,
l'op\'erateur de Casimir $C$ op\`ere par un scalaire $\rho(C)$ sur 
les vecteurs $C^\infty$ de $\rho$ (\S 1.1). On v\'erifie alors
que $\rho(C)$ est une fonction continue de $\rho$ pour la topologie 
de $\widehat G$ \cite{BD}. Pour $\pi _n\rightarrow \pi $ on a donc
$\pi _n(C)\rightarrow \pi (C)=0$ ; pour $n>>0$, $\pi _n(C)$ est non nul car 
l'ensemble des repr\'esentations unitaires $\pi$ v\'erifiant
(ii)-(iii) $\Longleftrightarrow$ (iv) est fini. On a donc~:

\begin{thm} \label{suite approx}
Soient $\Gamma\subset G$ un sous-groupe de congruence, $\pi $ une 
repr\'esenta\-tion de $G$ v\'erifiant (i--iii) et $\pi _n\rightarrow \pi
$ dans $\widehat G$. Si les repr\'esentations $\pi _n$ appartiennent 
\`a $\widehat G_{\rm Aut}$, la Conjecture $A^-(i)$ est en
d\'efaut pour les degr\'es $i$ tels que $H^i(\mathfrak{g},K;V(\pi ))\not =0$.
\end{thm}

Pour un groupe donn\'e, l'\'etude de $A^-(i)$ se divise donc en deux 
questions~:

\begin{prob} \label{pbm1}
D\'ecrire les repr\'esentations isol\'ees dans $\widehat G$ parmi les 
repr\'esenta\-tions cohomologiques.
\end{prob}

\begin{prob} \label{pbm2}
Pour $G/{\Bbb Q}$ donn\'e, soit $\pi \in \widehat G$ une 
repr\'esentation cohomologique non triviale et non isol\'ee. La
repr\'esentation $\pi
$ est-elle isol\'ee {\bf dans} $\widehat G_{\rm Aut} \cup \{ \pi \}$ ?
\end{prob}

Il existe une classe simple de repr\'esentations pour lesquelles les 
deux questions ont une solution n\'egative.

\bigskip

Rappelons qu'une repr\'esentation irr\'eductible de $G$ appartient 
\`a la s\'erie discr\`ete \index{s\'erie discr\`ete} si ses coefficients (\S 2.3) sont de
carr\'e int\'egrable. D'apr\`es Harish-Chandra, $G$ a une s\'erie 
discr\`ete si, et seulement si, $G$ contient un sous-groupe de
Cartan compact \cite{Knapp}. Parmi nos groupes ``hyperboliques'', $SU(n,1)$ a 
toujours une s\'erie discr\`ete et $SO(n,1)$ a une s\'erie
discr\`ete si et seulement si $n$ est {\bf pair} -- ainsi, 
$SO(2,1)\approx SL(2,{\Bbb R})$ a une s\'erie discr\`ete alors que
$SO(3,1)\approx SL(2,{\Bbb C})$ n'en a pas.

\bigskip

Rappelons qu'un sous-groupe parabolique $P=M^0AN$ de $G$ est {\bf 
cuspidal} \index{cuspidal, sous-groupe parabolique} si $M^0$ contient un sous-groupe de Cartan
compact (donc, poss\`ede une s\'erie discr\`ete). Une 
repr\'esentation $\pi $ appartient \`a la s\'erie discr\`ete si, et 
seulement
si, elle appara\^{\i}t comme sous-module {\bf discret} dans $L^2(G)$. 
Plus g\'en\'eralement, $\pi \in\widehat G$ est
temp\'er\'ee \index{temp\'er\'ee} si elle v\'erifie les conditions \'equivalentes qui suivent~:

\begin{eqnarray} \label{ds support}
\pi  \mbox{ {\it appartient au  support de} } L^2(G)
\end{eqnarray}
\index{support}

\begin{eqnarray} \label{ss module}
\begin{array}{c}
\pi \mbox{ {\it est un sous-module de} } {\rm 
ind}_{P=M^0AN}^G(\delta\otimes e^\nu \otimes 1)
\mbox{ {\it o\`u} } P
\mbox{ {\it est cuspidal,} }  \index{cuspidal} \\
\delta\in\widehat {M^0} \mbox{ {\it appartient \`a la s\'erie 
discr\`ete, et} \index{s\'erie discr\`ete} } \nu\in i\mathfrak{a}^* \mbox{ {\it est unitaire.}}
\end{array}
\end{eqnarray}

Supposons que $G$ n'a pas de s\'erie discr\`ete, et soit $P_f\subset 
G$ un parabolique tel que $A$ soit de dimension minimale.
Alors $P_f$ est n\'ecessairement cuspidal ; si $P_f=M_f\ N_f$, $M_f$ 
est uniquement d\'etermin\'e \`a conjugaison pr\`es dans
$G$. On dit que $P_f$ est un parabolique fondamental.

\begin{thm}[Borel-Wallach] \label{borelwallach}
Soit $\pi \in\widehat G$ une repr\'esentation temp\'er\'ee et 
supposons que $H^\bullet (\mathfrak{g},K;V(\pi ))\not =0$. Alors
\begin{description}
\item[{\rm (i)}] $\pi $ est un sous-module de ${\rm 
ind}_{P_f}^G(\delta\otimes e^0\otimes 1)$, $\delta\in\widehat {M^0}$
\'etant une repr\'esentation de la s\'erie discr\`ete.
\item[{\rm (ii)}] La cohomologie de $\pi $ est concentr\'ee dans 
l'intervalle $]q-e,q+e]$ o\`u $q=\frac12 \dim_{{\Bbb R}}(G/K)$ et
$e=\frac12(rg\ G-rg\ K)$.
\end{description}
\end{thm}
\index{Th\'eor\`eme de Borel-Wallach}

En fait la repr\'esentation induite de (i) n'est cohomologique que 
pour un nombre fini de repr\'esentations $\delta$,
explicitement d\'ecrites \cite[Thm. 5.1]{BorelWallach}. On v\'erifie que $q\pm e$ est 
entier \cite[p. 98]{BorelWallach}; tous les degr\'es contenus dans $]q-e,q+e]$
apparaissent dans $H^\bullet (\mathfrak{g},K;V)$ o\`u $V$ est l'espace des 
vecteurs $K$-finis de l'induite totale donn\'ee par (i). Noter
que si $G$ a une s\'erie discr\`ete, le Th\'eor\`eme reste vrai avec 
$P_f=G$ -- les seules repr\'esentations cohomologiques
temp\'er\'ees sont dans la s\'erie discr\`ete.

Nous allons d\'eduire de ce th\'eor\`eme~:

\begin{thm} \label{contre ex}
Soit $G$ arbitraire, et supposons que $G$ n'a pas de s\'erie 
discr\`ete. Alors la propri\'et\'e $A(i)$ est en d\'efaut pour
$i\in ]q-e,q+e]$.
\end{thm}
\index{repr\'esentations cohomologiques non isol\'ees}

Soit en effet $\pi $ une repr\'esentation v\'erifiant les conditions 
du Th\'eor\`eme \ref{borelwallach} (i) et telle que $H^i(\mathfrak{g},K,V(\pi ))\not
=0$. Alors $\pi $ est temp\'er\'ee {\bf et appartient donc au ``dual 
automorphe''} $\widehat G_{\rm Aut}$. Si $\pi _n\rightarrow \pi $
et $\pi _n$ appara\^{\i}t dans la d\'ecomposition (continue si $G$ 
est isotrope) de $L^2(\Gamma_n\backslash G)$ pour un groupe de
congruence, la valeur propre associ\'ee $\lambda_n$ de $C$ v\'erifie 
$\lambda_n\rightarrow 0$, et appara\^{\i}t dans le spectre
(peut-\^etre continu) du laplacien $\omega_i$.

La d\'emonstration a utilis\'e le th\'eor\`eme suivant~:

\begin{thm} \label{luck}
Toute repr\'esentation temp\'er\'ee $\pi \in\widehat G$ est limite de 
repr\'esentations $\pi _n$ apparaissant (discr\`etement) dans
$L^2(\Gamma_n\backslash G)$.
\end{thm}
\index{Th\'eor\`eme de Delorme}

Il y a plusieurs d\'emonstrations du Th\'eor\`eme \ref{luck}. Supposons 
d'abord $G$ anisotrope, donc les quotients $\Gamma\backslash G$
compacts. Si $\Gamma_n\subset \Gamma_0$ est une famille de 
sous-groupes distingu\'es (en fait $\Gamma_{n+1}$
distingu\'e dans $\Gamma_n$) d'intersection $\{1\}$, le th\'eor\`eme 
est d\^u \`a Delorme \cite{Delorme}; on sait m\^eme que la suite des
mesures spectrales (sur $\widehat G$) des $L^2(\Gamma_n\backslash G)$, 
pond\'er\'ees par $[\Gamma_0:\Gamma_n]$, tend vers celle de
$L^2(G)$.

\medskip
Une autre approche du th\'eor\`eme est due \`a Burger, Li et Sarnak 
\cite{BurgerLiSarnak}. Cet article ne donne pas la d\'emonstration, qui est expos\'ee
dans \cite{CL}. Pour $G$ isotrope, elle implique que $\pi $ est limite 
de repr\'esentations dans le support de $L^2(\Gamma_n\backslash G)$.

Enfin, le fait que $\pi $ est limite de repr\'esentations $\pi _n$ 
apparaissant {\bf discr\`etement} dans $L^2(\Gamma_n\backslash G)$
r\'esulte de la d\'emonstration donn\'ee dans \cite{CL}. Nous expliquons 
l'argument suppl\'ementaire au lecteur familier avec les m\'ethodes
ad\'eliques, en supposant $G$ simplement connexe comme dans \cite{CL}~: 
choisissons une place finie $p$, et rempla\c cons $G=
G({\Bbb R})$ par $G_{\infty,p}=G({\Bbb R})\times G ({\Bbb Q}_p)$. Alors \cite[\S 3]{CL}, le Th\'eor\`eme \ref{luck} reste vrai, en rempla\c cant $G$ par
$G_{\infty,p}$ et les $\Gamma_ n$ par des sous-groupes 
$S$-arithm\'etiques pour $S=\{\infty,p\}$. Si on remplace $\pi 
\in\widehat
G$ par $\pi \otimes \pi _p$ o\`u $\pi _p\in\widehat G ({\Bbb Q}_p)$ 
est supercuspidale, on voit que $\pi \otimes \pi _p$ est limite
de repr\'esentations ``automorphes'' $\pi ^{(n)}\otimes \pi 
_p^{(n)}$. Mais alors $\pi _p^{(n)}=\pi _p$ pour tout $n>>0$, et cette
repr\'esentation ne peut appara\^{\i}tre que dans le spectre discret.

\medskip
{\bf Exemple} \ref{luck}. Si $G=SO(n,1)$ avec $n$ impair $= 2m+1$, on a 
$q=\frac{n}{2}=m+\frac12$, $e=\frac12$, et la
cohomologie de $\pi $ appara\^{\i}t en degr\'es $\{m,m+1\}$.

\medskip
Terminons en remarquant que le Th\'eor\`eme \ref{contre ex} est reli\'e, par le 
th\'eor\`eme de L\"uck \cite{Luck}, au calcul des invariants de
Novikov-Shubin de l'espace $G/K$ (Lohou\'e et Mehdi \cite{MehdiLohoue}).

\markboth{CHAPITRE 2. SPECTRE DU LAPLACIEN}{2.5. PERSPECTIVES}

\section{Perspectives : contraintes locales et contraintes automorphes}

Nous revenons aux deux Probl\`emes (\ref{pbm1} et \ref{pbm2}) \'enonc\'es dans la
section pr\'ec\'edente et nous allons d\'ecrire, en utilisant toute la force de
la th\'eorie des repr\'esentations, deux cas oppos\'es o\`u ils 
peuvent \^etre r\'esolus, amenant \`a une solution satisfaisante du
probl\`eme $A^-(i)$.

Le Probl\`eme \ref{pbm1} a en fait \'et\'e compl\`etement r\'esolu par 
Vogan, dans un article longtemps clandestin \cite{Vogan2}. Sa solution est
difficile, et nous n'exposerons les r\'esultats que dans les cas qui 
nous int\'eressent. 
Le Probl\`eme \ref{pbm2} est encore plus profond~: l'un
des buts de ces notes est d'expliquer comment, pour les groupes associ\'es aux espaces hyperboliques, sa solution est implicitement
donn\'ee par les Conjectures d'Arthur.

Commen\c cons par le groupe $G=SU(n,1)$. Le sous-groupe compact 
maximal est $K\cong U(n)$. La d\'ecomposition de Cartan est
$\mathfrak{g_0}=\mathfrak{k}_0 \oplus \mathfrak{p}_0 $ avec $\mathfrak{p}_0 \cong {\Bbb C}^n$. On a $\mathfrak{p} = \mathfrak{p}_0 \otimes 
{\Bbb C} \cong {\Bbb C}^n \oplus {\Bbb C}^n$, le premier facteur s'identifiant \`a
l'espace tangent holomorphe en $o=1\cdot K$ \`a $X=G/K$ et le second 
\`a l'espace tangent antiholomorphe. Donc~:
\begin{eqnarray} \label{dec de Lp}
\Lambda^i(\mathfrak{p}) \cong \bigoplus_{p+q=i} \Lambda^p ({\Bbb C}^n)\otimes\Lambda^q ({\Bbb C}^n) .
\end{eqnarray}

Le groupe $K$ op\`ere sur $\mathfrak{p}_0$ par $u\cdot X=\det(u)^{-1}uX$ ($u\in 
K,X\in {\Bbb C}^n$). Il pr\'eserve la d\'ecomposition (\ref{dec de Lp}).

On sait que les repr\'esentations $\Lambda^p({\Bbb C}^n)$ sont 
irr\'eductibles sous $K$ ; chaque facteur de (\ref{dec de Lp}) contient alors une
unique repr\'esentation irr\'eductible contenant le vecteur 
$e_p\otimes e_q$, $e_p$ et $e_q$ \'etant respectivement un vecteur de
plus haut poids de $\Lambda^p$ et $\Lambda^q$ (sous un tore maximal 
de $K$). Notons-la $\tau_{pq}$. On a alors~:

\begin{thm}[Kumaresan, Vogan-Zuckerman \cite{VoganZuckerman}, Krajlevic \cite{Krajlevic}] \label{VZ}
\
\begin{description}
\item[{\rm (i)}] Pour tous $p,q$ tels que $p+q\leq n$ il existe une 
unique repr\'esentation unitaire irr\'eductible $V_{pq}$ de $G$ telle
que
\begin{description}
\item[{\rm (a)}] Hom$_K(\tau_{pq},V_{pq})\not =0$.
\item[{\rm (b)}] L'op\'erateur de Casimir op\`ere trivialement 
sur $V_{pq}$.
\end{description}
\item[{\rm (ii)}] Aucune repr\'esentation $V_{pq}$ n'est isol\'ee 
dans $\widehat G$.
\item[{\rm (iii)}] Les $V_{pq}$ constituent toutes les 
repr\'esentations unitaires cohomologiques de $G$.
\end{description}
\end{thm}
\index{Th\'eor\`eme de Kumaresan, Vogan-Zuckerman, Krajlevic}

Noter que d'apr\`es (a); (b),
$$
H^{p+q} (\mathfrak{g},K,V_{pq})=\mbox{Hom}_K (\Lambda^{p+q} \mathfrak{p}  ,V_{pq})\not =0\ .
$$

D'apr\`es \cite{VoganZuckerman} on sait en fait que les autres repr\'esentations de 
$\Lambda^{p+q} \mathfrak{p}$ n'interviennent pas dans $V_{pq}$, et que
$\tau_{pq}$ appara\^{\i}t avec multiplicit\'e~1. En particulier
\begin{eqnarray} \label{cohom1}
H^{p+q} (\mathfrak{g},K;V_{pq})\cong {\Bbb C}\ .
\end{eqnarray}

Par ailleurs (\ref{dec de Lp}) implique une d\'ecomposition de Hodge sur la 
cohomologie \cite[II, \S 4]{BorelWallach} et~:
\begin{eqnarray} \label{cohom2}
H^{p,q}(\mathfrak{g},K;V_{pq})\cong {\Bbb C}
\end{eqnarray}
et plus pr\'ecis\'ement
\begin{eqnarray} \label{cohom3}
H^{p+r,q+r}(\mathfrak{g},K;V_{pq})\cong {\Bbb C} \quad(r=0,\ldots n-p-q),
\end{eqnarray}
les autres composantes de Hodge de $H^*(\mathfrak{g},K;V_{pq})$ \'etant nulles 
\cite[Prop. 3.6]{VoganZuckerman}.

Il est clair que ces r\'esultats ne laissent aucun espoir quant au 
Probl\`eme {\bf local} \ref{pbm1}. En revanche nous allons voir (Ch. 6) que 
le
Probl\`eme global \ref{pbm2} devrait admettre une solution positive, 
d'apr\`es les Conjectures d'Arthur (Conjecture A de l'Introduction).

\newpage 

\markboth{CHAPITRE 3. GL(n)}{3.1. CLASSIFICATION DE LANGLANDS}

\chapter{GL(n)}

\section{Classification de Langlands}

La classification de Langlands pour $G=GL(n,\mathbb{R})$ d\'ecrit toutes les repr\'esentations admissibles
de $G$ \`a \'equivalence infinit\'esimale pr\`es. Nous notons 
$K=O(n)$ le sous-groupe compact maximal de $G$.

Une repr\'esentation admissible irr\'eductible $(\rho , V)$ admet un {\it caract\`ere central} \index{caract\`ere central} {\it i.e.}
un homomorphisme  $\omega_{\rho} :\mathbb{R}^* \rightarrow \mathbb{C}^*$ tel que 
$$\rho ( t I_n ) = \omega_{\rho} (t) Id_V , $$
pour tout $t \in \mathbb{R}^*$.

Soit $SL^{\pm} (m,\mathbb{R})$ le sous-groupe des \'el\'ements $g$ de $GL(m,\mathbb{R})$ tels que $|$det$(g)|=1$.
On va d'abord d\'ecrire certaines repr\'esentations irr\'eductibles de $SL^{\pm} (m,\mathbb{R})$ pour $m=1$ et 
$m=2$. Pour $m=1$ il n'y a que deux repr\'esentations, toutes les deux de dimension un; on note la repr\'esentation
triviale $1$ et ${\rm sgn}$ l'autre. Pour $m=2$, les repr\'esentations qui nous int\'eressent sont celles qui sont dans la
``s\'erie discr\`ete'' et que l'on note $D_l$ pour $l \in \mathbb{N}^*$. Pour $l \geq 1$, la repr\'esentation $D_l$
est induite de $SL(2,\mathbb{R} )$ :
\begin{eqnarray}
D_l = \mbox{ind}_{SL(2,\mathbb{R} )}^{SL^{\pm} (2,\mathbb{R} )} (D_l^+ ).
\end{eqnarray} \index{$D_l$} \index{s\'erie discr\`ete}
O\`u la repr\'esentation $D_l^+$ est donn\'ee par l'action de $SL(2,\mathbb{R} )$ sur l'espace des fonctions analytiques
$f$ sur le demi-plan de Poincar\'e de norme
$$||f|| =\left( \int \int |f(z)|^2 y^{l-1} dx dy \right)^{\frac12} $$
finie, l'action d'un \'el\'ement $g= \left( 
\begin{array}{cc}
a & b \\
c & d 
\end{array} \right)$ \'etant donn\'ee par
\begin{eqnarray}
D_l^+ (g) f(z) = (bz+d)^{-(l+1)} f \left( \frac{az+c}{bz+d} \right) .
\end{eqnarray} \index{$D_l^+$}
Les repr\'esentations $D_l$ de $SL^{\pm} (2,\mathbb{R} )$ sont irr\'eductibles et unitaires (cf. \cite{Knapp}, \cite{Lang})
\footnote{Le caract\`ere infinit\'esimal de $D_l$ est $l \rho$ o\`u $\rho \left( \begin{array}{cc} 
1 & 0 \\
0 & -1 
\end{array} \right) =1$, voir \cite[Probl\`eme 2 p. 276]{Knapp}.}.

Chaque repr\'esentation de $G$ est construite \`a partir de ``blocs \'el\'ementaires'' : les repr\'esentations de 
$GL(1,\mathbb{R})$ et de $GL(2,\mathbb{R} )$ que l'on obtient en tensorisant les repr\'esentations 
de $SL^{\pm}$ ci-dessus par un morphisme $a \mapsto |$det$a|_{\mathbb{R}}^t$ 
du groupe des matrices scalaires strictement positives de taille $1$ ou $2$ vers $\mathbb{C}^*$.
Ci-dessous, nous notons $|.|_{\mathbb{R}}$ la valeur absolue ordinaire et $t\in \mathbb{C}$.
On obtient donc les ``blocs \'el\'ementaires'' suivant :
\begin{eqnarray} \label{33}
\left. 
\begin{array}{r}
1 \otimes |.|_{\mathbb{R}}^t \\
{\rm sgn} \otimes |.|_{\mathbb{R}}^t 
\end{array} 
\right\} & \mbox{   pour  } GL(1,\mathbb{R} ), 
\end{eqnarray}
\begin{eqnarray} \label{34}
D_l \otimes | \mbox{det}(.)|_{\mathbb{R}}^t & \mbox{   pour  } GL(2,\mathbb{R}).
\end{eqnarray}

\`A chaque partition de $n$ en somme de $1$ et de $2$, not\'ee $(n_1 , \ldots , n_r )$, on associe le 
sous-groupe diagonal par blocs 
$${\Bbb M} = MA= GL(n_1 , \mathbb{R} ) \times \ldots \times GL(n_r , \mathbb{R} )$$
de $G$, o\`u comme d'habitude $A$ d\'esigne le sous-groupe des matrices diagonales positives dans le centre de ${\Bbb M}$ et 
$M$ est un produit de groupes $\{ \pm 1\}$ ou $SL^{\pm} (2,{\Bbb R})$. Pour chaque entier $j$ compris entre $1$ et $r$, soit $\sigma_j$ une repr\'esentation de 
$GL(n_j , \mathbb{R} )$ de la forme (\ref{33}) ou (\ref{34}) suivant que $n_j =1$ ou $2$; nous noterons $t_j$ le nombre
complexe $t$ correspondant. Alors $(\sigma_1 , \ldots , \sigma_r )$ d\'efinit, par produit tensoriel, 
une repr\'esentation du sous-groupe diagonal par blocs $MA$ que l'on peut \'etendre en une repr\'esentation
du sous-groupe parabolique correspondant $P=MAU$ (constitut\'e des matrices triangulaires sup\'erieures par blocs)
par l'identit\'e sur $U$ (les matrices strictement  triangulaires sup\'erieures par blocs). Nous notons alors :
\begin{eqnarray} \label{induite reelle}
I( \sigma_1 , \ldots , \sigma_r ) = \mbox{ind}_P^G (\sigma_1 , \ldots , \sigma_r ) ,
\end{eqnarray} \index{$I(\sigma_1 , \ldots , \sigma_r )$}
o\`u ind d\'esigne l'induction unitaire (cf. \cite{Knapp}). 

D'apr\`es \cite[Proposition 8.22]{Knapp}, la repr\'esentation d\'efinit par (\ref{induite reelle}) admet un caract\`ere infinit\'esimal \index{caract\`ere infinit\'esimal}
\begin{eqnarray} \label{chi reel}
(\lambda_{\sigma_1} , \ldots , \lambda_{\sigma_r } ) \in {\Bbb C}^n,
\end{eqnarray}
o\`u $\lambda_{\sigma_j}$ est \'egal \`a $t_j$ lorsque $n_j =1$ et est \'egal au couple $(t_j + \frac{l_j }{2} , 
t_j - \frac{l_j}{2} )$ lorsque $n_j=2$ et o\`u le vecteur (\ref{chi reel}) de ${\Bbb C}^n$ obtenu est vu comme 
une forme lin\'eaire sur la sous-alg\`ebre complexe de Cartan $(\mathfrak{a} \oplus \mathfrak{m} )^{{\Bbb C}}$
naturellement identifi\'ee \`a ${\Bbb C}^n$ via l'identification de $\mathfrak{g}^{{\Bbb C}}$ \`a l'alg\`ebre 
$\mathfrak{gl}_{n} ({\Bbb C} )$ \footnote{Par exemple si $n=n_1 =2$, on identifie 
$\mathfrak{a}^{{\Bbb C}}$ \`a ${\Bbb C} \left( \begin{array}{cc}
1 & 0 \\
0 & 1
\end{array} \right)$, $\mathfrak{m}^{{\Bbb C}}$ \`a ${\Bbb C}  \left( \begin{array}{cc}
1 & 0 \\
0 & -1 \end{array} \right)$ et on identifie un vecteur de ${\Bbb C}^2$ avec une forme lin\'eaire sur l'alg\`ebre des 
matrices diagonales de $M_2 ({\Bbb C})$.}.  

Le th\'eor\`eme qui suit r\'esulte des travaux de Langlands \cite{Langlands} et de la th\`ese de Speh (cf. aussi \cite{KnappZuckerman}).

\begin{thm}[Classification de Langlands pour $GL(n,\mathbb{R})$] \label{class reel}
Pour $G=GL(n,\mathbb{R} )$.
\begin{enumerate}
\item Si les param\`etres $t_j$ de $(\sigma_1 , \ldots , \sigma_r )$ v\'erifient
\begin{eqnarray} \label{37}
n_1^{-1} \mbox{Re} \, t_1 \geq n_2^{-1} \mbox{Re} \, t_2 \geq \ldots \geq n_r^{-1} \mbox{Re} \, t_r ,
\end{eqnarray}
alors $I(\sigma_1 , \ldots , \sigma_r )$ a un unique quotient irr\'eductible, son {\bf quotient de Langlands}, \index{quotient de Langlands} \index{$J(\sigma_1 , \ldots , \sigma_r )$}
que l'on note $J(\sigma_1 , \ldots , \sigma_r )$.
\item Toute repr\'esentation admissible irr\'eductible de $G$ est infinit\'esimalement \'equivalente \`a une 
repr\'esentation $J(\sigma_1 , \ldots , \sigma_r )$.
\item Deux repr\'esentations $J(\sigma_1 , \ldots , \sigma_r )$ et $J(\sigma_1 ' , \ldots , \sigma_r ')$ sont 
infinit\'esimalement \'equivalentes si et seulement si $r=r'$ et il existe une permutation $j(i)$ de 
$\{ 1, \ldots , r \}$ telle que $\sigma_i ' = \sigma_{j(i)}$ pour $1\leq i \leq r$.
\end{enumerate}
\end{thm}
\index{Classification de Langlands pour $GL(n, {\Bbb R} )$, Th\'eor\`eme de}

\medskip

On peut de la m\^eme mani\`ere classifier les repr\'esentations irr\'eductibles admissibles de 
$GL(n, \mathbb{C} )$. Dans la suite de cette section, nous notons $G=GL(n,\mathbb{C} )$ et $K=U(n)$ un 
sous-groupe compact maximal.

Rappelons que tout caract\`ere (non n\'ecessairement unitaire) de $\mathbb{C}^*$ s'\'ecrit 
$z \mapsto z^p \overline{z}^q$ o\`u $p,q \in \mathbb{C}$ et $p-q \in \mathbb{Z}$. Le caract\`ere est unitaire si, et 
seulement si, Re$(p+q)=0$. Posons $t_i = \frac{p_i + q_i }{2}$. Pour tout $j=1, \ldots , n$, soit $\sigma_j$ un 
caract\`ere de $\mathbb{C}^*$.
Alors $(\sigma_1 , \ldots , \sigma_n )$ d\'efinit, par produit tensoriel, une repr\'esentation de dimension
un du sous-groupe des matrices diagonales de $G$, que l'on \'etend, par l'identit\'e, en une repr\'esentation
de dimension un du sous-groupe $B$ des matrices triangulaires sup\'erieures dans $G$. On pose alors~:
\begin{eqnarray} \label{induite complexe}
I(\sigma_1 , \ldots , \sigma_n ) = \mbox{ind}_B^G ( \sigma_1 , \ldots , \sigma_n ) , 
\end{eqnarray} \index{$I(\sigma_1 , \ldots , \sigma_n )$}
o\`u l'induction ci-dessus est  unitaire.

D'apr\`es \cite[Proposition 8.22]{Knapp}, la repr\'esentation d\'efinit par (\ref{induite complexe}) admet un caract\`ere infinit\'esimal \index{caract\`ere infinit\'esimal}
\begin{eqnarray} \label{chi complex}
(p_1 , \ldots , p_n ) \times (q_1 , \ldots , q_n ) \in {\Bbb C}^n \times {\Bbb C}^n,
\end{eqnarray}
o\`u le vecteur (\ref{chi complex}) de ${\Bbb C}^n \times {\Bbb C}^n$ obtenu est vu comme 
une forme lin\'eaire sur la complexification de la sous-alg\`ebre de Cartan constitu\'ee des matrices complexes 
diagonales que l'on consid\`ere comme une alg\`ebre r\'eelle par restriction des scalaires de ${\Bbb C}$
\`a ${\Bbb R}$ \footnote{Autrement dit dans notre cas 
${\Bbb C}$ se plonge dans ${\Bbb C} \times {\Bbb C}$ par $p: z \mapsto (z, \overline{z} )$, puis l'application
$(u,v) \mapsto (\frac{u+v}{2} , \frac{u-v}{2} )$ r\'ealise un isomorphisme entre 
${\Bbb C} \times {\Bbb C}$ et une ${\Bbb R}$-alg\`ebre dont les points r\'eels s'identifient \`a l'image de $p$. On obtient ainsi 
des coordonn\'ees naturelles pour la complexification de ${\Bbb C}$. Dans ces coordonn\'ees le caract\`ere infinit\'esimal,
de la repr\'esentation $z \mapsto z^p \overline{z}^q$ sont $(p,q)$.}.

\begin{thm}[Classification de Langlands pour $GL(n,\mathbb{C})$] \label{class complexe}
Pour $G= GL(n,\mathbb{C})$. 
\begin{enumerate}
\item Si les param\`etres $t_j$ de $(\sigma_1 , \ldots , \sigma_n )$ v\'erifient
\begin{eqnarray} \label{310}
\mbox{Re} \, t_1 \geq \mbox{Re} \, t_2 \geq \ldots \geq \mbox{Re} \, t_n , 
\end{eqnarray}
alors $I(\sigma_1 , \ldots , \sigma_n )$ a un unique quotient irr\'eductible, son {\bf quotient de Langlands}, \index{quotient de Langlands}
que l'on note $J(\sigma_1 , \ldots , \sigma_n )$. \index{$J(\sigma_1 , \ldots , \sigma_n )$}
\item Toute repr\'esentation irr\'eductible admissible de $G$ est infinit\'esimalement \'equivalente \`a une
repr\'esentation $J(\sigma_1 , \ldots , \sigma_n )$.
\item Deux repr\'esentations $J(\sigma_1 , \ldots , \sigma_n )$ et $J(\sigma_1 ' , \ldots , \sigma_n ')$ sont 
infinit\'esimalement \'equivalentes si et seulement s'il existe une permutation $j(i)$ de $\{ 1 , \ldots , n \}$
telle que $\sigma_i ' = \sigma_{j(i)}$ pour $1\leq i \leq n$.
\end{enumerate}
\end{thm}
\index{Classification de Langlands pour $GL(n, {\Bbb C} )$, Th\'eor\`eme de}

Nous renvoyons au livre de Knapp \cite{Knapp} pour une d\'emonstration de ces Th\'eor\`emes.

\markboth{CHAPITRE 3. GL(n)}{3.2. CORRESPONDANCE DE LANGLANDS LOCALE}

\section{Correspondance de Langlands locale}

Le {\it groupe de Weil de $\mathbb{R}$}, \index{groupe de Weil r\'eel} not\'e $W_{\mathbb{R}}$, \index{$W_{\mathbb{R}}$} est l'extension de $\mathbb{C}^*$ par 
$\mathbb{Z} /2\mathbb{Z}$ (le groupe de Galois de $\mathbb{C}$ sur $\mathbb{R}$) :
$$W_{\mathbb{R}} = \mathbb{C}^* \cup j \mathbb{C}^* ,$$
o\`u $j^2 =-1$ et $jcj^{-1} = \overline{c}$. 
Dans cette section nous d\'ecrivons la {\it correspondance de Langlands locale pour $GL(n,\mathbb{R})$} {\it i.e.} \index{correspondance de Langlands locale pour $GL(n, {\Bbb R})$}
une bijection entre l'ensemble des classes d'\'equivalence de repr\'esentations complexes semi-simples de dimension $n$ 
de $W_{\mathbb{R}}$ et l'ensemble des classes
d'\'equivalence de repr\'esentations admissibles irr\'eductibles de $GL(n,\mathbb{R})$. On va donc 
s'int\'eresser aux repr\'esentations semi-simples du groupe de Weil de $\mathbb{R}$.

Il existe une suite exacte naturelle 
$$1 \rightarrow U \rightarrow W_ {{\Bbb R}} \stackrel{u}{\rightarrow} {\Bbb R}^* \rightarrow 1$$
donn\'ee par 
$$\begin{array}{l}
z \mapsto |z |_{{\Bbb C}} = z \overline{z} \  (z \in {\Bbb C}^*) \\
j \mapsto -1, 
\end{array}$$
et donc le noyau est $\{ z \in {\Bbb C}^* \, : \, |z|=1 \}$. Si $|z|=1$, on a $z=w/\overline{w} = w(j w j^{-1} )^{-1}$
qui est \'egal au commutateur $[w,j] = wjw^{-1} j^{-1}$. Donc $U$ s'identifie au groupe d\'eriv\'e de $W_{{\Bbb R}}$.

Par cons\'equent, les caract\`eres ab\'eliens de $W_{{\Bbb R}}$ sont de la forme 
$w \mapsto \chi ( u (w))$ o\`u $\chi$ est un caract\`ere de ${\Bbb R}^*$. 
Les repr\'esentations 
de dimension un sont donc param\`etr\'ees par un signe et un param\`etre complexe $t$ comme ci-dessous :
\begin{eqnarray} \label{rep de dim 1 de wr}
\begin{array}{cccc}
(+,t): & \varphi (z) = |z|_{\mathbb{R}}^t & \mbox{ et } & \varphi (j) = +1 ,\\
(-,t):  & \varphi (z) = |z|_{\mathbb{R}}^t & \mbox{ et } & \varphi (j) = -1.
\end{array}
\end{eqnarray}

Consid\'erons maintenant une repr\'esentation irr\'eductible $r$ de $W_{\mathbb{R}}$ sur un espace $V$, 
de degr\'e $\geq 2$. Alors $r_{| {\Bbb C}^*}$ est compl\`etement r\'eductible, donc somme directe de 
caract\`eres. Si $\chi$ est un caract\`ere de ${\Bbb C}^*$, soit 
$$V(\chi ) = \{ v \in V \, : \, z.v = \chi (z) v \} .$$
Puisque $jzj^{-1} =  \overline{z}$ ($z \in {\Bbb C}^*$), la somme directe $V(\chi ) \oplus V(\chi^{\sigma} )$, o\`u 
$\chi^{\sigma} (z) = \chi (\overline{z} )$, est stable par $W_{{\Bbb R}}$. Par cons\'equent, il existe au plus un 
couple de caract\`eres $\{ \chi ,\chi^{\sigma} \}$ dans la d\'ecomposition de $V$. 

Si $\chi = \chi^{\sigma}$, d'apr\`es le paragraphe au-dessus de (\ref{rep de dim 1 de wr}) $\chi$ s'\'etend en un caract\`ere ab\'elien 
$\chi '$ de $W_{{\Bbb R }}$. Alors $r\otimes (\chi ')^{-1}$ est une repr\'esentation irr\'eductible de $W_{{\Bbb R}}$, 
triviale sur ${\Bbb C}^*$, donc repr\'esentation de $W_{{\Bbb R}} / {\Bbb C}^* \cong \{ \pm 1 \}$. Elle est donc de 
dimension $1$, et $r$ est un caract\`ere, d\'ecrit dans (\ref{rep de dim 1 de wr}). 

Si $\chi \neq \chi^{\sigma}$, soit $v \in V(\chi )$. Alors $r(j)v \in V(\chi^{\sigma} )$ donc $v$ et $r(j)v$ sont ind\'ependants.
Alors, l'espace $\langle v,r(j)v \rangle$ est stable par $W_{{\Bbb R}} = \langle j, {\Bbb C}^* \rangle$. Donc $r$ est de dimension $2$, donn\'ee dans 
cette base par les matrices 
$$z \mapsto \left(
\begin{array}{cc}  
\chi (z) & \\
 & \chi ( \overline{z} ) 
\end{array}
\right), \,  z \in {\Bbb C}^*$$
$$j \mapsto \left(
\begin{array}{cc}
 & \chi (-1 ) \\
1 &  
\end{array}
\right) .$$
C'est la repr\'esentation de $W_{{\Bbb R}}$ {\it induite} \`a partir du caract\`ere $\chi$ de ${\Bbb C}^*$.
Nous avons donc obtenu que les classes d'\'equivalence de repr\'esentations 
irr\'eductibles $\varphi : W_{\mathbb{R}} \rightarrow GL(2,\mathbb{C})$ sont classifi\'ees par les paires
$(l,t)$ avec $l$ un entier $\geq 1$ et $t$ dans $\mathbb{C}$; \`a chaque 
paire $(l,t)$ correspond la classe d'\'equivalence du morphisme suivant~:
\begin{eqnarray} \label{rep de dim 2 de wr}
(l,t) : \left\{
\begin{array}{l}
\varphi (r e^{i\theta} ) = \left(
\begin{array}{cc}
r^{2t} e^{il\theta} & \\
 & r^{2t} e^{-il\theta} 
\end{array}
\right) \\
 \mbox{et } \\
  \varphi (j) = \left(
\begin{array}{cc}
 & (-1)^l \\
1 & 
\end{array} 
\right) .
\end{array} \right.
\end{eqnarray}

\medskip

Soit maintenant $\varphi$ une repr\'esentation complexe semi-simple de dimension $n$ de $W_{\mathbb{R}}$.
La liste des dimensions des composantes 
irr\'eductibles de $\varphi$ donne une partition de $n$ comme somme de $1$ et de $2$ que l'on note
$(n_1 , \ldots ,n_r )$. Soit $\varphi_j$ la composante irr\'eductible de $\varphi$ correspondante \`a l'entier $n_j$.
\`A $\varphi_j$ on peut associer une repr\'esentation $\sigma_j$ de $GL(n_j , {\Bbb R})$ de la mani\`ere suivante~:
\begin{eqnarray}
\begin{array}{lcl}
(+,t) \mbox{  dans (\ref{rep de dim 1 de wr})} & \mapsto & 1 \otimes |.|_{\mathbb{R}}^t \mbox{  dans (\ref{33})}, \\
(-,t) \mbox{  dans (\ref{rep de dim 1 de wr}) }  & \mapsto & {\rm sgn} \otimes |.|_{\mathbb{R}}^t \mbox{  dans (\ref{33})}, \\
(l,t) \mbox{  dans (\ref{rep de dim 2 de wr}) }  & \mapsto & D_l \otimes  |\mbox{det}(.)|_{\mathbb{R}}^t \mbox{  dans (\ref{34})}.
\end{array}
\end{eqnarray}
De cette fa\c{c}on, on associe \`a la repr\'esentation $\varphi$ un $r$-uplet $(\sigma_1 , \ldots ,\sigma_r )$.
De plus, quitte \`a permuter les $\sigma_j$, on peut supposer que (\ref{37}) est v\'erifi\'ee. Le Th\'eor\`eme \ref{class reel} permet 
donc de d\'efinir l'application~:
\begin{eqnarray} \label{314}
\varphi \mapsto \rho_{\mathbb{R}} (\varphi )=J(\sigma_1 , \ldots , \sigma_r )
\end{eqnarray}
et implique le~:

\begin{thm}[Correspondance de Langlands locale pour $GL(n,\mathbb{R})$] \label{CLLR}
L'application (\ref{314}) d\'efinit une bijection entre l'ensemble des classes d'\'equivalence de repr\'esentations complexes
semi-simples de dimension $n$ de $W_{\mathbb{R}}$ et l'ensemble des classes d'\'equivalence de repr\'esentations
admissibles irr\'eductibles de $GL(n,\mathbb{R})$.
\end{thm}
\index{correspondance  de Langlands locale pour $GL(n,\mathbb{R})$}

Remarquons de plus que la classification de Langlands implique l'\'equivalence des trois assertions suivantes :
\begin{enumerate}
\item le param\`etre $\varphi$ est {\it temp\'er\'e}, \index{tempr\'er\'e, param\`etre} {\it i.e.} est born\'e;
\item Re$(t_1) = \ldots =$Re$(t_r ) =0$;
\item la repr\'esentation $\rho_{\mathbb{R}} (\varphi )$ est temp\'er\'ee.
\end{enumerate}

\medskip

Dans la suite nous d\'ecrivons l'analogue  du Th\'eor\`eme \ref{CLLR} pour $GL(n,\mathbb{C})$. 
Le {\it groupe de Weil de $\mathbb{C}$}, \index{groupe de Weil complexe} not\'e $W_{\mathbb{C}}$, \index{$W_{\mathbb{C}}$} est donn\'e par $W_{\mathbb{C}} = \mathbb{C}^*$.

Une repr\'esentation semi-simple de dimension $n$ de $W_{\mathbb{C}}$, $\chi$, est donc donn\'ee par $n$ (quasi-)caract\`eres \index{caract\`eres (quasi-)}
$\chi_1 , \ldots , \chi_n$ o\`u $\chi_i (z) = z^{p_i} (\overline{z})^{q_i}$. Comme plus haut posons
$t_i = \frac{p_i + q_i}{2}$. Quitte \`a permuter
les $\chi_i$, on peut supposer que (\ref{310}) est v\'erifi\'ee. Le Th\'eor\`eme \ref{class complexe} permet donc de d\'efinir l'application~:
\begin{eqnarray} \label{315}
\chi \mapsto \rho_{\mathbb{C}} (\chi ) =J(\chi_1 , \ldots , \chi_n ) 
\end{eqnarray}
et implique le~:

\begin{thm}[Correspondance de Langlands locale pour $GL(n,\mathbb{C} )$] \label{correspondance locale}
L'application (\ref{315}) d\'efinit une bijection entre l'ensemble des classes d'\'equivalence de repr\'esentations complexes
semi-simples de dimension $n$ de $W_{\mathbb{C}}$ et l'ensemble des classes d'\'equivalence de repr\'esentations
admissibles irr\'eductibles de $GL(n,\mathbb{C} )$.
\end{thm}
\index{correspondance de Langlands locale pour $GL(n,{\Bbb C})$}

Enfin, remarquons que l\`a encore la classification de Langlands implique que les trois assertions 
suivantes sont \'equivalentes :
\begin{enumerate}
\item le param\`etre $\chi$ est {\it temp\'er\'e}, \index{temp\'er\'e, param\`etre} {\it i.e.} born\'e;
\item Re$(t_1)= \ldots =$Re$(t_n ) =0$;
\item la repr\'esentation $\rho_{\mathbb{C}} (\chi )$ est temp\'er\'ee.
\end{enumerate}

\medskip

\markboth{CHAPITRE 3. GL(n)}{3.3. UN PEU DE FONCTIONS $L$}

\section{Un peu de fonctions L}

\subsection{Caract\`eres de Dirichlet et ad\`eles}

Pour simplifier nous commen\c{c}ons par quelques rappels sur ces notions en se pla\c{c}ant sur 
le corps $\mathbb{Q}$, mais tous ces r\'esultats admettent un analogue dans 
le cas d'un corps de nombres quelconque.

Soit $n$ un entier sup\'erieur ou \'egal \`a $1$. 
Soit $\chi_0$ un caract\`ere du groupe multiplicatif fini $(\mathbb{Z} / n \mathbb{Z} )^*$, {\it i.e.} 
un morphisme de groupe \`a valeurs dans le cercle unit\'e.
On appelle {\it caract\`ere de Dirichlet de module $n$} \index{caract\`ere de Dirichlet} toute application $\chi$ obtenue  
\`a partir d'un caract\`ere $\chi_0$ comme ci-dessus par l'extension \`a $\mathbb{Z}$ :
$$\chi (a  ) = \left\{ 
\begin{array}{ll}
0 & \mbox{ si } (a, n) \neq 1 ,\\
\chi_0 ( a \mbox{ mod } n ) & \mbox{ si } (a , n) =1 .
\end{array}
\right.$$
Nous dirons de plus que $\chi$ est {\it primitif  de conducteur $n$} \index{caract\`ere primitif} \index{conducteur (d'un caract\`ere)} si $\chi$ ne provient pas par composition 
d'un caract\`ere de $(\mathbb{Z} / d\mathbb{Z} )^*$ o\`u $d$ divise $n$.

\medskip

Dans la suite nous allons ramener en partie l'\'etude du spectre des formes diff\'erentielles des vari\'et\'es hyperboliques
complexe arithm\'etiques \`a l'\'etude du spectre des formes diff\'erentielles des rev\^etements de congruence 
de 
$$SL(n,{\Bbb Z}) \backslash SL(n,{\Bbb R}) / SO(n) .$$
Ce lien est pr\'edit par la fonctorialit\'e de Langlands, dont nous d\'ecrirons quelques cas particuliers dans les prochains 
chapitres.
  
Nous voulons donc \'etudier les formes diff\'erentielles sur l'espace sym\'etrique
$$X= SL(n,{\Bbb R}) / SO(n)$$
invariantes par rapport \`a un sous-groupe de congruence de $SL(n,{\Bbb Z})$ et ce simultan\'ement pour 
diff\'erents sous-groupes de congruence.
Rappelons qu'un sous-groupe de congruence de $SL(n,{\Bbb Z})$ est un sous-groupe de $SL(n,{\Bbb Z})$ 
qui contient un sous-groupe $\Gamma_N = \{ M \in SL(n, {\Bbb Z}) \; : \; M \equiv I_n \; (\mbox{mod} \; N) \}$, pour 
un certain entier $N \geq 1$. 
\`A chaque inclusion $\Gamma ' \subset \Gamma$ de sous-groupes de congruence correspond une projection de 
rev\^etement $\Gamma ' \backslash X \rightarrow \Gamma \backslash X$. Puisque le nombre de classes 
d'id\'eaux de ${\Bbb Q}$ est $1$, le th\'eor\`eme qui suit montre que la limite projective de ce syst\`eme,
$$\mbox{lim proj}_{\Gamma} \Gamma \backslash X = SL(n, {\Bbb Q}) \backslash SL(n, {\Bbb A}) / SO(n,{\Bbb R}) ,$$
o\`u ${\Bbb A}$ est l'anneau des ad\`eles de ${\Bbb Q}$.
Ainsi l'introduction des ad\`eles permet de transformer l'\'etude du spectre des quotients de congruence de $X$ en l'\'etude d'un 
certain espace de fonctions sur un espace de classes \`a droite et \`a gauche du groupe $SL(n,{\Bbb A})$.

\medskip

Dans la suite nous notons $F$ un corps de nombres et ${\Bbb A}$ son anneau des ad\`eles. 
Soit $S=S_{\infty} \cup S_f$ l'ensemble des places de $F$, r\'eunion des places archim\'ediennes et des 
places finies de $F$. Si $v\in S$, nous notons $F_v$ la compl\'etion de $F$ en la place $v$; si $v$ est non-archim\'edienne,
nous notons $\mathfrak{o}_v $ l'anneau des entiers de $F_v$. Nous notons ${\Bbb A}_f$ l'anneau des ad\`eles finis et 
$F_{\infty}$ le produit de toutes les compl\'etions archim\'ediennes de $F$. 
Nous supposerons une certaine familiarit\'e avec la d\'efinition
de l'anneau des ad\`eles, sa topologie, et le fait que $F \subset {\Bbb A}$ est un sous-groupe discret, avec un quotient
${\Bbb A} /F$ compact.

\begin{thm}[Th\'eor\`eme d'approximation forte] \label{approxforte}
Soit $F$ un corps de nombres.
\begin{enumerate}
\item $SL(n,F_{\infty} ) SL(n,F)$ est dense dans $SL(n,{\Bbb A} )$.
\item Soit $K_0$ un sous-groupe compact ouvert de $GL(n,{\Bbb A}_f )$. Supposons que l'image de
$K_0$ dans ${\Bbb A}_f^*$ par le d\'eterminant soit $\prod_{v\notin S_{\infty}} \mathfrak{o}_v^*$. Alors 
le cardinal de 
$$GL(n,F)GL(n,F_{\infty} ) \backslash GL(n,{\Bbb A} ) / K_0$$
est \'egal au nombre de classes d'id\'eaux de $F$.
\end{enumerate}
\end{thm}
\index{Th\'eor\`eme d'approximation forte}

Dans le premier chapitre on a ramen\'e l'\'etude du spectre des formes diff\'erentielles sur les quotients 
$\Gamma \backslash X$ \`a l'\'etude de certaines repr\'esentations irr\'eductibles unitaires du groupe $SL(n, {\Bbb R})$,
celles apparaissant dans le spectre automorphe. Les ad\`eles permettent de voir chacune de ces repr\'esentations comme
la ``repr\'esentation \`a l'infini'' d'un seul et m\^eme groupe, le groupe $SL(n, {\Bbb A})$. Il sera en fait plus commode
de consid\'erer le groupe $GL(n,{\Bbb A})$. 

Nous dirons d'une repr\'esentation unitaire irr\'eductible de $GL(n,{\Bbb A})$ qu'elle est {\it automorphe} \index{repr\'esentation automorphe} si elle 
apparait comme sous-repr\'esentation de la repr\'esentation r\'eguli\`ere droite dans 
$L^2 (GL(n,F) \backslash GL(n,{\Bbb A}) )$. Remarquons qu'en g\'en\'eral une repr\'esentation automorphe n'apparait
pas discr\`etement, ce probl\`eme peut-\^etre \'eviter en ne consid\'erant que des repr\'esentations {\it cuspidales} \index{repr\'esentation cuspidale}
pour la d\'efinition desquelles nous renvoyons \`a \cite{Bump}.

\medskip

On appelle {\it caract\`ere de Hecke} \index{caract\`ere de Hecke} un caract\`ere continu $\chi$
de ${\Bbb A}^* / F^*$, o\`u $F$ est un corps de nombres et ${\Bbb A}$ son anneau des ad\`eles. 
Si $v$ est une place non-archim\'edienne de $F$, on dit que $\chi$ est {\it non ramifi\'e} \index{ramifi\'e, caract\`ere non} en $v$ si 
$\chi_v$, la compos\'ee de $\chi$ et de l'inclusion naturelle $F_v \hookrightarrow {\Bbb A}$, est triviale sur
$\mathfrak{o}_v^*$. Dans le cas contraire on dit que $\chi$ est ramifi\'e en $v$. 
Il est facile de v\'erifier qu'un caract\`ere de Hecke est non ramifi\'e en presque toutes les places. 
Lorsque $F={\Bbb Q}$, il d\'ecoule de la description de ${\Bbb A}^* / {\Bbb Q}^*$ que  
les caract\`eres de Hecke correspondent aux caract\`eres primitifs de Dirichlet de la fa\c{c}on suivante.  

\begin{prop} \label{caracteredeHeckevsDirichlet}
\begin{enumerate}
\item Supposons que $F={\Bbb Q}$ et que $\chi$ est un caract\`ere de ${\Bbb A}^* / F^*$. Il existe alors un unique 
caract\`ere $\chi_1$ d'ordre fini de ${\Bbb A}^* /F^*$ et un unique nombre imaginaire pur $\lambda$ tel que 
$\chi (x) = \chi_1 (x) |x|^{\lambda}$.
\item Supposons que $F={\Bbb Q}$ et que $\chi$ est un caract\`ere d'ordre fini de ${\Bbb A}^* / F^*$. Il existe alors un 
entier $N$ dont les diviseurs premiers sont  exactement les places finies de ${\Bbb Q}$ en lesquelles $\chi$ est ramifi\'e,
et un caract\`ere primitif de Dirichlet $\chi_0$ modulo $N$ tel que si $p$ est un nombre premier ne divisant pas $N$, 
alors $\chi_0 (p) = \chi_p (p)$. Cette correspondance $\chi \mapsto \chi_0$ est une bijection entre les caract\`eres 
d'ordre fini de ${\Bbb A}^* / F^*$ et les caract\`eres primitifs de Dirichlet.
\end{enumerate}
\end{prop}

\subsection{Fonctions $L$}

Soit $F$ un corps de nombre et ${\Bbb A}$ son anneau des ad\`eles.
La th\'eorie de Langlands veut associer \`a chaque repr\'esentation automorphe de $GL(n, \mathbb{A})$ 
une fonction $L$ donn\'ee comme produit eul\'erien de facteurs $L$ \'el\'ementaires, un pour chaque 
place de $F$. D'apr\`es un th\'eor\`eme de Jacquet-Langlands \cite{JacquetLanglands} et Flath \cite{Flath},
toute repr\'esentation irr\'eductible admissible de $GL(n, {\Bbb A} )$ est un ``produit tensoriel restreint'' 
$\pi = \otimes \pi_v$ de repr\'esentations de $GL(n, F_v)$ o\`u $v$ d\'ecrit les places de $F$ (y compris les 
places archim\'ediennes). D'apr\`es Godement-Jacquet \cite{GodementJacquet} on peut associer \`a $\pi$
une fonction 
$$L(\pi , s) = \prod_v L(\pi_v , s) .$$ \index{fonction $L$}
Pour $v$ finie, $L(\pi_v ,s)$ est un {\bf facteur eul\'erien} \index{facteur eul\'erien} du type usuel. Pour $v$ infinie, c'est un produit de 
fonctions $\Gamma$. (Godement-Jacquet consid\'eraient les $\pi$ cuspidales, le cas g\'en\'eral est trait\'e par 
Jacquet \cite{Jacquet}.)

En fait, on peut associer \`a $\pi$ une {\bf famille} de fonctions $L$, associ\'ees \`a des choix de vecteurs dans 
l'espace de la repr\'esentation et \`a des fonctions auxiliaires. Nous aurons besoin de consid\'erer {\bf la} fonction 
$L$ ayant \`a  toutes les places le {\bf bon} facteur pr\'edit par la fonctorialit\'e de Langlands \cite{K}. Ceci est obtenu 
par Jacquet \cite{Jacquet}. Pour les places archim\'ediennes, d'apr\`es le \S 3.3, la repr\'esentation $\pi_v$ 
est associ\'ee \`a une repr\'esentation de degr\'e $n$ de $W_{{\Bbb C}}$ ou $W_{{\Bbb R}}$. Il suffit donc d'\'ecrire
les facteurs locaux associ\'ees \`a celle-ci, et ceci est d\'ecrit par la proposition suivante. 

\begin{prop}[Facteurs $L_{\infty}$ pour $GL(n,\mathbb{R} )$] \index{facteurs $L$ r\'eels}
\label{facteurLreel}
Soit $\pi = \rho_{\mathbb{R}} (\varphi )$ une repr\'esentation de $GL(n,\mathbb{R} )$.
Le facteur $L_{\infty}$, $L(\pi , s) = L(\varphi ,s)$ est le produit des $L(\varphi_j ,s)$ o\`u 
les $\varphi_j$ sont les composantes irr\'eductibles de $\varphi$. Et si $\varphi$ est irr\'eductible, on a :
$$L(\varphi ,s) = \left\{
\begin{array}{ll}
\pi^{-\frac{s+t}{2}} \Gamma \left( \frac{s+t}{2} \right) & \mbox{ si } \varphi \mbox{ est donn\'ee par } (+,t) \mbox{ dans (\ref{rep de dim 1 de wr})}, \\
\pi^{-\frac{s+t+1}{2}} \Gamma \left( \frac{s+t+1}{2} \right) & \mbox{ si } \varphi \mbox{ est donn\'ee par } (-,t) \mbox{ dans (\ref{rep de dim 1 de wr})}, \\
2(2\pi )^{-(s+t+\frac{l}{2})} \Gamma \left( s+t+\frac{l}{2} \right) & \mbox{ si } \varphi \mbox{ est donn\'ee par } (l,t) \mbox{ dans (\ref{rep de dim 2 de wr})}.
\end{array}
\right. $$
\end{prop}

\begin{prop}[Facteurs $L_{\infty}$ pour $GL(n,\mathbb{C} )$] \index{facteurs $L$ complexes}
\label{facteurLcomplex}
Soit $\pi = \rho_{\mathbb{C}} (\chi )$ une repr\'esentation de $GL(n,\mathbb{C} )$.
Le facteur $L_{\infty}$, $L(\pi ,s) = L(\chi ,s)$ est le produit des $L(\chi_j ,s)$ o\`u 
les $\chi_j$ sont les caract\`eres composant $\chi$. Et si $\chi$ est un caract\`ere du type
$\chi (z) = z^p (\overline{z})^q$, on a~:
$$L(\chi ,s ) =  2(2\pi )^{-(s+\max (p,q))} \Gamma \left( s+\max (p,q) \right) .$$
\end{prop}

\medskip
 
\markboth{CHAPITRE 3. GL(n)}{3.4. DUAL UNITAIRE}

\section{Dual unitaire}

\subsection{Repr\'esentations de Speh}

Soit $F$ un corps \'egal \`a $\mathbb{R}$ ou $\mathbb{C}$. 
Soit $r$ un entier \'egal \`a $1$ si $F=\mathbb{C}$ et \'egal \`a $1$ ou $2$ si $F=\mathbb{R}$.
Soit $\delta$ une repr\'esentation de la s\'erie discr\`ete unitaire de $GL(r,F)$.
Ceci \'equivaut \`a ce que 
\begin{itemize}
\item $\delta (z) = z^p (\overline{z} )^q$ avec Re$(\frac{p+q}{2})=0$ si $F=\mathbb{C}$,  
\item $\delta = 1 \otimes |.|_{\mathbb{R}}^t$ ou ${\rm sgn} \otimes |.|_{\mathbb{R}}^t$ avec dans les deux cas 
Re$(t)=0$ si $F=\mathbb{R}$ et $r=1$ et,
\item $\delta = D_l \otimes |$det$(.)|_{\mathbb{R}}^t$ si $F=\mathbb{R}$ et $r=2$.
\end{itemize}
Soit $n=mr$ un multiple de $r$.
On note $\delta |.|^s$ ($s\in \mathbb{C}$) la repr\'esentation de $GL(r,F)$ donn\'ee par $\delta (g) |$det$(g)|_{F}^s$.
Consid\'erons la repr\'esentation de $GL(n,F)$ unitairement induite \`a partir de 
$$(\delta |.|^{\frac{m-1}{2}} , \delta |.|^{\frac{m-3}{2}} , \ldots , \delta |.|^{\frac{1-m}{2}}).$$
D'apr\`es les Th\'eor\`emes \ref{class reel} et \ref{class complexe}, elle admet un unique quotient de Langlands. On le note
$Sp (\delta ,m)$.

\begin{thm}[Speh \cite{Speh}]
La repr\'esentation $Sp(\delta , m)$ est unitaire.
\end{thm}
\index{Th\'eor\`eme de Speh} \index{$Sp(\delta , m)$}

\subsection{S\'eries compl\'ementaires}

Si ${\cal P}$ : $n=n_1 + \ldots +n_r$ ($n_i >0$) est une partition de $n$, on note $P$ le parabolique
associ\'e, form\'e des matrices triangulaires sup\'erieures par blocs. Son sous-groupe de Levi est 
$M= GL(n_1 , F) \times \ldots \times GL(n_r , F)$. Si $\pi_i$ est une repr\'esentation irr\'eductible de 
$GL(n_i , F)$, on note (de fa\c{c}on coh\'erente avec les notations de la premi\`ere section) 
$(\pi_1 , \ldots , \pi_r )$ la repr\'esentation de $P$ obtenue par extension triviale \`a $P$ de la repr\'esentation
tensorielle $\pi_1 \otimes \ldots \otimes \pi_r$ de $M$.

Soit $Sp ( \delta ,m)$ une repr\'esentation de Speh pour $GL(n,F)$, $n=mr$.
Dans $GL(2n , F)$, 
$$(Sp (\delta ,m) |.|^{\alpha} , Sp (\delta , m) |.|^{\beta } )$$
d\'efinit une repr\'esentation du sous-groupe parabolique associ\'e \`a la d\'ecomposition 
$2n = n+n$. \index{$(Sp (\delta ,m) |.|^{\alpha} , Sp (\delta , m) |.|^{\beta } )$}

\begin{thm}[Vogan \cite{Vogan}]
Pour $0 \leq \alpha < \frac12$, la repr\'esentation de $GL(2n,F)$ unitairement induite \`a partir de 
$$(Sp ( \delta , m) |.|^{\alpha} , Sp (\delta , m) |.|^{-\alpha} )$$
est irr\'eductible et unitaire.
\end{thm}
\index{repr\'esentation de Vogan}

Soit $V(\delta ,m,\alpha )$ la repr\'esentation ainsi d\'efinie. \index{$V(\delta , m, \alpha)$}

\subsection{Classification}

\begin{thm}[Vogan \cite{Vogan}] \label{Unit}
\begin{enumerate}
\item Soit $\pi$ une repr\'esentation unitaire de $GL(n,F)$. Alors, il existe une partition
$n=m_1 r_1 +\ldots +m_s r_s +2(m_{s+1} r_{s+1} +\ldots +m_{s+t} r_{s+t} )$, des repr\'esentations 
$\delta_i$ de la s\'erie discr\`ete de $GL(r_i , F)$ ($i=1, \ldots , s$) et des r\'eels $0<\alpha_i < \frac12$ ($i=s+1, 
\ldots , s+t$) tels que $\pi$ soit \'equivalente \`a l'induite unitaire de 
$$(Sp(\delta_1 ,m_1 ) , \ldots , Sp ( \delta_s , m_s ), V(\delta_{s+1} , m_{s+1} , \alpha_{s+1} ), \ldots , 
V(\delta_{s+t} , m_{s+t} , \alpha_{s+t} ) ).$$
\item Cette expression est unique, aux permutations pr\`es des $(\delta_i , m_i )$ ($i\leq s$) et des 
$(\delta_j , m_j , \alpha_j )$ ($j>s$).
\end{enumerate}
\end{thm}
\index{Th\'eor\`eme de Vogan}

Le Th\'eor\`eme ci-dessus implique que si les $t_j$ sont les param\`etres complexes, 
comme dans (\ref{37}) ou (\ref{310}), d'une repr\'esentation unitaire de $GL(n,F)$ alors~:
\begin{eqnarray}
\{ \overline{t_j} \} = \{ - t_k \}.
\end{eqnarray}

Enfin, remarquons que le Th\'eor\`eme de classification de Vogan permet (cf. \cite{BR}) de red\'emontrer que 
pour une repr\'esentation unitaire g\'en\'erique (ce qui est le cas de toute repr\'esentation automorphe 
cuspidale) les $t_j$ comme ci-dessus v\'erifient :
\begin{eqnarray}
|\mbox{Re}(t_j ) | < \frac12 .
\end{eqnarray}

\subsection{Retour sur les repr\'esentations de Speh}

Il est int\'eressant de noter que la construction de Speh est assez g\'en\'erale dans le sens qu'elle ne 
n\'ecessite pas r\'eellement que $\delta$ soit une repr\'esentation de la s\'erie discr\`ete. 

Pla\c{c}ons nous, pour simplifier, sur le corps ${\Bbb C}$. Soit $a$ un entier $\geq 1$ et $\tau$ une repr\'esentation
(unitaire) de $GL(a,{\Bbb C})$. Enfin soit $n=ab$ un multiple de $a$. Comme au-dessus notons $\tau |.|^s$ $(s \in 
{\Bbb C} )$ la repr\'esentation de $GL(a,{\Bbb C} )$ donn\'ee par $\tau (g) |\mbox{det} (g) |_{{\Bbb C}}^s$. On 
introduit alors la repr\'esentation de $GL(n,{\Bbb C})$ unitairement induite \`a partir de 
\begin{eqnarray} \label{J}
(\tau |.|^{\frac{b-1}{2}} , \tau |.|^{\frac{b-3}{2} }, \ldots , \tau |.|^{\frac{1-b}{2}} ).
\end{eqnarray}
D'apr\`es les Th\'eor\`emes \ref{class reel} et \ref{class complexe}, elle admet un unique quotient de Langlands, on le note 
$$J(\tau , b).$$ \index{$J(\tau ,b)$}

Ces repr\'esentations, nous le verrons, jouent (conjecturalement au moins) le r\^ole de pav\'es \'el\'ementaires
dans la description du spectre automorphe. Concluons ce chapitre par le calcul du caract\`ere infinit\'esimal d'une 
telle repr\'esentation. 

Supposons que le caract\`ere infinit\'esimal de la repr\'esentation $\tau$ soit  \index{caract\`ere infinit\'esimal}
$$(\lambda , \mu ) = ((\lambda_1 , \ldots , \lambda_a ), (\mu_1 , \ldots , \mu_a )) \in {\Bbb C}^a \times {\Bbb C}^a$$
toujours dans la param\'etrisation d'Harish-Chandra, cf. \S 1.3.
Alors, d'apr\`es \cite[Proposition 8.22]{Knapp}, le caract\`ere infinit\'esimal de $J(\tau, b)$ est \index{caract\`ere infinit\'esimal}
\begin{eqnarray} \label{caractere de J}
((\lambda +\frac{b-1}{2} , \lambda +\frac{b-3}{2} , \ldots , \lambda +\frac{1-b}{2}),(\mu +\frac{b-1}{2},  \ldots , \mu +\frac{1-b}{2} ))
\end{eqnarray}
(vecteur de ${\Bbb C}^n \times {\Bbb C}^n$).  

\newpage

\thispagestyle{empty}

\newpage

\markboth{CHAPITRE 4. REPR\'ESENTATIONS DE $U(n,1)$}{4.1. CLASSIFICATION DE LANGLANDS}

\chapter{Repr\'esentations de $U(n,1)$}

Soit $G=U(n,1)$ le sous-groupe de $GL(n+1,\mathbb{C} )$ laissant invariante la \index{$U(n,1)$}
forme hermitienne :
$$|z_1 |^2 + \ldots +|z_n |^2 -|z_{n+1} |^2 .$$
Soit $K=U(n+1) \cap G = U(n) \times U(1)$. 
L'alg\`ebre de Lie $\LG_0$ de $G$ est :
$$\LG_0 = \{ M \in M(n,\mathbb{C} )/  \; ^t \overline{M} I_{n,1} + I_{n,1} M = 0 \} ,$$
o\`u 
$$I_{n,1} = \left( 
\begin{array}{cc}
I_n & 0 \\
0    & -1 
\end{array}
\right) .$$ 
L'alg\`ebre de Lie $\LK_0$ de $K$ est :
$$\LK_0 = \left\{ \left( 
\begin{array}{cc}
A & 0 \\
0 & i\theta 
\end{array}
\right)
/ \theta \in \mathbb{R}  , \; ^t \overline{A} +A =0 \right\} .$$

La d\'ecomposition de Cartan correspondante est :
$$\LG_0 = \LK_0 \oplus \LM_0 ,$$
o\`u 
$$\LM_0 = \left\{ \xi (z) = 
\left( 
\begin{array}{cc}
0 & z \\
^t \overline{z} & 0 
\end{array}
\right)  / z\in M_{n,1} (\mathbb{C} ) \right\} .$$
Cette d\'ecomposition est orthogonale relativement \`a la forme de Killing \index{forme de Killing}
\begin{eqnarray} \label{forme de killing}
B(X,Y) = \frac12 Tr(XY)
\end{eqnarray}
sur $\LG_0$. Pour cette normalisation de la forme de Killing les 
courbures sectionnelles de l'espace sym\'etrique associ\'e (l'espace hyperbolique complexe) sont comprises 
entre $-4$ et $-1$, cf. deuxi\`eme partie.

Soit ${\rm Ad}$ la repr\'esentation adjointe de $G$. Puisque tout \'el\'ement de $K$ peut s'\'ecrire 
$$k= \left( 
\begin{array}{cc}
U & 0 \\
0 & v 
\end{array} \right), \; \; U \in U(n), \; v \in U(1)$$
le groupe $K$ agit sur $\mathfrak{p}_0 \cong {\Bbb C}^n$ par $Ad(k) X =UXv^{-1}$, et cette action pr\'eserve la 
structure complexe naturelle $J$. On en d\'eduit donc la d\'ecomposition
\begin{eqnarray} \label{structure complexe}
\mathfrak{p} = \mathfrak{p}_0 \otimes {\Bbb C} = \mathfrak{p}_+ \oplus \mathfrak{p}_-  \; \mbox{  avec  } \mathfrak{p}_{\pm}^* \cong \overline{
\mathfrak{p}_{\pm} } \cong \mathfrak{p}_{\mp},
\end{eqnarray}
comme $K$-modules. Pour $0\leq r \leq 2n$, soit
$$\tau_r := \Lambda^r ({\rm Ad}_+ \oplus {\rm Ad}_- ) \cong \Lambda^r ({\rm Ad} \oplus \overline{{\rm Ad}} )$$ \index{$\tau_r$}
la repr\'esentation de $K$ sur $\Lambda^r \mathfrak{p} = \Lambda^r (\mathfrak{p}_+ \oplus \mathfrak{p}_- )$.
Il est bien connu \cite{BorelWallach} (d\'ecomposition de Hodge) que chaque repr\'esentation $\tau_r$ se d\'ecompose en 
$$\tau_r = \oplus_{p+q = r} \tau_{p,q},$$
o\`u 
$$\tau_{p,q} := \Lambda^p \overline{{\rm Ad}} \otimes \Lambda^q {\rm A}d$$ \index{$\tau_{p,q}$}
est la repr\'esentation de $K$ sur $\Lambda^{p,q} \mathfrak{p} = \Lambda^p \mathfrak{p}_- \otimes \Lambda^q \mathfrak{p}_+$.
Via la formule de Matsushima, la d\'ecomposition de Lefschetz de la cohomologie d'une vari\'et\'e hyperbolique
complexe d\'ecoule de la d\'ecomposition de la repr\'esentation $\tau_{p,q}$ en irr\'eductibles. D'apr\`es 
\cite{BorelWallach}, on obtient~:
\begin{eqnarray} \label{lefschetz}
\tau_{p,q} = \oplus_{k=0}^{\min (p,q)} \tau_{p-k , q-k} ' , \index{$\tau_{a,b}$}
\end{eqnarray}
o\`u pour chaque couple d'entiers positifs $(i,j)$ de somme $i+j \leq n$, $\tau_{i,j} '$ est une repr\'esentation 
irr\'eductible. 

\section{Classification de Langlands}

Commen\c{c}ons par quelques notations.
Soit 
\begin{eqnarray}
H_0 = \left(
\begin{array}{ccc}
0 & 0 & 1 \\
0 & 0_{n-1} & 0 \\
1 & 0 & 0 
\end{array} \right) \in \LM_0 .
\end{eqnarray}
Alors $\la :={\Bbb R} H_0$ est un sous-espace de Cartan ({\it i.e.} ab\'elien maximal) \index{sous-espace de Cartan} dans $\LM_0$, et le 
sous-groupe de Lie correspondant de $G$ est param\`etr\'e par les \'el\'ements~:
$$a_t := \exp (tH_0 ) = \left( 
\begin{array}{ccc}
\cosh t & 0 & \sinh t \\
0 & I_{n-1} & 0 \\
\sinh t & 0 & \cosh t 
\end{array} \right) $$
o\`u $t \in {\Bbb R}$.

Soit $\alpha \in \la_0^*$ d\'efinie par $\alpha (tH_0 ) = t$. Alors $R(\LG_0 , \la_0 )=\{ \pm \alpha , \pm 2\alpha \}$
est un syst\`eme restreint de racines de $(\LG_0 , \la_0 )$ avec un sous-syst\`eme positif $R^+(
\LG_0 ,\la_0 ) = \{ \alpha , 2\alpha \}$. On utilisera dor\'enavant l'identification :
$$\begin{array}{ccc}
\la^* & \stackrel{\sim}{\rightarrow} & \mathbb{C} , \\
s \alpha & \mapsto & s . 
\end{array} $$
Soient $\lln_0 =(\LG_0 )_{\alpha} \oplus (\LG_0 )_{2\alpha }$ la somme des espaces propres des racines strictement positives,
$N$ le sous-groupe de Lie de $G$ correspondant et $\rho$ la demi-somme des racines dans $R^+ (\LG_0 , \la_0
)$, compt\'ees avec multiplicit\'es. Alors $\rho = n\alpha $ (dim $(\LG_0 )_{\alpha }=2(n-1)$ et dim $(\LG_0 )_{2\alpha }=1$).
Soit $M$ le centraliseur de $A$ dans $K$. Alors :
$$M= \left\{ \left( 
\begin{array}{ccc}
e^{i\theta} & 0 & 0 \\
0 & U & 0 \\
0 & 0 & e^{i\theta} 
\end{array} \right) \; : \;  \theta \in \mathbb{R} \mbox{  et  } U \in U(n-1) \right\} .$$
Nous notons $P=MAN$ le sous-groupe parabolique (minimal) usuel de $G$.
\'Etant donn\'es $\sigma \in \widehat{M}$ et $s \in \mathbb{C} \cong \la^*$, l'application :
$$(\sigma \otimes e^s \otimes 1) 
(ma_t n) = e^s \sigma (m) $$
d\'efinit une repr\'esentation de $P$ dans $V_{\sigma }$.
Nous notons $\pi_{\sigma , s}$ l'{\it induite unitaire} \index{induite unitaire} de cette repr\'esentation {\it i.e.} l'action de $G$
sur l'espace 
\begin{eqnarray*}
{\cal H}_{\sigma ,s} & = & L^2 (G, MAN , \sigma \otimes e^s \otimes 1 ) \\
                                & := & \{ f:G\rightarrow V_{\sigma } \; : \; f(xma_t n ) = e^{-(s+\rho)t} \sigma (m)^{-1} f(x) \mbox{ et }
f_{|K} \in L^2 (K) \} 
\end{eqnarray*} \index{${\cal H}_{\sigma , s}$}
par translation \`a gauche : $\pi_{\sigma ,s} (g) f(h) = f(g^{-1} h)$.
Remarquons que, comme $K$-module, chaque ${\cal H}_{\sigma , s}$ est isomorphe (pour tout $s$) \`a l'espace 
$L^2 (K,M, \sigma )$ des fonctions $f$
de carr\'e int\'egrable sur $K$ telles que $f(km)=\sigma (m^{-1} ) f(k)$.

Puisque $\LG_0$ et $\LK_0$ ont m\^eme rang $n$, le groupe $G$ admet une {\it s\'erie discr\`ete} \index{s\'erie discr\`ete} de repr\'esentations {\it i.e.}
des repr\'esentations irr\'eductibles unitaires dont les coefficients matriciels sont dans $L^2$.

\subsection{Repr\'esentations induites de $P$}

Dans cette sous-section nous revenons sur les repr\'esentations induites de $P$ pour d\'eterminer leurs $K$-types.

Puisque ${\cal H}_{\sigma , s}$ est $L^2 (K,M,\sigma )$, vue comme $K$-module, le th\'eor\`eme de r\'eciprocit\'e
de Frobenius nous dit que 
$${\rm Hom}_K ({\cal H}_{\sigma , s} , V_{\tau_{p,q} '} ) \cong {\rm Hom}_M (V_{\sigma} , V_{\tau_{p,q} '} ) \;  \mbox{  pour tout } s \in {\Bbb C}.$$
\index{Th\'eor\`eme de r\'eciprocit\'e de Frobenius}
On doit donc comprendre comment la repr\'esentation $\tau_{p,q} '$ se restreint \`a $M$. Ceci est classique, nous suivons ici 
l'article \cite{Pedon} de Pedon.

Soit $\mathfrak{h}_0$ (resp. $\mathfrak{t}_0$) la sous-alg\`ebre de Cartan de $\mathfrak{k}_0$ (resp. $\mathfrak{m}_0$) constitu\'ee
des \'el\'ements diagonaux. Pour chaque entier $1\leq i \leq n+1$, soit $\varepsilon_i$ la forme lin\'eaire sur 
$\mathfrak{h}$ d\'efinie par $\varepsilon_i ( \mbox{diag} (h_1 , \ldots , h_{n+1} )) =h_i$. Nous conservons
les m\^emes notations pour leurs restrictions \`a $\mathfrak{t}$. Les racines pour les paires $(\mathfrak{k} , 
\mathfrak{h} )$ et $(\mathfrak{m} , \mathfrak{t} )$ sont, respectivement,
\begin{eqnarray} \label{racines}
R_K := R( \mathfrak{k} , \mathfrak{h} ) = \{ \varepsilon_i - \varepsilon_j \; : \; 1 \leq i \neq j \leq n \}, \\
R_M := R(\mathfrak{m} , \mathfrak{t} ) = \{ \varepsilon_i - \varepsilon_j \; : \; 2 \leq i \neq j \leq n \},
\end{eqnarray}
et les syst\`emes positifs correspondants (pour l'ordre lexicographique) sont~:
\begin{eqnarray} \label{racines +}
R_K^+ = \{ \varepsilon_i - \varepsilon_j \; : \; 1 \leq i < j \leq n \}, \\
R_M^+ = \{ \varepsilon_i - \varepsilon_j \; : \; 2 \leq i<j \leq n \}. 
\end{eqnarray}

Les poids de la repr\'esentation adjointe ${\rm Ad}$ de $K$ sont les $\varepsilon_i -\varepsilon_{n+1}$, avec comme vecteur
de poids correspondant $e_i$ pour $1\leq i \leq n$. On en d\'eduit que les plus hauts poids des repr\'esentations 
irr\'eductibles $\Lambda^p \overline{{\rm Ad}}$ et $\Lambda^q {\rm Ad}$ de $K$ sont, respectivement, 
$$\mu_{\Lambda^p \overline{{\rm Ad}}} = -\sum_{k=n-p+1}^n \varepsilon_k +p \varepsilon_{n+1} \mbox{  et  }
\mu_{\Lambda^q {\rm Ad}} = \sum_{k=1}^q \varepsilon_k -q \varepsilon_{n+1} .$$
D'o\`u il d\'ecoule, \`a l'aide de \cite[Lemme 4.9, Chapitre VI]{BorelWallach}, que le plus haut poids de 
$\tau_{p,q}'$ est~:
\begin{eqnarray} \label{poids de tau}
\mu_{\tau_{p,q} ' } = \sum_{k=1}^q \varepsilon_k - \sum_{k=n-p+1}^n \varepsilon_k + (p-q) \varepsilon_{n+1} .
\end{eqnarray}
D'apr\`es le \cite[Th\'eor\`eme 4.4]{BaldoniSilva}, on a g\'en\'eriquement \footnote{Ce qui signifie que 
l'on demande que $\sigma_{a,b}  =0$ si $\min (a,b) <0$ ou si $\max (a,b) >n-1$.}  
$$(\tau_{p,q} ')_{|M}  = \sigma_{p,q}  \oplus \sigma_{p-1 ,q}  \oplus \sigma_{p, q-1}  \oplus  \sigma_{p-1,q-1} ,$$ \index{$\sigma_{a,b}$}
o\`u chaque $\sigma_{a,b}$ est une repr\'esentation irr\'eductible de $M$ de plus haut poids
\begin{eqnarray} \label{poids de sigma}
\mu_{\sigma_{a,b} } = \frac{a-b}{2} (\varepsilon_1 + \varepsilon_{n+1} ) + \sum_{k=2}^{b+1} \varepsilon_k - \sum_{k=n-a+1}^n
\varepsilon_k .
\end{eqnarray}

\subsection{S\'erie discr\`ete}

Revenons maintenant sur les repr\'esentations de la s\'erie discr\`ete. 

Comme \`a la section pr\'ec\'edente, soit $\mathfrak{h}_0 \subset \mathfrak{k}_0 \subset \mathfrak{g}_0$ la sous-alg\`ebre de Cartan
constitu\'ee des matrices diagonales. Avec les notations ci-dessus,
$$R_G := R(\mathfrak{g} , \mathfrak{h} ) = \{ \varepsilon_i - \varepsilon_j \; : \; 1 \leq i\neq j \leq n+1 \},$$
alors que $R_K$ et $R_K^+$ sont respectivement d\'efinis en (\ref{racines}) et (\ref{racines +}). Le groupe 
de Weyl $W_G$ (resp. $W_K$) est  le groupe  des permutations de $n+1$ (resp. $n$) \'el\'ements. Il y a donc 
$n+1$ sous-syst\`emes positifs de $R_G$ qui sont compatibles avec $R_K^+$. On choisit 
$$R_G^+ = \{ \varepsilon_i - \varepsilon_j \; : \; 1\leq i<j  \leq n+1 \}.$$
On obtient alors chaque syst\`eme positif compatible de la fa\c{c}on suivante. Soit $s_{\beta}$ la r\'eflexion relative \`a la 
racine $\beta$, on pose 
$$w_j = \prod_{k=j+1}^n s_{\varepsilon_k - \varepsilon_{n+1} } \; \mbox{ pour chaque } 0 \leq j \leq n-1, \; \mbox{ et } w_n =id .$$
Alors $W_K \backslash W_G = \{ W_K  w_j \; : \; 0 \leq j \leq n \}$ et les $n+1$ sous-syst\`emes positifs de $R_G$ compatibles
avec $R_K^+$ sont exactement les $w_j . R_G^+$, avec $0 \leq j \leq n$. Les sommes de racines positives sont~:
\begin{eqnarray} \label{delta G}
2 \delta_G = \sum_{k=1}^{n+1} (n+2-2k ) \varepsilon_k ,
\end{eqnarray}
\begin{eqnarray} \label{delta K}
2 \delta_K = \sum_{k=1}^n (n+1 -2k) \varepsilon_k .
\end{eqnarray}

Un r\'esultat classique d'Harish-Chandra, cf. \cite[Th\'eor\`emes 9.20 et 12.21]{Knapp} par exemple, affirme que les repr\'esentations
de la s\'erie discr\`ete \index{s\'erie discr\`ete} de $G$ sont, \`a \'equivalence pr\`es, uniquement d\'etermin\'ees par leur param\`etre 
d'Harish-Chandra $w_j \lambda$, o\`u $\lambda \in (i \mathfrak{h}_0^+ )'$ est tel que $\lambda + \delta_G$
est analytiquement int\`egre et $i \mathfrak{h}_0^+$ est la chambre de Weyl positives dans $i\mathfrak{h}_0$ 
correspondant \`a $R_G^+$. Nous notons $\pi_{w_j. \Lambda}$ la repr\'esentation de la s\'erie discr\`ete correspondante. \index{s\'erie discr\`ete}

Toujours d'apr\`es le \cite[Th\'eor\`eme 9.20]{Knapp}, la restriction \`a $K$ de la repr\'esentation 
$\pi_{w_j  \lambda}$ contient avec multiplicit\'e un la repr\'esentation de plus haut poids 
\begin{eqnarray} \label{poids de la sd}
w_j \lambda + w_j  \delta_G - 2 \delta_K .
\end{eqnarray}
De plus, tout plus haut poids d'un $K$-type de $\pi_{w_j \lambda}$, est de la forme
$$w_j  \lambda + w_j  \delta_G - 2 \delta_K + \sum_{\alpha \in R_G^+ } n_{\alpha} \alpha ,$$
o\`u les $n_{\alpha}$ sont des entiers $\geq 0$. 

Le fait suivant fait partie du "folklore" mais n'est \`a notre connaissance pas facile \`a trouver explicitement 
\'ecrit dans la litt\'erature. Dans notre cas le fait suivant peut-\^etre v\'erifi\'e \`a la main, la d\'emonstration que l'on propose,
valable en toute g\'en\'eralit\'e, est tir\'ee de la th\`ese de Pedon \footnote{On remercie Parthasarathy et Gaillard de nous avoir envoy\'e des 
d\'emonstrations diff\'erentes de ce fait.}.  

\medskip

\noindent
{\bf Fait.} Les repr\'esentations 
$\pi_{w_j  \delta_G}$ sont exactement les repr\'esentations de la s\'erie discr\`ete de $G$ qui contiennent 
un $K$-type intervenant comme sous-repr\'esentation de la repr\'esentation $\tau_r$ de $K$ pour un certain entier
$0 \leq r \leq n$. De plus $r$ doit alors \^etre \'egal \`a $n$.

\medskip

En effet, soit $\pi_{w_j  \lambda}$ une repr\'esentation de la s\'erie discr\`ete de $G$. Puisque $\lambda +\delta_G$ est 
int\`egre et dans la chambre de Weyl positive $(i\mathfrak{h}^+ )'$, on a $\lambda + \delta_G =  \sum_{\alpha \in R_G^+ } 
m_{\alpha} \alpha$ avec chaque $m_{\alpha} \in {\Bbb N}^*$. D'un autre c\^ot\'e, on a par d\'efinition, 
$\sum_{\alpha \in R_G^+ } \alpha = 2\delta_G$. Donc, $\lambda ':= \lambda - \delta_G = \sum_{\alpha \in R_G^+ } 
(m_{\alpha} -1 ) \alpha$ appartient au r\'eseau positif engendr\'e par $R_G^+$. 

Rappelons que tout plus haut poids d'un $K$-type de $\pi_{w_j \lambda}$ est de la forme 
$w_j \lambda ' + w_j  2\delta_G -2\delta_K + \sum_{\alpha \in R_G^+ } n_{\alpha} \alpha$, o\`u les 
$n_{\alpha}$ sont des entiers $\geq 0$.

Consid\'erons maintenant la repr\'esentation $\tau_r $ pour un certain entier $0 \leq r \leq n$. Puisque $\mathfrak{p}$
se d\'ecompose en sous-espaces de racines $\mathfrak{g}_{\alpha}$ de dimension $1$, les poids de $\tau_r$ sont toutes
les sommes possibles de $p$ racines distinctes non compactes $\alpha \in R_G \setminus R_K$. En particulier
$\mu_j = \sum_{\alpha \in (w_j . R_G^+ ) \setminus R_K^+} \alpha = w_j 2\delta_G -2\delta_K$ n'apparait que comme plus
haut poids d'une sous-repr\'esentation de $\tau_n$, et tous les autres plus hauts poids de $\tau_p$ appartiennent \`a 
$\mu_j - \sum_{\alpha \in R_G^+ } k_{\alpha} \alpha$, o\`u les $k_{\alpha}$ sont des entiers $\geq 0$.

Il d\'ecoule de tout ceci que $\pi_{w_j \lambda}$ et $\tau_r$ n'ont pas de K-type en commun, sauf lorsque 
$\lambda = \delta_G$ et $r=n$.

\bigskip

Si l'on note $\tau_{\mu_j}$ la repr\'esentation irr\'eductible de $K$ de plus haut
poids 
$$\mu_j = w_j  2\delta_G - 2 \delta_K .$$
Celle-ci intervient \`a la fois comme $K$-type minimal de $\pi_{w_j  \delta_G}$ et comme sous-repr\'esentation 
de la repr\'esentation $\tau_n$ de $K$. 

Concluons cette sous-section, en remarquant qu'un calcul simple montre que 
$$\mu_j = \sum_{k=1}^j \varepsilon_k - \sum_{k=j+1}^n \varepsilon_k + (n-2j) \varepsilon_{n+1}$$
et donc que 
$$\tau_{\mu_j} = \tau_{n-j,j} ' \; \mbox{  pour } 0\leq j \leq n.$$

\subsection{Caract\`eres infinit\'esimaux et action du Casimir}

Commen\c{c}ons par remarquer que le \cite[Th\'eor\`eme 9.20]{Knapp} affirme qu'une repr\'esentation 
$\pi_{w_j  \Lambda}$ dans la s\'erie discr\`ete de $G$ admet un caract\`ere infinit\'esimal \'egal \`a \index{caract\`ere infinit\'esimal}
$\Lambda$ (dans la param\'etrisation d'Harish-Chandra rappel\'ee au premier chapitre). En particulier, les
$\pi_{w_j  \delta_G}$ sont exactement les repr\'esentations de la s\'erie discr\`ete de $G$ qui ont un 
caract\`ere infinit\'esimal trivial. On a donc~:
$$\pi_{w_j  \delta_G} (C) =0,$$
o\`u $C$ d\'esigne le Casimir. 

La formule de Plancherel, d\'emontr\'ee par Harish-Chandra, sp\'ecialis\'ee dans le cas du groupe $G$ 
implique que le laplacien de Hodge-de Rham n'a pas de valeurs propres {\bf discr\`etes} autre 
que la valeur propre $0$ (avec multiplicit\'e infinie) sur les formes diff\'erentielles de degr\'e $n$, cf. \cite{Pedon}.
On peut en particulier en d\'eduire que le seul degr\'e o\`u le groupe de cohomologie $L^2$ r\'eduite de $X_G$ 
est non trivial est le degr\'e $n$. Ce th\'eor\`eme peut aussi se d\'emontrer \`a l'aide de m\'ethodes purement
g\'eom\'etriques comme celles d\'evelopp\'ees dans le Chapitre 14 de la deuxi\`eme partie, cette approche est 
due \`a Donnelly et Fefferman.

\medskip

Nous allons maintenant d\'eterminer les caract\`eres infinit\'esimaux des repr\'esentations induites du paraboliques
$P=MAN$. Soit $\mathfrak{t}$ la sous-alg\`ebre de Cartan de $\mathfrak{m}$ constitu\'ee des matrices diagonales.
D'apr\`es la \cite[Proposition 8.22]{Knapp}, si $\sigma$ admet un caract\`ere infinit\'esimal 
$\Lambda_{\sigma}$ la repr\'esentation $\pi_{\sigma , s}$ admet un caract\`ere 
infinit\'esimal \'egal \`a 
$$\Lambda_{\sigma} + s$$ \index{caract\`ere infinit\'esimal}
par rapport \`a la sous-alg\`ebre de Cartan $\mathfrak{a} \oplus \mathfrak{t}$.

Soit $\mu_{\sigma}$ le plus haut poids
de $\sigma$ et $2\delta_M = \sum_{k=2}^n (n+2-2k) \varepsilon_k$ la somme des racines dans $R_M^+$. 
Alors $\Lambda_{\sigma} = \mu_{\sigma} + \delta_M$. 

Le Lemme 12.28 de \cite{Knapp} implique que 
$$\pi_{\sigma ,s} (C) = (\langle \Lambda_{\sigma} +s , \Lambda_{\sigma} +s \rangle - \langle \delta_G ,\delta_G \rangle ) Id ,$$
o\`u $\langle .,. \rangle$ est le produit scalaire induit par la forme de Killing $B$ d\'efinie en (\ref{forme de killing}).
Un calcul facile montre alors que 
$$\pi_{\sigma ,s} (C) = - (n^2 -s^2 - \langle \mu_{\sigma} , \mu_{\sigma} + 2 \delta_M \rangle ) Id.$$
En particulier, \`a l'aide de (\ref{poids de sigma}), on obtient~: 
\begin{eqnarray} \label{valeur du casimir}
\pi_{\sigma_{p,q}  ,s} (C) = -((n-p-q)^2 -s^2 ) Id .    \index{action du Casimir}
\end{eqnarray}

\medskip

Remarquons (voir \cite{Pedon} pour plus de d\'etails) qu'\`a l'aide de la formule de Plancherel pour $G$, 
on peut alors d\'emontrer le th\'eor\`eme classique suivant.

\begin{thm} \label{spectre l2 de X}
Si pour $0\leq p+q \leq 2n$, on d\'esigne par $spec \Delta_{p,q}$ le spectre $L^2$ du laplacien de Hodge-de Rham
sur les formes diff\'erentielles de type $(p,q)$ sur l'espace hyperbolique complexe de dimension $n$. Alors,
$${\rm spec} \Delta_{p,q} = 
\left\{ 
\begin{array}{ll}
[(n-p-q)^2 , + \infty [  & \mbox{  si } p+q \neq n , \\
\{ 0 \} \cup [1 , +\infty [ & \mbox{  si } p+q = n .
\end{array} \right. $$
\end{thm} 

\bigskip

\`A l'aide des deux familles de repr\'esentations de $G$ que l'on a d\'ecrite ci-dessus, on peut d\'ecrire toutes les 
repr\'esentations admissibles de $G$, plus pr\'ecisemment on peut montrer (cf. \cite{Knapp}) :

\begin{thm}[Langlands, Knapp et Zuckerman] \label{classification de L pour unitaire}
Soit $\sigma$ une repr\'esentation irr\'eductible de $M$ et soit $s$ un nombre complexe de partie r\'eelle
$>0$. Alors la repr\'esentation induite $\pi_{\sigma , s }$ admet un unique quotient irr\'eductible,
son {\bf quotient de Langlands} : $J_{\sigma , s}$. \index{quotient de Langlands} \index{$J_{\sigma , s}$}

Chaque repr\'esentation irr\'eductible admissible de $G$ est soit une repr\'esentation de la s\'erie 
discr\`ete de $G$, soit une limite non d\'eg\'en\'er\'ee de s\'erie discr\`ete, soit une induite $\pi_{\sigma , s}$ o\`u $s \in i {\Bbb R}$ (et l'induite est irr\'eductible),
soit une repr\'esentation de la forme $J_{\sigma , s}$ comme ci-dessus. Et ces repr\'esentations sont deux \`a deux non \'equivalentes.
\end{thm}

\section{Groupe dual}

\markboth{CHAPITRE 4. REPR\'ESENTATIONS DE $U(n,1)$}{4.2. GROUPE DUAL}

Dans cette section nous expliquons la construction du {\bf groupe dual} de Langlands associ\'e \`a $G$. \index{groupe dual}

Tout d'abord, le groupe $G\times_{{\Bbb R}} {\Bbb C}$ obtenu \`a partir de $G$ par extension des scalaires \`a 
${\Bbb C}$ n'est autre que $GL(n+1,{\Bbb C})$. Son groupe dual est le groupe complexe $\widehat{G} = GL(n+1,{\Bbb C} )$,
nous le notons aussi $^L (G/{\Bbb C} )$. (Pour les motivations de cette d\'efinition voir \cite{Langlands}; voir aussi le cas des groupes
orthogonaux, Chapitre 6.)

Le groupe dual $^L G = {}^L (G/{\Bbb R} )$ du groupe {\bf r\'eel} $G$ est un produit semi-direct $\widehat{G} \rtimes {\rm Gal} ({\Bbb C}
/{\Bbb R} )$; $\widehat{G}$ est le groupe $GL(n+1, {\Bbb C})$, ${\rm Gal} ({\Bbb C} /{\Bbb R} )$ op\`ere par automorphismes 
{\bf holomorphes} ($\equiv$ alg\'ebriques) et la construction tient compte de la structure de $G$ comme groupe r\'eel.

Nous consid\'erons d'abord un groupe diff\'erent, le groupe $H=U(p+1,p)$ ($n=2p$ pair) ou $H=U(p+1, p+1)$ ($n=2p+1$
impair). Le groupe r\'eel $H$ est {\it quasi-d\'eploy\'e}, \index{quasi-d\'eploy\'e} c'est-\`a-dire qu'il existe un sous-groupe de Borel $B \subset 
H \times_{{\Bbb R}} {\Bbb C} = GL(n+1 , {\Bbb C} )$ d\'efini sur ${\Bbb R}$. Si par exemple $n$ est impair, la matrice
$$J= \left( 
\begin{array}{cccccc}
1 & & & & & \\
& \ddots & & & & \\
& & 1& & & \\
& & & -1 & & \\
& & & & \ddots & \\
& & & & & -1
\end{array} \right)$$
d\'efinissant la forme hermitienne est semblable \`a 
la matrice $J_*$ n'ayant que des $1$ sur l'anti-diagonale.
Le groupe $H=U(J)$ est donc isomorphe \`a 
$$U(J_* ) = \{ g\in GL(n+1 ,{\Bbb C} ) \; : \; {}^t \overline{g} 
J_* g = J_* \}.$$
La conjugaison complexe associ\'ee est 
$$\tau :  g \mapsto J_* {}^t \overline{g}^{-1} J_*  \; ;$$
$J_*$ \'etant antidiagonale, $\tau$ laisse globalement invariante le sous-groupe de Borel $B\subset 
GL(n+1 , {\Bbb C} )$ form\'e des matrices triangulaires sup\'erieures. M\^eme argument pour $n$ pair (avec la 
m\^eme matrice $J_*$).

Soit $R (H\times_{{\Bbb R}} {\Bbb C} , T)$ l'ensemble des racines positives de $GL(n+1 , {\Bbb C} )$ d\'efinies par $B$, 
et $\Delta$ l'ensemble des racines simples; $T=({\Bbb C}^* )^{n+1}$ est le tore maximal diagonal. La conjugaison complexe
$\tau$ op\`ere sur $R$ et $\Delta$ par ${}^{\tau} \alpha (t) = \overline{\alpha (\tau t)}$, $(t\in T)$. Si 
$\{ \alpha_1 , \ldots , \alpha_{n} \}$ sont les racines simples, $\tau$ op\`ere par $( \alpha_1 , \ldots , \alpha_n ) 
\mapsto ( \alpha_n^{-1} , \ldots , \alpha_1^{-1} )$.

Puisque $H \times_{{\Bbb R}} {\Bbb C} = G \times_{{\Bbb R}} {\Bbb C}$, le groupe $\widehat{H}$ est toujours $GL(n+1 , {\Bbb C})$.
On fait op\'erer la conjugaison complexe $\sigma \in {\rm Gal} ({\Bbb C} /{\Bbb R} )$ de fa\c{c}on holomorphe, sur $\widehat{H}$, 
de la fa\c{c}on suivante~:
\begin{description}
\item[(a)] $\sigma$ laisse globalement invariants $\widehat{T}$, $\widehat{B}$, o\`u $\widehat{T}$ est le tore diagonal, $\widehat{B}$ est 
le sous-groupe de Borel des matrices triangulaires sup\'erieures.
\item[(b)] Soient $\widehat{R}$, $\widehat{\Delta}$ d\'efinis comme ci-dessus mais relativement \`a $\widehat{T}$. Alors $\sigma$ op\`ere
par l'action pr\'ec\'edente sur $\widehat{R}$, $\widehat{\Delta}$.
\end{description}

La condition suivante est plus d\'elicate. Choisissons un {\it \'epinglage} \index{\'epinglage} de $\hat{B}$, c'est-\`a-dire un choix de matrices
nilpotentes dans $M_{n+1} ({\Bbb C}) = {\rm Lie} ( \widehat{G})$ associ\'ees aux racines simples. On choisira 
$$X_1 = \left( 
\begin{array}{cccc}
0 & 1 & &  \\
& 0 & & \\
& & \ddots & \\
& & & 0 
\end{array} \right) \; X_2 = \left(
\begin{array}{ccccc}
0 & & & & \\
& 0 & 1 &  & \\
& & 0 & & \\
& & & \ddots & \\
& & & & 0 
\end{array} \right) \; \ldots \; X_n = \left( 
\begin{array}{cccc}
0 & & & \\
& \ddots & & \\
& & 0 & 1\\
& & & 0
\end{array} \right) .$$
Alors
\begin{description}
\item[(c)] $\sigma$ pr\'eserve (globalement) $\{ X_1 , \ldots , X_n \}$.
\end{description}

On peut v\'erifier (cf. Borel \cite{Borel}) que $\sigma$ est uniquement d\'efini par le choix des $X_i$. Pour notre choix, 
$\sigma$ est alors de la forme suivante~:

\begin{lem} \label{action de galois}
Pour $g\in \widehat{G} = GL(n+1, {\Bbb C} )$, 
$$\sigma (g) = w_0 {}^t g^{-1} w_0^{-1}$$
o\`u 
$$w_0 = \left(
\begin{array}{ccccc}
& & & & (-1)^{n} \\
& & & \cdots & \\
& & 1 & & \\
& -1 & & & \\
1 & & & & 
\end{array} 
\right) .$$
\end{lem}

Noter que 
\begin{eqnarray} \label{w0 carre}
w_0^2 = (-1)^n .
\end{eqnarray}

Revenons enfin \`a $G$. Le groupe $G$ est {\it forme int\'erieure} \index{int\'erieure, forme} de $H$, {\it i.e.}, la conjugaison complexe 
$\sigma_G$ de $G({\Bbb C}) = GL(n+1 , {\Bbb C})$ d\'efinie par $G$ est conjugu\'ee \`a $\tau = \sigma_H$.
(Si $J_{1,n}$ est la matrice de la forme hermitienne d\'efinissant $G$, la premi\`ere est $g\mapsto J_{1,n} 
{}^t \overline{g}^{-1} J_{1,n}$ et la seconde est $g\mapsto J {}^t \overline{g}^{-1} J$.) Par 
d\'efinition (Langlands \cite{Langlands}, Borel \cite{Borel}), le groupe dual $^L G$ s'identifie \`a 
${}^L H= GL(n+1, {\Bbb C}) \rtimes {\rm Gal} ( {\Bbb C} / {\Bbb R} )$ o\`u $\sigma$ op\`ere comme dans le 
Lemme \ref{action de galois}. 

Nous aurons besoin de quelques notions simples relatives aux groupes duaux. Notons que $^L G$
est par construction un groupe r\'eductif complexe (non connexe). Un sous-groupe parabolique
$^L P$ de $^L G$ est un sous-groupe $\widehat{P} \rtimes {\rm Gal} ({\Bbb C}/{\Bbb R})$ o\`u
$\widehat{P} \subset \widehat{G}$ est un sous-groupe parabolique. Il revient au m\^eme de dire que c'est le normalisateur 
dans $^L G$ de $\widehat{P}$ o\`u $\widehat{P} \subset \widehat{G}$ est un sous-groupe parabolique globalement 
invariant par $\sigma$.

Rappelons qu'un sous-groupe parabolique $\widehat{P}$ de $\widehat{G}$ est conjugu\'e \`a un unique parabolique {\it
standard}, \index{sous-groupe parabolique standard} {\it i.e.}, contenant $\widehat{B}$. Il en r\'esulte qu'un sous-groupe parabolique $^L P$ est conjugu\'e
\`a un unique parabolique contenant $^L B= \widehat{B} \rtimes {\rm Gal} ({\Bbb C} /{\Bbb R} )$.

Enfin soit $P_0$ un sous-groupe parabolique de $G$ (donc d\'efini sur ${\Bbb R}$). Alors $P=P_0 \times_{{\Bbb R}} {\Bbb C}$
est un sous-groupe parabolique de $GL(n+1,{\Bbb C})$ d\'efini sur ${\Bbb R}$. La classe de conjugaison de $P$, en 
particulier, est d\'efinie sur ${\Bbb R}$. On en d\'eduit \cite{Borel} qu'on associe \`a $P$ un unique sous-groupe parabolique 
{\bf standard} $^L P$ contenant $^L B$. Pour notre groupe $G=U(n,1)$ on v\'erifie alors~:

\begin{lem} \label{paraboliques pertinents}
Les paraboliques standard $^L P \supset ^L B$ provenant de $G$ sont~:
\begin{enumerate}
\item Le groupe dual $^L G$.
\item Le normalisateur de $\widehat{P} \cong {\Bbb C}^* \times GL(n-1,{\Bbb C}) \times {\Bbb C}^* \ltimes \widehat{N}$ 
(o\`u $\widehat{N}$ est le radical unipotent).
\end{enumerate}
\end{lem}
\index{sous-groupe parabolique pertinent}

Dans (2) $\widehat{P}$ est l'ensemble des matrices triangulaires sup\'erieures par blocs de la taille indiqu\'ee. Un parabolique
$^L P$ {\bf provient de} $G$ si le parabolique standard associ\'e a cette propri\'et\'e.
Le Lemme se d\'eduit de Borel \cite[\S 3]{Borel}. Langlands appelle de tels paraboliques ``pertinents'' ({\it relevant}).

\markboth{CHAPITRE 4. REPR\'ESENTATIONS DE $U(n,1)$}{4.3. PARAM\`ETRES DE LANGLANDS}

\section{Param\`etres de Langlands}

Un {\it param\`etre de Langlands} \index{param\`etre de Langlands} pour $G$ est un homomorphisme de groupes $\varphi : W_{{\Bbb R}} \rightarrow
{}^L G$ tel que
\begin{eqnarray} \label{33a}
\begin{array}{c}
\mbox{Le morphisme } \varphi \mbox{ rend commutatif le diagramme} \\
\begin{array}{ccccc}
W_{{\Bbb R}} & & \stackrel{\varphi}{\rightarrow} & & ^L G \\
 & \searrow & & \swarrow & \\
& & {\rm Gal} ({\Bbb C} /{\Bbb R} ) && 
\end{array} 
\end{array}
\end{eqnarray}
\begin{eqnarray} \label{33b}
\begin{array}{l}
\mbox{L'image } \varphi (W_{{\Bbb C}} )= \varphi ({\Bbb C}^* ) \subset \widehat{G} = GL(n+1 , {\Bbb C}) \\
\mbox{ est form\'ee d'\'el\'ements semi-simples}.
\end{array}
\end{eqnarray}
\begin{eqnarray} \label{33c}
\begin{array}{l}
\mbox{Si } ^L P \subset {}^L G \mbox{ est un parabolique contenant } \varphi (W_{{\Bbb R}}), \\ 
^L P \mbox{ provient de } G.
\end{array}
\end{eqnarray}

Noter que (\ref{33b}) est \'equivalent \`a 
\begin{eqnarray} \label{33b'}
\varphi_{|W_{{\Bbb C}}} : W_{{\Bbb C}} \rightarrow GL (n+1,{\Bbb C} ) \mbox{ est semi-simple}.
\end{eqnarray}

On identifie deux param\`etres conjugu\'es par $\widehat{G}$.

Nous allons d\'ecrire tous les param\`etres de Langlands pour $G$. Il sera utile d'abord d'expliciter le changement 
de base dans cette situation. 

Soit $\varphi : W_{{\Bbb R}} \rightarrow {}^L G$ un param\`etre de Langlands. Alors $\varphi$ donne par restriction un 
param\`etre $\varphi_0 : W_{{\Bbb C}} = {\Bbb C}^* \rightarrow \widehat{G} = GL(n+1 , {\Bbb C} )$ qui d\'efinit donc (Chapitre 3)
une repr\'esentation de $GL(n+1 , {\Bbb C} )$.

\begin{lem} \label{ppte du chgt de base}
Si $\varphi_0 :  W_{{\Bbb C}} \rightarrow \widehat{G}$ provient par restriction d'un param\`etre pour $G$. 
\begin{enumerate}
\item La repr\'esentation $\varphi_0$ de ${\Bbb C}^*$ est isomorphe \`a $z\mapsto \widetilde{\varphi_0 (\overline{z} )} =
{}^t \varphi_0 (\overline{z} )^{-1} $.
\item La repr\'esentation $\pi_0$ de $GL(n+1 , {\Bbb C} )$ associ\'ee est isomorphe \`a sa conjugu\'ee $\pi_0 \circ
\sigma_G$, o\`u $\sigma_G$ est la conjugaison complexe d\'efinie par $G$.
\end{enumerate}
\end{lem}
{\it D\'emonstration.} La propri\'et\'e 2. se d\'eduit de 1. Rappelons que $\sigma_G$ et $\sigma_H$ sont conjugu\'es, $H$ \'etant la forme 
int\'erieure d\'ecrite dans la section pr\'ec\'edente. Mais $H$ a un sous-groupe de Borel, d\'efini par 
$$B_H = \{ g\in B = B_{GL(n+1,{\Bbb C} )} \; : \; {}^t \overline{g} J_* g = J_* \}$$
(section 4.2). Son tore maximal s'identifie \`a 
$$\{ z= (z_1 , \ldots , z_{n+1} ) \in {\Bbb C}^* \; : \; z_1 \overline{z}_{n+1} = z_2 \overline{z}_n = \ldots =1 \; (\mbox{et } 
z_{p+1} \overline{z}_{p+1} =1 \mbox{ si } n=2p+1) \} .$$
et la conjugaison complexe sur le tore maximal $T \cong ({\Bbb C}^* )^{n+1}$ de $B$ est donn\'ee par 
$(z_1 , \ldots , z_{n+1} ) \mapsto ( \overline{z}_{n+1}^{-1} , \ldots , \overline{z}_1^{-1} )$.
Si $\varphi_0 = \chi_1 \oplus \ldots \oplus \chi_{n+1}$ v\'erifie 1., on peut r\'eordonner les caract\`eres de sorte que 
$\chi_{n+1} (z) = \chi_1 (\overline{z} )^{-1} , \ldots , \chi_1 (z) = \chi_{n+1} (\overline{z} )^{-1}$. Alors $\pi_0 \cong 
\pi_0 \circ \sigma_G$ par transport de structure.

La propri\'et\'e 1. r\'esulte d'un calcul simple. \'Ecrivons $\varphi (j) = (h,\sigma ) \in \widehat{G} \times \{ \sigma \}$.
La d\'efinition de $^L G$ comme produit semi-direct donne 
$$(h,\sigma )^{-1} = (w_0 {}^t h w_0^{-1} , \sigma ) ; $$
on utilise les identit\'es
\begin{eqnarray} \label{w0 2}
{}^t w_0 = (-1)^n w_0 = w_0^{-1} .
\end{eqnarray}
Alors,
$$\varphi (\overline{z} ) =\varphi (jzj^{-1} ) = (hw_0 {}^t \varphi (z)^{-1} w_0^{-1} h^{-1} , 1)$$
comme il r\'esulte de l'expression du produit semi-direct et des \'egalit\'es pr\'ec\'edentes. D'o\`u le point 1.

\medskip

\begin{prop} \label{Lparam1}
Soit $\varphi$ une param\`etre de Langlands dont l'image n'est contenue dans aucun parabolique propre de $^L G$.
Alors $\varphi$ est \`a conjugaison pr\`es de la forme 
$$z \mapsto \left( (z/\overline{z} )^{p_1} , \ldots , (z/\overline{z} )^{p_{n+1}} \right)$$
o\`u $p_1 , \ldots , p_{n+1} \in \frac{n}{2} + {\Bbb Z}$, $p_i \neq p_j$; 
$$j \mapsto (w_0^{-1} = (-1)^n w_0 , \sigma ).$$
\end{prop}
{\it D\'emonstration.} Soit $\varphi_0 = \varphi_{|{\Bbb C}^*}$. On peut supposer que $\varphi_0 ({\Bbb C}^* ) \subset \widehat{T}$.
D'apr\`es le premier point du Lemme \ref{ppte du chgt de base} on peut supposer 
$$\varphi_0 (z) = (\chi_1 (z) , \ldots , \chi_{n+1} (z) )$$
o\`u $\chi_{n+2-i} (z) = \chi_i (\overline{z} )^{-1}$. Notons comme ci-dessus $\varphi (j) = (h, \sigma )$, et consid\'erons 
la restriction de $\varphi_0$ \`a ${\Bbb R}_+^*$. Nous notons $u$ un \'el\'ement de $\widehat{T} = ({\Bbb C}^* )^{n+1}$.
Quitte \`a r\'eordonner les $\chi_i$, on peut supposer que pour $x\in {\Bbb R}_+^*$
$$\varphi_0 (x) =(x^{s_1} , x^{s_1} , \ldots , x^{s_2} , x^{s_2} , \ldots  ,\ldots  x^{s_r} , x^{-s_r} , \ldots , x^{-s_1} ) =: u$$
o\`u $s_1 > s_2 > \ldots > s_r \geq 0$ (si $s_r$ est nul, on ne s\'epare pas les occurences de $x^{s_r} =1$ et de 
$x^{-s_r}$).

On a $\sigma u \sigma^{-1}  = w_0 {}^t u^{-1} w_0^{-1} = w_0 u^{-1} w_0^{-1} = u$. La conjugaison par 
$(1, \sigma ) \in ^L G$ fixe donc $u$. Puisque $x$ est r\'eel, $\varphi (j) \varphi_0 (x) \varphi (j)^{-1} = \varphi_0 (x)$.
Donc $(h,1)$ commute \`a $u$. Si $r>1$, il en r\'esulte que $h$ appartient au parabolique propre de 
$\widehat{G}$ de type $(n_1 , n_2 , \ldots , n_r , n_r , \ldots , n_1)$ (un seul bloc m\'edian si $s_r =0$). Puisque celui-ci est 
normalis\'e par l'action de $\sigma$, on voit que $\varphi_0 ({\Bbb C}^* )$, et $\varphi (j)$ appartiennent \`a un parabolique
propre de $^L G$. On en d\'eduit donc que tous les caract\`eres $\chi_i$ sont de la forme $(z/ \overline{z} )^{p_i}$ avec 
$p_i \in \frac12 {\Bbb Z}$.

Soit alors $z\in {\Bbb C}^*$ et $u=\varphi_0 (z)$. Puisque ${}^t u^{-1} = u^{-1}$, l'argument pr\'ec\'edent montre
que ${\rm Ad} (h w_0 )u^{-1} = u^{-1}$. Supposons par exemple que $p_1 = p_2 = \ldots =p_r$, les autres $p_i$ \'etant diff\'erents.
On peut supposer $\varphi_0$ de la forme $\left( \ldots , (z /\overline{z} )^{p_1} , \ldots , (z / \overline{z} )^{p_1 } , \ldots \right)$
les occurences de $(z/\overline{z} )^{p_1}$ \'etant sym\'etriques par rapport \`a $\frac{n+2}{2}$. Alors ${\rm Ad} (h w_0 )$ normalise 
un parabolique de $\widehat{G}$ qui s'\'etend en un parabolique de $^L G$, et on en d\'eduit de nouveau que $\varphi$ passe 
par un parabolique propre.

V\'erifions la condition de parit\'e sur les $p_i$. Rappelons qu'avec la relation de Langlands pour les caract\`eres complexes
(Chapitre 3) $(z/\overline{z} )^p$ ($p\in \frac12 {\Bbb Z}$) est \'egal \`a $(-1)^{2p}$ pour $z=-1$, et que $\varphi (j) = (h,\sigma )$ op\`ere sur
$\widehat{T}$ par ${\rm Ad} (hw_0 )$. Soit $u=\varphi_0 (z) =
\left( (z/\overline{z} )^{p_1} , \ldots , (z/\overline{z} )^{p_{n+1}} \right)$. On a vu que ${\rm Ad} (hw_0 ) u=u$; puisque $\varphi_0$ est un caract\`ere r\'egulier, ceci 
implique maintenant que ${\rm Ad}(hw_0 )u=u$ pour tout $u\in \widehat{T}$, donc $hw_0 \in \widehat{T}$; \'ecrivons $h= vw_0^{-1}$ o\`u
$v\in \widehat{T}$. Alors $\varphi (j) = (vw_0^{-1} , \sigma )$ et 
\begin{eqnarray*}
\varphi (j^2 ) & = & (vw_0^{-1} , \sigma )(v w_0^{-1} , \sigma ) = (vw_0^{-1} w_0 v {}^t w_0 w_0^{-1} , 1) \\
 & = & ((-1)^n , 1)
\end{eqnarray*}
- on utilise de nouveau (\ref{w0 2}). Mais $j^2 =-1$ et la remarque pr\'ec\'edente implique $p_i \equiv \frac{n}{2} \; (\mbox{mod }1)$.

Enfin l'\'egalit\'e, vraie pour $u\in \widehat{T}$~:
$$(u,1)(w_0^{-1} , \sigma ) (u^{-1} ,1) = (u,1) (w_0^{-1} w_0 u w_0^{-1} , \sigma ) = (u^2 w_0^{-1} , \sigma )$$
implique qu'\`a conjugaison pr\`es (par $\widehat{T}$, donc sans changer $\varphi_0$) $\varphi (j)$ est de la forme indiqu\'ee.

\medskip

La proposition suivante compl\`ete la description des param\`etres de Langlands pour $G$.

\begin{prop} \label{Lparam2}
Soit $\varphi$ un param\`etre de Langlands pour $G$ dont l'image est contenue dans le parabolique de type $(1,n-1,1)$.
Alors $\varphi$ est \`a conjugaison pr\`es de la forme 
$$z\mapsto \left( \chi (z) , (z/\overline{z})^{p_1} , \ldots , (z/\overline{z} )^{p_{n-1}} , \chi (\overline{z} )^{-1} \right)$$
$$j \mapsto ( \left( 
\begin{array}{ccc}
\varepsilon_1 & & \\
& w_0 & \\
& & \varepsilon_2 
\end{array} \right) , \sigma ),  \; \; \varepsilon_1 \varepsilon_2 = \chi (-1) $$
o\`u $\chi$ est un caract\`ere arbitraire de ${\Bbb C}^*$ et les $p_i$ sont tous disjoints et appartiennent \`a 
$\frac{n}{2} + {\Bbb Z}$.
\end{prop}
{\it D\'emonstration.} En effet l'image de $\varphi_0$ est contenue \`a conjugaison pr\`es dans ${\Bbb C}^* \times
GL(n-1, {\Bbb C} ) \times {\Bbb C}^*$, sur lequel $\sigma$ op\`ere de fa\c{c}on analogue \`a l'action pr\'ec\'edente
(sur le bloc m\'edian) et en permutant les deux facteurs ${\Bbb C}^*$. On est alors ramen\'e au cas pr\'ec\'edent. Les
d\'etails sont laiss\'es au lecteur.

\bigskip

On dira d'un param\`etre de Langlands $\varphi$ qu'il est {\it discret} \index{param\`etre discret}
s'il est comme dans la Proposition \ref{Lparam1}. Et on dira d'un param\`etre de Langlands qu'il est {\it temp\'er\'e} \index{param\`etre temp\'er\'e}
s'il est born\'e. Remarquons que les param\`etres de Langlands temp\'er\'es sont 
\begin{enumerate}
\item les param\`etres discrets et 
\item les param\`etres de la Proposition \ref{Lparam2} avec $\chi$ unitaire.
\end{enumerate}

\markboth{CHAPITRE 4. REPR\'ESENTATIONS DE $U(n,1)$}{4.4. CLASSIFICATION ET PARAM\`ETRES DE LANGLANDS}

\section{Classification et param\`etres de Langlands}

La classification de la s\'erie discr\`ete de $G$ a \'et\'e rappel\'e au \S 4.1.2. On va en d\'eduire~:

\index{classification de Langlands, Th\'eor\`eme de}
\begin{thm}[Langlands] \label{thm de L}
\begin{enumerate}
\item \`A tout param\`etre de Langlands $\varphi$ pour $G$ est associ\'e un ensemble fini $\Pi (\varphi )$ de 
repr\'esentations admissibles irr\'eductibles de $G$.
\item Les $\Pi (\varphi )$ correspondant \`a tous les param\`etres (\`a conjugaison pr\`es) sont disjoints, et leur
r\'eunion est le dual admissible de $G$.
\end{enumerate}
\end{thm}
Soit en effet $\varphi$ un param\`etre discret (Proposition \ref{Lparam1}). Notons $T_c$
le tore diagonal $U(1)^{n+1}$ de $G$. Soit $\Lambda$ le r\'eseau des caract\`eres de $T_c$; donc $\Lambda =
{\Bbb Z}^{n+1}$; on peut consid\'erer $(p_1 , \ldots , p_{n+1} )$ comme un \'el\'ement de $\frac12 \Lambda$.
De plus $\delta_G = \left( \frac{n}{2} , \frac{n}{2} -1 , \ldots , -\frac{n}{2} \right)$ (cf. (\ref{delta G})); toujours d'apr\`es la 
Proposition \ref{Lparam1}, $p\in \Lambda +\delta_G$. D'apr\`es le \S 4.1.2, on peut lui associer une repr\'esentation de la s\'erie
discr\`ete.

Par ailleurs $(p_1 , \ldots , p_{n+1} )$ n'est d\'efini qu'\`a l'ordre pr\`es, donc modulo le groupe de Weyl $W_G =
\Sigma_{n+1}$; et deux choix d'ordre conjugu\'es par $W_K = \Sigma_n$ d\'efinissent la m\^eme repr\'esentation.
On associe donc dans ce cas \`a $\varphi$ l'ensemble de $(n+1)$ repr\'esentations correspondantes \`a l'orbite de 
$W_G$ (modulo conjugaison par $W_K$). C'est le {\it $L$-paquet} \index{$L$-paquet} associ\'e \`a $\varphi$.

Si au contraire $\varphi$ est comme dans la Proposition \ref{Lparam2}, il d\'efinit de m\^eme une repr\'esentation 
$\tau$ du groupe $U(n-1)$ dont le param\`etre d'Harish-Chandra est $(p_1 , \ldots , p_{n-1} )$. Choisissons $\chi$
(en le rempla\c{c}ant au besoin par $\chi (\overline{z} )^{-1}$) de la forme 
$$\chi (z) = z^{\alpha} (\overline{z})^{\beta} \mbox{ avec Re} (\alpha +\beta )\geq 0.$$

Le groupe $U(n-1) \times {\Bbb C}^*$ est le sous-groupe de Levi du parabolique minimal $P$ de $G$. On distingue alors
deux cas~:
\begin{enumerate}
\item Si Re$(\alpha +\beta )>0$, la repr\'esentation ind$_P^G (\tau \otimes \chi )$ a un unique quotient irr\'eductible, son
quotient de Langlands (section 4.1) $J(\tau , \chi )$. On d\'efinit alors $\Pi (\varphi ) = \{ J(\tau , \chi ) \}$.
\item Supposons Re$(\alpha +\beta )=0$ : $\chi$ est donc unitaire, ainsi que la repr\'esentation induite. Celle-ci est en 
fait de longueur $1$ ou $2$ selon des r\'esultats g\'en\'eraux de Knapp et Zuckerman \cite{KnappZuckerman} et ses 
composantes sont de multiplicit\'e $1$. L'ensemble $\Pi (\varphi )$, de cardinal $1$ ou $2$, est l'ensemble de ces 
composantes.
\end{enumerate}

\bigskip 

\noindent
{\bf Remarque.}  Dans le cas 2, et si ind$_P^G (\tau \otimes \chi )$ est r\'eductible, elle se scinde en deux
repr\'esentations qui sont des limites de s\'eries discr\`etes et qui ont une description analogue aux repr\'esentations de la s\'erie discr\`ete
(\S 4.1.2). Nous ne donnons pas de description plus explicite de celles-ci, pour la raison suivante. Dans la section suivante,
nous calculons la valeur propre de l'op\'erateur de Casimir dans une repr\'esentation induite contenant un $K$-type 
$\Lambda^p \mathfrak{p}^+ \otimes \Lambda^q \mathfrak{p}^-$ associ\'e aux formes diff\'erentielles. Dans le Chapitre 8
nous v\'erifions pour les induites {\bf unitaires} que ces valeurs propres v\'erifient trivialement les bornes de la Conjecture
A(1). Par cons\'equent, pour d\'emontrer celle-ci, nous n'avons pas \`a nous soucier des induites unitaires ou de leurs
sous-modules.

\bigskip

D\'ecrivons maintenant les caract\`eres infinit\'esimaux des repr\'esentations obtenues. Dans le cas des s\'eries discr\`etes
le caract\`ere infinit\'esimal, d'apr\`es Harish-Chandra est simplement $\lambda = (p_1 , \ldots , p_{n+1} ) \in {\Bbb C}^n 
\cong \mathfrak{a}+\mathfrak{t} $.

Consid\'erons maintenant les repr\'esentations induites; soit $\chi (z) = z^{\alpha} (\overline{z})^{\beta}$, avec 
Re$(\alpha +\beta )\geq 0$. Soit $M=U(n-1) \times U(1)$ le groupe d\'efini au \S 4.1.3 et $\sigma$ la repr\'esentation
induisante de $M$. Le caract\`ere infinit\'esimal de $\sigma$ est alors \index{caract\`ere infinit\'esimal}
$$\lambda_{\sigma} = (p_1 , \ldots , p_{n-1} , \alpha -\beta )$$
et son plus haut poids
$$\mu_{\sigma} =  \lambda_{\sigma} - \delta_M = \left( p_1 - \frac{n-2}{2} , p_2 - \frac{n-4}{2} , \ldots , p_{n-1} +\frac{n-2}{2} , \alpha -\beta \right)$$
o\`u on a ordonn\'e les $p_i$ par $p_1 > p_2 > \ldots > p_n$.

\`A conjugaison pr\`es par $G$, l'alg\`ebre de Lie $\mathfrak{a}+\mathfrak{t}$ s'identifie aux matrices de la forme
$$H = \left( 
\begin{array}{ccccc}
X+iY & & & & \\
& i X_1 & & & \\
& & \ddots & & \\
& & & iX_{n-1} & \\
& & & & -X+iY 
\end{array} \right) .$$

La valeur du caract\`ere infinit\'esimal sur une telle matrice est alors
$$\lambda_{\pi} (H) = i\sum_{j=1}^{n-1} p_j X_j + i(\alpha - \beta )Y + sX$$
o\`u $s=\alpha +\beta$. Pour la forme duale \`a la forme de Killing (\ref{forme de killing}), on a alors
\begin{eqnarray} \label{421}
\langle \lambda_{\pi} , \lambda_{\pi} \rangle = s^2 + (\alpha -\beta )^2 +2\sum_{j=1}^{n-1} p_j^2 .
\end{eqnarray}

\bigskip

Remarquons que les param\`etres de Langlands discrets correspondent exactement aux $L$-paquets contenant une 
repr\'esentation discr\`ete (et dans ce cas toutes les repr\'esentations dans le $L$-paquet sont discr\`etes).
De m\^eme, les param\`etres de Langlands temp\'er\'es correspondent exactement aux $L$-paquets contenant une 
repr\'esentation temp\'er\'ee (et dans ce cas toutes les repr\'esentations dans le $L$-paquet sont temp\'er\'ees).

\bigskip

On peut maintenant poser la question de la fonctorialit\'e. Un param\`etre de Langlands 
$\varphi : W_{\mathbb{R}} \rightarrow {}^L G$ pour $G$ induit un param\`etre de Langlands 
$\varphi_{|\mathbb{C}^*} : \mathbb{C}^* \rightarrow GL(n+1 , \mathbb{C} )$ pour $GL(n+1, \mathbb{C})$.
Or les param\`etres de Langlands de $\mathbb{C}^*$ dans $GL(n+1, \mathbb{C} )$ param\`etrent les repr\'esentations 
admissibles. On obtient donc une application naturelle $\phi$ de l'ensemble des $L$-paquets de $G$ 
dans l'ensemble des repr\'esentations admissibles de $GL(n+1, \mathbb{C} )$. Et le probl\`eme de fonctorialit\'e, ici de changment de base,
s'\'enonce comme suit.

\medskip
\noindent
{\bf Question (Fonctorialit\'e)} {\it L'application $\phi$ envoie-t-elle tout $L$-paquet contenant une repr\'esentation
automorphe vers une repr\'esentation automorphe de $GL(n+1,\mathbb{C})$ ?}
\index{fonctorialit\'e} \index{changement de base}

\medskip

Une r\'eponse positive \`a cette question permettrait d'exploiter toute approximation de la Conjecture de Ramanujan 
sur $GL(n+1, \mathbb{C})$. 
N\'eanmoins, la r\'eponse \`a la question ci-dessus est fausse en general comme nous le verrons notamment 
au Chapitre 6 (Th\'eor\`eme 6.51).

Malgr\'e cela, on peut esp\'erer que l'application $\phi$ envoie de nombreux $L$-paquets automorphes de $G$ sur 
des repr\'esentations automorphes de $GL(n+1,\mathbb{C})$. On reviendra sur cette question au Chapitre 8.

\bigskip

Enfin, et bien que nous n'en aurons pas besoin, remarquons que l'on conna\^{\i}t le dual unitaire de $G$ (cf. \cite{Krajlevic}). 
Premi\`erement  tout $L$-paquet contenant une repr\'esentation unitaire n'est constitu\'e que de repr\'esentations 
unitaires. On peut donc parler de {\it $L$-paquet unitaire} \index{$L$-paquet unitaire}. La description du dual unitaire de $G$ est difficile, nous utiliserons
au Chapitre 6 la description qu'en font Knapp et Speh dans \cite{KnappSpeh}. Il en r\'esulte en particulier que les $L$-paquets 
$\Pi_{\varphi}$ o\`u le param\`etre $\varphi$ est 
\begin{itemize}
\item temp\'er\'e,  ou 
\item contient la repr\'esentation $J_{\sigma_{a,b} , s}$ pour $0<s \leq \frac{n-(a+b)}{2}$, 
\end{itemize}
sont unitaires. {\it A contrario}, si $s> \frac{n-(a+b}{2}$, la repr\'esentation $J_{\sigma_{a,b} , s}$ n'est pas unitaire. Nous donnons la 
description compl\`ete du dual unitaire de $U(2,1)$ au \S 4.6.

\section{$K$-types des repr\'esentations induites, action du Casimir}

\markboth{CHAPITRE 4. REPR\'ESENTATIONS DE $U(n,1)$}{4.5. $K$-TYPES, ACTION DU CASIMIR}

Soit $\varphi$ un param\`etre de Langlands d\'efinissant une repr\'esentation induite, celle-ci est donc param\'etr\'ee par
les donn\'ees $(p_j )_{j\leq n-1}$, $\alpha$, $\beta$. Soit $I_{\varphi}$ la repr\'esentation induite et $J_{\varphi}$
son quotient de Langlands. La repr\'esentation $I_{\varphi}$ ne nous int\'eresse
que si 
$$\mbox{Hom}_K (I_{\varphi} , \Lambda^p \mathfrak{p}_+ \otimes \Lambda^q \mathfrak{p}_- ) 
= \mbox{Hom}_M (\sigma , \Lambda^p \mathfrak{p}_+ \otimes \Lambda^q \mathfrak{p}_- ) \neq 0 .$$
D'apr\`es le \S 4.1.1 il existe alors un entier $k\in [0, \mbox{min} (p,q)]$ tel que 
\begin{eqnarray*}
\mu_{\sigma} & = & (p_1 -\frac{n-2}{2} , \ldots , p_{n-1} + \frac{n-2}{2} , \alpha -\beta ) \\
                      & (= & (\mu_1 , \ldots , \mu_{n-1} , \alpha - \beta ) )
\end{eqnarray*} 
soit de la forme
$$( \underbrace{1, \ldots 1}_{b \; {\rm termes}} , 0 \ldots 0, \underbrace{-1 , \ldots  -1}_{a \; {\rm termes}}, a-b )$$
et $(a,b)$ de la forme $(\pi , \rho )$, $(\pi -1 , \rho )$, $(\pi , \rho -1)$ ou $(\pi -1 , \rho -1)$ avec 
$$(\pi , \rho ) = (p-k , q-k).$$
On a donc 
$$\alpha - \beta = a-b .$$
En utilisant les indentit\'es 
$$p_{\sigma} = \mu +\delta_K \; \; (\mbox{o\`u } p=(p_1 , \ldots p_{n-1} , \alpha -\beta ) )$$
et
$$\langle \delta_K , \delta_K \rangle - \langle \delta_G , \delta_G \rangle =  -n^2$$
un calcul simple donne alors~:

\begin{prop} \label{K type de Jphi} 
La repr\'esentation $J_{\varphi}$ intervient dans le spectre des $(p,q)$-formes si et seulement s'il existe 
un entier $k \in [0, \min (p,q)  ]$ tel que 
\begin{enumerate} 
\item $\mu_1=   \ldots =\mu_{q-k} =1$, $\mu_{q-k+1} = \ldots =\mu_{n-p+k-1} = 0$, 
$\mu_{n-p+k} = \ldots = \mu_{n-1} =-1$ et $\alpha - \beta = p-q$ et dans ce cas 
$$J_{\varphi} (C) = -((n-p-q+2k)^2 - (\alpha + \beta )^2 ),$$
ou bien
\item  $k \leq p-1$, $\mu_1=  \ldots = \mu_{q-k} =1$, $\mu_{q-k+1} = \ldots = \mu_{n-p+k} =0$, 
$\mu_{n-p+k+1} = \ldots = \mu_{n-1} =-1$ et $\alpha - \beta = (p-q-1)$ et dans ce cas 
$$J_{\varphi} (C) = -((n-p-q+2k+1)^2 - (\alpha + \beta )^2 ),$$
ou bien
\item  $k \leq q-1$, $\mu_1= \ldots = \mu_{q-k-1} =1$, $\mu_{q-k} = \ldots = \mu_{n-p+k-1} =0$, 
$\mu_{n-p+k} = \ldots = \mu_{n-1} =-1$ et $\alpha - \beta = (p-q+1)$ et dans ce cas
$$J_{\varphi} (C) =-((n-p-q+2k+1)^2 - (\alpha +\beta )^2 ),$$
ou bien
\item  $k \leq p-1$, $k\leq q-1$, $\mu_1= \ldots = \mu_{q-k-1} =1$, $\mu_{q-k} = \ldots = \mu_{n-p+k} =0$, 
$\mu_{n-p+k+1} = \ldots = \mu_{n-1} =-1$ et $\alpha - \beta = (p-q)$ et dans ce cas
$$J_{\varphi} (C) = -((n-p-q+2k+2)^2 - (\alpha +\beta )^2 ).$$
\end{enumerate}
\end{prop} \index{$K$-type} \index{action du Casimir}

\markboth{CHAPITRE 4. REPR\'ESENTATIONS DE $U(n,1)$}{4.6. REPR\'ESENTATIONS DE $U(2,1)$}

\section{Repr\'esentations de $U(2,1)$}

Dans cette section, nous d\'etaillons la classification des repr\'esentations de $U(2,1)$, le cas dont nous aurons le 
plus besoin dans la suite. Dans cette section $G=U(2,1)$. Si $\chi \in {\rm Hom} (MA, {\Bbb C}^* )$, il existe $u$, $v \in 
{\Bbb C}^*$ uniques v\'erifiant $u-v \in {\Bbb Z}$ et un unique entier $r\in {\Bbb Z}$ tels que~:
$$\chi \left( \left(
\begin{array}{ccc}
e^t e^{i\theta} & & \\
& e^{i\eta } & \\
& & e^{-t} e^{i\theta }
\end{array} \right) \right) = e^{(u+v)t} e^{i(u-v)\theta} e^{ir(\eta+2\theta)} .$$
Nous noterons $\chi = (u,v,r)$. \footnote{La param\'etrisation que l'on adopte ici est l\'eg\`erement diff\'erente de celle des sections pr\'ec\'edentes. Nous adoptons en effet la param\'etrisation 
de Rogawski dans \cite{Rogawski} dont nous utiliserons les r\'esultats par la suite. Le lecteur constatera qu'il n'est n\'eanmoins pas difficile de se raccrocher aux r\'esultats des sections pr\'ec\'edentes.} 

\medskip

Cette fois l'ensemble des morphismes admissibles se scinde en deux familles (Propositions \ref{Lparam1} et \ref{Lparam2}).
\begin{enumerate}
\item Les morphismes de la forme :
$$\varphi (z) = \left( 
\begin{array}{ccc} 
\left( \frac{z}{\overline{z}} \right)^a &&\\
& \left( \frac{z}{\overline{z}} \right)^b & \\
&& \left( \frac{z}{\overline{z}} \right)^c 
\end{array} \right) $$
et 
$$\varphi (j) = \left( \left( 
\begin{array}{ccc}
&& 1 \\
& -1 & \\
1 && 
\end{array}
\right) , \sigma \right) $$
avec $a,b,c \in \mathbb{Z}$ et $a\geq b \geq c$.
\item Les morphismes de la forme :
$$\varphi (z) = \left( 
\begin{array}{ccc}
z^u \overline{z}^v & & \\
& \left( \frac{z}{\overline{z}} \right)^{\mu} & \\
&& z^{-v} \overline{z}^{-u} 
\end{array}
\right) $$
et 
$$\varphi (j) = \left( \left( 
\begin{array}{ccc}
1 &&\\
& 1 &\\
&& (-1)^{u-v} 
\end{array}
\right) , \sigma \right) $$
avec $\mu \in \mathbb{Z}$, $u,v \in \mathbb{C}$, $u-v \in \mathbb{Z}$, Re$(u+v) \geq 0$ et si $u+v =0$ alors 
$u \in \frac{1}{2} + \mathbb{Z}$. \footnote{Le cas o\`u $u\in {\Bbb Z}$ est d\'ej\`a d\'ecrit dans 1.}
\end{enumerate}

\medskip

D\'ecrivons la correspondance de Langlands. 
\begin{description}
\item[Premier cas :] Le param\`etre $\varphi$ est dans la premi\`ere famille. 
\begin{description}
\item[Premier sous-cas :] $a>b>c$. Alors (nous adoptons ici les notations de \cite{Rogawski}) :
$$\Pi_{\varphi} = \{ D_{\varphi} , D_{\varphi}^+ , D_{\varphi}^- \} $$
est le $L$-paquet {\bf discret} (et donc {\bf unitaire}) constitu\'e des repr\'esentations de la s\'erie discr\`ete 
ayant pour caract\`ere infinit\'esimal 
$$\chi_{\varphi} = (a-b,b-c,b ) .$$
\item[Deuxi\`eme sous-cas :] $a>b=c$. Alors :
$$\Pi_{\varphi} = \{ \pi_{\varphi}^1 , \pi_{\varphi}^2 \} $$
est le $L$-paquet {\bf temp\'er\'e} (et donc {\bf unitaire}) compos\'e des constituants irr\'eductibles de 
l'induite unitaire du caract\`ere 
$$\chi_{\varphi}^+ =(b-a , a-c , a) .$$
On a num\'erot\'e ces constituants de fa\c{c}on \`a ce que $\pi_{\varphi}^1$ soit un constituant de l'induite 
unitaire du caract\`ere 
$$\chi_{\varphi}^- = (a-c , c-b , c) .$$
\item[Troisi\`eme sous-cas :] $a=b>c$. Alors :
$$\Pi_{\varphi} =\{ \pi_{\varphi}^1 , \pi_{\varphi}^2 \} $$
est le $L$-paquet {\bf temp\'er\'e} (et donc {\bf unitaire}) compos\'e des constituants irr\'eductibles de l'induite 
unitaire du caract\`ere $\chi_{\varphi}^-$ et num\'erot\'es de fa\c{c}on \`a ce que $\pi_{\varphi}^1$ soit
un constituant de l'induite unitaire du caract\`ere $\chi_{\varphi}^+$.
\item[Quatri\`eme sous-cas :] $a=b=c$. Alors :
$$\Pi_{\varphi} =\{ i_G (\chi_{\varphi} ) \} $$
est le $L$-paquet {\bf temp\'er\'e} (et donc {\bf unitaire}) constitu\'e de l'induite unitaire (qui est irr\'eductible) du 
caract\`ere $\chi_{\varphi}$.
\end{description}
\item[Deuxi\`eme cas :] Le param\`etre $\varphi$ est dans la seconde famille.
\begin{description}
\item[Premier sous-cas :] $u,v \in \mathbb{Z}$ et $u\leq 0 < v$. Alors 
l'induite unitaire du caract\`ere 
$$\chi_{\varphi} =(u,v, \mu ) $$
a un unique quotient irr\'eductible que l'on note $J_{\varphi}$ 
(attention ce n'est pas tout \`a fait la notation de \cite{Rogawski}) et :
$$\Pi_{\varphi} = \{ J_{\varphi} \}.$$
Ce $L$-paquet n'est {\bf pas temp\'er\'e} et est {\bf unitaire ssi $u+v=1$}.
\item[Deuxi\`eme sous-cas :] $u,v \in \mathbb{Z}$ et $v\leq 0<u$. Alors 
l'induite unitaire du caract\`ere $\chi_{\varphi}$ a un unique quotient irr\'eductible que l'on note
$J_{\varphi}$ et :
$$\Pi_{\varphi} = \{ J_{\varphi} \} .$$
Ce $L$-paquet n'est {\bf pas temp\'er\'e} et est {\bf unitaire ssi $u+v=1$}.
\item[Troisi\`eme sous-cas :] $u,v \in \mathbb{Z}$ et $u,v >0$. Alors 
l'induite unitaire du caract\`ere $\chi_{\varphi}$ a un unique quotient irr\'eductible que l'on note 
$F_{\varphi}$ et :
$$\Pi_{\varphi} = \{ F_{\varphi} \} .$$
Ce $L$-paquet n'est {\bf pas temp\'er\'e} et {\bf de dimension finie}. (La repr\'esentation $F_{\varphi}$ est triviale si et seulement si $(u,v,\mu )= (1,1, 0)$.)
\item[Quatri\`eme sous-cas :] $(u,v) \notin \mathbb{Z}^2$. Alors 
l'induite unitaire du caract\`ere $\chi_{\varphi}$ est irr\'eductible et :
$$\Pi_{\varphi} =\{ i_G (\chi_{\varphi} ) \} .$$
Ce $L$-paquet est {\bf temp\'er\'e ssi Re$(u+v) =0$} et {\bf unitaire ssi (Re$(u+v)=0$) ou ($u=v$ et $|u+v|<2$) ou ($u-v$ impair 
et $|u+v| <1$)}.
\end{description}
\end{description}

\bigskip

Concluons cette section par la description d'une autre fonctorialit\'e reliant d'une part 
certaines repr\'esentations de $U(1,1) $ et de $GL(2,{\Bbb C})$ et d'autre part les repr\'esentations
de $U(1,1) \times U(1)$, et celles de $U(2,1)$.

Le groupe de Weil $W_{\mathbb{R}}$ se projette sur le groupe de Galois de $\mathbb{C}$ sur $\mathbb{R}$
et, via cette projection, agit sur :
$$\widehat{U(1,1)} = GL(2,{\Bbb C})$$
et donc sur 
$$\widehat{U(1,1) \times U(1)} = GL(2,\mathbb{C} ) \times GL(1,\mathbb{C} ) $$ 
et sur
$$\widehat{G} = GL(3,\mathbb{C} ).$$
Pour chaque entier $n$, on peut alors former les morphismes injectifs qui suivent~:
\begin{eqnarray} \label{betan}
\beta_n : 
GL(2,{\Bbb C}) \rtimes W_{{\Bbb R}} \rightarrow  GL(2,{\Bbb C}),  
\end{eqnarray}
dont la restriction \`a $GL(2,{\Bbb C})$ est l'identit\'e et la restriction au groupe 
de Weil $W_{\mathbb{R}}$ est donn\'ee par :
$$\beta_n (z) = \left( 
\begin{array}{cc}
\left( \frac{z}{\overline{z}} \right)^{\frac12 +n} &  \\
 & \left( \frac{z}{\overline{z}} \right)^{\frac12 +n} 
\end{array} \right) \times z \mbox{  si } z\in \mathbb{C}^* $$ \index{$\beta_n$}
et 
$$\beta_n (j) = \left( 
\begin{array}{cc}
1 &  \\
& -1 
\end{array} \right) \times j ;$$
et par ailleurs 
\begin{eqnarray} \label{xin} \index{$\xi_n$}
\xi_n : \widehat{H} \rtimes W_{\mathbb{R}} \rightarrow \widehat{G} \rtimes W_{\mathbb{R}}
\end{eqnarray}
o\`u~: 
\begin{enumerate}
\item la restriction de $\xi_n$ \`a $\widehat{H}$ est donn\'ee par le morphisme :
$$\begin{array}{ccc}
GL(2,\mathbb{C} ) \times GL(1, \mathbb{C} ) & \rightarrow & GL(3,\mathbb{C} ) \\
(\left( 
\begin{array}{cc}
a & b \\
c & d 
\end{array} \right) , \alpha ) & \mapsto & 
\left(
\begin{array}{ccc} 
a & 0 & b \\
0 & \alpha & 0 \\
c & 0 & d 
\end{array} \right) , 
\end{array} $$
\item la restriction de $\xi_n$ au groupe de Weil $W_{\mathbb{R}}$ est donn\'ee par :
$$\xi_n (z) = \left( 
\begin{array}{ccc}
\left( \frac{z}{\overline{z}} \right)^{\frac12 +n} & &  \\
& 1 & \\
& & \left( \frac{z}{\overline{z}} \right)^{\frac12 +n} 
\end{array} \right) \times z \mbox{  si } z\in \mathbb{C}^* $$
et 
$$\xi_n (j) = \left( 
\begin{array}{ccc}
1 & & \\
 & 1 & \\
& & -1 
\end{array} \right) \times j .$$
\end{enumerate}
On v\'erifie facilement que ce sont bien des morphismes.

\medskip

Chaque $L$-param\`etre $W_{{\Bbb R}} \rightarrow {}^L U(1,1)$ induit bien \'evidemment un 
morphisme $W_{{\Bbb R}} \rightarrow \widehat{U(1,1)} \rtimes W_{{\Bbb R}}$ qui compos\'e avec l'un quelconque 
des morphismes $\beta_n$ induit une repr\'esentation semi-simple irr\'eductible $W_{{\Bbb R}} \rightarrow GL(2,{\Bbb C})$
et donc un $L$-param\`etre pour le groupe $GL(2,{\Bbb C})$. L\`a encore se pose la question de la fonctorialit\'e~: 
\`a une repr\'esentation automorphe de $U(1,1)$ correspond-il une repr\'esentation automorphe de $GL(2,{\Bbb C})$ via la correspondance
induite par les morphismes $\beta_n$ ? \index{fonctorialit\'e}

On a ainsi deux (si l'on oublie l'entier $n$) mani\`eres de
relever des repr\'esentations de $U(1,1)$ au groupe $GL(2,{\Bbb C})$, les deux seront importantes dans la 
suite.

\medskip
  
De m\^eme chaque $L$-param\`etre $W_{{\Bbb R}} \rightarrow {}^L (U(1,1) \times U(1))$ induit bien \'evidemment 
un morphisme $W_{{\Bbb R}} \rightarrow (GL(2,{\Bbb C}) \times GL(1,{\Bbb C})) \rtimes W_{{\Bbb R}}$ qui compos\'e avec 
l'un quelconque des morphismes $\xi_n$ induit un morphisme $W_{{\Bbb R}} \rightarrow GL (3,{\Bbb C}) \rtimes W_{{\Bbb R}}$ 
qui est un $L$-param\`etre pour le groupe $U(2,1)$. L\`a encore se pose la question de la 
fonctorialit\'e~: \`a une repr\'esentation automorphe de $U(1,1) \times U(1)$ correspond-il une repr\'esentation 
automorphe de $U(2,1)$ via la correspondance induite par les morphismes $\xi_n$ ? \index{fonctorialit\'e}

C'est faux en g\'en\'eral, mais nous verrons au Chapitre 8 que les deux types de fonctorialit\'e 
que nous avons \'evoqu\'es dans ce chapitre suffisent essentiellement \`a d\'ecrire tout le spectre automorphe de $U(2,1)$.

\newpage

\thispagestyle{empty}  

\newpage

\markboth{CHAPITRE 5. REPR\'ESENTATIONS DE $U(a,b)$ ($a,b>1$)}{5.1. GROUPE DUAL}

\chapter{Repr\'esentations de $U(a,b)$ $(a,b>1)$}
Dans ce chapitre, nous d\'ecrivons le groupe dual de $U(a,b)$ pour $a,b>1$, puis ses repr\'esentations cohomologiques. Pour celles-ci, nous d\'ecrivons explicitement les donn\'ees de Langlands, d'abord comme r\'ealisations explicites par des quotients de Langlands de repr\'esentations standard (\S 5.2) puis comme param\`etres dans le groupe dual (\S 5.3). Enfin, nous explicitons pour ce groupe le Th\'eor\`eme de Vogan caract\'erisant les repr\'esentations cohomologiques isol\'ees dans le dual unitaire.

\section{Groupe dual}

La description du groupe dual a d\'ej\`a \'et\'e donn\'ee dans le chapitre~4, puisqu'il ne d\'epend que de la forme int\'erieure quasi-d\'eploy\'ee de $G$, qui est la m\^eme que pour le groupe $U(n-1,1)$ associ\'e (on pose $n=a+b$). Ainsi
$$
\widehat G= GL(n,{\Bbb C})
$$
et $^L \! G=\widehat G\rtimes\mbox{Gal}({\Bbb C}/{\Bbb R})$, l'\'el\'ement non-trivial $\sigma\in\mbox{Gal}({\Bbb C}/{\Bbb R})$ op\'erant par
$$
g\mapsto w_{0} \, ^{t} \! g^{-1} w_{0}^{-1}Ê\quad (g\in\widehat G),
$$
$$w_{0}=\left( 
\begin{array}{ccccc}
&&&&(-1)^{n+1} \\
&&&\adots & \\
&&1 && \\
&-1 &&& \\
1 &&&&
\end{array}\right) .
$$

Nous notons $G=U(a,b)$ le groupe unitaire d\'efini par la forme hermitienne de matrice

$$\left(
\begin{array}{cc}
1_{a} & 0 \\
0 & -1_{b}
\end{array} \right) .$$

et $K=U(a)\times U(b)$ le sous-groupe compact maximal associ\'e; soient $\mathfrak{g}_{0}$ et $\mathfrak{k}_0$ leurs alg\`ebres de Lie r\'eelles respectives. On a
$$
\mathfrak{g}_{0}=\mbox{Lie} (G)=\mathfrak{k}_{0} \oplus \mathfrak{p}_{0},
$$
$$
\mathfrak{p}_{0}=\left\{ \left( 
\begin{array}{cc}
0 & Y \\
- {}^{t}\overline Y & 0
\end{array}
\right) \; : \; Y\in M_{a,b}({\Bbb C} )\right\} .
$$

Enfin $\mathfrak{p}_{0} \otimes {\Bbb C} ={\Bbb C}^{ab} \otimes {\Bbb C} \cong {\Bbb C}^{ab} \oplus {\Bbb C}^{ab}=\mathfrak{p}^{+} \oplus \mathfrak{p}^{-}$ le premier facteur \'etant l'espace tangent holomorphe \`a l'origine \`a $X=G/K$ et le second l'espace tangent antiholomorphe. Noter que $\mathfrak{g}_{0} \otimes {\Bbb C}$ s'identifie naturellement \`a $\mathfrak{g} \mathfrak{l}(n,{\Bbb C})$,
que $\mathfrak{p}_{0}\otimes{\Bbb C}=\left\{ \left( 
\begin{array}{cc}
0 &Y \\
Z  &0
\end{array} \right) \in \mathfrak{g}\mathfrak{l}(n,{\Bbb C}) \right\} $
et qu'alors $\mathfrak{p}^{+}$ peut \^etre d\'efini par $Z=0$, la structure holomorphe sur $\mathfrak{p}_{0}$ \'etant d\'efinie par l'action adjointe de
$$
\iota=\begin{pmatrix}
\sqrt{-1} &\hfill 0\cr
0  & 1
\end{pmatrix} \in  K\ .
$$

\markboth{CHAPITRE 5. REPR\'ESENTATIONS DE $U(a,b)$ ($a,b>1$)}{5.2. REPR\'ESENTATIONS COHOMOLOGIQUES}

\section{Repr\'esentations cohomologiques}

D'apr\`es Vogan et Zuckerman \cite{VoganZuckerman}, celles-ci sont associ\'ees aux alg\`ebres paraboliques $\theta$-stables.

Soit $T\subset G$ le tore diagonal compact,
$$
\mathfrak{t}=\mbox{Lie}(T)=\left\{ \begin{pmatrix}
iH_{1}\cr
&\ddots\cr
&&iH_{n}
\end{pmatrix}\ ,\ H_{i}\in{\Bbb R}\right\}.
$$

Modulo $W_{K}=\mathfrak{S}_{a}\times \mathfrak{S}_{b} $ on peut supposer
$$
H_{1}\geq\ldots \geq H_{a}\ ,\ H_{a+1}\geq\ldots \geq H_{n}.
$$

Il sera commode de noter
$$
X=(H_{1},\ldots H_{a})\in{\Bbb R}^{a} \mbox{  et  }  Y=(H_{a+1} , \ldots , H_{n} )\in {\Bbb R}^{b}.
$$

On associe \`a $H$ l'alg\`ebre de Lie $\mathfrak{q}=\mathfrak{l}+\mathfrak{u} \subset \mathfrak{g}=\mathfrak{g}_{0} \otimes {\Bbb C}$ donn\'ee par
\begin{eqnarray} \label{l}
\mathfrak{l} = \mathfrak{g}^H 
\end{eqnarray}
\begin{eqnarray} \label{u}
\begin{array}{l}
\mathfrak{u} \mbox{ {\it est la somme des espaces radiciels de} } (\mathfrak{g},\mathfrak{t}) \\
\mbox{\it associ\'es aux racines } \alpha \mbox{ {\it telles que }} \langle \alpha,H \rangle >0 .
\end{array}
\end{eqnarray}
\index{sous-alg\`ebre parabolique}

Rappelons quelques propri\'et\'es de la repr\'esentation $A_{\mathfrak{q}}$ (\S 1.). Soit
$$
\begin{array}{l}
\mathfrak{u} \cap \mathfrak{p} = \mathfrak{u} \cap \mathfrak{p}^{+} \oplus \mathfrak{u} \cap \mathfrak{p}^{-} \\
R^{\pm}=\dim(\mathfrak{u} \cap \mathfrak{p}^{\pm}) \\
\mu=2\rho(\mathfrak{u} \cap \mathfrak{p})=\mbox{ somme des racines de } \mathfrak{u}\cap\mathfrak{p} \mbox{ par rapport \`a } \mathfrak{t} . 
\end{array}
$$

\begin{prop} \label{rep cohom}
\begin{enumerate}
\item $A_{\mathfrak{q}}$ contient, avec multiplicit\'e $1$, le $K$-type de plus haut poids $\mu$.
\item On a
$$
H^{p+R^{+},q+R^{-}}(\mathfrak{g},\mathfrak{k}; A_{\mathfrak{q}})\cong\mbox{Hom}_{\mathfrak{l}\cap\mathfrak{k}}(\Lambda^{2p}(\mathfrak{l}
\cap \mathfrak{p}),{\Bbb C}).
$$
En particulier $H^{R^{+},R^{-}}\cong{\Bbb C}$, et la cohomologie n'intervient qu'en degr\'es sup\'erieurs.
\end{enumerate}
\end{prop}
\index{repr\'esentation cohomologique} \index{$A_{\mathfrak{q}}$}

Pour 2. voir \cite[Prop. 6.19]{VoganZuckerman}.

Nous allons maintenant d\'ecrire la param\'etrisation de $A_{\mathfrak{q}}$ \`a l'aide de l'induction parabolique. (La construction de Vogan et Zuckerman est par induction
{\bf cohomologique}, cf. \cite{VoganZuckerman}). Nous aurons besoin des donn\'ees combinatoires suivantes. Soit $r$ le nombre de valeurs distinctes des coordonn\'ees $H_{i}$, et soit
$$
Z_{1}>\ldots>Z_{r} \quad (Z_{j}\in{\Bbb R})
$$
ces valeurs. Alors
\begin{eqnarray*}
X & = & (\underbrace{Z_{1},\ldots Z_{1}}_{a_{1}},\underbrace{Z_{2},\ldots Z_{2}}_{a_{2}},\ldots) \\
Y & = & (\underbrace{Z_{1},\ldots Z_{1}}_{b_{1}},\underbrace{Z_{2},\ldots Z_{2}}_{b_{2}},\ldots)
\end{eqnarray*}
donc $a=\sum a_{j}$\ ,\  $b=\sum b_{j}$.

L'alg\`ebre $\mathfrak{l} \subset \mathfrak{g}$ est en fait r\'eelle, $= \mathfrak{l}_{0} \otimes {\Bbb C}$ avec
\begin{eqnarray} \label{l0}
\mathfrak{l}_0 \cong \prod_{j} \mathfrak{u} (a_{j},b_{j}).
\end{eqnarray}

Soit $d_{j}=\inf(a_{j},b_{j})\geq 0\ ,\ c_{j}=\sup(a_{j} , b_{j})-\inf(a_{j},b_{j})\geq 0\ ,\ n_{j}=a_{j}+b_{j}$.

\medskip

{\bf Remarque.} (cf. \cite[Thm. A8]{Vogan2}) Pour obtenir toutes les repr\'esentations cohomologiques unitaires, il suffit de consid\'erer les sous-alg\`ebres $\mathfrak{q}$ telles que le sous-groupe connexe $L\subset G$ d'alg\`ebre de Lie $\mathfrak{l}_{0}$ soit sans facteur compact non-ab\'elien. On peut donc supposer $a_{j}=1$ (resp. $b_{j}=1$) si $b_{j}=0$ (si $a_{j}=0$).

\medskip

Nous allons associer \`a $\mathfrak{q}$ un sous-groupe parabolique $P=MAN$ de $G$, ainsi qu'une repr\'esentation induisante de $MA$.

La composante d\'eploy\'ee $A$ est la composante neutre d'un tore maximal d\'eploy\'e de $L$. On prend
$$
A=\prod_{d_{j}>0}A_{j}\ ,\ A_{j}\subset L_{j}\cong U(a_{j},b_{j}).
$$ 
Si $0< a_{j}\leq b_{j}$, $\mathfrak{a}_{j}=\mbox{Lie}(A_{j})$ est donn\'ee par
\begin{eqnarray} \label{aj1}
\mathfrak{a}_{j} = \left(
\begin{array}{c|c}
0 & 
\begin{array}{cc}
\begin{array}{ccc}
t_1 & & \\
& \ddots & \\
& & t_a
\end{array} & 0 
\end{array} \\ \hline
\begin{array}{c}
\begin{array}{ccc}
t_1 & & \\
& \ddots & \\
& & t_a 
\end{array} \\
0 
\end{array} & 0 
\end{array} \right) , \ a=a_{j} .
\end{eqnarray}
Si $0< b_{j}\leq a_{j}$,
\begin{eqnarray} \label{aj2}
\mathfrak{a}_{j} = \left(
\begin{array}{c|c}
0 & 
\begin{array}{c}
\begin{array}{ccc}
t_1 & & \\
& \ddots & \\
& & t_b
\end{array} \\
 0 
\end{array} \\ \hline
\begin{array}{cc}
0 &
\begin{array}{ccc}
t_1 & & \\
& \ddots & \\
& & t_b 
\end{array} 
\end{array} & 0 
\end{array} \right) , \ b=b_{j} .
\end{eqnarray}

Soit ${\Bbb M}$ le centralisateur de $A$ dans $G$. On a ${\Bbb C}^{n}={\Bbb C}^{a}\oplus{\Bbb C}^{b}$ ; on v\'erifie que si $a_{j}\leq b_{j}$ (\ref{aj1}) ${\Bbb M}$ doit pr\'eserver le sous-espace ${\Bbb C}^{2a_{j}}$ de ${\Bbb C}^{a}$ d\'etermin\'e par (\ref{aj1}) ; de m\^eme si $b_{j}\leq a_{j}$. Sur ce sous-espace, ${\Bbb M}$ op\`ere pour le centralisateur de $A_{j}\cong({\Bbb R}_{+}^* )^{d_{j}}$, d'o\`u un facteur $({\Bbb C}^{*})^{d_{j}}$ ; sur l'orthogonal de $\bigoplus_{j}{\Bbb C}^{2d_{j}}$, ${\Bbb C}^{2d_j }$ op\`ere comme 
$U(a-\sum d_{j}\ ,\ b-\sum d_{j})$.

Donc
$$
{\Bbb M}=MA ,\ A\cong({\Bbb R}^{*})^{\sum d_{j}}  ,
$$
\begin{eqnarray} \label{M}
M\cong\prod_{j}U(1)^{d_{j}} \times U(a-\sum d_{j},b-\sum d_{j}).
\end{eqnarray}
On v\'erifie que ${\Bbb M}$ est bien le sous-groupe de Levi d'un parabolique (r\'eel) de $G$.

Nous devons d\'eterminer une repr\'esentation $\sigma\otimes\nu $ de ${\Bbb M}=MA$. On a $\nu\in\mathfrak{a}^{*}$, $\nu=\bigoplus_{j}\nu_{j}$. Dans le cas (\ref{aj1}) on pose pour $T=(t_{1},\ldots t_{{a}_{j}})\in \mathfrak{a}_{j}$~:
\begin{eqnarray} \label{nuj1}
\nu_{j}(T)=(n_{j}-1)t_{1}+\ldots+(n_{j}+1-2a_{j})t_{{a}_{j}}\quad (n_{j}=a_{j}+b_{j}). 
\end{eqnarray}

Dans le cas (\ref{aj2}) on posera de m\^eme~:
\begin{eqnarray} \label{nuj2}
\nu_{j}(T)=(n_{j}-1)t_{1}+\ldots+(n_{j}+1-2b_{j})t_{{b}_{j}}  .
\end{eqnarray}

On renvoie le lecteur \`a Vogan-Zuckerman \cite[p.~82]{VoganZuckerman} pour la justification de ce choix. 
($\nu$ est la demi-somme des racines de $A$ dans une alg\`ebre d'Iwasawa de $\mathfrak{l}$.)

La repr\'esentation $\sigma$ de $M$ appartient \`a la s\'erie discr\`ete. On d\'ecrit donc son param\`etre d'Harish-Chandra (cf. \S 4.1.2). Soit $T^{+}$ la composante compacte d'un tore maximal $\theta$-stable, contenant $A$, de $L$. On prendra de fa\c con naturelle 
\begin{eqnarray} \label{T+}
T^{+}=\prod_{j} U(1)^{d_{j}}\times\prod_{j} U(1)^{c_{j}} , 
\end{eqnarray} 
les facteurs $U(1)^{d_{j}}$ sont plong\'es dans $U(a_{j},b_{j})$ de fa\c con \`a commuter \`a $\mathfrak{a}_{j}$ et les facteurs $U(1)^{c_{j}}$ sont contenus dans $U(a_{j})$ si $a_{j}>b_{j}$ et dans  $U(b_{j})$ si $a_{j}<b_{j}$. Alors $T^{+}$ est un tore maximal de $M$.

En fait $T^{+}$ est un tore maximal de $C=M\cap L=\displaystyle\prod_{j}U(1)^{d_{j}}\times\displaystyle\prod_{j}U(c_{j}) $.

Sur chaque facteur $U(c_{j})$ (contenu dans $U(a_{j})$ si $a_{j}>b_{j}$ ou dans $U(b_{i})$ si $a_{j}<b_{i})$ choisissons le syst\`eme de racines positives naturel par rapport au tore diagonal~:
$$
\Delta=\left\{ u_{k}u_{l}^{-1}\; : \;  k>l ,\; u=(u_{k})\in U(1)^{c_{j}}\right\}.
$$ 
Soit $\rho_{c}$ la demi-somme des racines positives dans $\frac12 (\Lambda_{T^{+}})$, $\Lambda_{T^{+}} \subset it^{*}$ \'etant le r\'eseau des caract\`eres.

Alors~:
\begin{eqnarray}
\mbox{{\it Le param\`etre d'Harish-Chandra de }} \index{Harish-Chandra, param\`etre d'} \sigma \mbox{ {\it est \'egal \`a }} \lambda_{\sigma}=\rho_{c}+\rho(\mathfrak{u}).
\end{eqnarray}

Il nous sera n\'ecessaire de conna\^{\i}tre explicitement $\lambda_{\sigma}$ selon la param\'etrisation (\ref{T+}) de $T^{+}$. On note un \'el\'ement $t\in T^{+}$ selon (\ref{T+})~:
$$
t=(v_{ik},u_{jk})\ (k\leq d_{j}, l \leq c_{j}).
$$

De m\^eme si $X\in\mbox{Lie}(T^{+})$,
$$
X=(V_{jk},U_{jl}).
$$

\begin{lem} \label{rho}
\begin{enumerate}
\item $\rho_{c}(X) = \displaystyle\sum_{j} \left( \frac{c_{j}-1}{2} \right) U_{j,1}+\ldots+\left( \frac{1-c_{j}}{2}\right) U_{j,c_{j}}$,
\item $\rho_{\mathfrak{u}}(X)=\displaystyle\sum_{j} \frac{m_{j}}{2}\left( \displaystyle\sum_{k}2V_{jk}+\displaystyle\sum_{l} U_{jl}\right)$ 
o\`u $m_{j}=-n_{1} - \ldots - n_{j-1}+n_{j+1}+\ldots+n_{r}  $.
\end{enumerate}
\end{lem}

L'assertion 1. r\'esulte de la description pr\'ec\'edente de $\rho_{c}$ ; quant \`a 2. noter que
$L_{{\Bbb C}}=\prod_{j} GL(n_{j},{\Bbb C})$ est le sous-groupe de Levi d'un sous-groupe parabolique de $G_{{\Bbb C}}=GL(n,{{\Bbb C}})$; $L_{{\Bbb C}}$ op\`ere sur 
$\mathfrak{u}$ avec pour d\'eterminant
$$
2\rho_{\mathfrak{u}}=(\det g_{1})^{m_{1}}(\det g_{2})^{m_{2}}\cdots(\det g_{r})^{m_{r}}\ ; 
$$
la formule 2. s'en d\'eduit par restriction.

V\'erifions que $\lambda_{\sigma}$ est bien le param\`etre d'Harish-Chandra pour une s\'erie discr\`ete de $M$ (d\'efini par (\ref{M})). Tout d'abord, $\lambda_{\sigma}$ est int\'egral sur $\displaystyle\prod_{j}U(1)^{d_{j}}\subset T^{+}$. Sur $T^{+}\cap U(a-\sum d_{j}, b- \sum d_j )$,  variables $u_{j,l}$, on doit v\'erifier d'apr\`es le \S 4.1.2 que les coordonn\'ees de $\lambda_{\sigma}$ sont $\equiv \frac{a+b-2\sum d_{j}-1}{2} [1]$. Ceci r\'esulte du Lemme (noter que $c_{j}\equiv n_{j}$ [2]). Enfin, on v\'erifie que $\lambda_{\sigma}$ est r\'egulier dans $X^{*}(T^{+})$ (par rapport aux racines de $M$).

Soit $P$ un sous-groupe parabolique de $G$ de radical de Levi $MA$ : $P=MAN$.

\begin{prop}[Vogan-Zuckerman \cite{VoganZuckerman}]
\begin{enumerate} 
\item La repr\'esentation unitairement induite 
$$
I(\sigma,\nu)={\rm ind}_{P}^{G} (\sigma\otimes e^{\nu})
$$
admet un unique quotient irr\'eductible $J(\sigma,\nu)$.
\item $J(\sigma,\nu)\cong A_{\mathfrak{q}}$.
\end{enumerate}
\end{prop}
\index{repr\'esentation cohomologique}

\markboth{CHAPITRE 5. REPR\'ESENTATIONS DE $U(a,b)$ ($a,b>1$)}{5.3. PARAM\`ETRES DES $A_{\mathfrak{q}}$ EN DUALIT\'E DE LANGLANDS}

\section{Param\`etres des $A_{\mathfrak{q}}$ en dualit\'e de Langlands}

Rappelons que ${}^{L}G=GL(n,{\Bbb C})\rtimes {\rm Gal}({\Bbb C}/{\Bbb R}) $, $\sigma$ op\'erant par
$$
g\mapsto w_{0}{}^{t}g^{-1}w_{0}^{-1}
$$ 
$$
w_{0}= \left(
\begin{array}{ccccc}
&&&&(-1)^{n+1} \\
&&&\adots & \\
&&1 && \\
&-1 &&& \\
1 & & & & 
\end{array} \right).
$$

Soit ${\Bbb M}$ le groupe de Levi associ\'e \`a $\mathfrak{q}$ : on v\'erifie facilement que
$$
{\Bbb M}\cong\prod_{j}({\Bbb C}^{\times })^{d_{j}} \times U(A,B)
$$
o\`u $A=a-\sum d_{j}$, $B=b-\sum d_{j}$ ; soit $N=A+B$, $D=\sum d_{i}$.

Soit $\widehat  {\Bbb M}=({\Bbb C}^{* })^{D}\times GL(N,{\Bbb C})\times ({\Bbb C}^{* })^{D}\subset\widehat G$ (plongement par blocs diagonaux). Alors $\widehat  {\Bbb M}$ est stable par l'action de $w_{0}$, et ${}^{L}\widehat  {\Bbb M}=\widehat  {\Bbb M}\times {\rm Gal}({\Bbb C}/{\Bbb R})$ s'identifie naturellement au groupe dual de ${\Bbb M}$.

D'apr\`es Langlands \cite{Langlands}, le param\`etre $\varphi:W_{\Bbb R} \rightarrow {}^{L}G$ associ\'e \`a $A_{\mathfrak{q}}$ est obtenu par composition \`a partir de celui de $\sigma\otimes e^{\nu}$. Pour simplifier, nous ne d\'ecrivons que $\varphi_{|_{W_{\Bbb C}}}$ o\`u $W_{\Bbb C}={\Bbb C}^{*}$.

La repr\'esentation $\sigma\otimes e^{\nu}$ de ${\Bbb M}$ d\'etermine un caract\`ere de $({\Bbb C}^{*})^{D} $, d\'eduit de (\ref{nuj1}) (ou (\ref{nuj2})) et du Lemme \ref{rho}. Rappelons que $D=\sum d_{j}$ ; les variables de $({\Bbb C}^{*})^{D}$ peuvent \^etre index\'ees par $(j,k)$, $k\leq d_{j}$ (\S 5.3). Le caract\`ere associ\'e est donn\'e par
\begin{eqnarray} \label{z}
z=z_{jk}\mapsto(z/\overline z)^{\frac{m_{j}}{2}}(z\overline z)^{\frac{n_{j}+1-2k}{2}}
\end{eqnarray}
(Lemme \ref{rho}). Par dualit\'e de Langlands, la partie relative au facteur $({\Bbb C}^{*})^{D}\times({\Bbb C}^{*})^{D} $ de $\widehat {\Bbb M}$ est alors~:
$$
\varphi_{1}:z\mapsto\left(
\begin{array}{cc}
(z/ \overline z)^{\frac{m_{j}}{2}}(z\overline z)^{\frac{n_{j}+1-2k}{2}}  & \\
&(z/ \overline z)^{\frac{m_{j}}{ 2}}(z\overline z)^{-\frac{n_{j}+1-2k}{2}} 
\end{array} \right)
$$
(matrices diagonales ; l'ordre des entr\'ees $j,k$ dans le second facteur doit \^etre inverse de l'ordre du premier facteur, pour que $\varphi_{1}$ soit compatible \`a l'action de $w_{0}$).

Par ailleurs $\sigma$ donne par restriction une repr\'esentation (discr\`ete) de $U(A,B)$, dont le param\`etre de Langlands se d\'eduit du param\`etre d'Harish-Chandra comme on l'a expliqu\'e dans le Chapitre~4 pour $U(n,1)$ ; toujours d'apr\`es le Lemme \ref{rho}, on en d\'eduit
$$
\begin{array}{cccl}
\varphi_{2} : &  {\Bbb C}^{*} & \rightarrow & GL(N,{\Bbb C}) \\
                       &  z                    & \mapsto       & \left( (z/\overline z)^{\frac{m_{j}+c_{j}+1-2l}{2}} \right) \ (j=1,\ldots r, \  l\leq c_{j}) .
\end{array}
$$   
On remarquera que $N\equiv n \ [2]$, que $c_{j}\equiv n_j \ [2]$ et qu'il r\'esulte alors de l'expression de $m_{j}$ que 
$m_{j}+c_{j}+1-2l \equiv N-1 \ [2]$, ce qui est n\'ecessaire pour que $\varphi_{2}$ d\'efinisse une s\'erie discr\`ete pour $U(A,B)$ (cf. la Prop. 4.3.2 en rang~1).

R\'ecapitulons le r\'esultat obtenu~:

\begin{prop} \label{phi}
Le param\`etre de Langlands de $A_{\mathfrak{q}}$~:
$$
\varphi_{|W_{\Bbb C}} : {\Bbb C}^{*} \rightarrow \widehat {\Bbb M} \subset\widehat G 
$$
est donn\'e par
$$
z\mapsto \left(
\begin{array}{ccc}
(z/\overline z)^{\frac{m_{j}}{2} }(z \overline z)^{\frac{n_{j}+1-2k}{2}} \\
& (z/\overline z)^{\frac{m_{j}+c_{j}+1-2l}{2}} \\
&& (z/\overline z)^{\frac{m_{j}}{2} }(z\overline z)^{-\frac{n_{j}+1-2k}{2}}
\end{array} \right)
$$
$(k\leq d_{j},\ l\leq c_{j})$.
\end{prop}
\index{param\`etres de Langlands des repr\'esentations cohomologiques}

\markboth{CHAPITRE 5. REPR\'ESENTATIONS DE $U(a,b)$ ($a,b>1$)}{5.4. REPR\'ESENTATIONS COHOMOLOGIQUES ISOL\'EES}

\section{Repr\'esentations cohomologiques isol\'ees}

Dans ce paragraphe nous explicitons pour $U(a,b)$ un th\'eor\`eme de Vogan \cite[Th\'eor\`eme A10]{Vogan2} caract\'erisant les repr\'esentations cohomologiques {\bf isol\'ees}. 

Nous conservons les notations du \S 5.2, donc~:
$$
Z_{1}>\ldots>Z_{r}  . 
$$ 

Nous ferons sur les donn\'ees l'hypoth\`ese de non-d\'eg\'en\'erescence expliqu\'ee dans la Remarque suivant la Proposition \ref{rep cohom}, soit
$$
a_{j} =0\Rightarrow b_{j}=1,\ \;  b_{j}=0\Rightarrow a_{j}=1 .
\leqno{\rm(H0)}$$

Le facteur correspondant de $L=\prod U(a_{j},b_{j})$ est $U(1)$. Notons
$$\begin{array}{ccc}
j_{1}<\ldots <j_{t} & \quad & (0\leq t\leq r) \\
Z_{j1}>\cdots>Z_{jt} & &  
\end{array}$$
les autres valeurs de $j$. Pour de tels $j$ le facteur $L_{j}\cong U(a_{j},b_{j})$ est non compact. Consid\'erons l'hypoth\`ese suivante~:
$$
\mbox{{\it Si} } j=j_{i} \; (i \leq t ) \mbox{ {\it le rang semi-simple de} } L_{j} \mbox{ {\it est} } >1  \Longleftrightarrow a_{j}, b_{j}\geq 2 . 
\leqno{\rm(H1)}$$

Supposons enfin que $t<r$. Notons $I=\{ j,j+1,\ldots , j+k \}$ un intervalle non vide d'entiers contenu dans $\{1,\ldots ,r \}-\{j_{1},\ldots , j_{t}\}$. Pour $j\in I$ on a donc $\{a_{j},b_{j}\}=
\{0,1\}$.

$$
\mbox{{\it Si} } t<r, \; a_{j} \mbox{ {\it (et donc} } b_{j}\mbox{{\it ) est constant sur tout intervalle} } I \mbox{ {\it du compl\'ementaire de} } \{j_{t}\} .
\leqno{\rm(H2)}$$

\begin{prop}[Vogan \cite{Vogan2}] \label{Vogan}
Supposons $\mathfrak{q}$ associ\'ee aux donn\'ees $(a_{j},b_{j})$ v\'erifiant {\rm (H0)}. Alors la repr\'esentation $A_{\mathfrak{q}}$ est isol\'ee dans le dual de $SU(a,b)$ si et seulement si {\rm(H1)} et {\rm(H2)} sont v\'erifi\'ees.
\end{prop}\index{Th\'eor\`eme de Vogan}

Rappelons (\S 5.2) que la cohomologie de $A_{\mathfrak{q}}$ n'appara\^{\i}t qu'en degr\'es $\geq R=R({\mathfrak{q}})$.

\begin{cor} \label{cor de vogan}
$({\rm rg}\ G\geq 2)$ Si  $A_{\mathfrak{q}}$ n'est pas isol\'ee, la cohomologie de  $A_{\mathfrak{q}}$ n'appara\^{\i}t qu'en degr\'es $i\geq a+b-2$. 
\end{cor} \index{repr\'esentations cohomologiques isol\'ees}

Noter que si $t=0$, (H2) veut dire que $a_{j} $ (et donc $b_{j}$) est constant. Le cas o\`u $t=0$ correspond \`a $L=U(1)^{n}$~; la repr\'esentation $A_{\mathfrak{q}}$  est alors une s\'erie discr\`ete. Sous (H0) la condition (H2) impliquerait alors que $a_{j}\equiv 1$ et $b_{j}\equiv 0$ (ou l'inverse), ce qui est impossible puisque $a,b>0$. Donc (H2) est viol\'ee si $A_{\mathfrak{q}}$  est une s\'erie discr\`ete~: celles-ci ne sont pas isol\'ees.

\medskip
{\bf Remarque.} Du point de vue de la th\'eorie des repr\'esentations de $G$, ce r\'esultat ne peut \^etre am\'elior\'e. Consid\'erons par exemple $U(2,2)$. On peut choisir $\mathfrak{q}$ de sorte que $\mathfrak{l}$ est l'alg\`ebre diagonale par blocs $\mathfrak{u}(2,1)\times \mathfrak{u}(1)$. Elle viole (H2) donc $A_{\mathfrak{q}}$ n'est pas isol\'ee~; sa cohomologie appara\^{\i}t en degr\'es $\geq R=2$. On trouve ais\'ement de tels exemples en dimensions sup\'erieures.

D\'emontrons le Corollaire. Tout d'abord, un calcul simple donne
$$
R=ab-\sum_{j=1}^{r}a_{j}b_{j}=ab-\sum_{i=1}^{t}a^{i}b^{i}
$$
o\`u on a \'ecrit pour simplifier $a^{i},b^{i}=a_{j_i},b_{j_i}$.

Supposons que $\mathfrak{q}$ viole (H1). Alors (\`a l'ordre pr\`es) on peut supposer $a^{1}=1$, $b^{1}\geq 1$. Alors 
$$
\sum_{i=2}^{t}a^{i}b^{i}\leq(\sum a^{i})(\sum b^{i})\leq(a-1)(b-b^{1})  
$$ 
donc
\begin{eqnarray*}
R & \geq & ab-b^{1}-(a-1)(b-b^{1}) \\
    & = & b+(a-2)b^{1} \geq a+b-2 
\end{eqnarray*}
($a\geq 2$ puisque ${\rm rg} G \geq 2$).

Supposons que $\mathfrak{q}$ viole (H2). Donc $t<r$ et $j\mapsto a_{j}$ prend les valeurs $1$ et $0$ (donc $b_{j}=1$) sur $\{1,\ldots , n\}-\{j_t \}$. Alors $\sum a^{i}\leq a-1$, $\sum b^{i}\leq b-1$ et dans ce cas
$$
R\geq ab-(a-1)(b-1)=a+b-1>a+b-2 .
$$

Il nous reste \`a d\'eduire la Proposition \ref{Vogan} des r\'esultats de Vogan. D'apr\`es le Th\'eor\`eme A.10 de \cite{Vogan2}, $A_{\mathfrak{q}}$ est isol\'ee si et seulement si ${\mathfrak{q}}$ v\'erifie certaines conditions (0)-(3). Ici (0) est notre hypoth\`ese (H0)~; (1) est automatiquement v\'erifi\'ee si $G=U(a,b)$~; (2) est (H1). Il nous reste \`a exprimer (3).

\medskip
La construction-classification de Vogan-Zuckerman \cite{VoganZuckerman} pour les $A_{\mathfrak{q}}$ est la suivante. Tout d'abord on a $A_{\mathfrak{q}}=R_{\mathfrak{q}}({\Bbb C})$ o\`u ${\Bbb C}$ est la repr\'esentation triviale de $\mathfrak{l}$ et le foncteur $R_{\mathfrak{q}}=R_{\mathfrak{q}}^{\dim({\mathfrak{u}}\cap \mathfrak{k})}$ est d\'efini dans les r\'ef\'erences cit\'ees par \cite{VoganZuckerman}. Supposons donn\'e un syst\`eme de racines positives $\Delta^{+}({\mathfrak{g}},{\mathfrak{t}})$ pour $({\mathfrak{g}},{\mathfrak{t}})$ tel que les racines de $\mathfrak{u}$ soient positives. Alors, avec les notations usuelles~:
$$
\rho_{\mathfrak{g}}=\rho_{\mathfrak{u}}+\rho_{\mathfrak{l}} .
$$ 
Le caract\`ere infinit\'esimal de ${\Bbb C}_{\mathfrak{l}}$ est $\rho_{\mathfrak{l}}$~; $\lambda=\rho_{\mathfrak{g}}$ v\'erifie les hypoth\`eses du Th\'eor\`eme A.10 de \cite{Vogan2} pour $\Delta^{+}({\mathfrak{g}},{\mathfrak{t}})$.

Soit $\prod\subset\Delta^{+}({\mathfrak{g}},{\mathfrak{t}})$ l'ensemble des racines simples et $\prod({\mathfrak{l}})$ le sous-ensemble form\'e des racines simples de $\mathfrak{l}$. Alors la condition (3) de Vogan s'\'ecrit~:
$$
\langle \beta^{\vee},\lambda \rangle = \langle \beta^{\vee},\rho_{\mathfrak{g}} \rangle \neq 1
\leqno(H2')$$
{\it pour toute racine \(imaginaire\) non compacte $\beta\in\prod$ orthogonale \`a $\prod({\mathfrak{l}})$}.

\medskip
Il s'agit d'expliciter (H2'). Rappelons (\S 5.2) que $H=(X,Y)$. Les valeurs des coordonn\'ees sont
$$\begin{array}{ccl}
& & Z_{1}>Z_{2}>\ldots>Z_{r} \\ 
\mbox{et} & \quad & Z_{j_1}>Z_{j_2}>\ldots>Z_{j_{t}}  
\end{array}
$$
sont les valeurs de multiplicit\'es $>1$, donc $\geq 4$ d'apr\`es (H1). \'Ecrivons alors
\begin{eqnarray} \label{X}
X=(X_{1}>X_{2} > \ldots >X_{\alpha_{1}}>\underbrace{Z_{j_{1}}}_{a_{1}}>\ldots>X_{\alpha_{2}}> 
\underbrace{Z_{j_{2}}}_{a_{2}}>\ldots>X_{\alpha}); 
\end{eqnarray}
les $X_{i}$ apparaissant avec multiplicit\'e $1$ d'apr\`es (H0)~; de m\^eme
\begin{eqnarray} \label{Y}
Y=(Y_{1}>\ldots>Y_{\beta_{1}}>\underbrace{Z_{j_{1}}}_{b_{1}}>\ldots>Y_{\beta_{2}}> 
\underbrace{Z_{j_{2}}}_{b_{2}}>\ldots>Y_{\beta}). 
\end{eqnarray}

Enfin soit $H'$ \'egal \`a $H=(X,Y)$, r\'eordonn\'e de fa\c con positive~:
\begin{eqnarray} \label{H'}
H'=\sigma H=(Z_{1}>Z_{2}>\ldots>\underbrace{Z_{j_{1}}}_{n_{1}}>\ldots> 
\underbrace{Z_{j_{2}}}_{n_{2}}>\ldots>Z_{r}),
\end{eqnarray}
$\sigma\in{\mathfrak{S}}_{n}$. 

La base $\Delta^{+}({\mathfrak{g}},{\mathfrak{t}})$ est un ensemble de racines positives $\varepsilon$ telles que
$\langle \varepsilon ,H \rangle \geq 0$~; $\sigma$ l'envoie sur un ensemble de racines telles que $\langle \varepsilon ,H' \rangle \geq 0$, que l'on prendra \'egal \`a l'ensemble usuel puisque $H'$ est dominant. Alors
\begin{eqnarray*}
\sigma\prod & = & \{\varepsilon_{i}=(0,\ldots ,0,1,-1,0, \ldots), \; i=1,\ldots n-1 \} \\
\sigma\prod & = & \{\varepsilon_{i} \; : \; H_{i}'=H_{i+1}' \},
\end{eqnarray*}
et $\sigma$ envoie l'orthogonal de $\prod({\mathfrak{l}})$ sur l'ensemble des racines $\varepsilon_{i}$ telles que
\begin{eqnarray} \label{formule}
\{i,i+1\}\subset\{1,\ldots , n\}-\displaystyle\cup_{i=1}^{t} I_{i},
\end{eqnarray}
$I_{i}$ \'etant le support de $Z_{j_i}$ dans l'expression de $H'$ (\ref{H'}).

Par ailleurs, le syst\`eme de racines \'etant de type $A_{n-1}$, $\beta^{\vee}=\beta$ si
$\beta\in\prod$ et donc $\langle \beta^{\vee},\rho_{\mathfrak{g}} \rangle =1$. La condition (H2') est donc \'equivalente \`a

$$
\mbox{{\it Il n'y a pas de racine imaginaire non compacte dans} } \prod \mbox{ {\it orthogonale \`a} } \prod({\mathfrak{l}}).
\leqno(H2'')$$

\medskip
Soit donc $\varepsilon_{i}$ v\'erifiant (\ref{formule}). Ceci implique donc que $i$ et $i+1$ appartiennent \`a une composante connexe de cardinal $\geq 2$ de $\{1,\ldots ,n\}-\cup I_{i}$. Supposons par exemple que $\alpha_{1}$ ou $\beta_{1}>1$ et que $\alpha_{1}+\beta_{1}\geq 2$. Alors dans l'expression (\ref{H'}) $i,i+1$ sont deux indices associ\'es \`a $Z_{i}>Z_{i+1}>Z_{j_1}$. Alors $\sigma^{-1}\varepsilon_{i}$ est associ\'ee dans l'expression (\ref{X}), (\ref{Y}) \`a deux couples de la forme $(X_{i},X_{i+1})$, $(X_{j},Y_{j'})$ ou $(Y_{j},Y_{j+1})$  \`a gauche de $Z_{j_1}$. Dans le second cas, la racine de valeur $X_{j}-Y_{j'}$ est non compacte.

Donc $Z_{1}$ et $Z_{2}$ sont tous deux (par exemple) de la forme $(X_{1},X_{2})$~; il en est de m\^eme pour $X_{2}$ et $X_{3}$, etc. Ceci veut dire que pour tout $j<j_{1}$ on a $a_j =1$ (ou $b_j=1$), conform\'ement \`a (H2). Le m\^eme argument s'applique \`a toute composante connexe, et il est clair que (H2) est en fait \'equivalente \`a (H2'').

\newpage

\thispagestyle{empty}

\empty

\markboth{CHAPITRE 6. CONJECTURES D'ARTHUR}{6.1. PARAM\`ETRES D'ARTHUR}

\chapter{Cons\'equences des Conjectures d'Arthur}

Dans deux articles fondamentaux, Arthur a donn\'e une description conjecturale des repr\'esentations des groupes r\'eductifs qui peuvent appara\^{\i}tre dans $L^{2}(\Gamma \backslash G)$ pour un sous-groupe de congruence. Le but de ce chapitre est d'expliquer les cons\'equences de ces conjectures pour la th\'eorie spectrale des formes diff\'erentielles.

Les Conjectures d'Arthur reposent elles-m\^emes sur la construction hypoth\'etique de ``paquets locaux'' (cf. \cite[\S 4]{Arthur}), familles finies de repr\'esentations de $G$ jouissant de propri\'et\'es de stabilit\'e. Arthur lui-m\^eme ne dit rien de la construction de ces ``paquets'', sauf dans le cas des repr\'esentations cohomologiques \cite[\S 5]{Arthur}. En particulier, sa formulation ne pr\'ecise pas quels ``param\`etres'' (\pa 6.1) devraient appara\^{\i}tre quand $G$ n'est pas quasi-d\'eploy\'e. Pour les groupes r\'eels, ces paquets locaux sont a priori construits par Adams, Barbasch et Vogan \cite{AdamsBarbaschVogan}. Mais la construction, de nature g\'eom\'etrique, est tr\`es difficile. 

Pour des groupes de rang $1$ assez simples tels que ceux qui nous int\'eressent, nous avons pr\'ef\'er\'e utiliser la th\'eorie d'Arthur  \textbf{a minima}. Une cons\'equence non ambigu\"e de celle-ci est la description d'une famille de \textbf{caract\`eres infinit\'esimaux}, \index{caract\`ere infinit\'esimal}
les seuls possibles pour les repr\'esentations apparaissant dans $L^{2}(\Gamma \backslash G)$. Pour les groupes de rang $1$, ces restrictions sur le caract\`ere infinit\'esimal imposent des limitations s\'ev\`eres aux repr\'esentations.
 
\section{Param\`etres d'Arthur}

Soit $G$ un groupe r\'eductif r\'eel, nous noterons $G= G(\mathbb{R})$, $\widehat G$ le groupe dual (groupe r\'eductif complexe) connexe, ${}^{L}G=\widehat G\times \mathrm{Gal}(\mathbb{C}/\mathbb{R}) $ le groupe dual. (Pour les groupes unitaires v. ch.~4 ; pour les groupes orthogonaux v. \S 3-4).

\begin{defi} \index{param\`etre d'Arthur}
Un param\`etre d'Arthur pour $G$ est un homomorphisme
\begin{eqnarray}
\psi:W_{\mathbb{R}}\times\mathrm{SL}(2,\mathbb{C})\rightarrow{}^{L}G
\end{eqnarray}
tel que
\begin{itemize}
\item[$(i)$] Le diagramme
$$
\begin{array}{rcccl}
W_{\mathbb{R}} & &\rightarrow & & {}^{L}G\\
&\searrow & \quad & \swarrow &\\
& & \mathrm{Gal}(\mathbb{C}/\mathbb{R}) & &
\end{array}
$$
est commutatif.

\item[$(ii)$] La restriction $\psi\vert_{ \mathrm{SL}(2,\mathbb{C})}$ est holomorphe $(\equiv$ alg\'ebrique$)$.

\item[$(iii)$] L'image de $\psi\vert_{W_{\mathbb{R}}}$ est d'adh\'erence compacte.
\end{itemize}
\end{defi}

On d\'eduit de $\psi$ un ``param\`etre de Langlands''
\begin{eqnarray}
\begin{array}{cl}
\varphi_{\psi}: & W_{\mathbb{R}}\rightarrow{}^{L}G\\
& w\mapsto\psi\left(w,\left(
\begin{array}{cc}
\vert w\vert^{1/2} & 0\\
0 & \vert w\vert^{-1/2} 
\end{array}\right) \right)
\end{array}
\end{eqnarray}
o\`u, pour $w\in W_{\mathbb{R}}$, $\vert w\vert$ est la valeur absolue de l'image de $w$ dans $\mathbb{R}^{\times }$ (Ch.~3).

\vskip0,2cm
Noter que $\varphi_{\psi}$ ne d\'efinit pas toujours une repr\'esentation de $G$ car il ne v\'erifie pas n\'ecessairement la condition (4.3.3) de la D\'efinition du \S 4.3.

Soit $\mathfrak{g}_{0}=\mathrm{Lie}(G)$, $\mathfrak{g}=\mathfrak{g}_{0}\otimes \mathbb{C}$, $\mathfrak{h}$ une sous-alg\`ebre de Cartan de $\mathfrak{g}$. Rappelons que le centre de
$\mathfrak{Z}$ de $U(\mathfrak{g})$ s'identifie \`a $S(\mathfrak{h})^{W}$, $W$ \'etant le groupe de Weyl\ $W(\mathfrak{g},\mathfrak{h})$. On peut v\'erifier alors que $\varphi_{\psi}$ d\'efinit un caract\`ere infinit\'esimal, {\it i.e.} un \'el\'ement de
$\mathfrak{h}^{\ast}/W$ (pour tout choix de $\mathfrak{h}$). Nous le ferons explicitement pour les groupes qui nous int\'eressent. Nous utiliserons alors la formulation tr\`es faible suivante des Conjectures d'Arthur~:

\begin{conj} \index{Conjecture d'Arthur}
Si une repr\'esentation irr\'eductible $\pi$ de $G$ appara\^{\i}t $($faiblement$)$ dans $L^{2}(\Gamma\backslash G)$ pour un sous-groupe de congruence, son caract\`ere infinit\'esimal $\lambda_{\pi}$ est associ\'e \`a un param\`etre d'Arthur $\varphi_{\psi}$.
\end{conj}

\markboth{CHAPITRE 6. CONJECTURES D'ARTHUR}{6.2. $G=U(n,1)$}

\section{$G=U(n,1)$}

Si $G=U(n,1)$ toute alg\`ebre de Cartan de $\mathfrak{g}_{\mathbb{C}}$ s'identifie \`a
$\mathbb{C}^{n+1}$, de fa\c con unique modulo l'action de $\mathfrak{S}_{n+1}=W$. Les racines sont donn\'ees par
\begin{eqnarray*}
(x_{1},\ldots x_{n+1})\mapsto x_{i}-x_{j}\quad(i\neq j). 
\end{eqnarray*}

Soit ${}^{L}G=\mathrm{GL}(n+1,\mathbb{C})\rtimes \mathrm{Gal}(\mathbb{C}/\mathbb{R})$. 

Rappelons que ${}^{L}G$ a \'et\'e construit dans ce cas au \S 4.2 ; l'\'el\'ement non trivial $\sigma\in\mathrm{Gal}(\mathbb{C}/\mathbb{R})$ op\`ere par
\begin{eqnarray*}
g\mapsto w_{0}{}^{t}g^{-1}w_{0}^{-1}\ ,
\end{eqnarray*}

\begin{eqnarray*}
w_{0} = \left(
\begin{array}{cccc}
&&&(-1)^{n}\\
&&\adots &\\
&-1& & \\
1 & & & 
\end{array}\right) .
\end{eqnarray*}

Soit
\begin{eqnarray*}
\psi:W_{\mathbb{R}}\times\mathrm{SL}(2,\mathbb{C})\rightarrow{}^{L}G 
\end{eqnarray*}
un param\`etre d'Arthur, consid\'erons sa restriction

\begin{eqnarray*}
\psi_{0}:(W_{\mathbb{C}}=\mathbb{C}^{\times })\times\mathrm{SL}(2,\mathbb{C})\rightarrow\widehat G=\mathrm{GL}(n+1,\mathbb{C}). 
\end{eqnarray*}
D'apr\`es la d\'efinition 6.1.1 (iii) $\psi_{0}$ est semi-simple et s'\'ecrit donc :

\begin{eqnarray}
\psi_{0}=\bigoplus_{j=1}^{r}r_{j}\otimes \chi_{j} 
\end{eqnarray}
$r_{j}$ \'etant une repr\'esentation irr\'eductible de degr\'e $n_{j}$ $(\sum n_{j}=n+1)$ de $\mathrm{SL}(2,\mathbb{C})$ et $\chi_{j}$ un caract\`ere de $\mathbb{C}^{\times }$,
\textbf{unitaire} d'apr\`es (iii). Par ailleurs d\'efinissons $\psi_{0}^{\sigma}$ par

\begin{eqnarray*}
\psi_{0}^{\sigma}(z,s)=\psi_{0}(\bar z,s)\quad(z\in\mathbb{C}^{\times },s\in\mathrm{SL}(2,\mathbb{C})).
\end{eqnarray*}
La condition (i) implique que $\psi_{0}^{\sigma}$ est \'equivalente \`a la duale $\widetilde \psi_{0}$ de $\psi_{0}$, soit

\begin{equation}
\{r_{j},\chi_{j}^{\sigma}\}=\{r_{j},\chi_{j}^{-1}\}
\end{equation}
o\`u $\chi^{\sigma}(z)=\chi(\bar z)$.

\vskip0,2cm
Chaque bloc de (6.2.1) contribue une somme de caract\`eres \`a $\varphi_{\psi}\vert_{{}_{\mathbb{C}^{\times }}}$, de la forme

\begin{equation}
z\mapsto
\begin{pmatrix}
(z\bar z)^{\frac{n_{j}-1}{2}}\\
&\ddots\\
&&(z\bar z)^{\frac{1-n_{j}}{2}} 
\end{pmatrix}
\otimes \chi_{j}(z).
\end{equation}

\'Ecrivons, de la fa\c con usuelle,
\begin{equation*}
\chi_{j}=z^{p_{j}}\bar z^{q_{j}}\quad(p_{j}-q_{j}\in\mathbb{Z},\ p_{j}+q_{j}\in\sqrt{-1}\mathbb{R}).
\end{equation*}

D\'efinissons $(P_{1},\ldots P_{n+1}, Q_{1},\ldots Q_{n+1})$, modulo l'action diagonale de $\mathfrak{S}_{n+1}:\tau(P,Q)=(\tau P,\tau Q)$, par
\begin{equation}
(P_{j,k},Q_{j,k})=\Bigl(p_{j}+\frac{n_{j}+1-2k}{2},q_{j}+\frac{n_{j}+1-2k}{2}\Bigr)\ ,\ k=1,\ldots,n_{j}\ . 
\end{equation}

Nous allons d\'eduire de la Conjecture d'Arthur :

\begin{lem}
Selon l'identification naturelle avec $\mathbb{C}^{n+1}$ d'une sous-alg\`ebre de Cartan de $\mathfrak{g}$, le caract\`ere infinit\'esimal de toute repr\'esentation de $G$ associ\'ee \`a $\psi$ est \'egal \`a $P\in\mathbb{C}^{n+1}$ $($modulo $\mathfrak{S}_{n+1})$.
\end{lem}

Soit en effet $\varphi_{\psi}$ le param\`etre de Langlands, donn\'e par ses blocs (6.2.4). Noter que puisque $\sigma$ op\`ere sur $\widehat G$ par $g\mapsto w_{0}{}^{t}g^{-1}w_{0}^{-1}$, on a (cf. 6.2.3) :
\begin{equation*}
\{(P_{i},Q_{i})\}=\{(-Q_{i},-P_{i})\}
\end{equation*}
o\`u les $(P_{i},Q_{i})$ $(i=1,\ldots n+1)$ sont les couples $(P_{j,k},Q_{j,k})$.

Soit $G^{\ast}=U(a,b)$, o\`u $a+b=n+1$ et $a-b=0$ ou $1$, la forme int\'erieure quasi-d\'eploy\'ee de $G$. D'apr\`es la classification de Langlands (\pa 4.3-4.4 pour $U(n,1)$ ; \cite{Langlands}), $\varphi_{\psi}$ d\'efinit \textbf{toujours} une famille finie de repr\'esentations de   
$G^{\ast}$. (La condition (4.3.3) qui restreignait fortement les param\`etres pour $G=U(n,1)$ n'intervient pas ici, car tous les paraboliques ${}^{L}P$ de ${}^{L}G^{\ast}={}^{L}G$ proviennent de $G^{\ast}$.) Par ailleurs $G$ et $G^{\ast}$ ont la m\^eme complexification, et leurs sous-alg\`ebres de Cartan complexifi\'ees s'identifient donc canoniquement (modulo l'action de $\mathfrak{S}_{n+1}$). Il est alors implicite dans l'article d'Arthur \cite[\pa 4]{Arthur}, et il r\'esulte explicitement de la description par Adams, Barbasch et Vogan des paquets d'Arthur (\cite{AdamsBarbaschVogan}: voir en particulier Thm.~22.7 , Cor.~19.16 et les d\'efinitions pr\'ec\'edant celui-ci) que le caract\`ere infinit\'esimal d'une repr\'esentation \textbf{de} $G$ associ\'ee \`a $\psi$ est \'egal, modulo cette identification, \`a celui d'une repr\'esentation de $G^{\ast}$ associ\'ee par Langlands \`a $\varphi_{\psi}$. Nous calculons donc dans $G^{\ast}$.

Il suffit alors de suivre les constructions de Langlands \cite{Langlands}. Soit ${}^{L}M\subset {}^{L}G^{\ast}={}^{L}G$ un sous-groupe de Levi minimal contenant l'image de $\varphi_{\psi}$. Alors ${}^{L}M$ est le groupe dual d'un sous-groupe de Levi cuspidal $M^{\ast}$ de $G^{\ast}$, donc de la forme
\begin{equation*}
M^{\ast}\cong(\mathbb{C}^{\times })^{s}\times U(A,B)
\end{equation*} 
o\`u $A-B=a-b=0$ ou $1$, et $n+1=N+2s$, avec $N=A+B$ ; si $G^{\ast}$ est d\'efini par la matrice $J^{\ast}$ du \S 4.2, $M^{\ast}$ s'\'ecrit sous forme diagonale par blocs :
\begin{equation*}
M^{\ast}=\left\{\begin{pmatrix}
z_{1}\\
&\ddots\\ &&z_{s}\\&&&&u\\&&&&&\bar z_{s}^{-1}\\&&&&&&\ddots\\&&&&&&& \bar z_{1}^{-1} 
\end{pmatrix}\right\}\ ,
\end{equation*}
$u\in U(A,B)$.

D'apr\`es des calculs analogues \`a ceux du \pa 4.3 (et r\'esultant de \cite{Langlands}), un param\`etre $\varphi=\varphi_{\psi}:W_{\mathbb{R}} \rightarrow{}^{L}M$, minimal, se restreint alors \`a $W_{\mathbb{C}}\subset W_{\mathbb{R}}$ en
\begin{equation*}
\varphi^{0}:z\mapsto\begin{pmatrix}
\eta_{1}(z)\\
&\ddots\\ &&\eta_{s}(z)\\&&&(z/\bar z)^{\alpha_{1}}\\&&&&\ddots \\&&&&&(z/\bar z)^{\alpha_{N}}\\&&&&&& \eta_{s}(\bar z)^{-1}\\ &&&&&&&\ddots\\&&&&&&&& \eta_{1}(\bar z)^{-1} 
\end{pmatrix}\ .
\end{equation*}

Les caract\`eres qui figurent dans cette expression sont ceux de (6.2.3) -- pour tous les blocs -- et leurs exposants $z^{p}\bar z^{q}$ sont donc donn\'es par (6.2.4). Vu la condition de minimalit\'e, $\varphi^{0}$ d\'efinit une repr\'esentation de la s\'erie discr\`ete de $M^{\ast}$ (modulo le centre) et on a donc comme au \S 4.3: $\alpha_{j}\in\frac{N+1}{2}+\mathbb{Z}$ ; le bloc central de $\varphi^{0}$ d\'efinit une repr\'esentation $\delta $ de la s\'erie discr\`ete de $U(A,B)$ et les repr\'esentations associ\'ees de $G^{\ast}$ sont des sous-quotients de

\begin{equation*}
\ind_{P^{\ast}}^{G^{\ast}}(\eta \otimes \delta \otimes \1) =:I
\end{equation*}
o\`u $P^{\ast}=M^{\ast}N^{\ast}$ est le parabolique associ\'e et $\eta\otimes \delta $ d\'esigne
\begin{equation*}
\eta_{1}(z_{1})\cdots\eta_{s}(z_{s})\delta (u),\ (z,u)\in M^{\ast}. 
\end{equation*}

Une sous-alg\`ebre de Cartan (r\'eelle) de $M^{\ast}$ est donn\'ee par

\begin{equation*}
\mathfrak{h}_{0}=\left\{X=\begin{pmatrix}
z_{1}\\ & \ddots\\ && z_{s}\\ &&&\sqrt{-1}\ \theta_{1}\\ &&&&\ddots\\ &&&&&\sqrt{-1}\ \theta_{N}\\&&&&&& 
-\bar z_{s}\\ &&&&&&& \ddots\\ &&&&&&&& -\bar z_{1}  
\end{pmatrix}\right\}\ .
\end{equation*}
($z_i \in {\Bbb C}$, $\theta_j \in {\Bbb R}$.)

Posons $\eta_i = z^{a_i} \bar z ^{b_i}$ ($a_i - b_i \in {\Bbb Z}$) pour $i=1, \ldots , s$.
Le calcul usuel du caract\`ere infinit\'esimal d'une induite donne
\begin{equation*}
\lambda_{I}(X)=\sum_{i=1}^{s}a_{i}z_{i}+\sum_{i=1}^{s}b_{i}\bar z_{i}+\sqrt{-1}\sum_{j=1}^{N}\alpha_{j}\theta_{j}\ .
\end{equation*}
En fonction des valeurs propres $X_{i}$ de $X$, ceci s'\'ecrit
$$\lambda_{I}(X) = \sum_{i=1}^{s}a_{i}X_{i}-\sum_{i=1}^{s} b_{i}X_{i+s}+\sum_{j=1}^{N}\alpha_{j}X_{2r+j} .$$
Mais, revenant \`a l'expression de $\varphi^0$, on voit que le param\`etre $P$ associ\'e est donn\'e par
$$P= (a_1 , \ldots , a_s , \alpha_1 , \ldots , \alpha_N , -b_1 , \ldots , -b_s ).$$
D'o\`u le Lemme.

\vskip0,2cm
\textbf{Remarque.} Des calculs analogues s'appliquent aux groupes orthogonaux, et nous les admettrons (Lemme 6.3.1, Lemme 6.4.1).

\vskip0,2cm
Consid\'erons alors une repr\'esentation irr\'eductible unitaire $\pi$ de $G$ et son param\`etre de Langlands $\varphi_{\pi}:W_{\mathbb{R}}\rightarrow{}^{L}G$.

\vskip0,2cm
\textbf{Cas A. $\pi$ est une s\'erie discr\`ete ou une limite de s\'erie discr\`ete}.

\vskip0,2cm
Dans ce cas (cf. Ch.~4) son caract\`ere infinit\'esimal est

\begin{equation*}
P=(p_{1},\ldots p_{n+1})\ ,\ p_{i}\in\frac{1}{2}\mathbb{Z}\ ,\ p_{i}\equiv \frac{n}{2}\ \ [1]. 
\end{equation*}
Alors $\pi$ est toujours associ\'ee \`a un param\`etre d'Arthur temp\'er\'e : on prend $r=n+1$, $r_{j}$ triviale et
\begin{equation*}
\chi_{j}=z^{p_{j}}(\bar z)^{-p_{j}}\quad j=1,\ldots n+1\ ; 
\end{equation*} 

\vskip0,2cm
\textbf{Cas B. $\pi$ est un quotient de Langlands $J(\sigma,\chi)$ o\`u $\sigma$ est une repr\'esentation de $U(n-1)$ et $\chi$ un caract\`ere de} $\mathbb{C}^{\times }$ (cf. \pa 4.3).

\vskip0,2cm
\'Ecrivons $\chi=z^{\alpha}(\bar z)^{\beta}$, $\alpha-\beta\in\mathbb{Z}$.

Si la repr\'esentation $J(\sigma,\chi)$ est unitaire on sait que la donn\'ee $(\sigma,\chi)$ doit \^etre \textbf{hermitienne} \cite{KnappSpeh}. Le groupe de Weyl $W(G,M)$ s'identifie \`a $\mathbb{Z}/2\mathbb{Z}$, op\'erant trivialement sur $U(n-1)$ et par $z\mapsto\bar z^{-1}$ sur $\mathbb{C}^{\times }$. On a alors $\overline\chi=\chi^{-1}$ ou $\overline\chi=w\cdot\chi^{-1} $ $(w\neq 1)$ soit : $\chi$ unitaire ou $\overline\chi(z)=\chi(\bar z)$, {\it i.e} : $\alpha+\beta\in i\mathbb{R}$ ou $\alpha,\beta\in\mathbb{R}$.

\vskip0,2cm
Le second cas (si $\alpha+\beta\neq 0)$ correspond aux repr\'esentations non temp\'er\'ees.

Par ailleurs le caract\`ere infinit\'esimal de $\sigma$ est repr\'esent\'e par $P_{\sigma}\in\mathbb{C}^{n-1}$,

\begin{equation*}
P_{\sigma}=(m_{1},\ldots m_{n-1})\ ,\ m_{i}\in\frac{1}{2}\mathbb{Z}\ ,\ m_{i}\equiv\frac{n}{2}\ \ [1]\ ,\ m_{i+1}>m_{i}\ .
\end{equation*}
Celui de $J(\sigma,\chi)$ est donn\'e par
\begin{equation}
P=(P_{\sigma},\alpha,-\beta)
\end{equation}
qui doit \^etre de la forme (6.2.4).

On en d\'eduit donc que tous les coefficients $P_{i}$ de (6.2.4) doivent \^etre dans $\frac{1}{2}\mathbb{Z}$, sauf au plus deux.

Supposons alors que $\alpha\notin\tfrac{1}{2}\mathbb{Z}$ ($\Longleftrightarrow\beta\notin
\tfrac{1}{2}\mathbb{Z}$). Il y a deux possibilit\'es que nous examinons maintenant.

La premi\`ere possibilit\'e est que
\begin{equation}
\alpha=p_{j} \mbox{ et } -\beta=p_{j'} \mbox{ avec } j\neq j' \mbox{ et } n_{j}=n_{j'}=1\ . 
\end{equation}
D'apr\`es (6.2.2) $\chi_{j}=z^{p_{j}}(\bar z)^{q_{j}}$ appara\^{\i}t avec $\chi_{j}^{-\sigma}=z^{-q_{j}}(\bar z)^{-p_{i}}$. Donc $p_{j'}=-q_{j}$, d'o\`u $\beta=q_{j}$ et $\alpha+\beta=p_{j}+q_{j}\in i\mathbb{R}$ ; $\chi$ est alors unitaire.

La seconde possibilit\'e est que
\begin{equation}
\begin{split}
\alpha &=p_{j}+\tfrac{1}{2}\\
-\beta &=p_{j}-\tfrac{1}{2} .
\end{split}
\end{equation}

Alors $\alpha+\beta=1$, contrairement \`a l'hypoth\`ese puisque $\alpha-\beta\in\mathbb{Z}$. On en d\'eduit : 

\begin{lem}
Si $J(\sigma,\chi)$ est associ\'ee \`a un param\`etre d'Arthur, $\chi$ est unitaire ou de la forme
\begin{equation*}
z^{\alpha}\bar z^{\beta}\ ,\ \alpha,\beta\in\tfrac{1}{2}\mathbb{Z}\ .
\end{equation*}
\end{lem}

Nous allons en fait d\'emontrer un r\'esultat bien plus pr\'ecis. Rappelons (\pa 4.5) que nous ne nous int\'eressons qu'aux repr\'esentations $\sigma$ apparaissant dans $\Lambda^{\ast}\mathfrak{p}$. D'apr\`es le \pa 4.5, on a alors $\sigma=\sigma_{a,b}$ avec $a+b\leq n-1$, et

\begin{equation*}
P_{\sigma}=\Bigl(\underbrace{\frac{n}{2},\ldots,\frac{n}{2}-(b-1)}_{b};\frac{n}{2}-(b+1),\ldots,-\frac{n}{2}+(a+1); \underbrace{-\frac{n}{2}+(a-1),\ldots,-\frac{n}{2} }_{a}\Bigr)
\end{equation*}
o\`u chacune des suites d'entiers s\'epar\'es par ``;'' est une suite d'entiers cons\'ecutifs et d\'ecroissants.

De plus, si $J(\sigma,\chi)$ rencontre $\Lambda^{\ast}\mathfrak{p}$, on doit avoir, avec les notations pr\'ec\'edentes~:
\begin{equation*}
a-b=\alpha-\beta
\end{equation*}
(\pa 4.5). Donc $\chi$ est de la forme
\begin{equation}
\chi(z)=(z/\bar z)^{\frac{a-b}{2}}(z\bar z)^{s}
\end{equation}
et l'on s'int\'eresse au cas o\`u $J$ n'est pas temp\'er\'ee, donc $s\in\mathbb{R}-\{0\}$.

\begin{prop}
On suppose que $\sigma=\sigma_{ab}$ et que $\chi$ est de la forme (6.2.8) avec $s$ r\'eel $>0$. Alors $J(\sigma,\chi)$ est unitaire et associ\'ee \`a un param\`etre d'Arthur si, et seulement si :
\begin{equation}
s=\frac{n-(a+b)}{2}-k\ ,\ 0\leq k\leq\Bigl[\frac{n-(a+b)}{2}\Bigr]\ .
\end{equation}
\end{prop} 

Nous ne d\'emontrons que l'implication directe, la seule utilis\'ee dans ces notes. Soit
$s_{0}=\tfrac{n-(a+b)}{2}$. On sait que le dual unitaire de $G=U(n,1)$ est compl\`etement classifi\'e. Nous utilisons la description des r\'esultats par Knapp \cite[\pa 3]{Knapp} -- voir aussi \cite{Krajlevic}. Il en r\'esulte que $I(\sigma,\chi)$ est irr\'eductible et unitaire pour 
$s\in]0,s_{0}[$ et que $J(\sigma,\chi)$ est (irr\'eductible) unitaire pour $s=s_{0}$. Pour $s>s_{0}$ la repr\'esentation n'est pas unitaire. Si le caract\`ere infinit\'esimal de $J(s,\chi)$ provient d'un param\`etre d'Arthur, nous devons v\'erifier que $s$ v\'erifie (6.2.9).

Nous revenons \`a la description du caract\`ere infinit\'esimal $P$ donn\'ee avant le Lemme 6.2.2. On sait d\'ej\`a que $\alpha,\beta\in\tfrac{1}{2}\mathbb{Z}$, avec $\alpha\equiv\beta[1]$.

\begin{lem}
La relation $(6.2.9)$ est \'equivalente \`a
\begin{equation*}
\alpha,\beta\equiv\frac{n}{2}\ \ [1].
\end{equation*}
\end{lem}

En effet, puisque $\alpha-\beta=a-b$ et $a,b$ sont entiers :

\begin{equation*}
s=\frac{1}{2}(\alpha+\beta) =\frac{1}{2}(\alpha-\beta)+\beta 
=\frac{1}{2}(a-b)+\beta\equiv\frac{1}{2}(a+b)+\beta\ \ [1]
\end{equation*}
donc $s\equiv\tfrac{n-(a+b)}{2}\Longleftrightarrow\alpha,\beta\equiv\tfrac{n}{2}$.

Nous d\'emontrons donc que $\alpha,\beta\equiv\tfrac{n}{2}$ si $J(\sigma,\chi)$ est associ\'ee \`a un param\`etre d'Arthur. Pour ceci nous allons reformuler la param\'etrisation des repr\'esentations du groupe unitaire.

Les param\`etres de Langlands (\pa 4.3) ou d'Arthur (\pa 6.1) sont de la forme

\begin{eqnarray*}
\begin{array}{rccccl}
\psi : & W & & \rightarrow & & {}^{L}G\\
 &  & \searrow & \quad & \swarrow & \\
 & &  & \mathrm{Gal}(\mathbb{C}/\mathbb{R}) & & 
\end{array}
\end{eqnarray*}
avec $W=W_{\mathbb{R}}$ (Langlands) ou $W_{\mathbb{R}}\times \SL(2,\mathbb{C})$ (Arthur). (On notera simplement $S$ le groupe $\SL(2,\mathbb{C})$. Dans le cas de Langlands, on d\'esigne donc par $\psi$ le param\`etre not\'e $\varphi$ au \S 4.3). Le groupe $W$ est de la forme $W^{0}\amalg W^{0}j$, o\`u $j$ s'envoie sur $\sigma\in\Gal(\mathbb{C}/\mathbb{R})$ ; il contient un \'el\'ement, not\'e $-1\in W^{0}$, tel que $j^{2}=-1$. Enfin, on notera $z\mapsto\bar{z}$ l'action de $j$, par conjugaison, sur $W^{0}$.

On suppose que $G$ est un groupe unitaire de rang $m$ ; ${}^{L}G$ est d\'ecrit au \S 4.2. La donn\'ee de $\psi$ se r\'eduit aux donn\'ees suivantes : 

\vskip0,2cm
\noindent (Aa)\qquad $\psi_{0}:W_{0}\rightarrow\widehat {G}=\GL(m,\mathbb{C})$

\vskip0,2cm
\noindent(Ab)\qquad $\psi(j)=(g,\sigma)\ ,\ g\in\GL(m,\mathbb{C})$

\noindent avec les restrictions suivantes :

\vskip0,2cm
\noindent(A1)\qquad $(g,\sigma)^{2}=\psi_{0}(-1)$ soit :

\quad \quad\ \  $gw_{0}{}^{t}g^{-1}w_{0}^{-1}=\psi_{0}(-1)$

\noindent(A2)\qquad  $(g,\sigma)\psi_{0}(z)(g,\sigma)^{-1}=\psi_{0}(\bar{z})$, soit :

\qquad
\ \ $gw_{0}{}^{t}\psi_{0}(z)^{-1}w_{0}^{-1}= g^{-1}\psi_{0}(\bar{z}) \ (z\in W_{0}) .$

Posant $h=gw_{0}$, et en tenant compte de l'identit\'e ${}^t w_{{_{0}}} =(-1)^{m+1}w_{0}$, les conditions se r\'e\'ecrivent :

\noindent(A1)\qquad $h\,{}^{t}h^{-1}=(-1)^{m+1}\psi_{0}(-1)$

\vskip0,2cm
\noindent(A2)\qquad $h\,{}^{t}\psi_{0}(z)^{-1}h^{-1}=\psi_{0}(\bar{z})$\ $(z\in W_{0}).$

Notons $V=\mathbb{C}^{m}$ et soit $V^{\ast}$ l'espace dual, lui aussi identifi\'e \`a $\mathbb{C}^{m}$ par la base duale. On peut alors consid\'erer $\psi_{0}$, $h$ comme les donn\'ees suivantes :

\vskip0,2cm
\noindent(Ba)\qquad $\psi_{0}$ = repr\'esentation de $W_{0}$ sur $V$.

\vskip0,2cm
\noindent(Bb)\qquad $H$ = isomorphisme $V^{\ast}\rightarrow V$\ .

\vskip0,2cm
Soit $H^{\ast}:V^{\ast}\rightarrow V$ l'adjoint. Les contraintes sont alors :

\vskip0,2cm
\noindent(B1)\qquad $H\psi_{0}^{\ast}(z)H^{-1}=\psi_{0}(\bar{z})$ \quad$(z\in W_{0})$

\vskip0,2cm
\noindent(B2)\qquad $H(H^{\ast})^{-1}=(-1)^{m+1}\psi_{0}(-1)$.

\vskip0,2cm
R\'eciproquement, une donn\'ee (B) d\'etermine, apr\`es choix d'une base de $V$, une donn\'ee (A).

Nous appliquons ceci \`a une donn\'ee d'Arthur, avec $m=n+1$. Donc $W=W_{\mathbb{R}}\times S$ ; la repr\'esentation $\psi_{0}$ de $W_{0}=\mathbb{C}^{\times }\times S$ est semi-simple, unitaire en restriction \`a $\mathbb{C}^{\times }$. Les consid\'erations du d\'ebut du \pa 6.2 s'appliquent. On supposera

\begin{eqnarray}
\begin{array}{l}
\textbf{Les coefficients}\ P_{i}\ \mathbf{de}\  P\in\mathbb{C}^{n+1 } - \mathrm{cf. (6.2.4) } - \textbf{appartiennent \`a }\\ 
\frac{1}{2}\mathbb{Z},
 \mathbf{et}\ (n-1)\ \textbf{d'entre eux sont } \equiv\frac{n}{2} \mathrm{[1] }.  
\end{array}
\end{eqnarray}

Cette condition est vraie dans notre cas, cf. (6.2.5). La Proposition 6.2.3 sera donc d\'emontr\'ee modulo le

\begin{lem}
Il n'existe pas de donn\'ee $(B)$ v\'erifiant $(6.2.10)$, et telle que les deux autres coefficients $\alpha,-\beta$ de $P$ v\'erifient $\alpha,\beta\equiv\frac{n+1}{2}$ $[1]$ et soient distincts.
\end{lem}

(Remarquer que $\alpha=-\beta$ correspond \`a $s=0$, donn\'ee permissible pour une repr\'esentation d'Arthur.)

Soient $\chi_{j}$ les caract\`eres apparaissant dans l'expression (6.2.1) de $\psi_{0}$. Alors $\chi_{j}=z^{p}\bar{z}^{q}$, avec $p\in\frac{1}{2}\mathbb{Z}$, et $p-q\in\mathbb{Z}$ et $p+q\in i\mathbb{R}$. Donc $\chi_{j}=(z/\bar{z})^{p_{j}}$, $p_{j}\in\frac{1}{2}\mathbb{Z}$.

Soit par ailleurs $\varphi_{\psi}$ la repr\'esentation de $W_{\mathbb{R}}$ d\'eduite de $\psi$ et $\varphi_{\psi}^{0}$ sa restriction \`a $\mathbb{C}^{\times }$. On a donc \textbf{sous} $\varphi_{\psi}^{0}$ :

\begin{equation*}
V=\bigoplus_{\eta}V(\eta)
\end{equation*}
o\`u $\eta=z^{p}(\bar z)^{q}$ d\'ecrit les caract\`eres de $\mathbb{C}^{\times }$. Soit

\begin{equation*}
T=\bigoplus_{\eta'}V(\eta')
\end{equation*}
o\`u $\eta'$ d\'ecrit les caract\`eres tels que $p\notin\tfrac{n}{2}+\mathbb{Z}$. Par hypoth\`ese $T$ est de dimension $2$.

\begin{lem}
$T$ est stable par $W_{0}=\mathbb{C}^{\times }\times S$ $($pour la repr\'esentation $\psi_{0})$
\end{lem}

En effet $V$ se d\'ecompose canoniquement sous $W_{0}$ en somme d'espaces isotypiques, eux-m\^emes multiples de repr\'esentations $\chi\otimes \rho$ o\`u $\chi=(z/\bar z)^{p}$ et $\rho$ est une repr\'esentation irr\'eductible de $S$. Sous 
$\varphi_{\psi}^{0}$, $\chi\otimes \rho$ donne des caract\`eres d'exposants $p'=p+\tfrac{r-1}{2}$ + (entier) o\`u $r=\dim\rho$. Ils ont tous la m\^eme parit\'e (mod 1), d'o\`u le Lemme~6.2.6.

Soit $T^{\ast}\subset V^{\ast}$ d\'efini de fa\c con analogue. D'apr\`es (B1) $H$ envoie $T$ sur $T^{\ast}$, et c'est alors un isomorphisme. Donc la donn\'ee $B$ d\'efinie par $(T,\psi_{0},H)$ v\'erifie les conditions (B) ; dans (B2) $\psi_{0}(-1)$ est bien s\^ur un endomorphisme de $T$. Il faut alors distinguer deux cas :

\vskip0,2cm
\noindent\textbf{Cas} I : $T\cong\chi\otimes \rho_{2}$ o\`u $\rho_{2}$ est irr\'eductible.

On peut consid\'erer $H:T\rightarrow T^{\ast}$ comme une forme bilin\'eaire sur $T$. D'apr\`es (B1) $H$ est invariante par $S$, donc antisym\'etrique. D'apr\`es (B2) on a donc :

\begin{equation*}
-1=(-1)^{n}\psi_{0}(-1)=(-1)^{n}\chi(-1)
\end{equation*}    
ce qui implique que $\chi=(z/\bar z)^{p}$ avec $p\equiv\tfrac{n+1}{2}$ [1]. Mais alors les exposants $\{p+\tfrac{1}{2},p-\tfrac{1}{2}\}$ de $\varphi_{\psi}$ apparaissant dans $T$ appartiennent \`a $\tfrac{n}{2}+\mathbb{Z}$, contradiction.

\vskip0,2cm
\noindent\textbf{Cas} II : $T\cong\chi_{1}\oplus\chi_{2}$, $S$ op\'erant trivialement.

Dans ce cas $\chi_{1}=(z/\bar z)^{\alpha}$, $\chi_{2}=(z/\bar z)^{\beta}$ et donc $\chi_{1}\neq\chi_{2}$ par hypoth\`ese. Soient $(e,f)$ la base de $T$ associ\'ee et $e^{\ast}$, $f^{\ast}$ la base duale de $T^{\ast}$. D'apr\`es (B1) $H:V^{\ast}\rightarrow V$ s'\'ecrit dans ces bases sous forme diagonale, $H=\begin{pmatrix}a\\&b
\end{pmatrix}$. Alors $H(H^{\ast})^{-1}=1$ d'o\`u (B2) $\psi_{0}(-1)=(-1)^{n}$. Or
$\psi_{0}(-1)=(-1)^{2\alpha}(=(-1)^{2\beta})$ donc $\alpha\equiv\beta\equiv\tfrac{n}{2}$ [1].

Ceci ach\`eve la d\'emonstration du Lemme 6.2.5.

En combinant la Proposition 6.2.3, le calcul des repr\'esentations $\sigma$ apparaissant dans $\Lambda^{pq}$ (\pa 4.5) et le calcul des valeurs propres correspondantes (Prop. 4.5.1), on en d\'eduit :

\begin{cor}
Soit $\pi$ une  repr\'esentation unitaire irr\'eductible de $G$ telle que Hom$_K (\Lambda^p \mathfrak{p}^+ \otimes \Lambda^q \mathfrak{p}^- , \pi ) \neq 0$ est associ\'ee \`a un 
param\`etre d'Arthur. Alors la valeur propre $\lambda$ de l'op\'erateur de Casimir dans ${\cal H}_{\pi}$ v\'erifie~:
$$\lambda = [n-a-b]^2 - [n-a-b-2k]^2 , $$
pour $a \leq p$, $b \leq q$, $(p-q)-(a-b) \in \{ -1 , 0 ,1\}$ et $0 \leq k \leq \left[ \frac{n-(a+b)}{2} \right]$.
\end{cor}

\markboth{CHAPITRE 6. CONJECTURES D'ARTHUR}{6.3. $G=SO(n,1)$, $n$ IMPAIR}

\section{$G=SO(n,1)$, $n$ impair}

Le groupe $G$ n'est pas connexe : son sous-groupe compact maximal est $K=S(O(n)\times \{ \pm 1\} )\cong O(n)$. Puisque les Conjectures d'Arthur ne s'appliquent \textbf{a priori} qu'aux groupes alg\'ebriques (et que nous pr\'ef\'erons \'eviter de consid\'erer le groupe des spineurs), nous consid\'erons n\'eanmoins, pour l'instant, $G$. Posons $n=2\ell -1$, de sorte que $\ell$ est le rang absolu de $G$. 

Pour $n$ impair, le groupe dual ${}^{L}G$ a pour composante neutre

\begin{equation}
\widehat {G}=SO(n+1,\mathbb{C}). 
\end{equation}

Par ailleurs, si $\ell$ est \textbf{pair}, $G$ est forme int\'erieure de $G^{\ast}=SO(\ell-1,\ell+1)$ ; si $\ell$ est \textbf{impair} $G$ est forme int\'erieure de $G^{\ast}=SO(\ell,\ell)$, le groupe $G^{\ast}$ \'etant toujours quasi-d\'eploy\'e, et d\'eploy\'e dans le second cas. On en d\'eduit que

\begin{align} \index{groupe dual}
&{}^{L}G=SO(n+1,\mathbb{C})\times \Gal(\mathbb{C}/\mathbb{R})\ \textrm{(produit direct)}, & \ell\ \textrm{impair}\\
&{}^{L}G=SO(n+1,\mathbb{C})\rtimes \Gal(\mathbb{C}/\mathbb{R}) , &\ell\ \textrm{pair.}\kern0,5cm
\end{align} 

Dans le second cas, l'action de $\sigma\neq 1$ respecte un sous-groupe de Borel, par exemple les matrices triangulaires sup\'erieures dans $SO(2m+2,\mathbb{C})\cong SO(2m+2,\widehat {J})$ o\`u $\widehat {J}$ est la forme de matrice

\begin{equation*}
\begin{pmatrix}
0 &1_{\ell}\\
1_{\ell} &0\\
\end{pmatrix}\ .
\end{equation*}
Elle doit aussi respecter un \'epinglage (voir \pa 4.2).

Supposons $G$ d\'efini pour la forme de matrice
 
\begin{equation*}
J=\begin{pmatrix}
\begin{tabular}{ c  | cc}
$1_{n-1} $\\ 
  \\ \hline 
 &&1 \\
&1
\end{tabular}
\end{pmatrix}\ .
\end{equation*}
Le sous-groupe parabolique propre (minimal) de $G$ s'\'ecrit

\begin{equation*}
P=MN
\end{equation*} 
o\`u $M\subset G$ stabilise la d\'ecomposition correspondante $\mathbb{R}^{n+1}=\mathbb{R}^{n-1}\oplus\mathbb{R}^{2} $ ; donc

\begin{equation*}
M=S(O(n-1)\times O(1,1))
\end{equation*}  
o\`u 
\begin{equation}
SO(1,1)=\left\{\begin{pmatrix}
a &0\hfill\\0 &a^{-1}
\end{pmatrix}, a\in\mathbb{R}^{\times }\right\}
\end{equation}

\noindent et 
\begin{equation}
O(1,1)=<SO(1,1),\varepsilon>, \varepsilon=\begin{pmatrix}&1\\1
\end{pmatrix}.
\end{equation}

\noindent On pose 
\begin{equation}
A=\left\{\begin{pmatrix}a\\&a^{-1}\end{pmatrix}:a>0\right\}\subset O(1,1)\subset M\ ,
\end{equation}

\noindent et 
\begin{equation}
M={}^{0}MA
\end{equation}

\vskip0,2cm
\noindent o\`u la composante neutre de ${}^0 M$ est ${}^{0}M^{0}\cong SO(n-1)$ et ${}^{0}M=<{}^{0}M^{0},\varepsilon,\eta>$,

\begin{equation*}
\eta=\begin{pmatrix}-1\\&-1\end{pmatrix}\in SO(1,1).
\end{equation*}

Si $\tau$ est une repr\'esentation irr\'eductible de ${}^{0}M$ et $s\in\mathbb{C}$ soit

\begin{equation*}
I(\tau,s)=\ind_{{}^0MAN}^{G}(\tau\otimes e^{s}\otimes 1).
\end{equation*}
On utilise la classification de Langlands (cf. Ch. 4 pour $U(n,1)$, \cite{Langlands}). Pour $n$ impair $G$ n'a pas de s\'erie discr\`ete. Alors toute repr\'esentation admissible irr\'eductible de $G$ est de la forme suivante~:

\begin{equation}
\pi\ \textrm{est un sous-module de}\ I(\tau,s); \tau\in\widehat {{}^{o}M}, s\in i\mathbb{R}\ .
\end{equation}

\begin{equation}
\begin{split}
&\pi\ \textrm{est l'unique quotient irr\'eductible}\ J(\tau,s)\ \textrm{de}\ I(\tau,s)\\ &\textrm{o\`u}\ \rho\in\widehat {{}^{o}M}\ \textrm{et}\ Re(s)>0\ .
\end{split}
\end{equation}

Dans le cas (6.3.8) on sait en fait que $I(\tau,s)$ est irr\'eductible ; sa restriction \`a $G^{0}$ est en fait irr\'eductible, cf. \cite{KnappZuckerman}. \footnote{Nous n'utiliserons pas ce fait, sauf en notant dans tous les cas $J(\tau,s)$ la repr\'esentation $\pi$ pour uniformiser les notations.}

Par ailleurs, une sous-alg\`ebre de Cartan $\mathfrak{h}_{0}$ de $\mathfrak{g}_{0}$ est donn\'ee par les matrices

\begin{equation}
X=\begin{pmatrix}
\begin{array}{ ccccc  | cc}
&-x_{1}&&&&&\\
x_{1}&&&&& \\
&&\ddots&&&\\
&&&&-x_{ \ell-1}&\\
&&&x_{\ell-1}&&\\
 \hline
&&&&&x_{\ell}\\
&&&&&&-x_{ \ell}\\
\end{array}
\end{pmatrix}(x_{i}\in\mathbb{R})\ .
\end{equation}

Posons $y_{i}=\sqrt{-1}\ x_{i}$ $(i\leq\ell-1)$, $y_{\ell}=x_{\ell}$. Les coordonn\'ees $y$ donnent un isomorphisme
\begin{equation*}
\mathfrak{h}=\mathfrak{h}_{0}\otimes \mathbb{C}^{\ell}
\end{equation*} 
pour lequel les racines dans $\Delta(\mathfrak{g},\mathfrak{h})$ sont r\'eelles. Le groupe de Weyl $W$ est $\mathfrak{S}_{\ell}\ltimes \{ \pm 1\}^{\ell-1}$, $\{\pm1\}^{\ell-1}$ \'etant le sous-groupe de $\{\pm 1\}^{\ell}$, op\'erant diagonalement, d\'efini par $\Pi s_{i}=1$.

Soit ${}^{0}\mathfrak{h}=\mathbb{C}^{\ell-1}$, naturellement plong\'e dans $\mathfrak{h}$. C'est une alg\`ebre de Cartan pour ${}^{0}\mathfrak{m}$. Soit $\lambda_{\tau}\in{}^{0}\mathfrak{h}^{\ast}$ le caract\`ere infinit\'esimal de $\tau$, d\'efini modulo $\mathfrak{G}_{\ell-1}\times (\pm 1)^{\ell-2}$. Celui de $I(\tau,s)$ est alors

\begin{equation}
\lambda_{\tau,s} =(\lambda_{\tau},s)\in\mathbb{C}^{\ell}\ .
\end{equation} 
Si $J(\tau,s)$ est unitaire on a $s\in i\mathbb{R}$ ou $s\in\mathbb{R}$ (ceci n'\'etant possible que pour certaines valeurs de $\tau$). 

Soit alors
\begin{equation*}
\psi:W_{\mathbb{R}}\times \SL(2,\mathbb{C})\rightarrow{}^{L}G
\end{equation*}
un param\`etre d'Arthur pour $G$, et consid\'erons

\begin{equation*}
\psi_{0}:\mathbb{C}^{\times }\times \SL(2,\mathbb{C})\rightarrow\widehat G\subset \GL(2\ell ,\mathbb{C}).
\end{equation*}
On peut d\'efinir alors, comme dans le \pa 6.2,

\begin{equation*}
P=(P_{1},\ldots,P_{2\ell} )\in\mathbb{C}^{2\ell}\ ;
\end{equation*}
il est sp\'ecifi\'e (modulo $\mathfrak{S}_{2\ell}$) par le fait que $\varphi_{\psi|_{\mathbb{C}^{\times }}}$ (\pa 6.1) s'\'ecrit

\begin{equation*}
z\mapsto(z^{P_{i}}(\bar z)^{Q_{i}}).
\end{equation*}

Puisque $\psi$ se factorise par $\widehat G \subset \GL(2\ell ,\mathbb{C})$, on a $\varphi_{\psi}\cong\widetilde\varphi_{\psi} $ ; donc $P$ est conjugu\'e \`a $(-P)$ par $\mathfrak{S}_{2\ell}$ et on en d\'eduit que modulo $\mathfrak{S}_{2\ell}$~:

\begin{equation*}
P=(P_{1},\ldots, P_{\ell},-P_{1},\ldots, -P_{\ell})=(P',P'').
\end{equation*}
Alors $P'$ est uniquement d\'efini modulo $W$, et l'on a, en supposant que $\pi$ est dans le paquet d'Arthur associ\'e \`a $\psi$ :

\begin{lem}
Le caract\`ere infinit\'esimal de $\pi$ est param\'etr\'e par $P'$.
\end{lem}

On l'admettra (voir la d\'emonstration du Lemme 6.2.1).

Nous supposons maintenant $\pi$ donn\'ee par (6.3.8) ou (6.3.9). Puisque $\tau$ est une repr\'esentation irr\'eductible de ${}^{0}M$, son caract\`ere infinit\'esimal est celui de sa restriction \`a ${}^{0}M^{0}=SO(n-1)$ ; on a alors, le tore maximal de $SO(n-1)$ \'etant param\'etr\'e par (6.3.10) :

\begin{align*}
&\lambda_{\tau}=\mu_{\tau}+\rho\\
&\rho=(\ell-2,\ell-1,\ldots, 0)\in\mathbb{C}^{\ell-1}\\
&\mu_{\tau}=(\mu_{1}, \ldots,\mu_{\ell-1}),\ \mu_{1}\in\mathbb{Z},\ \mu_{1}\geq\cdots\geq \mu_{\ell-1}, \mu_{\ell-2}+\mu_{\ell-1}\geq 0 
\end{align*}
de sorte que $\lambda_{\tau}$ v\'erifie :

\begin{equation}
\begin{split}
&\lambda_{\tau}=(\lambda_{1}, \ldots,\lambda_{\ell-1})\\
&\lambda_{i}\in\mathbb{Z},\ \lambda_{1}>\cdots>\lambda_{\ell-1}\ ,\ \lambda_{\ell-2}+\lambda_{\ell-1}>0\ . 
\end{split}
\end{equation}

Par ailleurs le caract\`ere infinit\'esimal de $\pi$ est param\'etr\'e par

\begin{equation}
P'=(\lambda_{\tau},s).
\end{equation}

Or le param\`etre $(P',-P')$ associ\'e \`a

\begin{equation*}
\psi:W_{\mathbb{C}}\times \SL(2,\mathbb{C})\rightarrow\widehat G\rightarrow\GL(n+1,\mathbb{C})
\end{equation*}
v\'erifie les conditions du \pa 6.2 ; on a~:

\begin{equation}
(P',-P')=(p_{j}+\tfrac{n_{j}+1-2k}{2})\ \mod\ \mathfrak{S}_{n+1}\ .
\end{equation}

Nous allons utiliser ces relations pour d\'emontrer :

\begin{lem}
Si le caract\`ere infinit\'esimal d'une repr\'esentation unitaire $\pi=J(\tau,s)$ est associ\'e \`a un param\`etre d'Arthur, $s\in i\mathbb{R}$ ou $s\in \mathbb{Z}$.
\end{lem}

La d\'emonstration est un peu compliqu\'ee. Nous d\'emontrons d'abord sous les m\^emes hypoth\`eses~:

\begin{lem}
$s\in i\mathbb{R}$ ou $s\in \tfrac{1}{2}\mathbb{Z}$\ . 
\end{lem}

Supposons en effet $s\notin \tfrac{1}{2}\mathbb{Z}$. D'apr\`es (6.3.13) et (6.3.14),

\begin{equation}
s=p+\tfrac{m+1-2k}{2}\ \textrm{pour un}\ p=p_{j}\ ,\ m=n_{j}\ .
\end{equation}
Donc $p\notin \tfrac{1}{2}\mathbb{Z}$. Par ailleurs la repr\'esentation $\chi\otimes r_{j}$ de $\mathbb{C}^{\times }\times \SL(2,\mathbb{C})$, o\`u $\chi=z^{p}(\bar z)^{q}$, appara\^{\i}t avec sa duale. Puisque $p\neq 0$, $\chi\neq\widetilde\chi=\chi^{-1}$. Si $m>1$, (6.3.14) contient au moins quatre coordonn\'ees $\notin \tfrac{1}{2}\mathbb{Z}$, contrairement \`a (6.3.13). Donc $m=1$. Enfin, on a $p+q\in i\mathbb{R}$ ($\chi$ est unitaire), $p-q\in\mathbb{Z}$ donc $p\in\tfrac{1}{2}\mathbb{Z}+i\mathbb{R}$. Or $p=s$ ou $-s$ et $s\in i\mathbb{R}\cup\mathbb{R}$ car $\pi$ est unitaire. D'o\`u $s\in i\mathbb{R}$.

\medskip

Pour d\'emontrer le Lemme 6.3.2, nous pouvons maintenant supposer que $s\in \tfrac{1}{2}\mathbb{Z}$ ; nous allons d\'eriver une contradiction. Noter que si $s\in \tfrac{1}{2}\mathbb{Z}$, $p\in \tfrac{1}{2}\mathbb{Z}$ ; puisque $p+q\in i\mathbb{R}$ et
$p-q\in \mathbb{Z}$ on a $p=-q$. Par ailleurs l'argument pr\'ec\'edent montre toujours que $m=1$.

Supposons d'abord $\ell$ \textbf{impair}, de sorte (6.3.2) que

\begin{equation*}
{}^{L}G=SO(2\ell,\mathbb{C})\times \Gal(\mathbb{C}/\mathbb{R}).
\end{equation*}

Vu comme repr\'esentation de $\mathbb{C}^{\times }\times \SL(2)$, on a

\begin{align*}
&\psi=\chi\oplus\chi^{-1}\oplus\sum_{i}\chi_{i}\otimes r_{i},\\
&\chi(z)=(z/\bar z)^{p}.
\end{align*}

Pour tout $i$, on a de plus $p_{i}+\tfrac{n_{i}+1}{2}\in\mathbb{Z}$, o\`u $\chi_{i}(z)=\chi^{p_{i}}(\bar z)^{q_{i}}$, toujours d'apr\`es (6.3.13).

Soit $j$ l'\'el\'ement ext\'erieur de $W_{\mathbb{R}}$ et $\psi(j)=(x,\sigma)\in{}^{L}G$. Alors

\begin{equation}
\psi(\bar z,s)=x\psi(z,s)x^{-1}\quad(z\in\mathbb{C}^{\times },s\in\SL(2)).
\end{equation}

Notons $\psi_{f}$ la partie de $\psi$ d'exposants fractionnaires~:

\begin{equation*}
\psi_{f}(z)=(\chi(z),\chi^{-1}(z),1,\ldots,1 ).
\end{equation*}

On d\'eduit facilement de (6.3.16), en consid\'erant les classes d'isotypie de $\psi$ vue comme repr\'esentation de degr\'e $2\ell$, que

\begin{equation}
\psi_{f}(\bar z)=x\, \psi_{f}(z)x^{-1}.
\end{equation}

En utilisant la forme explicite de $\widehat G$ donn\'ee apr\`es (6.3.6), on peut \'ecrire \`a conjugaison pr\`es $\psi_{f}$ sous la forme

\begin{equation*}
z\mapsto\begin{pmatrix}
\chi\\&1\\&&\ddots\\&&&1\\&&&&\ddots\\&&&&&1\\&&&&&&\chi^{-1}
\end{pmatrix}
\end{equation*}
o\`u on a abr\'eg\'e $\chi(z)$ en $\chi$. L'\'equation (6.3.17) montre que $x$ pr\'eserve la d\'ecomposition isotypique associ\'ee de $\mathbb{C}^{2\ell }$, donc s'\'ecrit $(x',x'')$ o\`u $x'\in O(2)$, $x''\in O(2\ell-2)$ et $\det(x')\det(x'')=1$.

Alors l'\'equation (6.3.17) implique

\begin{equation*}
x'\begin{pmatrix}
\chi\\&\chi^{-1}
\end{pmatrix}(x')^{-1}=\begin{pmatrix}
\chi^{-1}\\&\chi
\end{pmatrix}
\end{equation*}  
d'o\`u $x'=\begin{pmatrix}
&u^{-1}\\u
\end{pmatrix}\in O(2)-SO(2)$ ; $u\in\mathbb{C}^{\times }$. Mais alors $(x')^{2}=1$, contrairement \`a la condition

\begin{equation*}
\psi(j)^{2}=\psi(j^{2})=\psi(-1)=\begin{pmatrix}
-1\\&1\\&&\ddots\\&&&1\\&&&&\ddots\\&&&&&1\\&&&&&&-1
\end{pmatrix}.
\end{equation*}
(Car $p \in \frac12 {\Bbb Z}$, $p \notin {\Bbb Z}$ donc $\chi (-1) = -1$.)

Le calcul est similaire mais plus compliqu\'e si $\ell$ est pair. Dans ce cas $\widehat G=SO(2\ell,\mathbb{C})$,

\begin{equation*}
{}^{L}G=SO(2\ell,\mathbb{C})\rtimes\Gal(\mathbb{C}/\mathbb{R}) 
\end{equation*}
et nous devons pr\'eciser l'action de $\sigma$ donn\'ee par Langlands. Avec le m\^eme choix de la forme quadratique $\widehat J$, $\widehat G$ a pour sous-groupe de Borel $\widehat B$ ses matrices triangulaires sup\'erieures ; l'action de $\sigma$ doit pr\'eserver $\widehat G$, $\widehat B$, le tore diagonal $\widehat T$ ainsi qu'un \'epinglage (\pa 4.2), et \^etre non triviale. Sur $\widehat T$, on a alors

\begin{equation*}
\sigma:\begin{pmatrix}
x_{1}\\&\ddots\\&&x_{\ell}\\&&&x_{\ell}^{-1}\\&&&&\ddots\\&&&&&x_{1}^{-1}
\end{pmatrix}\mapsto\begin{pmatrix}
x_{1}\\&\ddots\\&&x_{\ell}^{-1}\\&&&x_{\ell}\\&&&&x_{\ell-1}^{-1}\\&&&&&\ddots\\&&&&&&x_{1}^{-1}
\end{pmatrix}
\end{equation*} 
de sorte que $\sigma$ \'echange les deux derni\`eres racines $\alpha_{\ell-1}=x_{{\ell-1}}$ $x_{\ell}^{-1}$ et $\alpha_{\ell}=x_{\ell-1}x_{\ell}$ du diagramme de Dynkin. Un calcul simple montre que

\begin{equation*}
\sigma(g)=w\ g\ w^{-1}\quad(g\in\widehat G)
\end{equation*}

\begin{equation*}
w=\begin{pmatrix}
1_{\ell-1}\\&\begin{matrix}0&1\\1&0\end{matrix}\\&&1_{\ell-1}
\end{pmatrix}
\end{equation*}
est l'automorphisme cherch\'e. Noter que $w\in O (\widehat J)$.

Consid\'erons alors, comme pour $\ell$ impair, l'homomorphisme $\psi_{f}$ :

\begin{eqnarray*}
\begin{array}{ccccc}
W_{\mathbb{R}}&  & \rightarrow & & {}^{L}G\\
&  \searrow & \quad& \swarrow & \\
& &\mathrm{Gal}(\mathbb{C}/\mathbb{R}) & &  
\end{array}
\end{eqnarray*}

Quitte \`a conjuguer $\psi_{f}$ par un \'el\'ement de $\widehat G$, on peut supposer que
$\psi_{f}(\mathbb{C}^{\times })\subset \widehat T$, puis que

\begin{equation*}
\psi_{f}(z)=\begin{pmatrix}
\chi\\ &1\\&&\ddots\\ &&&1\\&&&&\chi^{-1}
\end{pmatrix}
\end{equation*}
avec la notation pr\'ec\'edente. Soit $\psi_{f}(j)=(x,\sigma)$. Alors

\begin{equation*}
\begin{split}
\psi_{f}(\bar z)=\begin{pmatrix}
\chi^{-1}\\&1\\&&\ddots\\ &&&1\\&&&&\chi
\end{pmatrix}
&=(x,\sigma)\psi_{f}(z)(x,\sigma)^{-1}\\
&=x\ w\ \psi_{f}(z)w^{-1}x^{-1}\\
&=x\ \psi_{f}(z)x^{-1}\ .
\end{split}
\end{equation*}

Ceci implique comme dans le cas impair que $x\in O(2)\times O(2\ell-2)$ et que sa composante dans $O(2)$ est de la forme $\begin{pmatrix}
&z^{-1}\\z\end{pmatrix}$. Comme pr\'ec\'edemment, ceci contredit le fait que $\chi(j^{2})=\chi(-1)=-1$.

\medskip

D'apr\`es la classification du dual unitaire de $G$, on d\'eduit alors du Lemme~6.3.2~:

\begin{prop} \index{caract\`ere infinit\'esimal} \index{param\`etre d'Arthur}
Toute repr\'esentation unitaire $\pi$ de $G$ dont le caract\`ere infinit\'esimal est associ\'e \`a un param\`etre d'Arthur est de la forme $J(\tau,s)$ o\`u

\begin{itemize}
\item[$(i)$] $s\in i\mathbb{R}$
\end{itemize}
ou 
\begin{itemize}
\item[$(ii)$] $s\in \mathbb{Z}$\ .
\end{itemize}
\end{prop}

Pour $\tau$ fix\'e, et dans le cas (ii), la longueur du param\`etre $s\in\mathbb{Z}$ est bien s\^ur control\'ee par la description, connue, des s\'eries compl\'ementaires (voir par exemple \cite{KnappSpeh} pour un guide sur la litt\'erature concernant les duaux unitaires). Si $\tau=\1$ (repr\'esentation triviale), $\pi_{\tau,s}$ est unitaire si $|s|\leq\ell-1$ ; pour les repr\'esentations d'Arthur, on a donc $s\in i\mathbb{R}$ ou $\pm s=1,2,\ldots,\ell-1=\tfrac{n-1}{2} $ ; le point $s=\ell-1$ est la demi-somme des racines $\rho$ de $N$ dans $A$, et correspond \`a la repr\'esentation triviale. La Conjecture 6.1.2 implique donc une conjecture de Burger et Sarnak \cite{BurgerSarnak, BurgerLiSarnak}~:

\begin{conj} \index{Conjecture de Burger et Sarnak}
\begin{itemize}
\item[$(i)$] Si une repr\'esentation sph\'erique $J(1,s)$ appara\^{\i}t dans $L^{2}(\Gamma\ba G)$ pour un sous-groupe de congruence, $s\in i\mathbb{R}$ ou $s\in\{1,\ldots, \tfrac{n-1}{2}\}$.

 \item[$(ii)$] En particulier les valeurs propres du laplacien associ\'ees v\'erifient
 
\quad\    $\lambda\geq\Bigl(\tfrac{n-1}{2}\Bigr)^{2}$

  ou $\lambda=\Bigl(\tfrac{n-1}{2}\Bigr)^{2}-k^{2}$, $k=1,\ldots,\tfrac{n-1}{2}$.
 
 La premi\`ere valeur propre $\neq 0$ est donc $\lambda_{1}=n-2$. 
 \end{itemize}
\end{conj}

La partie (ii) r\'esulte de (i) et des calculs suivant le Th\'eor\`eme~2.3.2.

\markboth{CHAPITRE 6. CONJECTURES D'ARTHUR}{6.4. $G=SO(n,1)$, $n$ PAIR}

\section{$G=SO(n,1)$, $n$ pair}

On supposera $n\geq 4$. Les remarques du \pa 3 sur la connexit\'e restent exactes, ainsi que la description du parabolique minimal. Par ailleurs $G$ est un groupe de type $C_{\ell}$ o\`u $\ell=\tfrac{n}{2}$ et son dual (connexe) est

\begin{equation*} \index{groupe dual}
\widehat G=Sp(\ell)=Sp\Bigl(\tfrac{n}{2}\Bigr),
\end{equation*}
le groupe symplectique de la forme altern\'ee sur $\mathbb{C}^{2\ell}$ de matrice

\begin{equation*}
\widehat J=\begin{pmatrix}
&&&&&-1\\
&&&&\adots\\
&&&-1\\
&&1\\
&\adots\\1 
\end{pmatrix}\ .
\end{equation*}

On a dans $\widehat G$ un tore maximal $\widehat T$ form\'e des matrices

\begin{equation}
x=\begin{pmatrix}
x_{1}\\&\ddots\\&&x_{\ell}\\&&&x_{\ell}^{-1}\\ &&&&\ddots\\ &&&&&x_{1}^{-1}
\end{pmatrix}
\end{equation}
et un sous-groupe de Borel $\widehat B$ form\'e des matrices triangulaires sup\'erieures de $\widehat G$. Enfin

\begin{equation*}
{}^{L}G=\widehat G\times \Gal(\mathbb{C}/\mathbb{R}).
\end{equation*}

Une repr\'esentation admissible irr\'eductible $\pi$ de $G$ est de la forme suivante~:

\begin{align}
&\pi\ \textrm{appartient \`a la s\'erie discr\`ete}.\\
&\pi\  \textrm{est un sous-module de}\ I(\tau,s),\ \tau\in\widehat {{}^{0}M},\ s\in i\mathbb{R}.\\
\begin{split}
&\pi=J(\tau,s)\ \textrm{est l'unique quotient irr\'eductible de}\ I(\tau,s)\\
&\textrm{o\`u}\ \tau\in\widehat {{}^{0}M}\ \textrm{et}\ Re(s)>0\ .
\end{split}
\end{align}

Comme dans (6.3.10) une sous-alg\`ebre de Cartan $\mathfrak{h}_{0}$ de $\mathfrak{g}_{0}$ est donn\'ee par les matrices

\begin{equation*}
X=\begin{pmatrix}
\begin{array}{ ccccc  | cc}
&-x_{1}&&&&&\\
x_{1}&&&&& \\
&&\ddots&&&\\
&&&&-x_{ \ell-1}&\\
&&&x_{\ell-1}&&\\
&&&&0\\
 \hline
&&&&&x_{\ell}\\
&&&&&&-x_{ \ell}\\
\end{array}
\end{pmatrix}(x_{i}\in\mathbb{R})\ .
\end{equation*}
On a de nouveau $\mathfrak{h}\otimes\mathbb{C}\cong \mathbb{C}^{\ell}$, le groupe de Weyl \'etant maintenant $\mathfrak{S}_{\ell}\ltimes \{\pm 1\}^{\ell}=W $.

Soit $\psi:W_{\mathbb{R}}\times \SL(2)\rightarrow{}^{L}G$ un param\`etre d'Arthur, et
$P\in\mathbb{C}^{2\ell}$ le param\`etre holomorphe associ\'e (cf. \pa 6.3). Alors
$P=(P',-P')$ o\`u $P'\in\mathbb{C}^{\ell}$ est bien d\'efini modulo l'action de $W$. Donc $P'$ d\'efinit naturellement un \'el\'ement de $\mathfrak{h}^{\ast}$, et on peut v\'erifier :

\begin{lem}
$P'$ est le caract\`ere infinit\'esimal de la repr\'esentation associ\'ee \`a~$\varphi_{\psi}$. 
\end{lem} 

La d\'emonstration est laiss\'ee au lecteur (cf. Lemme~6.2.1).

Par ailleurs ${}^{0}\mathfrak{m}\cong \mathfrak{so}(n-1)$, l'alg\`ebre de Lie ${}^{0}\mathfrak{m}$ s'identifie \`a
$\mathbb{C}^{\ell-1}\subset \mathbb{C}^{\ell}$ et $\lambda_{\tau}$ v\'erifie

\begin{equation*}
\kern-6,6cm\lambda_{\tau}=\mu_{\tau}-\rho
\end{equation*}
\begin{align}
&\rho\ =(\ell-\tfrac{3}{2},\ell-\tfrac{5}{2},\ldots , \tfrac{1}{2})\\ 
&\mu_{\tau}=(\mu_{1},\ldots , \mu_{\ell-1})\ ,\  \mu_{i}\in\mathbb{Z}\ ,\ \mu_{1}\geq\cdots\geq\mu_{\ell-1}\geq 0\ . 
\end{align}

On a donc
\begin{equation}
\lambda_{\tau}=(\lambda_{i})\ ,\ \lambda_{i}\in\tfrac{1}{2}+\mathbb{Z}\ ,\ \lambda_{1}>\cdots> \lambda_{\ell-1}>0\ .
\end{equation}

Le caract\`ere infinit\'esimal de $I(\tau,s)$ est alors
\begin{equation}
\lambda_{I}=(\lambda_{\tau},s).
\end{equation}

\begin{lem}
Si le caract\`ere infinit\'esimal d'une repr\'esentation unitaire $\pi=J(\tau,s)$ est associ\'e \`a un param\`etre d'Arthur, $s\in i\mathbb{R}$ ou $s\in\tfrac{1}{2}+\mathbb{Z}$.
\end{lem}

Tout d'abord, le Lemme 6.3.3 s'applique et montre que $s\in i\mathbb{R}$ ou $s\in\tfrac{1}{2}\mathbb{Z}$. Supposons alors $s=p\in\mathbb{Z}$. Comme apr\`es le Lemme 6.3.3, $\psi:\mathbb{C}^{\times }\times \SL(2)\rightarrow\GL(2\ell,\mathbb{C})$ est de la forme
\begin{equation*}
\psi=\chi\oplus\chi^{-1}\oplus\sum_{i}\chi_{i}\otimes r_{i}
\end{equation*}
o\`u $\chi(z)=(z/\bar z)^{p}$. Les exposants holomorphes $P_{j}$ apparaissant dans
$\displaystyle\sum_{i}\chi_{i}\otimes r_{i}$ doivent \^etre dans $\tfrac{1}{2}+\mathbb{Z}$
d'apr\`es (6.4.8) et (6.4.7). On peut consid\'erer comme dans le \pa 6.3 la partie de $\psi$ d'exposants $P$ \textbf{entiers} :
\begin{equation*}
\psi_{e}(z)=(\chi(z),\chi^{-1}(z),1)
\end{equation*}
et
\begin{equation*}
\psi_{e}(\bar z)=x\ \psi_{e}(z)x^{-1}
\end{equation*}
o\`u $\psi(j)=(x,\sigma)\in\widehat G\times \Gal(\mathbb{C}/\mathbb{R})$.

Revenant au choix (6.4.1) de $\widehat T$ on peut supposer
\begin{equation*}
\psi_{e}(z)=\begin{pmatrix}
\chi\\ & 1 \\ && \ddots\\&&& 1\\&&&&\chi^{-1} 
\end{pmatrix}\ ,\ 
\psi_{e}(\bar z)=\begin{pmatrix}
\chi^{-1}\\ &1 \\&& \ddots\\&&&1\\&&&&\chi 
\end{pmatrix}.
\end{equation*}
Alors $x$ pr\'eserve la d\'ecomposition correspondante $\mathbb{C}^{2\ell}=\mathbb{C}^{2}\oplus\mathbb{C}^{2\ell-2}$ et s'\'ecrit $x=(x',x'')$, $x'\in SL(2,\mathbb{C})$, $x''\in Sp(\ell-1,\mathbb{C})$. La relation $x\psi_{e}(z)x^{-1}=\psi_{e}(\bar z)$ implique $x'=\begin{pmatrix}&b\\a\end{pmatrix}$, $-ab=1$, d'o\`u $x'{}^{2}=-1$ contrairement \`a $\chi(j^{2})=\chi(-1)=1$\ .

\medskip

Comme pr\'ec\'edemment, on en d\'eduit :
\begin{prop} \index{caract\`ere infinit\'esimal} \index{param\`etre d'Arthur}
Toute repr\'esentation unitaire $\pi$ de $G$ dont le caract\`ere infinit\'esimal est associ\'e \`a un param\`etre d'Arthur est de la s\'erie discr\`ete, ou \'egale \`a un sous-module de $I(\tau,s)$ ou \`a $J(\tau,s)$ avec
\begin{itemize}
\item[$(i)$] $s\in i\mathbb{R}$
 \item[$(ii)$] $s\in\tfrac{1}{2}+ \mathbb{Z}$\ .
\end{itemize} 
\end{prop}
  
Pour les repr\'esentations sph\'eriques, on a alors :

\begin{conj} \textrm{(Burger-Sarnak)} \index{Conjecture de Burger et Sarnak}
\begin{itemize}
\item[$(i)$] Si une repr\'esentation $J(1,s)$ $(s\in\mathbb{R}$, $s\neq 0)$ appara\^{\i}t dans $L^{2}(\Gamma\ba G),\pm s\in\{\tfrac{1}{2},\ldots\tfrac{n-1}{2}\} $.
\item[$(ii)$] En particulier les valeurs propres du laplacien associ\'ees v\'erifient 

 \quad\  $\lambda\geq\Bigl(\tfrac{n-1}{2}\Bigr)^{2}$

 ou $\lambda=\Bigl(\tfrac{n-1}{2}\Bigr)^{2}-\Bigl(\tfrac{n-1}{2}-k\Bigr)^{2}$ $k=0,1,\ldots\tfrac{n-2}{2} $.
\end{itemize} 
La premi\`ere valeur propre $\neq 0$ est $\lambda_{1}=n-2$ \footnote{ Pour $n>2$ ! c'est $\tfrac{1}{4}$ pour $n=2$.}.
\end{conj}

\markboth{CHAPITRE 6. CONJECTURES D'ARTHUR}{6.5. REPR\'ESENTATIONS EXCEPTIONNELLES}

\section{Existence des repr\'esentations exceptionnelles}
 
Comme nous l'avons rappel\'e les Conjectures 6.3.5 et 6.4.4 sont dues \`a Burger et Sarnak (cf. \cite{BurgerSarnak}). Il montrent de plus que si $G$ est un groupe
alg\'ebrique sur ${\Bbb Q}$ obtenu par restriction des scalaires \`a partir d'un groupe orthogonal sur un corps de nombres et tel que $G^{{\rm nc}} \cong 
SO(n,1)$, les valeurs propres
$$ \left( \frac{n-1}{2} \right)^2 - k^2 $$
pour $k=\frac{n-1}{2} , \ldots , \frac{n-1}{2} - \left[ \frac{n-1}{2} \right] $ apparaissent dans le spectre automorphe de $G$. 
Ce qui montre que les Conjectures 6.3.5 et 6.4.4 sont optimales.

Le Corollaire 6.2.7 et les Conjectures d'Arthur impliquent une conjecture beaucoup plus forte que la Conjecture A pour tout groupe $G$ tel que $G^{{\rm nc}} \cong SU(n,1)$.
Nous pensons que cette conjecture est optimale dans le m\^eme sens que la conjecture de Burger et Sarnak est optimale. La grande diff\'erence est que nous consid\'erons
cette fois tout le spectre automorphe et non seulement le spectre sph\'erique. Bien que nous ne sachions pas montrer que cette conjecture est optimale, nous 
pouvons montrer le r\'esultat suivant.

\begin{thm} 
Soit $G$ un groupe alg\'ebrique sur ${\Bbb Q}$ obtenu par restriction des scalaires \`a partir d'un groupe unitaire sur un corps de nombres et  tel que 
$G^{{\rm nc}} \cong SU(n,1)$. 
Soient $a,p,q,k$ des  entiers v\'erifiant 
\begin{itemize}
\item $0 \leq a \leq p, q \leq n$;
\item $p-q \in \{ -1, 0 ,1\}$;
\item $0 \leq k \leq \left[ \frac{n-2a}{2} \right]$.
\end{itemize}   
Alors, le dual automorphe $\widehat{G}_{{\rm Aut}}$ de $G$ contient une repr\'esentation unitaire irr\'eductible $\pi$ 
de $G$ telle que Hom$_K (\Lambda^p \mathfrak{p}^+ \otimes \Lambda^q \mathfrak{p}^- , \pi ) \neq 0$ et dont la valeur propre de l'op\'erateur de Casimir dans 
${\cal H}_{\pi}$ est \'egale \`a 
$$(n-2a)^2 - (n-2a-2k)^2 .$$
\end{thm}
{\it D\'emonstration.} Nous conservons les notations du \pa 2. On suppose que  $\sigma = \sigma_{a,b}$ et que $\chi$ est de la forme (6.2.8) avec 
$s = \frac{n-2a-2k}{2}$. La Proposition 4.5.1 implique que  
Hom$_K (\Lambda^p \mathfrak{p}^+ \otimes \Lambda^q \mathfrak{p}^- , J(\sigma , \chi) ) \neq 0$.

D'apr\`es \cite[Theorem 7.7]{Schlichtkrull} la repr\'esentation $J(\sigma , \chi )$ apparait discr\`etement comme sous repr\'esentation de $L^2 (H \backslash G)$
pour $H=S(U(a+k) \times U(n-a-k , 1))$. 

Puisque $G$ provient (par restriction des scalaires) d'un groupe unitaire sur un corps de nombre, il contient un ${\Bbb Q}$-sous-groupe alg\'ebrique $H$ tel que 
$H^{{\rm nc}} \cong S(U(a+k) \times U(n-a-k,1))$. Un th\'eor\`eme de Burger, Li et Sarnak \cite{BurgerLiSarnak} affirme que si $\rho \in \widehat{H}_{{\rm Aut}}$ et
si $\pi$ est une repr\'esentation de $G$ faiblement contenue dans l'induite (unitaire) de $\rho$ de $H$ \`a $G$, alors $\pi \in \widehat{G}_{{\rm Aut}}$. La repr\'esentation
triviale $1_H$ de $H$ appartient bien \'evidemment \`a $\widehat{H}_{{\rm Aut}}$ et l'induite de $H$ \`a $G$ de $1_H$ est $L^2 (H \backslash G)$. La repr\'esentation 
$J(\sigma , \chi )$ appartient donc \`a $\widehat{G}_{{\rm Aut}}$.

\medskip 

Remarquons que d'apr\`es les Propositions 6.3.4 et 6.4.3, les Conjectures d'Arthur impliquent la Conjecture A pour tout groupe $G$ v\'erifiant $G^{{\rm nc}} \cong SO(n,1)$.  
On peut de m\^eme facilement  d\'eduire des Propositions 6.3.4 et 6.4.3 et des Conjectures d'Arthur, une g\'en\'eralisation de la Conjecture A. Cette conjecture 
n'est certainement pas optimal. Un (fastidieux) exercice d'alg\`ebre lin\'eaire permettrait surement d'obtenir un r\'esultat optimal.

\newpage

\thispagestyle{empty}

\newpage

\markboth{CHAPITRE 7. TH\'EOR\`EME DE LUO-RUDNICK-SARNAK}{7.1. FONCTION DE $L$ DE RANKIN-SELBERG}

\chapter{Th\'eor\`eme de Luo-Rudnick-Sarnak}

Dans ce chapitre nous d\'emontrons la l\'eg\`ere g\'en\'eralisation suivante d'un th\'eor\`eme de Luo, Rudnick et 
Sarnak \cite{LRS}.

\begin{thm}[Approximation de Ramanujan] \label{LRS}
Soit $F$ un corps de nombre et ${\Bbb A}$ son anneau des ad\`eles. 
Soit $n$ un entier $\geq 2$ et $\pi = \otimes_v \pi_v$ une repr\'esentation automorphe cuspidale 
de $GL(n,{\Bbb A})$. Alors, en chaque place archim\'edienne $v$ la repr\'esentation $\pi_v$ est
un quotient de Langlands 
$$\pi_v = J(\sigma_{1,v} , \ldots , \sigma_{r,v} ),$$
o\`u les $\sigma_{j,v}$ sont comme dans les Th\'eor\`emes \ref{class reel} et \ref{class complexe} et tels que pour tout $j=1 , \ldots , r$, 
$$|\mbox{Re}(t_{j ,v})| \leq \frac12 - \frac{1}{n^2 +1} .$$
\end{thm}

\medskip

Remarquons que depuis l'annonce de nos r\'esultats dans \cite{BergeronClozel} W. M\"uller nous a inform\'e qu'il a, conjointement
avec B. Speh, lui aussi \'etendu le r\'esultat de Luo, Rudnick et Sarnak aux repr\'esentations non ramifi\'ees. 

Fixons un corps de nombres $F$, ${\Bbb A}$ son anneau des ad\`eles et $n$ un entier $\geq 2$.

\section{Fonction $L$ de Rankin-Selberg}

Dans cette section nous faisons quelques rappels concernant la th\'eorie de Rankin-Selberg. 
De la m\^eme mani\`ere que l'on a associ\'e \`a une repr\'esentation automorphe une fonction $L$,
la th\'eorie de Rankin-Selberg associe une fonction $L$ (dite de Rankin-Selberg) \`a toute paire de 
repr\'esentations irr\'eductibles unitaires cuspidales $\pi = \otimes_v \pi_v$ et $\pi ' = \otimes_v \pi_v '$ de 
$GL_n ({\Bbb A} )$ et $GL_{n'} ({\Bbb A})$ respectivement. Cette fonction $L$ est encore donn\'ee par un produit eul\'erien
$$L(s, \pi \times \pi ') = \prod_v L(s, \pi_v \times \pi_v ' ) $$
o\`u $v$ d\'ecrit l'ensemble des places de $F$.

Nous allons commencer par d\'ecrire les facteurs locaux $L(s, \pi_v \times \pi_v ')$ lorsque $v$ est une place 
finie, puis lorsque $v$ est une place archim\'edienne; enfin nous nous int\'eresserons \`a la fonction $L$ globale.

\medskip

Fixons $v$ une place de $F$. Dans les trois prochaines sous-sections nous nous int\'eressons aux facteurs locaux, 
pour simplifier les notations on notera $F=F_v$ et $\pi = \pi_v$. Pour l'instant $F$ est donc un corps local
de caract\'eristique z\'ero. 

Si $\pi$ est une repr\'esentation admissible irr\'eductible de $GL_n (F)$, nous notons $\tilde{\pi}$ {\it la 
repr\'esentation contragr\'ediente} \index{repr\'esentation contragr\'ediente} que l'on peut r\'ealiser dans le m\^eme espace que $\pi$ par l'action 
$\tilde{\pi} (g) = \pi ( ^t g^{-1} )$.

Soient $\pi_i$, $i=1, \ldots , r$, des repr\'esentations admissibles irr\'eductibles des groupes respectifs $GL_{n_i} (F)$.
Alors $\pi = \pi_1 \otimes \ldots \otimes \pi_r$ est une repr\'esentation admissible irr\'eductible de 
$$M(F) =GL_{n_1} (F) \times \ldots \times GL_{n_r} (F) .$$
Pour $\underline{s} = (s_1 ,\ldots , s_r ) \in {\Bbb C}^r$ soit $\pi_i [s_i ]$ la repr\'esentation de $GL_{n_i} (F)$ d\'efinie
par
$$\pi_i [s_i ] (g) = |{\rm det} |^{s_i} \pi_i (g) , \; \; \; g\in GL_{n_i} (F) .$$ \index{$\pi[s]$}
Nous notons
$$I_P^G ( \pi , \underline{s} )= {\rm ind}_{P(F)}^{G(F)} ( \pi_1 [s_1 ] \otimes \ldots \otimes \pi_r [s_r] )$$ \index{$I_P^G (\pi , \underline{s})$}
l'induite unitaire.

Remarquons que Shalika a d\'emontr\'e que la composante 
locale d'une repr\'esentation automorphe cuspidale de $GL_n ({\Bbb A})$ est g\'en\'erique \cite{Shalika}.
Et Jacquet et Shalika \cite{JacquetShalika2} ont montr\'e que toute repr\'esentation irr\'eductible unitaire g\'en\'erique $\pi$ 
de $GL_n (F)$ est \'equivalente \`a une repr\'esentation induite (et non seulement \`a un quotient)
$$\pi = I_P^G (\sigma , \underline{s} ),$$
o\`u $P$ est un sous-groupe parabolique standard de type $(n_1 , \ldots , n_r)$ de $GL_n$, $\underline{s} \in
{\Bbb C}^r$ et $\sigma$ est une repr\'esentation de carr\'e int\'egrable de $M_P (F)$, le sous-groupe de Levi de $P$
sur le corps local $F$ \footnote{C'est la seule propri\'et\'e des repr\'esentations g\'en\'eriques que nous utiliserons, 
le lecteur non familier peut donc accepter comme ``bo\^{\i}te noire'' le fait que toute repr\'esentation $\pi$ comme ci-dessus
est obtenue comme induite, r\'esultat qui d\'ecoule de \cite{Shalika} et \cite{JacquetShalika2}}. 

Dans les trois sous-sections suivantes nous d\'ecrivons les facteurs locaux de la fonction $L$ de Rankin-Selberg lorsque 
$F$ est non-archim\'edien, ou bien \'egal \`a ${\Bbb R}$ ou \`a ${\Bbb C}$.

\subsection{$F$ non-archim\'edien}

Comme nous l'avons rappel\'e, toute repr\'esentation irr\'eductible unitaire g\'en\'erique $\pi$ de $GL_n (F)$ est isomorphe  
\`a une repr\'esentation induite 
\begin{eqnarray} \label{pi comme induite}
\pi \cong I_P^G ( \sigma , \underline{s} ),
\end{eqnarray}
o\`u $P$ est d\'efini par une partition $(n_1 , n_2 , \ldots , n_r ) $ avec $n_{i+1} \leq n_i$,
$\underline{s} = (s_1 , \ldots , s_r ) \in {\Bbb C}^r$ et $\sigma = \otimes_i \sigma_i$
o\`u chaque $\sigma_i$ est une repr\'esentation de la s\'erie discr\`ete de $GL_{n_i} (F)$. 
Soit $\pi$ (resp. $\pi '$) une repr\'esentation irr\'eductible unitaire g\'en\'erique de $GL_n(F)$ (resp. $GL_{n'} (F)$). On d\'efinit le 
facteur local de Rankin-Selberg pour la paire $(\pi , \pi ')$ par r\'ecurrence. Tout d'abord, 
$L(s , \pi \times \pi  ') =L( s, \pi ' \times \pi )$. Si $\pi$ est comme dans (\ref{pi comme induite}),
\begin{eqnarray} \label{facteur L non archimedien} \index{facteur $L$ de Rankin-Selberg non archim\'edien}
L(s, \pi \times \pi ') = \prod_{i=1}^r L(s + s_i , \sigma_i \times \pi ' ). 
\end{eqnarray}

\medskip

Il reste \`a d\'efinir les facteurs $L$ pour des repr\'esentations de la s\'erie discr\`ete. 
De la m\^eme mani\`ere que les repr\'esentations de la s\'erie discr\`ete forment les blocs pour construire 
toutes les repr\'esentations unitaires g\'en\'eriques du groupe lin\'eaire, ce sont les repr\'esentations supercuspidales
qui forment les blocs \'el\'ementaires permettant de construire toutes les repr\'esentations de la s\'erie discr\`ete.
Rappelons donc qu'une repr\'esentation $\rho$ de $GL_d (F)$ est dite {\it supercuspidale} \index{supercuspidale} si elle n'appara\^{\i}t comme 
sous-quotient d'aucune repr\'esentation induite d'un sous-groupe parabolique propre. Les repr\'esentations 
supercuspidales sont donc exactement les repr\'esentations qui ne sont pas accessibles par le proc\'ed\'e 
d'induction. Les repr\'esentations supercuspidales peuvent aussi \^etre caract\'eris\'ees par le fait que leurs 
coefficients matriciels sont \`a support compact modulo le centre de $G$. On renvoie \`a \cite{BernsteinZelevinski}
pour plus de d\'etails sur les repr\'esentations supercuspidales.
 
Soit maintenant $\sigma$ une 
repr\'esentation de carr\'e int\'egrable de $GL_m (F)$.
D'apr\`es Bernstein et Zelevinski \cite{BernsteinZelevinski}, il existe un diviseur $d | m$, un sous-groupe parabolique 
standard $P$ de type $(d, \ldots ,d)$ de $GL_m (F)$, et une repr\'esentation $\rho$ de $GL_d (F)$ 
supercuspidale, irr\'eductible et unitaire tels que $\sigma$ soit l'unique quotient irr\'eductible de la 
repr\'esentation induite ${\rm ind}_P^G (\rho_1 \otimes \ldots \otimes \rho_r )$, o\`u $r = m/d$ et $\rho_j = \rho [j- (r+1)/2]$.
On d\'esigne cette repr\'esentation par $\sigma = \Delta (r, \rho )$. \index{$\Delta (r, \rho )$}

Soit donc $\sigma = \Delta (r , \rho)$ (resp. $\sigma ' = \Delta (r' , \rho ' )$) une repr\'esentation de la s\'erie discr\`ete
de $GL_m (F)$ (resp. $GL_{m'} (F)$), avec $m\leq m'$. On pose alors :
\begin{eqnarray} \label{facteur L non archimedien 2}
L(s, \sigma  \times \sigma ') = \prod_{j=1}^r L(s+\frac{r+r'}{2} -j , \rho \times \rho ' ) .
\end{eqnarray}

La description des facteurs $L$ de Rankin-Selberg est donc r\'eduite au cas de deux repr\'esentations supercuspidales.
Soient $\rho_1$ et $\rho_2$ deux repr\'esentations supercuspidales de respectivement $GL_{k_1} (F)$ et 
$GL_{k_2} (F)$. Lorsque $\rho_2$ n'est pas un twist de $\rho_1$, {\it i.e.} s'il n'existe pas de nombre complexe 
$t\in {\Bbb C}$ tel que $\rho_2 \cong \rho_1 [t]$, en particulier lorsque $k_1 \neq k_2$, on pose 
$L(s, \rho_1 \times \rho_2) =1$. On pose enfin  
\begin{eqnarray} \label{facteur L non archimedien 3}
\begin{array}{ccl}
L(s, \rho_1 \times \tilde{\rho}_1 [t] ) & = & L(s+t , \rho_1 \times \tilde{\rho}_1 ) \\
                                                         & = & (1-q^{-a(t+s)} )^{-1} ,
\end{array} 
\end{eqnarray}
o\`u $a | k_1$ est l'ordre du groupe cyclique des caract\`eres nonramifi\'es $\chi = |{\rm det}|^u$ tels que  
$\rho_1 \otimes \chi \cong \rho_1$. Si $\sigma_i = \Delta (r_i , \rho_i )$, (\ref{facteur L non archimedien 2}) et 
(\ref{facteur L non archimedien 3}) impliquent
\begin{eqnarray*}
L(s+2 {\rm Re}(s_i ), \sigma_i \times \tilde{\sigma}_i ) & = & \prod_{j=1}^{r_i} L(s+ 2 {\rm Re} (s_i )+r_i -j, \rho \times \tilde{\rho} ) \\
& = & \prod_{j=1}^{r_i} (1 -q^{-a_i (s+ 2 {\rm Re} (s_i) +r_i -j) } )^{-1} , 
\end{eqnarray*}
o\`u $a_i$ est l'ordre du groupe cyclique des caract\`eres nonramifi\'es $\chi$ tels que $\rho_i \otimes \chi$ est isomorphe
\`a $\rho_i$. On en d\'eduit que si $\pi = I_P^G (\sigma , \underline{s} )$ est une repr\'esentation irr\'eductible unitaire
de la s\'erie principale de $GL_n (F)$ et si $i$ est un indice tel que ${\rm Re}(s_i )>0$, il existe un facteur de la fonction 
$L(s, \pi \times \tilde{\pi} )$ ayant un p\^ole en $s=-2{\rm Re}(s_i )$.

\subsection{$F = {\Bbb R}$}

Soit $\pi$ (resp. $\pi '$) une repr\'esentation irr\'eductible unitaire de $GL_n ({\Bbb R} )$ (resp. $GL_{n'} ({\Bbb R})$).
On d\'efinit les facteurs locaux r\'eels de la fonction $L$ de Rankin-Selberg attach\'ee \`a la paire de repr\'esentation
$(\pi , \pi ')$ \`a partir des repr\'esentations semisimples $\varphi$ et $\varphi '$ du groupe de Weil $W_{{\Bbb R}}$ 
de degr\'es respectifs $n$ et $n'$ associ\'ees par la correspondance de Langlands locale aux repr\'esentations 
$\pi$ et $\pi '$.

Nous avons d\'eja d\'ecrit les facteurs $L$ associ\'es \`a une repr\'esentation semisimple du groupe de Weil $W_{{\Bbb R}}$.
Remarquons que le produit tensoriel $\varphi \otimes \varphi '$ d\'efinit une repr\'esentation semisimple du groupe 
de Weil $W_{{\Bbb R}}$ de degr\'e $nn'$. On d\'efinit alors :
\begin{eqnarray} \label{facteur L reel}
L(s, \pi \times \pi ' ) & = & L(s , \varphi \otimes \varphi ' ). \index{facteur $L$ de Rankin-Selberg r\'eel}
\end{eqnarray}

Supposons que la repr\'esentation
$\pi$ est isomorphe au quotient de Langlands $J(\sigma_1 , \ldots , \sigma_r )$ o\`u les $\sigma_j$ sont comme dans 
le Th\'eor\`eme \ref{class reel}. D'apr\`es le Th\'eor\`eme \ref{Unit} on d\'eduit de l'unitarit\'e de $\pi$ que pour chaque entier $j \in [1,r]$
soit
\begin{itemize}
\item $t_j$ est imaginaire pure, soit
\item Re$(t_j )>0$ et il existe $k \neq j$ tel que $n_k = n_j$, $l_k = l_j$ (ceux-ci pouvant \^etre \'egaux \`a un entier comme
\`a l'un des signes $+$ ou $-$) et $t_j = - \overline{t}_k$.
\end{itemize}

En utilisant la description des facteurs $L$ donn\'ee par la Proposition \ref{facteurLreel}, on d\'eduit de (\ref{facteur L reel}) que $L(s, \pi \times 
\tilde{\pi} )$ est un produit de facteurs Gamma de la forme
$$4 (2 \pi )^{-2s -2t_i +2t_j -\max ( l_i ,l_j )} \Gamma (s+t_i -t_j + \frac{|l_i -l_j |}{2} )  \Gamma (s+t_i -t_j + \frac{l_i +l_j}{2} )$$
ou
$$4 (2\pi )^{-2s - l_i} \Gamma (s+t_i - t_j + \frac{l_i}{2} ) \Gamma (s+t_j -t_i + \frac{l_i}{2} ) $$
ou
$$\pi^{-\frac{s+t_i - t_j + \varepsilon_{i,j} }{2}} \Gamma (\frac{s+t_i - t_j + \varepsilon_{i,j} }{2} ),$$
o\`u $\varepsilon_{i,j} \in \{ 0,1\}$ avec $\varepsilon_{i,j} = \varepsilon_i + \varepsilon_j $ mod $2$, $1 \leq i,j \leq r$ et les $\varepsilon_i$
correspondent aux signes dans la Proposition \ref{facteurLreel}.
Puisque $\pi$ est unitaire on d\'eduit de ces formules que pour tout indice $i \leq r$ tel que Re$(t_i ) \neq 0$, le 
$L$-facteur local $L(s, \pi \times \tilde{\pi} )$ contient un facteur Gamma du type 
$$\Gamma (s+2 {\rm Re}(t_i) )$$
ou 
$$\Gamma (\frac{s+2 {\rm Re}(t_i)}{2} ).$$
Puisque la fonction Gamma $\Gamma (s)$ n'a pas de z\'ero, on conclut que $L(s, \pi \times \pi ')$ a son premier
p\^ole en $s = 2 \max |$Re$(t_i )|$. Rappelons enfin que le Th\'eor\`eme \ref{Unit} implique que Re$(t_i ) < \frac12$ puisque 
$\pi$ est g\'en\'erique.

\subsection{$F={\Bbb C}$}

Le cas $F={\Bbb C}$ se traite de la m\^eme mani\`ere que le cas $F={\Bbb R}$. 
Si $\pi$ est une repr\'esentation irr\'eductible unitaire de $GL_n ({\Bbb C})$ isomorphe au quotient de 
Langlands $J(\sigma_1 , \ldots , \sigma_n )$, o\`u les $\sigma_j$ sont comme dans le Th\'eor\`eme \ref{classcomplex}.
Comme dans le cas r\'eel, le fait que $\pi$ est unitaire et g\'en\'erique se traduit par le fait que la partie r\'eelle de chaque 
$t_i$ est strictement inf\'erieure \`a $\frac{1}{2}$ et que pour chaque $i$ tel que Re$(t_i ) \neq 0$ il existe $k\neq i$ tel que 
$p_k - q_k =p_i -q_i$ et $t_i = - \overline{t}_k$.

Le groupe de Weil est cette fois \'egal \`a ${\Bbb C}^*$ et les facteurs locaux de la fonction $L$ de Rankin-Selberg
se d\'efinissent comme dans le cas r\'eel. On obtient en particulier que 
$$L(s, \pi \times \tilde{\pi} ) = \prod_{i,j =1}^n 2 (2\pi )^{-\frac{s+t_i -t_j + \frac{|p_i-q_i -p_j -q_j |}{2}}{2}} \Gamma 
(s+t_i -t_j + \frac{|p_i-q_i -p_j -q_j |}{2} ).$$ \index{facteur $L$ de Rankin-Selberg complexe}
Comme dans le cas r\'eel, on en d\'eduit que si $\pi$ est unitaire, le facteur $L$ $L(s, \pi \times \tilde{\pi} )$ 
contient pour tout indice $i$ tel que Re$(t_i ) \neq 0$ la fonction 
$$\Gamma (s+2 {\rm Re} (t_i ) )$$
comme facteur. Puisque la fonction Gamma $\Gamma (s)$ n'a pas de z\'ero, on conclut que 
$L(s, \pi \times \pi ')$ a son premier p\^ole en $s = 2 \max |$Re$(t_i )|$.

\subsection{La fonction $L$ globale}

Les sous-sections pr\'ec\'edentes permettent de d\'efinir au moins formellement la fonction $L$ de Rankin-Selberg globale d'une
paire de repr\'esentations unitaires cuspidales automorphes $\pi = \otimes_v \pi_v$ et $\pi ' =\otimes_v \pi_v '$ de 
$GL_n ({\Bbb A} )$ et $GL_{n' } ({\Bbb A})$, o\`u ${\Bbb A}$ est l'anneau des ad\`eles d'un corps de nombre 
$F$. La fonction $L$ de Rankin-Selberg est d\'efinie par le produit eul\'erien :
\begin{eqnarray} \label{def fonction L globale}
L(s, \pi \times \pi ' ) = \prod_v L(s, \pi_v \times \pi_v '),
\end{eqnarray}
o\`u $v$ d\'ecrit l'ensemble des places de $F$. 

Comme pour la fonction z\'eta de Riemann l'int\'er\^et de cette fonction dans l'\'etude des repr\'esentations 
automorphes provient de ce qu'elle d\'efinit une vraie fonction analytique lorsque $s$ est dans 
un demi-plan et que cette fonction se prolonge m\'eromorphiquement au plan complexe tout entier et 
v\'erifie une \'equation fonctionnelle. 
Il est devenu r\'ecemment possible 
de d\'emontrer toutes ces propri\'et\'es, ainsi que d'autres, de la fonction $L$ de 
Rankin-Selberg en adoptant le point de vue naturel des repr\'esentations int\'egrales d\'evelopp\'e
par Jacquet, Shalika, Piatetski-Shapiro (\cite{JacquetPiatetskiiShapiroShalika}, \cite{JacquetShalika}) et, plus r\'ecemment, Cogdell et Piatetski-Shapiro \cite{Cogdell}.
C'est le point de vue que nous adoptons ici. Soulignons n\'eanmoins que les r\'esultats que nous utiliserons
pouvaient \^etre obtenus en combinant les m\'ethodes de Jacquet, Shalika et Piatetski-Shapiro \cite{JacquetPiatetskiiShapiroShalika}, Shahidi 
\cite{Shahidi}, \cite{Shahidi2}, \cite{Shahidi3} et Moeglin et Waldspurger \cite{MoeglinWaldspurger}. Cette derni\`ere approche est celle utilis\'ee par M\"uller et Speh.

On r\'esume les propri\'et\'es de la fonction $L$ de Rankin-Selberg dont nous aurons besoin dans 
le th\'eor\`eme suivant.

\begin{thm} \label{fonction L globale}
Chaque produit (\ref{def fonction L globale}) est convergent lorsque la partie r\'eelle de $s$ est suffisamment grande.
Ces produits d\'efinissent des fonctions holomorphes de $s$ qui se prolongent m\'eromorphiquement au plan complexe tout 
entier. Les fonctions ainsi obtenues, et toujours not\'ees $L(s, \pi \times \pi ')$, sont sympathiques dans le sens que
\begin{enumerate}
\item chaque fonction $L(s, \pi \times \pi ')$ reste born\'ee dans les bandes verticales (loin de ses p\^oles),
\item elles satisfont \`a l'\'equation fonctionnelle
$$L(s , \pi \times \pi ') = \varepsilon (s, \pi \times \pi ') L(1-s, \tilde{\pi} \times \tilde{\pi}'),$$
o\`u le facteur $\varepsilon$ est une fonction enti\`ere de $s$ de la forme 
$$\varepsilon (s, \pi \times \pi ' ) = W(\pi \times \pi ') (D_{F/{\Bbb Q}}^{nn'} N(\pi \times \pi '))^{\frac{1}{2} -s},$$
o\`u $D_{F/{\Bbb Q}}$ est le discriminant du corps de nombre $F$, $W(\pi \times \pi ')$ est un nombre 
complexe de module  $1$, et $N(\pi \times \pi ' )$ un nombre entier,
\item si $n' <n$, la fonction $L(s, \pi \times \pi ')$ est enti\`ere,
\item si $n'=n$, la fonction $L(s, \pi \times \pi ')$ a au plus des p\^oles simples qui interviennent si et seulement si
$\pi \cong \tilde{\pi} ' [i \sigma]$ avec $\sigma$ r\'eel, et qui dans ce cas sont $s=-i\sigma$ et $s=1-i\sigma$.
\end{enumerate}
\end{thm} \index{fonction $L$ de Rankin-Selberg}

Il n'est pas question de donner une preuve compl\`ete de ce r\'esultat ici, une r\'ecente pr\'epublication de Cogdell et 
Piatetski-Shapiro \cite{Cogdell} contient la preuve compl\`ete de ce Th\'eor\`eme. Afin d'en faciliter la reconstitution, nous 
donnons les grandes id\'ees de la preuve lorsque $n=n'$, le cas qui nous int\'eressera par la suite.

Soient donc $(\pi , V_{\pi } )$ et $(\pi ' , V_{\pi '})$ deux repr\'esentations cuspidales unitaires de $GL_n ({\Bbb A})$.
Soient $\omega = \otimes \omega_v$ et $\omega '= \otimes \omega_v '$ leurs caract\`eres centraux.
Fixons $\psi = \otimes \psi_v $ un caract\`ere additif continu non trivial de ${\Bbb A}$ qui soit 
trivial sur $F$ et $\varphi \in V_{\pi}$ et $\varphi ' \in V_{\pi '}$ deux formes cuspidales.

\bigskip

Commen\c{c}ons par d\'emontrer la convergence de la fonction $L(s, \pi \times \pi ' )$ d\'efinie comme 
produit eul\'erien.

Soit $S$ un ensemble fini de places de $F$, contenant toutes les places \`a l'infini et tel que 
pour toute place $v \notin S$ on ait $\pi_v$ et $\pi_v '$ non ramifi\'ees, {\it i.e.} ayant un vecteur fixe 
sous l'action du compact maximal.

Pour les places $v \notin S$, pour lesquelles $\pi_v$ et $\pi_v '$ sont donc nonramifi\'ees, on 
v\'erifie que 
$$L(s, \pi_v \times \pi_v ') = \mbox{det} (I- q_v^{-s} A_{\pi_v} \otimes A_{\pi_v '} )^{-1}$$
o\`u $A_{\pi_v }$ et $A_{\pi_v '}$ sont les param\`etres de Satake associ\'es \`a $\pi_v$ et $\pi_v '$ et 
dont les valeurs propres sont toutes de valeur absolue inf\'erieure \`a $q_v^{\frac{1}{2}}$. Donc, cf. \cite{JacquetShalika4}
pour plus de d\'etails, la fonction $L$ partielle
$$L^S (s, \pi \times \pi ') = \prod_{v\notin S} L(s, \pi_v \times \pi_v ' )= \prod_{v\notin S} \mbox{det}(I-q_v^{-s} A_{\pi_v} \otimes A_{\pi_v '} )^{-1}$$
est absolument convergente pour Re$(s)>>0$. On en d\'eduit la m\^eme chose pour la fonction $L$ globale.

\bigskip

Il nous faut maintenant d\'emontrer les propri\'et\'es analytiques de cette fonction.
Comme dans le cas classique de la fonction z\'eta de Riemann, il va s'agir de repr\'esenter 
la fonction $L$ de Rankin-Selberg par une int\'egrale. C'est l\`a qu'inerviennent des s\'eries d'Eisenstein. \index{s\'eries d'Eisenstein}
Commen\c{c}ons donc par d\'ecrire ces s\'eries d'Eisenstein.

Pour construire ces s\'eries d'Eisenstein on peut suivre \cite{JacquetShalika4}. 
Dans $GL_n$ soit $P_n$ le sous-groupe stabilisant le vecteur $(0,\ldots , 0,1)$. On a $P_n \backslash GL_n 
\cong F^n -\{ 0 \}$. Si on d\'esigne par ${\cal S} ({\Bbb A})$ l'espace des fonctions de Schwartz-Bruhat sur 
${\Bbb A}^n$, alors chaque $\Phi \in {\cal S} ({\Bbb A})$ d\'efinie une fonction lisse sur $GL_n ({\Bbb A})$, invariante
\`a gauche sous $P_n ({\Bbb A})$, par $g \mapsto \Phi ((0, \ldots , 0,1)g) = \Phi (e_n g)$. 
Consid\'erons la fonction 
$$F(g,\Phi ; s) = | \mbox{det} (g) |^s \int_{{\Bbb A}^*} \Phi (ae_n g) |a|^{ns} \omega (a) \omega '(a) d^* a .$$
Si $P_n ' = Z_n P_n$ est le sous-groupe parabolique de $GL_n$ associ\'e \`a la partition $(n-1,1)$ alors 
on prouve \cite{JacquetShalika4} que l'on peut former les s\'eries d'Eisenstein
$$E(g, \Phi ; s) = \sum_{ \gamma \in P_n '(F) \backslash GL_n (F)} F(\gamma g, \Phi ; s).$$
Si l'on remplace dans cette somme la fonction $F$ par sa d\'efinition et que l'on permute le signe somme 
et l'int\'egrale, on obtient que ces s\'eries convergent absolument dans le demi-plan Re$(s) > 1$ \cite{JacquetShalika4}
et que
$$E(g, \Phi ;s) = |\mbox{det}(g)|^s \int_{F^* \backslash {\Bbb A}^*} \Theta_{\Phi} '(a,g) |a|^{ns} \omega (a) \omega '(a) d^* a $$
qui est essentiellement l'expression de la s\'erie d'Eisenstein comme transform\'ee de Mellin de la s\'erie Theta
$$\Theta_{\Phi} (a,g) = \sum_{\xi \in F^n } \Phi (a\xi g),$$
o\`u au-dessus on a not\'e $\Theta_{\Phi} ' (a,g) = \Theta_{\Phi} (a,g) -\Phi (0)$. En appliquant la transform\'ee de
Poisson \`a $\Theta_{\Phi}$, on obtient \cite[Section 4]{JacquetShalika4}~:

\begin{prop}
La s\'erie d'Eisenstein $E(g, \Phi;s)$ admet un prolongement m\'eromorphe \`a tout le plan ${\Bbb C}$ avec au plus 
des p\^oles simples en $s=-i\sigma, 1-i\sigma$ quand $\omega \omega '$ d\'efinit un caract\`ere non ramifi\'e de la 
forme $\omega (a) \omega ' (a) = |a|^{in\sigma}$. Comme fonction de $g$ elle est lisse et \`a croissance lente  et 
comme fonction de $s$ elle reste born\'ee dans les bandes verticales (loin des p\^oles possibles), uniform\'ement 
pour $g$ dans un compact. Elle v\'erifie de plus une \'equation fonctionnelle.
\end{prop}

Pour chaque paire de formes cuspidales $\varphi \in V_{\pi}$ et $\varphi ' \in V_{\pi '}$ on peut 
consid\'erer l'int\'egrale
$$I(s; \varphi , \varphi ' , \Phi ) = \int_{Z_n ({\Bbb A}) GL_n (F) \backslash GL_n ({\Bbb A})} \varphi (g) \varphi ' (g) E(g,\Phi; s) dg.$$
Celle-ci est bien d\'efinie puisque les formes cuspidales sont rapidement d\'ecroissantes. Elle d\'efinit donc une fonction
m\'eromorphe de $s$, born\'ee dans les bandes verticales loin des p\^oles et elle v\'erifie l'\'equation fonctionnelle
$$I(s;\varphi, \varphi ' ,\Phi ) = I(1-s; \tilde{\varphi }, \tilde{\varphi} ' , \hat{\Phi} ),$$
provenant de l'\'equation fonctionnelle de la s\'erie d'Eisenstein, o\`u $\tilde{\varphi} (g) = \varphi ((g^t )^{-1})$ et de 
m\^eme pour $\tilde{\varphi} '$.

Ces int\'egrales sont enti\`eres sauf si $\omega (a) \omega ' (a) = |a|^{in\sigma}$ est non ramifi\'e. Dans ce cas, les 
r\'esidus de cette int\'egrale en $s=-i\sigma$ et $s=1-i\sigma$ d\'efinissent des couplages invariants entre 
$\pi$ et $\pi ' [-i\sigma ]$ ou de mani\`ere \'equivalente entre $\tilde{\pi}$ et $\tilde{\pi} ' [i\sigma]$. Ainsi un 
r\'esidu ne peut \^etre non nul que si $\pi \cong \tilde{\pi}[i\sigma]$ et r\'eciproquement, on peut dans ce cas trouver 
$\varphi$, $\varphi '$ et $\Phi$ tel que ce r\'esidu ne soit pas nul.

Pour relier nos int\'egrales aux produits eul\'eriens, on remplace $\varphi$ et $\varphi '$ par leur d\'ecomposition 
en s\'eries de Fourier qui est dans ce cas de la forme
$$\varphi (g) = \sum_{\gamma \in N_n (F) \backslash P_n (F)} W_{\varphi} (\gamma g) ,$$
o\`u $W_{\varphi}$ est la fonction de Whittaker sur $GL_n ({\Bbb A})$ associ\'ee \`a $\varphi$ d\'efinie par $\psi$~:
$$W_{\varphi} (g) = \int_{N_n (F) \backslash N_n ({\Bbb A}) } \varphi (ng) \psi (n) dn,$$
$\psi$ \'etant un caract\`ere non d\'eg\'en\'er\'e de $N_n (F) 
\backslash N_n ({\Bbb A})$ (cf. \cite{Cogdell}).
En d\'eveloppant l'int\'egrale, on obtient~:
\begin{eqnarray*}
I(s;\varphi , \varphi ', \Phi ) & = & \int_{N_n ({\Bbb A}) \backslash GL_n ({\Bbb A})} W_{\varphi} (g) W_{\varphi '} (g) \Phi (e_n g) |\mbox{det}(g)|^s dg \\
                                            & =: & \Psi (s; W_{\varphi} , W_{\varphi '}  , \Phi ).
\end{eqnarray*}
Cette derni\`ere expression converge pour Re$(s)>>0$ par un estim\'e de jauge de Jacquet, Shalika et Piatetski-Shapiro, cf. \cite{Cogdell}.

Pour arriver \`a une expression sous la forme d'un produit Eul\'erien, on suppose que $\varphi$, $\varphi '$ et $\Phi$
se d\'ecomposent en produit tensoriels via $\pi= \otimes \pi_v$, $\pi '  =\otimes \pi_v '$ et ${\cal S} ({\Bbb A}^n ) = 
\otimes {\cal S} (F_v^n ) $ de telle mani\`ere que l'on ait 
$W_{\varphi} (g) = \prod_v W_{\varphi_v} (g_v )$, $W_{\varphi'} (g) = \prod_v W_{\varphi_v '} (g_v )$ et 
$\Phi (g) = \prod_v \Phi_v (g_v )$. Alors,
$$\Psi (s; W_{\varphi} , W_{\varphi '}  , \Phi ) = \prod_v \Psi_v (s; W_{\varphi_v } , W_{\varphi_v '}  , \Phi_v ),$$
o\`u 
$$\Psi_v (s; W_{\varphi_v } , W_{\varphi_v '}  , \Phi_v ) = \int_{N_n (F_v ) \backslash GL_n (F_v )} W_{\varphi_v } (g_v ) W_{\varphi_v '} (g_v ) \Phi_v (e_n g_v ) |\mbox{det}(g_v )|^s dg_v ,$$
qui converge lorsque Re$(s)>>0$ d'apr\`es des estim\'ees de jauges locales de Jacquet, Shalika et Piatetski-Shapiro, cf. \cite{Cogdell}.

\bigskip

Nous allons maintenant pouvoir d\'emontrer que la fonction $L$ globale de Rankin-Selberg peut-\^etre prolong\'ee 
m\'eromorphiquement au plan complexe tout entier. Jacquet et Shalika \cite{JacquetShalika} 
montrent que pour un choix appropri\'e de $\varphi \in V_{\pi }$, $\varphi '\in V_{\pi '}$ et $\Phi \in {\cal S} ({\Bbb A}^n )$
on a 
$$\mbox{det} (I-q_v^{-s} A_{\pi_v} \otimes A_{\pi_v '} )^{-1} = \Psi_v (s; W_{\varphi_v} ,W_{\varphi_v '}, \Phi_v ) .$$
On obtient alors :
\begin{eqnarray*}
I(s;\varphi , \varphi ' , \Phi ) & = & \left( \prod_{v\in S} \Psi_v (s; W_{\varphi_v} , W_{\varphi_v '}, \Phi_v ) \right) L^S (s, \pi \times \pi ') \\
& = & \left( \prod_{v \in S} \frac{ \Psi_v (s; W_{\varphi_v} , W_{\varphi_v '}, \Phi_v ) }{L(s,\pi_v  \times \pi_v ')} \right) L(s, \pi \times \pi ') \\
& = & \left( \prod_{v \in S} e_v (s; W_{\varphi_v} , W_{\varphi_v '} , \Phi_v ) \right) L(s, \pi \times \pi ').
\end{eqnarray*}
Chaque fonction $e_v (s; W_{\varphi_v} , W_{\varphi_v '} , \Phi_v )$ est enti\`ere. Pour les places non-archim\'ediennes
cela d\'ecoule du Th\'eor\`eme 2.7 de \cite{JacquetPiatetskiiShapiroShalika} et pour les places archim\'ediennes cela d\'ecoule
du Th\'eor\`eme 1.2 (et de son Corollaire) de \cite{Cogdell}. On en d\'eduit que la fonction $L(s, \pi \times \pi ')$
admet un prolongement m\'eromorphe \`a tout le plan complexe.

\bigskip

Attaquons-nous maintenant \`a l'\'equation fonctionnelle. Partons de l'\'equation fonctionnelle de l'int\'egrale 
globale :
$$I(s; \varphi , \varphi ' , \Phi ) = I (1-s; \tilde{\varphi} , \tilde{\varphi} ', \hat{\Phi} ).$$
Comme plus haut, pour un choix appropri\'e de $\varphi$, $\varphi '$ et $\Phi$,  
$$I(s; \varphi , \varphi ' , \Phi ) = \left( \prod_{v\in S} e_v (s; W_{\varphi_v } , W_{\varphi_v '} , \Phi_v ) \right) L(s, \pi \times \pi ')$$
alors que d'un autre c\^ot\'e 
$$I(1-s; \tilde{\varphi} , \tilde{\varphi} ' , \hat{\Phi} ) = \left( \prod_{v\in S} e_v (1-s; W_{\tilde{\varphi}_v } , W_{\tilde{\varphi}_v '} , \hat{\Phi}_v ) \right) L(1-s, \tilde{\pi} \times \tilde{\pi} ').$$
Mais d'apr\`es le Th\'eor\`eme 2.7 de \cite{JacquetPiatetskiiShapiroShalika} (pour les places non-archim\'ediennes) et le 
Th\'eor\`eme 5.1 de \cite{JacquetShalika} (pour les places archim\'ediennes), pour chaque place $v \in S$ on a 
$$e_v (1-s; W_{\tilde{\varphi}_v } , W_{\tilde{\varphi}_v '} , \hat{\Phi}_v ) = \left(  \omega_{\pi_v '} (-1)\right)^{n-1} \varepsilon (s, \pi_v \times \pi_v ' , \psi_v ) e_v (s;W_{\varphi_v },W_{\varphi_v '} , \Phi_v ),$$
o\`u le facteur epsilon est d\'ecrit explicitement dans les r\'ef\'erences ci-dessus
et $\omega_{\pi_v '}$ est le caract\`ere central de la repr\'esentation $\pi_v '$.

En combinant tout cela, on obtient 
$$L(s,\pi \times \pi ') = \left( \prod_{v\in S} \left( \omega_{\pi_v '} (-1)\right)^{n-1} \varepsilon (s, \pi_v \times \pi_v ' , \psi_v ) \right) L(1-s, \tilde{\pi} \times \tilde{\pi} ').$$
Le produit intervenant ci-dessus est classiquement not\'e $\varepsilon (s, \pi \times \pi' )$ (il est bien s\^ur ind\'ependant 
de $\psi$). Sa valeur est ais\'ement calculable \`a l'aide des r\'ef\'erences cit\'ees plus haut.

\bigskip

Nous admettrons que la fonction $L$ reste born\'ee dans les bandes verticales, propri\'et\'e d\'emontr\'ee dans 
\cite{Cogdell}. Concluons cette section par l'\'etude des p\^oles de la fonction $L$. Ceux-ci sont \'evidemment 
reli\'es aux p\^oles de l'int\'egrale globale $I(s; \varphi , \varphi ', \Phi )$. On a montr\'e que pour un choix 
appropri\'e de $\varphi$, $\varphi '$ et $\Phi$, 
$$I(s; \varphi , \varphi ' , \Phi ) = \left( \prod_{v\in S} e_v (s; W_{\varphi_v } , W_{\varphi_v '} , \Phi_v ) \right) L(s, \pi \times \pi ').$$
D'un autre c\^ot\'e, pour tout $s_0 \in {\Bbb C}$ et pour tout $v$ il existe un choix local $W_v$, $W_v '$ et $\Phi_v$ tel 
que $e_v (s_0 ;W_v , W_v ' ,\Phi_v ) \neq 0$. Pour les places $v$ archim\'ediennes c'est d\'emontr\'e 
dans \cite[Th\'eor\`eme 1.2 et son Corollaire]{Cogdell}. Pour les places $v$ non-archim\'ediennes cela d\'ecoule
de \cite[Th\'eor\`eme 2.7]{JacquetPiatetskiiShapiroShalika} qui d\'ecrit chaque facteur $L$ local comme g\'en\'erateur de
l'id\'eal fractionnaire engendr\'e par les int\'egrales locales correspondantes. Ceci implique l'existence d'un ensemble
fini de fonctions $W_{v,i}$, $W_{v,i} '$ et $\Phi_{v,i}$ telles que 
$$L(s, \pi_v \times \pi_v ' ) = \sum_i \Psi (s;W_{v,i} , W_{v,i} ' ,\Phi_{v,i} )$$
autrement dit 
$$1= \sum_i e(s;W_{v,i} , W_{v,i} ' , \Phi_{v,i} ).$$
D'o\`u l'on conclut que pour tout $s_0 \in {\Bbb C}$ l'un des $e (s_0 ; W_{v,i} , W_{v,i} ' ,\Phi_{v,i} )$ doit \^etre non nul. 

On d\'eduit des cas locaux que pour chaque $s_0 \in {\Bbb C}$, il existe des fonctions globales $\varphi$, $\varphi '$ 
et $\Phi$ telles que 
$$\frac{I(s_0 ; \varphi , \varphi ' , \Phi )}{L(s_0 ,\pi \times \pi ' ) } \neq 0 .$$
Les p\^oles de la fonction globale $L(s, \pi \times \pi ' )$ sont donc exactement ceux de la famille d'int\'egrales globales
$\{ I(s; \varphi , \varphi ' , \Phi ) \}$, o\`u $\varphi$, $\varphi '$ et $\Phi$ varient. Et l'\'etude de ces fonctions permet de 
conclure la d\'emonstration du Th\'eor\`eme \ref{fonction L globale}.

\markboth{CHAPITRE 7. TH\'EOR\`EME DE LUO-RUDNICK-SARNAK}{7.2. D\'EMONSTRATION DU TH\'EOR\`EME \ref{LRS}}

\section{D\'emonstration du Th\'eor\`eme \ref{LRS}}

Soit $F$ un corps de nombres de degr\'e $d$ d'anneau des entiers ${\cal O} = {\cal O}_F$.
Nous noterons $U=U_F$ son groupe des unit\'es et $D=D_{F/{\Bbb Q}}$ son discriminant. 
Enfin soit ${\Bbb A}$ l'anneau des ad\`eles de $F$.
Pour $n$ un entier $\geq 2$, soit $\pi =\otimes \pi_v $ une repr\'esentation automorphe cuspidale de $GL(n,{\Bbb A} )$ 
que l'on normalise de fa\c{c}on \`a ce qu'elle ait un caract\`ere central unitaire. Rappelons que la d\'ecomposition
$\pi = \otimes \pi_v$ est donn\'e par un Th\'eor\`eme de Flath \cite{Flath}, c'est d'ailleurs en fait un produit tensoriel 
restreint dans le sens que pour presque toute place $v$ de $F$, la representation $\pi_v$ de $GL_n (F_v )$ est 
{\it non ramifi\'ee} \index{non ramifi\'ee, repr\'esentation} {\it i.e.} admet un vecteur $GL_n ({\cal O}_v )$-invariant.

La fonction $L$ standard, $L(s,\pi )$, associ\'ee \`a $\pi$ (cf. \S 3.3) est de la forme :
$$L(s,\pi ) = \prod_v L(s,\pi_v ) .$$
Chaque $L(s, \pi_v )$ est le facteur d'Euler local (cf. \S 3.3).
Puisque $\pi$ est cuspidale, en chaque place archim\'edienne $v$ de $F$, on a rappel\'e que  la repr\'esentation $\pi_v$ de 
$GL(n,F_v )$ est g\'en\'erique et s'obtient comme une induite pleine :
$$\pi_v = J( \sigma_{1,v} , \ldots , \sigma_{r,v} ),$$
o\`u chaque $\sigma_{j,v}$ est comme dans les Th\'eor\`emes \ref{class reel} et \ref{class complexe}.
Le facteur $L$ local $L(s, \pi_v )$ est calcul\'e dans les Propositions \ref{facteurLreel} et \ref{facteurLcomplex}.

Sous ces notations, la Conjecture de Ramanujan g\'en\'eralis\'ee pour $\pi$ et pour les places archim\'ediennes 
postule que :
$$\mbox{Re}(t_{j ,v} )=0 \mbox{  pour tout } j.$$
Jacquet et Shalika \cite{JacquetShalika4,JacquetShalika3}, en utilisant la g\'en\'ericit\'e de $\pi$, ont montr\'e que :
\begin{eqnarray}
|\mbox{Re} (t_{j ,v} ) | < \frac12 .
\end{eqnarray}
Remarquons encore une fois (avec \cite{BR}) que cela d\'ecoule de la classification de Vogan du dual unitaire de 
$GL(n)$ (cf. Th\'eor\`eme \ref{Unit}). 

La preuve du Th\'eor\`eme \ref{LRS} repose en grande partie sur la th\'eorie de Rankin-Selberg dont 
on a rappel\'e les principaux r\'esultats dans la section pr\'ec\'edente. Soit $\pi$ comme ci-dessus et $\tilde{\pi}$ la 
repr\'esentation contragr\'ediente de $\pi$.
Soit $S$ un ensemble fini de places de $F$.

Si $\chi$ est un caract\`ere unitaire de ${\Bbb A}^{\times} / F^{\times}$ soit
$$L^S (s, (\pi \otimes \chi ) \times \tilde{\pi} ) = \prod_{v \notin S} L(s, (\pi_v \otimes \chi_v ) \times \tilde{\pi}_v ) .$$
On peut aussi former la fonction $L$ totale
$$L(s, (\pi \otimes \chi ) \times \tilde{\pi} )$$
qui jouit des propri\'et\'es donn\'ees par le Th\'eor\`eme \ref{fonction L globale}. On en d\'eduit que $L^S$ est holomorphe 
dans tout le plan complexe (sauf peut-\^etre en $s=0,1$ pour $\chi = 1$) car 
$$\prod_{v \in S} L(s, (\pi_v \otimes \chi_v ) \times \tilde{\pi}_v )^{-1}$$
est holomorphe.

Comme dans \cite{LRS} la d\'emonstration du Th\'eor\`eme \ref{LRS} va reposer sur le th\'eor\`eme 
suivant qui est le r\'esultat principal de \cite{LRS}~:
\begin{thm} \label{LRS2}
Soit $\beta$ un r\'eel  $> 1-\frac{2}{n^2 +1}$, il existe alors un ensemble infini ${\cal X}$ 
de caract\`eres $\chi$ de ${\Bbb A}^{\times} / F^{\times}$ tels que~:
\begin{enumerate}
\item $\chi_v = 1$ ($v \in S$)
\item $L^S (\beta , (\pi \otimes \chi )\times \tilde{\pi} ) \neq 0 .$
\end{enumerate}
\end{thm}

\medskip

Montrons maintenant que le Th\'eor\`eme \ref{LRS2} implique le Th\'eor\`eme \ref{LRS}.

Fixons dor\'enavant $v$ une place archim\'edienne de $F$ et prenons $S=\{ v \}$.
D'apr\`es le Th\'eor\`eme \ref{LRS2}, on peut choisir $\chi \in {\cal X}$ et tel que 
$L(s, (\pi \otimes \chi ) \times \tilde{\pi} )$ soit une fonction enti\`ere. Puisque $\chi_v =1$,
\begin{eqnarray}
L(s,(\pi \otimes \chi ) \times \tilde{\pi} ) = L (s, \pi_v \times \tilde{\pi}_v )  L^S (s,(\pi \otimes \chi ) \times \tilde{\pi} ).
\end{eqnarray}
Supposons que $s_0$, Re$(s_0 ) > 0$, soit un p\^ole de $L(s,\pi_v \times \tilde{\pi}_v )$. Alors $s_0$ doit \^etre 
un z\'ero de $L^S (s, (\pi \otimes \chi ) \times \tilde{\pi} )$. Donc la fonction $L(s, \pi_v \times \tilde{\pi}_v )$ doit 
\^etre holomorphe dans le demi-plan Re$(s) \geq 1 - \frac{2}{n^2 +1}$.

\medskip

Supposons d'abord que $v$ est une place r\'eelle. On \'ecrit $\pi_v = J(\sigma_{1,v} , \ldots , \sigma_{r,v} )$.
Fixons un entier $1\leq i\leq r$. Puisque $\pi_v$ est unitaire, on a vu que l'on peut supposer que Re$(t_{i,v} ) <0$. 
En utilisant les calculs des facteurs $L$ locaux rappel\'es \`a la section 
pr\'ec\'edente, on en d\'eduit que $L(s, \pi_v \times \tilde{\pi}_v )$ contient comme facteur 
$\Gamma (s +2{\rm Re}(t_{i,v} ))$ ou $\Gamma (\frac{s+ 2{\rm Re}(t_{i,v} )}{2} )$ et que les autres 
facteurs n'ont pas de z\'eros. On en d\'eduit que $-2$Re$(t_{i,v} ) >0$ est un p\^ole de ce facteur $L$ et donc d'apr\`es le 
paragraphe pr\'ec\'edent que $|$Re$(t_{i,v} ) |=-$Re$(t_{i,v} ) < \frac12 -\frac{1}{n^2 +1}$. 

\medskip

Supposons maintenant que $v$ est une place complexe et $\pi_v = J(\sigma_{1,v} ,\ldots , \sigma_{n,v} )$. Comme
dans le cas r\'eel, on peut choisir $t_{i,v}$ de partie r\'eelle strictement n\'egative. Le facteur local de la fonction
$L$ de Rankin-Selberg contient le facteur $\Gamma (s+ 2{\rm Re}(t_{i,v} ))$. L'argument pr\'ec\'edent s'applique de la 
m\^eme mani\`ere pour conclure que l\`a aussi $|$Re$(t_{i,v} )| < \frac12 -\frac{1}{n^2+1}$. Ce qui ach\`eve la d\'emonstration
du Th\'eor\`eme \ref{LRS}.

\markboth{CHAPITRE 7. TH\'EOR\`EME DE LUO-RUDNICK-SARNAK}{7.3. APPLICATION : LE TH\'EOR\`EME DE SELBERG}

\section{Application : le th\'eor\`eme de Selberg}

Le Th\'eor\`eme \ref{LRS} implique imm\'ediatement le th\'eor\`eme suivant.

\begin{thm}[Selberg] \label{selberg}
Le groupe $SL(2)_{|{\Bbb Q}}$ v\'erifie la Conjecture $A^- (0)$. 
\end{thm}
{\it D\'emonstration.} Une repr\'esentation du groupe $SL(2,{\Bbb R})$ qui intervient non trivialement dans la formule de Matsushima \'etendue (Th\'eor\`eme \ref{mat}) pour
un quotient $\Gamma \backslash {\Bbb H}$, o\`u $\Gamma $ est un sous-groupe de congruence de $SL(2)$, co\"{\i}ncide avec la composante \`a la place infinie d'une repr\'esentation
automorphe du groupe ${\Bbb Q}$-alg\'ebrique $SL(2)$.

Les repr\'esentations du groupe $SL(2,{\Bbb R})$ sont bien connues et se rel\`event en des repr\'esentations du groupe $GL(2,{\Bbb R})$ (c'est assez facile de le d\'eduire du \S 3.1).
Il en est de m\^eme des repr\'esentations automorphes nous appliqueront donc le Th\'eor\`eme \ref{LRS} au groupe $SL(2)$.

Le spectre du laplacien sur les fonctions d'un quotient $\Gamma \backslash {\Bbb H}$ est constitu\'e d'une partie continue $[\frac{1}{4} , +\infty[$ et d'une partie discr\`ete 
qu'il nous faut contr\^oler. Via la formule de Matsushima \'etendue la partie discr\`ete du spectre correspond aux repr\'esentations automorphes cuspidales.
Puisqu'\`a une repr\'esentation temp\'er\'ee de $SL(2, {\Bbb R})$, il correspond (toujours via la formule de Matsushima) 
un espace de fonctions propres du laplacien de valeur propre $\geq \frac{1}{4}$, il ne nous reste qu'\`a consid\'erer les repr\'esentations de la s\'erie compl\'ementaire
$$\pi_s = {\rm ind}_B^{SL(2,{\Bbb R})} (|x|^s )$$
o\`u $B= \left\{ \left( 
\begin{array}{cc}
x & * \\
0 & x^{-1} 
\end{array} \right) \right\}$ et $s \in [0, \frac12 [$. Mais d'apr\`es le Th\'eor\`eme \ref{LRS}, si une repr\'esentation automorphe de $SL(2)$ a pour composante \`a l'infini $\pi_s$ 
alors $s\leq \frac12 - \frac{1}{2^2 +1} = \frac{3}{10}$. Autrement dit la repr\'esentation $\pi_s$ est uniform\'ement (par rapport au choix du groupe $\Gamma$) isol\'ee
de la repr\'esentation triviale. Ce qui conclut la d\'emonstration du Th\'eor\`eme \ref{selberg}.

\medskip

La d\'emonstration ci-dessus donne l'estim\'ee $\varepsilon ( SL(2) , 0 )= \frac{4}{25}$ qui est moins bon que l'estim\'ee original de Selberg \'egal \`a $\frac{3}{16}$. Il y a 
n\'eanmoins deux grands avantages \`a la m\'ethode ci-dessus 1) elle est valable sur n'importe quel corps de nombres 2) elle est valable pour tous les groupes $GL(n)$.
Et cette m\'ethode intervient effectivement dans la d\'emonstration r\'ecente par Kim et Sarnak de l'estim\'ee $\varepsilon ( SL(2) , 0 )= \frac{1}{4} - \left( \frac{7}{64}\right)^2$.

\medskip

Gr\^ace \`a la correspondance de Jacquet-Langlands \cite{JacquetLanglands} entre les repr\'esentations automorphes du groupe des unit\'es d'une alg\`ebre de quaternions et 
celles du groupe $SL(2)$, le Th\'eor\`eme \ref{selberg} (\'etendu \`a un corps de nombre quelconque)
permet de v\'erifier la Conjecture $A^- (0)$ pour tout groupe $G$ comme dans la Conjecture $A$ avec $G^{{\rm}}$
localement isomorphe au groupe $SL(2,{\Bbb R})$ ou $SL(2 ,{\Bbb C} )$. Le principe de restriction de Burger et Sarnak (\S 2.3) permet alors de v\'erifier la 
Conjecture $A^- (0)$ pour tous les groupes $G$ tels que $G^{{\rm nc}}$ soit localement isomorphe au groupe $SO(n,1)$.

\newpage

\markboth{CHAPITRE 8. D\'EMONSTRATION DU TH\'EOR\`EME 1}{8.1. DESCRIPTION DU GROUPE $G$}

\chapter{D\'emonstration du Th\'eor\`eme~1}

On est maintenant amen\'e \`a d\'emontrer le Th\'eor\`eme 1 annonc\'e dans l'Introduction. Rappelons 
son \'enonc\'e~:

\medskip

\noindent
{\bf Th\'eor\`eme 1} {\it Soit $G$ un groupe alg\'ebrique connexe et presque simple sur ${\Bbb Q}$ tel que $G^{{\rm nc}}$ 
soit isomorphe au groupe $SU(2,1)$.
Alors, pour tout sous-groupe de congruence $\Gamma$ dans $G$,
$$\lambda_1^0 (\Gamma \backslash X_G ) \geq \frac{84}{25} $$
et 
$$\lambda_1^1 (\Gamma \backslash X_G ) \geq \frac{9}{25} .$$
}

\medskip

Rappelons avant d'en donner une d\'emonstration compl\`ete que ce r\'esultat est essentiellement contenu dans 
l'article \cite{HarrisLi} de Harris et Li.

\section{Description du groupe $G$}

Soit $G$ un groupe alg\'ebrique connexe et presque simple sur ${\Bbb Q}$ modulo son centre tel que $G^{{\rm nc}}$ 
soit isomorphe au groupe $U(2,1)$.
Remarquons que l'on pr\'ef\`ere consid\'erer le groupe $U(2,1)$ plut\^ot que le groupe $SU(2,1)$ pour des 
raisons techniques; il est n\'eanmoins bien \'evident que cela ne change pas 
l'\'enonc\'e du Th\'eor\`eme 1 puisque l'espace sym\'etrique associ\'e est le m\^eme (pour un argument plus pr\'ecis voir \cite{Clozel}).

Alors, il existe un entier $d\geq 1$ tel que 
\begin{eqnarray} \label{unitaire}
G ({\Bbb R} ) \cong U(2,1) \times U(3)^{d-1} .
\end{eqnarray} 

D\'ecrivons comment construire un tel groupe $G$.
Soit $F$ un corps de nombres totalement r\'eel de degr\'e $d$, et soit $K$ une extension quadratique totalement imaginaire de $F$.
On notera $x \mapsto \overline{x}$ la conjugaison de $K$ par rapport \`a $F$. 
Fixons un plongement  r\'eel $\sigma_1$ de $F$ et une extension $\tau_1$ de $\sigma_1$ en un plongement complexe de $K$. Soit $V=K^3$
\'equip\'e d'une forme hermitienne $h$; on 
suppose $h$ de signature $(2,1)$ en $\sigma_1$ (c'est \`a dire relativement au plongement $\tau_1$)
et d\'efinie positive en tous les autres plongements r\'eels de $F$. On notera $\Phi \in GL_3 (K)$ la 
matrice (hermitienne) de $h$, on a alors ${}^t \overline{\Phi} = \Phi$. Soit $G = U_{\Phi}$ le groupe des 
transformations unitaires de $V$~:
$$G = \{ g \in GL_3 (K) \; : \; g \Phi {}^t \overline{g} = \Phi \}.$$
On consid\`ere le groupe $G$ comme un groupe ${\Bbb Q}$-alg\'ebrique par restriction des scalaires \cite{Zimmer}
de $F$ \`a ${\Bbb Q}$. Alors le groupe $G ({\Bbb R})$ est un groupe de Lie semi-simple connexe et 
v\'erifie (\ref{unitaire}). 

Il existe une autre construction de groupes v\'erifiant (\ref{unitaire}).

Soit $D$ une alg\`ebre semi-simple de dimension $9$ sur $K$. Soit $\alpha$ une involution de $D$. 
On dit de la paire $(D, \alpha )$ qu'elle est {\it de seconde esp\`ece} \index{seconde esp\`ece, involution de} si l'involution $\alpha$ est de seconde 
esp\`ece {\it i.e.}, $\alpha$ est un anti-automorphisme de $D$ dont la restriction \`a $K$ est $\sigma$.
Si $D$ est une alg\`ebre simple de centre $K$, on peut d\'efinir le groupe alg\'ebrique sur $F$
$$U(F) = \{ g \in D^* \; : \; \alpha (g) g=1 \} .$$
Supposons de plus qu'\`a toutes les places r\'eelles $v$ sauf (exactement) une, le groupe $U(F)$ est compact
isomorphe \`a $U(3)$. Alors le groupe ${\Bbb Q}$-alg\'ebrique $G$, obtenu \`a partir de $U(F)$   
par restriction des scalaires de $F$ \`a ${\Bbb Q}$ v\'erifie 
(\ref{unitaire}). R\'eciproquement, tout groupe $G$ alg\'ebrique sur ${\Bbb Q}$, dont les points r\'eels forment
un groupe de Lie semi-simple connexe v\'erifiant $G({\Bbb R}) \cong SU(2,1) \times SU(3)^{d-1}$ s'obtient par la construction 
rappel\'ee ci-dessus (en prenant le groupe d\'eriv\'e),
cf. \cite{PlatonovRapinchuk}.

Remarquons que si $D=M_3 (K)$, alors $\alpha (g) = \Phi {}^t \overline{g} \Phi^{-1}$ pour une certaine forme hermitienne
$\Phi$ (voir \cite{PlatonovRapinchuk}). La construction ci-dessus contient donc comme cas particulier la 
premi\`ere construction d\'ecrite plus haut. On distinguera d'ailleurs deux cas dans la construction ci-dessus, 
le cas o\`u $D=M_3 (K)$ et le cas o\`u $D \neq M_3 (K)$.

\markboth{CHAPITRE 8. D\'EMONSTRATION DU TH\'EOR\`EME 1}{8.2. $D=M_3 (K)$}

\section{$D=M_3 (K)$}

Dans cette section nous supposons que $G$ est obtenu par une construction d\'ecrite dans la section
pr\'ec\'edente avec comme alg\`ebre $D=M_3 (K)$. Nous conservons les notations de la section 
pr\'ec\'edente. Nous d\'emontrons le Th\'eor\`eme 1 dans ce cas.

\medskip

Le groupe $G ({\Bbb R} )$ est isomorphe \`a $U(2,1) \times U(3)^{d-1}$, une repr\'esentation qui intervient 
non trivialement dans la formule de Matsushima \'etendue (Th\'eor\`eme \ref{mat}) pour un quotient $\Gamma \backslash X_G$, o\`u $\Gamma$
est un sous-groupe de congruence de $G$, est de la forme
\begin{eqnarray} \label{rep}
\sigma \otimes 1 \otimes \ldots \otimes 1,
\end{eqnarray}
o\`u $\sigma$ est une repr\'esentation unitaire du groupe (r\'eel) $U(2,1)$.
De plus elle co\"{\i}ncide avec le produit tensoriel des composantes aux places infinies d'une repr\'esentation automorphe 
du groupe $F$-alg\'ebrique $U(F)$.

Noter qu'alors, plus simplement, $\sigma$ apparait dans $L^2 (\Gamma ' \backslash U(2,1) )$ o\`u $\Gamma '$ se d\'eduit de fa\c{c}on \'evidente de 
$\Gamma$. On dira que $\Gamma'$ ``est un sous-groupe de congruence de $G^{{\rm nc}} = U(2,1)$''.

Les repr\'esentations du groupe $U(2,1)$ sont compl\`etement d\'ecrites au \S 4.6. 
Soit d'abord une repr\'esentation temp\'er\'ee $\sigma$ de $U(2,1)$. S'il correspond \`a la repr\'esentation 
$\sigma \otimes 1 \otimes \ldots \otimes 1$ un espace de fonctions $\lambda$-propres (resp. de $1$-formes $\lambda$-propres), 
pour le laplacien de Hodge, via la formule de Matsushima, alors $\lambda \geq 4$ (resp. $\lambda \geq 1$) (cf. Th\'eor\`eme 4.1.1).
Il nous suffit donc de consid\'erer les repr\'esentations (\ref{rep}) avec $\sigma$ non temp\'er\'ee (et unitaire).
Dans la classification du \S 4.6 une telle repr\'esentation est associ\'ee \`a un param\`etre $\varphi$ correspondant
\`a un triplet $(u,v,\mu )$ o\`u $\mu \in {\Bbb Z}$, $u,v \in {\Bbb C}$, $u-v \in {\Bbb Z}$, Re$(u+v) \geq 0$ et si $u+v=0$
alors $u\in \frac12 + {\Bbb Z}$. De plus d'apr\`es la description des $K$-types du \S 4.5,  
la repr\'esentation $\sigma \otimes 1 \otimes \ldots \otimes 1$ intervient non trivialement, via la formule de Matsushima,
dans la description du spectre sur les fonctions (resp. sur les $1$-formes) seulement si 
$u-v=0$ (resp. $u-v = \pm 3$). Finalement puisque seules les repr\'esentations unitaires et de dimension 
infinie sont \`a consid\'erer, il nous
reste \`a \'etudier le cas o\`u $\sigma$ est l'unique quotient irr\'eductible de l'induite unitaire du caract\`ere
$$\chi_{\varphi} = (u,v,\mu ),$$
o\`u le triplet $(u,v, \mu )$ est astreint \`a v\'erifier l'une des alternatives suivantes (\pa 4.6), la valeur propre du Casimir \'etant calcul\'ee par (4.1.11).
\begin{enumerate}
\item $(u,v,\mu)=(2,-1,-1)$, et alors (\ref{rep}) intervient dans le spectre sur les formes diff\'erentielles de bi-degr\'e 
$(1,0)$, pour la valeur propre $0$;
\item $(u,v,\mu)=(-1,2,1)$, et alors (\ref{rep}) intervient dans le spectre sur les formes diff\'erentielles de bi-degr\'e
$(0,1)$, pour la valeur propre $0$;
\item $(u,v,\mu )=(\frac{s}{2},\frac{s}{2},0)$ pour un nombre reel $s\in ]0,2[$, et alors (\ref{rep}) intervient dans le spectre sur les
fonctions, pour la valeur propre $4-s^2 $;
\item $(u,v,\mu )=(\frac{s+3}{2}, \frac{s-3}{2} ,-1)$ pour un nombre r\'eel $s\in ]0,1[$, et alors (\ref{rep}) intervient 
dans le spectre sur les formes diff\'erentielles de bi-degr\'e $(1,0)$, pour la valeur propre $1-s^2$;
\item $(u,v,\mu )=(\frac{s-3}{2}, \frac{s+3}{2} ,1)$ pour un nombre r\'eel $s\in ]0,1[$, et alors (\ref{rep}) intervient
dans le spectre sur les formes diff\'erentielles de bi-degr\'e $(0,1)$, pour la valeur propre $1-s^2$.
\end{enumerate} 

Il nous reste \`a d\'emontrer que si $\sigma$ intervient dans la repr\'esentation r\'eguli\`ere $L^2 (\Gamma \backslash 
G^{{\rm nc}})$ pour un certain sous-groupe de congruence $\Gamma$ de $G$ et si $\sigma$ est 
du type 3, 4 ou 5 ci-dessus alors $s\leq 4/5$.

Mais si $\sigma$ intervient dans la repr\'esentation r\'eguli\`ere $L^2 (\Gamma \backslash G^{{\rm nc}} )$ 
pour un certain sous-groupe de congruence $\Gamma$ de $G$ alors il existe une repr\'esentation unitaire
automorphe cuspidale $\pi = \otimes \pi_v $ de $U({\Bbb A}_F)$
telle que 
$$\otimes_{ v \ {\rm infinie}} \pi_{v} = \sigma \otimes 1 \otimes \ldots \otimes 1.$$
Seule la place infinie $v$ o\`u $\pi_v \cong \sigma$ nous interessera, pour simplifier nous la noterons $v= \infty$. 
D'apr\`es le Th\'eor\`eme 13.3.3 de \cite{Rogawski} et en tenant compte du fait que $\pi_{\infty}$ n'est pas temp\'er\'ee, 
\begin{itemize}
\item ou bien il correspond \`a la repr\'esentation $\pi$ une repr\'esentation automorphe cuspidale $\Pi$ de $GL(3)$ sur 
le corps $K$ telle que $\Pi_{\infty}$ soit obtenue par la fonctorialit\'e d\'ecrite au \S 4.4 (``changement de base''); \index{changement de base}
\item ou bien il existe une repr\'esentation automorphe $\rho$ de $U(1,1) \times U(1)$ sur le corps $K$ telle que 
$\pi_{\infty}$ soit obtenue \`a partir de $\rho_{\infty}$ par la fonctorialit\'e induite par $\xi_n$ pour un certain entier
$n$ (cf. \S 4.6).
\end{itemize}

\medskip

Dans le premier cas, la fonctorialit\'e associe \`a un param\`etre $\varphi : W_{{\Bbb R}} \rightarrow GL_3 ({\Bbb C} ) \rtimes W_{{\Bbb R}}$
un morphisme $\varphi_{|{\Bbb C}^*} \rightarrow GL_3 ({\Bbb C})$ admissible pour  $GL_3 ({\Bbb C})$. \`A ce dernier morphisme
est associ\'ee une repr\'esentation admissible de $GL_3 ({\Bbb C})$. Lorsque le param\`etre $\varphi$ correspond 
au triplet $(u,v,\mu )$, d'apr\`es le Th\'eor\`eme \ref{correspondance locale}, la repr\'esentation de $GL_3 ({\Bbb C})$ 
correspondant au param\`etre $\varphi_{|{\Bbb C}^*}$ est $J(\chi_1 , \chi_2 , \chi_3 )$ o\`u 
$\chi_1 (z) = z^u \overline{z}^v$, $\chi_2 (z) = z^{\mu} \overline{z}^{-\mu}$ et $\chi_3 (z) = z^{-v} \overline{z}^{-u}$.
Remarquons tout d'abord que si le param\`etre $\varphi$ correspond \`a la repr\'esentation triviale, {\it i.e.} 
si $(u,v,\mu)=(1,1,0)$, alors la repr\'esentation obtenue par fonctorialit\'e est la repr\'esentation triviale.
La repr\'esentation automorphe cuspidale $\Pi$ de $GL(3)$ obtenue par fonctorialit\'e \`a partir de $\pi$ v\'erifie 
le Th\'eor\`eme de Luo, Rudnick et Sarnak. Donc si $\pi_{\infty}$ est de la forme (\ref{rep}), avec $\sigma$ 
du type 3, 4 ou 5 ci-dessus, le Th\'eor\`eme  7.0.1 appliqu\'e \`a $\chi_1$ implique que le nombre r\'eel $s$ doit v\'erifier $\frac{s}{2} \leq \frac12 - \frac{1}{10}$, 
{\it i.e.} $s \leq 4/5$.

\medskip

Dans le deuxi\`eme cas, il existe un entier $n$ et un param\`etre $\eta : W_{{\Bbb R}} \rightarrow (GL_2 ({\Bbb C})
\times GL_1 ({\Bbb C}) )\rtimes W_{{\Bbb R}}$ correspondant soit \`a une repr\'esentation 
automorphe de $U(1,1) \times U(1)$ soit \`a une repr\'esentation unitaire de dimension un et tel que le param\`etre 
$\varphi  = \xi_n \circ \eta $ corresponde
au triplet $(u,v, \mu )$. On identifie les param\`etres $\eta$ possibles \`a l'aide de la Proposition 4.3.3. 
Un tel param\`etre $\eta$ associe \`a un \'el\'ement $z \in {\Bbb C}^*$, 
$$\eta (z) = \left(
\left(
\begin{array}{cc}
z^{\alpha} \overline{z}^{\beta} & \\
 & z^{-\beta} \overline{z}^{-\alpha} 
\end{array} \right) , \left( \frac{z}{\overline{z}} \right)^{\mu} \right) ,$$
o\`u $\alpha$ et $\beta$ sont deux nombres r\'eels tels que $\alpha - \beta \in {\Bbb Z}$,
$\alpha +\beta  \geq 0$ et si $\alpha + \beta =0$, $\alpha -\beta \in 2{\Bbb Z}$.
On a alors 
$$(u,v, \mu ) = (\alpha + \frac12 +n , \beta -\frac12 -n , \mu ).$$

D'abord si la repr\'esentation du groupe $U(1,1) \times U(1)$ correspondant au param\`etre $\eta$ 
est unitaire de dimension $1$, alors $\alpha$, $\beta \in  \frac12 + {\Bbb Z}$ et $\alpha + \beta =1$, donc $\alpha 
-\beta$ est paire. Or $u-v = \alpha - \beta +1 +2n$, et puisque $u-v=0$, $3$ ou $-3$, on doit avoir $u$ et $v$ 
entiers, cas que l'on a exclu.  

Puis, si la repr\'esentation du groupe $U(1,1) \times U(1)$ correspondant au param\`etre $\eta$ est la 
partie \`a l'infinie d'une repr\'esentation automorphe, le Th\'eor\`eme de Selberg implique que
$s = u+v = \alpha + \beta \leq \frac{1}{4}$ \footnote{On a remarqu\'e que le Th\'eor\`eme de Luo, Rudnick et 
Sarnak implique lui aussi une majoration sur $s$ : $s \leq \frac{3}{10}$. Celle-ci nous suffit pour d\'emontrer 
le Th\'eor\`eme 1.}.

\medskip

Dans tous les cas le nombre $s$ est inf\'erieur ou \'egal \`a $\frac{4}{5}$. Et le Th\'eor\`eme 1 est d\'emontr\'e.

\markboth{CHAPITRE 8. D\'EMONSTRATION DU TH\'EOR\`EME 1}{8.3. $D \neq M_3 (K)$}

\section{$D \neq M_3 (K)$}

Dans cette section nous supposons que $G$ est obtenu par la construction d\'ecrite dans la premi\`ere section
avec une alg\`ebre $D \neq M_3 (K)$. Nous conservons les notations de la premi\`ere section. 
Nous d\'emontrons le Th\'eor\`eme 1 dans ce cas.

\medskip

Comme dans la section pr\'ec\'edente, le groupe $G ({\Bbb R} )$ est isomorphe \`a $U(2,1) \times U(3)^{d-1}$, 
une repr\'esentation qui intervient non trivialement dans la formule de Matsushima d'un quotient $\Gamma \backslash X_G$, o\`u $\Gamma$
est un sous-groupe de congruence de $G$, est de la forme
\begin{eqnarray} \label{rep2}
\sigma \otimes 1 \otimes \ldots \otimes 1,
\end{eqnarray}
o\`u $\sigma$ est une repr\'esentation unitaire du groupe (r\'eel) $G^{{\rm nc}} = U(2,1)$.
De plus elle co\"{\i}ncide avec le produit tensoriel des repr\'esentations aux places infinies d'une repr\'esentation automorphe 
du groupe $F$-alg\'ebrique $U(F)$.

Nous passons en revue, comme dans le \pa 8.2, les possibilit\'es 1 -- 5.

Il nous reste \`a d\'emontrer que si $\sigma$ intervient dans la repr\'esentation r\'eguli\`ere $L^2 (\Gamma \backslash 
G^{{\rm nc}})$ et est du type 3, 4 ou 5 ci-dessus alors $s\leq 4/5$.

Mais, comme dans la section pr\'ec\'edente, si $\sigma$ intervient dans $L^2 (\Gamma \backslash G^{{\rm nc}})$ il existe alors une repr\'esentation unitaire
automorphe cuspidale $\pi = \otimes \pi_v $ du groupe $F$-alg\'ebrique $U(F)$ 
telle que pour une certaine place \`a l'infini $\pi_{\infty} = \sigma$.
Le groupe $U(F)$ est, dans cette section, une forme int\'erieure du groupe unitaire quasi-d\'eploy\'e associ\'e 
\`a l'alg\`ebre $M_3 (F)$. Le Th\'eor\`eme 14.6.3 et la Proposition 14.6.2 de \cite{Rogawski} \'etablissent une bijection 
entre les repr\'esentations automorphes (en fait plus exactement les $L$-paquets ce qui revient au m\^eme pour les
repr\'esentations que l'on consid\`ere) de $U(F)$ et certaines repr\'esentations automorphes de sa forme
quasi-d\'eploy\'ee. On est donc ramen\'e au cas o\`u $D=M_3 (F)$, cas que l'on a trait\'e dans la section 
pr\'ec\'edente. Le Th\'eor\`eme 1 est donc d\'emontr\'e.

\markboth{CHAPITRE 8. D\'EMONSTRATION DU TH\'EOR\`EME 1}{8.4. CONJECTURE DE CHANGEMENT DE BASE}

\section{La Conjecture de changement de base de Harris et Li}

Commen\c{c}ons par \'enoncer un affaiblissement des r\'esultats de Rogawski \cite{Rogawski}
qui synth\'etise en un seul \'enonc\'e les deux cas trait\'es dans les sections pr\'ec\'edentes.
Pour cela conservons les notations de la premi\`ere section. Pour simplifier les \'enonc\'es nous dirons
d'une repr\'esentation de $GL(n,{\Bbb C})$ qu'elle est {\it automorphe (cuspidale)} \index{repr\'esentation automorphe cuspidale} s'il existe un corps de nombre $K$, une 
repr\'esentation automorphe (cuspidale) $\Pi$ de $GL(n, {\Bbb A}_K )$ et une place archim\'edienne complexe $v$ de $K$ telle
que $\Pi_v \cong \pi$. Ce qui compte pour nous est qu'une repr\'esentation automorphe cuspidale de $GL(n,{\Bbb C})$ doit
donc v\'erifier les conclusions du Th\'eor\`eme de Luo, Rudnick et Sarnak (Chapitre 7). Rappelons que plus g\'en\'eralement
la Conjecture de Ramanujan g\'en\'eralis\'ee en les places archim\'ediennes (complexes) s'\'enonce 
de la fa\c{c}on suivante.

\medskip

\noindent
{\bf Conjecture de Ramanujan g\'en\'eralis\'ee} {\it Soit $\pi$ une repr\'esentation automorphe cuspidale de $GL(n,{\Bbb C})$.
Alors $\pi$ est une repr\'esentation induite
$$I(\sigma_1 , \ldots , \sigma_n ) = J(\sigma_1 , \ldots , \sigma_n ),$$
o\`u les $\sigma_j$ sont comme dans le Th\'eor\`eme \ref{class complexe} et tels que pour tout $j=1, \ldots , n$,
$t_j$ soit imaginaire pur.}
\index{Conjecture de Ramanujan}

\medskip 

Soit $G$ un groupe alg\'ebrique simple et connexe sur ${\Bbb Q}$ tel que $G^{{\rm nc}}$ soit isomorphe au groupe $U(2,1)$. Remarquons que toute repr\'esentation 
$\pi$ de $G({\Bbb R})$ se restreint en une repr\'esentation, que l'on notera, $\pi_1$ du groupe $U(2,1)$.

Le th\'eor\`eme suivant est une synth\`ese des particularisations de r\'esultats de Rogawski, voir les 
Chapitres 13 et 14 de \cite{Rogawski}, que 
nous avons utilis\'es dans les deux sections pr\'ec\'edentes. 

\begin{thm}[Rogawski] \label{roga}
Soit $\pi$ une repr\'esentation irr\'eductible de dimension infinie appartenant \`a $\widehat{G}_{\rm Aut}$.  
Alors, la repr\'esentation de $GL(3,{\Bbb C} )$ obtenue par changement de base (\S 4.4) \`a partir de $\pi_{1}$
a le m\^eme caract\`ere infinit\'esimal que la repr\'esentation 
$${\rm ind}_P^{GL(3,{\Bbb C} )} (J( \tau_1 , b_1 ) \otimes \ldots \otimes J(\tau_m , b_m ) ),$$
o\`u 
\begin{enumerate}
\item $P$ est le sous-groupe parabolique de $GL (3, {\Bbb C})$ associ\'e \`a la partition 
$3=r_1 + \ldots + r_m$;
\item chaque $r_i = a_i b_i$ avec $a_i, b_i \in {\Bbb N}$;
\item $\tau_i$ est une repr\'esentation automorphe cuspidale de $GL(a_i , {\Bbb C} )$ obtenue par 
changement de base \`a partir d'une repr\'esentation de $U(1)$ si $a_i =1$ et obtenue par le changement de 
base non standard $\beta_0 $ de (\ref{betan}) \`a partir d'une repr\'esentation de $U(1,1)$ si $a_i = 2$; et,
\item $J( \tau_i , b_i )$ d\'esigne la repr\'esentation de $GL (r_i , {\Bbb C} )$ associ\'ee, comme au \S 3.4.4.
\end{enumerate}
\end{thm} \index{Th\'eor\`eme de Rogawski}
{\it D\'emonstration.} Montrons que le Th\'eor\`eme \ref{roga} se d\'eduit des Chapitres 13 et 14 de \cite{Rogawski}.
Le Th\'eor\`eme 14.6.3 et la Proposition 14.6.2 de \cite{Rogawski} permettent, comme dans la section pr\'ec\'edente, de 
ramener la d\'emonstration aux cas o\`u $G$ est obtenue par restriction des scalaires \`a partir d'un groupe 
unitaire associ\'e \`a une alg\`ebre $D=M_3 (F)$. Dans ce cas la repr\'esentation $\pi_1$ est la repr\'esentation 
de $U(2,1)$ obtenue en la place archim\'edienne correspondante. Pour les besoins de la d\'emonstration notons 
$\pi$ une repr\'esentation automorphe du groupe unitaire $U({\Bbb A}_F )$ telle que $\pi_{\infty} =\pi_1$ (o\`u 
$\infty$ d\'esigne toujours la place archim\'edienne en laquelle le groupe unitaire est $U(2,1)$). On peut de plus
supposer que $\pi$ intervient discr\`etement dans $L^2 (U(F) \backslash U({\Bbb A}_F ) )$.
Puisque $\pi$ est de dimension infinie, les Th\'eor\`emes 13.3.3 et 13.3.5 de \cite{Rogawski} distinguent alors trois cas~:
\begin{enumerate}
\item La repr\'esentation $\pi$ est stable, et alors il existe une repr\'esentation automorphe cuspidale $\Pi $ de $GL(3, {\Bbb A}_K )$
telle que $\Pi_{\infty}$, vue comme repr\'esentation de $GL(3,{\Bbb C})$, soit obtenue \`a partir de $\pi_{\infty}$ par 
changement de base, comme au \S 4.4. Dans ce cas on peut prendre 
$m=b_1= 1$, $r_1 = a_1=3$ et $\tau_1 = \Pi_{\infty}$.
\item La repr\'esentation $\pi$ est endoscopique du premier type, et alors il existe une repr\'esentation automorphe cuspidale $\rho$ de 
$U(1,1) \times U(1)$ (vu comme groupe alg\'ebrique sur $F$) telle que $\pi_{\infty}$ s'obtienne 
par fonctorialit\'e \`a partir de $\rho_{\infty}$, comme au \S 4.6. 
D'apr\`es le changement de base pour $U(1,1)$ \cite[Proposition 11.4.1]{Rogawski}, il correspond \`a la repr\'esentation automorphe $\rho$
une repr\'esentation automorphe cuspidale $\tilde{\rho}$ de 
$GL(2, {\Bbb A}_K ) \times GL(1, {\Bbb A}_K )$ telle que $\tilde{\rho}_{\infty}$,
vue comme repr\'esentation de $GL(2,{\Bbb C}) \times GL(1,{\Bbb C})$, soit obtenue \`a partir de $\rho_{\infty}$ par
changement de base. Dans ce cas la repr\'esentation $\pi_{\infty}$ est associ\'ee \`a un param\`etre 
$\varphi : W_{{\Bbb R}} \rightarrow GL(3,{\Bbb C}) \rtimes W_{{\Bbb R}}$ obtenue comme  
$\xi_n \circ \eta$ pour un certain entier $n$ et un certain param\`etre $\eta : W_{{\Bbb R}} \rightarrow (GL(2,{\Bbb C})  
\times GL(1,{\Bbb C} ) ) \rtimes W_{{\Bbb R}}$. Autrement dit, la repr\'esentation de $GL(3,{\Bbb C})$  obtenue
\`a partir de $\pi_{\infty}$ par changement de base, est associ\'ee \`a un param\`etre 
$\psi : {\Bbb C}^* \rightarrow GL(3,{\Bbb C}) \times {\Bbb C}^*$ \'egal \`a la compos\'ee $\xi_n \circ \eta$ (pour un 
certain entier $n$) o\`u 
$\eta : {\Bbb C}^* \rightarrow (GL(2,{\Bbb C}) \times GL(1,{\Bbb C}) ) \times {\Bbb C}^*$ est un $L$-param\`etre associ\'e
\`a une repr\'esentation automorphe cuspidale de $GL(2,{\Bbb C}) \times GL(1,{\Bbb C})$. Il lui correspond alors
une repr\'esentation automorphe cuspidale $\tau_1$ (resp. $\tau_2$) de $GL(2,{\Bbb C})$ (resp. $GL(1,{\Bbb C})$).
La repr\'esentation $\tau_1$ est bien \'evidemment obtenue par changement de base non standard \`a partir d'une 
repr\'esentation de $U(1,1)$ et la repr\'esentation $\tau_2$ est obtenue par changement de base (standard) \`a partir
d'une repr\'esentation de $U(1)$. 
Puisque, d'apr\`es (3.1.9), le caract\`ere infinit\'esimal d'une repr\'esentation de $GL(n,{\Bbb C})$ est cod\'e par son 
$L$-param\`etre, on conclut facilement la d\'emonstration du Th\'eor\`eme \ref{roga} dans ce cas en prenant,
$m=2$, $r_1 =a_1= 2$, $r_2 =1$.
\item La repr\'esentation $\pi$ est endoscopique du deuxi\`eme type, et alors il existe une repr\'esentation automorphe de dimension un $\rho$ de 
$U(1,1) \times U(1)$ (vu comme groupe alg\'ebrique sur $F$) telle que $\pi_{\infty}$ s'obtienne 
par fonctorialit\'e \`a partir de $\rho_{\infty}$, comme d\'ecrit au \S 4.6. Mais, dans ce cas, Rogawski montre que 
$\rho_{\infty}$ s'obtient, par fonctorialit\'e, \`a partir d'une repr\'esentation de $U(1) \times U(1) \times U(1)$. 
On est dans la configuration $m=3$, $r_1 = r_2 = r_3 = 1$.  
\end{enumerate}

De mani\`ere g\'en\'eral on s'attend, voir par exemple \cite{HarrisLi} pour un \'enonc\'e l\'eg\`erement diff\'erent, 
\`a ce que la conjecture suivante soit v\'erifi\'ee.

\medskip

\noindent
{\bf Conjecture de changement de base} {\it Soit $G$ un groupe alg\'ebrique simple et connexe sur ${\Bbb Q}$ tel que 
$G^{{\rm nc}}$ soit isomorphe au groupe $U(p,q)$. Soit $\pi$ une repr\'esentation de dimension infinie appartenant \`a  $\widehat{G}_{\rm Aut}$ et $\pi_1$ 
la restriction de $\pi$ au groupe $U(p,q)$.
Alors, la repr\'esentation de $GL(p+q,{\Bbb C} )$ obtenue par changement de base (\S 4.4) \`a partir de $\pi_{1}$
a le m\^eme caract\`ere infinit\'esimal que la repr\'esentation 
$${\rm ind}_P^{GL(p+q,{\Bbb C} )} (J( \tau_1 , b_1 ) \otimes \ldots \otimes J(\tau_m , b_m )),$$
o\`u 
\begin{enumerate}
\item $P$ est le sous-groupe parabolique de $GL (p+q, {\Bbb C})$ associ\'e \`a la partition 
$p+q=r_1 + \ldots + r_m$;
\item chaque $r_i = a_i b_i$ avec $a_i, b_i \in {\Bbb N}$;
\item $\tau_i$ est une repr\'esentation automorphe cuspidale de $GL(a_i , {\Bbb C} )$; et,
\item $J( \tau_i , b_i )$ d\'esigne la repr\'esentation de $GL (r_i , {\Bbb C} )$ associ\'ee, comme au \S 3.4.4.
\end{enumerate}
}
\index{Conjecture de changement de base}

\medskip

Comme il a d\'ej\`a \'et\'e remarqu\'e par Harris et Li, la Conjecture de changement de base implique la 
Conjecture A$^-$ pour les groupes du type $U(n,1)$. Nous montrons plus g\'en\'eralement le th\'eor\`eme suivant.

\begin{thm} \label{csqce du chgt de base}
Soit $G$ un groupe comme dans la Conjecture de changement de base et v\'erifiant les conclusions de celle-ci.
Il existe alors un r\'eel strictement positif $\varepsilon = \varepsilon (p,q) >0$ tel que  
pour tout sous-groupe de congruence $\Gamma$ dans $G$ et pour tout entier $i$,
$$\lambda_1^i (\Gamma \backslash X_G ) \geq \varepsilon .$$
\end{thm}
{\it D\'emonstration.} Soit $(P,\sigma , \nu )$ un triplet constitu\'e d'un sous-groupe parabolique cuspidal $P$ de $U(p,q)$ 
avec une decomposition de Langlands $P=MAN$, une repr\'esentation $\sigma$ de la s\'erie discr\`ete de $M$ et 
un \'el\'ement $\nu \in \mathfrak{a}^*$ (o\`u, comme d'habitude, on a not\'e $\mathfrak{a}$ l'alg\`ebre de Lie complexe
de $A$). On d\'efinit
$$I(\sigma , \nu ) = {\rm ind}_P^G (\sigma \otimes e^{\nu} \otimes 1).$$
On a vu au Chapitre 5 que toute repr\'esentation cohomologique de $U(p,q)$ s'obtient comme quotient de Langlands 
d'un certain $I(\sigma , \nu )$. D'apr\`es \cite{Vogan2}, on se trouve n\'ecessairement dans l'un des deux cas suivant~: 
\begin{itemize}
\item soit la repr\'esentation cohomologique est isol\'ee dans le dual unitaire de $U(p,q)$;
\item soit il existe une suite $(s_i ) \in ( \mathfrak{a}^* )^{{\Bbb N}^*}$ (correspondant \`a des caract\`eres de $A$) qui converge
vers $\nu$ telle que que chaque $I(\sigma , s_i )$ poss\`ede un sous-quotient irr\'eductible unitaire.
\end{itemize} 
Pour conclure la d\'emonstration du Th\'eor\`eme \ref{csqce du chgt de base}, il nous reste \`a montrer que dans le deuxi\`eme
cas les repr\'esentations cohomologiques sont  isol\'ees  dans le dual automorphe. 

Soit $P=MAN$ un sous-groupe parabolique de $U(p,q)$ associ\'e, comme au Chapitre 5 (dont nous conservons les notations),
\`a une repr\'esentation cohomologique et soient $\sigma$ la repr\'esentation de la s\'erie discr\`ete de $M$ et $\nu$ l'\'el\'ement de 
$\mathfrak{a}^*$ associ\'es.
Le sous-groupe $A$ est isomorphe au groupe $({\Bbb R}^* )^{\sum_j d_j}$, on \'ecrit chaque \'el\'ement $s$ de 
$\mathfrak{a}^*$ comme $s= (s(j,k))_{j,k}$ avec $k\leq d_j$.

Il est clair que la Proposition \ref{csqce du chgt de base} d\'ecoule du lemme suivant.

\begin{lem} \label{chgt de base + LRS}
Supposons que $I(\sigma , s)$ admette un sous-quotient irr\'eductible $\pi (\sigma, s)$ appartenant au dual automorphe de $U(p,q)$. 
Alors, en admettant la Conjecture de changement de base, pour chaque couple d'entiers $(j,k)$ avec $k\leq d_j$, 
ou bien
\begin{enumerate}
\item $s(j,k)$ est un entier, ou bien
\item c'est un nombre complexe qui v\'erifie l'in\'egalit\'e 
$$| \mbox{Re} (s(j,k))|  \leq 1 - \frac{2}{(p+q)^2 +1} .$$
\end{enumerate}
\end{lem}
{\it D\'emonstration du Lemme \ref{chgt de base + LRS}.} Posons $G=U(p,q)$. On a d\'ecrit le $L$-groupe de $G$ 
au Chapitre 4. Soit 
$$\varphi : W_{{\Bbb R}} \rightarrow ^L G$$
un morphisme admissible et supposons que $\pi (\sigma ,s )$ appartiennent au $L$-paquet de $\varphi$. On peut 
supposer que l'image de $\varphi$ est contenue dans le $L$-groupe du sous-groupe de Levi $MA$ et que la restriction
de $\varphi$ au sous-groupe ${\Bbb C}^* \subset W_{{\Bbb R}}$ a son image contenue dans le $L$-groupe du sous-groupe 
de Cartan diagonal $T \subset MA$. Comme au Chapitre 5, on identifie $^L T^+$ avec $\prod_j ({\Bbb C}^* )^{d_j} \times
\prod_j ({\Bbb C}^* )^{c_j}$. Alors la restriction de $\varphi$ \`a ${\Bbb C}^* \subset W_{{\Bbb R}}$ doit \^etre de la forme
\begin{eqnarray}
\varphi (z) = \left( 
\begin{array}{ccc}
(z\overline{z} )^{\frac{s(j,k)}{2}} \left( \frac{z}{\overline{z}} \right)^{\frac{n(j,k)}{2}} & & \\
& \left( \frac{z}{\overline{z}} \right)^{\frac{m(j,l)}{2}} & \\
& & (z\overline{z} )^{\frac{-s(j,k)}{2}} \left( \frac{z}{\overline{z}} \right)^{\frac{n(j,k)}{2}} 
\end{array}
\right) ,
\end{eqnarray}
o\`u $k\leq d_j$, $l\leq c_j$ comme au Chapitre 5 et $n(j,k)$, $m(j,l) \in {\Bbb Z}$.

Comme dans \cite{Borel}, nous notons
$$\varphi (z) = z^{\lambda} \overline{z}^{\mu}  \; \; (z\in {\Bbb C}^* )$$
o\`u les poids $\lambda , \mu \in X^* (T) \otimes {\Bbb C}$ sont donn\'es (dans les coordonn\'ees du Chapitre 5) par 
les formules
\begin{eqnarray} \label{lamu}
\left\{ 
\begin{array}{ccl}
\lambda & = & \frac12 ( s(j,k) +n(j,k) , -s(j,k) +n(j,k) , m(j,l) ) \\
\mu & = & \frac12 ( s(j,k) - n(j,k) , -s(j,k) - n(j,k) , - m(j,l) ) 
\end{array} \right.
\end{eqnarray}
o\`u $\lambda, \mu \in \prod_j \prod_{k\leq d_j} ({\Bbb C}^* ) \times \prod_{l \leq c_j} ({\Bbb C}^* )$.
Soit $\varphi ' : W_{\Bbb R} \rightarrow GL(p+q, {\Bbb C}) \times GL(p+q , {\Bbb C})$ le morphisme obtenu par changement
de base. La restriction de $\varphi '$ \`a ${\Bbb C}$ est donn\'ee par $\varphi ' (z) = (\varphi (z) , \varphi (\overline{z}) )$.
D'apr\`es le Chapitre3, le caract\`ere infinit\'esimal de la repr\'esentation irr\'eductible de $GL(p+q , {\Bbb C})$ associ\'e
au morphisme $\varphi '$ est la paire $(\lambda , \mu )$ donn\'ee par (\ref{lamu}), \`a permutations des coordonn\'ees
pr\`es.

Supposons maintenant que pour un certain couple $(j,k)$, $s(j,k)$ ne soit pas un entier. Notons $\pi'$ la repr\'esentation
irr\'eductible de $GL(p+q , {\Bbb C} )$ associ\'ee au morphisme $\varphi '$. D'apr\`es la Conjecture de changement de base,
le caract\`ere infinit\'esimal de $\pi'$ est \'egal au caract\`ere infinit\'esimal d'une repr\'esentation de la forme
$${\rm ind}_P^{GL(p+q , {\Bbb C})} (J( \tau_1 , b_1 )\otimes \ldots \otimes J(\tau_m , b_m )).$$
\`A l'aide de la description au \S 3.4.4 du caract\`ere infinit\'esimal de cette derni\`ere, on obtient imm\'ediatement 
que $s(j,k)$ est \'egal \`a un certain param\`etre $t_j$ d'une certaine repr\'esentation automorphe cuspidale 
$\tau_i$ de $GL(a_i , {\Bbb C})$. Le Th\'eor\`eme \ref{LRS} implique donc que
$$|\mbox{Re} (s(j,k) )| \leq 1- \frac{2}{a_i^2 +1} \leq 1 - \frac{2}{(p+q)^2 +1} .$$
Ce qui ach\`eve la d\'emonstration du Lemme \ref{chgt de base + LRS}.

\medskip

Remarquons que la conclusion de la Proposition 8.4.2 d\'ecoule \'egalement des Conjectures d'Arthur d'apr\`es le chapitre 6, d'apr\`es ce m\^eme chapitre la 
m\^eme conclusion devrait \'egalement \^etre vraie pour le groupe $O(p,q)$ sauf pour $pq$ impair et $i=\frac{pq \pm 1}{2}$.

\newpage

\thispagestyle{empty}

\newpage

\markboth{CHAPITRE 9. D\'EMONSTRATION DU TH\'EOR\`EME 2}{9.1.FONCTIONS RADIALES}

\chapter{D\'emonstration du Th\'eor\`eme~2}

Le but de ce chapitre est la d\'emonstration du Th\'eor\`eme 2. Celle-ci repose sur l'\'etude de la d\'ecroissance \`a l'infini des fonctions 
sph\'eriques et sur le Th\'eor\`eme de Burger et Sarnak d\'ej\`a rencontr\'e au chapitre 2.

Nous commen\c{c}ons ce chapitre par des g\'en\'eralit\'es sur les fonctions sph\'eriques. Puis, et bien que notre m\'ethode fonctionne pour tous les
groupes de rang 1, nous nous sp\'ecialisons au groupe $U(n,1)$ dont l'on d\'ecrit pr\'ecis\'ement les fonctions sph\'eriques et leur comportement \`a 
l'infini.

On peut alors d\'emontrer le Th\'eor\`eme 2 dans le cas du groupe $U(n,1)$, laissant au lecteur le soin de v\'erifier que la d\'emonstration est identique 
dans le cas du groupe $O(n,1)$.

Fixons pour l'instant un groupe $G$ alg\'ebrique simple et connexe sur ${\Bbb Q}$; pour simplifier nous noterons \'egalement $G=G^{{\rm nc}}$ la 
partie non compacte (semi-simple) des points r\'eels de $G$. Nous noterons alors $K$ un sous-groupe compact maximal de $G$ et 
$X=G/K$ l'espace sym\'etrique associ\'e.

\section{Fonctions radiales, coefficients matriciels et fonctions sph\'eriques}

Soit $(\tau , V_{\tau})$ une repr\'esentation irr\'eductible de $K$. On appelle {\it fonction $\tau$-radiale} \index{fonction radiale}
toute fonction 
$$F : G \rightarrow {\rm End}(V_{\tau} ) $$
v\'erifiant la condition de double $K$-\'equivariance suivante :
\begin{eqnarray}
F(k_1 g k_2 ) = \tau (k_1 ) F(g) \tau (k_2 ) 
\end{eqnarray}
pour tout $g \in G$ et $k_1 , k_2 \in K$.

Soit $\pi$ une repr\'esentation continue de $G$ dans un espace de Hilbert $({\cal H},\langle.,.\rangle )$. On appelle
{\it coefficient matriciel} \index{coefficient matriciel} de $\pi$ toute fonction $G \rightarrow \mathbb{C}$ de la forme 
$$c_{v,v'} : g \mapsto \langle \pi (g) v, v' \rangle , $$
o\`u $v,v'\in {\cal H}$. Le lemme suivant est 
facile \`a v\'erifier.

\begin{lem} \label{coeff}
Soit $v' \in {\cal H}$.
\begin{enumerate}
\item L'application $v \mapsto c_{v,v'}$ est \'equivariante par rapport \`a la $\pi$-action de $G$ sur
${\cal H}$ et \`a l'action r\'eguli\`ere \`a droite $R$ de $G$ sur les fonctions de $G$ dans $\mathbb{C}$.
\item Si $v$ est dans le sous-espace ${\cal H}^{\infty}$ des vecteurs $C^{\infty}$, alors $c_{v,v'} \in C^{\infty} (G)$.
\item L'application $v \mapsto c_{v,v'}$ de ${\cal H}^{\infty} \rightarrow C^{\infty} (G)$ est \'equivariante par rapport
aux actions de $U(\LG )$, l'alg\`ebre enveloppante universelle sur ${\Bbb C}$ de l'alg\`ebre de Lie de $G$,  induites 
par $\pi$ et $R$.
\end{enumerate}
\end{lem}

Si la repr\'esentation $\pi$ a un caract\`ere infinit\'esimal, alors il d\'ecoule du Lemme ci-dessus que tout
coefficient matriciel $c_{v,v'}$, avec $v$ vecteur lisse, est une fonction dans $C^{\infty} (G)$ fonction propre
pour l'op\'erateur $R(C)$ (o\`u l'on note toujours $C$ le casimir de $G$).

\medskip

Dor\'enavant soit $\pi$ une repr\'esentation admissible de $G$.
Si $v$ et $v'$ sont deux vecteurs $K$-finis de l'espace de la representation $\pi$, le coefficient matriciel 
$c_{v,v'}$ est analytique (r\'eel) et se transforme finiment sous les actions \`a droite et \`a gauche de $K$.
 
On associe \`a $\pi$, l'id\'eal $I_{\pi}$ du centre ${\cal Z} (\mathfrak{g} )$ de l'alg\`ebre enveloppante, d\'efini par
$$I_{\pi} = \{ Z \in {\cal Z} (\mathfrak{g} ) \; : \; \pi (Z) =0 \} .$$
Il d\'ecoule alors du Lemme \ref{coeff} que 
\begin{eqnarray} \label{eqns}
R(Z) c_{v,v'} = 0 \; \; (Z \in I).
\end{eqnarray}
Rappelons, cf. \cite{Knapp} par exemple, le r\'esultat classique que chaque id\'eal $I_{\pi}$ comme ci-dessus
est cofini dans l'alg\`ebre ${\cal Z} (\mathfrak{g} )$. Enfin remarquons que si $\pi$ a un caract\`ere infinit\'esimal, 
l'id\'eal $I_{\pi}$ est de codimension $1$. Dans ce cas, (\ref{eqns}) est un syst\`eme d'\'equations propres. 

\medskip

Supposons que $\tau$ est un $K$-type de la repr\'esentation $\pi$.
Soient l'inclusion
$$i_{\tau} : V_{\tau }  \hookrightarrow {\cal H} $$
et la projection
$$p_{\tau} : {\cal H} \rightarrow V_{\tau } .$$
Alors, la fonction $F:G \rightarrow {\rm  End}(V_{\tau } )$ d\'efinie par
\begin{eqnarray} \label{fonction radiale}
F(g) = p_{\tau} \circ \pi (g) \circ i_{\tau}
\end{eqnarray}
est $\tau$-radiale.
De plus, si $v$ et $v'$ sont deux vecteurs dans $V_{\tau }$, le 
coefficient $c_{v,v'}$ de $\pi$ s'\'ecrit :
$$c_{v,v'} (g) = \langle F(g) v , v' \rangle .$$

\medskip

On dira d'une fonction 
$\Phi$ $\tau$-radiale et de classe $C^{\infty}$, qu'elle est {\it $\tau$-sph\'erique sur $G$} \index{fonction sph\'erique} si 
\begin{enumerate}
\item $\Phi (e) =Id$, et si 
\item $\Phi$ v\'erifie le syst\`eme d'\'equations propres
$$R(Z) \Phi =0,$$
o\`u $Z$ d\'ecrit un id\'eal $I$ de codimension $1$ dans ${\cal Z} (\mathfrak{g} )$. 
\end{enumerate} 
Remarquons que toute fonction $\tau$-sph\'erique est en fait analytique (r\'eelle).
On notera ${\cal A} (G, \tau )$ l'espace des fonctions $\tau$-sph\'eriques, et ${\cal A} (G, \tau , I)$ le 
sous-ensemble de celle v\'erifiant le 2. pour un id\'eal $I$ sp\'ecifi\'e.

Remarquons que si $\tau$ est un $K$-type d'une repr\'esentation admissible $\pi$ de $G$ admettant un caract\`ere 
infinit\'esimal, la fonction radiale $F$ du (\ref{fonction radiale}) est en fait $\tau$-sph\'erique et appartient plus 
pr\'ecis\'ement \`a ${\cal A} (G, \tau , I_{\pi})$.

Les fonctions sph\'eriques v\'erifient un syst\`eme d'\'equations diff\'erentielles d'apr\`es le deuxi\`eme point de la
d\'efinition ci-dessus. Ce syst\`eme d'\'equation diff\'erentielles a \'et\'e \'etudi\'e et r\'esolu par 
Harish-Chandra, cf. \cite{CasselmanMilicic}, \cite{Knapp}.

Nous nous int\'eressons en fait \`a l'asymptotique des fonctions sph\'eriques $\Phi \in {\cal A} (G, \tau )$. D'apr\`es la d\'ecomposition de Cartan
\begin{eqnarray} \label{dec cartan}
G= K \overline{A^+} K
\end{eqnarray}
il suffit de comprendre $\Phi (a)$ lorsque $a$ tend vers l'infini dans $\overline{A^+}$. Ici $\overline{A^+}$ d\'esigne 
la cl\^oture de $A^+ = \exp (\mathfrak{a}^+_0 )$, o\`u $\mathfrak{a}^+_0$ est une chambre de Weyl positive (ouverte)
dans $\mathfrak{a}_0$.

Remarquons que la restriction d'une fonction $\tau$-sph\'erique \`a $A^+$ prend ses valeurs dans 
$$E^M := \{ L \in {\rm End} (V_{\tau} ) \; : \; L= \tau (m) L \tau (m)^{-1} \mbox{  pour tout } m\in M \};$$
o\`u $M$ d\'esigne le centralisateur de $A$ dans $K$.

\markboth{CHAPITRE 9. D\'EMONSTRATION DU TH\'EOR\`EME 2}{9.2.FONCTIONS SPH\'ERIQUES DU GROUPE $U(n,1)$}

\section{Fonctions sph\'eriques du groupe $U(n,1)$}

Puisque la d\'emonstration du Th\'eor\`eme 2 est similaire (en fait plus simple) dans le cas hyperbolique r\'eel nous nous contenterons de traiter le 
cas hyperbolique complexe. Dor\'enavant nous supposons donc que
$G=U(n,1)$ et $K=U(n+1) \cap G = U(n) \times U(1)$.

Dans ce cas, la paire $(G,K)$ est une paire de Gelfand et 
l'alg\`ebre $C_c (G,K, \tau ,\tau )$ des fonctions continues $\tau$-radiales de support compact sur $G$ pour 
le produit de convolution :
$$(F*G )(x) = \int_{G} F(g^{-1} x ) H(g) dg = \int_{G} F(g) H(x g^{-1} ) $$
est une alg\`ebre commutative.
Cette propri\'et\'e va nous permettre de caract\'eriser plus facilement les fonctions sph\'eriques et de les 
classifier. Mais d'abord commen\c{c}ons par d\'ecrire les diff\'erents $K$-types possibles.

\medskip

Puisque tout \'el\'ement de $K$ peut s'\'ecrire comme :
$$k = \left( 
\begin{array}{cc}
U & 0 \\
0 & v 
\end{array} \right) , \; \; U \in U(n) , \; v\in U(1) , $$
$K$ agit sur $\LM \cong \mathbb{C}^n$  par ${\rm Ad}(k) X = v^{-1} UX$ et cette action pr\'eserve la structure
complexe $J$.  On a donc la d\'ecomposition en $K$-modules suivante~:
$$\LM_{\mathbb{C}} = \LM_+ \oplus \LM_- \mbox{  et  } \LM_{\pm}^* \cong \overline{\LM_{\pm}} \cong \LM_{\mp} .$$
Comme au chapitre 4, on note :
\begin{eqnarray}
\tau_{p,q} := \Lambda^p {\rm Ad}_+^* \otimes \Lambda^q {\rm Ad}_-^* \cong \Lambda^p \overline{{\rm Ad}} \otimes \Lambda^q {\rm Ad} .
\end{eqnarray}
C'est une repr\'esentation de $K$ r\'eductible qui se d\'ecompose en irr\'eductible de la fa\c{c}on suivante :
\begin{eqnarray}
\tau_{p,q} = \oplus_{k=0}^{\mbox{min}(p,q)}  \tau_{p-k , q-k} ' .
\end{eqnarray}
Remarquons que via la formule de Matsushima, la repr\'esentation $\tau_{p,q}$ 
correspond aux formes de type $(p,q)$  et la repr\'esentation $\tau_{p,q} '$ correspond aux 
formes {\it primitives} \index{primitives, formes} de type $(p,q)$, {\it i.e.} les formes de type $(p,q)$ qui ne peuvent s'\'ecrire comme
un multiple de la forme de Kaehler.

\medskip

Dor\'enavant soit $\tau$ une repr\'esentation de $K$ appartenant \`a l'ensemble des $\tau_{p,q} '$.
Soit $(\pi , H)$ une repr\'esentation irr\'eductible de $G$ telle que 
\begin{enumerate}
\item $\pi \subset L^2 (\Gamma \backslash G)$;
\item $\tau  \subset \pi_{|K}$ et 
\item $\pi (C) = -\lambda$.
\end{enumerate}

\medskip

Il d\'ecoule de la formule de r\'eciprocit\'e de Frobenius le fait suivant.

\medskip
\noindent
{\bf Fait.} La multiplicit\'e de $\tau$ dans $\pi_{|K}$ est \'egale \`a $1$.

\medskip

Notons $V(\tau )$ le sous-espace de la repr\'esentation $\tau$ dans $V_{|K}$.
On a $V(\tau ) \cong \mathbb{C}^n$ et on introduit 
$$E_{\tau} = \mbox{End}(V(\tau )).$$

\medskip

Une fonction 
$\Phi$ $\tau$-radiale de classe $C^{\infty}$ est $\tau$-sph\'erique sur $G$ si et seulement si \index{fonction sph\'erique}
\begin{enumerate}
\item $\Phi (e) =Id$, et 
\item $\Phi$ est une fonction propre pour la convolution par $C_c (G,K,\tau ,\tau )$, {\it i.e.} il existe un (unique)
caract\`ere $\lambda_{\Phi}$ de $C_c (G,K,\tau ,\tau )$ tel que, pour toute $F\in C_c (G,K,\tau ,\tau )$, 
$F*\Phi = \Phi *F =\lambda_{\Phi} (F) \Phi .$
\end{enumerate} 
Nous noterons $\Sigma (G,K,\tau ,\tau ) $ l'ensemble des fonctions $\tau$-sph\'eriques sur $G$.

\bigskip

\`A l'aide de la d\'ecomposition $G=KAN$ on d\'efinit la fonction $h$ sur $G$ par :
$$h(ka_t n) = t$$
et par $\underline{k}$ le projecteur sur $K$. Alors la repr\'esentation $\pi_{\sigma , s}$
agit sur $({\cal H}_{\sigma ,s} )_{|K} \cong L^2 (K,M,\sigma )$ par 
$$\pi_{\sigma ,s} (x) f(k) = e^{-(s+\rho )h(x^{-1} k)} f(\underline{k} (x^{-1} k)) $$
pour tout $x\in G$ et $k \in K$.

Soit $\sigma$ une repr\'esentation irr\'eductible de $M$ apparaissant dans $\tau_{|K}$ ce que nous notons~:
$\sigma \in \widehat{M} (\tau )$. D\'esignons par $P_{\sigma}^{\tau}$ le g\'en\'erateur de l'espace de dimension $1$ 
${\rm Hom}_K ({\cal H}_{\sigma , s} , V_{\tau } ) (\cong {\rm Hom}_K (L^2 (K,M, \sigma ) , V_{\tau })$ pour tout 
$s \in \mathbb{C}$) d\'efini par~:
$$P_{\sigma}^{\tau} (f) := \sqrt{\frac{{\rm dim} \tau}{{\rm dim} \sigma}} \int_K \tau (k) f(k) dk .$$
Posons $J_{\sigma}^{\tau} = (P_{\sigma}^{\tau} )^*$, {\it i.e.} $J_{\sigma}^{\tau}$ est le g\'en\'erateur de 
l'espace de dimension 1 ${\rm Hom}_K (V_{\tau} , {\cal H}_{\sigma ,s} ) ( \cong {\rm Hom}_K (V_{\tau} , L^2 (K,M,\sigma )))$
d\'efinit par 
$$J_{\sigma}^{\tau} \xi = \sqrt{\frac{{\rm dim} \tau}{{\rm dim} \sigma}} P_{\sigma} \circ \{ \tau (.) ^{-1} \xi \} ,$$
o\`u $P_{\sigma}$ d\'esigne la projection orthogonale de $V_{\tau}$ sur le sous-espace de $\sigma$ dans $V_{\tau}$.
Pour $\sigma \in \widehat{M} (\tau )$ et $s \in \mathbb{C}$, on a vu que l'application
$$\Phi_{\sigma ,s}^{\tau} : g \mapsto P_{\sigma}^{\tau} \circ \pi_{\sigma , s} (g) \circ J_{\sigma}^{\tau}$$
d\'efinit une fonction $\tau$-radiale sur $G$.

\begin{prop}[Classification des fonctions sph\'eriques] \index{classification des fonctions sph\'eriques}
Soit $\tau = \tau_{p,q} '$.
Pour tout $\sigma \in \widehat{M} (\tau )$ et $s \in \mathbb{C}$, la fonction $\Phi_{\sigma ,s}^{\tau}$ 
est $\tau$-sph\'erique. Elle admet la repr\'esentation suivante :
$$\Phi_{\sigma , s}^{\tau} (x)  = \frac{{\rm dim} \tau}{{\rm dim} \sigma} \int_K e^{-(s+\rho ) h(x^{-1} k)} 
\tau (k) \circ P_{\sigma} \circ \tau ( \underline{k} (x^{-1} k)^{-1} ).$$
En particulier, $\Phi_{\sigma , s}^{\tau}$ est holomorphe en la variable $s$.
Enfin, 
$$\Sigma (G,K,\tau ,\tau ) = \{ \Phi_{\sigma , s}^{\tau} \; : \; \sigma \in \widehat{M} (\tau ) , \; \lambda \in \mathbb{C} / \{ \pm 1 \} \}.$$
\end{prop}

D'apr\`es la d\'ecomposition de Cartan $G=KA^+ K$ et l'\'equation (9.1.1), 
une application $\tau$-radiale $F$ est compl\`etement d\'etermin\'ee par sa restriction \`a $A^+$.
De plus $F_{|A^+}$ a son image contenue dans
$$E_{\tau}^M := \{ T \in E_{\tau} \; : \; T=\tau (m) T \tau (m)^{-1} \mbox{  pour tout } m\in M \} ;$$
o\`u $M$ d\'esigne toujours le centralisateur de $A$ dans $K$.
Dans la suite, nous \'etudions les fonctions sph\'eriques le long de $A$.
Rappelons donc la d\'ecomposition en irr\'eductibles de la restriction de $\tau =\tau_{p,q} '$ \`a $M$ :
\begin{eqnarray}
(\tau_{p,q} ')_{|M} = \sigma_{p,q}  \oplus \sigma_{p-1,q}  \oplus \sigma_{p,q-1}  \oplus \sigma_{p-1,q-1}  
\end{eqnarray}
(avec les conventions du \pa 4.1.1 : certains facteurs disparaissent dans les cas ``non g\'en\'eriques'').

\medskip

\section{L'\'equation diff\'erentielle radiale}

Soit $\sigma = \sigma_{p,q} $ avec $0\leq p+q \leq n-1$. Soit $\tau \in \widehat{K}$ telle que $\sigma \in \widehat{M} (\tau )$
(g\'en\'eriquement $\tau$ est l'une des $\tau_{p,q} '$, $\tau_{p+1,q} '$, $\tau_{p,q+1} '$, $\tau_{p+1,q+1} '$).
Posons 
$$F (g)= \frac{{\rm dim} \tau}{{\rm dim} \sigma} \int_K e^{-(\overline{s} +n) h(g k)} \tau (\underline{k} (g k))
\circ P_{\sigma} \circ \tau (k^{-1} ) dk , $$
l'adjoint de $\Phi_{\sigma ,s}^{\tau} (g^{-1} )$ pour le produit scalaire de $V_{\tau}$. 
Dans cette section nous calculons la restriction de $R(C)F$ (o\`u $C$ d\'esigne toujours le casimir)
\`a $A^+$ en fonction de $F_{|A^+}$.
Nous allons suivre \cite[pp. 279-282]{Wallach} (c'est pour cette raison que nous sommes pass\'es \`a l'adjoint). 

Posons d'abord :
$$f_k (g ) = e^{-(\overline{s} +n )h(gk)} \tau (\underline{k} (gk)) \circ P_{\sigma} \circ \tau (k)^{-1} .$$
Alors $F(g) = \frac{{\rm dim} \tau}{{\rm dim} \sigma} \int_K f_k (g ) dk$ et 
$R(C)F (g) = \frac{{\rm dim} \tau}{{\rm dim} \sigma} \int_K ( R(C) f_k ) (g )$.  
D'un autre c\^ot\'e, $f_k (g) = f_e (gk) \tau (k)^{-1}$. Puisque $(R(C)f_e )(gk) = (R({\rm Ad}(k) C) f_k (g)) \tau (k)$, on obtient que
$R(C)f_k (g) = (R(C)f_e )(gk) \tau (k)^{-1}$. Il suffit donc de calculer $R(C)f_e$. Soit $f=f_e$.

\medskip

Notons $\theta$ l'involution de Cartan correspondant \`a la d\'ecomposition $\mathfrak{g} = \mathfrak{p} \oplus 
\mathfrak{k}$ de l'alg\`ebre de Lie $\mathfrak{g}$ et $B$ la forme de Killing de $\mathfrak{g}$, cf. (4.0.1).
Soit $\lln = \LG_{\alpha} + \LG_{2\alpha}$ la d\'ecomposition en espaces de racines. 
La dimension de $\LG_{\alpha}$ est $2n-2$ et 
celle de $\LG_{2\alpha}$ est $1$. Soit donc $X_1 , \ldots , X_{2n-2}, X_{2n-1}$ une base de $\lln$ telle que 
$X_1 , \ldots , X_{2n-2}$ soient dans $\LG_{\alpha}$, $X_{2n-1}$ soit dans $\LG_{2\alpha}$ et $B(X_i , \theta X_j ) =-\delta_{ij}$.
Alors $[X_i , \theta X_i ]=-H_0$ pour $i=1, \ldots , 2n-2$ et $[X_{2n-1} , \theta X_{2n-1} ]= -2H_0$.
Rappelons que $B(H_0 , H_0 )=1$. Enfin soit $U_1, \ldots ,U_r$ une base de $\LM$ telle que $B(U_i , U_j )=-\delta_{ij}$.
Alors :
$$C= -\sum X_i \theta X_i - \sum \theta X_i X_i - \sum U_i^2 + H_0^2 .$$
Soit $C_0 = \sum U_i^2$. Remarquons que si $m$ est dans $M$, alors ${\rm Ad} (m) C_0 = C_0$. Par d\'efinition de $f$, 
$f(gn)=f(g)$ pour $n$ dans $N$. Puisque $X_i \theta X_i = -H_0 +\theta X_i X_i$ pour $1\leq i \leq 2n-2$
et $X_{2n-1} \theta X_{2n-1} = -2H_0 + \theta X_{2n-1} X_{2n-1}$, on obtient que 
$$R(C)f = 2n R(H_0 )f -R(C_0 )f +R( H_0^2 )f .$$
D'un autre c\^ot\'e,
\begin{eqnarray*}
(R(C_0)f ) (g) & = & \sum_{i=1}^{m} \frac{d^2}{dt^2} f(ge^{t U_i}) _{|t=0} \\
                  & = & \sum_{i=1}^m e^{-(\overline{s} +n) h(g)} \tau (\underline{k} (g)) \tau (U_i )^2 P_{\sigma} \\
                  & = & e^{-(\overline{s}+n ) h(g)} \tau (\underline{k} (g)) \tau (C_0 ) P_{\sigma} .
\end{eqnarray*}
Mais $\tau ( C_0 )$ agit par multiplication par un scalaire sur chaque sous-espace $M$-irr\'eductible de $V_{\tau}$.
Donc $\tau ( C_0 ) P_{\sigma} = \lambda_{\sigma} P_{\sigma}$. 

Puis,
\begin{eqnarray*}
R(H_0 )f (g) & = & \frac{d}{dt} f(g e^{tH_0} )_{|t=0} \\
               & = & \frac{d}{dt} (e^{-(\overline{s} +n) (h(g) +t)})_{|t=0} \tau (\underline{k} (g) )  P_{\sigma} \\
               & = & -(\overline{s} +n) f(g) .
\end{eqnarray*}
Et donc : $R(H_0^2 )f(g) = (\overline{s} +n )^2 f(g)$. On obtient finalement la formule :
$$R(C)f = -2n(\overline{s} +n) f + (\overline{s}+n)^2 f -\lambda_{\sigma} f.$$
Et donc :
\begin{eqnarray} \label{action de C}
R(C) F(g) = (\overline{s}^2 -n^2 -\lambda_{\sigma} ) F(g) .
\end{eqnarray}

Si $\varphi : G \rightarrow {\rm End} (V_{\tau} )$ est une fonction $\tau$-radiale de classe $C^{\infty}$, 
on peut montrer, \cite[formule (4) p.282]{Wallach}, que :
\begin{eqnarray} \label{action radiale de C}
\begin{array}{cl}
(R(C)\varphi ) (a_t )& = \frac{d^2 \varphi (a_t) }{dt^2} + (2\coth (2t) +2(n-1) \coth (t) ) \frac{d\varphi (a_t  )}{dt} \\
& +Q_1 (t) \varphi (a_t ) -\varphi (a_t ) \tau ( C_0 ) ,
\end{array}
\end{eqnarray}
o\`u de l'expression de $Q_1$ donn\'e dans \cite{Wallach} nous ne retiendrons que la d\'ecroissance 
exponentielle (pour $t\rightarrow + \infty$) vers $0$.

\begin{prop}[\'Equation diff\'erentielle radiale]
 \index{\'equation diff\'erentielle radiale}
Pour simplifier notons $F(t)=F(a_t)$. Alors pour $t>0$, 
\begin{eqnarray} \label{equa diff radiale}
F''(t)+(2\coth (2t) + 2(n-1) \coth (t) ) F'(t) +Q_1 (t) F(t) = (\overline{s}^2 -n^2 ) F(t) .
\end{eqnarray}
\end{prop}
{\it D\'emonstration.} La Proposition d\'ecoule des \'equations (\ref{action de C}) et (\ref{action radiale de C}) en remarquant que 
$$F(a_t) \tau (C_0 )  = \lambda_{\sigma} F(a_t ) .$$

\medskip

\markboth{CHAPITRE 9. D\'EMONSTRATION DU TH\'EOR\`EME 2}{9.4. COMPORTEMENT ASYMPTOTIQUE}

\section{Comportement asymptotique}

Nous conservons les notations de la section pr\'ec\'edente. Notamment nous 
notons toujours $F(t)$ la fonction de la Proposition 9.3.1.
Puisque la fonction $Q_1 (t)$ tend exponentiellement vite vers $0$ lorsque 
$t$ tend vers l'infini, l'\'equation diff\'erentielle (\ref{equa diff radiale}) est exponentiellement asymptote 
(lorsque $t$ tend vers $0$) \`a l'\'equation diff\'erentielle :
\begin{eqnarray} \label{equa modele}
F''(t) + 2n F'(t) +(n^2 - \overline{s}^2 ) F(t) =0.
\end{eqnarray}
Un th\'eor\`eme classique nous assure alors que les solutions des \'equations diff\'erentielles
(\ref{equa diff radiale}) et (\ref{equa modele}) sont asymptotes lorsque $t$ tend vers l'infini.

Plus pr\'ecis\'ement on obtient la proposition suivante \footnote{Ici comme dans d'autres formules, on esp\`ere que le lecteur saura distinguer $n(\in {\Bbb N})$ et 
$n$ (\'el\'ement unipotent).}:
\begin{prop}[Asymptotique des fonctions sph\'eriques] \label{assymptot fonct spherique}
 \index{asymptotique des fonctions sph\'eriques}
Soit $F(t)$ l'adjoint de $\Phi_{\sigma ,s}^{\tau} (a_t^{-1} )$ agissant sur l'espace de Hilbert $V_{\tau}$. 
Si $s$ est dans $\mathbb{C}$ avec $n>$Re$(s) \geq \varepsilon >0 $, alors il existe un r\'eel $\delta >0$
ind\'ependant de $s$ et de $\varepsilon$  
et une constante $C_{\varepsilon}$ ne d\'ependant que de $\varepsilon$ tels que :
$$||e^{nt} F(t) - e^{\overline{s}t} P_{\sigma} B_{\tau} (\overline{s}) 
-e^{-\overline{s} t} \tau (k^* ) B_{\tau} (-\overline{s} )^* P_{\sigma} \tau ( k^* )^{-1} || \leq C_{\varepsilon} 
e^{-\delta  t} , $$
o\`u $k^*$ est un \'el\'ement de $K$ centralisant $A$ et tel que $k^* a_t (k^*)^{-1} =a_{-t}$ et 
$$B_{\tau} (z)= \int_{\overline{N}} e^{-(n+z)h(\overline{n})} \tau ( \underline{k} (\overline{n} ))^{-1} d\overline{n} .$$
(Ici $\overline{N}$ est le radical unipotent du parabolique oppos\'e \`a $P$.)
\end{prop}
{\it D\'emonstration.} Soit $\Delta (t) = ( \sinh t)^{2(n-1)} (\sinh 2t)$. Posons $G(t) = \Delta (t)^{\frac12} F(t)$ pour $t>0$.
L'\'equation (9.3.3) se r\'e\'ecrit :
\begin{eqnarray}
-G''(t) +Q(t)G(t)=-\overline{s}^2 G(t) ,
\end{eqnarray}
o\`u $Q(t)$ est une fonction qui d\'ecro\^{\i}t exponentiellement vite vers $0$ lorsque $t$ tend vers l'infini.
Plus pr\'ecis\'ement (cf. \cite{Wallach}) si $T$ est un r\'eel strictement positif, il existe une constante 
$C_1$ ne d\'ependant que de $T$ telle que $||Q(t)|| \leq C_1 e^{-t}$ ($t\geq T$). Donc, \'etant donn\'e un r\'eel 
$\delta$ strictement compris entre $0$ et $1$, il existe une constante $C_2$ telle que pour tout $u>T$ :
\begin{eqnarray}
\int_u^{+\infty} ||Q(t)|| (1+t)dt <C_2 e^{-\delta u} .
\end{eqnarray}
Les Lemmes A.8.2.12 et A.8.2.17 de l'appendice de \cite{Wallach} impliquent alors que si Re$(s) \geq \varepsilon >0$, 
il existe une constante $C_{\varepsilon}$ ne d\'ependant que de $\varepsilon$ et deux endomorphismes $A_1$ et $A_2$ de 
$V_{\tau}$ tels que :
\begin{eqnarray}
||G(t) - e^{\overline{s}t} A_1 -e^{-\overline{s}t} A_2|| < C_{\varepsilon} e^{-\delta t} .
\end{eqnarray}
Puisque $\Delta (t)^{\frac12}$ est exponentiellement asymptotique \`a $e^{nt}$, pour conclure la d\'emonstration 
il reste \`a d\'eterminer $A_1$ et $A_2$.
Encore une fois nous suivons \cite[pp.283-285]{Wallach} pour cela. Puisque Re$(s)>0$, d'apr\`es (9.4.4), on a :
$$\lim_{t \rightarrow +\infty} || e^{nt}F(t) - e^{\overline{s}t} A_1 || =0 .$$
Or, 
$$F(t) = \frac{{\rm dim} \tau}{{\rm dim} \sigma} \int_K e^{-(\overline{s} +n) h(a_t k)} \tau (\underline{k} (a_t k))
\circ P_{\sigma} \circ \tau (k^{-1} ) dk . $$
Les calculs de \cite[p.284]{Wallach} (et en particulier l'expression (14)), montrent alors que si $n>$Re$(s)>0$, 
$$\lim_{t\rightarrow +\infty} e^{(n-\overline{s})t} F(t) = P_{\sigma} B_{\tau} (\overline{s}) .$$
Le calcul de $A_2$ se fait de m\^eme (cf. \cite{Wallach}) et ach\`eve la d\'emonstration de la Proposition.

\medskip

\noindent
{\bf Remarques.}  Les fonctions $B_{\tau}$ ont \'et\'e \'etudi\'ees par Schiffmann \cite{Schiffmann} (cf. aussi \cite{Wallach}). Elles
sont explicites. Les seules propri\'et\'es dont nous aurons besoin sont~:
\begin{enumerate}
\item pour $n>$Re$(s) >0$, $B_{\tau} (s) \neq 0$ et,
\item les fonctions $B_{\tau}$ sont continues (en fait analytiques).
\end{enumerate}

La Proposition 9.4.1 n'est en fait qu'un cas particulier d'un r\'esultat plus g\'en\'eral de Langlands, cf. \cite{Langlands, CasselmanMilicic}.

\begin{cor}
Soient $\tau = \tau_{p,q}'$ et $v$ un vecteur dans $V_{\tau}$. 
Soit $\varepsilon$ un r\'eel strictement positif.
Il existe une constante $C_{\varepsilon}$ ne d\'ependant que de $\varepsilon$ et de $\tau$ et un r\'eel $\delta >0$ ind\'ependant  
de $\varepsilon$ tels que :
pour tout $\sigma \in \widehat{M} (\tau )$, 
\begin{enumerate}
\item si $s$ est un r\'eel dans $[\varepsilon , n]$, 
$$\langle \Phi_{\sigma , s}^{\tau} (e^{-tH} ) v,v \rangle_{V_{\tau}}  = B_{\tau} (s) ||P_{\sigma} (v)||^2_{V_{\tau}} e^{(-n+s)t} + o_{\varepsilon} (||v||_{V_{\tau}}^2 e^{(-n- {\rm min} (\varepsilon ,\delta))t} );$$
\item si $s$ est un r\'eel dans $]0, \varepsilon ]$,
$$\langle \Phi_{\sigma ,s}^{\tau} (e^{-tH} ) v,v \rangle_{V_{\tau}} \leq C_{\varepsilon} ||v||^2_{V_{\tau}} e^{(-n+\varepsilon )t} .$$
O\`u la notation $o_{\varepsilon}$ signifie que les constantes dans le $o$ ne d\'ependent que de $\varepsilon$.
\end{enumerate}
\end{cor}
{\it D\'emonstration.} Cela d\'ecoule de la Proposition ci-dessus pour le point 1 et de la d\'emonstration 
de la Proposition ci-dessus pour le point 2, en remarquant que $F(t) = F(a_t )$ est \'egal \`a l'adjoint
de $\Phi_{\sigma ,s}^{\tau } (e^{-tH} )$ pour le produit scalaire de $V_{\tau}$. 

\medskip

Rappelons que  
si $\tau = \tau_{p,q} '$, l'ensemble $\widehat{M} (\tau )$ contient, g\'en\'eriquement, quatre repr\'esentations.
Si $\sigma= \sigma_{a,b} $ est l'une de ces repr\'esentations et $s$ un nombre complexe, on a :
\begin{eqnarray}
\pi_{\sigma ,s} (C) = -((n-(a+b))^2 -s^2 ) Id .
\end{eqnarray}

Soit $\mathbb{D} (G,K,\tau )$ l'{\it alg\`ebre des op\'erateurs diff\'erentiels
invariants \`a gauche sur l'espace $C^{\infty} (G,K,\tau )$} 
\index{$\mathbb{D} (G,K,\tau )$} 
\index{op\'erateurs diff\'erentiels, alg\`ebre des} des fonctions $C^{\infty}$ $f:G\rightarrow V_{\tau}$ telles que 
$$f(xk)=f(x) \tau(k) .$$
Alors l'alg\`ebre $\mathbb{D} (G,K,\tau )$ est engendr\'ee par les op\'erateurs $\partial \partial^*$, 
$\partial^* \partial$, $\overline{\partial} \overline{\partial}^*$ et $\overline{\partial}^* \overline{\partial}$
cf. \cite{Pedon}. Elle contient notamment le laplacien de Hodge-de Rham : $\Delta = (\partial +\overline{\partial})( \partial^* + 
\overline{\partial}^* ) + (\partial^* +\overline{\partial}^* )(\partial +\overline{\partial} )$.

Et l'identit\'e (9.4.5) implique que 
la fonction sph\'erique $\Phi_{\sigma ,s}^{\tau}$ associ\'ee au $K$-type $\tau$ de $\pi_{\sigma , s}$
v\'erifie :
$$\Delta \Phi_{\sigma ,s}^{\tau} = ((n-(a+b))^2  -s^2 )\Phi_{\sigma ,s}^{\tau} .$$

Notons 
$$\Phi_{a,b}^{p,q} (s,x) = \Phi_{\sigma_{a,b}  , s}^{\tau_{p,q} '} (x)  ,$$
pour $\sigma_{a,b}  \in \hat{M} (\tau_{p,q} ')$.
On peut v\'erifier la caract\'erisation suivante des fonctions sph\'eriques.
\begin{thm} \index{fonctions sph\'eriques}
Soit $\Phi$ une fonction $\tau_{p,q} '$-radiale normalis\'ee, {\it i.e.} $\Phi (e) =Id$.
\begin{enumerate}
\item $\Phi =\Phi_{p,q}^{p,q} (s,.)$ si et seulement si $\Delta \Phi = \{ (n-(p+q))^2 -s^2 \} \Phi$, $\partial^* \Phi =0$ et 
$\overline{\partial}^* \Phi =0$;
\item $\Phi =\Phi_{p-1,q}^{p,q} (s,.)$ si et seulement si $\Delta \Phi = \{ (n-(p+q-1))^2 -s^2 \} \Phi$, $\partial \Phi =0$ et 
$\overline{\partial}^* \Phi =0$;
\item $\Phi =\Phi_{p,q-1}^{p,q} (s,.)$ si et seulement si $\Delta \Phi = \{ (n-(p+q-1))^2 -s^2 \} \Phi$, $\partial^* \Phi =0$ et 
$\overline{\partial} \Phi =0$;
\item $\Phi =\Phi_{p-1,q-1}^{p,q} (s,.)$ si et seulement si $\Delta \Phi = \{ (n-(p+q-2))^2 -s^2 \} \Phi$, $\partial \Phi =0$ et 
$\overline{\partial} \Phi =0$.
\end{enumerate}
\end{thm}

\medskip

\markboth{CHAPITRE 9. D\'EMONSTRATION DU TH\'EOR\`EME 2}{9.5. D\'EMONSTRATION DU TH\'EOR\`EME 2}

\section{D\'emonstration du Th\'eor\`eme 2}

Soit $G$ un groupe alg\'ebrique simple et connexe sur $\mathbb{Q}$ tel
que $G^{{\rm nc}}$ soit isomorphe au groupe $U(n,1)$. Soit $H$ un ${\Bbb Q}$-sous-groupe de $G$ tel que 
$H^{{\rm nc}}$ soit isomorphe au groupe $U(k,1)$. Dans ce cas le groupe (r\'eel) $H$ est stable par une involution de Cartan 
de $G$. Quitte \`a conjuguer $H$ dans $G$, on peut donc supposer que $K_H \subset K_G$ et $A_H =A=A_G$. 

\medskip

Nous allons maintenant d\'emontrer le Th\'eor\`eme 2. Mais
commen\c{c}ons par introduire  les d\'efinitions suivantes.

Soient $\varepsilon$ un r\'eel positif et $(p,q)$ un couple d'entiers.
Notons $H_{\varepsilon}^{(p,q)}$ l'hypoth\`ese sur le groupe $G$ suivante~: \index{$H_{\varepsilon}^{(p,q)}$}
$$H_{\varepsilon}^{(p,q)} : \left\{
\begin{array}{l}
\mbox{si $\lambda$ est dans le $(p,q)$-spectre automorphe de $G$, alors soit} \\
\mbox{1. }
\lambda = (n-(p+q)+k)^2 - (n-(p+q)+k-i)^2 \\
\mbox{  avec } k=0 , \ldots , p+q \mbox{ et } i= 0, \ldots , n-(p+q) +k , \\
\mbox{2. }  \lambda \geq (n-(p+q))^2 - \varepsilon^2 .
\end{array}
\right. $$
Les Conjectures d'Arthur pr\'evoient (cf. Chapitre 6) que l'hypoth\`ese $H_0^{(p,q)}$ est vraie pour 
tout couple d'entiers $(p,q)$ tel que $p+q \leq n-1$ (elles pr\'evoient m\^eme un peu plus, cf. Chapitre 6). Les hypoth\`eses $H_{\varepsilon}^{(p,q)}$ sont donc 
des approximations aux Conjectures d'Arthur pour les groupes unitaires.

Le but de cette section est la d\'emonstration du th\'eor\`eme suivant~:
\begin{thm}[Rel\`evement des Conjectures d'Arthur] \label{rel} \index{rel\`evement de Conjectures d'Arthur}
Soient $n$ et $n'$ deux entiers tels que $n\geq n'\geq 1$.
Soient $H \subset G$ deux groupes alg\'ebriques d\'efinis sur $\mathbb{Q}$ tels que :
$H^{{\rm nc}} \cong U(n',1)$ et $G^{{\rm nc}} \cong U(n,1)$. 
Supposons les hypoth\`eses $H_{\varepsilon}^{(p,q)}$ v\'erifi\'ees pour le groupe $H$ et pour tout couple 
d'entiers $(p,q)$ tel que $p+q \leq k \leq n'-1$. Alors, les hypoth\`eses $H_{n-n'+\varepsilon}^{(p,q)}$
sont vraies pour $G$ et pour tout couple d'entiers $(p,q)$ tel que $p+q  \leq k$.
\end{thm}

Le Th\'eor\`eme 2 est un corollaire imm\'ediat du Th\'eor\`eme \ref{rel}.
La d\'emonstration de ce dernier va essentiellement reposer 
sur un th\'eor\`eme de Burger et Sarnak \cite{BurgerSarnak} dont nous proposons la d\'emonstration suivante 
(emprunt\'ee \`a \cite{ClozelUllmo}).

\subsection*{Variation sur un th\`eme de Burger et Sarnak}

Si $\delta \in G({\Bbb Q})$, la double classe $\Gamma \delta \Gamma =T_{\delta}$ op\`ere de la fa\c{c}on suivante sur $L^2 (\Gamma \backslash G)$.
Si 
\begin{eqnarray} \label{orbitedehecke}
\Gamma \delta \Gamma = \coprod_i \Gamma \delta_i ,
\end{eqnarray}
alors
\begin{eqnarray} \label{opdehecke}
(T_{\delta})(g) = \sum_i f(\delta_i g) \quad (g \in G).
\end{eqnarray}

Soit deg$(T_{\delta})$ le degr\'e de $T_{\delta}$, \'egal au cardinal de $\Gamma \backslash \Gamma \delta \Gamma$. Nous appellerons {\it op\'erateur de Hecke}\index{op\'erateur de Hecke} une combinaison lin\'eaire finie \`a coefficients entiers $\geq0$ de $T_{\delta}$; son degr\'e est alors d\'efini par additivit\'e. Il n'est pas 
difficile de montrer (cf. \cite{BurgerSarnak}, \cite{ClozelUllmo}) que $T_{\delta}$ est un op\'erateur born\'e dans $L^2 (\Gamma \backslash G)$, de norme 
$$|| T_{\delta} || = \mbox{deg}(T_{\delta} ).$$ 
Si $T$ est un op\'erateur de Hecke, soit $\tilde{T} = \mbox{deg}(T)^{-1} T$ l'op\'erateur {\it normalis\'e}\index{op\'erateur de Hecke
normalis\'e} associ\'e.

Notons $L^2 (\Gamma \backslash G)^{\perp}$ l'orthogonal de l'espace des fonctions constantes. Le th\'eor\`eme suivant qui renforce une proposition de Burger et Sarnak est du \`a Clozel et Ullmo \cite{ClozelUllmo}.

\begin{thm} \label{CU}
Il existe un op\'erateur de Hecke autoadjoint $T$ tel que $\tilde{T} f=f$ pour $f$ constante sur $\Gamma \backslash G$ et $|| \tilde{T} _{|L^2 (\Gamma \backslash G)^{\perp}}|| <1$.
\end{thm}

La norme est la norme forte d'op\'erateur~:
\begin{eqnarray} \label{norme}
|| \tilde{T} f || \leq || \tilde{T}_{|L^2 (\Gamma \backslash G)^{\perp}} || \cdot || f || , \quad f \in L^2 (\Gamma \backslash G)^{\perp} .
\end{eqnarray}

La d\'emonstration s'esquisse comme  suit. Fixons une place $q$ telle que $G$ soit d\'eploy\'e sur ${\Bbb Q}_q$, que $G$ soit d\'efini sur ${\Bbb Z}_q$, et que 
l'intersection du sous-groupe compact-ouvert $K_f$, de $G({\Bbb A}_f )$ d\'efinissant le sous-groupe de congruence $\Gamma$, avec $G({\Bbb Q}_q )$ soit r\'eduite 
au sous-groupe hypersp\'ecial $G({\Bbb Z}_q )$. 

\`A l'exception du cas sp\'ecial (sans int\'eret pour nous) o\`u $G$ est obtenu par restriction des scalaires (pour une extension finie $F$ de ${\Bbb Q}$) de $SL(2)/F$, le groupe
$G({\Bbb Q}_q )$ est de rang $\geq 2$ et a donc la propri\'et\'e $(T)$ de Kazhdan. La repr\'esentation triviale est donc isol\'ee dans le dual unitaire de $G({\Bbb Q}_q )$. 
Clozel et Ullmo en d\'eduisent \cite{ClozelUllmo} le lemme suivant.

\begin{lem} \label{T}
Il existe une fonction $\varphi$ sur $G({\Bbb Q}_q )$, bi-$G({\Bbb Z}_q )$-invariante, positive, \`a coefficients entiers et auto-adjointe ($\varphi (g) = \varphi (g^{-1} )$)
et une constante $C<1$ tels que
\begin{eqnarray} \label{TT}
||\pi (\varphi ) || \leq C \mbox{deg}(\varphi ) 
\end{eqnarray}
si $\pi \in \widehat{G({\Bbb Q}_q )}$ est diff\'erente de la repr\'esentation triviale.
\end{lem}

Remarquons que si $G$ est obtenu par restriction des scalaires (pour une extension finie $F$ de ${\Bbb Q}$) de $SL(2)/F$, le Lemme \ref{T} reste vrai si l'on consid\`ere 
la composante $q$-adique d'une repr\'esentation automorphe de $G$, puisque l'on dispose d'une approximation de la Conjecture de Ramanujan.

Compl\`etons alors la d\'emonstration du Th\'eor\`eme \ref{CU}. \'Ecrivons 
\begin{eqnarray} \label{deco}
L^2 (\Gamma \backslash G) = {\Bbb C} \oplus L^2 (\Gamma \backslash G)^{\perp}
\end{eqnarray}
o\`u ${\Bbb C}$ d\'esigne l'espace des constantes; $L^2 (\Gamma \backslash G)$ est une repr\'esentation de $G$.
Le groupe $G$ est obtenu par restriction des scalaires (pour une extension finie $F$ de ${\Bbb Q}$) d'un groupe absolument quasi-simple $G^0 /F$. Sous nos hypoth\`eses,
la place $q$ est d\'ecompos\'ee dans $F$ et 
$$G({\Bbb Q}_q ) = \prod_{v|q} G^0 (F_v ),$$
chaque facteur \'etant isomorphe \`a $G^d ({\Bbb Q}_q )$ o\`u $G^d$ est le groupe d\'eploy\'e simplement connexe de m\^eme syst\`eme de racines que $G^0$.
Soit $v$ une place fix\'ee au-dessus de $q$. \'Ecrivons $K_f = K_v K^v$ ($K_v \subset G^0 (F_v )$) et soit $\Gamma_v = G({\Bbb Q} ) \cap K^v$. Alors,
\begin{eqnarray} \label{espce}
{\cal L}_v = L^2 (G({\Bbb Q}) \backslash G^0 ({\Bbb A}_F^v ) / K^v )
\end{eqnarray}
(notations \'evidentes) est une repr\'esentation de $G \times G^0 (F_v )$, et $L^2 (\Gamma \backslash G)$ est l'espace des $K_v$-invariants dans ${\cal L}_v$. 
Par approximation forte, $L^2 (\Gamma \backslash G)^{\perp} = ({\cal L}_v^{\perp} )^{K_v}$ o\`u ${\cal L}_v^{\perp}$ est l'orthogonal de l'espace des constantes.

La th\'eorie des s\'eries d'Eisenstein donne une d\'ecomposition 
\begin{eqnarray} \label{eisen}
{\cal L}_v = {\Bbb C} \oplus \int_{\widehat{G^0 (F_v )}} m(\pi_v ) \pi_v d\mu (\pi_v )
\end{eqnarray}
que nous n'expliciterons pas, mais o\`u l'int\'egrale porte sur l'espace des repr\'esentations automorphes non triviales dans $\widehat{G^0 (F_v )}$.
L'op\'erateur $T$ associ\'e \`a la fonction $\varphi$ d\'eduite du Lemme \ref{T} op\`ere alors sur ${\cal L}_v^{K_v}$ d\'ecompos\'e selon (\ref{eisen}) par 
\begin{eqnarray}
T = \mbox{deg} (T) \oplus \int_{\widehat{G^0 (F_v )}_{nr}} m(\pi_v ) \pi_v (\varphi ) d\mu (\pi_v );
\end{eqnarray}
o\`u $\pi_v (\varphi )$ est un scalaire de norme $\leq C \mbox{deg} (T)$. D'o\`u le Th\'eor\`eme \ref{CU}.

\medskip

\begin{lem}[Burger-Sarnak \cite{BurgerSarnak}]
Soit $f \in C_0 (\Gamma \backslash G)$.
Alors (\`a une constante strictement positive de normalisation pr\`es)~:
\begin{eqnarray}
\int_{\Gamma \backslash G} f(g) \overline{f(gh)} dg = \lim_{m\rightarrow +\infty} 
\frac{1}{d^m} \sum_{i=1}^{\lambda^{(m)}} \int_{\Lambda_i^{(m)} \backslash H} f(\eta_i^{(m)} h_1 ) \overline{
f(\eta_i^{(m)} h_1 h) } dh_1 ,
\end{eqnarray}
o\`u les $\Lambda_i^{(m)}$ sont des sous-groupes de congruences de $H$, $d={\rm deg}(T)$, $\eta_i^{(m)} \in G({\Bbb Q})$,
et la limite est uniforme sur les compacts de $H$.
\end{lem}
{\it Id\'ee de la d\'emonstration.} Soit $f \in C_0 (\Gamma \backslash G)$, et $f^0 (x) = \int_{\Gamma \backslash G} f(g) dg$ ($\forall x$) (mesure normalis\'ee)
le ``terme constant'' de $f$. Le Th\'eor\`eme 9.5.2 dit que $\tilde{T}^n f \rightarrow f^0$ dans $L^2$. Il n'est alors pas difficile de montrer (\cite{BurgerSarnak}) que 
$\tilde{T}^n f(x) \rightarrow f^0 (x)$ $\forall x \in \Gamma \backslash G$, la convergence \'etant uniforme sur les compacts. Par dualit\'e on a donc~:
$$\int_{\Gamma \backslash G} f(g) \overline{f(gh)} dg = \lim_{m \rightarrow \infty} \langle \tilde{T}^m  (\mu ) , R \rangle$$
o\`u 
$$R(g) =  f(g) \overline{f(gh)} $$
et $\mu$ est la mesure positive (normalis\'ee) sur $\Gamma \backslash G$ d\'efinie par
$$\langle \mu , F \rangle = \int_{(\Gamma \cap H) \backslash H} F(h) dh .$$

On \'ecrit 
$$T^m F(g) = \sum_{j=1}^{d^m} F(\delta_j^{(m)} g) .$$
L'ensemble des $\delta_j^{(m)}$ se d\'ecompose en une r\'eunion disjointe de $(\Gamma \cap H)$-orbites~:
$$\{ \delta_j^{(m)} \}_{j=1, \ldots , d^m} = \bigsqcup_{i=1}^{\lambda^{(m)}} \Gamma \eta_i^{(m)} (\Gamma \cap H) $$
o\`u chaque $\eta_i^{(m)} \in$ Comm$(\Gamma )$. 
Soit 
$$\Lambda_i^{(m)} = \{ h \in \Gamma \cap H \; : \; \Gamma \eta_i^{(m)} h = \Gamma \eta_i^{(m)} \} .$$
Le groupe $\Lambda_i^{(m)}$ est un sous-groupe de congruence de $H$.
Et,
$$\langle \tilde{T}^m (\mu ) , F \rangle = \frac{1}{d^m} \sum_{i=1}^{\lambda^{(m)}} \int_{\Lambda_i^{(m)} \backslash H} F(h_i^{(m)} h) dh.$$
Le Lemme s'en d\'eduit en prenant $F=R$.   

\bigskip

Le Lemme 9.5.4 implique imm\'ediatement le principe de restriction de Burger et Sarnak \cite{BurgerSarnak} que nous
avons rappel\'e au \pa 2.3 (Th\'eor\`eme 2.3.1).

\subsection*{D\'emonstration du Th\'eor\`eme \ref{rel}}

Nous supposons les hypoth\`eses $H_{\varepsilon}^{(p,q)}$ v\'erifi\'ees pour $H$
et pour tout couple d'entiers $(p,q)$ tel que $p+q \leq k$.
Fixons un couple d'entiers $(p,q)$ tel que $p+q \leq k$. Nous allons montrer que l'hypoth\`ese 
$H_{n-n'+\varepsilon}^{(p,q)}$ est v\'erifi\'ee pour le groupe $G$.

Soit donc $\lambda \in ]0 , (n-(p+q))^2 [$ dans le $(p,q)$-spectre automorphe de $G$. Soit $\tau = \tau_{p,q} '$.
D'apr\`es la formule de Matsushima, il existe une repr\'esentation $\sigma = \sigma_{a,b}  \in \hat{M} (\tau )$ et 
un nombre r\'eel strictement positif $s$ telle que~:
\begin{enumerate}
\item le quotient de Langlands  $J_{\sigma , s}$ de la repr\'esentation
$\pi_{\sigma , s}$ appartient \`a $\widehat{G}_{{\rm Aut}}$;
\item $\lambda =  (n-(a+b))^2 -s^2$.
\end{enumerate}

Nous noterons dor\'enavant $\pi = J_{\sigma ,s}$ (remarquons que les coefficients de $\pi$ sont des coefficients de $\pi_{\sigma ,s}$).

La restriction de la repr\'esentation $\pi$ \`a $H$ se d\'ecompose en somme directe :
\begin{eqnarray} \label{dec en somme directe}
\pi_{|H} = \int_{\hat{H}}^{\oplus} m_{\sigma} \sigma d\mu (\sigma ),
\end{eqnarray}
o\`u $\mu$ est une  mesure de probabilit\'e sur le dual unitaire $\widehat{H}$ de $H$ et $\sigma \mapsto m_{\sigma} \in \{ 1, 2, \ldots , \infty \}$
une fonction bor\'elienne.

Remarquons tout d'abord que le principe de restriction de Burger et Sarnak (Th\'eor\`eme 2.3.1) implique que :
\begin{eqnarray} \label{support de mu}
\mbox{support } \mu \subset \widehat{H}_{{\rm Aut}} .
\end{eqnarray}

Rappelons maintenant que $V_{\tau } \subset {\cal H}_{\pi}$. Soit $\tau'$ un $K_H$-type apparaissant dans $\tau_{|H}$. 
Fixons $v \in V_{\tau '}$. 
La d\'ecomposition (\ref{dec en somme directe}) implique la d\'ecomposition suivante :
\begin{eqnarray} \label{dec des coeff}
\langle \pi (a) v , v \rangle = \int_{\widehat{H} } m_{\sigma} \langle \sigma (a) v_{\sigma} , v_{\sigma} \rangle d\mu ( \sigma ) ,
\end{eqnarray}
o\`u $v  = \int_{\widehat{H}}^{\oplus} m_{\sigma} v_{\sigma}$ et $a$ est un \'el\'ement 
de $A_H$. Remarquons que dans la d\'ecomposition (\ref{dec des coeff}) la mesure $\mu$ ne charge que 
les repr\'esentations $\sigma \in \widehat{H}_{{\rm Aut}}$ contenant le $K_H$-type $\tau'$.

Nous allons appliquer le Corollaire 9.4.2 \`a $H$. Le Corollaire 9.4.2 donne pour chaque 
$\sigma \in \widehat{H}_{{\rm Aut}}$ contenant le $K_H$-type $\tau'$ et pour tout vecteur $v_{\sigma} \in V_{\tau'}$, 
\begin{eqnarray} \label{croissa}
\langle \sigma (e^{-tH}) v_{\sigma} , v_{\sigma} \rangle_{V_{\tau'}} 
\left\{
\begin{array}{cl}
= & B_{\tau '} (n'-j) ||v_{\sigma} ||_{V_{\tau'}}^2 e^{-jt} + o_{\varepsilon} (e^{-(n' + \delta)t} ) \\
\leq & C_{\varepsilon} ||v_{\sigma} ||_{V_{\tau'}}^2 e^{-(n' -\varepsilon)t}
\end{array}
\right.
\end{eqnarray} 
pour $t\in [0, + \infty [$, la constante implicite dans le $o_{\varepsilon}$, ainsi que $C_{\varepsilon}$, \'etant uniformes.

La premi\`ere ligne de (9.5.13) correspond aux repr\'esentations irr\'eductibles de $H$ (ind\'ex\'es par $j=0,1, \ldots ,k$) 
v\'erifiant la condition 1. de $H_{\varepsilon}^{(a,b)}$ ($a+b \leq k \leq n'-1$);
la seconde ligne correspond aux autres, v\'erifiant 2.

\medskip

Consid\'erons alors l'expression de $\langle \pi (e^{-tH}) v,v\rangle_{V_{\tau'}}$ donn\'ee par (\ref{dec des coeff}). Commen\c{c}ons par remarquer que 
\begin{eqnarray} \label{aze}
||v||_{V_{\tau '}}^2 = \int_{\widehat{H} } m_{\sigma} ||v_{\sigma} ||_{V_{\tau '}}^2 d\mu ( \sigma ) .
\end{eqnarray}
Selon le type de la repr\'esentation de $H$ dans le support de $\mu$, chaque 
coefficient (dans l'int\'egrale) s'exprime selon (\ref{croissa}); d'apr\`es l'uniformit\'e de $C_{\varepsilon}$ et des termes $o_{\varepsilon}$ et (\ref{aze}), 
on a l'expression~:
$$\langle \pi (e^{-tH}) v , v\rangle = \sum_{j=0}^{k} C_j \cdot e^{-jt} + O_{\varepsilon} (e^{-(n'-\varepsilon) t}) \quad (t\in [0, +\infty[ ) $$
avec les m\^emes propri\'et\'es des constantes implicites dans $O_{\varepsilon}$; mais toujours d'apr\`es le Corollaire 9.4.2 (cette fois appliqu\'e au groupe $G$)
ceci est en fait de la forme
$$C \cdot e^{-(n-s)t} + o(e^{-(n-s)t}) \quad (t\in[0 , + \infty [) ,$$
o\`u $C$ est une constante non nulle. 

Nous nous retrouvons alors avec deux possibilit\'es :
\begin{enumerate}
\item soit il existe un entier $j$ dans $[0,k]$ tel que $s=n-j$, 
\item soit $s\leq n-n' +\varepsilon$.
\end{enumerate}
En rempla\c{c}ant la valeur de $s$ dans l'expression de $\lambda$, on obtient le Th\'eor\`eme {\it i.e.}
que la propri\'et\'e $H_{\varepsilon}^{(p,q)}$ est v\'erifi\'ee.

\bigskip

\noindent
{\bf Remarque.} Sous l'hypoth\`ese $H_{n-n' + \varepsilon}$, la valeur propre $\lambda$ (apparaissant dans le $(p,q)$-spectre) est alors enti\`ere ou
$$\lambda \geq (n-n'+1)^2 - (n-n'+ \varepsilon)^2 >0 \quad {\rm si} \ \varepsilon <1$$ 
et la Conjecture $A^-$ est v\'erifi\'ee. 

\medskip

\newpage

\thispagestyle{empty}

\newpage

\markboth{CHAPITRE 10. D\'EMONSTRATION DU TH\'EOR\`EME 3}{10.1. VRAIS GROUPES UNITAIRES}

\chapter{D\'emonstration du Th\'eor\`eme~3}

Dans ce chapitre $G$ est un groupe anisotrope sur $\mathbb{Q}$ obtenu par restriction des scalaires \`a partir d'un groupe sp\'ecial unitaire $G_{F}$ sur un corps totalement r\'eel $F$. Alors $G(\mathbb{R})$ est un produit de groupes 
$SU(p,q)$, le produit portant sur les plongements r\'eels de $F$. Nous supposerons que
\begin{equation*}
G(\mathbb{R}) \cong SU (n,1)\times SU(n+1)^{d_1}\ .
\end{equation*} 

Notre but est de d\'emontrer les Conjectures A$^{-}(0)$ et A$^{-}(1)$ dans ce cas. Rappelons (\pa 8.1) qu'en g\'en\'eral un groupe unitaire provient d'une involution sur une alg\`ebre simple centrale. Si celle-ci est une alg\`ebre de matrices, on obtient simplement les groupes unitaires ``usuels'' des formes hermitiennes. Ce cas est trait\'e dans le \pa 10.1.

Quand $G$ est un groupe ``exotique'' (associ\'e \`a une alg\`ebre simple centrale qui n'est pas une alg\`ebre de matrices), des principes g\'en\'eraux impliquent que le spectre des $0$-formes et des $1$-formes devrait \^etre \textbf{plus restreint} que dans le cas pr\'ec\'edent, impliquant \textbf{a fortiori} les r\'esultats cherch\'es. Mais ceci -- qui repose sur des exemples de la fonctorialit\'e de Langlands -- n'est pas si facile \`a d\'emontrer. Nous obtenons le r\'esultat cherch\'e (Th\'eor\`eme~3 de l'Introduction, avec parfois de meilleures constantes spectrales) sous l'hypoth\`ese que le rang absolu du $G$ (c'est $(n+1)$ dans la description pr\'ec\'edente) n'est pas une puissance de $2$. Pour ceci nous sommes amen\'es \`a reprendre et \`a pr\'eciser les d\'emonstrations de \cite{Clozel} qui d\'emontraient la Conjecture $\tau$ pour de tels groupes.

Nous pensons que l'hypoth\`ese sur le rang peut \^etre \'evit\'ee, mais cela impose de renforcer significativement les r\'esultats de \cite{Clozel}.

\section{Vrais groupes unitaires}

Soient $F$ un corps totalement r\'eel de degr\'e $d$, $E/F$ une extension quadratique totalement imaginaire et $h$ une forme hermitienne sur $E^{n+1}$.

Pour toute place archim\'edienne $v$ de $F$, $E/F$ d\'efinit une extension quadratique $E\otimes F_{v}/F_{v}$ isomorphe \`a $\mathbb{C}/\mathbb{R}$ et $h$ d\'efinit donc une forme hermitienne $h_{v}$ sur $\mathbb{C}^{n+1}$. On suppose que les signatures sont $(n,1)$ en $v_{1}$ et $n$ en $v_{2},\ldots v_{d}$.

Le groupe $G_{F}$ est le groupe sp\'ecial unitaire, d\'efini sur $F$, $SU(E^{n+1},h)$ et $G$ est le groupe $G_{F}$ ``vu comme $\mathbb{Q}$-groupe''.

Dans ce cas le Th\'eor\`eme~3 r\'esulte des chapitres 8 et 9. On a en fait~:

\begin{thm}
Soit $G$ un vrai groupe unitaire, tel que $G^{\rm nc} \cong SU(n,1)$, $n\geq 2$.

$\mathrm{(1)}$ La plus petite valeur propre positive du laplacien sur les fonctions v\'erifie :
\begin{equation*}
\lambda_{1}^{0}\geq 2n-1
\end{equation*}  

$\mathrm{(2)}$ La plus petite valeur propre sur les $1$-formes v\'erifie :
\begin{equation*}
\lambda_{1}^{1}\geq \frac{2}{5}n-\frac{11}{25}.
\end{equation*}
\end{thm}
{\it D\'emonstration.} Consid\'erons d'abord les fonctions. On applique le Th\'eor\`eme 9.5.1 avec $n'=2$. D'apr\`es le Th\'eor\`eme~1
(Chapitre~8) on a alors (notations 9.5.1) $\varepsilon=\frac{4}{5}$. Alors le Th\'eor\`eme 9.5.1 donne

\begin{equation*}
\lambda=n^{2}-(n-i)^{2}\ ,\ i=0,\ldots n-1 
\end{equation*}
ou $\lambda\geq n^{2}-(n-2+\varepsilon)^{2}=\frac{12}{5}n-\frac{36}{25}$.

\noindent On v\'erifie ais\'ement que (pour $\lambda\neq 0$) la borne inf\'erieure obtenue est donn\'ee par $i=1$, $\lambda=2n-1$.

Si on consid\`ere les $1$-formes, on obtient de m\^eme, avec toujours $\varepsilon=\frac{4}{5}$ (Ch.~8)~:

\begin{equation*}
\lambda=(n-1+k)^{2}-(n-1+k-i)^{2}Ê\quad k=0,1\ ,\ i=0,\ldots,n-1\  
\end{equation*}
ou $\lambda\geq(n-1)^{2}-(n-\frac{6}{5})^{2}=\frac{2}{5}n-\frac{11}{25}$.

Dans ce cas la borne inf\'erieure est donn\'ee par la seconde estim\'ee.

\markboth{CHAPITRE 10. D\'EMONSTRATION DU TH\'EOR\`EME 3}{10.2. GROUPES EXOTIQUES : R\'EDUCTIONS}

\section{Groupes exotiques : r\'eductions}

En g\'en\'eral, un groupe alg\'ebrique $G$ sur $\mathbb{Q}$ tel que $G(\mathbb{R})\cong\mathrm{SU}(n,1)\times \mathrm{SU}(n)^{d-1}$ est obtenu de la fa\c con suivante (cf. \cite{PlatonovRapinchuk}).

On fixe comme auparavant un corps totalement r\'eel $F$, une extension quadratique totalement imaginaire $E$ de $F$.

Soit $\sigma$ le g\'en\'erateur de $\mathrm{Gal}(E/F)$. 

Soit $B$ une alg\`ebre simple sur $E$, dont le centre est $E$, et de rang r\'eduit $(n+1)$~: donc $B\otimes_{E}\overline{\mathbb{Q}}\cong M_{n+1}(\overline{\mathbb{Q}})$. On suppose donn\'ee sur $B$ une involution de seconde esp\`ece $\alpha$ (cf. \pa 8.1). Alors (notations du \pa pr\'ec\'edent)
\begin{equation*}
G(\mathbb{Q})=G_{F}(F)=\{g\in B^{\times }:g\alpha(g)=1\}\ .
\end{equation*}

Comme on l'a dit dans l'Introduction, nous utiliserons librement dans ce paragraphe les notions relatives aux formes automorphes sur les groupes ad\'eliques, ainsi que les r\'esultats de \cite{Clozel}.

Commen\c cons par reformuler la d\'efinition g\'en\'erale des groupes unitaires. Une alg\`ebre simple centrale $B$ sur $E$ s'\'ecrit $B=M_{r}(D)$ o\`u $D$ est une \textbf{alg\`ebre \`a division} de rang r\'eduit $d$ sur $E$. Le rang absolu du groupe unitaire associ\'e ($n+1$ avec nos notations usuelles) est $rd$. Notons $\ast$ une involution de seconde esp\`ece {\bf sur} $D$.
\index{involution de seconde esp\`ece} 
Il y a alors une notion naturelle de modules de rang $r$ sur $D$ (isomorphes \`a $D^r$); pour $V$ un tel module End$_D (V) \cong 
M_r (D)$. On peut d\'efinir des formes hermitiennes $h$ sur $V$ (\`a valeurs dans $D$) relatives \`a l'involution $\ast$ 
\cite{PlatonovRapinchuk}. Alors le groupe $G$ peut \^etre d\'efini de fa\c con naturelle par
\begin{equation*}
G(\mathbb{Q})=G_{F}(F)= SU(V,h).
\end{equation*}

Comme on l'a vu dans l'Introduction, nous supposons pour les groupes exotiques ({\it i.e.}, si $d>1$, ce que l'on supposera d\'esormais) 
que le rang $rd$ est impair. On a tout d'abord si $r>1$~:

\begin{prop}
Supposons $r\geq 3$ \(par exemple $r$ impair\). Alors les estim\'ees du Th\'eor\`eme 10.1.1 restent vraies pour $G$. 
\end{prop}

Ceci r\'esulte des arguments de \cite[\pa 1.3]{Clozel}. On peut supposer donn\'ee une base orthonormale de $D^{r}$ pour $h$. Celle-ci s'\'ecrit alors

\begin{equation*}
h(x,y)=\sum_{i=1}^{r}x_{i}^{\ast}f_{i}x_{i}
\end{equation*}
o\`u $x=(x_{1},\ldots x_{r}) $ et $y=(y_{1},\ldots y_{r})\in D^{r} $, $f_{i}\in D$ et $f_{i}=f_{i}^{\ast}$. La Proposition 1.3 de \cite{Clozel} montre que, sans changer $G$ (\`a isomorphisme pr\`es) on peut supposer que les $f_{i}$ commutent. Ils sont contenus alors dans un sous-corps totalement r\'eel maximal $L_{0}\subset D$, stable par l'involution. ($L_{0}$ est de degr\'e $d$ sur $F$). Le corps $EL_{0}=L\subset D$ est alors totalement imaginaire, quadratique sur $L_{0}$. L'argument donn\'e dans \cite[p.~303-304]{Clozel} montre alors que 
$G$ contient un sous-groupe $H$ sur
$\mathbb{Q}$ tel que avec $H(\mathbb{Q})=H_{L_{0}}(L_{0})$ ; ici 
$H_{L_{0}}$ est un groupe unitaire {\textbf{ordinaire} sur $L_{0}$, de rang $r$, dont la forme hermitienne est de matrice $(f_{1},\ldots f_{r})$.

Soient $v$ une place archim\'edienne de $F$ et $w_{1},\ldots w_{d}$ les places (r\'eelles) de $L_{0}$ \'etendant $v$. Si on note $(p,q)$ les signatures associ\'ees, on v\'erifie ais\'ement que

\begin{equation*}
p_{v}=\sum_{i}p_{w_{i}}\ ,\ q_{v}=\sum_{i}q_{w_{i}}\ .
\end{equation*} 
On en d\'eduit que $H_{L_{0}}(L_{0}\otimes \mathbb{R})=H(\mathbb{R})$ est isomorphe \`a $SU(1,r-1)\times SU(r)^{N-1} $ o\`u $N=d[F:\mathbb{Q}]$. Puisque c'est un groupe unitaire standard (et $r\geq 3$) il contient un groupe analogue de type $SU(2,1)\times SU(3)^{N-1}$. Enfin, la construction pr\'ec\'edente montre que le plongement $H\hookrightarrow G$, apr\`es ``complexification'' (extension des scalaires \`a $E$) provient d'un plongement $GL(r,L)(\subset GL(r,D) $ d\'eduit du plongement de $L$ dans $D$ comme sous-corps maximal. On en d\'eduit ais\'ement que le sous-groupe $SU(1,r-1)\subset SU(n,1)$ est donn\'e par le plongement usuel. Il en est de m\^eme pour $SU(2,1)$, et les arguments du Chapitre~9 s'appliquent.

Nous sommes maintenant r\'eduits \`a consid\'erer le cas o\`u $r=1$ ou $2$ et $d$ n'est pas une puissance de $2$. L'existence d'une base orthogonale montre alors que $G$ contient un sous-groupe de la forme

\begin{equation*}
H = SU(D,\ast)=\{x\in D:xx^{\ast}=1\}.
\end{equation*}
L'argument d\'ej\`a utilis\'e, relatif aux signatures, montre que l'on peut trouver $H$ de type $(d-1,1)$ en une place de $F$. Soit $\ell$ un diviseur premier impair de $d$. On a alors d\'emontr\'e dans \cite{Clozel} le fait suivant (\pa 1.3). Il existe une extension $L_{0}$ de $F$, totalement r\'eelle, de degr\'e $d/\ell$, une alg\`ebre \`a division $A$ sur $L=L_{0}\otimes E$ de degr\'e $\ell$, munie d'une involution de seconde esp\`ece (relativement \`a $L/L_{0})$ toujours not\'ee $\ast$, et telle que le groupe $SU(A,\ast)$ -- un groupe
sur $L_{0}$ -- se plonge dans $H = SU(D,\ast)$. Comme pr\'ec\'edemment, on v\'erifie ais\'ement que

\begin{equation*}
H_{1}(\mathbb{R})=SU(\ell-1,1)\times SU(\ell)^{N-1}
\end{equation*}
o\`u $N=(d/\ell)[F:\mathbb{Q}]$ se plonge dans $G(\mathbb{R})$ comme sous-groupe unitaire standard. 
(On a not\'e $H_{1}$ le $\mathbb{Q}$-groupe d\'eduit de $SU(A,\ast))$.

\markboth{CHAPITRE 10. D\'EMONSTRATION DU TH\'EOR\`EME 3}{10.3. GROUPES EXOTIQUES DE RANG PREMIER}

\section{Contr\^ole du spectre pour les groupes exotiques de rang premier}

Nous supposons maintenant que $G$ provient d'une alg\`ebre \textbf{\`a division}, de degr\'e premier impair $\ell$, sur $E$ o\`u $E/F$ est quadratique et totalement imaginaire. Comme on l'a annonc\'e dans l'Introduction \`a ce chapitre, on va obtenir dans ce cas des estim\'ees particuli\`erement fortes. (Ce ph\'enom\`ene est bien connu ; pour un analogue cohomologique voir \cite{Clozel2}). Comme on va le voir, on se trouve dans la situation optimiste du \pa 8.4 o\`u le changement de base est connu. De plus les repr\'esentations de $\mathrm{GL}(\ell,\mathbb{C})$ obtenus par changement de base sont plus s\'ev\`erement restreintes que dans le cas du \pa 8.4.

Supposons que $G(\mathbb{R})\cong SU(n,1)\times SU(n+1)^{d-1}$. Rappelons (Ch.~4, \pa~5)~\footnote{$\tau$ est not\'e $\sigma$ au Ch.~4. Dans ce chapitre $\sigma$ d\'esigne la conjugaison complexe\dots } 
que les repr\'esentations de $U(n,1)$ sont temp\'er\'ees ou de la forme $J(\tau,\chi)$ = quotient de

\begin{equation*}
J(\tau,\chi)=\mathrm{ind}_{U(n-1)\times\mathbb{C}^{\times }\times N }^{U(n,1)}\tau\otimes \chi
\end{equation*}
o\`u $\chi(z)=z^{\alpha}(\overline z)^{\beta}$ $(z\in \mathbb{C}^{\times })$. On veut borner les s\'eries compl\'ementaires, pour lesquelles $\alpha+\beta\in\mathbb{R}$.

\begin{thm}
Si $J(\tau,\chi)$ appara\^{\i}t dans $L^{2}(\Gamma\ba G)$ pour un sous-groupe de congruence, et si $\alpha+\beta\in\mathbb{R}$,
\begin{equation*}
\left| \frac{\alpha+\beta}{2} \right| \leq\frac{1}{2}-\frac{1}{\ell^{2}+1}\ . 
\end{equation*}
\end{thm}

Pour d\'emontrer ceci nous utilisons les m\'ethodes de \cite{Clozel}. Comme dans cet article, nous d\'esignons maintenant par 
$G$ le groupe de \textbf{similitudes} unitaires d\'efini par $(D,\ast)$ :

\begin{equation*}
G_{1}(F)= G(\mathbb{Q})=\{d\in D:dd^{\ast}\in F^{\times }\},
\end{equation*}
o\`u $G_{1}$ est un $F$-groupe. 

On a alors \cite[\pa 2.1]{Clozel} $G_{1}(E)=\mathbb{G}_{E}(\mathbb{Q})\cong D^{\times }\times E^{\times }$. Le groupe $G_{1}(E)$ est donc d\'eduit de
$D^{\times }\times E^{\times }$ par restriction des scalaires~\footnote{Le lecteur ne perdra rien \`a supposer que $F=\mathbb{Q}$.}.

Soient $A_{G}=\mathbb{R}_{+}^{\times }\subset G= G(\mathbb{R})$ (inclusion centrale). On a de m\^eme $A_{G_{E}}=\mathbb{R}_{+}^{\times }\subset G_{E}(\mathbb{R})$.

Si $f\in C_{c}^{\infty }(G(\mathbb{A}))$ et $\varphi\in C_{c}^{\infty }(G_{E}(\mathbb{A}))$ on consid\`ere la trace de $f$ dans la repr\'esentation $r$ sur $L^{2}(A_{G}G(\mathbb{Q})\ba G(\mathbb{A}))=\mathcal{A}_{G}$. De m\^eme on consid\`ere la trace de $\varphi$ dans la repr\'esentation $R$ sur $L^{2}(A_{G_{E}} G_{E}(\mathbb{Q})\ba G_{E}(\mathbb{A}))=\mathcal{A}_{G_{E}}$.

Supposons $f$ et $\varphi$ \textbf{associ\'ees}, {\it i.e.}, leurs int\'egrales orbitales (locales) se correspondent en toutes les places, cf. \cite[D\'ef.~2.7]{Clozel}. Soit $I_{\sigma}$ l'op\'erateur d'entrelacement de $R$ associ\'e \`a la conjugaison galoisienne de 
$G_{1}(E)$ par rapport \`a $G_{1}(F)$. Alors \cite[(2.17)]{Clozel}~:
\begin{equation*}
\mathrm{trace}(r(f))=\mathrm{trace}(R(\varphi)I_{\sigma}).
\end{equation*}

Si $\pi$ d\'ecrit les repr\'esentations de $G (\mathbb{A})$ -- avec leurs multiplicit\'es -- dans $\mathcal{A}_{G_{E}}$, ceci s'\'ecrit
\begin{equation} \label{10.1}
\sum_{\pi}\mathrm{trace}\ \pi(f)=\sum_{\Pi}\mathrm{trace}(\Pi(\varphi)I_{\sigma}).
\end{equation}

D'apr\`es un r\'esultat de \cite{AC} et \cite{V} (d\'emontr\'e par Harris et Taylor dans \cite[Ch.~VI]{HarrisTaylor}) 
il n'y a pas de multiplicit\'es dans la somme de droite; celle ci ne porte que sur les repr\'esentations $\Pi$ telles que $\Pi\cong\Pi\circ\sigma$. Les deux sommes convergent absolument pour des fonctions lisses.

Afin de poursuivre nous devons avoir un meilleur contr\^ole sur les repr\'esentations $\Pi$ de (\ref{10.1}) ainsi que sur les fonctions associ\'ees $f$ et $\varphi$ (aux places archim\'ediennes).

Pour toute place $v$ de $F$, $G_{1}(E\otimes F_{v})\cong \mathbb{C}^{\times }\times GL(\ell,\mathbb{C})$. Noter qu'en une place r\'eelle on a

\begin{equation*}
G_{1}(F_{v})\cong GU(\ell-1,1)(\mathbb{R})\cong \left( \mathbb{R}^{\times }\times U(\ell-1,1)(\mathbb{R})\right) /\pm 1\ .
\end{equation*}
On peut essentiellement n\'egliger la composante centrale dans les arguments locaux qui suivent (elle n'\'etait introduite que pour simplifier la d\'emonstration de (\ref{10.1}), cf. \cite{Clozel}). Consid\'erons alors $U(\ell-1,1)\subset GL(\ell,\mathbb{C})$, ou $U(\ell)\subset GL(\ell,\mathbb{C})$, aux places r\'eelles, et rempla\c cons $G_{E}$ par
$D^{\times }$. (Notation provisoire : $G\subset G_{\mathbb{C}}$ d\'esigne $U(\ell -1,1) \subset GL(\ell , {\Bbb C})$. On n\'eglige les
arguments, plus faciles, portant sur les autres places r\'eelles.) 

\begin{lem}
Soit $v$ est une place r\'eelle de $F$, $w$ la place complexe de $E$ associ\'ee. Si $\Pi$ intervient dans (\ref{10.1}) $\Pi_{w}$ est g\'en\'erique ou est un caract\`ere ab\'elien.
\end{lem}

Il r\'esulte en effet de \cite[Thm.~VI.1.1]{HarrisTaylor} que $\Pi_{w}$ co\"{\i}ncide avec la composante locale d'une repr\'esentation $\Pi'$ apparaissant dans le spectre discret de $GL(\ell,\mathbb{A}_{E})$. Puisque $\ell$ est premier, il y a d'apr\`es M\oe glin-Waldspurger deux possibilit\'es \cite{MoeglinWaldspurger}~: $\Pi'$ appartient aux formes cuspidales, et donc $\Pi_{w}'=\Pi_{w}$ est g\'en\'erique. Ou bien $\Pi'$ est un caract\`ere ab\'elien de $GL(\ell,\mathbb{A}_{E})$. Rappelons qu'une repr\'esentation (unitaire) g\'en\'erique de $GL(n,\mathbb{C})$ n'est autre qu'une repr\'esentations unitairement induite (i.e., \'egale \`a l'induite \textbf{totale}) \`a partir de caract\`eres (non n\'ecessairement unitaires) du sous-groupe de Borel : c'est une s\'erie principale \textbf{irr\'eductible}.

Consid\'erons une repr\'esentation $\pi$ apparaissant dans (\ref{10.1}). D'apr\`es le classification de Langlands (Ch.~4) il y a deux cas possibles, Dans le premier cas, $\pi$ est de la forme $J(\tau,\chi)$ (ou $\pi$ est une composante d'une repr\'esentation $I(\tau,\chi)$ avec $\chi$ \textbf{unitaire} si celle-ci est r\'eductible). La repr\'esentation $I(\tau,\chi)$ est elle-m\^eme d\'eduite d'un param\`etre de Langlands

\begin{align*}
&W_{\mathbb{C}}=\mathbb{C}^{\times }\rightarrow GL(\ell,\mathbb{C})\\
&z\mapsto((z/\overline z)^{p_{1}},\ldots(z/\overline z)^{p_{\ell-2}}\ ,\ z^{\alpha}\overline z^{\beta}\ ,\ z^{-\beta}(\overline z)^{-\alpha})
\end{align*}
avec $p_{i}\in\mathbb{Z}$, $p_{i}$ distincts.

Cette donn\'ee d\'efinit \`a son tour une repr\'esentation de la s\'erie principale pour $GL(\ell,\mathbb{C})$ (Ch.~3). En g\'en\'eral (si $\alpha+\beta\notin i\mathbb{R}$) cette repr\'esentation n'est pas irr\'eductible.

Notons $I_{\mathbb{C}}(\tau,\chi)$ la repr\'esentations de $GL(\ell,\mathbb{C})$ ainsi d\'efinie (s\'erie principale, peut-\^etre r\'eductible). Si $\sigma$ est la conjugaison complexe de $GL(\ell,\mathbb{C})$ par rapport \`a $U(\ell-1,1)$, on v\'erifie que $I_{\mathbb{C}}$ est $\sigma$-invariante : $I_{\mathbb{C}}\cong I_{\mathbb{C}} \circ\sigma$~\footnote{Dans un groupe de Grothendieck convenable si $I_{\mathbb{C}}$ est r\'eductible\dots }. Mieux, d'apr\`es \cite{BC} il existe un op\'erateur d'entrelacement $A_{\sigma}:I_{\mathbb{C}}\rightarrow I_{\mathbb{C}}$ entrela\c cant
$I_{\mathbb{C}}$ et $I_{\mathbb{C}}\circ\sigma$, et uniquement d\'efini au signe pr\`es si $A_{\sigma}^{2}=1$ (on le d\'efinit d'abord pour $\chi$ unitaire, puis par prolongement analytique pour tout $\chi$.)

Dans le second cas $\pi$ est une repr\'esentation de la s\'erie discr\`ete de $G$, associ\'ee (\pa 4.3) \`a

\begin{equation}\label{10.2}
z\mapsto((z/\overline z)^{p_{1}},\ldots,(z/\overline z)^{p_{\ell}}) 
\end{equation}
o\`u $p_{i}\in\mathbb{Z}$, $p_{i}$ distincts. Alors (\ref{10.2}) d\'efinit de m\^eme un param\`etre de Langlands pour $GL(\ell,\mathbb{C})$, donc une repr\'esentation $\pi_{\mathbb{C}}$ qui appartient \`a la s\'erie principale \textbf{unitaire}. Elle est $\sigma$-invariante, d'o\`u $A_{\sigma}$ comme ci-dessous.

Rappelons (\pa 4.4) que l'on associe au param\`etre (\ref{10.2}) $\ell$ repr\'esentations de $G$ de la s\'erie discr\`ete, qui forment un $L$-\textbf{paquet}. On le notera $\Pi$; $\pi_{\mathbb{C}}$ est donc d\'eduite du $L$-paquet. On pose
\begin{equation*}
\mathrm{trace}\ \pi(f)=\sum_{\pi'\in\Pi}\mathrm{trace}\ \pi'(f).
\end{equation*}

\begin{thm}
Soit $f\in C_{c}^{\infty }(G)$ une fonction $K$-finie. Il existe alors une fonction $K$ finie $\varphi\in C_{c}^{\infty }(G_{\mathbb{C}})$ telle que :

\item[$\mathrm{(i)}$]
\begin{itemize}
Si $\pi$ est une repr\'esentation de la s\'erie principale de $G$, ou si $\pi=\Pi$ est un $L$-paquet de repr\'esentations de la s\'erie discr\`ete, et $\pi_{\mathbb{C}}$ est la repr\'esentation $\sigma$-stable de $G_{\mathbb{C}}$ associ\'ee,
\begin{equation} \label{10.3}
\mathrm{trace}\ \pi(f)=\mathrm{trace}\ \pi_{\mathbb{C}}(\varphi)A_{\sigma}).
\end{equation}
\end{itemize}

\item[$\mathrm{(ii)}$]
\begin{itemize}
Si $\pi_{\mathbb{C}}$ est une repr\'esentation de la s\'erie principale de $G_{\mathbb{C}}$ qui ne provient pas de $G$ \(i.e. de $\pi$ discr\`ete ou de $I(\tau,\chi))$,
\begin{equation*}
\mathrm{trace}(\pi_{\mathbb{C}}(\varphi)A_{\sigma})=0
\end{equation*}
si $\pi_{\mathbb{C}}$ est $\sigma$-stable.
\end{itemize}

De plus, $f$ et $\varphi$ sont associ\'ees au sens de \cite{Clozel}.
\end{thm}

Ceci r\'esulte du travail de Delorme \cite{D} et de l'argument de \cite[Ch. 1, \S 7]{AC}.

Noter que dans les identit\'es pr\'ec\'edentes, $A_{\sigma}$ doit \^etre choisi convenablement ; par ailleurs les identit\'es sont vraies pour les s\'eries principales (unitaires) compl\`etes -- sans les r\'eduire -- et alors par prolongement analytique pour toutes les s\'eries principales (g\'en\'eralis\'ees) compl\`etes quels que soient leurs param\`etres.

\begin{cor}
Si $\pi_{\mathbb{C}}$ est une s\'erie principale $\sigma$-stable pour $G_{\mathbb{C}}$, $\mathrm{trace}(\pi_{\mathbb{C}}(\varphi)A_{\sigma})$ est nul ou de la forme $\displaystyle\sum_{i}\mathrm{trace}(\pi_{i}(f)$, les $\pi_{i}$ \'etant des repr\'esentations \(peut-\^etre non unitaires\) de $G$ en nombre fini.
\end{cor}

Ceci r\'esulte du Th\'eor\`eme, (i).

On aura besoin d'\'etendre ces identit\'es au cas des caract\`eres ab\'eliens de $G_{\mathbb{C}}$, suppos\'es $\sigma$-stables. Soit $\varepsilon$ un tel caract\`ere ; on l'\'ecrit par abus de langage $\varepsilon(g)=\varepsilon(\det g)$ o\`u $\varepsilon$ est un caract\`ere de $\mathbb{C}^{\times }$. Il est $\sigma$-stable si $\varepsilon = z^{p}\overline z^{q}$ avec $p=-q\in\frac{1}{2}\mathbb{Z}$. Si $p\in\mathbb{Z}$, $\varepsilon$ ``provient de $G$'' (ou de $\varepsilon_0$), {\it i.e.} est associ\'e naturellement au caract\`ere 
$\varepsilon_{0} : g\mapsto\det(g)^{p}$ de $G$. Sinon, on dira que $\varepsilon$ ne provient pas de $G$.

\begin{lem}
Soit $f$, $\varphi$ comme dans le Th\'eor\`eme 10.3.3. Alors

\item[$\mathrm{(i)}$]
\begin{itemize}
Si $\varepsilon$ provient de $\varepsilon_{0}$, 
\begin{equation*}
\langle \mathrm{trace}\ \varepsilon_{0},f \rangle =\eta \langle \mathrm{trace}\ \varepsilon,\varphi \rangle ,\eta=\pm 1\ .
\end{equation*}
\end{itemize}

\item[$\mathrm{(ii)}$]
\begin{itemize}
Si $\varepsilon$ ne provient de $G$, 
\begin{equation*}
\langle \mathrm{trace}\ \varepsilon,\varphi \rangle =0\ .
\end{equation*}
\end{itemize}
\end{lem}

Nous esquissons seulement l'argument. Pour (ii) noter que si $\pi$ et $\pi_{\mathbb{C}}$ sont associ\'ees il r\'esulte ais\'ement de 
(\ref{10.3}) que $\theta_{\pi}\circ N=\theta_{\pi_{\mathbb{C}}}$, $\theta_{\pi}$ et $\theta_{\pi_{\mathbb{C}}}$ (sur $U(1)$ et 
$\mathbb{C}^{\times })$ \'etant les caract\`eres centraux et $N:Z(G_{\mathbb{C}})=\mathbb{C}^{\times }\rightarrow Z(G)=U(1)$ \'etant donn\'ee par $z\mapsto z/\overline z$.

En particulier si $\pi_{\mathbb{C}}$ est une s\'erie principale g\'en\'eralis\'ee provenant de $G$, $\theta_{\pi_{\mathbb{C}}}(z)$ est de la forme $(z/\overline z)^{p}$ pour $p$ \textbf{entier}.

Si $\varepsilon$ ne provient pas de $G$, son caract\`ere central, \'egal \`a $\varepsilon^{\ell}$, n'a pas cette propri\'et\'e puisque $\ell$ est \textbf{impair}. Or d'apr\`es un th\'eor\`eme bien connu, $\varepsilon$ peut s'\'ecrire comme somme altern\'ee
\begin{equation*}
\varepsilon=\sum_{i=1}^{r}n_{i}\mathrm{ind}_{B_{\mathbb{C}}}^{GL(\ell,\mathbb{C})}\chi_{i}=\sum_{i=1}^{r}n_{i}I_{i}
\end{equation*}
o\`u les $\chi_{i}$ sont des caract\`eres du groupe de Borel, que l'on peut supposer en situation \textbf{positive} au sens de Langlands (Ch.~3), et les $n_{i}$ sont des entiers relatifs. Une telle d\'ecomposition est alors unique. Les induites et $\varepsilon$ ont le m\^eme caract\`ere central.

On peut choisir la conjugaison complexe $\sigma$, modulo conjugaison dans $G(\mathbb{C})$, telle qu'elle laisse $B_{\mathbb{C}}$ invariant (Ch.~4). Noter que $\varepsilon$ est $\sigma$-stable et que (en prenant $A_{\sigma}=1$ dans l'espace de $\varepsilon$) 
$\langle \textrm{trace}\ \varepsilon,\varphi \rangle = \langle \textrm{trace}\ \varepsilon,\varphi\times A_{\sigma} \rangle $. 
Il r\'esulte de l'unicit\'e de la d\'ecomposition que $\sigma$ fixe les $I_{i}$ ou les \'echange deux \`a deux sans point fixe. Alors

\begin{equation*}
\langle \textrm{trace}\ \varepsilon,\varphi \rangle =\sum_{i=1}^{r'}n_{i}\eta_{i} \langle \textrm{trace}\ I_{i},\varphi\times A_{\sigma}^{i} \rangle
\end{equation*}
o\`u les $\eta_{i}$ sont des signes, et la somme ne porte que sur les $I_{i}$ qui sont $\sigma$-stables. Si le caract\`ere central ne provient pas de $G$, la somme de droite est nulle, q.e.d. La partie (i) est d\'emontr\'ee dans \cite[Ch.~3]{BC}.

Revenons alors \`a l'\'egalit\'e (10.1) ; $v$ d\'esigne toujours la place archim\'edienne distingu\'ee de $F$. Notons simplement $\mathfrak{g}_{v}$ l'alg\`ebre de lie du groupe r\'eel $G_{1}(F_{v})$ et soit $\mathfrak{Z}$ le centre de son alg\`ebre enveloppante (complexe). Soit $\mathfrak{Z}_{\mathbb{C}}$ l'objet analogue pour $G_{1}(E\otimes F_{v})$. Donc $\mathfrak{Z}_{\mathbb{C}}\cong\mathfrak{Z}\otimes \mathfrak{Z}$ et il y a une application norme naturelle $N:\mathfrak{Z}_{\mathbb{C}}\rightarrow \mathfrak{Z}\otimes \mathfrak{Z}$  
\cite[\pa 4.2]{Clozel}. Si $\omega$ est un caract\`ere infinit\'esimal, $\omega:\mathfrak{Z}\rightarrow\mathbb{C}$ pour $G_{1}(F_{v})$ on en d\'eduit un caract\`ere $\Omega=\omega\circ N$ pour $\mathfrak{Z}_{\mathbb{C}}$. On dira que $\omega$ et $\Omega$ sont associ\'es.

D'apr\`es un r\'esultat fondamental d'Arthur (cf. \cite[p.~320]{Clozel}) on peut s\'eparer dans (10.1) les contributions des  caract\`eres infinit\'esimaux~:  l'\'egalit\'e reste vraie, $\omega$ \'etant fix\'e, quand $\pi$ d\'ecrit les repr\'esentations telles que $\omega(\pi_{v})=\omega$ et $\Pi$ celles v\'erifiant $\omega(\Pi_{v})=\Omega$. Si les fonctions $f^{v}$ et $\varphi^{v}$ sont fix\'ees en les places diff\'erentes de $v$ ($K$-finies aux places archim\'ediennes), les sommes portent alors sur un nombre \textbf{fini} de repr\'esentations.

Choisissons alors une repr\'esentation $\pi_{v}$ de $G_{1}(F_{v})$ apparaissant dans (10.1) et s\'eparons l'identit\'e selon les repr\'esentations de $G_{1}(F_{v})$ et $G_{1}(E\otimes F_{v})$ :

\begin{equation}
\begin{split}
\textrm{trace}\ \pi_{v}(f_{v}) &\sum_{\pi}\textrm{trace}\ \pi^{v}(f^{v})+\sum_{\rho_{v}\tilde\neq\pi_{v}}\textrm{trace}\ \rho_{v}(h)\sum_{\rho}\textrm{trace}\ \rho^{v}(f^{v})=\\
=&\sum_{\pi_{w}}\textrm{trace}(\pi_{w}(\varphi_{w})A_{\sigma})\sum_{\pi}\textrm{trace}(\pi^{w}(\varphi^{w})I_{\sigma}^{w}).
\end{split}
\end{equation}

Dans le membre de gauche, $\displaystyle\sum_{\pi}$ porte sur les repr\'esentations telles que $\pi_{v}$ soit la repr\'esentation fix\'ee;$\rho_{v}$ parcourt toutes les autres repr\'esentations de $G_{1}(F_{v})$; pour chacune, $\displaystyle\sum_{\rho}$ est d\'efini de m\^eme. A droite, $\Pi_{w}$ d\'ecrit les repr\'esentations de $G_{1}(E\otimes F_{v})$, et l'on a d\'ecompos\'e $I_{\sigma}$ (dans l'espace de $\Pi$) en un produit tensoriel $A_{\sigma}$ et $I_{\sigma}^{w}$.

Fixons $\pi^{o}$ telle que $(\pi^{o})_{v}=\pi_{v}$. Alors $\pi^{o}$ est fixe pour un sous-groupe compact $K_{f}\subset G_{1}(\mathbb{A}_{F,f})= G(\mathbb{A}_{\mathbb{Q},f})$; pour $K_{f}$ assez petit, la fonction caract\'eristique $f_{f}$ de $K_{f}$ admet une fonction
$\varphi_{f}$ associ\'ee (cf. \cite[\pa 2.5]{Clozel}). (On peut choisir les fonctions $f_{v'}$ aux places archim\'ediennes $\neq v$ \'egales \`a des caract\`eres des repr\'esentations $\pi_{v'}^{o}$). Alors $f^{v}=\displaystyle\bigotimes_{u\neq v}f_{u} $ admet une fonction associ\'ee $\varphi^{v}$, et la partie de (10.4) relative \`a $\pi_{v}$ est de la forme $c\ \textrm{trace}\ \pi_{v}(f_{v})$ o\`u $c$ est une constante $>0$.  

D'apr\`es le Lemme 10.3.2, le Corollaire 10.3.4 et le Lemme 10.3.5, le membre de droite s'\'ecrit
\begin{equation}
\sum_{\sigma_{v}}c(\sigma_{v})\textrm{trace}\ \sigma_{v}(f_{v})
\end{equation} 
o\`u $\sigma_{v}$ parcourt un ensemble fini de repr\'esentations de $G_{1}(F_{v})$. L'\'egalit\'e :

\begin{equation}
\begin{split}
c\ \textrm{trace}\ \pi_{v}(f_{v}) &+\sum_{\rho_{v}\not\cong \pi_{v}}c(\rho_{v})\textrm{trace}\ \pi^{v}(f^{v})=\\
&=\sum_{\sigma_{v}}c(\sigma_{v})\ \textrm{trace}\ \sigma_{v}(f_{v}),
\end{split}
\end{equation}
\textbf{qui porte sur un ensemble fini de repr\'esentations}, montre alors qu'il existe $\sigma_{v}$ \'egale \`a $\pi_{v}$.

Dans le membre de droite, $\sigma_{v}$ \'etait associ\'ee \`a une repr\'esentation $\Pi_{w}$ par le Corollaire 10.3.4 ou le Lemme 10.3.5. Dans le second cas $\sigma_{v}=\varepsilon_{0}$ est un caract\`ere ab\'elien. Dans le premier, $\sigma_{v}$ est l'une des composantes d'une s\'erie principale, ou appartient \`a la s\'erie discr\`ete (Thm.~10.3.3).

Nous pouvons maintenant d\'emontrer le Th\'eor\`eme 10.3.1. Supposons que $\pi_{v}=J(\tau,\chi)$ comme repr\'esentation de $U(\ell-1,1)$.
Donc $\pi_{v}$ est associ\'ee par le Cor. 10.3.4 \`a une repr\'esentation \textbf{g\'en\'erique} $\Pi_{w}$ qui appara\^{\i}t dans les formes automorphes sur $D^{\times }(\mathbb{A}_{E})$ ; d'apr\`es le th\'eor\`eme d'Harris et Taylor d\'ej\`a cit\'e, $\Pi_{w}$ est la composante locale d'une repr\'esentation cuspidale de $GL(\ell,\mathbb{A}_{E})$. Puisque la trace tordue de $\Pi_{w}$ \'evalu\'ee sur une fonction $\varphi_{w}$ provenant de $G_{v}$ est non-nulle, $\Pi_{w}=\pi_{\mathbb{C}}$ doit \^etre l'une des repr\'esentations d\'ecrites dans le Th\'eor\`eme~10.3.3.

Le cas des repr\'esentations $\pi_{\mathbb{C}}$ ``provenant des s\'eries discr\`etes'' est exclu car l'identit\'e de caract\`eres impliquerait que $\pi_{v}$ serait une s\'erie discr\`ete. Donc $\Pi_{w}$ provient d'une repr\'esentation $I(\sigma,\chi)$. Ecrivons

\begin{equation*}
\Pi_{w}=\ind_{B_{\mathbb{C}}}^{\mathrm{GL}(\ell,\mathbb{C})}((z/\bar z)^{p_{1}},\ldots(z/\bar z)^{p_{n-1}},z^{\alpha}\bar z^{\beta},z^{-\beta}(\bar z)^{-\alpha}). 
\end{equation*}
(Les donn\'ees $p$, $\alpha$, $\beta$ ne sont pas n\'ecessairement, pour l'instant, celles de $\pi_{v}$). Le Th\'eor\`eme 7.0.1 implique alors~:
\begin{equation*}
\Bigl|\frac{\alpha+\beta}{2}\Bigr|\leq\frac{1}{2}-\frac{1}{\rho^{2}+1}<\frac{1}{2}\ .
\end{equation*}

De plus $\Pi_{w}$ est irr\'eductible, et
\begin{equation*}
\textrm{trace}\ \Pi_{w}(\varphi_{w}A_{\sigma})=\sum_{i}\textrm{trace}\ \pi_{i}(f_{v})
\end{equation*}
o\`u la somme porte sur les composantes \'eventuelles de la repr\'esentation $I(\tau',\chi')$ associ\'ee \`a $\Pi_{w}$.

Il en r\'esulte que $J(\tau,\chi)$ est l'une de ces composantes. Mais la condition $|\alpha+\beta|<\frac{1}{2}$ implique que $I(\tau',\chi')$ est irr\'eductible, sauf peut-\^etre si $\alpha+\beta=0$. (Voir les r\'esultats de Knapp cit\'es dans le \pa 6.2. Pour $\alpha+\beta=0$, on peut avoir r\'eductibilit\'e des s\'eries principales unitaires, cf. \pa 4.4).

Si $\alpha+\beta\neq 0$, on voit donc que $I(\tau',\chi')$ est irr\'eductible et \'egale \`a $J(\tau,\chi)$, d'o\`u le Th\'eor\`eme.

\vskip0,2cm
\noindent
{\bf Remarque.} On a donc d\'emontr\'e que toute repr\'esentation non-temp\'er\'ee et non-ab\'elienne de $G(\mathbb{R})$ qui appara\^{\i}t dans les formes automorphes est une s\'erie principale. En particulier les modules de Vogan-Zuckerman de degr\'e primitif $\neq 0$, $d=\ell-1$ n'apparaissent pas. Si $\Gamma$ est un sous-groupe de congruences de $G$, on voit donc~:

\begin{thm}
Si $\Gamma\subset G(\mathbb{Q})$ est un sous-groupe de congruences et $X=\Gamma\ba SU(\ell-1,1)/K_{\infty }$, $X$ n'a pas de cohomologie primitive en degr\'es $0<i<\ell-1$\ .
\end{thm}

Ce r\'esultat \'etait d\'emontr\'e de fa\c con assez diff\'erente dans \cite{Clozel2}.

\markboth{CHAPITRE 10. D\'EMONSTRATION DU TH\'EOR\`EME 3}{10.4. D\'EMONSTRATION DU TH\'EOR\`EME 3}

\section{D\'emonstration du Th\'eor\`eme 3}

Si $G$ n'est pas un vrai groupe unitaire, on suppose $n+1\neq 2^{e}$. On utilise les r\'eductions du \pa 2 ; en particulier on peut supposer que $G$ contient un sous-groupe $H$ de rang premier impair $\ell$ du type consid\'er\'e dans le \pa 3. (On a $n+1=d$ ou $2d$, et $\ell$ est un diviseur premier de $d$). On utilise le Th\'eor\`eme 9.5.1 avec $n'=\ell-1$.

Consid\'erons d'abord le cas des fonctions. Pour le ``petit groupe'' $H$, la minoration $H_{(\varepsilon)}^{(0,0)}$ est ici, avec les notations du \pa 3, cf. aussi Prop.~4.5.1~:

\begin{equation*}
\lambda=n'{}^{2}-(\alpha+\beta)^{2}\geq(\ell-1)^{2}-1\ ,\ \textrm{soit}\ \varepsilon=1
\end{equation*} 
puisque $|\alpha+\beta|<1$ (Thm. 10.3.1).

(La \textbf{seule} autre valeur propre possible est nulle). D'apr\`es 9.5.1 on a donc pour $G$ l'hypoth\`ese $H_{n-n'+1}^{(0,0)}$. Ceci correspond dans la d\'efinition de $H_{\varepsilon}^{(0,0)}$ aux valeurs propres enti\`eres ou \`a
\begin{equation*}
\lambda\geq n^{2}-(n-n'+1)^{2}=2n(\ell-2)-(\ell-2)^{2}\ .
\end{equation*}

Si $\ell=3$ c'est l'in\'egalit\'e cherch\'ee $\lambda\geq 2n-1$.

En g\'en\'eral on a, en posant $r=\ell-2$,

\begin{equation*}
2nr-r^{2}\geq 2n-1
\end{equation*}
si $n\geq\frac{r+1}{2}$. Or $n+1 \geq\ell=r+2$, donc $n\geq r+1$.

Consid\'erons le cas des $1$-formes. La minoration $H_{(\varepsilon)}^{(p,q)}$ avec $p+q=1$ est donn\'ee pour $H$ par
\begin{align*}
&\lambda=(n'-1)^{2}-(\alpha+\beta)^{2}\geq(\ell-2)^{2}-1\\
\mathrm{ou}\kern 1cm &\lambda=(n'-1)^{2}-(\alpha+\beta)^{2}\geq\ell^{2}-1\ ,
\end{align*}
cf. Prop.~4.5.1. Puisque $n'=\ell-1$, on a donc $\varepsilon=1$ dans le premier cas ; la seconde majoration donne une borne sup\'erieure. Le Th\'eor\`eme 9.5.1 donne alors pour $G$ l'``hypoth\`ese'' $H_{n-n'+1}^{(1,0)\ \mathrm{ou}\ (0,1)}$ soit, pour les valeurs propres exceptionnelles,
\begin{equation*}
\lambda\geq(n-1)^{2}-(n-n'+1)^{2}=(n-1)^{2}-(n-r)^{2}
\end{equation*}
soit
\begin{equation}
\lambda\geq2(r-1)n-r^{2}+1
\end{equation} 
avec $r=\ell-2$.

On veut montrer que $\lambda\geq\frac{2}{5}n-\frac{11}{25}$. Ceci ne r\'esulte pas de (10.4.1) si $\ell=3$, mais dans ce cas on peut appliquer l'argument du \pa 10.1. Si $\ell\geq 5$, l'in\'egalit\'e cherch\'ee r\'esulte de (10.4.1) si
\begin{equation*}
2(r-\frac{6}{5})n\geq r^{2}- \frac{36}{25}=(r-\frac{6}{5})(r+\frac{6}{5}).
\end{equation*} 
Or $n+1\geq r+2$, donc $2n\geq r+2>r+\frac{6}{5}$.

\medskip

\newpage

\thispagestyle{empty}

\newpage


\part{Homologie des vari\'et\'es hyperboliques}

\markboth{CHAPITRE 11. L'ESPACE HYPERBOLIQUE COMPLEXE}{11.1. MOD\`ELE DE L'HYPERBOLO\"IDE}

\chapter{L'espace hyperbolique complexe}

Ce court chapitre est une br\`eve introduction \`a l'espace hyperbolique complexe, on peut l'omettre il n'est 
pas essentiel pour la suite. Il peut n\'eanmoins \^etre utile \`a la compr\'ehension des g\'eom\'etries
plus compliqu\'ees que nous \'etudions dans les chapitres suivants.
(Pour plus de d\'etails sur l'espace hyperbolique complexe, se reporter au livre de Goldman \cite{Goldman} ou \`a l'article \cite{Epstein})

\section{Mod\`ele de l'hyperbolo\"{\i}de et mod\`ele projectif} 

Soit ${\Bbb C}^{n,1}$ l'espace vectoriel complexe de dimension $n+1$ (sur ${\Bbb C}$) constitu\'e des 
$(n+1)$-uplets $Z=(Z_1 , \ldots ,Z_{n+1} )\in  {\Bbb C}^{n+1}$ et \'equipp\'e de la forme hermitienne~:
$$\langle Z,W \rangle = Z_1 \overline{W}_1 +\cdots +Z_n \overline{W}_n -Z_{n+1} \overline{W}_{n+1} .$$
Soit $U(n,1)$ le groupe des automorphismes (unitaires) de ${\Bbb C}^{n,1}$. L'\'equation $\langle Z,Z \rangle =-1$ d\'efinit une 
hypersurface r\'eelle $H$ dans ${\Bbb C}^{n,1}$. Le groupe $U(n,1)$ agit transitivement sur $H$. 
D'un autre cot\'e, le groupe ${\Bbb S}^1 = \{ e^{i\theta} \}$ agit librement sur $H$ par $Z \mapsto e^{i\theta} Z$; on 
appelle {\it espace hyperbolique complexe de dimension $n$} \index{espace hyperbolique complexe}la base ${\Bbb H}_{{\Bbb C}}^n$ du fibr\'e principal $H$ 
avec pour groupe ${\Bbb S}^1$. Notons $\pi$ l'application canonique de $H$ dans ${\Bbb H}_{{\Bbb C}}^n$. 
Dans la suite on d\'esignera par $P_0$ le point de ${\Bbb H}_{{\Bbb C}}^n$ au-dessous de $(0, \ldots ,0,1) \in H$ 
{\it i.e.} $$P_0 = \pi (0, \ldots ,0,1).$$

L'action de $G:= SU(n,1) = \{ A \in U(n,1) \; : \; \mbox{det}A=1 \}$  sur ${\Bbb H}_{{\Bbb C}}^n$ est transitive; 
le groupe d'isotropie du point $P_0$ est $K:= S(U(n) \times U(1))$. On a alors l'identification suivante~:
$${\Bbb H}_{{\Bbb C}}^n = G/K.$$

Un vecteur $Z \in {\Bbb C}^{n,1}$ est dit {\it n\'egatif} \index{vecteur n\'egatif} (resp. {\it nul, positif}) \index{vecteur nul} 
\index{vecteur positif} si le produit hermitien $\langle Z,Z \rangle$ est n\'egatif 
(resp. nul, positif). L'espace hyperbolique complexe ${\Bbb H}_{{\Bbb C}}^n$ s'identifie au sous-espace de 
${\Bbb P}({\Bbb C}^{n,1} )$ constitu\'e des droites n\'egatives dans ${\Bbb C}^{n,1}$. L'image $PU(n,1)$ de $U(n,1)$ 
dans $PGL( {\Bbb C}^{n,1} )$ est le groupe des biholomorphismes de ${\Bbb H}_{{\Bbb C}}^n$.

\markboth{CHAPITRE 11. L'ESPACE HYPERBOLIQUE COMPLEXE}{11.2. MOD\`ELE DE LA BOULE}

\section{Mod\`ele de la boule} 

Soit ${\Bbb C}^n = {\cal M}_{n,1} ({\Bbb C})$ l'espace vectoriel complexe de dimension $n$ munit du produit 
hermitien standard
$$\langle \langle z,w \langle \rangle = {}^t \overline{w} z = z_1 \overline{w}_1 + \cdots +z_n \overline{w}_n $$
et soit $U(n)$ son groupe (compact) d'automorphismes unitaires. On peut identifier ${\Bbb H}_{{\Bbb C}}^n$ avec la 
boule unit\'e 
$${\Bbb B}^n = \{ z\in {\Bbb C}^n / \langle \langle z,z \rangle \rangle = {}^t \overline{z} z <1 \}$$
par le plongement biholomophe suivant : 
$$\begin{array}{ccc}
{\Bbb C}^n & \longrightarrow & {\Bbb P}({\Bbb C}^{n,1} ) \\
\left( \begin{array}{c}
z_1 \\
\vdots \\
z_n 
\end{array}         
\right) & \longmapsto & \left[                                     
\begin{array}{c}                                    
z_1 \\                                    
\vdots \\                                    
z_n \\                                    1                                     
\end{array}                                               
\right] .
\end{array}$$
Plongement qui envoie l'origine de ${\Bbb C}^n$ sur le point $P_0$. Dans la suite nous travaillerons dans chacun de ces 
mod\`eles; les $z$ minuscules indiqueront que l'on se place dans le mod\`ele de la boule et les $Z$ majuscules que 
l'on se place dans le mod\`ele projectif.   

Soit $g\in G$, nous notons 
$$g= \left( 
\begin{array}{cc}
A & b \\
c & d 
\end{array}           
\right)$$
o\`u $A\in {\cal M}_{n,n} ({\Bbb C})$, $b , {}^t c \in {\Bbb C}^n$ et $d\in {\Bbb C}$.
L'action de $g$ sur ${\Bbb H}_{{\Bbb C}}^n $ est donn\'ee par :
$$gz = (Az +b) (cz +d)^{-1} ,$$
pour tout $z \in {\Bbb H}_{{\Bbb C}}^n$.

\markboth{CHAPITRE 11. L'ESPACE HYPERBOLIQUE COMPLEXE}{11.3. STRUCTURE KAEHL\'ERIENNE}

\section{Structure kaehl\'erienne} 

Dans cette section, nous \'equipons l'espace ${\Bbb H}_{{\Bbb C}}^n$ d'une structure kaehl\'erienne. Rappelons 
qu'une structure kaehl\'erienne sur un vari\'et\'e est \'equivalente \`a une structure complexe $J$ et une structure 
symplectique $\omega$ compatible dans le sens que $\omega$ est une $(1,1)$-forme positive par rapport \`a $J$. 

La structure complexe sur ${\Bbb H}_{{\Bbb C}}^n$ est induite, de mani\`ere \'equivalente, par celle de 
${\Bbb P}({\Bbb C}^{n,1} )$ ou par celle de ${\Bbb C}^n$. 

Soit $Z$ un \'el\'ement non nul de ${\Bbb C}^{n,1}$. L'espace 
tangent \`a $\pi (Z) \in {\Bbb H}^n_{{\Bbb C}}$ peut \^etre identifi\'e avec
$$Z^{\perp} = \{ W / \langle Z,W \rangle =0 \}.$$ 
Cet espace est identifi\'e avec l'espace $(\lambda Z)^{\perp}$ par multiplication par $\lambda \in {\Bbb C}^*$.
Plus g\'en\'eralement, il est souvent commode d'autoriser un vecteur quelconque $W\in {\Bbb C}^{n,1}$ \`a 
repr\'esenter un vecteur tangent \`a $Z$, en le projetant sur $Z^{\perp}$. En tenant compte du fait que $Z$ n'est 
d\'etermin\'e qu'\`a un multiple scalaire pr\`es, on d\'efinit
\begin{eqnarray*}
g_Z (W ,W) & = & \frac{\langle W- \frac{\langle W,Z \rangle Z}{\langle Z,Z \rangle }, W-
\frac{\langle W,Z \rangle Z}{\langle Z,Z \rangle }\rangle}{-\langle Z,Z \rangle } \\
                   & = & \frac{\langle Z,Z \rangle \langle W,W \rangle - \langle Z,W\rangle \langle W,Z \rangle }{-\langle Z,Z \rangle^2} .
\end{eqnarray*}
On obtient ainsi une m\'etrique $g$ sur ${\Bbb H}^n_{{\Bbb C}}$ qui est hermitienne par rapport \`a $J$.
Dans le mod\`ele de la boule, en appliquant les formules ci-dessus \`a $Z=({}^t z , 1)$ et $W= (d {}^t z ,0)$, on trouve que 
la m\'etrique sur ${\Bbb H}^n_{{\Bbb C}}$ est donn\'ee par :
\begin{eqnarray*}
ds^2 & =  & \mbox{tr}[ (I - z {}^t \overline{z} )^{-1} dz (1- {}^t \overline{z} z )^{-1} d {}^t \overline{z} ] \\     
& =  & \frac{(1-\sum_i z_i \overline{z}_i )(\sum_i dz_i d\overline{z}_i )+(\sum_i \overline{z}_i dz_i )(\sum_i z_i d\overline{z}_i )}{(1-\sum_i z_i \overline{z}_i )^2 } .
\end{eqnarray*} 

On peut maintenant montrer que ${\Bbb H}_{{\Bbb C}}^n$ est une vari\'et\'e kaehl\'erienne. Ce qui revient \`a montrer 
que si l'on d\'efinit une $2$-forme $\omega$ sur ${\Bbb H}_{{\Bbb C}}^n$ par la formule 
$$\omega (X,Y) = g (X,JY) ,$$
on obtient une $2$-forme ferm\'ee. On rappelle que $J$ d\'esigne la structure complexe, elle peut se voir comme 
l'application de l'espace tangent en un point donn\'e dans lui-m\^eme par multiplication par $\sqrt{-1}$. Pour montrer 
que $\omega$ est ferm\'ee on introduit 
\begin{eqnarray}
D_n & = & 1- {}^t \overline{z} z \\
\Phi_n & = & \partial \overline{\partial} \log D_n .
\end{eqnarray}
Un calcul simple montre que 
$$\omega = \sqrt{-1} \Phi_n .$$
En particulier la forme $\omega$ est ferm\'ee et ${\Bbb H}_{{\Bbb C}}^n$ est une vari\'et\'e k\"ahl\'erienne.

\markboth{CHAPITRE 11. L'ESPACE HYPERBOLIQUE COMPLEXE}{11.4. COURBURE}

\section{Courbure} 
Muni de sa m\'etrique k\"ahl\'erienne, l'espace ${\Bbb H}^n_{{\Bbb C}}$ est un espace sym\'etrique. 
Soient $\mathfrak{g}_0$ et $\mathfrak{k}_0$ les alg\`ebres de Lie de $G$ et $K$. Soit $\mathfrak{p}_0$ le suppl\'ementaire 
orthogonal de $\mathfrak{k}_0$ dans $\mathfrak{g}_0$ par rapport \`a la forme de Killing. 
\'Etant donn\'e $z\in {\cal M}_{n,1} ({\Bbb C})$, nous notons 
$$\xi (z) = \left(
                         \begin{array}{cc}
                          0 & z \\
                          ^t \overline{z} & 0                   
                         \end{array}                                
                                           \right) .$$
Alors $\mathfrak{p}_0 = \{ \xi (z) \; : \;  z \in {\cal M}_{n,1} ({\Bbb C} ) \}$. 
Nous identifions $\mathfrak{p}_0$ avec l'espace tangent $T_0 ({\Bbb H}_{{\Bbb C}}^n )$ \`a ${\Bbb H}_{{\Bbb C}}^n$ en $0$. 

Pour $z \in {\cal M}_{n,1} ({\Bbb C})$, soit $\tau_t$ la courbe $\tau_t =(\exp t \xi (z) )0$.
L'image de $\xi (z)$ dans $T_0 ({\Bbb H}_{{\Bbb C}}^n )$ est le vecteur tangent $\dot{\tau}_0$ \`a $\tau_t$ en  $t=0$. 
Sous cette identification, la m\'etrique riemannienne est d\'ecrite par :
$$g_0 (\xi (z) , \xi (w)) = \mbox{Re}( ^t \overline{w} z).$$
De plus d'apr\`es \cite[Th\'eor\`eme 3.2 chapitre XI]{KobayashiNomizu}, pour $X,Y,U \in \mathfrak{p}$ on a~:
\begin{eqnarray} \label{courbure}
R(X,Y)U=-[[X,Y],U], 
\end{eqnarray}
o\`u $R(.,.)$ est le tenseur de courbure de ${\Bbb H}_{{\Bbb C}}^n$. 

Enfin, concluons cette section en remarquant que si $X, Y \in {\Bbb C}^{n,1}$ sont deux \'el\'ements non nuls, 
la distance $d$ entre les points qu'ils repr\'esentent dans ${\Bbb H}_{{\Bbb C}}^n$ est donn\'ee par :
$$(\cosh d)^2 = \frac{\langle X,Y \rangle \langle Y,X  \rangle}{\langle X,X\rangle \langle Y,Y \rangle} .$$
En particulier dans le mod\`ele de la boule la distance $d$ d'un point $z \in {\Bbb H}^n_{{\Bbb C}}$ \`a $0$ est donn\'ee 
par :
\begin{eqnarray}
(\cosh d)^2 & = & \frac{1}{D_n} \\
            & = & \frac{1}{1- {}^t \overline{z} z} . 
\end{eqnarray}

\markboth{CHAPITRE 11. L'ESPACE HYPERBOLIQUE COMPLEXE}{11.5. VOLUMES DES BOULES}

\section{Volume des boules}

D\'eterminons enfin comment varie le volume d'une boule g\'eod\'esique de rayon $\rho$ en fonction de $\rho$. 
Au point ${}^t (0, \ldots ,0,r) \in {\Bbb B}^n$ la forme de Kaehler 
$$\omega = \sqrt{-1} \left( \sum_{j=1}^{n-1} \frac{1}{1-r^2} dz_j \wedge d\overline{z}_j + \frac{1}{(1-r^2)^2} dz_n \wedge d\overline{z}_n \right).$$
La forme volume vaut donc :
\begin{eqnarray*}
\frac{1}{n!} \omega^n & = & (\sqrt{-1})^n \frac{1}{(1-r^2)^{n+1}} dz_1 \wedge d\overline{z}_1 \wedge \ldots \wedge dz_n \wedge d\overline{z}_n \\
                                  & = & \frac{2^n}{(1-r^2)^{n+1}} \left( \frac{\sqrt{-1}}{2} \right) ^n dz_1 \wedge d\overline{z}_1 \wedge \ldots \wedge dz_n \wedge d\overline{z}_n \\
                                  & = & \frac{2^n}{(1-r^2)^{n+1}} dx_1 \wedge dy_1 \wedge \ldots \wedge dx_n \wedge dy_n \\
                                  & = & \frac{2^n r^{2n-1}}{(1-r^2 )^{n+1}} dr d\sigma 
\end{eqnarray*}
o\`u $d\sigma$ d\'esigne la forme volume sur la sph\`ere unit\'e $r=1$.

Puisque la distance euclidienne $r$ est reli\'ee \`a la distance hyperbolique $\rho$ par 
$$r = \tanh \left( \frac{\rho}{2} \right) , $$
le volume d'une boule de rayon $\rho$ est donn\'ee par 
\begin{eqnarray*}
\mbox{vol}(B(\rho)) & = & \int_{\tanh^{-1} (r) \leq \rho } \frac{2^n r^{2n-1}}{(1-r^2 )^{n+1}} dr d\sigma \\
                               & = & 2^n \sigma_{2n-1} \int_0^{\rho} (\sinh R)^{2n-1} (\cosh R ) dR \\
                               & = & \frac{2^n \sigma_{2n-1}}{2n} (\sinh \rho )^{2n} \sim \frac{2^n \sigma_{2n-1}}{2n} e^{2n\rho} 
\end{eqnarray*}
(o\`u $\sigma_{2n-1} = 2\pi^n /n!$ est le volume euclidien de la sph\`ere unit\'e ${\Bbb S}^{2n-1} \subset {\Bbb C}^n $).

Dans le chapitre suivant nous nous int\'eressons plus g\'en\'eralement aux espaces sym\'etriques associ\'es au groupe 
$SU(p,q)$ et \`a leurs sous-espaces totalement g\'eod\'esiques. Notons que le cas de l'espace hyperbolique 
complexe est particuli\`erement important pour nous, et peut servir de guide pour la compr\'ehension de ce qui suit.

\newpage

\thispagestyle{empty}

\newpage

\markboth{CHAPITRE 12. L'ESPACE ${\cal D}_{p,q}$}{12.1. PR\'ELIMINAIRES}

\chapter{Espaces sym\'etriques associ\'es aux groupes unitaires}

\section{Pr\'eliminaires}

Soient $p\geq q$ deux entiers strictement positifs.
Dans ce chapitre nous notons $G=SU(p,q)$, $K=G \cap U(p) \times U(q) = S(U(p) \times U(q) )$ et ${\cal D}_{p,q} =G/K$, 
l'espace sym\'etrique associ\'e. Remarquons que ${\cal D}_{n,1}$ s'identifie \`a l'espace hyperbolique 
complexe. Nous r\'ealisons plus g\'en\'eralement ${\cal D}_{p,q}$ comme un domaine complexe born\'e
$${\cal D}_{p,q} = \{ Z \in M_{p,q} ({\Bbb C}) \; : \; {}^t \overline{Z}  Z < I_q \} .$$
Nous noterons g\'en\'eralement ${\cal D} = {\cal D}_{p,q}$, \`a moins que le contexte ne soit pas clair.
\'Etant donn\'e $g\in G$, on \'ecrit 
$$g = \left(
\begin{array}{cc}
A & B \\
C & D 
\end{array} \right)$$
o\`u $A\in M_{p,p} ({\Bbb C})$, $B\in M_{p,q} ({\Bbb C})$, $C\in M_{q,p} ({\Bbb C})$ et $D \in M_{q,q} ({\Bbb C})$. 
L'action de $g$ sur ${\cal D}$ est donn\'ee par 
$$gZ = (AZ+B)(CZ+D)^{-1}, \; \; Z\in {\cal D}, \; g \in G.$$
Le groupe $G$ agit transitivement sur ${\cal D}$ et le groupe d'isotropie de $0$ est $K$. Sur ${\cal D}$,
on a une m\'etrique kaehl\'erienne d\'efinie par 
$${\rm tr} ((I_p -Z {}^t \overline{Z}  )^{-1} dZ (I_q -{}^t \overline{Z}  Z )^{-1} d {}^t \overline{Z} ).$$
L'espace ${\cal D}$ \'equip\'e de la m\'etrique riemannienne correspondante est un espace sym\'etrique hermitien.

Soit $\mathfrak{g}_0$ (resp. $\mathfrak{k}_0$) l'alg\`ebre de Lie de $G$ (resp. $K$). Soit $\mathfrak{p}_0$ le suppl\'ementaire
orthogonal de $\mathfrak{k}_0$ dans $\mathfrak{g}_0$ par rapport \`a la forme de Killing. Les alg\`ebres de Lie 
$\mathfrak{g}_0$ et $\mathfrak{k}_0$ ainsi que la forme de Killing ont \'et\'e d\'ecrites au Chapitre 4 dans le cas $q=1$. Rappelons 
que si, pour $Z \in M_{p,q} ({\Bbb C})$, nous notons 
$$\xi (Z) = \left( 
\begin{array}{cc}
0 & Z \\
{}^t \overline{Z} & 0 
\end{array}
\right) ,$$
alors $\mathfrak{p}_0 = \{ \xi (Z) \; : \; Z\in M_{p,q} ({\Bbb C}) \}$. La forme de 
Killing induit sur $\mathfrak{p}_0$ le produit scalaire Re$({\rm tr}(Z {}^t \overline{W} ))$, o\`u Re d\'esigne la partie r\'eelle 
d'un nombre complexe. 

Nous identifions $\mathfrak{p}_0$ avec l'espace tangent
$T_0 ({\cal D} )$ \`a ${\cal D}$ en $0$. Pour $Z \in M_{p,q} ({\Bbb C})$, soit $\tau_t$ la courbe 
$\tau_t = (\exp t\xi (Z) ) 0$. L'image de $\xi (Z)$ dans $T_0 ({\cal D} )$ est le vecteur tangent $\dot{\tau}_0$ 
\`a la courbe $\tau_t$ en $t=0$. Sous cette identification, la m\'etrique riemannienne $g$ de ${\cal D}$ est 
induite par la forme de Killing~:
$$g_0 (\xi (Z) , \xi (W)) = \mbox{Re} ({\rm tr}(Z {}^t \overline{W} )) .$$

\markboth{CHAPITRE 12. L'ESPACE ${\cal D}_{p,q}$}{12.2. SOUS-ESPACES TOTALEMENT G\'EOD\'ESIQUES}

\section{Sous-espaces totalement g\'eod\'esiques}

Si $v\in {\Bbb C}^n$, o\`u $n=p+q$, nous d\'ecomposons $v$ en 
$$v= \left( 
\begin{array}{c}
v_+ \\
v_-
\end{array} \right), \; v_+ \in {\Bbb C}^p , \; v_- \in {\Bbb C}^q .$$
Soit $g = \left(
\begin{array}{cc}
A & B \\
C & D 
\end{array} \right) \in G$ et $Z \in {\cal D}$, on introduit les facteurs d'automorphie~:
\begin{eqnarray}
J(g,Z) & = & \left( 
\begin{array}{cc}
l(g,Z) & 0 \\
0 & j(g,Z ) 
\end{array}
\right), \\
j(g,Z)  & = & CZ+D, \\
l(g,Z) & = & A -(gZ)C .
\end{eqnarray}
L'action de $g$ sur ${\cal D}$ peut alors prendre la forme suivante~:
$$g \left( \begin{array}{c} 
Z \\
I_q
\end{array} \right) = \left( \begin{array}{c}
gZ \\
I_q 
\end{array} \right) j(g,Z ) , \; Z\in {\cal D} .$$

\medskip

Dans la suite, $n=p+q$ et $Q$ est la forme hermitienne sur ${\Bbb C}^n$ de matrice (elle aussi not\'ee $Q$)~: 
$\left( \begin{array}{cc}
I_p & 0 \\
0 & -I_q 
\end{array} \right)$.

Soit $V$ un sous-espace complexe de ${\Bbb C}^n$ positif par rapport \`a $Q$. On associe \`a un tel espace 
un sous-groupe $G_V$ de $G$ et une sous-vari\'et\'e ${\cal D}_V$ de ${\cal D}$ d\'efinis par~:
$$G_V = \{ g \in G \; : \; g \mbox{ laisse invariant le sous-espace } V \},$$
$${\cal D}_V = \{ Z \in {\cal D} \; : \; {}^t \overline{Z} v_+ = v_- \mbox{ pour tout } v\in V \}.$$

\begin{lem} \label{sous-espace}
La sous-vari\'et\'e ${\cal D}_V$ et le sous-groupe $G_V$ ont les propri\'et\'es suivantes.
\begin{enumerate}
\item Pour tout $g\in G$, $g {\cal D}_V = {\cal D}_{gV}$.
\item Le groupe $G_V$ agit transitivement sur ${\cal D}_V$.
\item La sous-vari\'et\'e ${\cal D}_V$ est un sous-espace sym\'etrique totalement g\'eod\'esique de dimension 
complexe $(p-r)q$, o\`u $r={\rm dim}_{{\Bbb C}} V$. En tant qu'espace sym\'etrique ${\cal D}_V$ est 
isomorphe \`a ${\cal D}_{p-r,q}$.
\end{enumerate}
\end{lem}
{\it D\'emonstration.} Il d\'ecoule facilement des d\'efinitions que 
$${\cal D}_V = \left\{ Z \in {\cal D} \; : \; {}^t \overline{v} Q \left( 
\begin{array}{c}
Z \\
I_q 
\end{array}
\right) = 0, \mbox{ pour tout } v \in V \right\}.$$
Alors si $g\in G$, $Z \in {\cal D}_V$ et $v \in V$, on v\'erifie facilement que 
$$0 = {}^t \overline{v}  Q \left(
\begin{array}{c}
Z \\
I_q 
\end{array} \right) = {}^t \overline{v}  {}^t \overline{g}  Q g \left(
\begin{array}{c}
Z \\
I_q 
\end{array} \right) = {}^t \overline{gv} Q \left( 
\begin{array}{c}
gZ \\
I_q
\end{array} \right) j(g,Z).$$
Donc ${}^t \overline{v}  {}^t \overline{g}  Q g \left(
\begin{array}{c}
Z \\
I_q 
\end{array} \right) = {}^t \overline{gv} Q \left( 
\begin{array}{c}
gZ \\
I_q
\end{array} \right) =0$ et le premier point est d\'emontr\'e.

Puis, il d\'ecoule du Th\'eor\`eme de Witt et du premier point que l'on peut supposer $V={\Bbb C}^r$. Dans 
ce cas,
$${\cal D}_V = \left\{ \left( 
\begin{array}{c}
0 \\
W
\end{array} \right) \; : \; W \in M_{p-r,q} ({\Bbb C}), \; {}^t \overline{W} W < I_q \right\}$$
et 
$$G_V = \left\{ \left( 
\begin{array}{cc}
u & 0 \\
0 & h 
\end{array} \right) \in G \; : \; h \in U(p-r,q), \; u\in U(r) \right\} .$$
Et les points 2. et 3. du Lemme \ref{sous-espace} s'en d\'eduisent facilement.

\medskip

Soit $e_1 , \ldots , e_n$ la base standard de ${\Bbb C}^n$. Fixons un entier $1\leq r <p$ et soit $V$ le sous-espace 
engendr\'e par $e_{p-r+1} , \ldots , e_p$. Nous \'etudions maintenant la fonction distance $d(Z, {\cal D}_V )$ d'un \'el\'ement $Z \in {\cal D}$
\`a ${\cal D}_V$. \'Etant donn\'e $Z \in {\cal D}$, nous d\'ecomposons $Z$ en 
$$Z = \left( 
\begin{array}{c}
Z_1 \\
Z_2
\end{array} \right),$$
o\`u $Z_1 \in M_{p-r,q} ({\Bbb C})$ et $Z_2 \in M_{r,q} ({\Bbb C})$. Le sous-espace ${\cal D}_V$ est alors 
donn\'e par 
$${\cal D}_V = \{ Z \in {\cal D} \; : \; Z_2 =0 \}.$$
Un \'el\'ement $g\in G_V$ s'\'ecrit comme matrice par blocs
$$g = \left( 
\begin{array}{ccc}
A_1 & 0 & B_1 \\
0 & u & 0  \\
C_1 & 0 & D_1  
\end{array} \right).$$
Notons ${\cal D}_1$ l'espace 
$${\cal D}_1 = \{ W \in M_{p-r,q} ({\Bbb C}) \; : \; {}^t \overline{W} W < I_q \} .$$
Le groupe $G_V$ agit transitivement sur ${\cal D}_1$ par~:
$$gW
= (A_1 W + B_1 ) (C_1 W + D_1 )^{-1}, \; g\in G_V, \; W \in {\cal D}_1.$$
L'action de $G_V$ sur ${\cal D}$ s'\'ecrit~:
\begin{eqnarray} \label{action}
gZ = \left(
\begin{array}{c}
gZ_1 \\
uZ_2 j(g,Z)^{-1} 
\end{array} \right) , \; g \in G_V , \; Z \in {\cal D} .
\end{eqnarray}
D'apr\`es (\ref{action}), il existe un \'el\'ement $g\in G_V$ tel que $gZ=Z'$ avec $Z_1 ' =0$. La m\'etrique 
riemannienne de ${\cal D}$ \'etant $G$-invariante, $d(Z, {\cal D}_V )= d(Z' , {\cal D}_V) = d(0,Z')$. Il nous 
suffit donc d'\'etudier la fonction distance $d(0,Z)$ de $0$ \`a $Z$.

\begin{lem} \label{distance}
Soit $Z \in {\cal D}$, $d= d(0,Z)$ et $m$ le rang de $Z {}^t \overline{Z}$. Alors,
\begin{enumerate}
\item si $m=1$, $\cosh^2 d = (\det (I_p -Z{}^t \overline{Z} ))^{-1}$,
\item et en g\'en\'eral,
$$\frac{1}{2^m} e^d \leq (\det (I_p -Z {}^t \overline{Z} ))^{-1} \leq e^{\sqrt{m} d}.$$
\end{enumerate}
\end{lem}
{\it D\'emonstration.} Il existe un unique $Y \in M_{p,q} ({\Bbb C})$ v\'erifiant 
\begin{eqnarray} \label{expo}
\exp (\xi (Y)) 0 = Z .
\end{eqnarray}
La courbe $\exp (t \xi (Y))0$, $0 \leq t \leq 1$, est une g\'eod\'esique joignant $0$ \`a $Z$. On a donc
$$d^2 = {\rm tr} (Y {}^t \overline{Y}  ).$$
Soit $A$ (resp. $B$) une matrice hermitienne positive v\'erifiant 
$$A^2 = Y {}^t \overline{Y}  \  ({\rm resp.} \  B^2 = {}^t \overline{Y} Y ) .$$
Il d\'ecoule des d\'efinitions que 
$$\exp (\xi (Y)) = \left(
\begin{array}{cc}
\cosh A & \sum_{k=0}^{\infty} \frac{A^{2k}}{(2k+1)!} Y \\
\sum_{k=0}^{\infty} \frac{B^{2k}}{(2k+1)!} {}^t \overline{Y} & \cosh B
\end{array} \right) .$$
Puisque $A^{2k} Y= Y B^{2k}$, on d\'eduit de l'expression ci-dessus et de (\ref{expo}) que 
$Z {}^t \overline{Z} = \tanh^2 (A)$ et donc que~:
\begin{eqnarray} \label{expA}
e^A = \frac{I_p + \sqrt{Z {}^t \overline{Z} }}{(I_p -Z{}^t\overline{Z} )^{1/2}} .
\end{eqnarray}
Il d\'ecoule facilement du fait que $A$ est de rang $m$ que
$$d \leq {\rm tr} (A) \leq \sqrt{m} d .$$
Puisque $Z{}^t \overline{Z} < I_p$, si l'on applique le d\'eterminant \`a (\ref{expA}), on obtient le point 2. du Lemme. Si
$m=1$, $d={\rm tr} (A)$ et le point d\'ecoule encore de (\ref{expA}).

\medskip

\begin{lem} \label{h}
Soient $h_Z = (I_q - {}^t Z_1 \overline{Z}_1 )^{-1}$ et $\tilde{h}_Z =(I_q - {}^t Z \overline{Z} )^{-1}$. Alors, 
$$\tilde{h}_{gZ} = \overline{j(g,Z)} \tilde{h}_Z  {}^t j(g,Z) , \mbox{ pour tout } g\in G,$$
$$h_{gZ} = \overline{j(g,Z)} h_Z {}^t j(g,Z) , \mbox{ pour tout } g\in G_V .$$
\end{lem}
{\it D\'emonstration.} Soit $l_Z =(I_q - {}^t Z \overline{Z} )$. On a~:
\begin{eqnarray*}
l_Z & = & -( {}^t Z \   I_q) Q \left( 
\begin{array}{c}
\overline{Z} \\
I_p 
\end{array} \right) \\
& = & -({}^t Z \   I_q) {}^t g Q \overline{g} \left( 
\begin{array}{c}
\overline{Z} \\
I_q 
\end{array} \right) \\
& = & {}^t j(g,Z) l_{gZ} \overline{j(g,Z)} .
\end{eqnarray*}
Puisque $\tilde{h}_Z = l_Z^{-1}$, on obtient la premi\`ere propri\'et\'e annonc\'ee. Mais 
$h_Z = \tilde{h}_{\left(
\begin{array}{c}
Z_1 \\
0
\end{array} \right) }$ et pour tout $g\in G_V$, $j(g,Z) = j \left( g, \left( 
\begin{array}{c}
Z_1 \\
0
\end{array} \right) \right)$, d'o\`u la seconde propri\'et\'e annonc\'ee.

\bigskip

Pour $Z \in {\cal D}$, on introduit les fonctions $A$ et $B$ sur ${\cal D}$ d\'efinies par 
\begin{eqnarray} \label{A et B}
A & = & \mbox{det} (I_q - {}^t Z \overline{Z} ) \\
B & = & \mbox{det} (I_q - {}^t Z_1 \overline{Z}_1 ).
\end{eqnarray}
La fonction $B$ est obtenue en restreignant la fonction $A$ \`a ${\cal D}_V$, puis en l'\'etendant \`a ${\cal D}$
tout entier de fa\c{c}on constante dans la direction de $Z_2$.

\begin{lem} \label{B sur A}
La fonction $\frac{B}{A}$ est $G_V$-invariante.
\end{lem}
{\it D\'emonstration.} Cela d\'ecoule imm\'ediatement du Lemme \ref{h}.

\medskip

Nous pouvons maintenant estimer la fonction $d(Z, {\cal D}_V )$.

\begin{prop} \label{distance2}
Soient $Z \in {\cal D}$ et $m$ le rang de la matrice $Z_2 {}^t \overline{Z}_2$.
\begin{enumerate}
\item Si $m=1$, $(\cosh d(Z,{\cal D}_V ))^2 = \frac{B}{A}$.
\item En g\'en\'eral, on a~: 
$$4^m \left( \frac{B}{A} \right) \geq e^{2 d(Z,{\cal D}_V )} , $$
et 
$$e^{2 \sqrt{m} d(Z,{\cal D}_V )} \geq \frac{B}{A} .$$
\end{enumerate}
\end{prop}
{\it D\'emonstration.} Les fonctions $\frac{B}{A}$ et $d(.,{\cal D}_V )$ sont toutes deux $G_V$-invariantes. On a 
vu, cf. (\ref{action}), que l'on pouvait se ramener \`a ce que $Z_1 =0$ et donc $d(Z,{\cal D}_V )=d(0,Z_2) $.
Mais alors, $\frac{B}{A} = (\mbox{det}(I_q - Z_2 {}^t \overline{Z}_2 ))^{-1}$. Et la Proposition d\'ecoule alors 
du Lemme \ref{distance}.

\medskip

Nous aurons \'egalement besoin dans la suite des expressions suivantes.

\begin{lem} \label{outils}
On a l'\'egalit\'e 
$$\frac{B}{A} = \det \{ I_r + \overline{Z}_2 (I_q - {}^t Z \overline{Z} )^{-1} {}^t Z_2 \} $$
et l'in\'egalit\'e
$$1+ r^{-1} {\rm tr} \left( \overline{Z}_2 (I_q - {}^t Z \overline{Z}^{-1} {}^t Z_2 \right) \geq \left( \frac{B}{A} \right)^{1/r} .$$
\end{lem}
{\it D\'emonstration.} Remarquons d'abord que 
\begin{eqnarray*}
I_q - {}^t Z_1 \overline{Z}_1 & = & I_q - {}^t Z \overline{Z} + {}^t Z_2 \overline{Z}_2 \\
& = & (I_q - {}^t Z \overline{Z} )^{1/2} \left\{ I_q + (I_q - {}^t Z \overline{Z} )^{-1/2} {}^t Z_2 \overline{Z}_2 (I_q - {}^t Z \overline{Z} )^{-1/2} 
\right\} (I_q - {}^t Z \overline{Z} )^{1/2} .
\end{eqnarray*}
On en d\'eduit que 
\begin{eqnarray*}
\frac{B}{A} & = & \det \left\{ I_q + (I_q - {}^t Z \overline{Z} )^{-1/2} {}^t Z_2 \overline{Z}_2 (I_q - {}^t Z \overline{Z} )^{-1/2} \right\} \\
                 & = & \det \left\{ I_r + \overline{Z}_2 (I_q -{}^t Z \overline{Z} )^{-1} {}^t Z_2 \right\} .
\end{eqnarray*}
Ce qui d\'emontre la premi\`ere partie du Lemme. Remarquons maintenant que la matrice $Z_2 (I_q - {}^t Z \overline{Z} )^{-1} {}^t Z_2$
est positive. Notons $\lambda_1 , \ldots , \lambda_r$ ses valeurs propres (r\'eelles positives). Alors,
\begin{eqnarray*}
1+\frac{1}{r} {\rm tr} \left( \overline{Z}_2 (I_q - {}^t Z \overline{Z} )^{-1} {}^t Z_2 \right) & = &  
\frac{(1+\lambda_1 ) + \ldots + (1+\lambda_r )}{r} \\
& \geq & \left\{ (1+\lambda_1 ) \ldots (1+\lambda_r ) \right\}^{1/r} \\
& = & \det \left\{ I_r + \overline{Z}_2 (I_q - {}^t Z \overline{Z} )^{-1} {}^t Z_2 \right\}^{1/r} \\
& = & \left( \frac{B}{A} \right)^{1/r} .
\end{eqnarray*}

\medskip

Remarquons qu'\`a l'aide de la th\'eorie g\'en\'erale des espaces sym\'etriques (cf. \cite{Helgason}, \cite{KobayashiNomizu}), on peut 
montrer que les seuls sous-espaces totalement g\'eod\'esiques de l'espace hyperbolique complexe ${\Bbb H}_{{\Bbb C}}^n = {\cal D}_{n,1}$ 
sont soit des sous-espaces ${\cal D}_V$ comme ci-dessus, soit des sous-vari\'et\'es totalement g\'eod\'esiques totalement r\'eelles. On ne s'occupera 
ici que des premiers. En ce qui concerne les espaces sym\'etriques ${\cal D}_{p,q}$, il existe en g\'en\'eral d'autres 
sous-espaces totalement g\'eod\'esiques holomorphes que ceux consid\'er\'es ci-dessus.

\medskip

\markboth{CHAPITRE 12. L'ESPACE ${\cal D}_{p,q}$}{12.3. CROISSANCE DU VOLUME}

\section{Croissance du volume}

\'Examinons maintenant les champs de Jacobi \'emanant de ${\cal D}_V$.(Une r\'ef\'erence g\'en\'erale pour les champs
de Jacobi est \cite{Sakai}.)

\begin{lem} \label{champs de Jacobi}
Soient $Z\in {\cal D}_V$, $T_Z ({\cal D}_V )$ l'espace tangent \`a ${\cal D}_V$ en $Z$ et $T_Z ({\cal D}_V )^{\perp}$
le suppl\'ementaire orthogonal de $T_Z ({\cal D}_V )$ dans $T_Z ({\cal D})$. Soit $Y$ un vecteur dans 
$T_Z ({\cal D}_V )^{\perp}$ avec $g_Z (Y,Y)=1$. Alors, 
\begin{enumerate}
\item les espaces $T_Z ({\cal D}_V )$ et $T_Z ({\cal D}_V )^{\perp}$ sont invariants sous 
l'application $R(.,Y)Y$, 
\item il existe $\lambda_1 \geq \lambda_2 \geq \ldots \geq \lambda_l \geq 0$, $l=\max \{ r,q \}$ tels que 
\begin{itemize}
\item $\lambda_1^2 + \ldots + \lambda_l^2 =1$,
\item $\lambda_i =0$, si $i>\min \{ r, q \}$,
\item l'op\'erateur $R(.,Y)Y_{|{T_Z ({\cal D}_V )}}$ a pour valeurs propres 
$$\begin{array}{cccc}
\underbrace{-\lambda_1^2 , \ldots , -\lambda_1^2} , & \underbrace{-\lambda_2^2 , \ldots , -\lambda_2^2 }, &  \ldots  &, \underbrace{-\lambda_q^2 , \ldots , - \lambda_q^2 } ,\\
2(p-r) & 2(p-r) & & 2(p-r) 
\end{array}$$
\item l'op\'erateur $R(.,Y)Y_{|{T_Z ({\cal D}_V )^{\perp} }}$ a pour valeurs propres 
$$-(\lambda_i -\lambda_j )^2 , \; -(\lambda_i +\lambda_j )^2 , \; 1\leq i \leq r , \; 1\leq j \leq q .$$
\end{itemize}
\end{enumerate}
\end{lem}
{\it D\'emonstration}. D'apr\`es (\ref{action}), on peut supposer que $Z=0$ et
$Y=\xi \left(
         \begin{array}{c}
          0 \\
          M
          \end{array}  \right) $. 
D'apr\`es \cite[Theorem 3.2, Chap. XI]{KobayashiNomizu}, \'etant donn\'e $X \in T_0 ({\cal D} )$, le tenseur de 
courbure est donn\'e par~:
\begin{eqnarray} \label{courbure}
R(X,Y)Y= -[[X,Y],Y].
\end{eqnarray}
Si $X= \xi \left( 
\begin{array}{c}
N \\
0
\end{array} \right) \in T_0 ({\cal D}_V )$,
un calcul simple donne alors
\begin{eqnarray} \label{cas1}
R(X,Y)Y = \xi \left(
                  \begin{array}{c}
                       -N {}^t \overline{M} M \\
                        0
                  \end{array} \right) .
\end{eqnarray}
Si maintenant $X= \xi \left(
                          \begin{array}{c}
                           0 \\
                           L
                           \end{array} \right) \in T_0 ({\cal D}_V )^{\perp} $, un autre calcul simple 
\`a l'aide de (\ref{courbure}) donne
\begin{eqnarray} \label{cas2}
R(X,Y)Y = \xi \left(
                  \begin{array}{c}
                  0 \\
                  -L {}^t \overline{M} M +2M {}^t \overline{L} M - M {}^t \overline{M} L
                   \end{array}  \right) .
\end{eqnarray}
Il d\'ecoule de (\ref{cas1}) et (\ref{cas2}) que les espaces $T_0 ({\cal D}_V )$ et $T_0 ({\cal D}_V )^{\perp}$ sont invariants
sous l'application $R(.,Y)Y$. Remarquons que l'on peut toujours supposer que $M$ est 
de la forme 
\begin{eqnarray*}
M & = & \left( 
\begin{array}{c}
\begin{array}{ccc}
\lambda_1 & \ldots & 0 \\
\vdots & \ddots & \vdots \\
0 & \ldots & \lambda_q 
\end{array} \\
0 
\end{array} \right) \; \mbox{ si } r>q, \\
M & = & \left( 
\begin{array}{ccc}
\lambda_1 & \ldots & 0 \\
\vdots & \ddots & \vdots \\
0 & \ldots & \lambda_q 
\end{array} \right) \; \mbox{ si } r=q, \\
M & = & \left( 
\begin{array}{cc}
\begin{array}{ccc}
\lambda_1 & \ldots & 0 \\
\vdots & \ddots & \vdots \\
0 & \ldots & \lambda_r 
\end{array} & 0 
\end{array} \right) \; \mbox{ si } r<q,
\end{eqnarray*}
avec $\lambda_1 \geq \ldots \geq \lambda_l \geq 0$, $l=\max \{r,q \}$ et $\lambda_i =0$ si $i> \min \{ r,q \}$.
Puisque $\lambda_1^2 + \ldots + \lambda_l^2 = {\rm tr} (M {}^t \overline{M}  ) =g_0 (Y,Y) = 1$, le deuxi\`eme 
point du Lemme \ref{champs de Jacobi} d\'ecoule alors des formules (\ref{cas1}) et (\ref{cas2}).

\bigskip

Soit $\tau$ une g\'eod\'esique perpendiculaire \`a ${\cal D}_V$. Nous pouvons maintenant \'etudier les 
champs de Jacobi le long de $\tau$. Soit $\tau =\tau_t$, o\`u $t$ est la longueur d'arc de ${\cal D}_V$ \`a 
$\tau_t$ et $Y= \dot{\tau}_0$. Dans la suite nous d\'ecrivons les champs de Jacobi $X=X(t)$ le long de $\tau$ v\'erifiant 
\begin{eqnarray} \label{conditions initiales}
X(0) \in T_{\tau_0} ({\cal D}_V ) \mbox{  et  } \nabla_Y X \in T_{\tau_0} ({\cal D}_V )^{\perp} .
\end{eqnarray}

L'\'equation de Jacobi est donn\'ee par 
$$\nabla_{\tau_t}^2 X +R(X,\dot{\tau}_t )\dot{\tau}_t =0 .$$

D'apr\`es le Lemme \ref{champs de Jacobi}, les valeurs propres de $R(.,Y)Y$ sont n\'egatives.
Soit $X_0 \in T_{\tau_0} ({\cal D}_V )$ un vecteur propre de $R(.,Y)Y$ pour la valeur propre 
$-\lambda^2$ ($\lambda \geq 0$). Soit $X_t$ le transport parall\`ele de $X_0$ le long de $\tau$. On peut v\'erifier que 
\begin{eqnarray} \label{chp de Jacobi1}
X(t) & = & (\cosh \lambda t) X_t 
\end{eqnarray}
est un champ de Jacobi le long de $\tau$ v\'erifiant (\ref{conditions initiales}). Soit 
$L_0 \in T_{\tau_0 }({\cal D}_V )^{\perp} $ un vecteur propre de $R(.,Y)Y$ pour la valeur propre 
$-\lambda^2$ ($\lambda \geq 0$). Soit $L_t$ le transport parall\`ele de $L_0$ le long de $\tau$. On peut v\'erifier que 
\begin{eqnarray} \label{chp de Jacobi2}
L(t) & = & \left\{ \begin{array}{cc}
                   (\sinh \lambda t) L_t & \mbox{ si } \lambda \neq 0, \\
                   t L_t                 & \mbox{ sinon}
                   \end{array} \right.
\end{eqnarray}
est un champ de Jacobi le long de $\tau$ v\'erifiant (\ref{conditions initiales}). L'espace des champs de Jacobi v\'erifiant 
(\ref{conditions initiales}) a pour dimension $2pq$. Cet espace est engendr\'e par les champs de Jacobi construits 
ci-dessus.

\medskip

L'espace  ${\cal D} - {\cal D}_V$ se d\'ecompose en un produit :
$${\cal F} \times ]0, +\infty [ $$
o\`u ${\cal F}$ est l'hypersurface de ${\cal D}$ constitu\'ee des points \`a distance $1$ de 
${\cal D}_V$ et o\`u nous identifions un point $(W,t) \in {\cal F} \times ]0, +\infty [$ avec le point 
$Z \in {\cal D}$ \`a distance $t$ de ${\cal D}_V$ et tel que la g\'eod\'esique passant par 
$Z$ et $W$ soit perpendiculaire \`a ${\cal D}_V$.

Si $S$ est un sous-ensemble mesurable de ${\cal F}$, nous noterons $\omega (t,S)$ le volume
de $\{ Z \in {\cal F} \times \{ t \} \; : \; Z_1 \in S \}$. D'apr\`es (\ref{action}), $\omega (t,S) = \omega (t ,gS)$ pour 
tout $g\in G_V$. On peut donc voir $\omega (t,S)$, pour chaque $t$, comme une mesure invariante sur 
${\cal D}_V$. Il existe alors une fonction $f(t)$ telle que $\omega (t,S) = f(t) \mbox{vol} (S)$.
 
\begin{lem} \label{croissance du volume}
Il existe une constante $c$ telle que~:
\begin{enumerate}
\item si $r=1$, 
$$\omega (t,S) = c {\rm vol}(S) (\sinh 2t) (\sinh t)^{2(q-1)} (\cosh t)^{2(p-1)}   ,$$
\item en g\'en\'eral,
$$\omega (t,S) \leq c {\rm vol}(S) (1+t^{pq} ) e^{2(p+q-1)\sqrt{m} t)}  ,$$
o\`u $m=\min \{ r,q\}$.
\end{enumerate}
\end{lem}
{\it D\'emonstration.} Si l'on \'ecrit $\omega (t,S) = f(t) \mbox{vol} (S)$, il nous faut estimer $f(t)$, par exemple en la comparant
\`a la constante $f(1)$. Soit $x_s$ une courbe dans ${\cal F}$. On note $x_s^t$ le point $(x_s , t)$, $x_{(s)}^t$ la 
courbe \`a $s$ fix\'e et $x_s^{(t)}$ la courbe \`a $t$ fix\'e. Les courbes $x_{(s)}^t$ sont des g\'eod\'esiques et 
$\dot{x}_s^{(t)}$ est un champ de Jacobi le long de cette g\'eod\'esique qui v\'erifie (\ref{conditions initiales}).
Mais d'apr\`es le Lemme \ref{champs de Jacobi} et (\ref{chp de Jacobi1}), (\ref{chp de Jacobi2}), 
l'espace $T_{x_s} ({\cal  F})$ admet un base orthonorm\'ee r\'eelle
$$X_1 , \ldots , X_{2(p-r)} , X_{2(p-r)+1} , \ldots , X_{4(p-r)}, \ldots \ldots , X_{2(p-r)q} ,$$
$$Y_1 , \ldots , Y_{2rq-1}$$
et il existe des r\'eels 
$$\lambda_1 \geq \ldots \geq \lambda_l \geq 0 , \; l=\max \{ r, q\}$$
v\'erifiant~:
\begin{enumerate}
\item $\lambda_1^2 +\ldots + \lambda_l^2 =1$,
\item $\lambda_i =0$ pour tout $i> \min \{ r,q\}$,
\item au point $(x_s , t)$,
\begin{eqnarray*}
\begin{array}{l}
||X_1 ||= \ldots = ||X_{2(p-r)} || =\frac{\cosh \lambda_1 t}{\cosh \lambda_1} , \\
\  \  \  \vdots \\
||X_{2(p-r)(q-1)+1}|| = \ldots = ||X_{2(p-r)q}|| = \frac{\cosh \lambda_q t}{\cosh \lambda_q} ,
\end{array}
\end{eqnarray*}
\item au point $(x_s , t)$, l'ensemble des $||Y_j ||$ pour $1\leq j \leq 2rq-1$ (compt\'ees avec multiplicit\'es) co\"{\i}ncide,
\`a une permutation pr\`es, avec l'ensemble des $a_{ij}$ pour $1\leq i \leq r$, $1\leq j \leq q$ et des $b_{ij}$ pour 
$1\leq i\leq r$, $1\leq j\leq q$ et $(i,j) \neq (1,1)$ tels que
$$a_{ij} = \left\{
\begin{array}{lc}
\frac{\sinh (\lambda_i +\lambda_j )t}{\sinh (\lambda_i +\lambda_j )}, & \mbox{ si } \lambda_i + \lambda_j \neq 0, \\
t, & \mbox{ si } \lambda_i + \lambda_j =0, 
\end{array} \right. \; \mbox{ et }  \;
b_{ij} = \left\{
\begin{array}{lc}
\frac{\sinh |\lambda_i -\lambda_j |t}{\sinh |\lambda_i -\lambda_j |}, & \mbox{ si } \lambda_i \neq \lambda_j , \\
t, & \mbox{ si } \lambda_i = \lambda_j . 
\end{array} \right.$$
\end{enumerate}
On en d\'eduit alors facilement que 
$$f(t) = f(1) \frac{\sinh 2t \sinh^{2(q-1)} t \cosh^{2(p-1)} t}{\sinh 2 \sinh^{2(q-1)} 1 \cosh^{2(p-1)} 1} \  \   \mbox{  si  } r=1$$
et en g\'en\'eral qu'il existe une constante $c$ telle que
$$f(t) \leq c (1+t^{2pq}) e^{2(p+q-1)\sqrt{m} t} , $$
o\`u $m = \min \{ r,q\}$.

\bigskip

Remarquons que la fonction $\frac{B}{A}$ est plus naturelle que la fonction distance $d(.,{\cal D}_V )$. Dans 
la suite nous reprenons l'\'etude du volume \`a l'aide de la fonction $\frac{B}{A}$.

Notons d'abord que si $g\in G$ et $Z\in {\cal D}$,
$$d(gZ) = l (g,Z) dZ j(g,Z)^{-1} .$$
On sait que 
$$\det ( l (g,Z )) = \det (j(g,Z))^{-1} .$$
Donc si l'on pose
$$\{ dZ \} = \prod_{i=1}^p \prod_{j=1}^q dZ_{ij} d \overline{Z}_{ij} $$
et si $g\in G$,
$$\{ dgZ \} = |\det (j(g,Z)) |^{-2(p+q)} \{ dZ \} .$$
Puis d'apr\`es le Lemme \ref{h}, 
$$A (gZ) = |\det (j(g,Z)) |^{-2} A(Z) .$$
La forme volume invariante $dv_{{\cal D}}$ de ${\cal D}$ s'\'ecrit donc 
\begin{eqnarray} \label{forme volume de D}
dv_{{\cal D}} = (\sqrt{-1})^{pq} A^{-(p+q)} \{ dZ \} .
\end{eqnarray}

Si $\left( 
\begin{array}{c}
Z_1 \\
0
\end{array} \right) \in {\cal D}_V$, soit $F_{Z_1}$ la fibre au-dessus de ce point dans le fibr\'e ${\cal D} \rightarrow 
{\cal D}_V$. On a donc~:
$$F_{Z_1} = \{ Z\in {\cal D} \; : \; Z_1 \mbox{ fix\'e} \} .$$
Soit $g\in G_V$ l'\'el\'ement 
$$g = \left( 
\begin{array}{ccc}
(I_{p-r} - Z_1 {}^t \overline{Z}_1 )^{-1/2} & 0 & -( I_{p-r} - Z_1 {}^t \overline{Z}_1 )^{-1/2} Z_1 \\
0 & I_r & 0 \\
-(I_q - {}^t \overline{Z}_1  Z_1 )^{-1/2} {}^t \overline{Z}_1 & 0 & (I_q - {}^t \overline{Z}_1 Z_1 )^{-1/2}
\end{array} \right) .$$
Alors $g$ envoie $F_{Z_1}$ isom\'etriquement sur $F_0$, et 
$$ g \left( 
\begin{array}{c}
Z_1 \\
Z_2
\end{array} \right) = \left( 
\begin{array}{c}
0 \\
Z_2 (I_q - {}^t \overline{Z}_1 Z_1 )^{-1/2} 
\end{array} \right) .$$
Sur $F_0$, l'\'el\'ement de volume est $(\sqrt{-1})^{rq} \det (I_q - {}^t Z_2 \overline{Z}_2 )^{-(r+q)} \{ dZ_2 \}$, l'\'el\'ement de volume
sur $F_{Z_1}$ est donc
$$dv_F = (\sqrt{-1})^{rq} A^{-r} \left( \frac{B}{A} \right)^q \{ dZ_2 \} .$$
D'o\`u il d\'ecoule que 
\begin{eqnarray} \label{decomposition du volume}
dv_{{\cal D}} = \left( \frac{B}{A} \right)^{p-r} dv_{{\cal D}_V} dv_F ,
\end{eqnarray}
o\`u $dv_{{\cal D}_V} = (\sqrt{-1})^{(p-r)q} B^{-(p+q-r)} \{ dZ_1 \}$ est la forme volume invariante sur ${\cal D}_V$. 

\begin{lem} \label{formules d'integration}
On a les formules d'int\'egration~:
\begin{enumerate}
\item
$$\int_{{\cal D}} A^{s} \{ dZ \} = \left( \frac{2\pi}{\sqrt{-1}} \right)^{pq} \prod_{i=1}^p \frac{\Gamma (s+i)}{\Gamma (s+q+i)} ,$$
d\`es que Re$(s) >-1$; et 
\item 
$$\int_{\Gamma_V \backslash {\cal D}} \left( \frac{A}{B} \right)^s dv_{{\cal D}} = \left( \frac{2\pi}{\sqrt{-1}} \right)^{rq} 
\prod_{i=1}^r \frac{\Gamma (s-p-q+i)}{\Gamma (s-p+i)}  {\rm vol} (\Gamma_V \backslash {\cal D}_V ) ,$$
d\`es que ${\rm Re} (s) >p+q-1$.
\end{enumerate}
\end{lem}
{\it D\'emonstration.} On introduit tout d'abord $f(s,p,q) = \int_{{\cal D}} A^s \{ dZ \}$. 
On d\'eduit de (\ref{decomposition du volume}), avec $r=1$, la relation de r\'ecurrence
$$f(s,p,q) = f(s+1, p-1 , q)f(s,1,q) .$$
Mais,
\begin{eqnarray*}
f(s,1,q) & = & \int_{\sum_i |z_i |^2 \leq 1} (1-(|z_1 |^2 + \ldots + |z_q |^2 ))^s dz_1 d\overline{z}_1 \ldots dz_q d\overline{z}_q  \\
& = & \left( \frac{2}{\sqrt{-1}} \right)^q \int_{\sum_i (|x_i |^2 + |y_i |^2 ) \leq 1}  (1-(|x_1|^2 +|y_1|^2 + \ldots + |x_q|^2 +|y_q|^2 ))^s dx_1 dy_1 \ldots dx_q dy_q \\
& = & \frac{(2\pi)^q}{(\sqrt{-1})^q \Gamma (q)} \int_0^1 (1-t)^s t^{q-1} dt \\
& = & \frac{(2\pi )^q \Gamma (s+1)}{(\sqrt{-1})^q \Gamma (s+q+1)} \; \mbox{ d\`es que Re}(s)>-1 . 
\end{eqnarray*}
Alors le premier point du Lemme \ref{formules d'integration} d\'ecoule d'une simple r\'ecurrence.

Concernant le deuxi\`eme point, il d\'ecoule de (\ref{decomposition du volume}) que l'int\'egrale vaut
$$\int_{\Gamma_V \backslash {\cal D}} \left( \frac{A}{B} \right)^{s+r-p} dv_F dv_{{\cal D}_V} .$$
Puisque $\frac{A}{B}$ et $dv_F$ sont $G_V$-invariant, l'int\'egrale
$$\int_{F_{Z_1}} \left( \frac{A}{B} \right)^{s+r-p} dv_F$$
est ind\'ependante de $Z_1$. En $Z_1 =0$, sa valeur est 
$$\int_{{\cal D}_2} \det (I_q - {}^t Z_2 \overline{Z}_2 )^{s-p-q} \{ dZ_2 \} ,$$
o\`u ${\cal D}_2 = \{ Z_2 \; : \; {}^t \overline{Z}_2  Z_2 < I_q \}$. Le deuxi\`eme point du Lemme \ref{formules d'integration}
d\'ecoule donc directement du premier point.

\medskip

\markboth{CHAPITRE 12. L'ESPACE ${\cal D}_{p,q}$}{12.4. FONCTION DISTANCE \`A L'HYPERSURFACE}

\section{Fonction distance \`a l'hypersurface}

Soit $F$ la fonction distance g\'eod\'esique \`a la sous-vari\'et\'e ${\cal D}_V$. La fonction $Z \mapsto F(Z)$ est bien \'evidemment
lisse pour $Z\in {\cal D}-{\cal D}_V$. Nous notons $\nabla^2 F$ le hessien de $F$. 
Rappelons que le {\it hessien} \index{hessien} d'une fonction $C^2$ $F$ de ${\cal D}$ dans ${\Bbb R}$ est la seconde d\'eriv\'ee 
covariante $\nabla^2 F$ de $F$, {\it i.e.} 
$$\nabla^2 F (X,Y) = X(YF) - (\nabla_X Y)F ,$$
pour n'importe quels champs de vecteurs $X$, $Y$ sur ${\cal D}$ et o\`u $\nabla$ est la connexion de 
Levi-Civit\`a induite par la structure riemannienne de ${\cal D}$. Le hessien $\nabla^2 F$ d\'efinit donc  un
tenseur sym\'etrique de type $(0,2)$. Nous appelons  {\it valeurs propres du hessien} \index{valeurs propes du hessien} les fonctions qui \`a 
chaque point de ${\cal D}$ associe les valeurs propres de la matrice associ\'ee dans n'importe quelle base 
orthonorm\'ee de l'espace tangent \`a ${\cal D}$ au point $x$.
  
\begin{prop} \label{fonction distance a l'hypersurface}
Notons $\{ \gamma_i (Z) \}_{1 \leq i \leq 2pq}$ les valeurs propres du hessien $\nabla^{2} F$.
Si $Z \in {\cal D}$, $m= \min \{ r,q \}$ 
et $l=\max \{ r,q \}$, il existe alors $\lambda_1 \geq \lambda_2 \geq \ldots \geq \lambda_l \geq 0$ tels que 
\begin{itemize}
\item $\lambda_1^2 + \ldots + \lambda_l^2 =1$,
\item $\lambda_i =0$, pour tout $i>m$,
\item quitte \`a r\'eordonner les $\gamma_i (Z)$, 
$$
\begin{array}{l}
\gamma_1 (Z) = \lambda_1 \tanh (\lambda_1 F(Z)) , \ldots , \gamma_{2(p-r)} (Z) = \lambda_1 \tanh (\lambda_1 F(Z)) , \\
\ldots \ldots , \\
\gamma_{2(p-r)(q-1)+1} (Z) = \lambda_q \tanh (\lambda_q F(Z)) , \ldots , \gamma_{2(p-r)q} (Z) = \lambda_q \tanh (\lambda_q F(Z)) 
\end{array}
$$
et les $\gamma_k (Z)$ pour $2(p-r)q+1 \leq k \leq 2pq$, sont (\`a permutations pr\`es) les nombres 
$n_{i,j}$ et $m_{i,j}$, pour $1\leq i \leq r$, $1\leq j \leq q$, tels que
$$n_{i,j} = \left\{
\begin{array}{l}
(\lambda_i + \lambda_j ) \coth ((\lambda_i +\lambda_j ) F(Z)) \; \mbox{ si } \lambda_i + \lambda_j \neq 0 \\
\frac{1}{F(Z)} \;  \mbox{ si } \lambda_i + \lambda_j =0 
\end{array} \right. $$
et  
$$m_{i,j} = \left\{ 
\begin{array}{l}
|\lambda_i -\lambda_j | \coth (|\lambda_i - \lambda_j| F(Z)) \; \mbox{ si } \lambda_i \neq \lambda_j \\
\frac{1}{F(Z)} \; \mbox{ si } \lambda_i = \lambda_j \mbox{ et } (i,j) \neq (1,1) ,\\
0 \; \mbox{ si } (i,j)=(1,1) .
\end{array} \right. $$ 
\end{itemize}
\end{prop}
{\it D\'emonstration.} Soit toujours $\tau$ une g\'eod\'esique perpendiculaire \`a ${\cal D}_V$ avec 
$\tau =\tau_{t}$ o\`u $t$ est la longueur d'arc de ${\cal D}_V$ \`a $\tau_t$.  Soit $Y= \dot{\tau}$. 
D'apr\`es le Lemme \ref{champs de Jacobi}, il existe donc 
\begin{itemize}
\item $l$ r\'eels $\lambda_1 \geq \lambda_2 \geq 
\ldots \geq \lambda_l \geq 0$ tels que $\lambda_1^2 + \ldots + \lambda_l^2 =1$ et $\lambda_i =0$, pour tout $i>m$;
\item un champs de 
bases orthonorm\'ees le long de $\tau$~: $\{ e_{\alpha}, f_{i,j} , f_{i,j} ' \; : \; 1\leq \alpha \leq 2(p-r)q \mbox{ et } 1\leq i,j \leq 2rq \}$ tel
que pour tout entier $1 \leq \beta \leq q$ et pour tout entier $2(\beta-1)(p-r)+1 \leq \alpha \leq 2\beta(p-r)$ le vecteur 
$e_\alpha (0 ) \in T_{\tau_0 } ({\cal D})$ et soit un vecteur $(-\lambda_{\beta}^2)$-propre de $R(.,Y)Y$ et que pour toute paire d'entiers
$1\leq i,j \leq 2rq$ le vecteur $f_{i,j} (0 )$ (resp. $f_{i,j} ' (0)$) $\in T_{\tau_0} ({\cal D})^{\perp}$ et soit
un vecteur $(-(\lambda_i + \lambda_j )^2)$-propre (resp. $(-(\lambda_i -\lambda_j )^2)$-propre de $R(.,Y)Y$.
\end{itemize}
Nous supposerons de plus (ce que l'on peut bien \'evidemment faire) que le vecteur $Y$ est \'egal au vecteur 
$f_{1,1} ' (0)$. 

Alors d'apr\`es (\ref{chp de Jacobi1}) et pour tout couple d'entiers $(\beta ,\alpha )$ v\'erifiant $1\leq \beta \leq q$ et 
$2(\beta-1)(p-r)+1\leq \alpha \leq 2\beta (p-r)$,  les champs de vecteurs~:
\begin{eqnarray} \label{base1}
v_{\alpha} (t) = \cosh (\lambda_{\beta} t) e_{\alpha} (t)  
\end{eqnarray}
sont des champs de Jacobi le long de $\tau$ v\'erifiant (\ref{conditions initiales}).
Puis, d'apr\`es (\ref{chp de Jacobi2}) et pour tout couple d'entiers $(i,j)$ v\'erifiant $1\leq i \leq r$ et $1 \leq j \leq q$,
les champs de vecteurs~:
\begin{eqnarray} \label{base2}
w_{i,j} (t) = \left\{ 
\begin{array}{ll}
\sinh ((\lambda_i +\lambda_j )t) f_{i,j} (t), &  \mbox{ si } \lambda_i + \lambda_j \neq 0 \\
t f_{i,j} (t), & \mbox{ si } \lambda_i + \lambda_j =0,
\end{array} \right.
\end{eqnarray}
et
\begin{eqnarray} \label{base3} 
w_{i,j} ' (t) = \left\{
\begin{array}{ll}
\sinh (|\lambda_i - \lambda_j |t) f_{i,j} ' (t), &  \mbox{ si } \lambda_i \neq \lambda_j \\
t f_{i,j} ' (t), & \mbox{ si } \lambda_i = \lambda_j  ,\\
\end{array} \right.
\end{eqnarray}
sont des champs de Jacobi le long de $\tau$ v\'erifiant (\ref{conditions initiales}).
De plus, nous avons vu que les champs de vecteurs (\ref{base1}), (\ref{base2}) et (\ref{base3})
forment une base orthogonale de l'espace des champs de Jacobi le long de $\tau$ v\'erifiant (\ref{conditions initiales}). 
La formule de la variation seconde \cite{Sakai} nous dit alors que le Hessien 
$\nabla^2 t (= \nabla^2 F)$ se diagonalise dans la base 
$\{ e_{\alpha} \}_{1\leq \alpha \leq 2q(p-r)} \cup \{ f_{i,j} ,  f_{i,j} \}_{
\begin{array}{c}
1\leq i \leq r \\
1 \leq j \leq q
\end{array}}$. Et plus pr\'ecisemment permet de calculer par exemple
$$\nabla^2 t(y) (e_{\alpha} ,e_{\alpha} ) = \frac{d^2}{ds^2}_{|s=0} L(\tau^s ) , $$
o\`u si $\tau$ va de $x:= \tau_0 \in {\cal D}_V$ \`a $y:= \tau_{t(y)}$, $\tau^s$ d\'esigne la g\'eod\'esique minimisante
joignant ${\cal D}_V$ au point $\exp_y (s e_{\alpha} )$ et $L( \tau^s )$ sa longueur. Or, si $\alpha$ est un entier compris entre $2(\beta-1)(p-r)+1$ et 
$2\beta (p-r)$ pour un certain entier $1\leq \beta \leq q$, le champ de vecteur 
$\hat{v}_{\alpha} = \frac{v_\alpha}{\sinh ( \lambda_{\beta} t(y))}$ est un champ de Jacobi le long de $\tau$, perpendiculaire \`a $\dot{\tau}$ et 
v\'erifiant~: $\hat{v}_{\alpha} (t(y)) =e_{\alpha} (t(y))$ et (\ref{conditions initiales}). La formule de la variation seconde implique alors~:
\begin{eqnarray*}
\nabla^2 t(y) (e_{\alpha} ,e_{\alpha} ) & = & < \nabla \hat{v}_{\alpha} (t(y)), \hat{v}_{\alpha} (t(y)) > , \\
& = & \lambda_{\beta} \frac{\sinh (\lambda_{\beta} t(y))}{\cosh (\lambda_{\beta} t(y))} .
\end{eqnarray*}
De la m\^eme mani\`ere, si $(i,j)$ est un couple d'entiers v\'erifiant $1\leq i \leq r$ et $1 \leq j \leq q$, on obtient~:
\begin{eqnarray*}
\nabla^2 t(y) (f_{i,j} ,f_{i,j} )  = \left\{
\begin{array}{ll} 
(\lambda_i + \lambda_j ) \frac{\cosh ((\lambda_i +\lambda_j )t(y))}{\sinh ((\lambda_i + \lambda_j )t(y))}, & \mbox{ si } \lambda_i + \lambda_j \neq 0 \\
\frac{1}{t(y)}, & \mbox{ si } \lambda_i = \lambda_j , 
\end{array} \right.
\end{eqnarray*}
et 
\begin{eqnarray*}
\nabla^2 t(y) (f_{i,j}' , f_{i,j}' )  = \left\{
\begin{array}{ll} 
|\lambda_i - \lambda_j | \frac{\cosh (|\lambda_i - \lambda_j |t(y))}{\sinh (|\lambda_i -\lambda_j |t(y))}, & \mbox{ si } \lambda_i \neq \lambda_j  \\
\frac{1}{t(y)}, & \mbox{ si } \lambda_i = \lambda_j \mbox{ et } (i,j) \neq (1,1), \\
0, & \mbox{ si } (i,j)=(1,1). 
\end{array} \right.
\end{eqnarray*}
Ce qui conclut la d\'emonstration de la Proposition \ref{fonction distance a l'hypersurface}.

\medskip

La vari\'et\'e ${\cal D}$ est, en plus de sa structure riemannienne, naturellement \'equip\'ee d'une structure complexe
$J : T^* {\cal D} \rightarrow T^* {\cal D}$ induite par la multiplication par $\sqrt{-1}$ sur $\mathfrak{p}$.
Nous avons \'equip\'e ${\cal D}$ d'une m\'etrique hermitienne kaehl\'erienne. On associe naturellement \`a cette derni\`ere 
la $2$-forme r\'eelle ({\it i.e.} qui est de bidegr\'e $(1,1)$ et a des coefficients r\'eels dans des coordonn\'es r\'eelles)
$$-\sqrt{-1} \partial \overline{\partial} \log \det (I_q - {}^t \overline{Z} Z) .$$
On peut voir cette derni\`ere comme l'oppos\'ee de la partie imaginaire de la m\'etrique hermitienne de ${\cal D}$.
Remarquons que ce jonglage entre m\'etrique et $2$-forme diff\'erentielle est tout d'abord r\'eversible et, de mani\`ere
g\'en\'erale applicable \`a tout $2$-tenseur covariant hermitien, {\it i.e.} tout tenseur locallement de la forme 
$\sum_{\alpha , \beta} H_{\alpha  \beta} dz^{\alpha} d\overline{z}^{\beta}$ o\`u $H_{\alpha \beta}= \overline{H}_{\beta \alpha}$.
Comme d'habitude on confond la m\'etrique avec la forme diff\'erentielle dans les \'enonc\'es. Ainsi par exemple, 
si $f$ est une fonction \`a valeurs r\'eelles, l'assertion ``$\partial \overline{\partial} f$ est d\'efinie positive'' signifie 
r\'eellement que ``le tenseur hermitien $\sum_{\alpha , \beta} \frac{\partial^2 f}{\partial z^{\alpha} \partial \overline{z}^{\beta}}
dz^{\alpha} d\overline{z}^{\beta}$ associ\'e \`a $\frac{\sqrt{-1}}{2} \partial \overline{\partial} f$ est d\'efinie positif''. 

Soit $f$ une fonction r\'eelle sur la vari\'et\'e complexe ${\cal D}$; on d\'efinit sa {\it forme de Levi} \index{forme de Levi} $Lf$ par~:
$$Lf = 2\sum_{\alpha , \beta } \frac{\partial^2 f}{\partial z^{\alpha} \partial \overline{z}^{\beta}} dz^{\alpha} d\overline{z}^{\beta} ,$$
qui est juste le tenseur hermitien associ\'e \`a $\sqrt{-1} \partial \overline{\partial} f$.
Nous appelons {\it valeurs propres de la forme de Levi} \index{valeurs propres de la forme de Levi} $Lf$ les valeurs propres de la matrice hermitienne associ\'ee 
au tenseur $Lf$ dans une base orthonorm\'ee. 

Rappelons le lien bien connu entre la forme de Levi de $f$ et son hessien sur une vari\'et\'e kaehl\'erienne.
\begin{lem} \label{levi et hessien}
Soit $f$ une fonction r\'eelle sur une vari\'et\'e kaehl\'erienne et soit $X_0 = \frac12 (X - \sqrt{-1} JX)$ un vecteur 
de bidegr\'e $(1,0)$. Alors,
$$Lf (X_0 , X_0 ) = \frac12 \left( \nabla^2 f (X,X) + \nabla^2 f (JX,JX) \right) .$$
\end{lem}

\medskip

Revenons maintenant \`a notre fonction ``distance \`a ${\cal D}_V$'', la fonction $F$. Il est facile de v\'erifier que 
lorsque $r$ ou $q$ est \'egale \`a $1$ ({\it i.e.} $m=1$) le hessien de $F$ se diagonalise dans une 
base $J$-invariante. On d\'eduit alors du Lemme \ref{levi et hessien} que les valeurs propres de la forme de 
Levi de $F$ sont 
\begin{itemize}
\item si $r=1$, les valeurs propres~:  $\coth (2F(Z))$ avec multiplicit\'e $1$, $\tanh F(Z)$ avec multiplicit\'e $p-1$, 
$\coth F(Z)$ avec multiplicit\'e $q-1$ et $0$ avec multiplicit\'e $(p-1)(q-1)$;
\item si $q=1$, les valeurs propres~: $\coth (2F(Z))$ avec multiplicit\'e $1$, $\tanh F(Z)$ avec multiplicit\'e $p-r$,
$\coth F(Z)$ avec multiplicit\'e $r-1$.
\end{itemize} 

Les valeurs propres de la forme de Levi de $F$ sont, lorsque $F(Z)$ tend vers l'infini, plus proches 
les unes des autres que celles du hessien de $F$. Ceci sera important pour nous dans la suite, l'importance de cette
propri\'et\'e pour l'\'etude du spectre et de la cohomologie $L^2$ a \'et\'e mise en \'evidence pour la premi\`ere fois 
dans \cite{DonnellyFefferman} et les m\^emes raisons fondent son importance pour nous.

Lorsque $m>1$, on l'a vu la ``bonne'' fonction tenant compte de la structure complexe de ${\cal D}_V$ n'est plus 
la fonction $F$ mais la fonction $\frac{B}{A}$ (ou plut\^ot $\log \left( \frac{B}{A} \right)$). Dans la suite nous aurons besoin 
de connaitre les valeurs propres de sa forme de Levi.

\begin{prop} \label{vp de Levi}
Les valeurs propres de la forme de Levi de $\log \left( \frac{B}{A} \right)$ au point $Z \in {\cal D}$ sont la valeur propre $1$ avec multiplicit\'e
$qr$ et chacune des valeurs propres de $(I_q - {}^t \overline{Z}_1  Z_1 )^{-1/2} {}^t \overline{Z}_2  Z_2 (I_q - (\overline{Z}_1 )^t Z_1 )^{-1/2}$ avec 
multiplicit\'e $p-r$. En particulier, la valeur propre nulle intervient avec multiplicit\'e $(p-r)(q-m)$ o\`u $m \leq \min \{ q,r \}$ 
est le rang de la matrice $Z_2$. 
\end{prop}
{\it D\'emonstration.} Nous ordonnons les coordonn\'ees de $Z\in {\cal D}$ par 
$Z_{11}, \ldots , Z_{1q} ,$ $Z_{21} , \ldots , Z_{2q} , \ldots \ldots , Z_{p1} , \ldots , Z_{pq}$.  
Alors, au point $Z$, la matrice des coefficients de la m\'etrique kaehl\'erienne est la matrice 
\begin{eqnarray} \label{matrice de la metrique}
\left( 
\begin{array}{ccc}
(I_q -{}^t \overline{Z} Z )^{-1} & & \\
& \ddots & \\
& & (I_q -{}^t \overline{Z} Z)^{-1} 
\end{array} \right) \widetilde{(I_p -  \overline{Z} {}^t Z )^{-1}} ,
\end{eqnarray}
et son inverse est 
\begin{eqnarray} \label{inverse}
\left( 
\begin{array}{ccc}
(I_q - {}^t \overline{Z} Z ) & & \\
& \ddots & \\
& & (I_q -{}^t \overline{Z} Z) 
\end{array} \right) \widetilde{(I_p -\overline{Z} {}^t Z )} ,
\end{eqnarray}
o\`u \'etant donn\'e $X = (x_{ij }) \in M_{p} ({\Bbb C})$,  nous notons $\widetilde{X}$ la matrice $(x_{ij} I_q ) \in M_{pq} ({\Bbb C})$.

Calculons maintenant la matrice de $\sqrt{-1} \partial \overline{\partial} \log \left( \frac{B}{A} \right)$. Remarquons tout d'abord
que 
$$\sqrt{-1} \partial \overline{\partial} \log \left( \frac{B}{A} \right) = \sqrt{-1} \partial \overline{\partial} \log B - \sqrt{-1} 
\partial \overline{\partial} \log A .$$
La fonction $\frac{B}{A}$ est $G_V$-invariante, d'apr\`es (\ref{action}) il nous suffit donc de d\'eterminer la 
matrice de $\sqrt{-1} \partial \overline{\partial} \log \left( \frac{B}{A} \right)$ aux points $Z$ tels que $Z_1 =0$. 
Il est bien connu \cite{KobayashiNomizu} que $-\log A$ (resp. $-\log B$) est un potentiel pour la m\'etrique kaehl\'erienne de ${\cal D}$
(resp. ${\cal D}_V$). Autrement dit $-\sqrt{-1} \partial \overline{\partial} \log A$ s'identifie \`a la m\'etrique de Kaehler
de ${\cal D}$ et sa matrice est donc (\ref{matrice de la metrique}); et $-\sqrt{-1} \partial \overline{\partial} \log B$ 
s'identifie \`a la m\'etrique de Kaehler de ${\cal D}_V$ et sa matrice en $Z_1 =0$ est donc $I_{(p-r)q}$.
La matrice 
\begin{eqnarray} \label{chgt de base}
\left( 
\begin{array}{ccc}
(I_q -{}^t \overline{Z} Z )^{1/2} & & \\
& \ddots & \\
& & (I_q - {}^t \overline{Z} Z)^{1/2} 
\end{array} \right) \widetilde{(I_p -\overline{Z} {}^t Z )^{1/2}} ,
\end{eqnarray}
est hermitienne et de carr\'e la matrice (\ref{inverse}), elle r\'ealise donc un changement de base de la 
base donn\'ee par les coordonn\'ees canonique vers une base orthonorm\'ee pour la m\'etrique. 
La matrice de $\sqrt{-1} \partial \overline{\partial} \log \left( \frac{B}{A} \right)$ dans cette nouvelle base et en $Z_1 =0$
est donc 
$$\left(
\begin{array}{cccc}
{}^t \overline{Z}_2 Z_2 & & & \\
& \ddots & &  \\
& & {}^t \overline{Z}_2  Z_2 & \\
& & & I_{qr} 
\end{array} \right) .$$
Ce qui conclut la d\'emonstration de la Proposition \ref{vp de Levi}.

\medskip

Concluons cette section par un corollaire dont nous aurons besoin dans la suite du texte.

\begin{cor} \label{vp asymptotique}
Les valeurs propres de la forme de Levi de $\log \left( \frac{B}{A} \right)$ au point $Z \in {\cal D}$ sont toutes 
positives, inf\'erieures (ou \'egales) \`a $1$, et parmi celles-ci au moins $p+qr-r$ tendent vers $1$ lorsque $\frac{B}{A} (Z)$ tend vers l'infini.
\end{cor}
{\it D\'emonstration.} On peut encore se ramener au cas $Z_1 =0$. Le Corollaire \ref{vp asymptotique} d\'ecoule 
alors facilement de la Proposition \ref{vp de Levi}, puisque lorsque $\frac{B}{A} (Z)$ tend vers l'infini, $\det (I_q -
{}^t Z_2 \overline{Z}_2 )$ tend vers $0$ et donc la plus grande valeur propre de ${}^t \overline{Z}_2  Z_2$ tend vers $1$.

\medskip

\markboth{CHAPITRE 12. L'ESPACE ${\cal D}_{p,q}$}{12.5. S\'ERIES DE POINCAR\'E}

\section{S\'eries de Poincar\'e}

Soit $\phi$ une forme diff\'erentielle de degr\'e $l$ sur ${\cal D}$. Nous notons $|| \phi ||$ (resp. $|| \phi ||_0 $) la 
norme ponctuelle induite par la m\'etrique $g$ (resp. la m\'etrique euclidienne).

\begin{lem} \label{comparaison des normes}
On a les in\'egalit\'es suivantes~:
$$||\phi ||_0 \geq ||\phi || \geq A^l || \phi ||_0 ,$$
o\`u $A=\det (I_q - {}^t \overline{Z} Z )$.
\end{lem}
{\it D\'emonstration.} La forme de Kaehler de ${\cal D}$, s'\'ecrit $\kappa = {\rm tr} ((I_p -Z {}^t \overline{Z} )^{-1} dZ (I_q -{}^t \overline{Z}
Z)^{-1} d {}^t \overline{Z} )$. Il est donc imm\'ediat que 
$${\rm tr} (dZ d{}^t \overline{Z} ) \leq \kappa \leq A^{-2} {\rm tr} ( dZd{}^t \overline{Z} ).$$
Le Lemme \ref{comparaison des normes} d\'ecoule trivialement de ces derni\`eres in\'egalit\'es.

\bigskip

\begin{cor} \label{comparaison}
Soit $\phi$ une forme diff\'erentielle $G_V$-invariante de degr\'e $l$. Supposons que chaque coefficient 
de $\theta_1 \wedge \ldots \wedge \theta_l$, avec $\theta_1 , \ldots , \theta_l \in \{ dZ_{ij} , d \overline{Z}_{ij} \; : \; 
1\leq i \leq p \; 1\leq j \leq q \}$, soit born\'e en $\left(
\begin{array}{c}
0 \\
Z_2 
\end{array} \right)$. Il existe alors deux constantes $C_1$, $C_2 >0$ telles que 
$$C_1 \geq || \phi || \geq C_2 \left( \frac{A}{B} \right)^l .$$
\end{cor}
{\it D\'emonstration.} D'apr\`es (\ref{action}), il suffit de le v\'erifier en $\left(
\begin{array}{c}
0 \\
Z_2 
\end{array} \right)$. Mais en $\left(
\begin{array}{c}
0 \\
Z_2 
\end{array} \right)$, $B=1$ et le Corollaire \ref{comparaison} d\'ecoule alors du Lemme \ref{comparaison des normes}.

\bigskip

Soit $\Gamma \subset G$ un sous-groupe discret sans torsion de type fini. Alors 
$M= \Gamma \backslash {\cal D} $ est une vari\'et\'e complexe hermitienne compl\`ete orient\'ee de dimension (complexe) $pq$. 
Une telle vari\'et\'e sera appel\'ee {\it $G$-vari\'et\'e} \index{$G$-vari\'et\'e} et, lorsque $q=1$, {\it vari\'et\'e hyperbolique complexe} \index{vari\'et\'e hyperbolique complexe}. 
Soit $\Gamma_V = \Gamma \cap G_V$, et soit 
$C_V =\Gamma_V \backslash {\cal D}_V $. On obtient alors le diagramme commutatif suivant~:
$$\begin{array}{ccccc}
   & {\cal D}_V          & \hookrightarrow           & {\cal D}         & \\
   & \downarrow                         &                           & \downarrow                    & \\
C_V = & \Gamma_V \backslash {\cal D}_V  & \stackrel{i}{\rightarrow} & \Gamma \backslash {\cal D} & =M
\end{array}$$
o\`u l'application $i$ est induite par l'inclusion de ${\cal D}_V$ dans ${\cal D}$. 
En g\'en\'eral, le groupe $\Gamma_V$ est r\'eduit \`a l'identit\'e. Dans la suite nous supposons que $C_V$ est de 
volume fini (plus loin nous supposerons m\^eme que $C_V$ est compacte). Soit 
$M_V =\Gamma_V \backslash {\cal D}$.
Remarquons que la fibration naturelle
$$\pi : \left\{ 
\begin{array}{rcl}
{\cal D} & \rightarrow & {\cal D}_V \\
Z & \mapsto & \left( 
\begin{array}{c}
Z_1 \\
0
\end{array} \right) 
\end{array} \right. $$ 
induit une fibration, nous la notons \'egalement $\pi$: $M_V = \Gamma_V \backslash {\cal D} 
\rightarrow \Gamma_V \backslash {\cal D}_V =C_V$.

\medskip

La Proposition 3 de \cite{MargulisSoifer} (ou le Lemme principal de \cite{EnseignMath}) implique(nt) le lemme suivant.

\begin{lem} \label{effeuillage}
Il existe une suite $\{ \Gamma_m \}$ de sous-groupes d'indices finis dans $\Gamma$, d\'ecroissante pour l'inclusion, 
telle que 
$$\Gamma_V = \bigcap_{m\in {\Bbb N}} \Gamma_m \mbox{  et  } \Gamma_0 =\Gamma .$$
Si de plus $\Gamma$ est un sous-groupe de congruence, on peut choisir les $\Gamma_m$ de congruence.
\end{lem}

Le Lemme \ref{effeuillage} implique que lorsque $\Gamma$ est de type fini, la vari\'et\'e $M$ admet une suite 
croissante $\{ M_m \}$ de rev\^etements finis telle que la suite $\{ M_m \}$ converge uniform\'ement sur tout compact 
vers la vari\'et\'e $M_V$ (il suffit de poser $M_m =\Gamma_m \backslash {\cal D}$). Nous appelons une telle 
suite de rev\^etements finis, une {\it tour d'effeuillage autour de $C_V$} \index{tour d'effeuillage}. Dans la suite, nous supposons que $M$ poss\`ede 
une telle tour et notons $\Gamma_m$ le groupe fondamental de $M_m$.

Nous allons travailler tout au long de cette partie avec des formes diff\'erentielles sur ${\cal D}$, $M_V$ ou 
$M_m$. Il sera plus commode de consid\'erer toutes ces formes diff\'erentielles comme d\'efinies sur 
${\cal D}$ et invariantes sous l'action des groupes $\{ e \}$, $\Gamma_V$ ou $\Gamma_m$. \'Etant donn\'e 
un entier $m_0$, un \'el\'ement $\gamma \in \Gamma_{m_0}$ et une forme diff\'erentielle $\omega$ sur $M_m$ (avec 
$m\in {\Bbb N} \cup \{ \infty \}$, $m\geq m_0$), nous pourrons notamment parler de la forme diff\'erentielle 
$\gamma^* \omega$. 

\medskip

Concluons ce chapitre par l'\'etude promise des s\'eries de Poincar\'e.

Soient $Z_1 ,Z_2 \in {\cal D}$, $t \in {\Bbb R}$, $t>0$. On introduit~: 
\begin{eqnarray} 
\nu  (Z_1 ,Z_2 ,t) := | \{ \gamma \in \Gamma \; : \; d(Z_1 , \gamma Z_2 ) \leq t \} |,
\end{eqnarray}
et
\begin{eqnarray} \label{denombrement}
N (Z, t) := |\{ \gamma \in \Gamma_V \backslash \Gamma \; : \; d(\gamma Z , {\cal D}_V ) \leq t \} |.
\end{eqnarray}

\begin{lem} \label{comptage}
Il existe une constante $c_1 (Z) >0$ (qui d\'epend de $\Gamma$) telle que pour tout $t>0$ on ait~:
$$N (Z,t)  \leq c_1 (Z)  \int_0^{t+1} (1+t^{2pq}) e^{2(p+q-1) \sqrt{m} t} dt .$$
De plus on peut choisir $c_1 (Z)$ de mani\`ere \`a ce qu'elle soit born\'ee sur les compacts de ${\cal D}$.
\end{lem}
{\it D\'emonstration.}
Soit $\varepsilon$ un nombre r\'eel strictement compris entre $0$ et $1$ et suffisamment petit pour que 
$$B(Z ,\varepsilon ) \cap B(\gamma Z,\varepsilon ) \neq \emptyset \Rightarrow \gamma = e ,$$
o\`u $B(Z, \varepsilon)$ d\'esigne la boule de rayon $\varepsilon$ autour du point $Z$ et 
$e$ d\'esigne l'\'el\'ement neutre du groupe $\Gamma$. 
Dans la suite \'etant donn\'ee une sous-vari\'et\'e ${\cal V}$ de ${\cal D}$, nous noterons $B({\cal V}, r)$ l'ensemble des points de ${\cal D}$
\`a distance plus petite que $r$ de ${\cal V}$.
On a alors : 
\begin{eqnarray*}
N (Z,t)            & \leq & | \{ [\gamma] \in \Gamma_V \backslash \Gamma \; : \; \gamma (B(Z,\varepsilon )) \subset B({\cal D}_V ,t+\varepsilon ) \} |.
\end{eqnarray*}
Mais, d'apr\`es \ref{action}, si $\gamma \in \Gamma$ v\'erifie que $\gamma (B(Z, \varepsilon )) \subset B({\cal D}_V , \varepsilon )$
quitte \`a translater $\gamma$ par un \'el\'ement de $\Gamma_V$, on peut supposer que $\gamma (B(Z, \varepsilon ))
\subset B(S , \varepsilon )$, o\`u $S$ est un domaine fondamental mesurable pour l'action de $\Gamma_V$ sur 
${\cal D}_V$. On d\'eduit alors du Lemme \ref{croissance du volume}~:
\begin{eqnarray*}
N(Z,t)             & \leq & \frac{\mbox{vol}(B(S,t+\varepsilon ))}{\mbox{vol}(B(Z,\varepsilon ))} \\
                      & \leq & \frac{c}{\mbox{vol}(B(P,\varepsilon ))} \int_0^{t+\varepsilon} (1+t^{2pq}) e^{2(p+q-1) \sqrt{m}t} dt. \\
\end{eqnarray*}  
Ce qui ach\`eve la d\'emonstration du Lemme \ref{comptage}.

\bigskip

Remarquons que pour $r=p$, la d\'emonstration du Lemme \ref{comptage} permet d'estimer $\nu (Z_1, Z_2 ,t)$
uniform\'ement par rapport \`a $Z_2$. On obtient, en effet, que
pour tout $Z_2 \in {\cal D}$ et pour $t>0$, 
\begin{eqnarray} \label{nu}
\nu (Z_1 , Z_2 , t) \leq c_1 (Z_1 ) \int_0^{t+1} (1+t^{2pq} )e^{2(p+q-1) \sqrt{q} t} dt .
\end{eqnarray}
On en d\'eduit la proposition suivante.

\begin{prop} \label{serie de Poincare}
Soit $K$ un compact de ${\cal D}$. Alors il existe une constante $c_2 (K)$
(qui d\'epend de $\Gamma$) telle que pour tout point $Z_1 \in K$, tout point $Z_2 \in {\cal D}$
et $t \geq 0$, on ait~:
$$\sum_{\begin{array}{c}
\gamma \in \Gamma \\
d(Z_1 , \gamma Z_2) \leq t
\end{array}} e^{-(2(p+q-1+s) d(Z_1 ,Z_2 ) } \leq c_2 (K)\left( 1+ \frac{1}{s} \left( 1 + \frac{1}{s^{2pq}} \right) \right) , $$
pour tout $s>0$.
\end{prop}
{\it D\'emonstration.} D'apr\`es (\ref{nu}), il existe une constante $c_1 (K)$ telle que 
$$d\nu (Z_1 , Z_2 , t ) \leq c_1 (K) (1+(t+1)^{2pq} )e^{2(p+q-1) \sqrt{q} t} dt ,$$
pour tout $Z_1 \in K$, $Z_2 \in {\cal D}$ et $t>0$. On a donc~: 
\begin{eqnarray*}
\sum_{\begin{array}{c}
\gamma \in \Gamma \\
d(Z_1 , \gamma Z_2) \leq t
\end{array}} e^{-(2(p+q-1)\sqrt{q} + s) d(Z_1 ,Z_2 ) } & = & \int_0^t e^{-(2(p+q-1)\sqrt{q} +2s-2)t}  d \nu (Z_1 ,Z_2 ,t)  \\ 
                                                                                   & = & c_1 (K) \int_0^t (1+(t+1)^{2pq} ) e^{-st} dt .  
\end{eqnarray*}  
Et la Proposition \ref{serie de Poincare} d\'ecoule d'un calcul simple et d'approximations grossi\`eres.                                                      

\bigskip

De mani\`ere analogue on d\'emontre la proposition suivante.

\begin{prop} \label{serie de formes}
Soit $\phi$ une forme diff\'erentielle $\Gamma_V$-invariante de degr\'e $l$ sur ${\cal D}$.
Si $||\phi || \leq c \left( \frac{A}{B} \right)^{(p+q-1) \sqrt{m} + \varepsilon }$, $m= \min \{ r,q \}$
pour un r\'eel strictement positif $\varepsilon >0$, alors la s\'erie
$$\sum_{\Gamma_V \backslash \Gamma} \gamma^* \phi $$
converge uniform\'ement sur les compacts de ${\cal D}$.
\end{prop}
{\it D\'emonstration.} Commen\c{c}ons par remarquer que la norme $||\gamma^* \phi ||$ au point $Z$ est 
\'egale \`a la norme $|| \phi ||$ au point $\gamma Z$. D'apr\`es l'hypoth\`ese faite sur la norme de $\phi$,
\begin{eqnarray*}
\sum_{\Gamma_V \backslash \Gamma} ||\gamma^* \phi || & \leq & c \sum_{\Gamma_V \backslash \Gamma} \left( \frac{A}{B} (\gamma Z) \right)^{(p+q-1)\sqrt{m} +\varepsilon} \\
                                                                                             & \leq & \frac{c}{4^m} \sum_{\Gamma_V \backslash \Gamma} e^{-2((p+q-1)\sqrt{m} +\varepsilon )d(\gamma Z , {\cal D}_V )},
\end{eqnarray*}
d'apr\`es la Proposition \ref{distance2}. On conclut alors facilement comme pour la 
Proposition \ref{serie de Poincare}.

\medskip

\newpage

\thispagestyle{empty}

\newpage

\markboth{CHAPITRE 13. CONSTRUCTION DE LA FORME DUALE}{13.1. FORMES SINGULI\`ERES DE BOTT ET CHERN}

\chapter{Construction de la forme duale}

\section{Formes singuli\`eres de Bott et Chern}

Soit $X$ une vari\'et\'e complexe de dimension complexe $m$ et $E\rightarrow X$ un fibr\'e vectoriel holomorphe
de rang $q$. Supposons fix\'ee une section holomorphe $v : X \rightarrow E$ telle que 
\begin{enumerate}
\item l'ensemble $X_v = \{ x\in X \; : \; v(x) =0 \}$ des z\'eros de $v$ soit non singulier de codimension $q$, et
\item la section nulle de $X$ dans $E$ lui soit transverse.
\end{enumerate}
Soit $C_q (E)$ la forme de Chern maximale associ\'ee \`a une structure hermitienne fix\'ee sur $E$. Dans \cite{BottChern},
Bott et Chern montrent 

\begin{prop}  \label{BC}
 \index{Proposition de Bott et Chern}
Il existe une forme diff\'erentielle $\tau$ de type $(q-1,q-1)$ sur $X-X_v$ et \`a valeurs r\'eelles telle que 
$$\overline{\partial} \partial \tau = C_q (E).$$
\end{prop}

La construction de $\tau$ que l'on va rappeler ci-dessous est explicite, remarquons que $\tau$ doit n\'ecessairement 
avoir des singularit\'es le long de $X_v$. Il est  bien connu \cite{KobayashiNomizu} que lorsque $X$ est compacte $C_q (E)$ repr\'esente la 
classe de cohomologie duale \`a la classe d'homologie $[X_v ]$. Nous verrons que dans le cas non compact (celui qui
nous int\'eresse dans la suite) cette dualit\'e subsiste en un sens plus faible.

\medskip

Rappelons le formalisme de Bott et Chern. En un point $x \in X$, soient $E_x$ et $T_x^*$ respectivement 
la fibre de $E$ au-dessus de $x$ et la fibre de l'espace cotangent holomorphe $T^*$ au-dessus de $x$. On a alors
les fibr\'es $\Lambda_{rs}^{pq} \rightarrow X$, o\`u 
$$\Lambda_{rs}^{pq} = \Lambda^{pq} (E_x \oplus \overline{E}_x ) \otimes \Lambda^{rs} (T_x^* \oplus \overline{T}_x^* ).$$
Notons $\tilde{\Lambda}_{rs}^{pq}$ l'espace des sections $C^{\infty}$ de $\Lambda_{rs}^{pq}$. Une connexion affine 
sur $E$ est d\'efinies par la donn\'ee de deux op\'erateurs 
$$\begin{array}{ll}
\partial : \tilde{\Lambda}_{00}^{10} \rightarrow \tilde{\Lambda}_{10}^{10} , & \overline{\partial} : \tilde{\Lambda}_{00}^{10} \rightarrow \tilde{\Lambda}^{10}_{01} \\
d= \partial + \overline{\partial} &
\end{array}$$
v\'erifiant
$$d(a+b) = da +db , \; a,b \in \tilde{\Lambda}_{00}^{10} $$
$$d(fa) = df \otimes a + f da , \; f \in \tilde{\Lambda}_{00}^{00} .$$
Une connexion affine $d$ est dite de type $(1,0)$ si $\overline{\partial} (a) = 0$ pour toute section holomorphe $a$.
Toute connexion affine d\'efinit bien \'evidemment de mani\`ere unique des op\'erateurs
$$ \begin{array}{ll}
\partial : \tilde{\Lambda}_{rs}^{pq} \rightarrow \tilde{\Lambda}_{(r+1)s}^{pq} , & \overline{\partial} : \tilde{\Lambda}_{rs}^{pq} \rightarrow \tilde{\Lambda}^{pq}_{r(s+1)} \\
d= \partial + \overline{\partial} &
\end{array}$$
v\'erifiant les propr\'et\'es usuelles.

Soit $e=(e_1 , \ldots , e_q )$ un champs de base local. Alors,
\begin{eqnarray} \label{forme de courbure}
de = e \omega = (e_1 , \ldots , e_q ) (\omega_{ij} ) \\
d^2 e = e \Omega = (e_1 , \ldots , e_q ) (\Omega_{ij} ) , 
\end{eqnarray}
o\`u $\omega$, $\Omega$ sont respectivement les matrices de connexion et de courbure.

Il d\'ecoule de (\ref{forme de courbure}) que
\begin{eqnarray} \label{omega}
\Omega = d\omega + \omega \wedge \omega.
\end{eqnarray}
Remarquons que les \'el\'ements de chaque fibre sont vus comme des vecteurs colonnes. La $q$-i\`eme forme de
Chern de la connexion affine est la forme diff\'erentielle $\left( \frac{\sqrt{-1}}{2\pi} \right)^q \mbox{det } \Omega$.

Fixons maintenant une structure hermitienne sur $E$, {\it i.e.}  un produit scalaire hermitien $\langle.,.\rangle_x$ sur chaque 
fibre $E_x$ ($x\in X$) qui soit $C^{\infty}$ en $x$. \'Etant donn\'e un champs de base local $e=(e_1 , \ldots , e_q )$,
notons 
\begin{eqnarray}  \label{H}
H=(h_{ij} ), \; h_{ij} = \langle e_i , e_j \rangle.
\end{eqnarray}

La structure hermitienne d\'efinit une unique connexion de type $(1,0)$ pour laquelle on a 
$$\omega =H^{-1} \partial H \; \mbox{ et } \; \Omega = \overline{\partial} (H^{-1} \partial H) .$$
Afin de construire $\tau$ on se fixe une section holomorphe $v$ du fibr\'e $E$. Soient alors 
\begin{eqnarray*}
\alpha & = & dv d\overline{v} = \overline{\alpha} \in \tilde{\Lambda}_{11}^{11} \\
K & = & -e \Omega H^{-1} {}^t \overline{e} = -(e_1 \ldots e_q ) \Omega H^{-1} \left( 
\begin{array}{c}
\overline{e}_1 \\
\vdots \\
\overline{e}_q 
\end{array} \right) = \overline{K} \in \tilde{\Lambda}_{11}^{11} \\
y_k & = & - \overline{y}_k = v \overline{v} \alpha^{k-1} K^{q-k} \in \tilde{\Lambda}_{q-1, q-1}^{qq}, \; 1\leq k \leq q ,\\
s_k & = & (\overline{v} dv) \alpha^{k-1} K^{q-k} \in \tilde{\Lambda}^{qq}_{q,q-1} , \; 1\leq k \leq q ,\\
w_k & = & \alpha^k K^{q-k} \in \tilde{\Lambda}_{qq}^{qq} , \; 0\leq k \leq q,
\end{eqnarray*}
et
\begin{eqnarray*}
s_k ' & = & |v|^{-2k} s_k , \; y_k ' = |v|^{-2k} y_k , \; w_k ' = |v|^{-2k} w_k \\
s_k '' & = & s_k ' +ky_k ' \partial \log |v|^2 , \; 1\leq k \leq q, \\
u_k & = & w_k ' +ks_k ' \overline{\partial} \log |v|^2 , \; 0 \leq k \leq q ,\\
\chi & = & (\det H)^{-1/2} e_1 \ldots e_q , 
\end{eqnarray*}
o\`u $|v|$ est la norme de $v$.

\begin{lem} \label{formules recursives}
On a les formules de r\'ecurrence~:
\begin{eqnarray} \label{y}
\partial y_k ' = -s_k '' - \frac{k-1}{q-k+1} s_{k-1} '' , 
\end{eqnarray}
\begin{eqnarray} \label{s}
\overline{\partial} s_k ' = u_k + \frac{k}{q-k+1} u_{k-1} .
\end{eqnarray}
\end{lem}
{\it D\'emonstration.} D\'ecomposons tout d'abord $dv$ et $d^2 v$ dans la direction de $v$ et de son suppl\'ementaire
ortogonal~:
\begin{eqnarray*}
dv = \theta v + \beta , \; \theta \in \tilde{\Lambda}_{10}^{00} , \; \beta \in \tilde{\Lambda}_{10}^{10} , \\
d^2 v =  |v|^2 ( \phi v + \gamma ), \; \phi \in  \tilde{\Lambda}_{11}^{00} , \; \gamma \in \tilde{\Lambda}_{11}^{10} , \\ 
\theta = \partial \log |v|^2 , \; \phi + \overline{\phi} =0 .
\end{eqnarray*}
Par r\'ecurrence, on obtient
\begin{eqnarray*}
\alpha^{k-1}&  = & \alpha_1^{k-1} + (k-1) \alpha_1^{k-2} (v\theta \overline{\beta} + \overline{v} \overline{\theta} \beta ) \\
& & + (k-1)^2 \alpha_1^{k-2} v\overline{v} \theta \overline{\theta} , \\
K^{q-k} & = & K_1^{q-k} + (q-k) (v \overline{\gamma} + \overline{v} \gamma ) K_1^{q-k-1}  \\
& & + (q-k) v\overline{v}  K_1^{q-k-2} \left\{ -K_1 \phi + (q-k-1) \gamma \overline{\gamma}  \right\} , \\
\alpha_1 & = & \beta \overline{\beta} .
\end{eqnarray*}
D'o\`u l'on d\'eduit que 
\begin{eqnarray*}
s_k & = & -k \theta y_k + (q-k) v\overline{v} \beta \overline{\gamma} \alpha_1^{k-1} K_1^{q-k-1} , \\
\frac{1}{q-k} w_k & = & v\overline{v} \alpha_1^k K_1^{q-k-2} \left\{ -K_1 \phi + (q-k-1) \gamma \overline{\gamma} \right\} \\
& & +kv\overline{v} (-\theta \overline{\beta} \gamma + \overline{\theta} \beta \overline{\gamma} ) \alpha_1^{k-1} K_1^{q-k-1} \\
& & + \frac{k^2}{q-k} v\overline{v} \theta \overline{\theta} \alpha_1^{k-1} K_1^{q-k} .
\end{eqnarray*}
Puisque $dK=0$ et $d^2 v \in \tilde{\Lambda}^{10}_{11}$,
$$\partial y_k = -s_k - \frac{k-1}{q-k+1} |v|^2 s_{k-1} - \frac{(k-1)^2}{q-k+1} |v|^2 y_{k-1} \partial \log |v|^2 .$$
D'o\`u il d\'ecoule l'\'egalit\'e (\ref{y}).
De m\^eme, on a 
$$\frac{1}{|v|^2} (\overline{\partial} s_k -w_k) = \frac{k}{q-k+1} w_{k-1} + \frac{k(k-1)}{q-k+1} s_{k-1} \overline{\partial} \log |v|^2 .$$
Ce qui implique (\ref{s}) et conclut la d\'emonstration du Lemme \ref{formules recursives}.

\medskip

Remarquons qu'il d\'ecoule de la d\'emonstration du Lemme \ref{formules recursives} que
\begin{eqnarray} \label{zero en q}
s_q '' =0 \; \mbox{ et } \; u_q =0.
\end{eqnarray}

\medskip

On peut maintenant d\'efinir la forme singuli\`ere $\tau$ par~:
\begin{eqnarray} \label{to}
c \tau \chi \overline{\chi} & = &  -q \sum_{k=1}^{q-1} (-1)^k \left( 
\begin{array}{c}
q-1 \\
k-1
\end{array} \right) \left( \frac{1}{k} + \ldots + \frac{1}{q-1} \right) y_k ' \\
& & -q \sum_{k=1}^{q-1} (-1)^k \left(
\begin{array}{c}
q-1 \\
k-1
\end{array} \right) y_k ' \log |v|^2 , 
\end{eqnarray}
o\`u $c= (-1)^{q^2 /2 -1} q! (2\pi )^q$. On introduit \'egalement les formes singuli\`eres $\psi_k$, $1\leq k \leq q$, et $\psi$
d\'efinies par~:
\begin{eqnarray} \label{psik}
c \psi_k \chi \overline{\chi} = s_k ' , \; 1\leq k \leq q
\end{eqnarray}
et 
\begin{eqnarray} \label{psi1}
\psi = \sum_{k=1}^q (-1)^k \left( 
\begin{array}{c}
q \\
k 
\end{array} \right) \psi_k .
\end{eqnarray}

\begin{prop} \label{bott et chern}
Sur $X-X_v$, les formes $\tau$ et $\psi$ sont lisses (non singuli\`eres) et v\'erifient
$\partial \tau = \psi$ et $\overline{\partial} \psi = C_q (E)$.
\end{prop}
{\it D\'emonstration.} On sait que 
$$\partial (\chi \overline{\chi} )  = \overline{\partial} (\chi \overline{\chi}) =0 .$$
L'identit\'e $\partial \tau = \psi$ sur $X-X_v$ d\'ecoule donc de (\ref{y}), (\ref{zero en q}) et (\ref{psi1}).

L'\'equation (\ref{s}) implique alors
$$\overline{\partial} \sum_{k=1}^q (-1)^k \left(
\begin{array}{c}
q \\
k 
\end{array} \right) s_k ' = -K^q .$$
Puisque $C_q (E) = \left( \frac{i}{2\pi} \right)^q \det (\Omega_{ij} )$, 
$$K^q = (-1)^{q^2/2} q! (2\pi )^q C_q (E) \chi \overline{\chi}.$$
Et l'identit\'e $\overline{\partial} \psi = C_q (E)$ sur $X-X_v$ s'en d\'eduit.

\medskip

La signification g\'eom\'etrique de la forme diff\'erentielle $\psi$ est contenue dans la proposition 
suivante.

\begin{prop} \label{psi duale}
Soient $T_{\varepsilon} \subset X$ un voisinage tubulaire de rayon $\varepsilon$ autour de $X_v$
et $\eta$ une $(m-q , m-q )$-forme sur $X$. Alors,
$$\lim_{\varepsilon \rightarrow 0} \int_{\partial T_{\varepsilon} } \psi \wedge \eta = - \int_{X_v} \eta .$$
\end{prop}
{\it D\'emonstration.} Ce r\'esultat est local, on peut donc supposer que $X= {\Bbb C}^m$, $E={\Bbb C}^m \times {\Bbb C}^q \rightarrow X$
est le fibr\'e produit et que $v(Z_1 , \ldots , Z_m ) = (Z_1 , \ldots , Z_q )$. Supposons maintenant que dans ces 
coordonn\'ees la structure hermitienne soit donn\'ee par la matrice $H= (h_{ij} )$. Choisissons un champ local de bases orthonorm\'ees
$(e_1 , \ldots , e_q )$. D'apr\`es la construction ci-dessus, la forme diff\'erentielle $\psi$ a un p\^ole d'ordre $2q-1$
le long de $X_v$ qui appara\^{\i}t dans le terme $\psi_q$. Exprimons la section $v$ dans notre champ de base locale~:
$$v = \sum_{\lambda =1}^q f_{\lambda} e_{\lambda} $$
$$dv = \sum_{\lambda =1}^q (df_{\lambda}) e_{\lambda} + \sum_{\lambda=1}^q \sum_{j=1}^q (\omega_{\lambda j} f_j )e_{\lambda}$$
et il d\'ecoule de (\ref{psik}) que la forme $\psi_q$ est \'egale \`a 
$$(-1)^{-q^2 /2 +1+q } \frac{(q-1)!}{(2\pi)^q (f)^{2q}} (df_1 \wedge \ldots \wedge df_q ) \left( \sum_{\lambda} 
(-1)^{\lambda -1} \overline{f}_{\lambda} d\overline{f}_1 \wedge \ldots \wedge \hat{d \overline{f}_{\lambda} } \wedge \ldots \wedge d\overline{f}_q \right)$$
plus des termes d'ordres inf\'erieurs. Il s'ensuit que 
$$\lim_{\varepsilon \rightarrow 0} \psi_q \wedge \eta = (-1)^{1+q} \int_{X_v} \eta $$
et donc que 
$$\lim_{\varepsilon \rightarrow 0} \int_{\partial T_{\varepsilon}} \psi \wedge \eta = (-1)^q \lim_{\varepsilon \rightarrow 0} 
\int_{\partial T_{\varepsilon}} \psi_q \wedge \eta = -\int_{X_v} \eta .$$

\markboth{CHAPITRE 13. CONSTRUCTION DE LA FORME DUALE}{13.2. CONSTRUCTION DE LA FORME DUALE}

\section{Construction de la forme duale} 

On revient aux notations du chapitre pr\'ec\'edent, soit $e_1 , \ldots , e_n$ la base standard de ${\Bbb C}^n$, $n=p+q$.
\'Etant donn\'e un entier $1\leq r <p$, soit $V$ un sous-espace positif de dimension $r$ pour la forme hermitienne $Q$
et ${\cal D}_V$ et $G_V$ comme au chapitre pr\'ec\'edent. Soit encore $\Gamma$ un sous-groupe discret et sans torsion
de $G$ et $\Gamma_V = \Gamma \cap G_V$. On suppose de plus~:
\begin{enumerate}\item $\Gamma \backslash {\cal D}$ compacte, et
\item $\Gamma_V \backslash {\cal D}_V$ compacte.
\end{enumerate}

\medskip

Comme au chapitre pr\'ec\'edent, il sera plus commode dans la suite de supposer que le sous-espace $V$ est le
sous-espace engendr\'e par $e_{p-r+1} , \ldots , e_p$. 
Le groupe $G_V$ agit sur ${\cal D} \times M_{qr} ({\Bbb C})$ par~:
$$g(z,M) =(gz,  {}^t j(g,z)^{-1} M {}^t u ) .$$
En quotientant par le groupe $\Gamma_V$, on obtient un fibr\'e vectoriel $E_V = {\cal D} \times_{\Gamma_V} M_{qr} ({\Bbb C})$
au-dessus de $\Gamma_V \backslash {\cal D}$. Ce fibr\'e est naturellement muni de deux m\'etriques hermitiennes, \`a 
savoir~:
$$\langle M_1 , M_2 \rangle_Z = {\rm tr} ({}^t \overline{M}_1 h_Z M_2 ) $$
et 
$$\langle M_1 , M_2 \rangle^{\sim}_Z = {\rm tr} ({}^t \overline{M}_1  \tilde{h}_Z M_2 ),$$
o\`u $h_Z$ et $\tilde{h}_Z$ sont les matrices du Lemme 12.2.3. D'apr\`es ce dernier, la m\'etrique $\langle .,. \rangle_Z^{\sim}$ est 
$G$-invariante
et la m\'etrique $\langle .,.\rangle_Z$ est $G_V$-invariante. Dans cette section nous ne nous servirons que de la  
m\'etrique $\langle .,.\rangle^{\sim}$. Enfin, le fibr\'e $E_V$ nous arrive \'equip\'e d'une section holomorphe
canonique $v: \Gamma_V \backslash {\cal D} \rightarrow E_V$
$$v(Z) = (Z, {}^t Z_2 ).$$
Il est clair que $\Gamma_V \backslash {\cal D}_V = \{ Z \in \Gamma \backslash {\cal D} \; : \; v(Z)=0 \}$. On peut
donc appliquer les r\'esultats de la section pr\'ec\'edente. 

\medskip

On verra $M_{qr} ({\Bbb C})$ comme $M_r (M_q ({\Bbb C}))$. \'Etant donn\'e une matrice carr\'ee $A \in M_q ({\Bbb C})$, nous notons 
$A^{[r]}$ la matrice carr\'ee diagonale par bloc constitu\'ee de $r$ fois le bloc $A$. 
Alors la m\'etrique hermitienne $\langle .,.\rangle^{\sim}_Z$, de la fibre au-dessus de $Z$, est donn\'ee par la matrice hermitienne d\'efinie positive
$H =\tilde{h}_Z^{[r]}$ par rapport au rep\`ere mobile standard $e=(e_1 , \ldots , e_r )$, o\`u chaque $e_i$ est un 
vecteur $(e_{1i} , e_{2i} , \ldots , e_{qi} )$. Un calcul facile montre~:
$$d\tilde{h}_Z = \tilde{h}_Z \omega_1 + {}^t \overline{\omega}_1  \tilde{h}_Z , \; \omega_1 = d{}^tZ \overline{Z} (I_q - {}^t Z \overline{Z} )^{-1} .$$
Soit $\omega = \omega_1^{[r]}$. On choisit maintenant une connexion m\'etrique telle que $de=e\omega$. 
La forme de courbure est alors donn\'ee par 
$$\Omega = \Omega_1^{[r]} , \; \Omega_1 = -d{}^t Z (I_p - \overline{Z} {}^t Z )^{-1} d\overline{Z} (I_q - {}^t Z \overline{Z} )^{-1} .$$

\begin{lem} \label{action sur connexion}
Si $g$ est un \'el\'ement de $G$, on a~:
$$g^* \omega_1 = {}^t j(g,Z)^{-1}  \omega_1 {}^t j(g,Z) + {}^t j(g,Z)^{-1}  d {}^t j(g,Z) $$
et 
$$g^* \Omega_1 = {}^t j(g,Z)^{-1}  \Omega_1 {}^t j(g,Z) .$$
\end{lem}
{\it D\'emonstration.} Rappelons que 
$$I_q - {}^t Z \overline{Z} = {}^t j(g,Z) (I_q - {}^t (gZ) \overline{gZ} ) \overline{j(g,Z)} .$$
En diff\'erentiant cette \'egalit\'e, on obtient~:
$$-d{}^t Z \overline{Z} = d{}^t j(g,Z) (I_q -{}^t (gZ) \overline{gZ} ) \overline{j(g,Z)}  - j(g,Z)^t d(gZ)^t \overline{gZ} \overline{j(g,Z)},$$
puis
$$d{}^t (gZ) \overline{gZ} = {}^t j(g,Z)^{-1} d{}^t Z \overline{Z} \overline{j(g,Z)} ^{-1} + {}^t j(g,Z)^{-1} d{}^t j(g,Z) (I_q - {}^t (gZ) \overline{gZ} ),$$ 
et donc
$$d{}^t (gZ) \overline{gZ} (I_q - {}^t (gZ) \overline{gZ} )^{-1} = {}^t j(g,Z)^{-1} d{}^t Z \overline{Z} ( I_q - {}^t Z \overline{Z})^{-1}
{}^t j(g,Z) + {}^t j(g,Z)^{-1} d{}^t j(g,Z) .$$
Ce qui d\'emontre la premi\`ere identit\'e annonc\'ee par le Lemme \ref{action sur connexion}.
La deuxi\`eme identit\'e s'obtient alors en diff\'erenciant la premi\`ere et en utilisant que 
$\Omega_1 = d\omega_1 + \omega_1 \wedge \omega_1$.

\medskip

Dans la suite, pour simplifier les notations, notons
$$E = \left(
\begin{array}{ccc}
e_{11} & \ldots & e_{q1} \\
\vdots &           & \vdots \\
e_{1r} & \ldots & e_{qr} 
\end{array} \right).$$
Dans le rep\`ere mobile standard, on a~:
\begin{eqnarray*}
v & = & {\rm tr} \left( E {}^t Z_2 \right),\\
dv & = & {\rm tr} \left( E d {}^t Z \overline{Z} (I_q - {}^t Z \overline{Z} )^{-1} {}^t Z_2 + E d{}^t Z_2 \right) , \\
K & = & {\rm tr} \left( E d{}^t Z (I_p -\overline{Z} {}^t Z )^{-1} d\overline{Z}  {}^t \overline{E} \right) , \\
|v|^2 & = & {\rm tr} \left( Z_2 (I_q - {}^t Z \overline{Z} )^{-1} {}^t \overline{Z}_2  \right) .
\end{eqnarray*}
Les formes $\tau$, $\psi_k$ et $\psi$ de la section pr\'ec\'edente sont bien d\'efinies par les formules 
(\ref{to}), (\ref{psik}) et (\ref{psi}). Pour se rappeler dans la suite que ces formes sont obtenues en consid\'erant le 
produit hermitien $\langle .,.\rangle^{\sim}$ dans la fibre, nous pr\'ef\`erons les noter $\tilde{\tau}$, $\tilde{\psi}_k$ et $\tilde{\psi}$.

\begin{prop} \label{action sur psi}
Les formes $\tilde{\psi}_k$ et $\tilde{\psi}$ sont $G_V$-invariantes.
\end{prop}
{\it D\'emonstration.} Notons $v(E,Z)$, $dv(E,Z)$, $K(E,Z)$ et $\chi (E,Z)$ les fonctions $v$, $dv$, $K$ et $\chi$
dont on pointe la d\'ependance en les variables $E$ et $Z$. D'apr\`es (\ref{action}) et le Lemme 12.2.3, il est clair que, si $g \in G_V$,
$$v(E,gZ ) = v({}^t u E {}^t j(g,Z)^{-1} ,Z)$$
et
$$\chi (E,gZ) = | \det j(g,Z) |^{-r} \chi (E,Z) .$$
Puis il d\'ecoule du Lemme \ref{action sur connexion} que
$$dv (E,gZ)= dv ({}^t u E {}^t j(g,Z)^{-1}  ,Z)$$
et
$$K(E, gZ) = K({}^t u E {}^t j(g,Z)^{-1}  ,Z).$$
Puisque $|v|$ est $G_V$-invariant, d'apr\`es la formule pour $s_k '$ donn\'ee dans la section pr\'ec\'edente et les
identit\'es ci-dessus, on obtient~:
\begin{eqnarray*}
s_k ' (E ,gZ) & = & s_k '({}^t u E {}^t j(g,Z)^{-1}  , Z) \\
                    & = & |\det u|^{2q} |\det j(g,Z) |^{-2r} s_k '(E,Z) \\
                    & = & |\det j(g,Z) |^{-2r} s_k '(E,Z);
\end{eqnarray*}
et donc
$$g^* \tilde{\psi}_k = \tilde{\psi}_k$$
pour tout $g\in G_V$. Puisque par d\'efinition, $\tilde{\psi}$ est une combinaison lin\'eaire de $\tilde{\psi}_k$, la Proposition 
\ref{action sur psi} est d\'emontr\'ee.

\bigskip

Rappelons maintenant que, d'apr\`es le Lemme \ref{outils}, 
$$\frac{B}{A} = \det \{ I_r + Z_2 (I_q -{}^tZ \overline{Z} )^{-1} {}^t Z_2 \}  .$$
La matrice $Z_2 (I_q - {}^t Z \overline{Z} )^{-1} {}^t Z_2$ est positive (au sens large),
notons $\lambda_1 , \ldots , \lambda_r$ ses valeurs propres (r\'eelles et positives). Remarquons 
que $|v|^2 = {\rm tr} (Z_2 (I_q - {}^t Z \overline{Z} )^{-1} {}^t Z_2) = \lambda_1 + \ldots + \lambda_r $. On en d\'eduit que 
$$\frac{A}{B} = \prod_{i=1}^r (1+\lambda_i )^{-1} \geq \prod_{i=1}^r (1-\lambda_i ) \geq 1 -\sum_{i=1}^r \lambda_i  = 1-|v|^2 ,$$
autrement dit
\begin{eqnarray} \label{module de v}
\frac{C}{B} \leq |v|^2 ,  \; \; C=B-A .
\end{eqnarray}
De plus, la fonction $(C/B)/|v|^2$ tend vers $1$ lorsque $v$ tend vers $0$. 

D'apr\`es la Proposition \ref{action sur psi}, pour calculer les normes de $\tilde{\psi}_k$ et $\tilde{\psi}$ on peut 
supposer que $Z_1 =0$. Il d\'ecoule alors du Lemme \ref{comparaison des normes} et de (\ref{module de v}) que 
$$||\tilde{\psi}_k || \prec \left( \frac{B}{C} \right)^{k-1/2} \left( \frac{B}{A} \right)^{r+rq+k-1}$$
ce qui implique que 
\begin{eqnarray} \label{norme de psi}
|| \tilde{\psi} || \prec \left( \frac{B}{C} \right)^{rq -1/2} \left( \frac{B}{A} \right)^{r+2rq -1},
\end{eqnarray}
o\`u dans les deux derni\`eres expressions, le signe $\prec$ signifie que l'on a une in\'egalit\'e $\leq$ 
\`a une constante positive pr\`es. 

\medskip

\'Etant donn\'e un nombre complexe $s$, soit $h_s (t)$ la fonction d\'efinie par
$$h_s (t) = -\int_t^{\infty} x^{-s} (x-r)^{qr -1} dx \; \; (\mbox{Re} (s) > qr ) .$$
On v\'erifie facilement que la fonction $h_s (t)$ v\'erifie~:
\begin{enumerate}
\item $h_s '(t) = t^{-s} (t-r)^{qr -1}$, et
\item $h_s (r) = -r^{qr -s} \frac{\Gamma (s-qr) \Gamma (qr)}{\Gamma (s)}$.
\end{enumerate}
On peut maintenant d\'efinir la forme diff\'erentielle $\omega_s$ par
\begin{eqnarray} \label{omegas}
\omega_s & = & \frac{-1}{h_s (r)} d\left( h_s (r+ |v|^2 ) \tilde{\psi} \right) .
\end{eqnarray}
D'apr\`es l'expression de $h_s '$, 
$$\omega_s = \frac{-1}{h_s (r)} \left\{ 2(r+|v|^2 )^{-s} |v|^{2rq-1} d|v| \wedge \tilde{\psi} + h_s (r+|v|^2 ) d\tilde{\psi} \right\} .$$
Il est clair que les formes $|v|^{2rq-1} \tilde{\psi}_k \wedge d|v|$ ($k\leq rq-1$) sont lisses (non singuli\`eres) et d'apr\`es (\ref{zero en q}),
la forme $|v|^{2rq -1} d|v| \wedge \tilde{\psi}_{rq}$ est, elle aussi, lisse. Par construction $d\tilde{\psi}$ est lisse aussi. On en 
d\'eduit donc que la forme $\omega_s$ est bien lisse (non singuli\`ere).  De plus, il d\'ecoule du Lemme \ref{outils} et 
de (\ref{norme de psi}) que 
\begin{eqnarray} \label{norme de omegas}
|| \omega_s || \prec \left( \frac{B}{A} \right)^{r+q+2rq +1 - r^{-1} {\rm Re}(s)} .
\end{eqnarray}

\medskip

D'apr\`es la Proposition \ref{action sur psi}, la forme $\omega_s$ est $G_V$-invariante. On peut montrer que 
la forme $\omega_s$, pour Re$(s)>>0$, peut se voir comme la forme duale \`a $\Gamma_V \backslash {\cal D}_V$ dans 
$\Gamma_V \backslash {\cal D}$. Nous allons pr\'eciser un peu ce r\'esultat.

\medskip

Soit d'abord $\mu$ une $k$-forme {\it harmonique} \index{harmonique, forme} sur $C_V := \Gamma_V \backslash {\cal D}_V$. On note $*_0$ 
l'op\'erateur $*$ de Hodge de la vari\'et\'e $C_V$ relativement \`a sa m\'etrique riemannienne. La forme $\pi^* (*_0 \mu )$
est alors une $(2(p-r)q - k)$-forme ferm\'ee sur ${\cal D}$. On va en fait construire la forme duale (dans un sens $L^2$)
\`a la forme $\pi^* (*_0 \mu )$.

Afin de d\'ecrire ce que cela signifie, soit $\phi$ une forme lisse quelconque sur $\Gamma_V \backslash {\cal D}$ de 
degr\'e $k$. On suppose que 
\begin{eqnarray} \label{condition sur phi}
||\phi || \prec \left( \frac{A}{B} \right)^N
\end{eqnarray}
pour un certain entier $N$. Remarquons que la condition (\ref{condition sur phi}) est v\'erifi\'ee par tout forme born\'ee pour
$N=0$. D'apr\`es (\ref{norme de omegas}) et le Lemme \ref{formules d'integration}, l'int\'egrale 
$\int_{\Gamma_V \backslash {\cal D}} (\omega_s \wedge \pi^* (*_0 \mu )) \wedge \phi$ est absolument convergente pour 
Re$(s) >>0$, la constante ne d\'ependant que de $N$. 

\begin{thm} \label{dualite}
Soit $\phi$ une forme ferm\'ee, lisse, de degr\'e $k$ sur $\Gamma_V \backslash {\cal D}$ et v\'erifiant la 
condition (\ref{condition sur phi}). Alors,
$$\int_{\Gamma_V \backslash {\cal D}} (\omega_s \wedge \pi^* (*_0 \mu )) \wedge \phi = 
\int_{C_V} (*_0 \mu) \wedge \phi \; \; (\mbox{Re}(s)>>0).$$
\end{thm}
{\it D\'emonstration.} La vari\'et\'e $C_V$ est une vari\'et\'e riemannienne compl\`ete. Soit $x_0$ un 
point fixe de $C_V$. Pout $t>0$, soit $B_t$ la boule ferm\'ee de centre $x_0$ et de rayon $t$. Soit $\partial B_t$ son 
bord et vol$(\partial B_t )$ son volume par rapport \`a la m\'etrique induite. Il est clair que 
$$\int_0^{\infty}  \mbox{vol} (\partial B_t ) dt = \mbox{vol} (C_V ) < \infty .$$
D'o\`u il d\'ecoule que 
\begin{eqnarray} \label{vol du bord}
\lim_{t \rightarrow \infty} \mbox{vol} (\partial B_t ) =0 .
\end{eqnarray}
\'Etant donn\'e deux r\'eels strictement positifs $t$, $\varepsilon>0$ et un r\'eel $l>>0$, soit
$N(t,\varepsilon , l)$ le sous-ensemble de $\Gamma_V \backslash {\cal D}$ constitu\'e des projet\'es des points 
$Z = \left(
\begin{array}{c}
Z_1 \\
Z_2
\end{array} \right)$ tels que 
\begin{enumerate}
\item $\left( \begin{array}{c}
Z_1 \\
0
\end{array} \right) \in B_t$,
\item $|v(Z)|\geq \varepsilon $,
\item la distance, $d(Z,{\cal D}_V )$ de $Z$ \`a ${\cal D}_V$ est inf\'erieure ou \'egale \`a $l$.
\end{enumerate}

La condition (\ref{condition sur phi}) assure la convergence absolue de l'int\'egrale et donc que 
\begin{eqnarray} \label{double limite}
\int_{\Gamma_V \backslash {\cal D}} (\omega_s \wedge \pi^* (*_0 \mu )) \wedge \phi = 
\lim_{t\rightarrow \infty} \lim_{
\begin{array}{c}
\varepsilon \rightarrow 0 \\ 
l \rightarrow \infty
\end{array} } \int_{N(t,\varepsilon ,l)} (\omega_s \wedge
\pi^* (*_0 \mu )) \wedge \phi .
\end{eqnarray}
Il est clair que 
$$\partial N(t,\varepsilon , l) = {\cal F}_{\varepsilon} \cup {\cal F}_l \cup {\cal F}_{\partial} ,$$
o\`u 
\begin{eqnarray*}
{\cal F}_{\varepsilon} & = & \{ Z\in N(t,\varepsilon , l) \; : \; |v(Z)| =\varepsilon \} ,\\
{\cal F}_l & = & \{ Z \in N(t,\varepsilon ,l) \; : \; d(Z,{\cal D}_V ) =l \} , \\
{\cal F}_{\partial} & = & \{ Z \in N(t, \varepsilon , l) \; : \; \left(
\begin{array}{c}
Z_1 \\
0
\end{array} \right) \in \partial B_t \} .
\end{eqnarray*}
Remarquons maintenant que 
$$(\omega_s \wedge \pi^* (*_0 \mu ) ) \wedge \phi = d \left( \frac{-h_s (r+|v|^2 )}{h_s (r)} (\tilde{\psi} \wedge \pi^* (*_0 \mu ))\wedge \phi \right) .$$
Le Th\'eor\`eme de Stokes implique donc~:
\begin{eqnarray*}
\int_{N(t,\varepsilon ,l)} (\omega_s \wedge \pi^* (*_0 \mu )) \wedge \phi & = & -\int_{{\cal F}_l} \frac{h_s (r+|v|^2 )}{h_s (r)} 
\tilde{\psi} \wedge \pi^* (*_0 \mu ) \wedge \phi \\
& & - \int_{{\cal F}_{\varepsilon}}  \frac{h_s (r+|v|^2 )}{h_s (r)} 
\tilde{\psi} \wedge \pi^* (*_0 \mu ) \wedge \phi \\
& & -\int_{{\cal F}_{\partial}} \frac{h_s (r+|v|^2 )}{h_s (r)} \tilde{\psi} \wedge \pi^* (*_0 \mu ) \wedge \phi .
\end{eqnarray*}
Notons $I_l$, $I_{\varepsilon}$ et $I_{\partial}$ les trois int\'egrales du membre de droite de l'\'egalit\'e ci-dessus. 
D'apr\`es (\ref{norme de psi}), on peut supposer que 
$$||h_s (r+|v|^2 ) \tilde{\psi} \wedge \pi^* (*_0 \mu ) \wedge \phi || \prec \left( \frac{B}{C} \right)^{rq-1/2} \left( \frac{A}{B} \right)^a$$
avec $a>>0$. Puis, d'apr\`es le Lemme 12.3.2, il existe une constante $b>0$ telle que vol$({\cal F}_l ) \prec e^{bl}$vol$(B_t)$.
Enfin, il d\'ecoule de la Proposition 12.2.5 que $\frac{A}{B} \prec e^{-2d(Z,{\cal D}_V )}$. On obtient donc que
$$I_l \prec \mbox{vol} (B_t ) e^{(b-2a)l} .$$
Puisque $a>>0$, il en d\'ecoule que 
\begin{eqnarray} \label{il}
\lim_{l\rightarrow \infty} I_l =0 .
\end{eqnarray}
La d\'emonstration de la Proposition 13.1.4 implique 
\begin{eqnarray} \label{iepsilon}
\lim_{\varepsilon \rightarrow 0} I_{\varepsilon} = - \int_{B_t } (*_0 \mu) \wedge \phi .
\end{eqnarray}
Afin d'estimer $|I_{\partial} |$, nous commen\c{}ons par int\'egrer le long de la fibre. Soit $\eta$ la forme volume de $\partial B_t$.
D'apr\`es (12.3.8), $||\eta \wedge dv_F || \succ \left( \frac{A}{B} \right)^c$ pour une certaine constante $c$. Puisque 
Re$(s)>>0$, il d\'ecoule de (\ref{norme de psi}) et de la d\'emonstration du Lemme 12.3.3 que
\begin{eqnarray} \label{idelta}
|I_{\partial} | \prec \mbox{vol} (\partial B_t ) .
\end{eqnarray}
Finalement, il d\'ecoule facilement de (\ref{vol du bord})-(\ref{idelta}) que 
$$\int_{\Gamma_V \backslash {\cal D}} (\omega_s \wedge \pi^* (*_0 \mu ))\wedge \phi = \lim_{t \rightarrow \infty} \int_{B_t} (*_0 \mu ) \wedge \phi 
= \int_{C_V} (*_0 \mu ) \wedge \phi .$$

\bigskip

Nous notons~:
\begin{eqnarray} \label{Omega1}
\Omega (s) = \omega_s \wedge \pi^* (*_0 \mu ) ,
 \index{$\Omega (s)$}
\end{eqnarray}
que l'on appellera {\it forme ``duale''} \index{forme duale} associ\'ee \`a $\mu$ dans $M_V$.

\medskip

Dans le cas de l'espace hyperbolique complexe nous aurons besoin de nous assurer que la forme $\Omega (s)$
peut-\^etre choisie $\partial$ et $\overline{\partial}$ ferm\'ee. C'est l'objet de la section suivante.

\markboth{CHAPITRE 13. CONSTRUCTION DE LA FORME DUALE}{13.3. LE CAS HYPERBOLIQUE COMPLEXE}

\section{Pr\'ecisions dans le cas hyperbolique complexe}

Dans la section pr\'ec\'edente, on a construit un fibr\'e vectoriel $E_V$ au-dessus de $M_V = \Gamma_V \backslash {\cal D}$. 
On a de plus \'equip\'e ce fibr\'e de deux m\'etriques hermitiennes $\langle.,.\rangle$ et $\langle.,.\rangle^{\sim}$.

La classe de Chern du fibr\'e $E_V$ 
est un \'el\'ement ${\bf c} \in H^{2r} (M_V )$ dual au cycle $C_V$ dans $M_V$ (au sens de la premi\`ere section).
\'Etant donn\'e un rep\`ere mobile holomorphe local $e=(e_1 , \ldots ,e_r )$, soient 
\begin{eqnarray*}
H = (h_{ij} ) , \; h_{ij} = \langle e_i ,e_j \rangle 
\end{eqnarray*}
et
\begin{eqnarray*}
\tilde{H} = (\tilde{h}_{ij} ) , \; \tilde{h}_{ij} = \langle e_i ,e_j \rangle^{\sim} .
\end{eqnarray*}
Chaque structure hermitienne d\'etermine de mani\`ere unique une connexion de type $(1,0)$, dont les 
formes de connexion et de courbure sont, d'apr\`es la premi\`ere section, donn\'ees par~:
\begin{eqnarray*}
\omega = H^{-1} \partial H  \mbox{  et  } \Omega= \overline{\partial} (H^{-1} \partial H ) , \\
\tilde{\omega} = \tilde{H}^{-1} \partial \tilde{H} \mbox{  et  } \tilde{\Omega} = \overline{\partial} (\tilde{H}^{-1} \partial \tilde{H} ) .
\end{eqnarray*}

\medskip

Nous supposons dor\'enavant $q=1$ (autrement dit que ${\cal D}$ est l'espace hyperbolique complexe de dimension $p$). 
Les formes de connexions et de courbures pour ces deux m\'etriques \'evalu\'ees au point 
$\left( \begin{array}{c} 
               0 \\
               Z_2
               \end{array} \right) $ sont donc :
$$\begin{array}{ccc}
\left\{ 
\begin{array}{ccc}
        \omega & = & (0) I_r \\
        \Omega & = & (\partial \overline{\partial} \log B ) I_r
         \end{array} \right. & \mbox{ et } \left\{ \begin{array}{ccc}
                                                  \tilde{\omega} & = & \left\{ - \left( \frac{B}{A} \right) \partial \left( \frac{B}{A} \right) \right\} I_r \\
                                                  \tilde{\Omega} & = &(\partial \overline{\partial} \log A) I_r .                                                   
\end{array} \right.
\end{array}$$

Il est connu (cf. \cite{GriffithHarris}) que les $2r$-formes 
$$\left( \frac{\sqrt{-1}}{2\pi} \right)^r \mbox{det}\Omega $$
et 
$$\left( \frac{\sqrt{-1}}{2\pi} \right)^r \mbox{det}\tilde{\Omega} $$
repr\'esentent toutes deux la classe de Chern ${\bf c} \in H^{2r} (M_V )$. 

En suivant la m\'ethode de Bott et Chern rappel\'ee dans la premi\`ere section, Tong et Wang \cite{TongWang} construisent deux 
formes $\psi$ et $\tilde{\psi}$ sur $M_{V} -C_V$ telles que
\begin{enumerate}
\item on a :
$$\overline{\partial} \psi = \left( \frac{\sqrt{-1}}{2\pi} \right)^r \mbox{det}\Omega $$
et 
$$\overline{\partial} \tilde{\psi} = \left( \frac{\sqrt{-1}}{2\pi} \right)^r \mbox{det}\tilde{\Omega} ,$$
\item et en cohomologie : 
$$\overline{\partial} [\psi ] = {\bf c} -[F]$$
et 
$$\overline{\partial} [\tilde{\psi} ] = {\bf c} -[F].$$
\end{enumerate}

En suivant Tong et Wang, on remarque que si 
$$a_r =1-\left( 1-\frac{B}{A} \right)^r \mbox{  et  } C_r =\sqrt{-1} (2\pi )^r $$
alors
\begin{eqnarray} \label{1-2}
\begin{array}{ll}
\pi^* (*_0 \mu ) \wedge \psi -a_r \pi^* (*_0 \mu )\wedge \tilde{\psi}  & = 
 -\frac{1}{C_r } \sum_{t=0}^{r-1} \sum_{\lambda =0}^{r-1-t} (-1)^{r+t+\lambda +1} \left(      
\begin{array}{c}     
r \\     
t+\lambda      
\end{array}  \right) 
\left( \frac{A}{B} \right)^{t+\lambda -1-r} \\
& \pi^* (*_0 \mu ) \wedge \partial \left( \frac{A}{B} \right) (\partial \overline{\partial} \log A)^{t-1} (\partial \overline{\partial} \log B)^{r-t} 
\end{array}
\end{eqnarray}
est une forme non singuli\`ere.

\medskip

Pour pouvoir former la s\'erie th\'eta de cette forme nous voulons obtenir une forme $\partial$-ferm\'ee cohomologue \`a cette 
derni\`ere multipli\'ee par un poids $\left( \frac{A}{B} \right)^s $. Comme dans \cite{TongWang} on introduit la 
forme suivante~:
\begin{eqnarray} \label{phi (s)}
\begin{array}{ll}
\phi (s) & =   \frac{(\frac{A}{B})^s}{s} \pi^* (*_0 \mu )\wedge \psi \\
 &  + \left[ \sum_{\lambda =1}^r \frac{(\frac{A}{B} )^{s-\lambda }}{s-\lambda } (-1)^{\lambda } \left( 
\begin{array}{c}
       r \\
       \lambda        
\end{array}  \right) \right] \pi^* (*_0 \mu) \wedge \tilde{\psi}  \\
& +\frac{(-1)^{r-1}}{C_r} \sum_{t=1}^{r-1} \sum_{\lambda =0}^{r-1-t} (-1)^{t+\lambda } 
\left( \begin{array}{c}
       r \\
       t+\lambda
        \end{array}  \right) \\
&  \frac{\left( \frac{A}{B} \right)^{s+t+\lambda -1-r}}{s+t+\lambda -1-r} \pi^* (*_0 \mu ) \wedge 
\partial \left( \frac{A}{B} \right) (\partial \overline{\partial} \log A)^{t-1} (\partial \overline{\partial} \log B)^{r-t} .
\end{array}
\end{eqnarray}

La forme $\phi (s)$ est $\partial$-ferm\'ee et sa $\overline{\partial}$-d\'eriv\'ee dans $M_V - C_V$ s'\'etend de 
mani\`ere lisse \`a $M_V$. On a :
\begin{eqnarray} \label{delbar de phi}
\begin{array}{ll}
\overline{\partial} \phi (s)  & =  \frac{1}{C_r} \left\{ \frac{(\frac{A}{B} )^s}{s} \pi^* (*_0 \mu ) \wedge (\partial \overline{\partial} \log A)^r \right. \\
& +\left[ \sum_{\lambda =1}^r (-1)^{\lambda} \left(
                                               \begin{array}{c}
                                               r \\
                                               \lambda
                                                \end{array}
                                                   \right)  \frac{(\frac{A}{B} )^{s-\lambda}}{s-\lambda} \right] \pi^* (*_0 \mu ) \wedge (\partial \overline{\partial} \log B)^r  \\
& + (-1)^r \sum_{t=1}^{r-1} \sum_{\lambda =0}^{r-1-t} (-1)^{t+\lambda }\left(   
\begin{array}{c}
   r \\
   t+\lambda    
\end{array}      \right) \\
& \left.   \frac{(\frac{A}{B} )^{s+t+\lambda -1-r}}{s+t+\lambda-1-r} \pi^* (*_0 \mu) \wedge \partial \overline{\partial} \left( \frac{A}{B} \right) (\partial \overline{\partial} \log A)^{t-1} (\partial \overline{\partial} \log B)^{r-t} \right\} .
\end{array}
\end{eqnarray}

La d\'emonstration de la Proposition 12.4.3 implique le lemme suivant.

\begin{lem} \label{norme de phi n}
La norme $||\partial \overline{\partial} \log A ||$ est constante et $||\partial \overline{\partial} \log B ||$ est born\'ee par une constante. En particulier, on a~:
$$||\overline{\partial} \phi (s) ||, \ ||\phi (s) || \prec \left( \frac{A}{B} \right)^{{\rm Re}(s) -r+\frac{2p-2r-k}{2}} .$$
\end{lem}

\begin{prop} \label{forme duale}
Soit $\eta$ une $k$-forme ferm\'ee et born\'ee sur $M_V$.
Alors pour Re$(s)>2r+\frac{k}{2}$, la forme $\eta \wedge \overline{\partial} \phi (s)$ est int\'egrable sur $M_V$, et 
$$\int_{M_V} \eta \wedge \overline{\partial} \phi (s) = \kappa (s) \int_{C_V } \eta \wedge *_0 \mu $$
avec 
$$\kappa (s) = \frac{(-1)^r r!}{s(s-1) \cdots (s-r)}.$$
\end{prop}
{\it D\'emonstration.} On note toujours ${\cal F}_t$ l'image de ${\cal F}_t$ dans $M_V$. Sur $M_V -C_V$, 
$\eta \wedge \overline{\partial} \phi (s) = \overline{\partial} \{\eta \wedge \phi (s) \}$. La formule de Stokes implique 
donc : 
\begin{eqnarray*}
\int_{M_V} \eta \wedge \overline{\partial} \phi (s) & = & \lim_{\begin{array}{c}r  \rightarrow  0 \\R  \rightarrow  
+\infty \end{array} } \left\{ \int_{{\cal F} (R)} \eta \wedge \phi (s) - \int_{{\cal F} (r)} \eta \wedge \phi (s) \right\} .
\end{eqnarray*}

Or sur $C_V$, $\left( \frac{A }{B} \right) \equiv 1$ donc l'\'egalit\'e (\ref{phi (s)}) et la propri\'et\'e 2 de $\psi$ et 
$\tilde{\psi}$ impliquent~:
\begin{eqnarray*} 
\lim_{r \rightarrow 0} \int_{{\cal F} (r)} \eta \wedge  \phi (s) & = & -\left\{ \frac{1}{s} + \sum_{\lambda =1}^r (-1)^{\lambda} \frac{
\left( \begin{array}{c}
       r \\
       \lambda
       \end{array}  \right) }{s-\lambda} \right\} \int_{F} \eta \wedge *_0 \mu \\
& = & - \frac{(-1)^r r!}{s(s-1) \cdots (s-r)} \int_{F} \eta \wedge *_0 \mu .
\end{eqnarray*}

Il reste donc \`a montrer que $\lim_{R \rightarrow +\infty} \int_{{\cal F} (R)} \eta \wedge \phi (s) = 0$.
Mais pour Re$(s) \geq 2r+\frac{k}{2} +\varepsilon$ ($\varepsilon >0$) d'apr\`es la Proposition 12.2.5 et le Lemme 13.3.1, 
en un point de ${\cal F} (R)$ on a : 
$$|| \eta \wedge \phi (s) || \prec e^{-2(p+\varepsilon ) R}.$$
On conclut alors gr\^ace au Lemme \ref{comptage}. 

\bigskip

On d\'efinit donc :
\begin{eqnarray} \label{psi}
\Psi (s) & = & \frac{(-1)^r s(s-1) \cdots (s-r)}{r!} \overline{\partial} \phi (s) .
\end{eqnarray}

On peut montrer (cf. \cite{TongWang}) qu'\`a un cobord pr\`es dans la formule de 
$C_r \overline{\partial} \phi (s)$ donn\'ee par (\ref{delbar de phi}), on peut remplacer la double somme par 
$$(-1)^r \sum_{t=1}^{r-1} \sum_{\lambda =0}^{r-1-t} (-1)^{t+\lambda +1} \left(   
\begin{array}{c}
  r  \\
  t+\lambda   
\end{array}  \right) \frac{(\frac{A}{B} )^{s+t+\lambda -r}}{s+t+\lambda -r} \pi^* (*_0 \mu ) \wedge (\partial \overline{\partial} \log A)^{t-r} (\partial \overline{\partial} \log B)^{r-t+1} $$
$$+ (-1)^r \sum_{t=1}^{r-1} \sum_{\lambda =0}^{r-1-t} 
(-1)^{t+\lambda } \left(
                     \begin{array}{c}
                     r   \\
                     t+\lambda
                      \end{array}  \right) \frac{(\frac{A}{B})^{s+t+\lambda -r}}{s+t+\lambda-r} \pi^* (*_0 \mu ) \wedge (\partial \overline{\partial} \log A)^t (\partial \overline{\partial} \log B)^{r-t} .$$
On obtient ainsi une forme cohomologue \`a $\Psi (s)$~: 
\begin{eqnarray} \label{Omega2} \index{$\Omega (s)$}
\begin{array}{ll}
\Omega (s) &  =  \frac{(-1)^r s(s-1) \cdots (s-r)}{C_r r!} \sum_{\lambda =0}^r  \\
& (-1)^{\lambda} \left(  
\begin{array}{c}
  r  \\
  \lambda
\end{array}  \right) \frac{(\frac{A}{B} )^{s-\lambda }}{s-\lambda} \pi^* (*_0 \mu )\wedge (\partial \overline{\partial} \log A)^{r-\lambda} (\partial \overline{\partial} \log B)^{\lambda} ,
\end{array}
\end{eqnarray}
qui reste $\partial$ et $\overline{\partial}$ ferm\'ee.
Dans la suite nous appellerons $\Omega(s)$ la {\it forme ``duale''} associ\'ee \`a $\mu$ dans $M_V$.        \index{forme duale}

\bigskip

Concluons cette section en remarquant qu'il est \'egalement possible de construire une famille de formes duales 
$\partial$ et $\overline{\partial}$ ferm\'ees dans le cas g\'en\'eral des espaces sym\'etriques associ\'es au groupe
$SU(p,q)$. Esquissons cette construction. Revenons aux notations de la section pr\'ec\'edente. Notons $\psi_k$, $\psi , \ldots$ (resp. $\tilde{\psi}$,
$\tilde{\psi} , \ldots$) les formes construites en suivant la premi\`ere section et pour la m\'etrique hermitienne $\langle .,. \rangle$ 
(resp. $\langle .,.\rangle^{\sim}$). Il est facile de v\'erifier que~:
\begin{eqnarray} \label{psi - psi tilde}
\tilde{\psi} = \psi - \sum_{k=1}^q  \left( 
\begin{array}{c}
q \\
k 
\end{array} \right) \left( \frac{C}{A} \right)^k \psi_k .
\end{eqnarray}

En d\'eveloppant la formule (\ref{psi - psi tilde}), on obtient~:
\begin{eqnarray} \label{devpt de psi}
\begin{array}{ll}
\tilde{\psi} & =\left[ \sum_{j=1}^q \left(
\begin{array}{c}
q \\
j
\end{array} \right) (-1)^j \left( \frac{B}{A} \right)^j \right] \psi  \\
& - \sum_{k=1}^{q-1} \left( \left( 
\begin{array}{c}
q \\
k
\end{array} \right) \sum_{l=1}^{q-k} \left( 
\begin{array}{c}
q-k \\
l
\end{array} \right) (-1)^l \left( \frac{B}{A} \right)^{l+k} \right) \left( 
\left( \frac{C}{B} \right)^k \psi_k \right) .
\end{array}
\end{eqnarray}
Ce qui nous conduit \`a d\'efinir la forme d\'ependant d'un param\`etre $s \in {\Bbb C}$~:
\begin{eqnarray} \label{Phis}
\begin{array}{ll}
\Phi (s) & = \frac{ (A/B)^{s+q}}{s+q} \tilde{C} (E)  - \left[ \sum_{j=1}^q \left( 
\begin{array}{c} 
q \\
j 
\end{array} \right) (-1)^j \frac{ (A/B)^{s+q-j}}{s+q-j} \right] C (E)  \\
& + \sum_{k=1}^{q-1} \left( \left( 
\begin{array}{c}
q \\
k
\end{array} \right) \sum_{l=1}^{q-k} \left( 
\begin{array}{c}
q-k \\
l 
\end{array} \right) (-1)^l \frac{ (A/B)^{s+q-l-k}}{s+q-l-k}  \right) \overline{\partial} \left( \left( 
\frac{C}{B} \right)^k \psi_k \right) .
\end{array}
\end{eqnarray}
Pour Re$(s)\geq 0$, la forme $\Phi (s)$ est lisse sur ${\cal D}$ et $\partial$ et $\overline{\partial}$ ferm\'ee. 
Nous montrons finalement~:

\begin{prop} \label{forme duale 3}
Soit $\eta$ une $(q(p-1), q(p-1))$-forme $\partial$ et $\overline{\partial}$ ferm\'ee sur $\Gamma_V \backslash {\cal D}$
de norme born\'ee. Alors, pour Re$(s) > >1$,
$$\int_{\Gamma_V \backslash {\cal D}} \Phi (s) \wedge \eta = \kappa (s) \int_{\Gamma_V \backslash {\cal D}_V} \eta ,$$
o\`u $\kappa (s)=  \frac{1}{s+q} - \sum_{j=1}^q \left( 
\begin{array}{c}
q \\
j
\end{array} \right) \frac{(-1)^j}{s+q-j}$.
\end{prop}
{\it D\'emonstration.} Soit
$$\psi (s) = \frac{(A/B)^{s+q}}{s+q} \tilde{\psi} - \left[ \sum_{j=1}^q \left( 
\begin{array}{c}
q \\
j
\end{array} \right) (-1)^j \frac{(A/B)^{s+q-j}}{s+q-j} \right] \psi $$
$$+ \sum_{k=1}^{q-1} \left( \left( 
\begin{array}{c}
q \\
k
\end{array} \right) \sum_{l=1}^{q-k} \left( 
\begin{array}{c}
q-k \\
l
\end{array} \right) (-1)^l \frac{(A/B)^{s+q-l-k}}{s+q-l-k} \right) \left( \left(
\frac{C}{B} \right)^k \psi_k \right) .$$
D'apr\`es (\ref{devpt de psi}) on a $\overline{\partial} \psi (s) =\Phi (s)$ sur ${\cal D} \setminus {\cal D}_V$. Soit 
$T(r)$ le voisinage tubulaire de rayon $r$ de ${\cal D}_V$ dans ${\cal D}$. On a~:
$$\int_{\Gamma_V \backslash {\cal D}} \Phi (s) \wedge \eta = \lim_{
\begin{array}{c}
r \rightarrow +\infty \\
\varepsilon \rightarrow 0
\end{array}} \left[ \int_{\partial T(r)} \psi (s) \wedge \eta - \int_{\partial T(\varepsilon) } \psi (s) \wedge \eta \right] .$$
Mais $||\psi (s) \wedge \eta || \prec (A/B)^{{\rm Re}(s)}$ et, d'apr\`es la Proposition 12.2.5, pour Re$(s)>>1$, 
$$\lim_{r \rightarrow +\infty} \int_{\partial T(r)} \psi (s) \wedge \eta = 0 .$$
Enfin, d'apr\`es la d\'emonstration de la Proposition \ref{forme duale}, 
$$\lim_{\varepsilon \rightarrow 0} \int_{\partial T(\varepsilon )} \psi (s) \wedge \eta = -\kappa (s) \int_{\Gamma_V 
\backslash {\cal D}_V} \eta .$$
Ce qui conclut la d\'emonstration de la Proposition \ref{forme duale 3}.

\bigskip

On pourrait l\`a encore former une s\'erie de Poincar\'e et obtenir une famille de formes duales $\Omega (s)$
$\partial$ et $\overline{\partial}$ ferm\'ees m\^eme dans le cas de l'espace sym\'etrique associ\'e au groupe 
$SU(p,q)$.

\markboth{CHAPITRE 13. CONSTRUCTION DE LA FORME DUALE}{13.4. TOURS DE REV\^ETEMENTS FINIS}

\section{Tours de rev\^etements finis}

Nous conservons les notations des sections pr\'ec\'edentes. Nous supposons de plus la vari\'et\'e $M$ compacte. 
D'apr\`es la Proposition 12.5.6 et le Lemme 13.3.1, pour Re$(s) >>1$ la s\'erie  
\begin{eqnarray} \label{omegas}
\omega_s^m         & = & \sum_{\gamma \in \Gamma_V \backslash \Gamma_m  } \gamma^* \Omega (s) 
\end{eqnarray}
converge uniform\'ement sur tout compact de ${\cal D}$ et d\'efinit une $(2pq-k)$-forme ferm\'ee sur 
$M_m$. Le statut de ``forme duale'' de $\Omega (s)$ implique que pour toute $k$-forme ferm\'ee $\eta $ sur $M_m$ on a~:
$$\int_{M_m} \omega_s^m \wedge \eta = \int_{M_V} \Omega (s) \wedge \eta = \int_{C_V} *_0 \mu \wedge \eta .$$

On obtient donc le th\'eor\`eme suivant.

\begin{thm} \label{forme duale 2}
Soit $c$ un cycle de dimension $k$ dans $C_V$. Soit $\mu$ la $k$-forme harmonique sur $C_V$ dont l'\'etoile de Hodge est 
duale \`a $c$. Alors pour Re$(s) >>1$, les formes $\omega_s^m$ sur $M_m$ forment une famille de formes 
ferm\'ees duales au cycle $c$ dans $M_m$. 
\end{thm}

Soit $\Delta = \delta d +d \delta$, o\`u $\delta$ est l'adjoint de $d$, l'op\'erateur laplacien que l'on \'etend en un 
op\'erateur, toujours not\'e $\Delta$, agissant sur l'espace $L^2 \Omega^{2pq-k} ({\cal D} )$ des 
$(2pq-k)$-formes de carr\'e int\'egrable sur ${\cal D}$ de fa\c{c}on essentiellement auto-adjointe. Alors le 
Th\'eor\`eme spectral s'applique et il existe une famille spectrale $\{ P_{\lambda} \; : \; \lambda \in [0, +\infty [\}$ associ\'ee 
\`a $\Delta$.

Notons $P_{\lambda} (x,y)$ le noyau de Schwartz de $P_{\lambda}$. 
On a $\Delta = \int_0^{+\infty} \lambda dP_{\lambda}$. \`A toute fonction $f\in C_0 ([0,+\infty [)$, on associe l'op\'erateur 
$$f(\Delta ) = \int_0^{+\infty} f(\lambda )dP_{\lambda }.$$
Le laplacien est un op\'erateur elliptique. Soit $\omega \in L^2 \Omega^{2n-p} ({\cal D} )$.
On a :
$$\Delta (f(\Delta ) \omega ) = F(\Delta )\omega $$
o\`u $F$ est la fonction qui \`a $x$ associe $xf(x)$. D'apr\`es le Th\'eor\`eme de r\'egularit\'e sur les op\'erateurs 
elliptiques, la forme $f(\Delta )\omega$ est lisse. De plus, pour tout $x \in {\cal D}$, il existe une 
constante $C(x,f,{\cal D} )$ telle que~:
$$|f(\Delta )\omega |(x) \leq C(x,f,{\cal D} ) ||\omega ||_{L^2 ({\cal D} )} .$$
En particulier, l'application 
$$\left\{
\begin{array}{ccc}
L^2 \Omega^{2pq-k} ({\cal D} ) & \rightarrow & L^2 \Omega^{2pq-k}_x ({\cal D} ) \\
\omega & \mapsto & f(\Delta )\omega (x) 
\end{array}
\right. $$
est continue. D'apr\`es le Th\'eor\`eme de Riesz, il existe donc 
$f(\Delta )(x,.) \in L^2 \Omega^{2pq-k} ({\cal D} )$ tel que 
$$f(\Delta )\omega (x) = \int_{{\cal D}} f(\Delta )(x,y) \omega (y) dy .$$
De plus, pour tout compact $K$ de ${\cal D}$, 
$\int_{K \times {\cal D}} ||f(\Delta ) (x,y)||^{2} dxdy \leq \int_K C(x,f,{\cal D})^2 dx $. 
Donc $f(\Delta )(.,.) \in L^2_{loc} ({\cal D} \times {\cal D} )$. Or 
$$(\Delta_x + \Delta_y )f(\Delta )(x,y) = 2F(\Delta ) (x,y)$$
et l'op\'erateur $\Delta_x + \Delta_y$ est elliptique. Le Th\'eor\`eme de r\'egularit\'e elliptique implique donc que 
$f(\Delta ) (x,y)$ est une fonction $C^{\infty}$ en $x$ et $y$.

\medskip

Sur les vari\'et\'es $M_m$ et $M_{V}$, nous notons le laplacien respectivement $\Delta_m$ et $\Delta_{\infty}$; 
de m\^eme nous notons respectivement $P_{\lambda}^m$ et $P_{\lambda}^{\infty}$ les familles spectrales associ\'ees. 
Si l'on d\'esigne, de mani\`ere coh\'erente avec les notations pr\'ec\'edentes,
le noyau de la chaleur sur les $(2pq-k)$-formes de ${\cal D}$ par $e^{-t\Delta} (x,y)$, il est connu \cite{Donnelly} 
que pour tout $T>0$, il existe une constante $\alpha >0$ et une constante $C_T$ (d\'ependante de $T$) telles que~:
\begin{eqnarray} \label{noyau de la chaleur}
|e^{-t\Delta} (x,y)| & \leq & C_T e^{-\alpha d(x,y)^2 /t } ,
\end{eqnarray}
pour tout $t\in ]0,T]$. 

\begin{lem} \label{noyau de la chaleur 2}
Pour $t>0$ fix\'e, la s\'erie 
$$\sum_{\gamma \in \Gamma} |e^{-t\Delta} (x,\gamma y)|$$
converge uniform\'ement pour $x \in {\cal D}$ et $y$ dans un compact vers une fonction born\'ee.
\end{lem}
{\it D\'emonstration.} Soit $K$ un compact de ${\cal D}$ et soit $t$ un r\'eel strictement positif. D'apr\`es 
le Lemme 12.5.4, il existe une constante $c_1 (K)$ telle que 
$$\nu(x,y,R ) := |\{ \gamma \in \Gamma \; : \;  d(x,\gamma y) \leq R \} |\leq c_1 (K) \int_0^{R+1} (1+u^{2pq} ) e^{2(p+q-1)\sqrt{q} u} du$$
pour tout $x \in {\cal D}$ et $y \in K$. 

Soient $C=C_t$ et $\beta = \alpha /t$. Alors, d'apr\`es l'in\'egalit\'e (\ref{noyau de la chaleur}), pour 
$x \in {\cal D}$ et $y \in K$, on a~:
\begin{eqnarray*}
\sum_{ 
\begin{array}{c}
\gamma \in \Gamma \\
d(x,\gamma y) \leq R
\end{array} 
}
|e^{-t\Delta } (x,\gamma y )| & \leq & 
C \sum_{ 
\begin{array}{c}
\gamma \in \Gamma \\
d(x,\gamma y) \leq R
\end{array} 
}
e^{-\beta d(x,\gamma y )^2 } \\
& \leq & C \left( \int_0^R e^{-\beta r^2} d\nu (x,y,r) \right) \\
& \leq & C \left( [e^{-\beta r^2 } \nu (x,y,r) ]_0^R + 2 \beta \int_0^R re^{-\beta r^2} \nu (x,y,r) dr \right) \\
& \leq & C c_2 (K) \left( e^{-\beta R^2} \int_0^{R+1} (1+u^{2pq} )e^{2(p+q-1) \sqrt{q} u} du  \right. \\
&        & \left. +2\beta \int_0^R re^{-\beta r^2} \int_0^{r+1} (1+u^{2pq} ) e^{2(p+q-1) \sqrt{q} u} du dr \right) .
\end{eqnarray*}
Le Lemme d\'ecoule de ces in\'egalit\'es en faisant tendre $R$ vers l'infini et du fait que la vari\'et\'e $M$ est compacte.  

\bigskip

On obtient alors que pour tout $t>0$,
$$e^{-t\Delta_m } (x,y) = \sum_{\gamma \in \Gamma_m} (\gamma_y )^* e^{-t\Delta } (x,y) $$
$$e^{-t\Delta_{\infty}} (x,y) = \sum_{\gamma \in \Gamma_V } (\gamma_y )^* e^{-t\Delta } (x,y) .$$
Et la convergence est absolue et uniforme pour $x$, $y$ dans un compact. En particulier, le noyau de la chaleur 
$e^{-t\Delta_m } (x,y)$ est $\Gamma_m$-invariant. De plus, puisque d'apr\`es le Lemme \ref{effeuillage}, 
$\cap_m \Gamma_m = \Gamma_V$ et $\Gamma_{m+1} \subset \Gamma_m$, on en d\'eduit que si $t>0$ est fix\'e, 
$$e^{-t\Delta_{\infty}} (x,y) = \lim_{m\rightarrow +\infty} e^{-t\Delta_m} (x,y) $$
uniform\'ement pour $x$ et $y$ dans un compact de ${\cal D}$. 

\medskip

Le lemme suivant est une cons\'equence du Th\'eor\`eme d'approximation de Weierstrass.

\begin{lem}[cf. \cite{Donnelly2}]  \label{approx par des polynomes}   
Soit $f \in C_0 ([0, +\infty [)$. Alors $f$ peut-\^etre uniform\'ement approch\'ee sur $[0, +\infty[$ par une combinaison     
lin\'eaire finie d'exponentielles $e^{-tx}$, $t>0$.
\end{lem}
{\it D\'emonstration.} On se restreint d'abord \`a l'intervalle $]0, +\infty [$ et on effectue le changement de variable 
$y=e^{-x}$. Soit $g(y) = f(x)$. Alors $g \in C_0 (]0,1 [)$.

Soit $\varepsilon >0$, le Th\'eor\`eme d'approximation de Weierstrass donne un polyn\^ome proche de $g$ :
$$|g(y) - \sum_{j=0}^n a_j y^j | < \varepsilon /2 .$$
Puisque $g(0)=0$, on a $|a_0| <\varepsilon /2$. Ainsi $|g(y) - \sum_{j=1}^n a_j y^j | < \varepsilon $. Le changement de 
variable inverse de $y$ \`a $x$ implique que $|f(x) -\sum_{j=1}^n a_j e^{-jx} |<\varepsilon $.  

\bigskip

Comme Donnelly dans \cite{Donnelly}, montrons que ce lemme implique que~:
$$f(\Delta_{\infty} )(x,y) = \lim_{m \rightarrow +\infty } f(\Delta_m ) (x,y) ,$$
pour tout $f \in C_0 ([0,+\infty [)$ et uniform\'ement pour $x$, $y$ dans un compact. 

Soit $f \in C_0 ([0, +\infty [)$. Soit $\varepsilon >0$. Posons $g(x) =f(x) e^x$ pour $x \in [0, +\infty [$. Bien s\^ur 
$g \in C_0 ([0, +\infty [)$. D'apr\`es le Lemme \ref{approx par des polynomes}, il existe une suite $g_l (x)$ de 
polyn\^omes en $e^{-x}$ qui approche uniform\'ement $g(x)$. Alors pour $x$, $y \in {\cal D}$ et 
$m \in {\Bbb N}$, on a~:
$$|f(\Delta_{\infty} )(x,y) -f(\Delta_m ) (x,y) |\leq A_1 +A_2 +A_3 $$
o\`u 
$$A_1 = |f(\Delta_{\infty} )(x,y) -(g_l (\Delta_{\infty} )e^{-\Delta_{\infty} }) (x,y) | ,$$
$$A_2 = |(g_l (\Delta_{\infty} )e^{-\Delta_{\infty}} ) (x,y) -(g_l (\Delta_{m} )e^{-\Delta_{m} }) (x,y) | $$ 
et
$$A_3 = |(g_l (\Delta_{m} )e^{-\Delta_{m} }) (x,y) - f(\Delta_m )(x,y) |.$$

Le Th\'eor\`eme spectral pour $\Delta_{\infty}$ implique que $g_l (\Delta_{\infty} )$ converge fortement vers 
$g(\Delta_{\infty} )$ ({\it i.e.} $g_l (\Delta_{\infty} ) \omega \rightarrow g (\Delta_{\infty} )\omega $ pour tout 
$\omega \in L^2 \Omega^{2pq-k} (M_{V})$). Les noyaux de Schwartz $g_l (\Delta_{\infty} )e^{-\Delta_{\infty}} (x,y)$ 
convergent donc dans $L^2_{loc}$ vers $g(\Delta_{\infty } )e^{-\Delta_{\infty} } (x,y) = f(\Delta_{\infty}) (x,y)$. Et le 
Th\'eor\`eme de r\'egularit\'e elliptique implique que l'on peut choisir $l$ suffisamment grand pour que 
$A_1 \leq \varepsilon /3$. 

Notons maintenant $\{ f_i^m \} _{i \geq 0}$ une base orthonorm\'ee de $L^2 \Omega^{2pq-k} (M_{m})$ constitu\'ee de 
formes propres pour $\Delta_m$ et $\{ \lambda_i^m \}_{i \geq 0}$ la suite des valeurs propres (compt\'ees avec 
multiplicit\'ee) associ\'ees. Alors~:
$$h(\Delta_m ) = \sum_i h(\lambda_i^m )f_i^m (x) \otimes f_i^m (y) ,$$
pour toute fonction continue $h$ \`a support compact dans ${\Bbb R}_+$. Donc :
\begin{eqnarray*}    
A_3 & \leq & (\sup_{\lambda \in {\Bbb R}_+ } |g(\lambda ) - g_l     (\lambda ) |) \times \sum_i e^{-\lambda_i^m } |f_i^m (x)||f_i^m     (y)| \\
    & \leq & (\sup_{\lambda } |g(\lambda ) - g_l     (\lambda ) |) \times \sqrt{tr(e^{-\Delta_m} (x,x) )}     \sqrt{tr(e^{-\Delta_m} (y,y) )} \\
    & & \mbox{(d'apr\`es l'in\'egalit\'e de Cauchy-Schwarz)} \\
    & \leq & K_0 (\sup_{\lambda } |g(\lambda ) - g_l     (\lambda ) | ) 
\end{eqnarray*}
o\`u $K_0$ est une constante ind\'ependante de $m$ que l'on obtient gr\^ace \`a l'estim\'ee (\ref{noyau de la chaleur}) comme dans la preuve 
du Lemme \ref{noyau de la chaleur 2} (on renvoie \`a \cite{Donnelly2} pour plus de d\'etails). Ainsi, on peut choisir $l$ suffisamment grand  
de mani\`ere \`a ce que $A_3 \leq \varepsilon /3 $ pour tout $m$. 

Fixons maintenant $l$ de mani\`ere \`a ce que 
$A_1 +A_3 \leq \varepsilon /3$. D'apr\`es le Lemme \ref{effeuillage} et le Lemme \ref{noyau de la chaleur 2}, on a 
$A_2 \leq \varepsilon /3$ pour $m$ suffisamment grand.

En conclusion, comme Donnelly dans \cite{Donnelly}, on a montr\'e que :
$$f(\Delta_{\infty} )(x,y) = \lim_{m \rightarrow +\infty } f(\Delta_m ) (x,y) ,$$
pour tout $f\in C_0 ([0,+\infty [)$ et uniform\'ement pour $x$, $y$ dans un compact.

\begin{lem} \label{convergence des noyaux}
Pour tout $f\in C_0 ([0,+\infty [)$, la suite $ f(\Delta_m )(x,y)$ converge vers $f(\Delta_{\infty} )(x,y)$ 
uniform\'ement pour $x$ et $y$ dans un compact. Et l'expression
$$| f(\Delta_m )(x, y)) -f(\Delta_{\infty} )(x,y) |$$
est uniform\'ement born\'ee (ind\'ependamment de $m$) pour $x\in {\cal D}$ et $y$ dans un compact.
\end{lem}
{\it D\'emonstration.} En effet, pour tout $x$, $y \in {\cal D}$ et $t >0$, on a~:
$$|e^{-t\Delta_{\infty}} (x,y) - e^{-t\Delta_m } (x,y)| \leq \sum_{\gamma \in \Gamma_m -\Gamma_V } |e^{-t\Delta } (x, \gamma y)| .$$

Puisque $\cap_m \Gamma_m =\Gamma_V$, le Lemme \ref{convergence des noyaux} pour la fonction $f(.) =e^{-t.}$ d\'ecoule du Lemme \ref{noyau de la chaleur 2}. 
On conclut la preuve (comme ci-dessus ou dans \cite{Donnelly}) \`a l'aide du Lemme \ref{approx par des polynomes}.  

\bigskip

Remarquons que, puisque $ e^{-t\Delta_m } (.,.)$ est $\Gamma_m$-bi-invariant, il en est de m\^eme pour 
$f(\Delta_m )(.,.)$.
\begin{prop} \label{convergence des noyaux 2}
Soient $f \in C_0 ([0,+\infty [)$ et $s\in {\Bbb C}$, ${\rm Re}(s) >>1$. Alors, la suite 
$f (\Delta_{m} ) \omega^m_s  $ converge vers $f(\Delta_{\infty } ) \Omega (s)$ uniform\'ement sur les compacts.
\end{prop}
{\it D\'emonstration.} Soit $s\in {\Bbb C}$, Re$(s) >>1$. Si ${\cal F}_m$ est un domaine fondamental pour 
l'action de $\Gamma_m$ sur ${\cal D}$, on a~:
\begin{eqnarray*}    
f(\Delta_m )\omega_s^m (.) & = & \int_{M_m} \omega_s^m (x) \wedge     *f(\Delta_m ) (x,.) dx \\
    & = & \int_{M_m} \left( \sum_{\gamma \in \Gamma_V \backslash     \Gamma_m } \gamma^* \Omega (s) (x) \right) \wedge *f(\Delta_m     )(x,.) dx \\
    & = & \sum_{\gamma \in \Gamma_V \backslash     \Gamma_m } \int_{{\cal F}_m}  \gamma^* \Omega (s) (x) \wedge *f(\Delta_m     )(x,.) dx \\
      & = & \sum_{\gamma \in \Gamma_V \backslash     \Gamma_m } \int_{\gamma {\cal F}_m} \Omega (s) (x) \wedge *f(\Delta_m     )(x,.) dx \\
    & & \mbox{(car $f(\Delta_m )(.,.)$ est $\Gamma_m$-bi-invariant)} \\
    & = & \int_{M_{V}} \Omega (s) (x) \wedge *f(\Delta_m     )(x,.) dx .
\end{eqnarray*}

On obtient donc sur ${\cal D}$~:
\begin{eqnarray*}    
f(\Delta_m )\omega_s^m (.) - f(\Delta_{\infty } ) \Omega (s) (.) & = & \int_{M_{V}} \Omega (s) (x) \wedge *f(\Delta_m ) (x,.) dx \\
&  & - \int_{M_{V}} \Omega (s) (x) \wedge *f(\Delta_{\infty} ) (x,.) dx \\
& = & \int_{M_{V}} \Omega (s) (x) \wedge *(f(\Delta_m ) (x,.) -f(\Delta_{\infty} ) (x,.)) dx .
\end{eqnarray*}

De plus d'apr\`es le Lemme \ref{convergence des noyaux}, l'expression 
$$|f(\Delta_m ) (x,y) -f(\Delta_{\infty} ) (x,y)|$$
est uniform\'ement born\'ee (ind\'ependamment de $m$) pour $x\in M_{V}$ et $y$ dans un compact.
Puisque la forme $\Omega (s)$ est dans $L^1$, le Th\'eor\`eme de convergence domin\'ee et le Lemme \ref{convergence des noyaux} impliquent 
que pour tout r\'eel $\varepsilon >0$ et pour tout compact $K$, il existe un entier $m_0$ tel que pour tout $m\geq m_0$, 
les applications $f(\Delta_m )\omega_s^m$ et $f(\Delta_{\infty } ) \Omega (s)$ sont $\varepsilon$-proches sur $K$. 
Ce qui ach\`eve la d\'emonstration de la Proposition \ref{convergence des noyaux 2}.  

\bigskip

\newpage 

\thispagestyle{empty}

\newpage

\markboth{CHAPITRE 14. COHOMOLOGIE $L^2$ R\'EDUITE}{14.1. RAPPELS SUR LA COHOMOLOGIE $L^2$}

\chapter{Cohomologie $L^2$ r\'eduite}

Nous conservons dans ce chapitre les notations des chapitres pr\'ec\'edents. Soient donc toujours $G = SU(p,q)$ ($p\geq q$), ${\cal D}$
l'espace sym\'etrique associ\'e, ${\cal D}_V$ le sous-espace totalement g\'eod\'esique de ${\cal D}$ associ\'e \`a un 
sous-espace vectoriel de dimension (complexe) $r$ de ${\Bbb C}^{p+q}$, $G_V$ le sous-groupe de $G$ pr\'eservant 
${\cal D}_V$ et $\Gamma_V$ un sous-groupe discret sans torsion et cocompact dans $G_V$. Nous notons $C_V = \Gamma_V 
\backslash {\cal D}_V$ et $M_V = \Gamma_V \backslash {\cal D}$. Dans la suite, nous notons $H^k_c (M_V )$, $H^k_2 (M_V )$ et ${\cal H}^k_{2} (M_V )$ 
d\'esignent respectivement le groupe de cohomologie \`a support compact de degr\'e $k$ de $M_V $, le 
groupe de cohomologie $L^2$ r\'eduite de degr\'e $k$ de $M_V$ et l'espace des $k$-formes harmoniques $L^2$ de $M_V$.
Le but de ce chapitre est la d\'emonstration du th\'eor\`eme suivant.

\begin{thm} \label{cohom l2}
Pour tout entier, $k < p+qr-r$, on a les isomorphismes naturels suivants~:
$$H^{k-2qr} (C_V ) \cong H_c^k (M_V ) \stackrel{\simeq}{\rightarrow} H_{2}^k (M_V ) \cong {\cal H}_2^k (M_V ) .$$
Si de plus, $q=1$, l'espace $H_{2}^p (M_V ) \cong {\cal H}_2^{p} (M_V )$ est de dimension infini et 
l'application naturelle $(H^{p-r} (C_V ) \cong ) H_c^p (M_V ) \rightarrow H_2^k (M_V)$ est injective.
\end{thm}

La d\'emonstration pour $q=1$ annonc\'ee dans \cite{BergeronClozel} reposait sur une proposition de Donnelly et 
Fefferman. La d\'emonstration de cette derni\`ere telle qu'esquiss\'ee dans \cite{DonnellyFefferman} est fausse comme
nous l'ont signal\'es Gilles Carron et Nader Yeganefar \`a qui nous devont \'egalement la r\'ef\'erence \cite{OhsawaTakegoshi}
qui nous permet dans ce chapitre de d\'emontrer compl\`etement le Th\'eor\`eme \ref{cohom l2}. Remarquons 
n\'eanmoins que le r\'esultat principal de \cite{OhsawaTakegoshi} repose encore sur un \'enonc\'e 
(compl\`etement d\'emontr\'e cette fois) tir\'e de \cite{DonnellyFefferman}. Remarquons enfin que le cas $q=1$ peut maintenant 
\^etre d\'eduit d'un r\'esultat de Yeganefar \cite{Yeganefar}.
 
Commen\c{c}ons par des rappels sur la cohomologie $L^2$.

\section{Rappels sur la cohomologie $L^2$}

Soit $M$ une vari\'et\'e complexe hermitienne compl\`ete. Nous notons $C_0^{\infty} (\Lambda^k T^* M)$ (resp. 
$L^2 (\Lambda^k T^* M)$, etc...) l'ensemble des $k$-formes lisses \`a support compact (resp. de carr\'e 
int\'egrable, etc...) dans $M$. De m\^eme nous notons $C_0^{\infty} (\Lambda^{a,b} T^* M)$ (resp. 
$L^2 (\Lambda^{a,b} T^* M)$, etc...) l'ensemble des formes lisse de bidegr\'e $(a,b)$ \`a support compact (resp.
de carr\'e int\'egrable, etc...) dans $M$. 

\medskip

Commen\c{c}ons par oublier la structure complexe de $M$ et par ne voir $M$ que comme une vari\'et\'e riemannienne
compl\`ete.
Le $k$-i\`eme espace de cohomologie $L^2$ (r\'eduite) de $M$ est d\'efini par 
$$H_2^k (M) = \{ \alpha \in L^2 (\Lambda^k T^* M) \; : \; d\alpha =0 \} / \overline{dC_0^{\infty} (\Lambda^{k-1} T^* M)} ^{L^2} .$$
Un autre espace tr\`es proche souvent consid\'er\'e est l'espace de cohomologie $L^2$ non r\'eduite, qui, en degr\'e $k$, est le 
quotient de $\{ \alpha \in L^2 (\Lambda^k T^* M) \; : \; d\alpha =0 \}$ par  $\{ d\alpha \; : \; \alpha \in L^2 (\Lambda^{k-1} T^* M), \; d\alpha \in L^2 \}$,
sans prendre d'adh\'erence. En g\'en\'eral, cohomologies $L^2$ r\'eduite et non r\'eduite sont diff\'erentes.
Il y a n\'eanmoins \'egalit\'e en degr\'e $k$ lorsque $0$ n'est pas dans le spectre essentiel du laplacien $\Delta$ sur 
les formes diff\'erentielles de degr\'e $k$. Dans la suite, ``cohomologie $L^2$'' voudra dire ``cohomologie $L^2$ r\'eduite''.

Il y a une interpr\'etation de la cohomologie $L^2$ en termes de formes harmoniques. En effet, notons ${\cal H}^k_{2}$ l'espace
des $k$-formes harmoniques $L^2$ de $M$~:
$${\cal H}^k_2 (M) = \{ \alpha \in L^2 (\Lambda^k T^* M) \; : \; d\alpha = \delta \alpha =0 \}$$
o\`u $\delta$ est l'op\'erateur d\'efini initialement sur les formes lisses \`a support compact comme l'adjoint de $d$. 
Comme $M$ est compl\`ete, ${\cal H}_2^* (M)$ est aussi le noyau $L^2$ du laplacien $\Delta = d\delta +\delta d$. Un fait important est la d\'ecomposition
de Hodge-de Rham-Kodaira \cite{GriffithHarris}~:
$$L^2 (\Lambda^k T^*M )={\cal H}^k_2 (M) \oplus \overline{dC_0^{\infty} (\Lambda^{k-1} T^*M )} \oplus \overline{\delta C_0^{\infty} (\Lambda^{k+1} T^* M)} ,$$
et de plus,
$$\{ \alpha \in L^2 (\Lambda^k T^*M ) \; : \; d\alpha =0 \} = {\cal H}^k_2 (M) \oplus \overline{dC_0^{\infty} (\Lambda^{k-1} T^*M )}.$$
On en d\'eduit que 
$$H_2^k (M) \cong {\cal H}_2^k (M) .$$

\medskip

Les r\'esultats pr\'ec\'edents admettent bien entendu des analogues complexes. Rappelons qu'une vari\'et\'e complexe 
hermitienne est dite {\it compl\`ete} \index{compl\`ete, vari\'et\'e hermitienne} si la vari\'et\'e riemannienne sous-jacente est compl\`ete. En rempla\c{c}ant le 
complexe de de Rham $(K^* ,d)$ d\'efini par 
$$K^k = \{ \alpha \in L^2 (\Lambda^k T^* M) \; : \; d\alpha \in L^2 \} $$
par le complexe de Dolbeault $(A^* , \overline{\partial})$ d\'efini par $A^* = \oplus A^{a,b}$ o\`u
$$A^{a,b} = \{ \alpha \in L^2 (\Lambda^{a,b} T^* M) \; : \; \overline{\partial} \alpha \in L^2 \} ,$$
on d\'efinit les groupes de cohomologie $L^2$ de bidegr\'e $(a,b)$. Autrement dit,
$$H_2^{a,b} (M) = \{ \alpha \in L^2 (\Lambda^{a,b} T^* M) \; : \; \overline{\partial}\alpha =0 \} / \overline{\overline{\partial}C_0^{\infty} (\Lambda^{a,b-1} T^* M)} ^{L^2} .$$
Il y a l\`a encore une interpr\'etation de la cohomologie $L^2$ en termes de formes harmoniques. Mais cette fois on
utilise le laplacien complexe $\square$. Notons ${\cal H}_2^{a,b}$ l'espace des formes harmoniques $L^2$ de bidegr\'e
$(a,b)$ de $M$~:
$${\cal H}_2^{a,b} (M) = \{ \alpha \in L^2 (\Lambda^{a,b} T^* M) \; : \; \overline{\partial} \alpha = \overline{\partial}^* \alpha =0 \} ,$$
o\`u $\overline{\partial}^*$ est l'op\'erateur d\'efini intitialement sur les formes lisses \`a support compact
comme l'adjoint de $\overline{\partial}$. Comme $M$ est compl\`ete, ${\cal H}^*_2$ est aussi le noyau $L^2$ du laplacien
complexe $\square = \overline{\partial} \overline{\partial}^* + \overline{\partial}^* \overline{\partial}$. L\`a encore~:
$$L^2 (\Lambda^{a,b} T^*M )={\cal H}^{a,b}_2 (M) \oplus \overline{\overline{\partial} C_0^{\infty} (\Lambda^{a,b-1} T^*M )} \oplus \overline{\overline{\partial}^* C_0^{\infty} (\Lambda^{a,b+1} T^* M)} ,$$
et de plus,
$$\{ \alpha \in L^2 (\Lambda^{a,b} T^*M ) \; : \; d\alpha =0 \} = {\cal H}^{a,b}_2 (M) \oplus \overline{\overline{\partial}C_0^{\infty} (\Lambda^{a,b-1} T^*M )}.$$
On en d\'eduit que 
$$H_2^{a,b} (M) \cong {\cal H}_2^{a,b} (M) .$$

\medskip

Rappelons \cite{GriffithHarris} que lorsque $M$ est une vari\'et\'e kaehl\'erienne, le laplacien de Hodge-de Rham $\Delta$ est \'egal \`a deux
fois le laplacien complexe $\square$. On d\'eduit donc des rappels ci-dessus la proposition suivante.

\begin{prop} \label{k=a+b}
Soit $M$ une vari\'et\'e kaehl\'erienne. 
\begin{enumerate}
\item Sans autre hypoth\`ese, on a pour tout $k$ une d\'ecomposition orthogonale
$${\cal H}_2^k (M) =\bigoplus_{a+b=k} {\cal H}_2^{a,b} (M), \; \overline{{\cal H}_2^{a,b} (M)} = {\cal H}_2^{b,a} (M) .$$
\item Si de plus $M$ est compl\`ete, il y a des isomorphismes canoniques
$$H_2^k (M) \cong \bigoplus_{a+b=k} H_2^{a,b} (M), \; \overline{H_2^{a,b} (M)} = H_2^{b,a} (M) .$$
\end{enumerate}
\end{prop}

\medskip

\markboth{CHAPITRE 14. COHOMOLOGIE $L^2$ R\'EDUITE}{14.2. TH\'EORIE DE HODGE}

\section{Th\'eorie de Hodge des vari\'et\'es kaehl\'eriennes faiblement pseudoconvexes}

Les vari\'et\'es kaehl\'eriennes compl\`etes que nous consid\`ererons seront faiblement pseudoconvexes. 
Rappelons qu'une vari\'et\'e complexe $X$ est dite {\it faiblement pseudoconvexe} \index{faiblement pseudoconvexes} s'il existe une fonction 
d'exhaustion psh $\psi$ de classe $C^{\infty}$ sur $X$ \footnote{Rappelons qu'une fonction $\psi$ est dite psh (pluri-sous-harmonique)
\index{fonction psh} 
si sa forme de Levi est semi-positive et est dite exhaustive si $\psi (z)$ tend vers $+\infty$ quand $z$ tend vers 
l'infini suivant le filtre des compl\'ementaires de parties compactes de $X$.}.

\medskip

En particulier la fonction $\psi (Z) = \log \frac{B}{A} (Z)$ ($Z \in {\cal D}$) qui descend sur $M_V$ fait de 
$M_V$ une vari\'et\'e kaehl\'erienne faiblement convexe. Le Corollaire \ref{vp asymptotique} montre de plus que cette derni\`ere vari\'et\'e a un certain 
nombre de directions strictement pseudoconvexes. 

\medskip

Nous allons dans cette section d\'ecrire un th\'eor\`eme de 
d\'ecomposition de Hodge pour des vari\'et\'es kaehl\'eriennes faiblement pseudoconvexes ayant justement ``suffisamment
de directions strictement pseudoconvexes''. Suivant \cite{AndreottiGrauert}, une vari\'et\'e complexe $X$ sera dite 
{\it absolument $l$-convexe} \index{absolument $l$-convexe} si $X$ poss\`ede une fonction d'exhaustion psh $\psi$ qui est fortement $l$-convexe sur \index{fortement $l$-convexe}
le compl\'ementaire $X\setminus K$ d'une partie compacte, {\it i.e.} telle que la forme de Levi de $\psi$ a au moins 
$n-l+1$ valeurs propres positives en tout point de $X \setminus K$, o\`u $n=\dim_{{\Bbb C}} X$.

\medskip

Remarquons imm\'ediatement que d'apr\`es le Corollaire \ref{vp asymptotique}, la vari\'et\'e $M_V$ est 
en fait absolument $((p-r)(q-1)+1)$-convexe.

\medskip

Mous pouvons maintenant \'enoncer le th\'eor\`eme de d\'ecomposition de Hodge pour les vari\'et\'es absolument 
$l$-convexe. Ce r\'esultat est d\^u \`a Ohsawa; dans \cite{Demailly} Demailly en donne une d\'emonstration
simplifi\'ee.

\begin{thm} \label{Ohsawa}
 \index{Th\'eor\`eme d'Ohsawa}
Soit $X$ une vari\'et\'e kaehl\'erienne et $n=\dim_{{\Bbb C}} X$. On suppose que $X$ est absolument $l$-convexe.
Alors, en des degr\'es convenables, il y a d\'ecomposition et sym\'etrie de Hodge~:
$$H^k (X) \cong \bigoplus_{a+b=k} H^{a,b} (X) , \; \overline{H^{a,b} (X)} \cong H^{b,a} (X) , \; k\geq n+l , $$
$$H_c^k (X) \cong \bigoplus_{a+b=k} H_c^{a,b} (X) , \; \overline{H^{a,b}_c (X)} \cong H^{b,a}_c (X) , \; k\leq n-l ,$$
tous ces groupes \'etant de dimension finie. ($H_c^k (X)$ et $H^{a,b}_c (X)$ d\'esignent ici les groupes de cohomologie 
\`a support compact). De plus, on a un isomorphisme de Lefschetz (induit par la multiplication par une puissance convenable 
de la forme de Kaehler)
$$H_c^{a,b} (X) \stackrel{\sim}{\rightarrow} H^{n-b, n-a} (X), \; a+b \leq n-l .$$
\end{thm}

Gr\^ace au Th\'eor\`eme \ref{Ohsawa} (et \`a la Proposition \ref{k=a+b}) la d\'emonstration du Th\'eor\`eme \ref{cohom l2}
se r\'eduit \`a l'\'etude des groupes $H_c^{a,b} (M_V)$ et $H_2^{a,b} (M_V )$. C'est ce que permet un th\'eor\`eme 
d'Ohsawa et Takegoshi que nous d\'ecrivons dans la section suivante. Cette r\'ef\'erence nous a \'et\'e fourni par 
Nader Yeganefar, nous l'en remercions.

\medskip

\markboth{CHAPITRE 14. COHOMOLOGIE $L^2$ R\'EDUITE}{14.3. UN TH\'EOR\`EME D'OHSAWA ET TAKEGOSHI}

\section{Un th\'eor\`eme d'Ohsawa et Takegoshi}

Dans  \cite[Theorem 4.1]{OhsawaTakegoshi}, Ohsawa et Takegoshi d\'emontrent le th\'eor\`eme suivant.

\begin{thm} \label{OT}
 \index{Th\'eor\`eme d'Ohaswa et Takegoshi}
Soit $X$ une vari\'et\'e kaehl\'erienne compl\`ete, $n=\dim_{{\Bbb C}} X$ et $l$ un entier naturel. Supposons qu'il existe
une fonction d'exhaustion $C^{\infty}$ $\psi :X \rightarrow {\Bbb R}$ v\'erifiant les 
conditions suivantes.
\begin{enumerate}
\item Soient $\gamma_1 \geq \ldots \geq \gamma_n$ les valeurs propres de la forme de Levi de $\psi$. Alors,
$$\lim_{c \rightarrow +\infty} \sup_{X\setminus X_c } \gamma_1 =1 \; \mbox{ and } \; \lim_{c\rightarrow +\infty} \inf_{X\setminus X_c} \gamma_{n-l+1} =1 .$$
\item La limite lorsque $c$ tend vers l'infini, de l'infimum de la plus petite valeur propre de 
$\partial \overline{\partial} \psi -8 \partial \psi \overline{\partial} \psi$ sur $X \setminus X_c$ est sup\'erieure ou 
\'egale \`a $-\frac{1}{100n}$.
\end{enumerate}
Dans les condition ci-dessus $X_c = \{ x \in X \; : \; \psi (x) <c \}$. 

Alors, $\dim H_2^{a,b} (X) < \infty $ et 
$$H^{a,b} (X) \cong H_2^{a,b} (X), \; \mbox{ pour } a+b \geq n+l .$$
\end{thm}

\medskip

Remarquons que beaucoup de constantes dans l'\'enonc\'e ci-dessus sont arbitraires. Quitte \`a remplacer la 
fonction $\psi$ par $x \mapsto \lambda^{-2} \psi ( \lambda x)$, on peut par exemple changer le $8$ de la condition 2
par n'importe constante positive. La condition 2 sera donc (en particulier) garantie d\`es que les valeurs propres de la forme de 
Levi de $\psi$ seront toutes positives ou nulles (autrement dit $\psi$ psh) et que la norme de $d\psi$ sera uniform\'ement
major\'ee. 

\medskip

L'hypoth\`ese importante dans l'\'enonc\'e du Th\'eor\`eme \ref{OT} est bien \'evidemment la condition 1.
\`A l'aide de celle-ci, la d\'emonstration du Th\'eor\`eme \ref{OT} repose de mani\`ere essentielle sur une proposition de 
Donnelly et Fefferman que nous d\'ecrivons maintenant. 

\medskip

Soit $M$ une vari\'et\'e complexe hermitienne compl\`ete de dimension complexe $n$. Notons 
$J : T^* M \otimes {\Bbb C} \rightarrow T^* M \otimes {\Bbb C}$ la structure presque complexe de $M$. Comme d'habitude
l'action de $J$ s'\'etend en une action sur les formes diff\'erentielles sur $M$ et $J \phi = i^{a-b} \phi$, pour 
$\phi$ de type $(a,b)$.

Supposons que $F$ soit une fonction r\'eelle de classe $C^2$ sur $M$. Comme dans les autres chapitres nous notons
$dF$ et $\nabla^2 F$ respectivement la diff\'erentielle ext\'erieure de $F$ et le hessien de $F$. Le hessien
d\'efinit une transformation lin\'eaire sym\'etrique $\nabla^2 F : T^* M \rightarrow T^* M$.

Soit $S : T^* M \rightarrow T^* M$ une transformation lin\'eaire sym\'etrique de valeurs propres r\'eelles $\gamma_1 , 
\ldots , \gamma_n$. On dit que $S$ est compatible avec $J$ si, pour chaque vecteur propre $v_i$, de valeur propre
$\gamma_i$, de $S$, il existe un indice $i^*$ tel que $v_{i^*} =J v_i$ soit aussi un vecteur propre, de valeur propre
associ\'ee $\gamma_{i^*}$. Soient $\mu_1, \ldots , \mu_n $ les nombres obtenus en moyennant  $\gamma_i$ et
$\gamma_{i^*}$ pour $i=1, \ldots , 2n$. Autrement dit, $\mu_1 = (\gamma_1 + \gamma_{1^*} )/2$, et si 
$v_{1^*} \neq v_2$, alors $\mu_2 = (\gamma_2 + \gamma_{2^*} )/2$, et ainsi de suite...

Si $\phi$ est une forme diff\'erentielle sur $M$, nous notons $|\phi |$ la norme ponctuelle de $\phi$ et $||\phi ||_{L^2}^2 =
\int_M |\phi |^2$ la norme $L^2$ globale.

\begin{prop} \label{DF}
Soit $\phi$ une forme diff\'erentielle dans $C_0^{\infty} (\Lambda^{a,b} T^* X)$.
Soit $F$ une fonction r\'elle $C^2$ sur le support de $\phi$ et telle que
la transformation lin\'eaire $\nabla F$ soit sym\'etrique et compatible avec $J$ avec pour valeurs 
propres r\'eelles $\gamma_1 , \ldots , \gamma_{2n}$. Si $|dF| \leq 1$, alors~:
$$[||d\phi ||_{L^2} + ||\delta \phi ||_{L^2} ] ||\phi ||_{L^2} \geq \int_M \left[ \sum \mu_i - (a+b) \max_i (\mu_i ) \right] |\phi |^2 ,$$
o\`u $\mu_1, \mu_2 , \ldots , \mu_n$ sont les nombres obtenus en moyennant $\gamma_i$ et $\gamma_{i^*}$, pour 
$i =1, \ldots ,2n$. 
\end{prop}

Remarquons que Donnelly et Fefferman \'enoncent une version \`a bord de cette proposition pour une forme diff\'erentielle $\phi$
g\'en\'erale de degr\'e $k$ (et non de bidegr\'e $(a,b)$). Nous ne savons pas si cette g\'en\'eralisation est vraie, si c'\'etait le cas un 
argument classique de suite exacte \`a la Cheeger permettrait de donner une d\'emonstration directe du Th\'eor\`eme \ref{cohom l2}. 
N\'eanmoins comme nous l'ont signal\'e Gilles Carron et Nader Yeganefar la d\'emonstration sugg\'er\'ee par Donnelly et Fefferman utilise
que la composante de type $(a,b)$ d'une forme diff\'erentielle de degr\'e $k=a+b$ v\'erifiant la condition absolue au bord
v\'erifie elle aussi cette condition, ce qui est faux en g\'en\'eral.

\medskip

La vertu essentielle de la Proposition \ref{DF} est de contr\^oler le spectre essentiel de $M_V$.

Soit $\psi$ la fonction r\'elle sur ${\cal D}$ d\'efinie par 
$$\psi (Z) = \log \left( \frac{B}{A} \right) $$
(nous pourrions \'egalement consid\'erer $\psi (Z) =d(Z, {\cal D}_V )$ si $\min \{ r,q \} =1$).
D'apr\`es le Lemme 12.4.2 et le Corollaire \ref{vp asymptotique}, le hessien 
$\nabla^2 \psi$ de $\psi$ est compatible avec $J$ et les valeurs propres moyenn\'ees 
$\mu_1 (Z) , \ldots , \mu_n (Z)$ v\'erifient qu'il existe une constante $c_0 >0$ telle que 
pour tout $Z$ tel que $\psi (Z) > c_0$, 
\begin{eqnarray} \label{1/10}
\sum_i \mu_i (Z) - (a+b) \max_i (\mu_i (Z)) \geq 1/10 ,
\end{eqnarray}
pour $a+b < p+qr -r$.

Remarquons que la fonction $\psi $ est $G_V$-invariante et descend donc en une fonction sur $M_V = \Gamma_V 
\backslash {\cal D}$. Dans la suite, $X= M_V $ et 
$X_c = \{ Z \in X \; : \; \psi (Z )< c\} $ pour tout r\'eel $c>0$. 

\begin{lem} \label{spectre isole de 0}
Pour tout bidegr\'e $(a,b)$ tel que $a+b <p+qr-r$, le spectre essentiel du laplacien sur les 
formes de bidegr\'e $(a,b)$ est isol\'e de $0$.
\end{lem}
{\it D\'emonstration.} Il est bien connu que le spectre essentiel du laplacien ne 
d\'epend que de la g\'eom\'etrie \`a l'infini. 
Or, la Proposition \ref{DF} et l'in\'egalit\'e (\ref{1/10}) impliquent que si $\omega$ est une forme diff\'erentielle de 
bidegr\'e $(a,b)$, $a+b < p+qr-r$, \`a support dans $X_{c_0}$ alors~:
\begin{eqnarray*}
[||d\omega ||_2 + ||\delta \omega ||_2 ] ||\omega ||_2 & \geq &\int_X \left[ \sum \mu_i -k \max_i (\mu_i ) \right] |\omega |^2 \\
& \geq & \frac{1}{10} ||\omega ||_2^2 .
\end{eqnarray*}
Le spectre du laplacien \`a l'infini (et donc le spectre essentiel) est donc isol\'e de $0$.
Le Lemme \ref{spectre isole de 0} est d\'emontr\'e.  

\bigskip

Le Lemme \ref{spectre isole de 0} est l'\'etape essentielle dans la d\'emonstration du Th\'eor\`eme \ref{OT}.
Il est naturel de se demander si dans l'\'enonc\'e du Th\'eor\`eme \ref{cohom l2} ou du Lemme \ref{spectre isole de 0}
le nombre $p+qr-r$ est optimal. En g\'en\'eral on ne sait pas r\'epondre \`a cette question. Peut-\^etre la formule de Plancherel
pour les espaces sym\'etriques pseudoriemannien $G/H$ peut-elle apporter une r\'eponse. Remarquons
n\'eanmoins le lemme suivant.

\begin{lem} \label{q=1}
Si $q=1$, le spectre essentiel du laplacien sur les formes de bidegr\'e $(a,b)$ avec $a+b=p$ est encore isol\'e de $0$.
\end{lem}
{\it D\'emonstration.} C'est \'evident puisque, le laplacien complexe commutant 
aux op\'erateurs $\overline{\partial}$ et $\overline{\partial}^*$, son spectre sur les formes de bidegr\'e $(a,b)$ est 
contenu dans la r\'eunion du spectre du laplacien sur les formes de bidegr\'e $(a,b-1)$, du spectre du laplacien sur 
les formes de bidegr\'e $(a,b+1)$ et (\'eventuellement) de la valeur propre $0$. Mais par dualit\'e de Hodge le 
spectre du laplacien sur les formes de bidegr\'e $(a,b+1)$ co\"{\i}ncide avec le spectre du laplacien sur les formes
de bidegr\'e $(b,a-1)$. Le Lemme \ref{q=1} d\'ecoule donc du Lemme \ref{spectre isole de 0}.

\medskip

\markboth{CHAPITRE 14. COHOMOLOGIE $L^2$ R\'EDUITE}{14.4. D\'EMONSTRATION DU TH\'EOR\`EME \ref{cohoml2}}

\section{D\'emonstration du Th\'eor\`eme \ref{cohom l2}}

Commen\c{c}ons par appliquer le Th\'eor\`eme \ref{OT} \`a la vari\'et\'e $M_V$ en utilisant la fonction d'exhaustion 
$\psi (Z) = \log \frac{B}{A} (Z) $. D'apr\`es le Corollaire \ref{vp asymptotique} la condition 1 du Th\'eor\`eme \ref{OT}
est v\'erifi\'ee pour $l=(p-r)(q-1)+1$, par ailleurs les valeurs propres de la forme de Levi de $\psi$ sont toutes 
positives (ou nulles) et la norme de la diff\'erentielle de $\psi$ est uniform\'ement born\'ee, d'apr\`es la 
remarque suivant le Th\'eor\`eme \ref{OT} celui-ci s'applique \`a $M_V$. On en d\'eduit que pour tout bidegr\'e 
$(a,b)$ tel que $a+b \geq pq + (p-r)(q-1)+1$, le groupe $H_2^{a,b} (M_V )$ est de dimension finie et 
$$H_2^{a,b} (M_V ) \cong H^{a,b} (M_V ) .$$

On d\'eduit alors du Th\'eor\`eme \ref{Ohsawa} et de la Proposition \ref{k=a+b} que 
$$H_2^k (M_V ) \cong H^k (M_V ) $$
pour tout $k\geq pq+(p-r)(q-1)+1$. Puis par dualit\'e, il d\'ecoule que l'application naturelle
$$H_c^k (M_V ) \rightarrow H_2^k (M_V )$$
est un isomorphisme pour tout entier $k \leq p+qr-r-1$.

On sait par ailleurs que $H_2^k (M_V )$ est (en tout degr\'e $k$) isomorphe \`a ${\cal H}_2^k (M_V )$ et puisque 
$M_V$ est hom\'eomorphe au produit $C_V \times {\Bbb R}^{qr}$, la formule de K\"unneth implique que 
$$H^{k-2qr} (C_V ) \cong H_c^k (M_V ) ,$$
pour tout degr\'e $k$.  La premi\`ere partie du Th\'eor\`eme \ref{cohom l2} est donc d\'emontr\'ee. 

\medskip

Supposons maintenant que $q=1$. On a toujours l'isomorphisme $H_2^p (M_V ) \cong {\cal H}_2^p (M_V )$. Le groupe
${\cal H}_2^p (M_V )$ ($p$ est ici la dimension r\'eelle moiti\'e de $M_V$) ne d\'epend que de la structure 
conforme de $M_V$. Il est donc facile de v\'erifier que l'espace ${\cal H}_2^p (M_V )$ est de dimension infinie et que 
l'application naturelle $H_c^p (M_V ) \rightarrow H_2^p (M_V )$ est injective.

\newpage

\thispagestyle{empty}

\newpage

\markboth{CHAPITRE 15. D\'EMONSTRATIONS DES TH\'EOR\`EMES 4, 5 ET 8}{15.1. D\'EMONSTRATION DU TH\'EOR\`EME 4}

\chapter{D\'emonstrations des Th\'eor\`emes 4, 5 et 8}

Comme l'indique le titre, le but de ce chapitre est la d\'emonstration des 
Th\'eor\`emes 4, 5 et 8. N\'eanmoins nous ne parlerons pas de l'espace hyperbolique r\'eel {\it i.e.} du cas 
o\`u ${\Bbb G}$ est un groupe alg\'ebrique du type $SO(n,1)$. Dans ce cas le Th\'eor\`eme 4 est en effet plus facile 
\`a d\'emontrer et est de toute fa\c{c}on compl\`etement trait\'e dans \cite{MathZ}.

\section{D\'emonstration du Th\'eor\`eme 4}

Nous conservons les notations des chapitres pr\'ec\'edents. Le but de cette section est la d\'emonstration du 
th\'eor\`eme suivant qui \`a l'aide du Lemme 12.5.3 implique imm\'ediatement le Th\'eor\`eme 4.

\begin{thm}    \label{sur l'homologie}
Soit $M = \Gamma \backslash {\cal D}$ une vari\'et\'e compacte localement sym\'etrique model\'ee sur l'espace sym\'etrique 
${\cal D}$ associ\'e au groupe semi-simple $G= SU(p,q)$. 
Supposons que l'espace ${\cal D}$ contienne un sous-espace ${\cal D}_V$ tel que la vari\'et\'e 
$C_V =\Gamma_V \backslash {\cal D}_V$, avec $\Gamma_V = \Gamma \cap G_V$, soit compacte. 
Soit $k$ un entier $> 2pq-p-qr+r$, ou \'egal \`a $p$ si $q=1$. Supposons l'hypoth\`ese suivante v\'erifi\'ee~:
\begin{description}	
\item[(H)] $M$ admet une tour d'effeuillage autour de $C_V$ dont la premi\`ere valeur propre non nulle du laplacien sur 	
les $(2pq-k)$-formes ferm\'ees est uniform\'ement minor\'ee.
\end{description}
Alors, il existe un rev\^etement fini $\widehat{M}$ de $M$ tel que 
\begin{enumerate}
\item l'immersion de $C_V$ dans $M$ se rel\`eve en un plongement de $C_V$ dans $\widehat{M}$,
\item l'application induite~:
$$H_k (C_V ) \rightarrow H_k (\widehat{M} )$$
est injective.
\end{enumerate}

De plus pour tout entier $N$ et tout cycle $c$ dans $H_k (C_V)$, il existe un rev\^etement fini $M_N$ de $M$ contenant 
$N$ pr\'eimages de $i(c)$ lin\'eairement ind\'ependantes dans $H_k (M_N )$.
\end{thm}

\bigskip

Notons $\{ M_m \}$ la tour d'effeuillage fournie par l'hypoth\`ese (H).

Remarquons imm\'ediatement que le point 1. du Th\'eor\`eme \ref{sur l'homologie} se d\'emontre de la m\^eme 
mani\`ere que le Th\'eor\`eme 1 de \cite{EnseignMath}.

Fixons maintenant un cycle $c$ non nul dans $H_k (C_V )$. Soit $\mu$ la forme harmonique sur $C_V$ telle que $*_0 \mu $ 
soit duale \`a $c$ dans $C_V$. Dans la suite nous conservons les notations des chapitres pr\'ec\'edents.

\medskip

\noindent
{\bf Fait 1.} Soit $s \in {\Bbb C}$, Re$(s)>>1$. Alors, la suite $P_0^m \omega_s^m$ converge uniform\'ement 
sur tout compact de ${\cal D}$ vers $P_0^{\infty} \Omega (s)$.

\medskip

En effet, soit $\lambda$ un r\'eel strictement positif tel que la premi\`ere valeur propre non nulle du laplacien sur les 
$(2pq-k)$-formes ferm\'ees de $M_m$ (resp. $M_{\infty}$) soit strictement sup\'erieure \`a $\lambda$ (un tel $\lambda$ 
existe d'apr\`es l'hypoth\`ese (H)). On introduit la fonction $h_{\lambda} \in C_0 ([0, +\infty [)$ qui vaut $1$ sur 
l'intervalle $[0, \frac{\lambda }{2} ]$, $0$ sur l'intervalle $[\lambda , +\infty [$ et qui d\'ecroit lin\'eairement sur 
$[\frac{\lambda}{2} ,\lambda ]$. Puisque~:
\begin{enumerate}
\item la seule valeur propre du laplacien sur les $(2pq-k)$-formes ferm\'ees de $M_m$ (resp. $M_{\infty}$) 
est strictement inf\'erieure \`a $\lambda$ est $0$,
\item l'espace des formes ferm\'ees est ferm\'e, et
\item les formes $\omega_s^m$ et $\Omega (s)$ sont ferm\'ees,
\end{enumerate}
on a~:
$$h(\Delta_m ) \omega_s^m = P_0^m \omega_s^m \mbox{  et  } h(\Delta_{\infty} ) \Omega (s) = P_0^{\infty} \Omega (s) .$$

Or, d'apr\`es la Proposition 13.5.5 et pour $s \in {\Bbb C}$, Re$(s)>>1$, la suite $h_{\lambda} \omega_s^m$ 
converge vers $h_{\lambda }\Omega (s)$ uniform\'ement sur les compacts.
Ce qui conclut la d\'emonstration du Fait 1.

\bigskip

\noindent
{\bf Fait 2.}  Soit $s \in {\Bbb C}$, Re$(s)>>1$. Alors, la forme harmonique 
$P_0^{\infty} \Omega (s)$ est non nulle (et ne d\'epend pas de $s$). 

\medskip

En effet, la forme $\Omega (s)$ repr\'esente une classe de cohomologie dans $H^{2pq-k}_c (M_V )$ 
(avec les notations du chapitre pr\'ec\'edent). Plus pr\'ecisement, d'apr\`es le Th\'eor\`eme 13.2.3, cette 
classe est ind\'ependante de $s$ et correspond, via la dualit\'e de Poincar\'e 
$H^{2pq-k}_c (M_V ) \cong H_k (V ) \cong H_k (C_V)$, au cycle $c$ dans $C_V$. Elle est donc non nulle. 
D'apr\`es le Th\'eor\`eme \ref{cohom l2} et puisque $k >2pq-p-qr+r$, on a les isomorphismes naturels
$H_c^{2pq-k} (M_V ) \cong H_2^{2pq-k} (M_V ) \cong {\cal H}_2^{2pq-k} (M_V )$,
le projet\'e harmonique $L^2$, $P_0^{\infty} \Omega (s)$, dans ${\cal H}_2^{2pq-k}$ est donc non nul (et ne d\'epend pas
de $s$). Remarquons que lorsque $q=1$ et $k=p$, l'application naturelle $H_c^p (M_V ) \rightarrow H_2^p (M_V) \cong 
{\cal H}_2^p (M_V )$ est encore injective et donc, l\`a encore, $P_0^{\infty} \Omega (s) \neq 0$.

\bigskip

On peut maintenant d\'emontrer la premi\`ere partie du Th\'eor\`eme \ref{sur l'homologie}. 

\medskip

Par la dualit\'e de Poincar\'e et le Th\'eor\`eme de Hodge-de Rham, l'application $H_k (C_V ) \rightarrow H_k (M_m )$ correspond 
\`a l'application $\mu \mapsto P_0^m \omega_s^m$ allant des $k$-formes harmoniques sur $C_V$ vers les 
$(2pq-k)$-formes harmoniques sur $M_m$. Mais, d'apr\`es le fait 1, cette derni\`ere converge simplement vers 
l'application, injective d'apr\`es le fait 2, des $k$-formes harmoniques sur $C_V$ vers ${\cal H}^{2pq-k}_2 (M_V )$ qui \`a 
$\mu$ associe $P_0^{\infty} \Omega (s)$. Or $H_k (C_V )$ est de dimension finie donc la convergence est uniforme et pour 
$m$ grand, l'application $H_k (C_V ) \rightarrow H_k (M_m )$ est injective. Comme $M_m$ est un rev\^etement fini de 
$M$, la premi\`ere partie du Th\'eor\`eme \ref{sur l'homologie} est d\'emontr\'ee.

\bigskip

Pour montrer la deuxi\`eme partie du Th\'eor\`eme \ref{sur l'homologie}, nous allons d'abord d\'eduire des faits pr\'ec\'edents un autre 
corollaire.

\begin{prop} \label{translate}   
On se place sous l'hypoth\`ese (H).  Soit $s \in {\Bbb C}$, Re$(s)>>1$. Supposons $k>2pq-p-qr+r$ (ou $k=p$ si $q=1$). Si 
$P_0^0 \omega_s^0 \neq 0$, il existe un entier $m \geq 0$ et un \'el\'ement $\gamma \in \Gamma$ tels que les formes 
harmoniques $P_0^m \omega_s^m$ et $\gamma^* P_0^m \omega_s^m$ soient lin\'eairement ind\'ependantes.
\end{prop}
{\it D\'emonstration.} Nous montrons d'abord par l'absurde qu'il existe un entier $m \geq 0$ tel que la forme 
$P_0^m \omega_s^m$ ne soit pas invariante sous l'action de $\Gamma$. Soit $m$ un entier $\geq 0$. Supposons que 
la forme $P_0^m \omega_{s}^m$  soit invariante sous l'action de $\Gamma$. Alors,
$$P_0^0 \omega_{s}^{0} = [\Gamma  :\Gamma_m ] P_0^m \omega_{s}^m .$$
Or $[\Gamma :\Gamma_m ]$ tend vers l'infini avec $m$, donc $P_0^m \omega_s^m$ tend vers $0$ avec $m$ ce qui 
contredit les faits 1 et 2. Il existe donc un entier $m \geq 0$ et un \'el\'ement $\gamma \in \Gamma$ tels que les formes 
$P_0^m \omega_s^m$ et $\gamma^* P_0^m \omega_s^m$ soient distinctes. Puisque 
$$\sum_{g \in \Gamma_m \backslash \Gamma } g^* P_0^m \omega_s^m = \sum_{g\in \Gamma_m \backslash \Gamma} g^* (\gamma^* P_0^m \omega_s^m )=P_0^0 \omega_s^{0} \neq 0, $$
les formes $P_0^m \omega_s^m$ et $\gamma^* P_0^m \omega_s^m$ sont en fait n\'ecessairement lin\'eairement 
ind\'ependantes. Ce qui ach\`eve la d\'emonstration de la Proposition \ref{translate}.  

\bigskip

\`A l'aide de la Proposition \ref{translate}, nous pouvons maintenant conclure la d\'emonstration du Th\'eor\`eme \ref{sur l'homologie}.

\medskip

Soit $k >2pq-p-qr+r$ (ou $k=p$ si $q=1$). Soit $\mu$ une $k$-forme harmonique sur $C_V$. Nous allons montrer par r\'ecurrence sur $N\geq 1$ qu'il 
existe un rev\^etement fini $M_N$ de $M$ et $N$ formes harmoniques de degr\'e $2pq-k$ sur $M_N$ lin\'eairement 
ind\'ependantes. Le Th\'eor\`eme \ref{sur l'homologie} en d\'ecoule imm\'ediatement. 

Supposons qu'il existe un tel rev\^etement $M_N$ pour un certain $N\geq 1$.
Notons $\omega_1 , \ldots , \omega_N$ les $N$ formes harmoniques ind\'ependantes.
On peut supposer que la forme $\omega_1$ est la forme harmonique associ\'ee \`a $\mu$ (qui est une forme 
harmonique sur une pr\'eimage de $C_V$ dans $M_N$). La Proposition \ref{translate} implique qu'il existe un 
rev\^etement fini $M_{N+1}$ de $M_N$, une forme harmonique $\hat{\omega}_1$ sur $M_{N+1}$ et un \'el\'ement 
$\gamma \in \pi_1 M_N$ tels que les formes harmoniques $\hat{\omega}_1$ et $\gamma^* \hat{\omega}_1$ soient 
lin\'eairement ind\'ependantes. Supposons qu'il existe $N+1$ r\'eels $\alpha_0 , \alpha_1,\ldots ,\alpha_N$ tels que 
$$\alpha_0 \hat{\omega}_1 + \alpha_1 \gamma^* \hat{\omega}_1 +\alpha_2 \omega_2 +\cdots +\alpha_N \omega_N =0.$$

Alors en moyennant par $\pi_1 M_{N+1} \backslash \pi_1 M_{N}$, on obtient~:
$$(\alpha_1 +\alpha_2) \omega_1 + [\pi_1 M_N  : \pi_1 M_{N+1} ] \{ \alpha_2 \omega_2 +\cdots +\alpha_N \omega_N \} =0.$$
L'hypoth\`ese de r\'ecurrence implique donc~:
$$\alpha_0 +\alpha_1 =\alpha_2 =\ldots =\alpha_N =0.$$
Et puisque les formes harmoniques $\hat{\omega}_1$ et $\gamma^* \hat{\omega}_1$ sont lin\'eairement 
ind\'ependantes, on obtient finalement que~:
$$\alpha_0 =\alpha_1 =0.$$
Enfin puisque le rev\^etement de $M_{N+1}$ sur $M$ est fini, on peut le supposer galoisien, la forme 
$\gamma^* \hat{\omega}_1$ est alors d\'efinie sur $M_{N+1}$. Ce qui ach\`eve la r\'ecurrence et la d\'emonstration du 
Th\'eor\`eme \ref{sur l'homologie}.

\bigskip

Remarquons que d'apr\`es la section 13.3 et la d\'emonstration du Th\'eor\`eme \ref{sur l'homologie}, lorsque 
$q=1$, on peut remplacer l'hypoth\`ese (H) par l'hypoth\`ese suivante plus faible~:
\begin{description}
\item[(H')] $M$ admet une tour d'effeuillage autour de $C_V$ dont la premi\`ere valeur propre non nulle du laplacien sur 	
les $(2pq-k)$-formes $\partial$ et $\overline{\partial}$ ferm\'ees est uniform\'ement minor\'ee.
\end{description}
C'est cette hypoth\`ese que nous nous attacherons \`a v\'erifier dans la d\'emonstration du Th\'eor\`eme 5. Mais avant 
de passer \`a la d\'emonstration de celui-ci, on fait un petit aparte concernant les vari\'et\'es arithm\'etiques.

\markboth{CHAPITRE 15. D\'EMONSTRATIONS DES TH\'EOR\`EMES 4, 5 ET 8}{15.2. VARI\'ET\'ES ARITHM\'ETIQUES}

\section{Vari\'et\'es arithm\'etiques}

En g\'en\'eral, la classe des vari\'et\'es arithm\'etiques
est plus large que celle des vari\'et\'es de congruences auxquelles les r\'esultats de la premi\`ere partie s'appliquent.

\medskip

Un sous-groupe discret $\Gamma$ de $SU(p,q)$ est dit {\it arithm\'etique} \index{groupe arithm\'etique} s'il existe un groupe ${\Bbb Q}$-alg\'ebrique 
$G \subset GL_N ({\Bbb R} )$ pour un certain entier naturel $N$ et un homomorphisme continue et surjectif 
$\rho : G \rightarrow SU(p,q)$ tels que~:
\begin{enumerate}    
\item le noyau de $\rho$ soit un sous-groupe {\bf compact} de $G$;    
\item l'image par $\rho$ de $G \cap GL_N ({\Bbb Z} )$ soit {\bf commensurable} avec $\Gamma$, 
{\it i.e.} l'intersection $\Gamma \cap \rho (G \cap GL_N ({\Bbb Z} ))$ est d'indice fini dans $\Gamma$ et dans 
$\rho (G \cap GL_N ({\Bbb Z} )$.
\end{enumerate}

\bigskip

\'Etant donn\'e un sous-groupe arithm\'etique de $\Gamma$ de $SU(p,q)$, on appelle {\it sous-groupe principal de congruence} \index{sous-groupe principal de congruence} de 
$\Gamma$ tout sous-groupe de la forme~:
$$\Gamma (m) = \Gamma \cap \rho (G \cap {\rm ker}( GL_N ({\Bbb Z}) \rightarrow GL_N ({\Bbb Z} /m {\Bbb Z} ) )) , $$
pour un certain entier naturel $m$. Plus g\'en\'eralement, on appelle {\it sous-groupe de congruence} \index{sous-groupe de congruence} de $\Gamma$ \footnote{La terminologie est dangereuse, voir la remarque qui suit.}  
tout sous-groupe de $\Gamma$ contenant un sous-groupe principal de congruence.
Remarquons qu'un sous-groupe de congruence d'un groupe arithm\'etique donn\'e $\Gamma$ n'est en g\'en\'eral pas un sous-groupe de 
congruence de $SU(p,q)$. Ceci n'est vrai que lorsque $\Gamma$ est d\'ej\`a un sous-groupe de congruence de $SU(p,q)$.
 
Enfin, on appelle {\it vari\'et\'e arithm\'etique} \index{vari\'et\'e arithm\'etique} tout quotient de l'espace ${\cal D}$ par un sous-groupe arithm\'etique sans torsion
$\Gamma$ de $SU(p,q)$. Remarquons qu'une telle vari\'et\'e est de volume fini. Les rev\^etements fini d'une vari\'et\'e 
aritm\'etique obtenues \`a l'aide de sous-groupes de congruence (au sens ci-dessus) seront appel\'es {\it rev\^etements de congruence}.\index{rev\^etement de
congruence}

\medskip

Il existe des vari\'et\'es hyperboliques complexes de volume fini non arithm\'etiques pour $n=1,2$ et $3$, construites par Mostow 
(cf. \cite{Mostow}). La question de l'existence de vari\'et\'es hyperboliques complexes non arithm\'etiques pour $n\geq 4$ reste 
ouverte. Pour $q>1$, toute vari\'et\'e localement sym\'etrique de volume fini model\'ee sur ${\cal D}$ est
arithm\'etique d'apr\`es le c\'el\`ebre Th\'eor\`eme d'arithm\'eticit\'e de Margulis. En revanche, la 
question de savoir si tous les sous-groupes arithm\'etiques sont de congruence est encore ouverte.(La question est, maintenant, 
essentiellement r\'eduite au cas des groupes ``exotiques'' du Chapitre 10. Pour une discussion d\'etaill\'ee voir Prasad \cite{Prasad}.) 

\medskip

Rappelons maintenant qu'il existe effectivement des vari\'et\'es arithm\'etiques quotients  de l'espace
${\cal D}$ associ\'e au groupe $SU(p,q)$. Soit $K$ un corps de nombre totalement r\'eel de degr\'e $m$ sur 
${\Bbb Q}$, soit ${\cal O}$ son anneau des entiers et soit $\Sigma = \{\sigma_1 , \ldots , \sigma_m \}$ l'ensemble des 
plongements de $K$ dans ${\Bbb R}$. Soit $K'=K(\sqrt{-1} )$. On \'etend chaque $\sigma \in \Sigma$ en un plongement 
de $K'$ dans ${\Bbb C}$ laissant $\sqrt{-1}$ fixe. Soit 
$$h(z_1 , \ldots , z_n , z_{n+1} )=a_1 |z_1 |^2 + \ldots + a_p |z_p |^2 -a_{p+1} |z_{p+1} |^2 - \ldots - a_{p+q} |z_{p+q} |^2 $$
une forme hermitienne diagonale avec $a_i \in K$. Supposons que $^{\sigma_1 } h$ a pour signature $(p,q)$ et que 
$^{\sigma_i } h$ est d\'efinie positive pour $i=2,3,\ldots ,m$. Le sous-groupe $\Gamma (h)$ de $SL_{n+1} ({\cal O} 
(\sqrt{-1} ))$ pr\'eservant $h$ s'identifie \`a un sous-groupe de $SU(p,q)$. Si $\Gamma < SU(p,q)$ est un 
sous-groupe commensurable \`a $\Gamma (h)$, alors $\Gamma$ est un groupe arithm\'etique 
(on obtient le morphisme $\rho$ et le groupe $G$ \`a l'aide d'une restriction des scalaires de $K$ \`a ${\Bbb Q}$ 
appliqu\'ee \`a la structure $K$-alg\'ebrique de $SU(p,q)$ donn\'ee par la forme hermitienne $h$).
On appelle un tel groupe~: {\it groupe arithm\'etique standard} \index{groupe arithm\'etique standard} (bien que cette appellation ne le soit pas elle m\^eme). 
Quitte \`a passer \`a un sous-groupe d'indice fini et gr\^ace \`a un lemme classique de Selberg, on peut supposer 
$\Gamma $ sans torsion, il agit alors librement sur ${\cal D}$ et l'espace quotient est une {\it vari\'et\'e arithm\'etique 
standard}. \index{vari\'et\'e arithm\'etique standard} Les vari\'et\'es arithm\'etiques sont toutes compactes sauf celles qui sont standard et telles que
$K={\Bbb Q}$ et $h$ repr\'esente z\'ero sur ${\Bbb Q}$; dans ce cas elles sont de volume fini.

\medskip

Remarquons que notre Th\'eor\`eme 5 est non vide en vertu du th\'eor\`eme suivant.

\begin{thm}[Kazhdan \cite{Kazhdan}, Shimura \cite{Shimura} et Borel-Wallach \cite{BorelWallach}] \label{b1}
Soit $M$ une vari\'et\'e hyperbolique complexe (donc associ\'ee au groupe $SU(n,1)$) 
arithm\'etique standard et soit $N$ un entier naturel quelconque. 
Alors, $M$ admet un rev\^etement (fini) de congruence avec un premier nombre de Betti $\geq N$.
\end{thm}
\index{Th\'eor\`eme de Kazhdan, Shimura et Borel-Wallach}

Autrement dit, si $F$ est une vari\'et\'e hyperbolique complexe arithm\'etique standard compacte de dimension (r\'elle) $d$, 
quitte \`a passer \`a un rev\^etement fini, on peut supposer que $H_{d-1} (F)$ est non nul (et m\^eme de rang 
arbitrairement grand). 

\bigskip

Montrons par ailleurs~:

\begin{prop} \label{arith}
Soit $M= \Gamma \backslash {\cal D}$ une vari\'et\'e arithm\'etique (resp. de congruence). Supposons que ${\cal D}$ contienne un 
sous-espace ${\cal D}_V $ tel que la vari\'et\'e $C_V =\Gamma_V \backslash {\cal D}_V$, avec 
$\Gamma_V = \Gamma \cap G_V$, soit compacte. Alors, la vari\'et\'e $C_V$ est arithm\'etique (resp. de congruence).
\end{prop}
{\it D\'emonstration.} Le groupe $\Gamma$ est arithm\'etique, notons $\rho$ et $G$ comme dans la d\'efinition d'un 
groupe arithm\'etique. Soit $H$ la composante connexe de l'\'el\'ement neutre du groupe~: 
$$\rho^{-1} (G_V ) \subset G \subset GL_{N} ({\Bbb R} ) .$$
Le groupe $H$ est r\'eductif et connexe, il est donc \'egal \`a la composante connexe de l'\'el\'ement neutre de sa 
cl\^oture de Zariski. De plus, le sous-groupe $H \cap GL_N ({\Bbb Z} )$ est un r\'eseau cocompact dans $H$ car~:
\begin{enumerate}    
\item $\rho (H \cap GL_N ({\Bbb Z})) = \rho (G \cap GL_N ({\Bbb     Z})) \cap G_V $,    
\item le groupe $\rho (G \cap GL_N ({\Bbb     Z}))$ est commensurable avec $\Gamma$,    
\item le groupe $\Gamma_V = \Gamma \cap G_V$ est cocompact dans $G_V$, et     
\item le noyau de $\rho$ est compact.
\end{enumerate}
Notons $H'$ l'adh\'erence de Zariski de $H \cap GL_N ({\Bbb Z} )$ dans $GL_N ({\Bbb R})$. 
La composante connexe de l'\'el\'ements neutre $(H')_0$ de $H'$ est contenue dans $H$ et, 
d'apr\`es le Th\'eor\`eme de densit\'e de Borel, elle se surjecte sur $H/K$ o\`u $K$ d\'esigne un compact distingu\'e 
maximal de $H$.

\medskip

\noindent
{\bf Fait.} Le groupe $H'$ est ${\Bbb Q}$-alg\'ebrique.

\medskip

En effet, soit $I_d$ l'ensemble des polyn\^omes de degr\'e $\leq d$ s'annulant sur $H'$. Puisque $H \cap GL_N ({\Bbb Z} )$ 
est Zariski-dense dans $H'$, un polyn\^ome $P$ de degr\'e $\leq d$ est dans $I_d$ si et seulement si $P(h) =0$ pour 
tout $h \in H \cap GL_N ({\Bbb Z} )$. On voit maintenant ces \'equations comme un syst\`eme d'\'equations 
homog\`enes lin\'eaires, une \'equation pour chaque $h \in H \cap GL_N ({\Bbb Z} )$, o\`u les inconnues sont les $r$ 
coefficients de $P$ et les coefficients de l'\'equation lin\'eaire sont dans ${\Bbb Z}$. Puisqu'il y a $r$ inconnus, on peut 
trouver $r$ \'equations telles que les \'equations repr\'esentent le syst\`eme. Mais le noyau d'une matrice $r\times r$ 
\`a coefficients dans ${\Bbb Z} \subset {\Bbb Q}$, vue comme transformation lin\'eaire de ${\Bbb R}^r$, a une base 
constitu\'ee de vecteurs dans ${\Bbb Q}^r$. Donc $I_d$ a une base constitu\'ee d'\'el\'ements \`a coefficients dans 
${\Bbb Q}$ et, puisque c'est vrai pour tout $d$, on obtient bien que $H'$ est d\'efini sur ${\Bbb Q}$. 

\medskip

Finalement, la vari\'et\'e $F$ est commensurable au quotient $K' \backslash H' / (H' \cap GL_N ({\Bbb Z} ))$, c'est donc une vari\'et\'e arithm\'etique. Enfin si $\Gamma = G({\Bbb Q}) \cap K_f$ pour un certain sous-groupe compact-ouvert $K_f \subset G({\Bbb A}_f )$, 
$\Gamma_V = H'({\Bbb Q}) \cap K_f$, la vari\'et\'e $F$ est donc de congruence lorsque $\Gamma$ est un sous-groupe de congruence. 
Ce qui conclut la d\'emonstration du Th\'eor\`eme \ref{arith}.  

\bigskip

Remarquons que la d\'emonstration ci-dessus montre en fait que toute sous-vari\'et\'e localement sym\'etrique d'une
vari\'et\'e localement sym\'etrique arithm\'etique est arithm\'etique. Ce fait est sans doute bien connu mais nous 
n'avons pas trouv\'e de r\'ef\'erences dans la litt\'erature.

\medskip

Rappelons que les vari\'et\'es arithm\'etiques forment en g\'en\'eral une classe plus grande de vari\'et\'es que les 
vari\'et\'es de congruence. Ainsi, il existe des surfaces hyperboliques arithm\'etiques dont la 
premi\`ere valeur propre non nulle du laplacien sur les fonctions est arbitrairement petite. De tels exemples peuvent
par exemple \^etre obtenus en consid\'erant un sous-groupe d'indice fini de $SL(2,{\Bbb Z})$ se surjectant sur ${\Bbb Z}$ et 
l'ensemble des ses sous-groupes de quotient fini. 
N\'eanmoins, une fois une vari\'et\'e arithm\'etique fix\'e on s'attend \`a ce que l'ensemble de ses rev\^etements de 
congruence soit soumis aux m\^emes genres de ph\'enom\`enes que l'ensemble des vari\'et\'es de congruence.
Illustrons par exemple ce principe vague par la cons\'equence suivante de la d\'emonstration de la ``Conjecture $\tau$'' (cf. \S 2.1).

\begin{thm} \label{clozel arith}
Soit $M$ une vari\'et\'e arithm\'etique. Alors, la premi\`ere valeur propre non nulle du laplacien sur les fonctions, 
reste uniform\'ement minor\'ee dans les rev\^etements de congruence de $M$.
\end{thm}
{\it D\'emonstration.} Lorsque $M$ est de congruence c'est le Th\'eor\`eme 2.1.1. Pour passer au cas 
arithm\'etique g\'en\'eral, on utilise un r\'esultat de Brooks \cite{Brooks} qui montre que la propri\'et\'e d'avoir une tour 
de rev\^etements fini $\{ M_m \}$ de $M$ ayant la propri\'et\'e d'avoir une premi\`ere valeur propre non nulle sur le spectre 
du laplacien sur les fonctions uniform\'ement minor\'ee est une propr\'et\'e combinatoire dans le sens qu'elle ne d\'epend
que de la suite des graphes de Schreier des quotients $\pi_1 M_m / \pi_1 M$ pour un syst\`eme de g\'en\'erateur fini et 
fix\'e de $\pi_1 M$. Mais $\pi_1 M$ est commensurable au groupe fondamental d'une vari\'et\'e de congruence, et la 
combinatoire de la tour $\{ M_m \}$ se comporte donc de la m\^eme mani\`ere que la combinatoire d'une tour de vari\'et\'es 
de congruence. D'o\`u l'on d\'eduit le Th\'eor\`eme \ref{clozel arith}. 

\bigskip

Remarquons que si dans le cas des vari\'et\'es de congruence les Conjectures d'Arthur pr\'evoient des valeurs explicites
pour la minoration du spectre, on ne peut en esp\`erer autant dans le cas des vari\'et\'es arithm\'etiques puisque, comme 
on l'a rappel\'e au-dessus, il se peut que la vari\'et\'e arithm\'etique $M$ avec laquelle on d\'emarre ait d\'ej\`a une premi\`ere 
valeur propre non nulle du laplacien sur les fonctions arbitrairement petite. Le Th\'eor\`eme \ref{clozel arith} n'est donc 
pas un corollaire imm\'ediat du Th\'eor\`eme 2.1.1. On ne sait d'ailleurs pas si le Th\'eor\`eme 3 qui traite des $1$-formes
diff\'erentielles peut s'\'etendre (sans minoration explicite) au cas des vari\'et\'es arithm\'etiques. Cette question nous para\^{\i}t 
int\'eressante. Reformulons-la de la mani\`ere suivante~:

\medskip

\noindent
{\bf Question.} La Conjecture A$^-$ implique-t-elle que pour toute vari\'et\'e hyperbolique 
r\'eelle ou complexe arithm\'etique $M$ de dimension (r\'elle $d$) et pour tout entier naturel $i\leq \frac{d}{2} -1$, il 
existe une constante strictement positive $\varepsilon (M,i)$ telle que pour tout rev\^etement de congruence $M'$
de $M$,
$$\lambda^i_1 (M' ) \geq \varepsilon (M,i) \  ?$$
 
\bigskip

Nous allons maintenant pouvoir passer \`a la d\'emonstration du Th\'eor\`eme 5. Celle-ci repose sur la v\'erification 
de l'hypoth\`ese (H') qui n\'ecessite quelques rappels sur les vari\'et\'es kaehl\'eriennes que nous d\'etaillons dans 
la section suivante.

\markboth{CHAPITRE 15. D\'EMONSTRATIONS DES TH\'EOR\`EMES 4, 5 ET 8}{15.3. RAPPELS SUR LES VARI\'ET\'ES KAEHL\'ERIENNES}

\section{Rappels sur les vari\'et\'es kaehl\'eriennes}

Dans cette section $M$ est une vari\'et\'e kaehlerienne compacte. Nous avons d\'ej\`a rappel\'e que 
le laplacien usuel et le laplacien complexe co\"{\i}ncident (\`a une constante multiplicative pr\`es).
Dans cette section nous allons utiliser le laplacien complexe~:
$$\square = -\sqrt{-1} \{ \overline{\partial} \partial \Lambda -\overline{\partial} \Lambda \partial +\partial \Lambda 
\overline{\partial} -\Lambda \partial \overline{\partial} \} ,$$
o\`u $\Lambda$ est l'op\'erateur adjoint \`a l'op\'erateur $L$ de multiplication par la forme de Kaehler. 
Plus pr\'ecisemment sur une vari\'et\'e complexe de dimension $n$, munie d'une m\'etrique hermitienne 
$$ds^2 = \sum_{\alpha , \beta =1}^n g_{\alpha \beta} dz^{\alpha} d\overline{z}^{\beta}, $$
o\`u l'on note $(g^{\alpha , \beta })$ l'inverse $(g_{\alpha , \beta} )^{-1}$ de la matrice $(g_{\alpha , \beta} )$. Soit 
$$\varphi = \frac{1}{p! q!} \sum \varphi_{\alpha_1 \ldots \alpha_p \beta_1 \ldots \beta_q} dz^{\alpha_1} \wedge \ldots 
\wedge d\overline{z}^{\beta_q}.$$
Alors 
\begin{eqnarray} \label{lambda}
\Lambda \varphi = \frac{1}{(p-1)! (q-1)!} \sum \sqrt{-1} g^{\beta \alpha} \varphi_{\alpha \beta \alpha_2 \ldots \alpha_p \beta_2 \ldots \beta_q} 
dz^{\alpha_2} \wedge \ldots \wedge d\overline{z}^{\beta_q} ,
\end{eqnarray}
c'est \`a dire, 
\begin{eqnarray*}
(\Lambda \varphi )_{\alpha_2 \ldots \alpha_p \beta_2 \ldots \beta_q} & = & \sum_{\alpha ,\beta} \sqrt{-1} g^{\beta \alpha} 
\varphi_{\alpha \beta \alpha_2 \ldots \alpha_p \beta_2 \ldots \beta_q} \\
                                                                     & = & (-1)^{p-1}  \sum_{\alpha ,\beta} \sqrt{-1} g^{\beta \alpha} \varphi_{\alpha \alpha_2 \ldots \alpha_p \beta \beta_2 \ldots \beta_q} .
\end{eqnarray*}
En particulier si $\varphi$ est une forme de degr\'e $j$, 
\begin{eqnarray} \label{eqn}
-L\Lambda \varphi +\Lambda L \varphi = (n-j) \varphi .
\end{eqnarray}
Un calcul simple montre alors que le laplacien $\Box$ commute \`a $L$.
Enfin, remarquons que si $\omega$ est une forme $\partial$ et $\overline{\partial}$ ferm\'ee, on a~: 
\begin{eqnarray} \label{laplacien}
\square \omega = \partial \overline{\partial} \sqrt{-1} \Lambda \omega .
\end{eqnarray}

\begin{prop} \label{astuce}    
Soit $M$ une vari\'et\'e k\"ahlerienne compacte de dimension complexe $n$. 
Notons $\lambda$ la premi\`ere valeur propre non nulle du laplacien sur les $1$-formes de $M$. Alors, 
la premi\`ere valeur propre non nulle du laplacien sur les $3$-formes $\partial$ et $\overline{\partial}$ ferm\'ees     
de $M$ est $\geq \lambda$. 
\end{prop}
{\it D\'emonstration.} Soit $\lambda$ la premi\`ere valeur propre non nulle du laplacien sur les $1$-formes de $M$.
Soit $\beta$ une $3$-forme $\mu$-propre, $\partial$ et $\overline{\partial}$ ferm\'ee de $M$. On peut d\'ecomposer 
l'espace des $3$-formes sur $M$ comme somme orthogonale de 
$$ \overline{\{ L\alpha \; : \; \alpha \mbox{ est une } 1 \mbox{-forme sur } M \} }$$
et de son suppl\'ementaire orthogonal. Quitte \`a prendre des suites, on peut donc \'ecrire~:
$$\beta =L\alpha \oplus^{\perp} \gamma .$$
Si $L\alpha= 0$, la forme $\gamma$ est orthogonale \`a l'image de $L$ restreinte aux $1$-formes et donc par dualit\'e, 
on obtient que $\Lambda \beta =0$. Et donc $\mu=0$ puisque $\beta$ est $\partial$ et $\overline{\partial}$ ferm\'ee. 
Enfin, si $\mu \neq 0$, puisque le laplacien pr\'eserve la d\'ecomposition orthogonale ci-dessus, on obtient que $\mu$ 
est une valeur propre du laplacien sur les $1$-formes, avec comme forme propre~: $\alpha$. 
Donc si $\mu \neq 0$, on a n\'ecessairement $\mu \geq \lambda$. Ce qui conclut la d\'emonstration de la Proposition \ref{astuce}.

\bigskip

\`A l'aide du Th\'eor\`eme 3, la Proposition \ref{astuce} va nous permettre de v\'erifier l'hypoth\`ese (H') dans certains cas 
int\'eressants et de d\'emontrer le Th\'eor\`eme 5.

\medskip

Concluons cette section de rappels sur les vari\'et\'es kaehl\'eriennes par un cas particulier \'el\'ementaire de la d\'ecomposition de 
Lefschetz.
\begin{prop} \label{lefschetz}    
Soit $M$ une vari\'et\'e k\"ahlerienne de dimension complexe $n$.     
Alors, l'application $L^k$ pour $1\leq k<n$ est injective des $1$-formes sur les $(2k+1)$-formes et envoie les formes     
harmoniques sur des formes harmoniques.
\end{prop}
{\it D\'emonstration.} La premi\`ere partie est bien connue \cite{GriffithHarris}. De plus, 
le laplacien $\square$ commute \`a $L$ et la Proposition \ref{lefschetz} est d\'emontr\'ee.  

\bigskip

Soit 
$$vb_i (M) = \mbox{sup} \{ b_i (\widehat{M} ) : \widehat{M} \mbox{ est un rev\^etement fini de } M \} ,$$
o\`u $i$ est un entier naturel inf\'erieur \`a la dimension de la vari\'et\'e $M$ et $b_i (M)$ d\'esigne le $i$-\`eme nombre 
de Betti de $M$. On appelle $vb_i (M)$ le {\it $i$-\`eme nombre de Betti virtuel de $M$}.\index{nombre de Betti virtuel}
On d\'eduit de la Proposition \ref{lefschetz} et du Th\'eor\`eme \ref{b1} le corollaire suivant.

\begin{cor}    
Les vari\'et\'es hyperboliques complexes standard ont tous leurs nombres    de Betti virtuels infinis.
\end{cor}

Le Th\'eor\`eme 5 permet de d\'ecrire comment certaines classes non triviales en homologie apparaissent g\'eom\'etriquement. 
Il est maintenant temps de passer \`a la d\'emonstration proprement dite du Th\'eor\`eme 5.

\bigskip

\newpage

\markboth{CHAPITRE 15. D\'EMONSTRATIONS DES TH\'EOR\`EMES 4, 5 ET 8}{15.4. D\'EMONSTRATION DU TH\'EOR\`EME 5}

\section{D\'emonstration du Th\`eor\`eme 5}

Soient $M$ une vari\'et\'e hyperbolique complexe compacte de congruence de dimension $d+2$ 
et $F$ une sous-vari\'et\'e complexe compacte connexe totalement g\'eod\'esique de dimension $d$ et immerg\'ee
dans $M$. On peut supposer, en concervant les notations des chapitres pr\'ec\'edents, que $M= \Gamma \backslash {\cal D}$, $F= \Gamma_V \backslash {\cal D}_V$  
o\`u $\Gamma_V = \Gamma \cap G_V$. 

D'apr\`es le Th\'eor\`eme \ref{sur l'homologie} pour conclure il nous faut construire une tour d'effeuillage autour de $F$ 
v\'erifiant l'hypoth\`ese (H') pour $k=d-1$. La vari\'et\'e $M$ est de congruence, donc d'apr\`es la Proposition \ref{arith}, 
la vari\'et\'e $F$ aussi. Et il existe deux groupes ${\Bbb Q}$-alg\'ebriques $H$ et $G$, avec $H$ un 
${\Bbb Q}$-sous-groupe de $G$, $G ({\Bbb R} ) \subset GL_N ({\Bbb R})$ (pour un certain entier $N$) et 
un morphisme continue \`a noyau compact $\rho : G \rightarrow SU(n,1)$ tels que~:
l'image par $\rho$ de $G\cap GL_N ({\Bbb Z} )$ (resp. $H \cap GL_N ({\Bbb Z} )$) soit $\Gamma$ (resp. $\Gamma_V $).
On introduit par r\'ecurrence, la suite de sous-groupes d'indices finis dans $\Gamma$~:
$$\Gamma_0 = \Gamma$$
et pour tout $m \geq 1$,
$$\Gamma_m = \rho ( p_m^{-1} (p_m ( H \cap GL_N ({\Bbb Z} )))) \cap \Gamma_{m-1} ,$$
o\`u $p_m$ est la projection de $GL_N ({\Bbb Z})$ sur $GL_N ({\Bbb Z} /m{\Bbb Z})$. 
Notons $M_m = \Gamma_m \backslash {\cal D}$. La suite de rev\^etements finis $\{ M_m \}$ de $M$ est une tour 
d'effeuillage autour de $F$ dans $M$; cette tour est, de plus, constitu\'ee de rev\^etements de congruences. 
Mais, d'apr\`es le Th\'eor\`eme 3 et la Proposition \ref{astuce}, la premi\`ere valeur non nulle du laplacien sur les $3$-formes 
$\partial$ et $\overline{\partial}$ ferm\'ees sur $M_m$ est uniform\'ement (par rapport \`a $m$) minor\'ee. 
En particulier, l'hypoth\`ese $(H')$ est v\'erifi\'ee pour cette tour d'effeuillage et pour $k =d-1$.
Le Th\'eor\`eme \ref{sur l'homologie} permet donc de conclure la d\'emonstration du Th\'eor\`eme 5.  

\bigskip

Le Th\'eor\`eme 5 permet la construction de $3$-classes d'homologie non triviales \`a partir du Th\'eor\`eme \ref{b1} de la mani\`ere 
suivante. Si $M$ est une vari\'et\'e hyperbolique complexe arithm\'etique standard  et de congruence d\'efinie par une forme hermitienne 
$h$ \`a $n$ variables (comme au \S 15.2), on obtient une sous-vari\'et\'e (immerg\'ee) $F$ holomorphe totalement g\'eod\'esique 
de codimension (complexe) $1$ (qui est donc de congruence) en restreignant la forme $h$ \`a un $(n-1)$-plan. 
Quitte \`a passer \`a un rev\^etement de congruence de $M$ (et donc de $F$), on peut supposer que $H_{2n-3} (F)$ est 
non trivial, d'apr\`es le Th\'eor\`eme \ref{b1}. Alors, le Th\'eor\`eme 5 permet de relever cette classe d'homologie dans $F$ 
en une classe d'homologie non triviale dans un rev\^etement fini de $M$. 

\bigskip

Bien s\^ur, par d\'efinition des groupes $H_* (Sh^0 G)$, le Corollaire 1 d\'ecoule imm\'ediatement du Th\'eor\`eme
5. Enfin, en rempla\c{c}ant le Th\'eor\`eme 3 par le Corollaire 5.4.2, la d\'emonstration du Th\'eor\`eme 5 s'\'etend au groupe 
$SU(p,q)$. On en d\'eduit le th\'eor\`eme suivant qui implique imm\'ediatement le Th\'eor\`eme 8.

\begin{thm}
Soit $M = \Gamma \backslash {\cal D}$ une vari\'et\'e compacte localement sym\'etrique de congruence model\'ee sur l'espace sym\'etrique 
${\cal D}$ associ\'e au groupe semi-simple $G= SU(p,q)$, avec $p\geq q\geq 2$. 
Supposons que l'espace ${\cal D}$ contienne un sous-espace ${\cal D}_V$ tel que la vari\'et\'e 
$C_V =\Gamma_V \backslash {\cal D}_V$, avec $\Gamma_V = \Gamma \cap G_V$, soit compacte. 
Soit $k$ un entier $> 2pq-p-q+1$. Alors, il existe un rev\^etement fini $\widehat{M}$ de $M$ tel que 
\begin{enumerate}
\item l'immersion de $C_V$ dans $M$ se rel\`eve en un plongement de $C_V$ dans $\widehat{M}$,
\item l'application induite~:
$$H_k (C_V ) \rightarrow H_k (\widehat{M} )$$
est injective.
\end{enumerate}

De plus, pour tout entier $N$ et tout cycle $c$ dans $H_k (C_V)$, il existe un rev\^etement fini $M_N$ de $M$ contenant 
$N$ pr\'eimages de $i(c)$ lin\'eairement ind\'ependantes dans $H_k (M_N )$.
\end{thm}

\newpage

\thispagestyle{empty}

\newpage

\markboth{BIBLIOGRAPHIE}{BIBLIOGRAPHIE} 

\bibliography{bibliographie}

\def\cprime{$'$} \def\cprime{$'$} \def\cprime{$'$} \def\cprime{$'$}
\begin{thebibliography}{100}

\bibitem{AdamsBarbaschVogan}
Jeffrey Adams, Dan Barbasch, and David~A. Vogan, Jr.
\newblock {\em The {L}anglands classification and irreducible characters for
  real reductive groups}, volume 104 of {\em Progress in Mathematics}.
\newblock Birkh\"auser Boston Inc., Boston, MA, 1992.

\bibitem{AndreottiGrauert}
Aldo Andreotti and Hans Grauert.
\newblock Th\'eor\`eme de finitude pour la cohomologie des espaces complexes.
\newblock {\em Bull. Soc. Math. France}, 90:193--259, 1962.

\bibitem{Arthur}
James Arthur.
\newblock Unipotent automorphic representations: conjectures.
\newblock {\em Ast\'erisque}, (171-172):13--71, 1989.
\newblock Orbites unipotentes et repr\'esentations, II.

\bibitem{AC}
James Arthur and Laurent Clozel.
\newblock {\em Simple algebras, base change, and the advanced theory of the
  trace formula}, volume 120 of {\em Annals of Mathematics Studies}.
\newblock Princeton University Press, Princeton, NJ, 1989.

\bibitem{BaldoniSilva}
M.~W. Baldoni~Silva.
\newblock The unitary dual of {${\rm Sp}(n,\,1)$}, {$n\geq 2$}.
\newblock {\em Duke Math. J.}, 48(3):549--583, 1981.

\bibitem{BR}
Laure Barthel and Dinakar Ramakrishnan.
\newblock A nonvanishing result for twists of {$L$}-functions of {${\rm
  GL}(n)$}.
\newblock {\em Duke Math. J.}, 74(3):681--700, 1994.

\bibitem{IRMN}
N.~Bergeron.
\newblock Lefschetz properties for arithmetic real and complex hyperbolic
  manifolds.
\newblock {\em Int. Math. Res. Not.}, (20):1089--1122, 2003.

\bibitem{EnseignMath}
Nicolas Bergeron.
\newblock Premier nombre de {B}etti et spectre du laplacien de certaines
  vari\'et\'es hyperboliques.
\newblock {\em Enseign. Math. (2)}, 46(1-2):109--137, 2000.

\bibitem{MathZ}
Nicolas Bergeron.
\newblock Asymptotique de la norme {$L\sp 2$} d'un cycle g\'eod\'esique dans
  les rev\^etements de congruence d'une vari\'et\'e hyperbolique
  arithm\'etique.
\newblock {\em Math. Z.}, 241(1):101--125, 2002.

\bibitem{BergeronClozel}
Nicolas Bergeron and Laurent Clozel.
\newblock Spectre et homologie des vari\'et\'es hyperboliques complexes de
  congruence.
\newblock {\em C. R. Math. Acad. Sci. Paris}, 334(11):995--998, 2002.

\bibitem{BD}
Pierre Bernat and Jacques Dixmier.
\newblock Sur le dual d'un groupe de {L}ie.
\newblock {\em C. R. Acad. Sci. Paris}, 250:1778--1779, 1960.

\bibitem{BernsteinZelevinski}
I.~N. Bernstein and A.~V. Zelevinsky.
\newblock Induced representations of reductive {$p$}-adic groups. {I}.
\newblock {\em Ann. Sci. \'Ecole Norm. Sup. (4)}, 10(4):441--472, 1977.

\bibitem{Borel}
A.~Borel.
\newblock Automorphic {$L$}-functions.
\newblock In {\em Automorphic forms, representations and $L$-functions (Proc.
  Sympos. Pure Math., Oregon State Univ., Corvallis, Ore., 1977), Part 2},
  Proc. Sympos. Pure Math., XXXIII, pages 27--61. Amer. Math. Soc., Providence,
  R.I., 1979.

\bibitem{BorelHarishChandra}
Armand Borel and Harish-Chandra.
\newblock Arithmetic subgroups of algebraic groups.
\newblock {\em Ann. of Math. (2)}, 75:485--535, 1962.

\bibitem{BorelWallach}
Armand Borel and Nolan~R. Wallach.
\newblock {\em Continuous cohomology, discrete subgroups, and representations
  of reductive groups}, volume~94 of {\em Annals of Mathematics Studies}.
\newblock Princeton University Press, Princeton, N.J., 1980.

\bibitem{BottChern}
Raoul Bott and Shiing~S. Chern.
\newblock Some formulas related to complex transgression.
\newblock In {\em Essays on Topology and Related Topics (M\'emoires d\'edi\'es
  \`a Georges de Rham)}, pages 48--57. Springer, New York, 1970.

\bibitem{Brooks}
Robert Brooks.
\newblock The spectral geometry of a tower of coverings.
\newblock {\em J. Differential Geom.}, 23(1):97--107, 1986.

\bibitem{Bump}
Daniel Bump.
\newblock {\em Automorphic forms and representations}, volume~55 of {\em
  Cambridge Studies in Advanced Mathematics}.
\newblock Cambridge University Press, Cambridge, 1997.

\bibitem{BurgerLiSarnak}
M.~Burger, J.-S. Li, and P.~Sarnak.
\newblock Ramanujan duals and automorphic spectrum.
\newblock {\em Bull. Amer. Math. Soc. (N.S.)}, 26(2):253--257, 1992.

\bibitem{BurgerSarnak}
M.~Burger and P.~Sarnak.
\newblock Ramanujan duals. {II}.
\newblock {\em Invent. Math.}, 106(1):1--11, 1991.

\bibitem{CasselmanMilicic}
William Casselman and Dragan Mili{\v{c}}i{\'c}.
\newblock Asymptotic behavior of matrix coefficients of admissible
  representations.
\newblock {\em Duke Math. J.}, 49(4):869--930, 1982.

\bibitem{CL}
Laurent Clozel.
\newblock Spectral theory of automorphic forms.
\newblock \`a para\^{\i}tre (Conf\'erence du Park City Math. Institute,
  AMS/IAS, 2002).

\bibitem{BC}
Laurent Clozel.
\newblock Changement de base pour les repr\'esentations temp\'er\'ees des
  groupes r\'eductifs r\'eels.
\newblock {\em Ann. Sci. \'Ecole Norm. Sup. (4)}, 15(1):45--115, 1982.

\bibitem{Clozel2}
Laurent Clozel.
\newblock On the cohomology of {K}ottwitz's arithmetic varieties.
\newblock {\em Duke Math. J.}, 72(3):757--795, 1993.

\bibitem{Clozel}
Laurent Clozel.
\newblock D\'emonstration de la conjecture {$\tau$}.
\newblock {\em Invent. Math.}, 151(2):297--328, 2003.

\bibitem{ClozelUllmo}
Laurent Clozel and Emmanuel Ullmo.
\newblock \'{E}quidistribution des points de {H}ecke.
\newblock In {\em Contributions to automorphic forms, geometry, and number
  theory}, pages 193--254. Johns Hopkins Univ. Press, Baltimore, MD, 2004.

\bibitem{LPSS}
J.~Cogdell, J.-S. Li, I.~Piatetski-Shapiro, and P.~Sarnak.
\newblock Poincar\'e series for {${\rm SO}(n,1)$}.
\newblock {\em Acta Math.}, 167(3-4):229--285, 1991.

\bibitem{Cogdell}
James~W. Cogdell and Ilya~I. Piatetski-Shapiro.
\newblock Remarks on {R}ankin-{S}elberg convolutions.
\newblock In {\em Contributions to automorphic forms, geometry, and number
  theory}, pages 255--278. Johns Hopkins Univ. Press, Baltimore, MD, 2004.

\bibitem{D}
P.~Delorme.
\newblock Th\'eor\`eme de {P}aley-{W}iener invariant tordu pour le changement
  de base {${\bf C}/{\bf R}$}.
\newblock {\em Compositio Math.}, 80(2):197--228, 1991.

\bibitem{Delorme}
Patrick Delorme.
\newblock Formules limites et formules asymptotiques pour les multiplicit\'es
  dans {$L\sp 2(G/\Gamma)$}.
\newblock {\em Duke Math. J.}, 53(3):691--731, 1986.

\bibitem{Demailly}
Jean-Pierre Demailly.
\newblock Th\'eorie de {H}odge {$L\sp 2$} et th\'eor\`emes d'annulation.
\newblock In {\em Introduction \`a\ la th\'eorie de Hodge}, volume~3 of {\em
  Panor. Synth\`eses}, pages 3--111. Soc. Math. France, Paris, 1996.

\bibitem{Dixmier}
Jacques Dixmier.
\newblock {\em Les {$C\sp *$}-alg\`ebres et leurs repr\'esentations}.
\newblock Les Grands Classiques Gauthier-Villars. [Gauthier-Villars Great
  Classics]. \'Editions Jacques Gabay, Paris, 1996.
\newblock Reprint of the second (1969) edition.

\bibitem{Donnelly2}
Harold Donnelly.
\newblock On the spectrum of towers.
\newblock {\em Proc. Amer. Math. Soc.}, 87(2):322--329, 1983.

\bibitem{Donnelly}
Harold Donnelly.
\newblock Elliptic operators and covers of {R}iemannian manifolds.
\newblock {\em Math. Z.}, 223(2):303--308, 1996.

\bibitem{DonnellyFefferman}
Harold Donnelly and Charles Fefferman.
\newblock {$L\sp{2}$}-cohomology and index theorem for the {B}ergman metric.
\newblock {\em Ann. of Math. (2)}, 118(3):593--618, 1983.

\bibitem{EGM}
J.~Elstrodt, F.~Grunewald, and J.~Mennicke.
\newblock Kloosterman sums for {C}lifford algebras and a lower bound for the
  positive eigenvalues of the {L}aplacian for congruence subgroups acting on
  hyperbolic spaces.
\newblock {\em Invent. Math.}, 101(3):641--685, 1990.

\bibitem{Epstein}
D.~B.~A. Epstein.
\newblock Complex hyperbolic geometry.
\newblock In {\em Analytical and geometric aspects of hyperbolic space
  (Coventry/Durham, 1984)}, volume 111 of {\em London Math. Soc. Lecture Note
  Ser.}, pages 93--111. Cambridge Univ. Press, Cambridge, 1987.

\bibitem{Flath}
D.~Flath.
\newblock Decomposition of representations into tensor products.
\newblock In {\em Automorphic forms, representations and $L$-functions (Proc.
  Sympos. Pure Math., Oregon State Univ., Corvallis, Ore., 1977), Part 1},
  Proc. Sympos. Pure Math., XXXIII, pages 179--183. Amer. Math. Soc.,
  Providence, R.I., 1979.

\bibitem{GelbartJacquet}
Stephen Gelbart and Herv{\'e} Jacquet.
\newblock A relation between automorphic representations of {${\rm GL}(2)$} and
  {${\rm GL}(3)$}.
\newblock {\em Ann. Sci. \'Ecole Norm. Sup. (4)}, 11(4):471--542, 1978.

\bibitem{GelfandFomin}
I.~M. Gel{\cprime}fand and S.~V. Fomin.
\newblock Unitary representations of {L}ie groups and geodesic flows on
  surfaces of constant negative curvature.
\newblock {\em Doklady Akad. Nauk SSSR (N.S.)}, 76:771--774, 1951.

\bibitem{GGPP}
I.~M. Gel{\cprime}fand, M.~I. Graev, and I.~I. Pyatetskii-Shapiro.
\newblock {\em Representation theory and automorphic functions}, volume~6 of
  {\em Generalized Functions}.
\newblock Academic Press Inc., Boston, MA, 1990.
\newblock Translated from the Russian by K. A. Hirsch, Reprint of the 1969
  edition.

\bibitem{GodementJacquet}
Roger Godement and Herv{\'e} Jacquet.
\newblock {\em Zeta functions of simple algebras}.
\newblock Springer-Verlag, Berlin, 1972.
\newblock Lecture Notes in Mathematics, Vol. 260.

\bibitem{Goldman}
William~M. Goldman.
\newblock {\em Complex hyperbolic geometry}.
\newblock Oxford Mathematical Monographs. The Clarendon Press Oxford University
  Press, New York, 1999.
\newblock Oxford Science Publications.

\bibitem{GriffithHarris}
Phillip Griffiths and Joseph Harris.
\newblock {\em Principles of algebraic geometry}.
\newblock Wiley Classics Library. John Wiley \& Sons Inc., New York, 1994.
\newblock Reprint of the 1978 original.

\bibitem{HarishChandra}
Harish-Chandra.
\newblock Representations of a semisimple {L}ie group on a {B}anach space. {I}.
\newblock {\em Trans. Amer. Math. Soc.}, 75:185--243, 1953.

\bibitem{HarrisLi}
Michael Harris and Jian-Shu Li.
\newblock A {L}efschetz property for subvarieties of {S}himura varieties.
\newblock {\em J. Algebraic Geom.}, 7(1):77--122, 1998.

\bibitem{HarrisTaylor}
Michael Harris and Richard Taylor.
\newblock {\em The geometry and cohomology of some simple {S}himura varieties},
  volume 151 of {\em Annals of Mathematics Studies}.
\newblock Princeton University Press, Princeton, NJ, 2001.
\newblock With an appendix by Vladimir G. Berkovich.

\bibitem{Helgason}
Sigurdur Helgason.
\newblock {\em Differential geometry, {L}ie groups, and symmetric spaces},
  volume~34 of {\em Graduate Studies in Mathematics}.
\newblock American Mathematical Society, Providence, RI, 2001.
\newblock Corrected reprint of the 1978 original.

\bibitem{HoweMoore}
Roger~E. Howe and Calvin~C. Moore.
\newblock Asymptotic properties of unitary representations.
\newblock {\em J. Funct. Anal.}, 32(1):72--96, 1979.

\bibitem{Iwaniec}
Henryk Iwaniec.
\newblock {\em Spectral methods of automorphic forms}, volume~53 of {\em
  Graduate Studies in Mathematics}.
\newblock American Mathematical Society, Providence, RI, second edition, 2002.

\bibitem{JacquetLanglands}
H.~Jacquet and R.~P. Langlands.
\newblock {\em Automorphic forms on {${\rm GL}(2)$}}.
\newblock Springer-Verlag, Berlin, 1970.
\newblock Lecture Notes in Mathematics, Vol. 114.

\bibitem{JacquetPiatetskiiShapiroShalika}
H.~Jacquet, I.~I. Piatetskii-Shapiro, and J.~A. Shalika.
\newblock Rankin-{S}elberg convolutions.
\newblock {\em Amer. J. Math.}, 105(2):367--464, 1983.

\bibitem{JacquetShalika3}
H.~Jacquet and J.~A. Shalika.
\newblock On {E}uler products and the classification of automorphic forms.
  {II}.
\newblock {\em Amer. J. Math.}, 103(4):777--815, 1981.

\bibitem{JacquetShalika4}
H.~Jacquet and J.~A. Shalika.
\newblock On {E}uler products and the classification of automorphic
  representations. {I}.
\newblock {\em Amer. J. Math.}, 103(3):499--558, 1981.

\bibitem{Jacquet}
Herv{\'e} Jacquet.
\newblock Principal {$L$}-functions of the linear group.
\newblock In {\em Automorphic forms, representations and $L$-functions (Proc.
  Sympos. Pure Math., Oregon State Univ., Corvallis, Ore., 1977), Part 2},
  Proc. Sympos. Pure Math., XXXIII, pages 63--86. Amer. Math. Soc., Providence,
  R.I., 1979.

\bibitem{JacquetShalika2}
Herv{\'e} Jacquet and Joseph Shalika.
\newblock The {W}hittaker models of induced representations.
\newblock {\em Pacific J. Math.}, 109(1):107--120, 1983.

\bibitem{JacquetShalika}
Herv{\'e} Jacquet and Joseph Shalika.
\newblock Rankin-{S}elberg convolutions: {A}rchimedean theory.
\newblock In {\em Festschrift in honor of I. I. Piatetski-Shapiro on the
  occasion of his sixtieth birthday, Part I (Ramat Aviv, 1989)}, volume~2 of
  {\em Israel Math. Conf. Proc.}, pages 125--207. Weizmann, Jerusalem, 1990.

\bibitem{Kazhdan}
David Kazhdan.
\newblock Some applications of the {W}eil representation.
\newblock {\em J. Analyse Mat.}, 32:235--248, 1977.

\bibitem{KimSarnak}
Henry~H. Kim.
\newblock Functoriality for the exterior square of {${\rm GL}\sb 4$} and the
  symmetric fourth of {${\rm GL}\sb 2$}.
\newblock {\em J. Amer. Math. Soc.}, 16(1):139--183 (electronic), 2003.
\newblock With appendix 1 by Dinakar Ramakrishnan and appendix 2 by Kim and
  Peter Sarnak.

\bibitem{KnappSpeh}
A.~W. Knapp and B.~Speh.
\newblock Status of classification of irreducible unitary representations.
\newblock In {\em Harmonic analysis (Minneapolis, Minn., 1981)}, volume 908 of
  {\em Lecture Notes in Math.}, pages 1--38. Springer, Berlin, 1982.

\bibitem{KnappZuckerman}
A.~W. Knapp and Gregg~J. Zuckerman.
\newblock Classification of irreducible tempered representations of semisimple
  groups.
\newblock {\em Ann. of Math. (2)}, 116(2):389--455, 1982.

\bibitem{Knapp}
Anthony~W. Knapp.
\newblock {\em Representation theory of semisimple groups}.
\newblock Princeton Landmarks in Mathematics. Princeton University Press,
  Princeton, NJ, 2001.
\newblock An overview based on examples, Reprint of the 1986 original.

\bibitem{KobayashiNomizu}
Shoshichi Kobayashi and Katsumi Nomizu.
\newblock {\em Foundations of differential geometry. {V}ol. {II}}.
\newblock Wiley Classics Library. John Wiley \& Sons Inc., New York, 1996.
\newblock Reprint of the 1969 original, A Wiley-Interscience Publication.

\bibitem{Krajlevic}
Hrvoje Kraljevi{\'c}.
\newblock Representations of the universal convering group of the group {${\rm
  SU}(n,\,1)$}.
\newblock {\em Glasnik Mat. Ser. III}, 8(28):23--72, 1973.

\bibitem{K}
Stephen~S. Kudla.
\newblock The local {L}anglands correspondence: the non-{A}rchimedean case.
\newblock In {\em Motives (Seattle, WA, 1991)}, volume~55 of {\em Proc. Sympos.
  Pure Math.}, pages 365--391. Amer. Math. Soc., Providence, RI, 1994.

\bibitem{Lang}
Serge Lang.
\newblock {\em {${\rm SL}\sb 2({\bf R})$}}, volume 105 of {\em Graduate Texts
  in Mathematics}.
\newblock Springer-Verlag, New York, 1985.
\newblock Reprint of the 1975 edition.

\bibitem{Langlands}
R.~P. Langlands.
\newblock On the classification of irreducible representations of real
  algebraic groups.
\newblock In {\em Representation theory and harmonic analysis on semisimple Lie
  groups}, volume~31 of {\em Math. Surveys Monogr.}, pages 101--170. Amer.
  Math. Soc., Providence, RI, 1989.

\bibitem{MehdiLohoue}
No{\"e}l Lohoue and Salah Mehdi.
\newblock The {N}ovikov-{S}hubin invariants for locally symmetric spaces.
\newblock {\em J. Math. Pures Appl. (9)}, 79(2):111--140, 2000.

\bibitem{Luck}
W.~L{\"u}ck.
\newblock Approximating {$L\sp 2$}-invariants by their finite-dimensional
  analogues.
\newblock {\em Geom. Funct. Anal.}, 4(4):455--481, 1994.

\bibitem{LRS}
Wenzhi Luo, Ze{\'e}v Rudnick, and Peter Sarnak.
\newblock On the generalized {R}amanujan conjecture for {${\rm GL}(n)$}.
\newblock In {\em Automorphic forms, automorphic representations, and
  arithmetic (Fort Worth, TX, 1996)}, volume~66 of {\em Proc. Sympos. Pure
  Math.}, pages 301--310. Amer. Math. Soc., Providence, RI, 1999.

\bibitem{MargulisSoifer}
G.~A. Margulis and G.~A. So{\u\i}fer.
\newblock Maximal subgroups of infinite index in finitely generated linear
  groups.
\newblock {\em J. Algebra}, 69(1):1--23, 1981.

\bibitem{Matsushima}
Yoz{\^o} Matsushima.
\newblock A formula for the {B}etti numbers of compact locally symmetric
  {R}iemannian manifolds.
\newblock {\em J. Differential Geometry}, 1:99--109, 1967.

\bibitem{MoeglinWaldspurger}
C.~M{\oe}glin and J.-L. Waldspurger.
\newblock Le spectre r\'esiduel de {${\rm GL}(n)$}.
\newblock {\em Ann. Sci. \'Ecole Norm. Sup. (4)}, 22(4):605--674, 1989.

\bibitem{Mostow}
G.~D. Mostow.
\newblock On a remarkable class of polyhedra in complex hyperbolic space.
\newblock {\em Pacific J. Math.}, 86(1):171--276, 1980.

\bibitem{Oda}
Takayuki Oda.
\newblock A note on the {A}lbanese variety of an arithmetic quotient of the
  complex hyperball.
\newblock {\em J. Fac. Sci. Univ. Tokyo Sect. IA Math.}, 28(3):481--486 (1982),
  1981.

\bibitem{OhsawaTakegoshi}
Takeo Ohsawa and Kensh{\=o} Takegoshi.
\newblock Hodge spectral sequence on pseudoconvex domains.
\newblock {\em Math. Z.}, 197(1):1--12, 1988.

\bibitem{Pedon}
Emmanuel Pedon.
\newblock Harmonic analysis for differential forms on complex hyperbolic
  spaces.
\newblock {\em J. Geom. Phys.}, 32(2):102--130, 1999.

\bibitem{PlatonovRapinchuk}
Vladimir Platonov and Andrei Rapinchuk.
\newblock {\em Algebraic groups and number theory}, volume 139 of {\em Pure and
  Applied Mathematics}.
\newblock Academic Press Inc., Boston, MA, 1994.
\newblock Translated from the 1991 Russian original by Rachel Rowen.

\bibitem{Prasad}
Gopal Prasad.
\newblock Semi-simple groups and arithmetic subgroups.
\newblock In {\em Proceedings of the International Congress of Mathematicians,
  Vol.\ I, II (Kyoto, 1990)}, pages 821--832, Tokyo, 1991. Math. Soc. Japan.

\bibitem{Rogawski}
Jonathan~D. Rogawski.
\newblock {\em Automorphic representations of unitary groups in three
  variables}, volume 123 of {\em Annals of Mathematics Studies}.
\newblock Princeton University Press, Princeton, NJ, 1990.

\bibitem{Sakai}
Takashi Sakai.
\newblock {\em Riemannian geometry}, volume 149 of {\em Translations of
  Mathematical Monographs}.
\newblock American Mathematical Society, Providence, RI, 1996.
\newblock Translated from the 1992 Japanese original by the author.

\bibitem{Sarnak}
P.~Sarnak.
\newblock The arithmetic and geometry of some hyperbolic three-manifolds.
\newblock {\em Acta Math.}, 151(3-4):253--295, 1983.

\bibitem{Schiffmann}
G{\'e}rard Schiffmann.
\newblock Int\'egrales d'entrelacement et fonctions de {W}hittaker.
\newblock {\em Bull. Soc. Math. France}, 99:3--72, 1971.

\bibitem{Schlichtkrull}
Henrik Schlichtkrull.
\newblock The {L}anglands parameters of {F}lensted-{J}ensen's discrete series
  for semisimple symmetric spaces.
\newblock {\em J. Funct. Anal.}, 50(2):133--150, 1983.

\bibitem{Selberg}
Atle Selberg.
\newblock On the estimation of {F}ourier coefficients of modular forms.
\newblock In {\em Proc. Sympos. Pure Math., Vol. VIII}, pages 1--15. Amer.
  Math. Soc., Providence, R.I., 1965.

\bibitem{Shahidi3}
Freydoon Shahidi.
\newblock Local coefficients and normalization of intertwining operators for
  {${\rm GL}(n)$}.
\newblock {\em Compositio Math.}, 48(3):271--295, 1983.

\bibitem{Shahidi}
Freydoon Shahidi.
\newblock Fourier transforms of intertwining operators and {P}lancherel
  measures for {${\rm GL}(n)$}.
\newblock {\em Amer. J. Math.}, 106(1):67--111, 1984.

\bibitem{Shahidi2}
Freydoon Shahidi.
\newblock Local coefficients as {A}rtin factors for real groups.
\newblock {\em Duke Math. J.}, 52(4):973--1007, 1985.

\bibitem{Shalika}
J.~A. Shalika.
\newblock The multiplicity one theorem for {${\rm GL}\sb{n}$}.
\newblock {\em Ann. of Math. (2)}, 100:171--193, 1974.

\bibitem{Shimura}
Goro Shimura.
\newblock Automorphic forms and the periods of abelian varieties.
\newblock {\em J. Math. Soc. Japan}, 31(3):561--592, 1979.

\bibitem{Speh}
Birgit Speh.
\newblock Unitary representations of {${\rm Gl}(n,\,{\bf R})$} with nontrivial
  {$(\mathfrak{g},\,K)$}-cohomology.
\newblock {\em Invent. Math.}, 71(3):443--465, 1983.

\bibitem{TongWang}
Y.~L. Tong and S.~P. Wang.
\newblock Harmonic forms dual to geodesic cycles in quotients of {${\rm
  SU}(p,\,1)$}.
\newblock {\em Math. Ann.}, 258(3):289--318, 1981/82.

\bibitem{Venky}
T.~N. Venkataramana.
\newblock Cohomology of compact locally symmetric spaces.
\newblock {\em Compositio Math.}, 125(2):221--253, 2001.

\bibitem{V}
M.-{F}. Vign\'eras.
\newblock {\em On the global correspondence between ${G}{L} (n)$ and division
  algebras}.
\newblock notes de l'{I}nstitute for {A}dvanced {S}tudy, {P}rinceton, 1984.

\bibitem{Vogan2}
D.~A. Vogan.
\newblock Isolated unitary representations.
\newblock to appear in the 2002 Park City summer school volume.

\bibitem{Vogan}
David~A. Vogan, Jr.
\newblock The unitary dual of {${\rm GL}(n)$} over an {A}rchimedean field.
\newblock {\em Invent. Math.}, 83(3):449--505, 1986.

\bibitem{VoganZuckerman}
David~A. Vogan, Jr. and Gregg~J. Zuckerman.
\newblock Unitary representations with nonzero cohomology.
\newblock {\em Compositio Math.}, 53(1):51--90, 1984.

\bibitem{Wallach}
Nolan~R. Wallach.
\newblock {\em Harmonic analysis on homogeneous spaces}.
\newblock Marcel Dekker Inc., New York, 1973.
\newblock Pure and Applied Mathematics, No. 19.

\bibitem{Yeganefar}
Nader Yeganefar.
\newblock Formes harmoniques {$L\sp 2$} sur les vari\'et\'es asymptotiquement
  hyperboliques complexes.
\newblock In {\em S\'eminaire de Th\'eorie Spectrale et G\'eom\'etrie. Vol. 21.
  Ann\'ee 2002--2003}, volume~21 of {\em S\'emin. Th\'eor. Spectr. G\'eom.},
  pages 55--59. Univ. Grenoble I, Saint, 2003.

\bibitem{Zimmer}
Robert~J. Zimmer.
\newblock {\em Ergodic theory and semisimple groups}, volume~81 of {\em
  Monographs in Mathematics}.
\newblock Birkh\"auser Verlag, Basel, 1984.

\end{thebibliography}

\bibliographystyle{plain}

\newpage

\thispagestyle{empty}

\newpage

\markboth{INDEX}{INDEX}

\printindex

\newpage

\thispagestyle{empty}

\noindent
N. \textsc{Bergeron}, L. \textsc{Clozel} \\
Universit\'e Paris-Sud \\
Unit\'e Mixte de Recherche 8628 du CNRS, \\
Laboratoire de Math\'ematiques, B\^at. 425, \\
91405 Orsay Cedex, France \\
{\it adresse electronique :} \texttt{Nicolas.Bergeron@math.u-psud.fr}

\end{document}